\documentclass[11pt,reqno]{amsart}

\pdfoutput=1

\usepackage{CL1_Macros}

\begin{document}

\title[Asymptotic stability of the sine-Gordon kink]{Asymptotic stability of the sine-Gordon kink}

\begin{abstract}
    We establish the full asymptotic stability of the sine-Gordon kink outside symmetry under small perturbations in weighted Sobolev norms. 
    Our proof consists of a space-time resonances approach based on the distorted Fourier transform to capture modified scattering effects combined with modulation techniques to take into account the invariance under Lorentz transformations and under spatial translations.
    A major challenge is the slow local decay of the radiation term caused by the threshold resonances of the non-selfadjoint linearized matrix operator around the moving kink. Our analysis crucially relies on two remarkable null structures in the quadratic nonlinearities of the evolution equation for the radiation term and of the modulation equations.

    The entire framework of our proof, including the systematic development of the distorted Fourier theory, is not specific to the sine-Gordon model and extends to many other asymptotic stability problems for moving solitons in relativistic scalar field theories on the line.
\end{abstract}

\author[G. Chen]{Gong Chen}
\address{School of Mathematics, Georgia Institute of Technology, Atlanta, GA 30332, USA}
\email{gc@math.gatech.edu}

\author[J. L\"uhrmann]{Jonas L\"uhrmann}
\address{Department of Mathematics \\ Texas A\&M University \\ College Station, TX 77843, USA}
\email{luhrmann@tamu.edu}

\thanks{
The first author was partially supported by NSF grant DMS-2350301 and by Simons foundation MP-TSM-0000225.
The second author was partially supported by NSF grant DMS-1954707 and by NSF CAREER grant DMS-2235233.
}

\maketitle 

\tableofcontents

\section{Introduction}

\subsection{The sine-Gordon model}

The sine-Gordon model is a classical field theory on the line for scalar fields $\phi \colon \bbR \times \bbR \to \bbR$ governed by the Lagrangian action functional
\begin{equation}
    \iint_{\bbR^{1+1}} \biggl( \frac12 (\pt \phi)^2 - \frac12 (\px \phi)^2 - \bigl( 1 - \cos(\phi) \bigr) \biggr) \, \ud x \, \ud t.
\end{equation}
The associated Euler-Lagrange equation is called the sine-Gordon equation 
\begin{equation} \label{equ:intro_sG} 
    \pt^2 \phi - \px^2 \phi + \sin(\phi) = 0, \quad (t,x) \in \bbR \times \bbR.
\end{equation}
Dating back to the 1860s the sine-Gordon equation first arose in the study of surfaces with constant negative curvature \cite{Bour1862}, and it has since found a broad range of applications in physics. 
We refer to \cite{DauxPey10, Lamb80, MantSut04, SG_Series} for more background.
It is also a prime example of a completely integrable nonlinear dispersive equation on the line in the sense that it admits a formulation as a Lax pair equation \cite{AKNS,Tal, Kaup, ZTF} and it admits B\"acklund transformations \cite{Lamb80,MSbook,APM,KY23}.

While the sine-Gordon equation admits infinitely many conservation laws, a distinct role is played by the conservation of energy
\begin{equation}
    E = \int_\bbR \biggl( \frac12 (\pt \phi)^2 + \frac12 (\px \phi)^2 + \bigl( 1 - \cos(\phi) \bigr) \biggr) \, \ud x
\end{equation}
and the conservation of momentum 
\begin{equation}
    P = \int_\bbR (\pt \phi) (\px \phi) \, \ud x.
\end{equation}
The Cauchy problem for \eqref{equ:intro_sG} is globally well-posed in the energy space $H^1_{\sin}(\bbR) \times L^2(\bbR)$, see \cite{BuckinghamMiller08, LaireGravejat18}, where 
\begin{equation}
    H^1_{\sin}(\bbR) := \Bigl\{ \phi \in \dot{H}^1(\bbR) \, \Big| \, \sin\Bigl( \frac{\phi}{2} \Bigr) \in L^2(\bbR) \Bigr\}.
\end{equation}
Moreover, the sine-Gordon equation~\eqref{equ:intro_sG} is invariant under space-time translations
\begin{equation*}
    \phi(t,x) \mapsto \phi(t-t_0, x-x_0), \quad (t_0,x_0) \in \bbR \times \bbR,
\end{equation*}
and under Lorentz transformations
\begin{equation*}
    \phi(t,x) \mapsto \phi\bigl( \gamma(t-\ell x), \gamma(x - \ell t) \bigr), \quad \ell \in (-1,1), \quad \gamma := \frac{1}{\sqrt{1-\ell^2}} \in [1,\infty).
\end{equation*}

Among the many intriguing features of the sine-Gordon equation are its soliton solutions called kinks and breathers. A static kink solution to \eqref{equ:intro_sG} is given by
\begin{equation} \label{eq:staticK}
    K(x) = 4 \arctan(e^x).
\end{equation}
Breathers are time-periodic, spatially localized solutions to~\eqref{equ:intro_sG} such as
\begin{equation*}
 B_\beta(t,x) = 4 \arctan \biggl( \frac{\beta}{\alpha} \frac{\sin(\alpha t)}{\cosh(\beta x)} \biggr), \quad \beta \in (-1,1) \backslash \{0\}, \quad \alpha := \sqrt{1-\beta^2}.
\end{equation*}
The invariance under Lorentz transformations then gives rise to families of moving kink solutions and of moving breather solutions to \eqref{equ:intro_sG}.

Kinks and breathers are building blocks of the long-time dynamics for the sine-Gordon equation in the sense that typical solutions to \eqref{equ:intro_sG} asymptotically decompose into finitely many moving kinks, anti-kinks, and breathers, along with a radiation term that behaves like a modified linear Klein-Gordon wave.
Such a soliton resolution result was rigorously proven in \cite{CLL20} using complete integrability techniques.

\subsection{Motivation: Asymptotic stability of moving kinks}

The main motivation for this paper is to develop a robust perturbative framework to study the full asymptotic stability of moving kink solutions that arise in many relativistic scalar field theories on the line.
Before stating our main result specifically for the sine-Gordon model, we provide a brief discussion of our work in this broader context.
The sine-Gordon model is a prime example of a scalar field theory on the line whose Lagrangian features a double-well self-interaction potential denoted by $W \colon \bbR \to [0,\infty)$,
\begin{equation} \label{equ:intro_lagrangian_general}
    \iint_{\bbR^{1+1}} \biggl( \frac12 (\pt \phi)^2 - \frac12 (\px \phi)^2 - W(\phi) \biggr) \, \ud x \, \ud t.
\end{equation}
Specifically, $W(\phi)$ is a sufficiently regular non-negative function that has (at least) two consecutive global minima denoted by $\phi_- < \phi_+$ with $W(\phi_\pm) = W'(\phi_\pm) = 0$ and $W''(\phi_\pm) > 0$. Often, the minima $\phi_-$ and $\phi_+$ are referred to as the vacuum states of the model.
Apart from the sine-Gordon model for which $W(\phi) = 1 -\cos(\phi)$, another prime example is the $\phi^4$ model with $W(\phi) = \frac14 (1-\phi^2)^2$. Further well-known examples include the more general $P(\phi)_2$ theories \cite{Lohe79} and the double sine-Gordon theories \cite{CampbellPeyrardSodano86}.

The Euler-Lagrange equations associated with ~\eqref{equ:intro_lagrangian_general} are the wave-type equations
\begin{equation} \label{equ:intro_wave_type_field_equ}
 \begin{aligned}
  (\pt^2 - \px^2) \phi + W'(\phi) = 0, \quad (t,x) \in \bbR \times \bbR,
 \end{aligned}
\end{equation}
which admit static solutions $H(x)$ called kinks that connect the two distinct vacuum states
\begin{equation}
 \left\{ \begin{aligned}
  -\px^2 H + W'(H) &= 0, \quad x \in \bbR, \\
  \lim_{x\to\pm\infty} H(x) &= \phi_\pm.
 \end{aligned} \right.
\end{equation}
The invariance of \eqref{equ:intro_wave_type_field_equ} under Lorentz transformations and under spatial translations leads to the following two-parameter family of moving kink solutions
\begin{equation} \label{equ:intro_general_moving_kink_family}
    H\bigl(\gamma_0(x-\ell_0 t - x_0)\bigr), \quad \ell_0 \in (-1,1), \quad x_0 \in \bbR, \quad \gamma_0 := \frac{1}{\sqrt{1-\ell_0^2}} \in [1,\infty).
\end{equation}

The asymptotic stability of kinks is a fundamental question concerning the long-time dynamics of the wave-type equations~\eqref{equ:intro_wave_type_field_equ}.
The interplay of weak dispersion in one space dimension, low-power nonlinearities, and intriguing spectral features of linearized operators around kinks — such as threshold resonances or internal modes — creates a rich and challenging class of asymptotic stability problems. Despite substantial progress in recent years, as detailed in Subsection 1.4 below, significant open questions remain.

It has become customary to distinguish two notions of asymptotic stability for solitons on the line that complement each other.
One typically seeks to derive a decomposition of the perturbed solution into a possibly modulated kink and a radiation term that asymptotically behaves like a (possibly modified) linear Klein-Gordon wave.
{\it Full asymptotic stability} refers to proving sharp decay estimates and asymptotics for the radiation term at the expense of strong assumptions such as requiring the initial perturbations to be small in weighted Sobolev norms. 
{\it Local asymptotic stability} involves establishing decay of the radiation term locally in space (without explicit decay rates) but under much weaker assumptions on the initial data, typically for small finite energy perturbations. 

Focusing now on full asympotitic stability,
in recent years there has been significant progress on the full asymptotic stability of kinks under symmetry assumptions (odd perturbations) that prevent the kink from starting to move.
Outside symmetry, the only full asymptotic stability result for a general class of scalar field models was obtained in \cite{KK11_2} under rather tailored assumptions about the double-well potential and about the linearized operator around the kink. These assumptions in particular bypass having to deal with modified scattering behavior of the radiation term. In the context of perturbations of kinks, modified scattering of the radiation term typically refers to logarithmic phase corrections with respect to the free flow caused by constant coefficient cubic nonlinearities in the evolution equation for the radiation term. However, many open questions are related to potential additional amplitude corrections such as slow-downs of the decay rate with respect to the free flow. These may result from the impact of threshold resonances or internal modes of the linearized operator on the dynamics, primarily affecting the radiation term through the quadratic nonlinearities.

In this work we develop a robust perturbative framework to establish the full asymptotic stability of the family of moving kink solutions to the sine-Gordon equation under arbitrary perturbations that are small with respect to a weighted Sobolev norm. This completes previous work for the specific case of odd perturbations \cite{LS1}.
Our proof is based on a space-time resonances approach using the distorted Fourier transform associated with the linearized operator around the moving kink to capture the modified scattering behavior of the radiation term and on modulation techniques to take into account the invariance under Lorentz transformations and under spatial translations.
Relying on complete integrability techniques such a full asymptotic stability result outside symmetry for the family of moving sine-Gordon kinks was previously obtained in \cite{CLL20, KY23}, in fact with slightly more refined information on the asymptotic behavior of the perturbed solution. 

\medskip

\noindent {\it Outlook.}
Beyond the intrinsic interest in a perturbative proof of the full asymptotic stability of moving kink solutions to the sine-Gordon equation, our main motivation stems from the potential applications of the framework and the techniques introduced in this paper to many other full asymptotic stability problems for moving solitons in relativistic scalar field theories on the line.
For example, combined with insights from \cite{GP20}, we expect to be able to prove a full asymptotic stability result for moving kinks for a fairly general class of double-well self-interaction potentials under the key assumption that the linearized operator around the kink does not exhibit any threshold resonances nor any internal modes.

In cases where the linearized operator does exhibit threshold resonances or internal modes, understanding the full asymptotic stability of kinks remains a major open challenge, even in the special case of symmetric perturbations. In this direction, the sine-Gordon model serves as a valuable intermediate step: while the linearized operator around the moving kink has threshold resonances, their most severe effects on the dynamics are mitigated by favorable null structures in the quadratic nonlinearities.
In the absence of such favorable structures, at least a proof of long-time decay estimates for the radiation term outside symmetry (up to time scales that are polynomial or exponential in terms of inverse powers of the size of the initial perturbations) is within reach of the techniques introduced in this paper.

For ground states of focusing Klein-Gordon equations on the line, similar remarks hold regarding the problem of their full asymptotic stability outside symmetry on the corresponding center-stable manifolds.

\subsection{Main result} \label{subsec:main_result}

In this work we study the sine-Gordon equation \eqref{equ:intro_sG} in terms of the equivalent first-order system 
\begin{equation} \label{equ:intro_sG_1st_order}
    \pt \bmphi = \begin{bmatrix} 0 & 1 \\ \px^2 & 0 \end{bmatrix} \bmphi + \begin{bmatrix} 0 \\ - W'(\phi) \end{bmatrix}, \quad \bm{\phi} = (\phi, \pt \phi) = (\phi_1, \phi_2),
\end{equation}
where the self-interaction potential is given by
\begin{equation}
    W(\phi) := 1 - \cos(\phi).
\end{equation}
The static kink~\eqref{eq:staticK} and the invariance of the sine-Gordon equation under Lorentz transformations as well as under spatial translations give rise to the following two-parameter family of (vectorial) moving kink solutions to \eqref{equ:intro_sG_1st_order},
\begin{equation} \label{equ:intro_vectorial_moving_kink_family}
    \begin{aligned}
        \begin{bmatrix}
            K\bigl(\gamma_0(x-\ell_0 t -x_0)\bigr) \\ -\gamma_0 \ell_0 K'\bigl(\gamma_0(x-\ell_0 t - x_0)\bigr)
        \end{bmatrix},
        \quad \ell_0 \in (-1,1), \quad x_0 \in \bbR, \quad \gamma_0 := \frac{1}{\sqrt{1-\ell_0^2}} \in [1,\infty).
    \end{aligned}
\end{equation}
Our main result is the full asymptotic stability of the family of moving sine-Gordon kinks \eqref{equ:intro_vectorial_moving_kink_family} under small perturbations in weighted Sobolev norms.
\begin{theorem} \label{thm:main}
 For any $\ell_0 \in (-1,1)$ there exist constants $0 < \varepsilon_0 \ll 1$, $0 < \delta \ll 1$, and $C \geq 1$ such that for any $x_0 \in \bbR$, and any $\bmu_0 = (u_{0,1}, u_{0,2}) \in H^3_x(\bbR) \times H^2_x(\bbR)$ with 
 \begin{equation} \label{equ:statement_theorem_smallness_data}
     \varepsilon := \bigl\| \jx (u_{0,1}, u_{0,2}) \bigr\|_{H^3_x \times H^2_x} \leq \varepsilon_0,
 \end{equation}
 the solution to the sine-Gordon equation \eqref{equ:intro_sG_1st_order} with initial data 
 \begin{equation} \label{equ:statement_theorem_initial_data}
     \bmphi(0,x) = \bmK_{\ell_0,x_0}(x) + \bmu_0(x-x_0), \quad \bmK_{\ell_0, x_0}(x) := \begin{bmatrix}
        K\bigl(\gamma_0 (x-x_0)\bigr) \\ -\gamma_0 \ell_0 K'\bigl(\gamma_0(x-x_0)\bigr)
    \end{bmatrix},
 \end{equation}
 exists globally in time and there exist continuously differentiable paths $\ell \colon [0,\infty) \to (-1,1)$ and $q \colon [0,\infty) \to \bbR$ with 
 \begin{equation*}
     |\ell(0) - \ell_0| + |q(0)-x_0| \leq C \varepsilon
 \end{equation*}
 such that the following holds:
 \begin{itemize}[leftmargin=1.8em]
        \item[(1)] Decomposition of the solution into a modulated kink and a radiation term:
        \begin{equation}
            \bmphi(t,x) = \bmK_{\ell(t),q(t)}(x) + \bmu\bigl(t,x-q(t)\bigr), \quad t \geq 0,
        \end{equation}
        where
        \begin{equation} \label{equ:statement_theorem_modulated_kink}
            \bmK_{\ell(t),q(t)}(x) := \begin{bmatrix}
                        K\bigl( \gamma(t) (x-q(t)) \bigr) \\
                        - \gamma(t) \ell(t) K'\bigl( \gamma(t) (x-q(t)) \bigr)
                     \end{bmatrix}, 
                     \quad \bmu = \begin{bmatrix} u_1 \\ u_2 \end{bmatrix},
                     \quad \gamma(t) := \frac{1}{\sqrt{1-\ell(t)^2}}.
        \end{equation}        

        \item[(2)] Dispersive decay of the radiation term:
        \begin{equation} \label{equ:statement_theorem_dispersive_decay_estimate}
            \|u_1(t)\|_{L^\infty_x(\bbR)} + \|u_2(t)\|_{L^\infty_x(\bbR)} + \|\px u_1(t)\|_{L^\infty_x(\bbR)} \leq C \varepsilon \jt^{-\frac12}, \quad t \geq 0.
        \end{equation}

        \item[(3)] Asymptotic behavior of the modulation parameters: 
        \begin{equation} \label{equ:statement_theorem_modulation_parameters_asymptotics}
            \bigl|\ell(t) - \ell_\infty\bigr| + \bigl|\dot{q}(t) - \ell_\infty\bigr| \leq C \varepsilon \jt^{-1+\delta}, \quad t \geq 0,
        \end{equation}
         for some $\ell_\infty \in (-1,1)$ with $|\ell_\infty - \ell_0| \leq C \varepsilon$.

        \item[(4)] The asymptotics of the radiation term exhibit logarithmic phase corrections with respect to the free flow, specifically there exists $g_{\infty} \in L^\infty_\xi(\bbR)$ such that 
 \begin{equation} \label{equ:statement_theorem_radiation_asymptotics}
     \bigl\| \bmu\bigl(t, \cdot - q(t)\bigr) - \bmu_\infty\bigl(t, \cdot \bigr) \bigr\|_{L^\infty_x(\bbR)} \leq C \varepsilon t^{-\frac12-2\delta}, \quad t \geq 1,
 \end{equation}
 where
 \begin{equation} \label{equ:statement_theorem_radiation_asymptotic_state}
     \begin{aligned}
         \bmu_\infty(t,x) &:= \frac{1}{t^{\frac12}} \Im \biggl( e^{i\frac{\pi}{4}} e^{i\rho(t,x)} e^{i \Gamma_\infty(\xi_\ast) \log(t)} e^{i\xi_\ast\theta_\infty(t)} \\
         &\qquad \qquad \qquad \times g_\infty(\xi_\ast) \bm{m}_\infty(t, x, \xi_\ast) \one_{(-t,t)}\bigl( x - \theta_\infty(t) - q(0) \bigr) \biggr),
     \end{aligned}
 \end{equation} 
 and where for $-t < x - \theta_\infty(t) - q(0) < t$,
 \begin{equation} \label{equ:statement_theorem_asymptotics_quantities_definitions}
 \begin{aligned}
     \rho(t,x) &:= \sqrt{t^2 - (x - \theta_\infty(t) - q(0))^2}, 
     \quad \quad \theta_\infty(t) := \int_0^t \bigl( \dot{q}(s) - \ell_\infty \bigr) \, \ud s, \\
     \Gamma_\infty(\xi_\ast) &:= \frac{1}{16} \jap{\xi_\ast}^{-3} |g_\infty(\xi_\ast)|^2, 
     \quad \quad \xi_\ast := - \frac{x - \theta_\infty(t) - q(0)}{\rho(t,x)}, \\
    \bm{m}_\infty(t,x, \xi_\ast) &:= \begin{bmatrix} m_{\infty,1}(t,x,\xi_\ast) \\ m_{\infty,2}(t,x,\xi_\ast) \end{bmatrix}, \\ 
     m_{\infty,1}(t,x,\xi_\ast) &:= \jap{\xi_\ast}^{-1} \frac{\gamma_\infty(\xi_\ast + \ell_\infty \jap{\xi_\ast}) + i \tanh\bigl(\gamma_\infty (x - \ell_\infty t - \theta_\infty(t) - q(0)) \bigr)}{|\gamma_\infty(\xi_\ast + \ell_\infty \jap{\xi_\ast})| - i}, \\
     m_{\infty,2}(t,x,\xi_\ast) &:= i \jap{\xi_\ast} m_{\infty,1}(t,x,\xi_\ast) - \ell_\infty \px m_{\infty,1}(t,x,\xi_\ast).
 \end{aligned}
 \end{equation}
 \end{itemize}
\end{theorem}

\begin{remark}
    Heuristically, bounds on the function $\theta_\infty(t)$ defined in \eqref{equ:statement_theorem_asymptotics_quantities_definitions} quantify how sharply the center of the moving kink is known, see the discussion at the end of Subsection~\ref{subsec:overview_modulation}. The decay estimates \eqref{equ:statement_theorem_modulation_parameters_asymptotics} imply that $|\theta_\infty(t)| \lesssim \varepsilon \jt^\delta$, which corresponds to a small uncertainty in the location of the center of the kink that has to be propagated throughout the entire proof of Theorem~\ref{thm:main}.
\end{remark}

\begin{remark}
We compare the statement of Theorem~\ref{thm:main} with the full asymptotic stability results for moving sine-Gordon kinks established  in \cite{CLL20, KY23} using complete integrability techniques.
    
In \cite[Corollary 1.6]{CLL20} the asymptotic stability of the family of moving sine-Gordon kinks is obtained as a by-product of proving soliton resolution for the sine-Gordon equation for generic initial data via the inverse scattering transform and the study of Riemann-Hilbert problems.
In comparison to Theorem~\ref{thm:main}, the assumptions on the initial data are slightly weaker in \cite[Corollary 1.6]{CLL20} and sharp control of the center of the moving kink is established there.
However, the formulas for the asymptotics of the radiation term in \cite{CLL20} are less transparent\footnote{The formulas for the asymptotics away from the centers of kinks, anti-kinks, and breathers are more explicit, see Proposition~8.5 in \cite{CLL20}.} than the formula \eqref{equ:statement_theorem_radiation_asymptotic_state} obtained in the statement of Theorem~\ref{thm:main}.

In \cite{KY23} the B\"acklund transformation is used to map a neighborhood of a moving sine-Gordon kink to a neighborhood of the zero solution. In this manner the proof of the asymptotic stability of the family of moving sine-Gordon kinks in \cite{KY23} is transformed to proving decay and asymptotics of small solutions to the sine-Gordon equation, which is achieved using the method of testing by wave packets \cite{IT15}. 
Under slightly weaker assumptions on the initial data than in this work, \cite[Theorem~1(b)]{KY23} establishes decay and asymptotics for the radiation term with limited control on the center of the moving kink. Imposing stronger weighted Sobolev norm assumptions on the initial data, \cite[Theorem~1(c)]{KY23} additionally gains sharp control of the center of the moving kink.
\end{remark}

\subsection{References} \label{subsec:references}

The study of the asymptotic stability of solitons is a rich and vast subject that we do not attempt to comprehensively review here.
Instead, we focus on the significant advancements made in recent years regarding the asymptotic stability of kinks in scalar field theories on the line and of solitary waves in nonlinear Schr\"odinger models in one space dimension. 
In this subsection we collect the corresponding references and we briefly describe how this work relates to closely connected articles.
For additional references to results in higher space dimensions and for other models, we refer to the reviews \cite{Martel_ICM, CuccMaeda20_Survey, Germain24_Review, KMM17_Review, Tao09}.

The orbital stability of kinks in relativistic scalar field theories on the line is classical \cite{HPW82}.
Local asymptotic stability results for kinks include \cite{KMM17, KMM17_short, KMMV20, KM22, CuccMaeda_kink_23, AMP20}, 
and full asymptotic stability results for kinks have been obtained in \cite{KK11_1, KK11_2, GermPusZhang22, CLL20, KY23, LS1, DelMas20}.
Closely related are the following local or full (co-dimension one) asymptotic stability results for ground states of focusing nonlinear Klein-Gordon equations on the line \cite{KNS12, KMM19, KairzhanPusateri22, LL1, LS2, PalaciosPusateri24, CuccMaedaMurgScrob23}.
We also refer to local asymptotic stability results for solitary waves in nonlinear Schr\"odinger models on the line \cite{CuccMaeda19, CuccMaeda22, CuccMaeda2404, CuccMaeda2405, Martel22, Martel24, Rialland23, Rialland24}
and to full asymptotic stability results for such solitary waves \cite{BusPerel92, BusSul03, KS06, Mizumachi08, Chen21, MasMurphSeg20, CG23, LL2}.

The study of the full asymptotic of kinks involves proving decay and asymptotics for small solutions to nonlinear Klein-Gordon equations on the line with variable coefficient low power nonlinearities and with potentials 
\cite{LLS1, LLS2, LLSS, GP20, Sterb16, LS15}. See \cite{Del01, LS05_1, LS05_2, HN08, HN10, HN12, Del16_KG, Stingo18, CL18} for earlier results on the modified scattering of small solutions to 1D Klein-Gordon equations with constant coefficient quadratic or cubic nonlinearities, but without potentials.
Related modified scattering results for small solutions to 1D nonlinear Schr\"odinger equations without potentials include \cite{HN98, LS06, KatPus11, IT15, GermPusRou18, CP21, CP22, Del16, Naum16, Naum18, MasMurphSeg19, NaumWed22}.
We also refer to the expository articles \cite{IfrimTataru_wave_packets_expository, Murphy_expository} on modified scattering for cubic nonlinear dispersive equations on the line.

\medskip 

\noindent {\it Closely related articles.}
A central component of our framework for proving Theorem~\ref{thm:main} is the development of a space-time resonances approach to capture the modified scattering behavior of the radiation term. This approach relies on the distorted Fourier transform associated with the non-selfadjoint matrix operator obtained from linearizing around a moving kink. While we build the distorted Fourier theory for this operator from the ground up, our functional framework and the concept of nonlinear spectral distributions draw from the setting of scalar linearized operators around static kinks such as in \cite{GP20, LS1} and from the corresponding Schr\"odinger settings \cite{GermPusRou18, CP21, CP22}.  
We of course also rely on insights from the perturbative proof of the asymptotic stability of the sine-Gordon kink under odd perturbations in \cite{LS1}, in particular we generalize a variable coefficient quadratic normal form from \cite{LS1} to our setting.

Another key element of our framework is the use of modulation techniques to account for the invariances of the sine-Gordon equation under Lorentz transformations and under spatial translations, as in the earlier local \cite{KMMV20} and full \cite{KK11_2} asymptotic stability results for kinks.

We point out that the only prior work \cite{KK11_2} on full asymptotic stability of kinks outside symmetry is under rather tailored assumptions on the double-well self-interaction potential. These ensure that the nonlinearities in the evolution equation for the radiation term are all spatially localized (at least up to tenth order). Another crucial assumption in \cite{KK11_2} is that the linearized operator does not exhibit threshold resonances. Consequently, the radiation term enjoys improved local decay and is not subject to modified scattering effects in the setting of \cite{KK11_2}. For this reason, in \cite{KK11_2} the modulation setup can be combined with a relatively simple functional framework based on pointwise-in-time local decay estimates to close the argument.

Instead, in this work we combine for the first time in the setting of moving kinks modulation techniques with a much more involved space-time resonances method based on the distorted Fourier transform of the associated matrix linearized operator.
Such an approach has previously been taken in the context of proving full asymptotic stability results for solitary wave solutions to nonlinear Schr\"odinger equations on the line in \cite{Chen21, CG23} in the generic setting (no threshold resonances, no internal modes) and in \cite{LL2} for the focusing cubic Schr\"odinger equation (with threshold resonances, but still under symmetry assumptions). 
We point out that the distorted Fourier theory applied in \cite{Chen21, CG23,LL2} had largely been developed earlier in \cite{BusPerel92,DT, KS06}.

\subsection{Organization of the paper}

In Section~\ref{sec:outline_argument} we provide an overview of the proof of Theorem~\ref{thm:main}.
In Section~\ref{sec:preliminaries} we introduce basic notational conventions and we recall the standard distorted Fourier theory for scalar linearized operators around kinks.
Then we develop from the ground up the spectral and distorted Fourier theory for the matrix operator obtained from linearizing around a moving kink in Sections~\ref{sec:spectral_theory} and~\ref{sec:distorted_fourier_theory}.
Linear decay estimates for the evolution are established in Section~\ref{sec:linear_decay_estimates}.
In preparation for the nonlinear analysis we compute quadratic and cubic spectral distributions in Section~\ref{sec:nonlinear_spectral_distributions}.
The proof of Theorem~\ref{thm:main} begins in earnest in Section~\ref{sec:setting_up}, where we lay out the modulational setup and where we prepare the evolution equation of the radiation term for the nonlinear analysis.
In Section~\ref{sec:bootstrap_setup_proof_thm} we formulate the entire bootstrap argument in terms of two central bootstrap propositions and, based on their conclusions, we provide a succinct proof of Theorem~\ref{thm:main}.
The remainder of the paper is then devoted to the proofs of the two bootstrap propositions.
In Section~\ref{sec:modulation_control} we obtain control of the modulation parameters. In Section~\ref{sec:energy_estimates} we establish all energy estimates for the distorted Fourier transform of the profile of the radiation term, and in Section~\ref{sec:pointwise_profile_bounds} we deduce the corresponding pointwise estimates for the profile.

\subsection{Acknowledgements}

The authors are grateful to Wilhelm Schlag for helpful discussions and to Herbert Koch, Yongming Li, Wilhelm Schlag, and Sohrab Shahshahani for valuable comments on the manuscript. 
Part of this work was conducted in the summer of 2024 while the authors participated in the thematic program ``Nonlinear Waves and Relativity'' at the Erwin-Schr\"odinger International Institute for Mathematics and Physics whose hospitality and support is gratefully acknowledged.
The second author also warmly thanks the Mathematical Institute of the University of Bonn for its hospitality and support in the fall of 2024.

\section{Overview of the Proof.} \label{sec:outline_argument}

In this section we outline the main ideas and concepts that enter the proof of Theorem~\ref{thm:main}.
Overall we develop a robust framework for the study of moving solitons in relativistic scalar field theories on the line that brings together a set of techniques related to {\it modified scattering}, {\it distorted Fourier theory}, and {\it modulation}.

\subsection{Linearized operator}

\begin{figure}[ht] \label{fig:spectrum_linearized_operator}
\centering
\begin{tikzpicture}[domain=-1.74:1.74,samples=50, scale=1.3]

\coordinate[label=left:{\footnotesize $i\gamma^{-1}$}] (pm) at ($(90:1)$);
\coordinate[label=left:{\footnotesize $-i\gamma^{-1}$}] (pp) at ($(270:1)$);
\coordinate (px) at ($(0:0)$);         

\draw[->] ($(180:2)$) --++ ($(0:4)$) node[right] {$\Re$};
\draw[->] ($(270:2)$) --++ ($(90:4)$)node[above] {$\Im$};

\fill[blue] (px) circle (2pt);         
\draw[line width=5pt,color=orange,draw opacity=0.45] (pp) --++ ($(270:1)$);
\draw[line width=5pt,color=orange,draw opacity=0.45] (pm) --++ ($(90:1)$);     
\fill[red] (pm) circle (2pt);
\fill[red] (pp) circle (2pt);


\end{tikzpicture}
\caption{Spectral features of the operator $\bfL_\ell$ defined in \eqref{equ:overview_definition_bfLell}: The orange bands indicate the essential spectrum, the blue dot corresponds to the zero eigenvalue, and the red dots represent the threshold resonances.}
\end{figure}
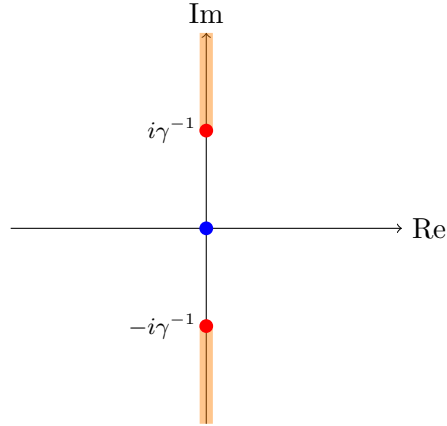

We begin with a discussion of the spectral features of the linearized operator around the moving sine-Gordon kink.
Starting from the first-order formulation \eqref{equ:intro_sG_1st_order} of the sine-Gordon equation, we linearize around a moving kink using the decomposition 
\begin{equation}
    \bmphi(t,x) = \begin{bmatrix} K\bigl(\gamma(x-\ell t)\bigr) \\ -\gamma \ell K'\bigl(\gamma(x-\ell t)\bigr) \end{bmatrix} + \bmu(t,x-\ell t), \quad \ell \in (-1,1), \quad \gamma := \frac{1}{\sqrt{1-\ell^2}}.
\end{equation}
Upon passing to the moving frame coordinate $y := x - \ell t$, we arrive at the following non-selfadjoint matrix operator
\begin{equation} \label{equ:overview_definition_bfLell}
    \bfL_\ell := \begin{bmatrix} \ell \py & 1 \\ - L_\ell & \ell \py \end{bmatrix}, \quad L_\ell := -\py^2 - 2\sech^2(\gamma y) + 1.
\end{equation}
The operator $\bfL_\ell$ is defined on $L^2(\bbR) \times L^2(\bbR)$ and it is closed on its domain $H^2(\bbR) \times H^1(\bbR)$. Its essential spectrum is $(-i\infty,-i\gamma^{-1}] \cup [i\gamma^{-1}, i\infty)$ with no embedded eigenvalues. The only eigenvalue is $0$ with a two-dimensional generalized kernel. The generalized kernel elements are generated by the invariances under spatial translations and under Lorentz transformations. They are explicitly given by
\begin{equation} \label{equ:overview_definition_generalized_eigenfunctions}
    \begin{aligned}
        \bmY_{0,\ell}(y) := \begin{bmatrix} -\gamma K'(\gamma y) \\ \gamma^2 \ell K''(\gamma y) \end{bmatrix}, \quad
        \bmY_{1,\ell}(y) := \begin{bmatrix} \gamma^3 \ell y K'(\gamma y) \\ - \gamma^3 K'(\gamma y) - \gamma^4 \ell^2 y K''(\gamma y) \end{bmatrix},
    \end{aligned}
\end{equation}
and satisfy
\begin{equation*}
    \bfL_\ell \bmY_{0,\ell} = 0, \quad \bfL_\ell \bmY_{1,\ell} = \bmY_{0,\ell}.
\end{equation*}
Additionally, the operator $\bfL_\ell$ exhibits threshold resonances on the two edges of its essential spectrum. Indeed, the functions
\begin{equation}
    \bmPhi_{\pm, \ell}(y) := \begin{bmatrix} \tanh(\gamma y) \\ \pm i \gamma \tanh(\gamma y) - \gamma^2 \ell \sech^2(\gamma y) \end{bmatrix} e^{\mp i\gamma \ell y}
\end{equation}
belong to $L^\infty(\bbR) \times L^\infty(\bbR) \, \backslash \, L^2(\bbR) \times L^2(\bbR)$ and satisfy
\begin{equation*}
    \bfL_\ell \bmPhi_{\pm,\ell} = \pm i \gamma^{-1} \bmPhi_{\pm,\ell}.
\end{equation*}  
All of the spectral features of the operator $\bfL_\ell$ are illustrated in Figure~\ref{fig:spectrum_linearized_operator}.

\subsection{Distorted Fourier theory and linear decay estimates} \label{subsec:overview_distorted_FT}

Since we adopt a space-time resonances approach to establish decay for the radiation term, the distorted Fourier transform associated with the linearized operator~\eqref{equ:overview_definition_bfLell} plays a central role for the entire analysis in this work.

In Sections~\ref{sec:spectral_theory}~and~\ref{sec:distorted_fourier_theory} we systematically develop the distorted Fourier theory for \eqref{equ:overview_definition_bfLell}.
We stress that all of our arguments apply to general matrix operators of the form~\eqref{equ:overview_definition_bfLell} that arise from linearizing around a moving soliton in relativistic scalar field theories on the line.
Since the potential in \eqref{equ:overview_definition_bfLell} belongs to the family of P\"oschl-Teller potentials \cite{PoschlTeller}, the expressions for the distorted Fourier basis elements turn out to be explicit. However, we do not make use of this specific feature of the sine-Gordon model in the systematic development of the distorted Fourier theory\footnote{The only minor place where we exploit this feature is in Subsection~\ref{subsec:mapping_properties} to give straightforward proofs of simple mapping properties of the distorted Fourier transform. These proofs can be replaced easily by general arguments.} for \eqref{equ:overview_definition_bfLell}.

The starting point for the construction of the distorted Fourier transform associated with~\eqref{equ:overview_definition_bfLell} is to express the evolution $e^{t\bfL_\ell} \Pe$ in terms of the jump of the resolvent across the essential spectrum
    \begin{equation} \label{equ:overview_evolution_jump_across_resolvent}
        e^{t\bfL_\ell} \Pe = \frac{1}{2\pi} \int_{\Sigma_\ell} e^{i\lambda t} \bigl[ \bfR(\lambda+i0)-\bfR (\lambda-i0) \bigr] \, \ud \lambda. 
    \end{equation}
Here, $\Pe$ denotes the projection to the essential spectrum relative to $\bfL_\ulell$, $\Sigma_\ell := \bigl(-\infty,-\gamma^{-1}\bigr] \cup \bigl[\gamma^{-1},\infty\bigr)$, and $\bfR(\mu) := \bigl(\bfL_\ell - i \mu \bbI_{2 \times 2}\bigr)^{-1}$ for $\mu \in \bbC \backslash \Sigma_\ell$ is the resolvent. 
While \eqref{equ:overview_evolution_jump_across_resolvent} is just Stone's formula in the selfadjoint setting, the derivation of \eqref{equ:overview_evolution_jump_across_resolvent} for the non-selfadjoint operator $\bfL_\ell$ is based on a suitable contour argument, generally proceeding as in \cite{ErSch}, see Lemma~\ref{lem:jumpformula}. However, unlike the Schr\"odinger case considered in \cite{ErSch}, the resolvent $\bfR(\mu)$ does not decay as $|\mu|\rightarrow\infty$, which reflects the relativistic features of wave-type equations. For this reason we have to carry out the contour deformation argument in a manner that exploits the local energy decay, which is consistent with wave-type equations. 

Next, we determine the integral kernel  of the resolvent $\bfR(\mu)$ in Lemma~\ref{lem:Rexpression}. In this computation the integral kernel  of the resolvent of the scalar operator corresponding to the determinant of $\bfL_\ell - i \mu \bbI_{2 \times 2}$ arises naturally. The latter can in turn be computed from the well-known integral kernel  of the resolvent of the scalar linearized operator $-\partial^2 - 2 \sech^2(\cdot) + 1$ around the static sine-Gordon kink.

Upon inserting the resulting expressions back into \eqref{equ:overview_evolution_jump_across_resolvent}, we obtain a representation formula for the evolution $e^{t\bfL_\ell} \Pe \bmf$ of a pair of real-valued Schwartz functions $\bm{f} = (f_1,f_2) \in \calS(\bbR) \times \calS(\bbR)$ at any time $t \in \bbR$,
\begin{equation} \label{equ:overview_representation_evolution_dist_FT}
    \begin{aligned}
        \big( e^{t\bfL_\ell} P_\mathrm{e} \bm{f} \big)_1(y) &= \Re \, \biggl( \calF_\ell^\ast\Bigl[ e^{i t (\jxi + \ell \xi)} i\jxi^{-1} \Bigl( \calF_{\ell,D}[f_1](\xi) - \calF_\ell[f_2](\xi) \Bigr) \Bigr](y) \biggr), \\
        \big( e^{t\bfL_\ell} P_\mathrm{e} \bm{f} \big)_2(y) &= \Re \, \biggl( \calF_{\ell,D}^\ast \Bigl[ e^{i t (\jxi + \ell \xi)} i{\jxi}^{-1} \Bigl( \calF_{\ell,D}[f_1](\xi) - \calF_\ell[f_2](\xi) \Bigr) \Bigr](y) \biggr).
    \end{aligned}
\end{equation}
The preceding identities feature the following scalar transforms and their adjoints
\begin{align}
    \calF_\ell[g](\xi) :=  \int_\bbR \overline{e_\ell(y,\xi)} g(y) \, \ud y, \quad
    \calF_{\ell,D}[g](\xi) := \int_\bbR \overline{e_\ell(y,\xi)} D g(y) \, \ud y, \quad D := \ell \partial_y - i (\jxi + \ell \xi),
\end{align}
where
\begin{equation} \label{equ:overview_definition_dist_Fourier_basis_element}
\begin{aligned}
 e_\ell(y,\xi) := \frac{1}{\sqrt{2\pi}} e^{i y \xi}
  \begin{cases}
   \frac{\gamma(\xi+\ell\jxi) + i \tanh(\gamma y)}{\gamma(\xi+\ell\jxi) - i} & \xi \geq -\gamma \ell,
   \\
   \\
   \frac{\gamma(\xi+\ell\jxi) + i \tanh(\gamma y)}{\gamma(\xi+\ell\jxi) + i} & \xi < -\gamma \ell.
  \end{cases}
\end{aligned}
\end{equation}
In the present moving kink setting, we invite the reader to think of the frequency $\xi = -\gamma \ell$ for which $\xi+\ell\jxi = 0$ as the ``zero frequency'' in the static case.

For the sine-Gordon model the distorted Fourier basis \eqref{equ:overview_definition_dist_Fourier_basis_element} has a jump discontinuity at the frequency $\xi = -\gamma \ell$. Indeed, we have
\begin{equation}
    \lim_{\xi \downarrow -\gamma \ell} e_\ell(y,\xi) = - \lim_{\xi \uparrow -\gamma \ell} e_\ell(y,\xi) = - \frac{1}{\sqrt{2\pi}} e^{-iy\gamma\ell} \tanh(\gamma y).
\end{equation}
This discontinuity is related to the threshold resonances of the operator $\bfL_\ell$. However, for the
nonlinear analysis it is preferable to have better regularity properties of the distorted Fourier basis.
Following \cite{CP22} we remove the discontinuity by a simple sign change for frequencies $\xi < -\gamma \ell$. Introducing the modified distorted Fourier basis element
\begin{equation} \label{equ:overview_definition_dist_FT_basis_element}
    e_\ell^\#(y,\xi) := \frac{1}{\sqrt{2\pi}} \frac{\gamma(\xi+\ell\jxi) + i \tanh(\gamma y)}{|\gamma(\xi+\ell\jxi)| - i} e^{i y \xi}
\end{equation}
and the corresponding scalar transforms 
\begin{equation} \label{equ:overview_scalar_transforms_modified}
        \calF^\#_\ell[g](\xi) := \int_\bbR \overline{e^{\#}_\ell(y,\xi)} g(y) \, \ud y, \quad
        \calF^{\#}_{\ell,D}[g](\xi) := \int_\bbR \overline{e^{\#}_\ell(y,\xi)} D g(y) \, \ud y, \quad D := \ell \partial_y - i (\jxi + \ell \xi),
\end{equation}
we conclude from \eqref{equ:overview_representation_evolution_dist_FT} the following representation formula in terms of the modified distorted Fourier transform for the evolution $e^{t\bfL} P_\mathrm{e} \bmf$ of a pair of real-valued Schwartz functions $\bm{f} = (f_1,f_2) \in \calS(\bbR) \times \calS(\bbR)$ at any time $t \in \bbR$, 
\begin{equation} \label{equ:overview_representation_evolution_modified_dist_FT}
    \begin{aligned}
        \big( e^{t\bfL} P_\mathrm{e} \bm{f} \big)_1(y) &= \Re \, \biggl( \calF_\ell^{\#,\ast}\Bigl[ e^{i t (\jxi + \ell \xi)} i\jxi^{-1} \Bigl( \calF_{\ell,D}^{\#}[f_1](\xi) - \calF_\ell^{\#}[f_2](\xi) \Bigr) \Bigr](y) \biggr), \\
        \big( e^{t\bfL} P_\mathrm{e} \bm{f} \big)_2(y) &= \Re \, \biggl( \calF_{\ell,D}^{\#,\ast} \Bigl[ e^{i t (\jxi + \ell \xi)} i{\jxi}^{-1} \Bigl( \calF_{\ell,D}^{\#}[f_1](\xi) - \calF_\ell^{\#}[f_2](\xi) \Bigr) \Bigr](y) \biggr).
    \end{aligned}
\end{equation}
The modified distorted Fourier transform associated with the operator $\bfL_\ell$ should be thought of as the transformation 
\begin{equation} \label{equ:overview_dist_FT_vector_to_scalar}
    \bmf = (f_1, f_2) \in H^1_y(\bbR;\bbR) \times L^2_y(\bbR; \bbR) \quad \mapsto \quad \calF_{\ell,D}^{\#}\bigl[f_1\bigr] - \calF_{\ell}^{\#}\bigl[f_2\bigr] \in L^2_\xi(\bbR;\bbC),
\end{equation}
which maps a \emph{vector} of two real-valued functions to a \emph{scalar} complex-valued function.
For convenience, we will often omit the adjective ``modified'' in what follows, referring to \eqref{equ:overview_dist_FT_vector_to_scalar} simply as the distorted Fourier transform associated with $\bfL_\ell$.

Evaluating \eqref{equ:overview_representation_evolution_modified_dist_FT} at time $t=0$ yields the following inversion formula for the distorted Fourier transform on the physical side
\begin{equation} \label{equ:overview_inversion_formula_dist_FT_physical_side}
    \Pe = \Re \begin{bmatrix}
                    \calF_\ell^{\#,\ast}\Bigl[ i\jxi^{-1} \Bigl( \calF_{\ell,D}^{\#}[\cdot](\xi) - \calF_\ell^{\#}[\cdot](\xi) \Bigr) \Bigr] \\
                    \calF_{\ell,D}^{\#,\ast} \Bigl[ i{\jxi}^{-1} \Bigl( \calF_{\ell,D}^{\#}[\cdot](\xi) - \calF_\ell^{\#}[\cdot](\xi) \Bigr) \Bigr]
               \end{bmatrix}.
\end{equation}
At first sight it is not obvious what the corresponding inversion formula looks like on the distorted Fourier side. In Lemma~\ref{lem:inversionXi} we determine that it reads 
\begin{equation} \label{equ:overview_inversion_formula_dist_FT_frequency_side}
        \mathrm{Id}_\xi = \Bigl( \calF_{\ell, D}^{\#} \, \Re \, \calF_{\ell}^{\#,\ast} - \calF_{\ell}^{\#} \, \Re \, \calF_{\ell, D}^{\#,\ast} \Bigr) i \jxi^{-1}.
\end{equation}
This identity is key to study the evolution equation for the radiation term on the distorted Fourier side, where almost all of the nonlinear analysis is carried out in this paper. 

The preceding properties of the distorted Fourier transform are established in Subsections~\ref{subsec:vectorial_dist_FT} and \ref{subsec:fourier_inversion}.
In Subsection~\ref{subsec:higher_order_sobolev} we show the equivalence between Sobolev norms on the physical side and weighted norms on the frequency side. Specifically for any integer $k \geq 1$ and any $\bmf \in H^k_y \times H^{k-1}_y$ with $\bmf = \Pe \bmf$, we have 
\begin{equation}
    \bigl\| \bmf \bigr\|_{H^k_y \times H^{k-1}_y} \simeq \bigl\| \jxi^k \bigl( \calF_{\ell,D}^{\#}\bigl[f_1\bigr] - \calF_{\ell}^{\#}\bigl[f_2\bigr] \bigr) \bigr\|_{L^2_\xi}.
\end{equation}
Finally, in Subsection~\ref{subsec:mapping_properties} we prove further mapping properties of the distorted Fourier transform that will be needed for the nonlinear analysis.

\medskip 

In Section~\ref{sec:linear_decay_estimates} the oscillatory integral representation \eqref{equ:overview_representation_evolution_modified_dist_FT} of the evolution in terms of the distorted Fourier transform is the point of departure for the derivation of linear decay estimates.
In Subsection~\ref{subsec:pointwise_decay} we establish pointwise estimates and asymptotics. In particular, we find that for any pair of real-valued functions $\bmf = (f_1, f_2) \in \calS \times \calS$ and for all times $t \geq 1$,
\begin{equation} \label{equ:overview_linear_evolution_dispersive_decay}
    \begin{aligned}
        &\bigl\|\bigl( e^{t\bfL_\ell} P_\mathrm{e} \bm{f} \bigr)_1\bigr\|_{L^\infty_y} + \bigl\|\bigl( e^{t\bfL_\ell} P_\mathrm{e} \bm{f} \bigr)_2\bigr\|_{L^\infty_y} + \bigl\|\py \bigl( e^{t\bfL_\ell} P_\mathrm{e} \bm{f} \bigr)_1\bigr\|_{L^\infty_y} \\
        &\quad \lesssim \frac{1}{t^{\frac12}} \bigl\| \jxi^{\frac32} g_\ell^{\#}(\xi) \bigr\|_{L^\infty_\xi} + \frac{1}{t^{\frac23}} \Bigl( \bigl\|\jxi^2 \partial_\xi g_\ell^{\#}(\xi) \bigr\|_{L^2_\xi} + \bigl\|\jxi^2 g_\ell^{\#}(\xi) \bigr\|_{L^2_\xi} \Bigr),
    \end{aligned}
\end{equation}
where we use the short-hand notation
\begin{equation} \label{equ:overview_dispersive_decay_effective_profile_naturally}
    g_\ell^{\#}(\xi) := \calF_{\ell,D}^{\#}\bigl[f_1\bigr](\xi) - \calF_{\ell}^{\#}\bigl[f_2\bigr](\xi).
\end{equation}
As we aim to develop a space-time resonances approach for the moving kink setting, the preceding dispersive decay estimate naturally identifies \eqref{equ:overview_dispersive_decay_effective_profile_naturally} as a suitable notion of ``profile'' on the distorted Fourier side.
Moreover, an inspection of the right-hand side of \eqref{equ:overview_linear_evolution_dispersive_decay} then points to the norms of this profile that must be propagated through a bootstrap argument.

Moreover, we establish several improved local decay estimates for the evolution in Subsection~\ref{subsec:improved_local_decay_estimates} and we prove an integrated local energy decay estimate with a moving center in Subsection~\ref{subsec:ILED}, which we hope to be of independent interest. In this work the latter is crucial later on to close certain weighted energy estimates for the contributions of spatially localized terms with cubic-type time decay in the proofs of Lemma~\ref{lem:pxi_cubic} and of Lemma~\ref{lem:pxi_calR}.

\subsection{Modulation and evolution equation for the effective profile} \label{subsec:overview_modulation}

Our study of the asymptotic stability of the family of moving kink solutions~\eqref{equ:intro_general_moving_kink_family} begins in earnest in Section~\ref{sec:setting_up}, where we prepare the coupled system of the evolution equation for the radiation term and of the modulation equations for the nonlinear analysis.

The proof of Theorem~\ref{thm:main} proceeds via a grand bootstrap argument. 
On a given bootstrap time interval $[0,T]$, we use standard modulation techniques to decompose the solution to \eqref{equ:intro_sG_1st_order} in a neighborhood of the kink into a modulated kink and a radiation term 
\begin{equation}
    \bmphi(t,x) = \bmK_{\ell(t),q(t)}(x) + \bmu(t, x-q(t)), \quad 0 \leq t \leq T,
\end{equation}
where 
        \begin{equation} \label{equ:overview_definition_modulated_kink}
            \bmK_{\ell(t),q(t)}(x) := \begin{bmatrix}
                        K\bigl( \gamma(t) (x-q(t)) \bigr) \\
                        - \gamma(t) \ell(t) K'\bigl( \gamma(t) (x-q(t)) \bigr)
                     \end{bmatrix}, 
                     \quad \bmu = \begin{bmatrix} u_1 \\ u_2 \end{bmatrix},
                     \quad \gamma(t) := \frac{1}{\sqrt{1-\ell(t)^2}}.
        \end{equation}     
The modulation parameters $\ell \colon [0,T] \to (-1,1)$ and $q \colon [0,T] \to \bbR$ are uniquely chosen so that the radiation term satisfies the orthogonality properties
        \begin{equation} 
            \bigl\langle \bfJ \partial_q \bm{K}_{\ell,q}, \bmu(t,\cdot-q) \bigr\rangle = \bigl\langle \bfJ \partial_\ell \bm{K}_{\ell,q}, \bmu(t,\cdot-q) \bigr\rangle = 0, \quad 0 \leq t \leq T, \quad \bfJ = \begin{bmatrix} 0 & 1 \\ -1 & 0 \end{bmatrix}.
        \end{equation}
Passing to the moving frame coordinate $y := x - q(t)$ and viewing the radiation term $\bmu$ from now on as a function of $t$ and $y$, we arrive at the following evolution equation for the radiation term
    \begin{equation} \label{equ:overview_evol_equ_radiation_term}
        \partial_t \bmu - \bfL_{\ell} \bmu = (\dot{q} - \ell) \py \bmu + {\mathcal Mod} + \calN(\bmu)
    \end{equation}
coupled to the first-order system of modulation equations 
    \begin{equation} \label{equ:overview_modulation_equations}
            \bbM_{\ell}[\bmu] \begin{bmatrix} \dot{\ell} \\ \dot{q}-\ell \end{bmatrix} =
        \begin{bmatrix}
            - \bigl\langle \bfJ \bmY_{0,\ell}, \calN(\bm{u}) \bigr\rangle \\
            \bigl\langle \bfJ \bmY_{1,\ell}, \calN(\bm{u}) \bigr\rangle 
        \end{bmatrix}.
    \end{equation}
Here, $\bfL_\ell$ denotes the linearized matrix operator defined in \eqref{equ:overview_definition_bfLell}, ${\mathcal Mod}$ are mild modulation terms, $\calN(\bmu)$ is the nonlinearity, and $\bbM_{\ell}[\bmu]$ is an invertible matrix of the schematic form
\begin{equation}
    \bbM_{\ell}[\bmu] = \gamma^3 \|K'\|_{L^2} \begin{bmatrix} 1 & 0 \\ 0 & 1 \end{bmatrix} + \calO\bigl( \|\bmu\|_{L^2} \bigr).
\end{equation}
It is important to keep in mind that in the preceding identities the modulation parameters are time-dependent.

In order to study the long-time behavior of the radiation term via spectral methods, we need to pass to a time-independent reference operator $\bfL_\ulell$ for a suitably chosen fixed value $\ulell$ of the path $\ell(t)$. For the purposes of the heuristic discussion in this section, the reader should think of $\ulell$ as the final value $\ell(T)$ of the path $\ell(t)$ on the given bootstrap interval $[0,T]$. 
Focusing from now on only on the essential features of the system \eqref{equ:overview_evol_equ_radiation_term}--\eqref{equ:overview_modulation_equations}, we disregard various milder terms and only consider the following schematic evolution equation for the radiation term
    \begin{equation} \label{equ:overview_evol_equ_radiation_term_essence}
        \pt \bmu - \bfL_\ulell \bmu = (\dot{q}-\ulell) \py \bmu + \calQ_\ulell(\bmu)  + \calC(\bmu)  + \ldots 
    \end{equation}
coupled to the following schematic first-order system of differential equations for the modulation parameters
    \begin{equation} \label{equ:overview_evol_equ_modulation_equations_essence}
        \begin{bmatrix} \ulg^3 \|K'\|_{L^2} & 0 \\ 0 & \ulg^3 \|K'\|_{L^2} \end{bmatrix}
        \begin{bmatrix} \dot{\ell} \\ \dot{q}-\ell \end{bmatrix} 
        =
        \begin{bmatrix}
            - \bigl\langle \bfJ \bmY_{0,\ulell}, \calQ_\ulell(\bmu) \bigr\rangle \\
            \bigl\langle \bfJ \bmY_{1,\ulell}, \calQ_\ulell(\bmu) \bigr\rangle 
        \end{bmatrix} + \ldots,
    \end{equation} 
where the leading order quadratic and cubic nonlinearities\footnote{We caution the reader that the notation introduced in \eqref{equ:overview_definition_quad_cubic_nonlin} is useful in this section for the purposes of providing an overview of the proof of Theorem~\ref{thm:main}, but that we introduce slightly different definitions for these nonlinearities in Subsection~\ref{subsec:evolution_equation_profile}.} are given by 
\begin{equation} \label{equ:overview_definition_quad_cubic_nonlin}
    \calQ_\ulell(\bmu) := \begin{bmatrix} 0 \\ \alpha(\ulg y) u_1^2 \end{bmatrix}, \quad \calC(\bmu) := \begin{bmatrix} 0 \\ {\textstyle \frac16} u_1^3 \end{bmatrix}, \quad \alpha(\ulg y) := - \sech(\ulg y) \tanh(\ulg y).
\end{equation}

At this stage our task is to simultaneously prove dispersive decay for the radiation term and to deduce sufficient decay for $\dot{\ell}$ as well as $\dot{q}-\ell$ to infer convergence of the modulation parameters to the final moving kink parameters.
Taking a space-time resonances approach, we define the profile of the radiation term relative to the reference operator $\bfL_\ulell$ on the physical side by
    \begin{equation}
        \bmf_\ulell(t) = \bigl(f_{\ulell,1}(t), f_{\ulell,2}(t) \bigr) := e^{-t\bfL_\ulell} \bigl( \ulPe \bmu(t) \bigr),
    \end{equation}
where $\ulPe$ denotes the projection to the essential spectrum relative to $\bfL_\ulell$. As anticipated in the preceding subsection, we then consider the distorted Fourier transform of the profile
    \begin{equation}
        \gulellsh(t,\xi) := \calFulellDsh\bigl[ f_{\ulell,1}(t) \bigr](\xi) - \calFulellsh\bigl[ f_{\ulell,2}(t) \bigr](\xi),
    \end{equation}
which from now on we refer to as the (complexified) \emph{effective profile} on the distorted Fourier side.

Using the inversion formula on the distorted Fourier side \eqref{equ:overview_inversion_formula_dist_FT_frequency_side} along with some basic properties of the scalar transforms \eqref{equ:overview_scalar_transforms_modified}, we conclude from \eqref{equ:overview_evol_equ_radiation_term_essence} the following \emph{scalar} evolution equation for the effective profile       
    \begin{equation} \label{equ:overview_gulellsh_evolution_equation}
        \begin{aligned}
        &\pt \gulellsh(t,\xi) \\
        &\quad =  \bigl( \dot{q} - \ulell \bigr) i\xi \gulellsh(t,\xi) - e^{-it(\jxi+\ulell\xi)} \Bigl( \calFulellsh\bigl[ \alpha(\ulg \cdot) u_1(t,\cdot)^2 \bigr](\xi) + \calFulellsh\bigl[ {\textstyle \frac16} u_1(t,\cdot)^3 \bigr](\xi) \Bigr) + \ldots 
        \end{aligned}
    \end{equation}
If we ignore the first term on the right-hand side of \eqref{equ:overview_gulellsh_evolution_equation}, this equation for the effective profile is remarkably similar to the corresponding evolution equation for the (scalar) distorted Fourier transform of the profile of (scalar) perturbations of \emph{static} kinks.
The evolution equation for such scalar perturbations of static kinks is of the schematic form
\begin{equation} \label{equ:overview_scalar_perturbation_NLKG}
    \bigl( \pt^2 - \px^2 + V(x) + m^2 \bigr) u = \alpha(x) u^2 + \beta u^3 + \ldots,
\end{equation}
where $V(x)$ is a smooth decaying potential, $\alpha(x)$ can often be thought of as a spatially localized variable coefficient, and $\beta \in \bbR$ is a constant coefficient.
Thanks to major advances in recent years, see for instance \cite{LLS1, LLS2, LLSS, LS1, LS2, GP20, GermPusZhang22, CP21, CP22}, the decay and the asymptotics of small solutions to the deceptively simple looking nonlinear Klein-Gordon equation~\eqref{equ:overview_scalar_perturbation_NLKG} is by now well-understood in the generic setting and for odd perturbations of the sine-Gordon kink. However, in the presence of threshold resonances or internal modes, many challenging open questions remain.

For the analysis of the evolution equation of the effective profile \eqref{equ:overview_gulellsh_evolution_equation}, we build on these developments and extend them to the moving kink setting.
In particular, we adopt ideas and techniques from \cite{GP20, LS1, CP22}.
A significant difference in our work in comparison to the preceding references is obviously the coupling of \eqref{equ:overview_gulellsh_evolution_equation} to the modulation equations~\eqref{equ:overview_evol_equ_modulation_equations_essence}. The dynamics of the modulation parameters impacts the radiation term especially through the first term on the right-hand side of \eqref{equ:overview_gulellsh_evolution_equation}.
For the nonlinear analysis, it is helpful to absorb this term into the phase via an integrating factor as in \cite[Proposition 9.5]{CG23}.
Introducing
\begin{equation} \label{equ:overview_definition_theta}
    \theta(t) := \int_0^t \bigl( \dot{q}(s) - \ulell \bigr) \, \ud s,
\end{equation}
we then arrive at the following renormalized evolution equation for the effective profile
\begin{equation} \label{equ:overview_evol_equ_radiation_term_renormalized}
    \begin{aligned}
    &\pt \Bigl( e^{-i\xi\theta(t)} \gulellsh(t,\xi) \Bigr) \\
    &\quad = - e^{-i\xi\theta(t)} e^{-it(\jxi+\ulell\xi)} \Bigl( \calFulellsh\bigl[ \alpha(\ulg \cdot) u_1(t,\cdot)^2 \bigr](\xi) + \calFulellsh\bigl[ {\textstyle \frac16} u_1(t,\cdot)^3 \bigr](\xi) \Bigr) + \ldots
    \end{aligned}
\end{equation}

Heuristically, the quantity \eqref{equ:overview_definition_theta} measures how sharply we control the center of the moving kink. A uniform-in-time bound on \eqref{equ:overview_definition_theta} would correspond to a sharp localization of the center, while growth bounds on \eqref{equ:overview_definition_theta} would quantify the uncertainty of the localization of the center in our analysis. 

Writing $\theta(t) = \int_0^t \bigl( \dot{q}(s) - \ell(s) \bigr) \, \ud s + \int_0^t \bigl( \ell(s) - \ulell \bigr) \, \ud s$ and keeping in mind that $\ulell = \ell(T)$ for the purposes of this discussion, in order to deduce uniform-in-time bounds on $\theta(t)$ without taking into account any favorable oscillations, we would need $\dot{q}(s) - \ell(s)$ and $\ell(s) - \ell(T)$ to decay faster than $\varepsilon \js^{-1}$ for $0 \leq s \leq T$.
The latter would in turn follow crudely from writing $\ell(s)-\ell(T) = -\int_s^T \dot{\ell}(\tau) \, \ud \tau$ if $\dot{\ell}(\tau)$ was twice integrable in time, i.e., if $\dot{\ell}(\tau)$ decayed faster than $\varepsilon \jap{\tau}^{-2}$ for $0 \leq \tau \leq T$.
Looking ahead to the next subsection, we remark that such decay rates are achievable in the generic setting, but they appear difficult to attain in the presence of threshold resonances.

\subsection{Slow local decay and the null structures} \label{subsec:overview_slow_local_decay}

If we had a uniform-in-time bound on \eqref{equ:overview_definition_theta} and if the contribution of the quadratic nonlinearity on the right-hand side of \eqref{equ:overview_evol_equ_radiation_term_renormalized} was milder than that of the cubic nonlinearity, the proof of Theorem~\ref{thm:main} would at this stage be reduced to capturing the well-understood long-range effects of the cubic nonlinearity.
To assess how close we are to this scenario, the key observation is that the quadratic nonlinearity on the right-hand side of \eqref{equ:overview_evol_equ_radiation_term_renormalized} features the spatially localized variable coefficient $\alpha(\ulg y)$ and that the leading order quadratic nonlinearities on the right-hand side of the modulation equations~\eqref{equ:overview_evol_equ_modulation_equations_essence} are also spatially localized.
This highlights the critical role of the \emph{local decay} of the radiation term for the entire nonlinear analysis in this paper.

However, the local decay of the radiation term is drastically affected by the presence of the threshold resonances of the linearized operator $\bfL_\ulell$.
They induce slow local decay, which is captured sharply by the following refined local decay estimate\footnote{In this work we establish versions of the refined local decay estimate \eqref{equ:overview_refined_local_decay_estimate} adapted to the limited regularity of the effective profile that we propagate in the nonlinear analysis, see Lemma~\ref{lem:improved_local_decay} and Lemma~\ref{lem:local_decay_resonance_subtracted_off}. Minor modifications of the proofs of these lemmas and some additional straightforward stationary phase arguments would yield a proof of \eqref{equ:overview_refined_local_decay_estimate}. We refer to \cite[Corollary 2.17]{LLSS} for a proof of the analogue of the refined local decay estimate \eqref{equ:overview_refined_local_decay_estimate} for (scalar) linearized operators around static kinks.} for the evolution
\begin{equation} \label{equ:overview_refined_local_decay_estimate}
    \begin{aligned}
        \biggl\| \jy^{-5} \biggl( \bigl( e^{t\bfL_\ulell} \ulPe \bmf \bigr)(y) 
        - c_{+, \ulell} \frac{e^{it\ulg^{-1}}}{t^{\frac12}} \bm{\Phi}_{+,\ulell}(y) c_{+}(\bmf) &- c_{-, \ulell} \frac{e^{-it\ulg^{-1}}}{t^{\frac12}} \bm{\Phi}_{-,\ulell}(y) c_{-}(\bmf) \biggr) \biggr\|_{L^2_y} \\
        &\qquad \qquad \qquad \qquad \quad \lesssim \frac{1}{t^{\frac32}} \bigl\| \jy^5 \bmf \bigr\|_{L^2_y},
    \end{aligned}
\end{equation}
where $c_{\pm}(\bmf)$ should be thought of as projections of $\bmf$ onto the threshold resonances $\bm{\Phi}_{\pm,\ulell}$ and where $c_{\pm,\ulell}$ are constants depending only on the value of $\ulell$.
In the absence of threshold resonances for the linearized operator, the entire evolution would enjoy the improved local decay $t^{-\frac32}$ on the right-hand side of \eqref{equ:overview_refined_local_decay_estimate}, which would significantly simplify the analysis of \eqref{equ:overview_evol_equ_radiation_term_renormalized} as explained above.
Instead, in our setting the threshold resonances of the linearized operator $\bfL_\ulell$ around the moving sine-Gordon kink lead to the slowly decaying terms on the left-hand side of \eqref{equ:overview_refined_local_decay_estimate}. These should be thought of as projections of the evolution onto the threshold resonances. 
Hence, to leading order the quadratic nonlinearities on the right-hand sides of \eqref{equ:overview_evol_equ_modulation_equations_essence} and \eqref{equ:overview_evol_equ_radiation_term_renormalized} exhibit slowly decaying source terms of the schematic form
\begin{equation} \label{equ:overview_source_terms}
    \biggl( \frac{e^{\pm i t \ulg^{-1}}}{t^{\frac12}} \bm{\Phi}_{\pm, \ulell}(y) \varepsilon \biggr)^2.
\end{equation}
At first sight these source terms are highly problematic, because they could potentially lead to additional corrections in the decay rate of the radiation term and in the dynamics of the modulation parameters.
However, in the case of the sine-Gordon model the worst effects of these source terms are suppressed owing to two remarkable null structures of the quadratic nonlinearities in the evolution equation for the effective profile as well as in the modulation equations.

\medskip 

\noindent {\it Null structure I: Evolution equation for the effective profile.}
In view of the distinct oscillations of the quadratic source terms \eqref{equ:overview_source_terms}, it turns out that the contribution of the quadratic nonlinearity to the evolution equation for the effective profile \eqref{equ:overview_evol_equ_radiation_term_renormalized} features critically decaying resonant terms of the schematic form
\begin{equation}
    e^{-it(\jxi+\ulell\xi)} e^{it2\ulg^{-1}} \, \calFulellsh\Bigl[\alpha(\ulg \cdot) \bigl( \bmPhi_{+,\ulell,1}(y)\bigr)^2 \Bigr](\xi) \, \frac{\varepsilon^2}{t},
\end{equation}
whose phase function vanishes at specific frequencies
\begin{equation}
    \jxi + \ulell \xi - 2\ulg^{-1} = 0 \quad \Leftrightarrow \quad \xi = \ulg (-2\ulell \pm \sqrt{3}).
\end{equation}
It turns out that the scalar distorted transform $\calFulellsh\bigl[\alpha(\ulg \cdot) \bigl( \bmPhi_{+,\ulell,1}(y)\bigr)^2 \bigr](\xi)$ vanishes exactly at these problematic frequencies $\xi = \ulg (-2\ulell \pm \sqrt{3})$, as we uncover in Lemma~\ref{lem:null_structure1}. This is a generalized version of the null structure observed in \cite[Lemma 3.1]{LS1} and \cite[Remark 1.2]{LLSS} for the special case of odd perturbations of the static sine-Gordon kink. In fact, this null structure is only the ``tip of the iceberg'' of an all-encompassing null structure of the quadratic nonlinearity in \eqref{equ:overview_evol_equ_radiation_term_renormalized}, as we explain in comment (ii) in Subsection~\ref{subsec:further_comments} and in Appendix~\ref{appendix:deeper_null_structure}.

Thanks to the remarkable null structure from Lemma~\ref{lem:null_structure1}, we can use a variable coefficient quadratic normal form to transform the contributions of the quadratic source terms on the right-hand side of \eqref{equ:overview_evol_equ_radiation_term_renormalized} into spatially localized terms with cubic-type time decay, which are slightly better behaved than the constant coefficient cubic terms on the right-hand side of \eqref{equ:overview_evol_equ_radiation_term_renormalized}. We refer to Subsection~\ref{subsec:normal_form_transformation} for the details. 

\medskip 

\noindent {\it Null structure II: Modulation equations.}   
As explained at the end of the preceding Subsection~\ref{subsec:overview_modulation}, the main challenge in the analysis of the dynamics of the modulation parameters is to deduce sufficiently fast convergence of the Lorentz boost parameter $\ell(t)$ to its final value $\ell(T)$ on the given bootstrap interval $[0,T]$. 
As in \cite{Chen21, CG23, LL2}, to this end we try to exploit the oscillations in the modulation equation for $\dot{\ell}$.
Using the fundamental theorem of calculus and inserting \eqref{equ:overview_evol_equ_modulation_equations_essence}, we obtain to leading order for $0 \leq t \leq T$ that
\begin{equation} \label{equ:overview_null_structure2_equ1}
    \ell(t) - \ell(T) = - \int_t^T \dot{\ell}(s) \, \ud s = \ulg^{-3} \|K'\|_{L^2}^{-2} \int_t^T \bigl\langle \bfJ \bmY_{0,\ulell}, \calQ_\ulell(\bmu) \bigr\rangle \, \ud s + \ldots 
\end{equation}
Inserting \eqref{equ:overview_representation_evolution_modified_dist_FT} to express the radiation term as an oscillatory integral in terms of the effective profile, the preceding equation reads schematically
\begin{equation} \label{equ:overview_null_structure2_equ2}
    \ell(t) - \ell(T) \approx \int_t^T \iint_{\bbR^2} e^{is(\jap{\xi_1}+\ulell\xi_1 - \jap{\xi_2} - \ulell\xi_2)} \gulellsh(s,\xi_1) \overline{\gulellsh(s,\xi_2)} \, \mu_{\ulell;+-}(\xi_1, \xi_2) \, \ud \xi_1 \, \ud \xi_2 \, \ud s + \ldots,
\end{equation}
where we only display the most delicate resonant quadratic contribution on the right-hand side.
The quadratic spectral distribution $\mu_{\ulell;+-}(\xi_1,\xi_2)$ in \eqref{equ:overview_null_structure2_equ2} results from computing the interaction of the distorted Fourier basis elements~\eqref{equ:overview_definition_dist_FT_basis_element}, the variable coefficient $\alpha(\ulg y)$ of the quadratic nonlinearity in \eqref{equ:overview_definition_quad_cubic_nonlin}, and the generalized kernel element $\bmY_{0,\ulell}(y)$. 
Clearly, the phase in \eqref{equ:overview_null_structure2_equ2} has time resonances. Remarkably, as we uncover by an explicit computation in Subsection~\ref{subsec:quadratic_spectral_distributions}, the quadratic spectral distribution $\mu_{\ulell;+-}(\xi_1,\xi_2)$ is divisible by the phase $(\jap{\xi_1}+\ulell\xi_1 - \jap{\xi_2} - \ulell\xi_2)$. Thus, we can use a normal form to turn the quadratic contribution \eqref{equ:overview_null_structure2_equ2} into a cubic contribution so that \eqref{equ:overview_null_structure2_equ2} is effectively of the form
\begin{equation} \label{equ:overview_null_structure2_equ3}
    \ell(t) - \ell(T) \approx \int_t^T \bigl\langle \bfJ \bmY_{0,\ulell}, \calC(\bmu) \bigr\rangle \, \ud s + \ldots 
\end{equation}
Then we observe that the cubic interactions of the leading order local decay behavior of the radiation term do not exhibit time resonances, whence we can use another normal form to transform \eqref{equ:overview_null_structure2_equ3} into quartic contributions. These can be easily shown to enjoy almost twice integrable time decay $\lesssim \varepsilon^4 \js^{-2}$, which just about suffices for the nonlinear analysis.
See the treatment of the contributions of the term $I(s)$ in Section~\ref{sec:modulation_control} for the details.

\subsection{Bootstrap setup} \label{subsec:overview_bootstrap_setup}

We now describe in more detail the bootstrap setup for the proof of Theorem~\ref{thm:main} in order to simultaneously establish dispersive decay for the radiation term and to deduce sufficiently fast convergence of the modulation parameters to those of the final moving kink.

Informed by the pointwise decay estimate~\eqref{equ:overview_linear_evolution_dispersive_decay} discussed at the end of Subsection~\ref{subsec:overview_distorted_FT}, we introduce the following collection of pointwise and weighted energy norms for the effective profile on a given time interval $[0,T]$,
\begin{equation} \label{equ:overview_definition_XT_norm}
    \begin{aligned}
        \bigl\| g_\ulell^\# \bigr\|_{X(T)} &:= \sup_{0 \leq t \leq T} \, \biggl( \bigl\| \jxi^{\frac32} g_\ulell^\#(t,\xi) \bigr\|_{L^\infty_\xi} + \jt^{-\delta} \bigl\| \jxi^2 g_\ulell^\#(t,\xi) \bigr\|_{L^2_\xi} \\
        &\qquad \qquad \quad + \jt^{-\delta} \bigl\| \pxi g_\ulell^\#(t,\xi) \bigr\|_{L^2_\xi} + \jt^{-2\delta} \bigl\| \jxi^2 \pxi g_\ulell^\#(t,\xi) \bigr\|_{L^2_\xi} \biggr).
    \end{aligned}
\end{equation}
Our asymptotic stability proof then boils down to establishing the following core bounds for the effective profile and for the modulation parameters on a given time interval $[0,T]$ via an intertwined bootstrap argument 
\begin{equation} \label{equ:overview_summary_bootstrap_bounds}
    \begin{aligned}
        \bigl\| g_\ulell^\# \bigr\|_{X(T)} \lesssim \varepsilon, \quad 
        \sup_{0 \leq t \leq T} \, \jt^{1-\delta} \bigl| \ell(t) - \ell(T) \bigr| \lesssim \varepsilon, \quad
        \sup_{0 \leq t \leq T} \, \jt \bigl| \dot{q}(t) - \ell(t) \bigr| \lesssim \varepsilon.
    \end{aligned}
\end{equation}
Here, $0 < \varepsilon \ll 1$ measures the size of the initial perturbation in the weighted Sobolev norm \eqref{equ:statement_theorem_smallness_data} and $0 < \delta \ll 1$ is a small absolute constant.  

Our framework is structured around two central bootstrap propositions, presented in Subsection~\ref{subsec:bootstrap_propositions} as Proposition~\ref{prop:modulation_parameter_control} and Proposition~\ref{prop:profile_bounds}.
At that stage, as direct consequences of the bootstrap assumptions, we also derive various decay estimates and other bounds in Corollary~\ref{cor:consequences_bootstrap_assumptions}, which are repeatedly used throughout the nonlinear analysis in all subsequent sections.
Based on the conclusions of the two core bootstrap propositions, we provide a succinct proof of Theorem~\ref{thm:main} in Subsection~\ref{subsec:proof_of_main_theorem}.
The purpose of the remainder of this paper, namely Sections~\ref{sec:modulation_control}--\ref{sec:pointwise_profile_bounds}, is then to prove these two propositions.
We conclude the description of our bootstrap setup by pointing out important aspects of their proofs.
In particular, we clarify the origin of the two distinct slow growth rates for the weighted energy norms in the second line of \eqref{equ:overview_definition_XT_norm}, which may initially seem surprising to readers familiar with the usual modified scattering techniques.

\medskip 
\noindent {\it Control of the modulation parameters (Section~\ref{sec:modulation_control}).}
While the decay estimate $|\dot{q}(t)-\ell(t)| \lesssim \varepsilon \jt^{-1}$ is a straightforward consequence of the dispersive decay $\varepsilon \jt^{-\frac12}$ of the radiation term and the fact that the nonlinearity on the right-hand side of the modulation equations \eqref{equ:overview_evol_equ_modulation_equations_essence} is at least quadratic,
establishing the asserted convergence rate $|\ell(t) - \ell(T)| \lesssim \varepsilon \jt^{-1+\delta}$ is much more subtle.
As already described at the end of Subsection~\ref{subsec:overview_slow_local_decay}, the proof of the latter exploits the oscillations of $\dot{\ell}(t)$ and hinges on a remarkable null structure of the quadratic nonlinearity of the modulation equations. 
Moreover, the analysis of the contributions of the cubic nonlinearity relies on a leading order local decay decomposition of the radiation term into a slowly decaying term with distinct oscillations (corresponding to the slowly decaying terms that get subtracted on the left-hand side of \eqref{equ:overview_refined_local_decay_estimate}) and a bulk term with improved local decay, see the decompositions \eqref{equ:consequences_decomposition_leading_order_local_decay}, \eqref{equ:consequences_decomposition_leading_order_local_decay_usubetwo}, and \eqref{equ:consequences_aux_KG_leading_order_local_decay_decomp}.

\medskip  
\noindent {\it Weighted energy estimates for the profile (Section~\ref{sec:energy_estimates}).}
In view of the non-integrable decay rates for the two bounds on the modulation parameters in \eqref{equ:overview_summary_bootstrap_bounds}, we propagate a slowly growing bound $|\theta(t)| \lesssim \varepsilon \jt^{\delta}$ on the quantity $\theta(t)$ defined in \eqref{equ:overview_definition_theta}.
This small uncertainty\footnote{In the generic setting (no threshold resonances, no internal modes) we expect to be able to propagate integrable decay rates for the modulation parameters in \eqref{equ:overview_summary_bootstrap_bounds} and thus obtain sharp control of the center of the moving kink.} in the center of the moving kink within our proof turns out to be the reason for the two differing slow growth rates for the weighted energy norms in the second line of \eqref{equ:overview_definition_XT_norm}. See the estimates \eqref{equ:pxi_profile_differing_slow_growth_rates1} and \eqref{equ:pxi_profile_differing_slow_growth_rates2} for the details.

The weighted energy norm of the profile with additional Sobolev weights $\|\jxi^2 \pxi \gulellsh(t,\xi)\|_{L^2_\xi} \lesssim \varepsilon \jt^{2\delta}$ is primarily needed to estimate the faster decaying remainder term in the dispersive decay estimate \eqref{equ:overview_linear_evolution_dispersive_decay}.
In order to be able to close the overall bootstrap argument in this paper, it is key that in the proof of the convergence rate estimate $|\ell(t) - \ell(T)| \lesssim \varepsilon \jt^{-1+\delta}$ for the Lorentz boost parameter,
we only need to make use of the improved local decay estimates \eqref{equ:consequences_Remusubeone_H1y_local_decay}, \eqref{equ:consequences_Remusubetwo_L2y_local_decay}, \eqref{equ:consequences_aux_KG_Remv_improved_H1y_local_decay},  
because for their proofs the weighted energy norm of the profile without additional Sobolev weights $\|\pxi \gulellsh(t,\xi)\|_{L^2_\xi} \lesssim \varepsilon \jt^\delta$ suffices.

There are essentially two types of terms on the right-hand side of the renormalized evolution equation for the effective profile \eqref{equ:overview_evol_equ_radiation_term_renormalized}, namely the contributions of various spatially localized terms with at least cubic-type time decay $\jt^{-\frac32+\delta}$ and the contribution of the constant coefficient cubic term $\frac16 u_1^3$. 
In order to estimate the contributions of the localized terms with cubic-type time decay to the weighted energy bounds, we adopt ideas from \cite{CP21, CP22, LLS2, LS1} to the moving kink setting and use integrated local energy decay estimates with a moving center. We prove the latter in Subsection~\ref{subsec:ILED} and hope that they are of independent interest. 
Instead, in order to estimate the contributions of the constant coefficient cubic nonlinearity to the weighted energy bounds, we rely on its decomposition into cubic interactions with a Dirac kernel, with a Hilbert-type kernel, or with a regular kernel (spatially localized terms), see Subsection~\ref{subsec:structure_cubic_nonlinearities}. 
While we also use the integrated local energy decay estimates with a moving center for the interactions with a regular kernel, we exploit the structure of the cubic interactions with a Dirac or a Hilbert-type kernel to transfer the frequency derivative onto the inputs up to better behaved terms.

\medskip 
\noindent {\it Pointwise estimates for the profile (Section~\ref{sec:pointwise_profile_bounds}).}
In order to establish the uniform-in-time pointwise bound on the effective profile corresponding to the first term in \eqref{equ:overview_definition_XT_norm}, we view the renormalized evolution equation \eqref{equ:overview_evol_equ_radiation_term_renormalized} for the effective profile as an ordinary differential equation for every fixed frequency.
Then all spatially localized terms with cubic-type time decay on the right-hand side of \eqref{equ:overview_evol_equ_radiation_term_renormalized} are integrable error terms so that only the constant coefficient cubic nonlinearity contributes to leading order.
In view of its fine structure determined in Subsection~\ref{subsec:structure_cubic_nonlinearities}, in fact only the cubic interactions with a Dirac kernel or a Hilbert-type kernel contribute to leading order.
Using a stationary phase analysis we extract their leading order contributions. Interestingly, the cubic interactions with a Hilbert-type kernel do not give a resonant critically decaying leading order term (no time oscillations) due to a vanishing property of the corresponding coefficient, see \eqref{equ:cubic_spectral_distributions_pv_vanishing}. 
Using a standard integrating factor argument, the resonant leading order term stemming from the cubic interactions with a Dirac kernel can be absorbed into the phase, which captures the logarithmic phase corrections in the asymptotics of the radiation term.

\subsection{Further comments} \label{subsec:further_comments}

We conclude the overview of the proof of Theorem~\ref{thm:main} with a few more mostly technical comments.

\begin{enumerate}[label=(\roman*)]
    
    \item It is an interesting question whether the strategy of the proof of Theorem~\ref{thm:main} can be enhanced to deduce integrable decay rates for the modulation parameters in \eqref{equ:statement_theorem_modulation_parameters_asymptotics}.

    \item During the work on this paper we identified a notable ``deeper null structure'' of the quadratic nonlinearities in the evolution equation for the radiation term. This structure allows for the removal of all quadratic nonlinearities in the evolution equation for the radiation term using a suitable normal form argument. 
    However, in proving Theorem 1.1, we do not rely on this feature, as it is sufficient for closing the overall bootstrap argument to only remove the worst quadratic source terms stemming from the threshold resonances via a variable coefficient normal form based on the corresponding null structure recorded in Lemma~\ref{lem:null_structure1}. 
    Nonetheless, we record this deeper null structure in Appendix~\ref{appendix:deeper_null_structure} in the (algebraically simpler) case of odd perturbations of the static sine-Gordon kink. 

    \item In the perturbative proof of the asymptotic stability of the sine-Gordon kink under odd perturbations in \cite{LS1}, a specific super-symmetric factorization property of the linearized operator around the sine-Gordon kink is exploited to transform it to the flat linear Klein-Gordon operator via a Darboux transformation. The corresponding transformed evolution equation for the radiation term can then be studied taking a simpler space-time resonances approach based on the flat Fourier transform. 
    Following \cite[Sections 4.1-4.2]{KMMV20}, we could exploit similar factorization properties for the matrix linearized operator around the moving sine-Gordon kink to pass to a transformed matrix operator without a potential. However, for the corresponding transformed evolution equation for the radiation term the bookkeeping of the number of terms in the differentiated transformed nonlinearities becomes quite involved. Moreover, the second crucial null structure in the modulation equations appears significantly more difficult to capture taking such a transformed approach. An analogous comment in fact applies to the deeper null structure discussed in the preceding item (ii). For these reasons, taking a space-time resonances approach based on the distorted Fourier transform turned out more advantageous in the moving sine-Gordon kink setting.

    \item Null structures similar to those discussed in Subsection~\ref{subsec:overview_slow_local_decay} also occur in the asymptotic stability problem for solitary waves of the focusing cubic Schr\"odinger equation on the line, as uncovered in \cite{LL2, Li23}.

   \item We made use of the Wolfram Mathematica software system in the computation of the quadratic spectral distributions in Subsection~\ref{subsec:quadratic_spectral_distributions} and in the computation of the null structure in Lemma~\ref{lem:null_structure1}.
    
\end{enumerate}

\section{Preliminaries} \label{sec:preliminaries}

\subsection{Notation and conventions}

We collect several notational conventions and basic definitions that will be used throughout this work.

\medskip 

\noindent {\it Standard conventions.}
An absolute constant whose value may change from line to line is denoted by $C > 0$. 
In many estimates the implicit constants will depend on the initial Lorentz boost parameter $\ell_0 \in (-1,1)$ specified in the statement of Theorem~\ref{thm:main}. In this regard we remark that although the Lorentz boost modulation parameter is time-dependent, the dependence of the implicit constants in all estimates in this paper is always in terms of the initial boost parameter $\ell_0$ because of the comparison estimate \eqref{equ:comparison_estimate} and Lemma~\ref{lem:setting_up_basic_comparison}.
For non-negative $X,Y$ we write $X \lesssim Y$ if $X \leq C Y$, and we use the notation $X \ll Y$ to indicate that the implicit constant should be regarded as small.
Additionally, for non-negative $X$ and arbitrary $Y$, we use the short-hand notation $Y = \calO(X)$ if $|Y| \leq C X$.
Throughout, we use the Japanese bracket notation          
\begin{equation*}
 \jt := (1 + t^2)^{\frac12}, \quad \jx := (1+x^2)^{\frac12}, \quad \jxi := (1+\xi^2)^{\frac12},
\end{equation*}
and we denote by $\bfJ$ the matrix
\begin{equation} \label{equ:bfJ_definition}
    \bfJ := \begin{bmatrix} 0 & 1 \\ -1 & 0 \end{bmatrix}.
\end{equation}
We use standard notations for the Lebesgue spaces $L^p$ as well as for the Sobolev spaces $H^k$ and $W^{k,p}$.
Moreover, we set $\bm{L}^2(\bbR) := L^2(\bbR) \times L^2 (\bbR)$.
We caution the reader that in this work we write $D_x := -i\partial_x$ or $D_y := -i\partial_y$, while the notation $D := \ell \partial - i (\jxi+\ell\xi)$ is reserved for the derivative expression in the differentiated scalar transform $\calF_{\ell,D}^{\#}$ defined in \eqref{eq:modFD}, which is one of the building blocks of the distorted Fourier transform associated with the matrix linearized operator $\bfL_\ell$ around the moving kink.

\medskip
\noindent {\it Flat Fourier transform.}
Our conventions for the flat Fourier transform of a Schwartz function on the line are
\begin{equation} \label{equ:flat_FT_definition}
 \begin{aligned}
  \whatcalF[f](\xi) &= \hat{f}(\xi) = \frac{1}{\sqrt{2\pi}} \int_{\bbR} e^{-i x \xi} f(x) \, \ud x, \\
  \whatcalF^{-1}[f](x) &= \check{f}(x) = \frac{1}{\sqrt{2\pi}}  \int_{\bbR}  e^{ix \xi} f(\xi) \, \ud \xi.
 \end{aligned}
\end{equation}
The convolution laws for $g,h \in \calS(\bbR)$ then read
\begin{equation*}
\widehat{\calF}\bigl[g \ast h\bigr] = \sqrt{2\pi} \hatg \hath, \quad \widehat{\calF}\bigl[g h\bigr] = \frac{1}{\sqrt{2\pi}} \hatg \ast \hath.
\end{equation*}
Moreover, we recall that in the sense of tempered distributions
\begin{align}
    \widehat{\calF}[1](\xi) &= \sqrt{2\pi} \delta_0(\xi), \label{equ:prelim_FT_one} \\
    \widehat{\calF}[\tanh(\cdot)](\xi) &= -i \sqrt{\frac{\pi}{2}} \pvdots \cosech \Bigl(\frac{\pi}{2} \xi\Bigr). \label{equ:prelim_FT_tanh}
\end{align}
While \eqref{equ:prelim_FT_one} is standard, we refer to \cite[Lemma 5.6]{LS1} for a proof of \eqref{equ:prelim_FT_tanh}.

\medskip 
\noindent {\it Inner products.}
In terms of the $L^2$ inner product of complex-valued functions, we use
\begin{equation}\label{eq:L2inner}
    \langle f, g\rangle = \Re \int_\bbR  f\overline{g}\,\ud x.
\end{equation} 
Given two pairs of real-valued vector functions $\bmf = (f_1, f_2)$ and $\bmg = (g_1, g_2)$, their $\bm{L}^2$ inner product is given by
\begin{equation}\label{eq:L2L2inner}
    \langle \bmf, \bmg\rangle:=  \int_\bbR  \bigl( f_1 g_1 + f_2 g_2 \bigr) \, \ud x.
\end{equation}

\subsection{Scalar scattering theory}
In this subsection, we recall basics on the spectral and scattering theory for the linear operator
\begin{equation}\label{eq:LK}
   L := -\px^2 + W''(K(x)) = -\px^2 - 2 \sech^2(x) + 1
\end{equation}
obtained from linearizing around the static sine-Gordon kink. Its potential belongs to the family of reflectionless P\"oschl-Teller potentials. 
We start with the following standard spectral results for~$L$.

\begin{lemma} \label{lem:spectrumL} 
The linearized operator $L$ defined in \eqref{eq:LK} is a selfadjoint operator with domain $H^2(\bbR)$. Its spectrum is given by
\begin{equation}
    \sigma(L) = \sigma_d \cup \sigma_c := \{0\} \cup [1,\infty). 
\end{equation}
The eigenfunction for the zero eigenvalue is related to the translation invariance,
\begin{equation}
    LY=0,\qquad Y(x) := K'(x).
\end{equation}
Moreover, $L$ exhibits a threshold resonance at the edge of the continuous spectrum given by
\begin{equation}\label{eq:Lres}
    L \psi= \psi, \qquad \psi(x):=\tanh(x).
\end{equation}
Finally, $L$ satisfies the following coercivity properties:
    \begin{itemize}[leftmargin=1.8em]
        \item[(a)] There exists $\mu_0 \in (0,1)$ such that for all $v \in H^1(\bbR)$,
        \begin{equation} \label{eq:scalar_L_coercivity_with_Y}
            \langle Y, v \rangle = 0 \quad \Rightarrow \quad \langle Lv, v \rangle \geq \mu_0 \|v\|_{H^1}^2.
        \end{equation}

        \item[(b)] There exists $\mu_1 \in (0,1)$ such that for all $\ell \in (-1,1)$ and all $v \in H^1(\bbR)$,
        \begin{equation} 
            \langle Lv, v \rangle \geq \mu_1 \|v\|_{H^1}^2 - \frac{1}{\mu_1} \bigl| \langle Z_\ell, v \rangle \bigr|^2,
        \end{equation}
        where $Z_\ell(x) := (1+\ell^2) Y(x) + 2 \ell^2 x Y'(x)$ for $\ell \in (-1,1)$.
    \end{itemize}        
\end{lemma}
\begin{proof}
    For the proof of part (a) we first recall from \cite[Problem 39]{Flugge} that the spectrum of the linear operator $L$ with domain $H^2$ is $\{0\} \cup [1,\infty)$ and that $\mathrm{ker}(L) = \mathrm{span}(Y)$ with $Y := K'$. 
    Thus, we have 
    \begin{equation*}
        \begin{aligned}
            \langle L v, v \rangle \geq \|v\|_{L^2_x}^2 - \|Y\|_{L^2_x}^{-2} \bigl( \langle Y, v \rangle \bigr)^2.
        \end{aligned}
    \end{equation*}
    On the other hand, we obtain directly from the definition of $L$ that 
    \begin{equation*}
        \begin{aligned}
            \langle Lv, v \rangle = \|v\|_{H^1_x}^2 - \int_\bbR  2 \sech^2(x)  v(x)^2 \, \ud x.
        \end{aligned}
    \end{equation*}
    Hence, for small $0 < c < 1$ we find for any $v \in H^1_x$ with $\langle Y, v \rangle = 0$ that 
    \begin{equation*}
        \begin{aligned}
            \langle Lv, v\rangle \geq c \|v\|_{H^1}^2 - c \int_\bbR2 \sech^2(x) v(x)^2 \, \ud x + (1-c) \|v\|_{L^2}^2,
        \end{aligned}
    \end{equation*}
    which implies \eqref{eq:scalar_L_coercivity_with_Y}.
    The proof of part (b) follows along the lines of the proof of \cite[Lemma 2.6]{KMMV20}. 
\end{proof}

In the moving kink setting, we also need versions of the coercivity estimates from Lemma~\ref{lem:spectrumL} above for a rescaled version of the linearized operator.
\begin{corollary}\label{cor:coerLell}
For any fixed $\ell \in (-1,1)$, the linear operator
\begin{equation}
    L_\gamma:=-\gamma^{-2}\px^2-2\sech^2(\gamma x)+1, \quad \gamma = \frac{1}{\sqrt{1-\ell^2}},
\end{equation}
is a selfadjoint operator with domain $H^2(\bbR)$. Its spectrum is given by
\begin{equation}
    \sigma(L_\ell)=\sigma_d\cup \sigma_c:= \{0\} \cup \bigl[1,\infty\bigr). 
\end{equation}
The eigenfunction associated with the discrete eigenvalue $0$ is given by
\begin{equation}
    L_\gamma Y_\ell=0,\qquad Y_\ell(x) := \gamma K'(\gamma x).
\end{equation}
Moreover $ L_\gamma$ has a threshold resonance at the edge of the continuous spectrum given by
\begin{equation}
     L_\gamma \psi_\ell= \psi_\ell,\qquad \psi_\ell(x):=\tanh(\gamma x).
\end{equation}
Finally, $L_\gamma$ satisfies the following coercivity properties: there exists $\mu_0 \in (0,1)$ such that for all $v \in H^1(\bbR)$,
        \begin{equation} 
            \langle Y_\ell, v \rangle = 0 \quad \Rightarrow \quad \langle  L_\gamma v, v \rangle \geq \mu_0 \|v\|_{H^1}^2.
        \end{equation}        
\end{corollary}
\begin{proof}
All results follow directly from Lemma \ref{lem:spectrumL} by a change of variable $y = \gamma x$.
\end{proof}

Next we recall the basic scattering theory for the Schr\"odinger operator
\begin{equation} \label{eq:Hk}
    H:=-\px^2 - 2\sech^2(x).
\end{equation}
We refer to \cite{Yafaev, DT} for more background.
The Jost solutions for $H$ with the property that
 \begin{equation}
    H f_\pm(x,\zeta) = \zeta^2 f_\pm(x,\zeta), \quad e^{\mp i x \zeta} f_{\pm}(x, \zeta) \to 1 \quad \text{as} \quad x \to \pm \infty,
 \end{equation}
 are given by
 \begin{equation}
     f_\pm(x,\zeta) := c_\pm(\zeta) \calD^\ast \bigl( e^{\pm i x \zeta} \bigr), \quad \calD^\ast = -\partial + \tanh(\cdot), \quad c_\pm(\zeta) := \mp \bigl( i \zeta - 1 \bigr)^{-1},
 \end{equation}
 that is,
 \begin{equation}\label{eq:JostH}
    f_\pm(x,\zeta) := c_\pm(\zeta) \bigl( \mp i \zeta + \tanh(x) \bigr) e^{\pm i x \zeta} , \quad c_\pm(\zeta) := \mp \bigl( i \zeta - 1 \bigr)^{-1}.
 \end{equation}
 The latter identities follow from $-\px^2 ( e^{\pm i x \zeta} ) = \zeta^2 e^{\pm i x \zeta}$ and from the adjoint of the beautiful conjugation identity
 \begin{equation}
     \calD \bigl( -\px^2 - 2\sech^2(x) \bigr) = (-\px^2 ) \calD, \quad \calD := \partial + \tanh(\cdot).
 \end{equation}
The asymptotics of $f_{+}(x,\xi)$ as $x\to-\infty$ are expressed via the scattering data. Indeed, one has
\begin{align}\label{eq:scattrela}
    T(\xi) f_+(\cdot,\xi) & = f_-(\cdot,-\xi) + R_-(\xi) f_-(\cdot,\xi), \\
T(\xi) f_-(\cdot,\xi) & = f_+(\cdot,-\xi) + R_+(\xi) f_+ (\cdot,\xi),
\end{align}
where the transmission coefficient is given by
\begin{equation}\label{eq:TW}
 T(\xi)W\bigl(f_+(\cdot,\xi), f_-(\cdot,\xi)\bigr)={-2i\xi}   
\end{equation} 
with $W = W(\xi)$ being the Wronskian,  and where $R(\xi)$ is the reflection coefficient. Here we also recall that the Schr\"odinger operator $H$ defined in \eqref{eq:Hk} is reflectionless, i.e., $R(\xi)=0$, and that the transmission coefficient is given by
\begin{equation}\label{eq:Htrans}
 T(\xi) = \frac{\xi^2+2i\xi-1}{\xi^2+1} = \frac{\xi+i}{\xi-i}.
\end{equation}
By the standard spectral and scattering theory, the Jost solutions $f_{\pm}$ give rise to the integral kernel of the resolvent on the real axis (approached from the upper half plane), and by Stone's formula, therefore also to the 
spectral measure. More precisely, we have 
\begin{equation}\label{eq:Hresol}
    (H-(\xi^2+i0))^{-1} (x,y) = \frac{f_+(x,\xi) f_-(y,\xi)\one_{[x>y]}+ f_+(y,\xi) f_-(x,\xi)\one_{[y>x]}}{W\bigl(f_+(\cdot,\xi), f_-(\cdot,\xi)\bigr)}, 
\end{equation}
and the density of the spectral resolution $E(\ud \xi^2)$ of $H$ restricted to the continuous spectrum $[0,\infty)$ has the kernel
\begin{equation*}
    \frac{E(\ud \xi^2)}{\ud\xi}(x,y) = \frac{|T(\xi)|^2}{2\pi } \bigl(f_+(x,\xi)f_+(y,-\xi)+f_-(x,\xi)f_-(y,-\xi)\bigr).
\end{equation*}
The preceding formula implies that the distorted Fourier basis associated with $H$ takes the form
\begin{equation*}
 e(x,\xi) := \frac{1}{\sqrt{2 \pi}} \left\{ \begin{aligned}
 &T(\xi) \frac{i \xi - \tanh(x)}{i \xi - 1} e^{ix \xi} &\text{for } \xi \geq 0, \\
&T(-\xi) \frac{i \xi - \tanh(x)}{i \xi + 1} e^{ix\xi} &\text{for } \xi < 0,
\end{aligned} \right. 
\end{equation*}
where $T(\xi)$ denotes the transmission coefficient associated with $H$ introduced in \eqref{eq:Htrans}.

\section{Spectral Theory} \label{sec:spectral_theory}

In this section, we develop the spectral theory and the scattering theory for the non-selfadjoint matrix operator $\bfL_\ell$ obtained from linearizing around the moving kink. 
Moreover, we derive a representation formula for the evolution generated by this matrix operator in terms of the associated distorted Fourier transform.
We emphasize that all arguments in this section do not rely on specific features of the sine-Gordon model and extend to general moving kink settings.

For the reader's orientation, we note that in the nonlinear analysis in the remainder of this paper, the results from this section will be applied in a moving frame with the spatial coordinate $y := x - q(t)$, as introduced in \eqref{equ:setting_up_moving_frame_coordinate}. 
However, since we hope that the results from this section are of independent interest, we use the more conventional notation $x$ for the spatial coordinate here.

\subsection{Basic setting} 

Fix $\ell \in (-1,1)$. Throughout this section we drop the subscript $\ell$ from the notation $\bfL_\ell$ for the matrix linearized operator around the moving kink and we just write 
\begin{equation}\label{eq:bfL}
    {\bf L} = \begin{bmatrix} \ell \px & 1 \\ - L_\ell & \ell \px \end{bmatrix} 
\end{equation}
with 
\begin{equation} \label{eq:Lell}
    L_\ell := -\px^2+V_\ell(x)+1, \quad V_\ell(x):=-2\sech^2(\gamma x).
\end{equation}
We observe that 
\begin{equation}
    \bfL = \bfJ \bfH,
\end{equation}
where 
\begin{equation}
     \bfJ:=\begin{bmatrix}
	0& 1\\
	-1 & 0 
	\end{bmatrix},
    \quad
    \bfH:=\begin{bmatrix}
	L_\ell & -\ell\px\\
	\ell\px & 1 
	\end{bmatrix}.
\end{equation}
The equation for the kink in vectorial form, see \eqref{eq:vecKeq} in Subsection~\ref{subsec:modulation_and_orbital_stability} for the derivation, implies that
\begin{equation} \label{equ:Y0Y1_definition}
  \bmY_{0}(x) := \begin{bmatrix} -\gamma K'(\gamma x) \\ \gamma^2 \ell K''(\gamma x) \end{bmatrix}, \quad
  \bmY_{1}(x) := \begin{bmatrix} \gamma^3 \ell x K'(\gamma x) \\ - \gamma^3 K'(\gamma x) - \gamma^4 \ell^2 x K''(\gamma x) \end{bmatrix}
\end{equation}
satisfy
\begin{equation} \label{eq:Y0Y1}
    \bfL \bmY_{0} =0, \quad \bfL \bmY_{1}= \bmY_{0}.
\end{equation}

Integration by parts shows that $\bfH$ is symmetric with respect to the $\bm{L}^2$ inner product \eqref{eq:L2L2inner},
\begin{equation}
  \langle \bfH \bm{\phi},\bm{\psi}\rangle=\langle \bm{\phi}, \bfH 
 \bm{\psi}\rangle.
\end{equation}
Moreover, $\bfH$ is a selfadjoint operator and we have the following conclusion.
\begin{lemma} \label{lem:Hself}
$\bfH$ is a closed operator on the  domain $\Dom(\bfH)= H^2(\bbR)\times H^1(\bbR)$. On this domain $\bfH$ is symmetric, i.e., 
\begin{equation}
  \langle \bfH \bm{\phi},\bm{\psi}\rangle=\langle \bm{\phi}, \bfH \bm{\psi}\rangle.
\end{equation}
The complexification of $\bfH$ on $$\Dom_\bbC (\bfH):=\Dom(\bfH)+i\Dom(\bfH)\subset \bm{L}^2+i \bm{L}^2=:\bm{L}^2_\bbC$$ is selfadjoint on this domain.
\end{lemma}
\begin{proof}
It is clear that $\bfH$ is closed on its domain.
Integrating by parts, one can directly show that the operator $\bfH$ is symmetric with respect to the $\bm{L}^2(\bbR)$ inner product.

For the self-adjointness over the complex Hilbert space $\bm{L}^2_\bbC$, considering the inner product over the complex Hilbert space $\bm{L}_\bbC^2$, one has
\begin{equation}
  \langle \bfH\bm{\phi},\bm{\phi}\rangle_{\bm{L}^2_\bbC}:=\int_\bbR \overline{\phi}_1L_\ell\phi_1+|\phi_2|^2-\ell\overline{\phi_1}\px \phi_2+\overline{\phi_2}\ell\px\phi_1\,\ud x.
\end{equation} 
Note that
\begin{equation}
    \langle \bfH\bm{\phi},\bm{\phi}\rangle_{\bm{L}^2_\bbC} \in \bbR,
\end{equation}
which implies $\bfH=\bfH^*$. 
\end{proof}
We also have the following coercivity estimate for $\bfH$.
\begin{lemma}\label{lem:coerH}
    There exist absolute constants $C \geq 1$ and $\mu \in (0,1)$ such that for any  $\bmu = (u_1, u_2)  \in \calH$, the following holds.
    \begin{itemize}[leftmargin=1.8em]
        \item[(a)] Upper bound:
            \begin{equation} \label{eq:H_upper_bound}
                \bigl\langle \bfH \bmu, \bmu \bigr\rangle \leq 2 \|\bmu\|_{H^1\times L^2}^2.
            \end{equation}
        \item[(b)] Coercivity:
            \begin{equation} \label{eq:H_lower_bound}
                \bigl\langle \bfJ \bmY_1, \bmu \bigr\rangle = 0 \quad \Rightarrow \quad \bigl\langle \bfH \bmu, \bmu \bigr\rangle \geq \mu \gamma^{-2} \|\bmu\|_{H^1\times L^2}^2.
            \end{equation}          
    \end{itemize}
\end{lemma}
\begin{proof}
The bound \eqref{eq:H_upper_bound} follows by direct computation. 
For the proof of the lower bound \eqref{eq:H_lower_bound}, we first observe that
\begin{equation*}    
    \bigl\langle \bfH  \bmu, \bmu \bigr\rangle = \langle L_{\ell} u_1, u_1 \rangle + \|u_2\|_{L^2}^2 + 2\ell P[\bmu],
\end{equation*}
where 
\begin{equation*}
    P[\bmu] := \int_\bbR u_2 (\px u_1) \, \ud x.
\end{equation*}
Here it is convenient to pass to the rescaled variables
\begin{equation*}
    v_1(\gamma x) := u_1(x), \quad v_2(\gamma x) := u_2(x).
\end{equation*}
Then we have 
\begin{equation*}
    \|v_1\|_{L^2}^2 = \gamma \|u_1\|_{L^2}^2, \quad \|\px v_1\|_{L^2}^2 = \gamma^{-1} \|\px u_1\|_{L^2}^2, \quad \|v_2\|_{L^2}^2 = \gamma \|u_2\|_{L^2}^2,
\end{equation*}
and a direct computation gives 
\begin{equation*}
    \langle L_{\ell} u_1, u_1 \rangle + \|u_2\|_{L^2}^2 + 2\ell P[\bmu] = \gamma^{-1} \Bigl( \langle L v_1, v_1 \rangle + \| v_2 + \ell \gamma \px v_1 \|_{L^2}^2 \Bigr).        
\end{equation*}
Thus, we find that 
\begin{equation*}
    \gamma \bigl\langle \bfH \bmu, \bmu \bigr\rangle = \langle L v_1, v_1 \rangle + \| v_2 + \ell \gamma \px v_1 \|_{L^2}^2
\end{equation*}where $L$ is given by \eqref{eq:LK}.

Next, we compute that 
\begin{equation*}
    \bigl\langle \bfJ \partial_\ell \bmK_{\ell,q}, \bmu \bigr\rangle = 0 \quad \Leftrightarrow \quad \langle Z_\ell, v_1 \rangle + \ell \gamma^{-1} \langle x Y, v_2 + \ell \gamma \px v_1 \rangle = 0
\end{equation*}
with $Y(x) := K'(x)$ and $Z_\ell(x) := (1+\ell^2) Y(x) + 2 \ell^2 x Y'(x)$ for $\ell \in (-1,1)$, whence
\begin{equation*}
    \bigl|\langle Z_\ell, v_1 \rangle\bigr| \leq C_1 \|v_2 + \ell \gamma \px v_1 \|_{L^2}^2
\end{equation*}
for some absolute constant $C_1 \geq 1$. 
Hence, by Lemma~\ref{lem:spectrumL}(b) we obtain that
\begin{equation*}
    \begin{aligned}
        \gamma \bigl\langle \bfH \bmu, \bmu \bigr\rangle = \langle L v_1, v_1 \rangle + \bigl\| v_2 + \ell \gamma \px v_1 \bigr\|_{L^2}^2 &\geq \frac{\mu_1}{2C_1^2} \langle L v_1, v_1 \rangle + \bigl\| v_2 + \ell \gamma \px v_1 \bigr\|_{L^2}^2 \\
        &\geq \frac{\mu_1^2}{2C_1^2} \|v_1\|_{H^1_x}^2 + \frac12 \bigl\| v_2 + \ell \gamma \px v_1 \bigr\|_{L^2}^2. 
    \end{aligned}
\end{equation*}
Finally, using that uniformly for all $\ell \in (-1,1)$, we have
\begin{equation*}
    \|\bmu\|_{H^1 \times L^2}^2 \leq 3 \gamma \Bigl( \|v_1\|_{H^1}^2 + \| v_2 + \ell \gamma \px v_1 \|_{L^2}^2 \Bigr),
\end{equation*}
the asserted lower bound \eqref{eq:H_lower_bound} follows.
\end{proof}

\subsection{Free operator} 

Before we turn to the spectrum of the operator $\bfL$, we first consider the corresponding free operator
\begin{equation} \label{equ:free_matrix_operator}
    {\bf L}_0 := \begin{bmatrix} \ell \px & 1 \\ - \bigl( - \px^2 + 1 \bigr) & \ell \px \end{bmatrix}, \quad \ell \in (-1,1).
\end{equation}
Consider the real Hilbert space $$\calH:=H^1(\bbR)\times L^2(\bbR).$$ Setting ${\bf L}_0= \bfJ {\bf A}_0+{\bf B}_0$, where 
 \begin{equation}\label{eq:A0B0}
 {\bf A}_0:=\begin{bmatrix}
	 -\partial_x^2+1& 0\\
	0 & 1 
	\end{bmatrix},\qquad 
{{\bf B}_0}:=\begin{bmatrix}
	 \ell\partial_x & 0\\
	0 &  \ell\partial_x 
	\end{bmatrix},
\end{equation}
the inner product on $\calH$ over $\bbR$ is given by
\begin{equation}\label{eq:innerprodH1}
    \langle \bm{\phi},\bm{\psi}\rangle_{\calH}:= \langle {\bf A}_0 \bm{\phi},\bm{\psi}\rangle.
\end{equation}
The norm induced by this inner product is denoted by $\|\bm{\phi}\|_{\calH}$.

\begin{lemma}\label{lem:calLsa}
    In the Hilbert space $\calH$ with the inner product given by \eqref{eq:innerprodH1},  $ {\bf L}_0$ is a closed operator on the  domain $\Dom( {\bf L}_0)=H^2(\bbR)\times H^1(\bbR)$. On this domain $ {\bf L}_0$ is anti-symmetric, i.e., 
     \begin{equation} \label{equ:free_operator_vanishing}
     \langle  {\bf L}_0 \bm{\phi},\bm{\phi}\rangle_{\calH} = 0\quad \forall \bm{\phi}\in \Dom( {\bf L}_0),
     \end{equation}
     which implies for all $\lambda\in\bbR$ that
     \begin{equation} \label{eq:LowerR}
         \|({\bf L}_0-\lambda)\bm{\phi}\|_{\calH}\ge |\lambda|\|\bm{\phi}\|_{\calH}.
     \end{equation}
    The complexification of ${\bf L}_0$ on $$\Dom_\bbC({\bf L}_0):=\Dom({\bf L}_0)+i\Dom({\bf L}_0)\subset \calH+i\calH=:\calH_\bbC$$ is anti-selfadjoint on this domain, i.e., $i{\bf L}_0$ is selfadjoint on $\Dom_\bbC({\bf L}_0)$. 
\end{lemma}
\begin{proof}
    The closedness is clear. To prove the vanishing \eqref{equ:free_operator_vanishing}, we note that for all $\bm{\phi}\in\Dom({\bf L}_0)$,
    \[
        \langle {\bf L}_0 \bm{\phi},\bm{\phi}\rangle_{\calH} = \langle (\bfJ {\bf A}_0+{\bf B}_0)\bm{\phi},{\bf A}_0\bm{\phi}\rangle = 0.
    \]
    The lower bound \eqref{eq:LowerR} then follows from the Cauchy-Schwarz inequality. To show the anti-selfadjointness of the complexification of $\bfL_0$, we work on $\calH_\bbC$. We compute the adjoint ${\bf L}_0^*$ whose domain consists of all $\widetilde{\bm{\phi}}\in \calH_\bbC$ so that
    there exists $\bm{\phi}\in\calH_\bbC$ with 
    \[
    \langle \bm{\phi}, \bm{\psi}\rangle_{\calH_\bbC}:= \langle {\bf A}_0 \bm{\phi},\bm{\psi}\rangle_{\bm{L}^2_\bbC}= \langle \widetilde{\bm{\phi}}, {\bf L}_0\bm{\psi}\rangle_{\calH_\bbC}
    \]
    for all $\bm{\psi}\in\Dom_\bbC({\bf L}_0)$ and the inner product is the one in~$\calH_\bbC$. Equivalently, again over $\bbC$, 
    \begin{equation} \label{eq:dom*}
        \langle {\bf A}_0\bm{\phi}, \bm{\psi}\rangle_{\bm{L}^2_\bbC} = \langle {\bf A}_0\widetilde{\bm{\phi}}, {\bf L}_0\bm{\psi}\rangle_{\bm{L}^2_\bbC}. 
    \end{equation}
    Now suppose $\bm{g}=\binom{g_1}{g_2}\in \calH_\bbC$. Then we solve 
    \[
    {\bf L}_0\bm{\psi} = \bm{g} \quad \Leftrightarrow \quad \binom{\ell\px \psi_2}{(\partial_x^2-1)\psi_1+\ell \partial_x \psi_2} = \binom{g_1}{g_2}
    \]
    uniquely for $\bm{\psi}\in \Dom_\bbC({\bf L}_0)$. Thus, \eqref{eq:dom*} is equivalent to  
    \[
    \langle {\bf A}_0\bm{\phi},  {\bf L}_0^{-1}\bm{g} \rangle_{\bm{L}^2_\bbC} = \langle {\bf A}_0\widetilde{\bm{\phi}}, \bm{g}\rangle_{\bm{L}^2_\bbC}. 
    \]
    This shows that $\widetilde{\bm{\phi}}={\bf A}_0^{-1}({\bf L}_0^{-1})^*{\bf A}_0\bm{\phi}$ with the adjoint in $\bm{L}^2_\bbC$,  in the sense of distributions. Now ${\bf A}_0\bm{\phi}\in H^{-1}\times L^2$, and $({\bf L}_0^{-1})^*{\bf A}_0\bm{\phi}\in L^2\times H^{1}$. Thus, finally, $ \widetilde{\bm{\phi}} \in H^2\times H^1$ whence  $\Dom_\bbC({\bf L}_0^*)=\Dom_\bbC({\bf L}_0)$ as claimed. 
\end{proof}

Next, we determine the spectrum of $\bfL_0$ explicitly.
\begin{lemma}\label{lem:spectrumL0}
The spectrum of ${\bf L}_0$ defined in \eqref{equ:free_matrix_operator} is given by 
\begin{equation*}
 \sigma({\bf L}_0) = \bigl(-i\infty,-i\gamma^{-1}\bigr] \cup \bigl[ i\gamma^{-1},i\infty\bigr), \quad \gamma := \frac{1}{\sqrt{1-\ell^2}}.
\end{equation*} 
\end{lemma}
\begin{proof}
First, it follows from Lemma \ref{lem:calLsa} that $\sigma({\bf L}_0)\in i\bbR$. To find the spectrum $\sigma({\bf L}_0)$, we solve 
\begin{equation*}
 {\bf L}_0 \begin{bmatrix} f \\ g \end{bmatrix} = z \begin{bmatrix} f \\ g \end{bmatrix} 
\end{equation*}
over $\bbC$ in  the form $\ell f'+g=zf$ and $(1-\ell^2)f''-(1+z^2)f+2\ell\gamma z f'=0$. The latter has a fundamental system $e^{k_\pm(z) x}$ with 
\[
    k_{\pm}(z) := \gamma \bigl( -\gamma \ell z \pm \sqrt{1+\gamma^2 z^2} \bigr).
\]
The spectrum is characterized by
\[
    \sigma({\bf L}_0) = \bigl\{ z\in\bbC\:|\: k_\pm(z)\in i\bbR \bigr\}.
\]
Solving for $k_\pm(z)=i\tau$, $\tau\in\bbR$, yields 
\begin{equation} \label{eq:zspec}
    z= i \bigl( \ell\tau\pm\sqrt{1+ \tau^2} \bigr).
\end{equation}
Analyzing the graphs of the functions on the right-hand side over $\tau\in\bbR$ provides the desired description of the spectrum. 
\end{proof}

\subsection{Spectrum of the linearized operator}

Next, we determine the spectrum of the operator~$\bfL$.

\begin{lemma} \label{lem:Lspec}
 The essential spectrum of $\bfL$ is $(-i\infty,-i\gamma^{-1}]\cup[i\gamma^{-1},i\infty)$ without embedded eigenvalues. The discrete spectrum of $\bfL$ contains only $0$.
 The geometric multiplicity of the zero eigenvalue is one and its algebraic multiplicity is two.
 It holds that $\ker \, (\bfL^n)= \ker \, (\bfL^2)$ for all $n \geq 2$, and $\ker \, (\bfL^2) = \mathrm{span} \, \{\bmY_{0},\bmY_{1}\}$.
\end{lemma}
\begin{proof}
The assertion about the essential spectrum of $\bfL$ follows from Weyl's theorem and Lemma~\ref{lem:spectrumL0}, because $\bfL$ is a compact perturbation of $\bfL_0$. Now we turn to the discrete spectrum.  Assume for some $\lambda \in \bbC$ that there exists $\bm{\phi} = (\phi_1, \phi_2) \in \bm{L}_{\bbC}^2$ such that
\begin{equation*}
    \bfL \bm{\phi}=\lambda \bm{\phi}.
\end{equation*}
Then we have 
\begin{align*}
    \ell\px \phi_1+\phi_2 &= \lambda \phi_1, \\
   - \bigl( - \px^2 - 2 \sech^2(\gamma x) + 1 \bigr) \phi_1+\ell\px\phi_2 &= \lambda \phi_2.
\end{align*}
Combining the two preceding equations gives 
\begin{equation}
 - \bigl( - \px^2 - 2 \sech^2(\gamma x) + 1 \bigr) \phi_1+\ell\px(\lambda\phi_1-\ell\px \phi_1) = \lambda (\lambda\phi_1-\ell\px \phi_1),
\end{equation}
which simplifies to 
\begin{equation}
    - \bigl( -(1-\ell^2) \px^2 - 2 \sech^2(\gamma x) + 1 \bigr) \phi_1 + 2\lambda\ell \px \phi_1 = \lambda^2 \phi_1.
\end{equation}
Setting $\rho(y) := e^{\ell \gamma \lambda y} \phi_1(\gamma^{-1} y)$, the equation for $\rho(y)$ is
\begin{equation}
     \bigl( - \py^2 - 2 \sech^2(y) + 1  \bigr) \rho(y)= -\lambda^2 \gamma ^2\rho(y).
\end{equation}
From the spectral information for the linear operator  $- \py^2 - 2 \sech^2(y) + 1  $ summarized in Lemma \ref{lem:spectrumL},  we conclude that $\lambda = 0$ and that $\rho(y)$ is a multiple of $Y(y) := K'(y)$.

We now analyze the generalized kernel of $\bfL$. 
Recall from \eqref{eq:Y0Y1} that
\begin{equation}
\bfL \bmY_{0} =0, \qquad\bfL  \bmY_{1}= \bmY_{0}.
\end{equation}We claim that there is no solution $\bm{\phi}$ such that
\begin{equation*}
    \bfL \bm{\phi} = \bmY_{1}.
\end{equation*}
Writing out both rows of the equation above explicitly, we get
\begin{equation*}
    \phi_2=\ell\gamma^3yK'(\gamma x)-\ell \px \phi_1
\end{equation*}
and
\begin{equation*}
    - \bigl( - \px^2 - 2 \sech^2(\gamma x) + 1 \bigr) \phi_1 = - \gamma^3 K'(\gamma x) - \ell^2 \gamma^4 x K''(\gamma x)-\ell\px \phi_2.
\end{equation*}
Putting the two equations together, we obtain
\begin{equation}
\bigl( -(1-\ell^2)\px^2 - 2 \sech^2(\gamma x) + 1 \bigr) \phi_1=2 \ell^2 \gamma^4 x K''(\gamma x)+\gamma^3(1+\ell^2)K'(\gamma x).
\end{equation}
Setting $\kappa := \frac{1}{2} \ell^2 x^2 \gamma^5 K'(\gamma x) + \phi_1$
and using the equation for $K'(\gamma x)$, we find
\begin{equation}
    \bigl( -(1-\ell^2)\px^2 - 2 \sech^2(\gamma x) + 1 \bigr) \kappa = -\gamma^3 K'(\gamma x).
\end{equation}
It follows that $\kappa$ is in the generalized kernel of the linear operator on the left-hand side since the right-hand side of the equation is in the kernel of the operator on the left-hand side. This implies that $\kappa$ does not exist since $ -(1-\ell^2)\px^2 - 2 \sech^2(\gamma x) + 1$ is selfadjoint. 
Thus, we must have $\ker \, (\bfL^n) = \ker \, (\bfL^2)$ for $n\geq 2$, which finishes the proof.
\end{proof}

Moreover, we observe that the threshold resonance \eqref{eq:Lres} of the scalar operator $L$ induces two threshold resonances $\bmPhi_{\pm}(x)$ for the matrix operator $\bfL$ at the two edges $\pm i\gamma^{-1}$ of the essential spectrum
\begin{equation}
    \bfL \bmPhi_\pm (x) = \pm i\gamma^{-1} \bmPhi_\pm
\end{equation}
with 
\begin{equation}
    \bmPhi_\pm(x) = \binom{\tanh(\gamma x)}{\pm i \gamma \tanh(\gamma x)-\ell\gamma^2\sech^2(\gamma x)} e^{\mp i\gamma \ell x}.
\end{equation}

Let $P_{\mathrm{d}}$ be the Riesz projection onto the discrete spectrum of $\bfL$ defined by
\begin{equation}\label{eq:rieszproj}
  P_\mathrm{d} := \frac{1}{2\pi i} \oint_\Gamma (\bfL - \lambda)^{-1} \, \ud \lambda,
\end{equation}
where $\Gamma$ is a simple closed curve within the resolvent set that encloses the zero eigenvalue.
We define the projection onto the essential spectrum of $\bfL$ by
\begin{equation}\label{eq:defPe}
    P_\mathrm{e}:=1-P_\mathrm{d}.
\end{equation}
Next, we show that the range of the projection $P_{\mathrm{d}}$ is the generalized kernel of $\bfL$.

\begin{lemma} \label{lem:ranPd} 
    We have $$\mathrm{ran} (P_\mathrm{d})=\calN:=\ker(\bfL^2).$$  
\end{lemma}
\begin{proof}
By construction, it is clear that $\ker(\bfL^n)\subset\mathrm{ran} (P_\mathrm{d})$ for all $n\geq1$. 

Now suppose $z \notin \sigma(\bfL_0)$. Since $(\bfL_0-z)$ is invertible, we write
\begin{equation}
    \bfL - z = \bigl( 1+{\bf V}(\bfL_0-z)^{-1} \bigr)(\bfL_0-z), \quad {\bf V} := \begin{bmatrix}0 & 0 \\  2 \sech ^2(\gamma x)  & 0 \end{bmatrix}.
\end{equation}
Note that ${\bf V}(\bfL_0-z)^{-1}$ is analytic for $z \notin \sigma_(\bfL_0)$ and compact. The analytic Fredholm theorem implies that $(\bfL-z)$ is a meromorphic function and that $(1+{\bf V}(\bfL_0-z)^{-1})$
is invertible for all but a discrete set of $\bbC \setminus \sigma({\bf L}_0)$.  Moreover, we also know that this discrete set of poles is the  discrete spectrum of $\bfL$, which is $\{0\}$ by Lemma \ref{lem:Lspec}. Therefore, one has
\begin{equation*}
    \|(\bfL-z)^{-1}\|_{\calH\rightarrow\calH}\lesssim |z|^{-m}
\end{equation*} for some $m\in \mathbb{N}$.  It follows that $\bfL^m P_\mathrm{d}=0$ and $\mathrm{ran} (P_\mathrm{d})\subset \ker (\bfL^m)$. By Lemma \ref{lem:Lspec}, we conclude that $\mathrm{ran} (P_\mathrm{d})= \ker (\bfL^2)$.
\end{proof}

We carry out a similar analysis for the adjoint operator $\bfL^{*}$ (with respect to the $\bm{L}^2(\bbR)$ inner product \eqref{eq:L2L2inner}). Since $\bfL = \bfJ \bfH$, Lemma~\ref{lem:Hself} gives $\bfL^{*}=-\bfH\bfJ$, and \eqref{eq:Y0Y1} implies
\begin{equation}
    \bfH  \bmY_{0}=0,\quad  \bfH  \bmY_{1} = -\bfJ \bmY_{0}.
\end{equation}  
Analogously to the proof of Lemma~\ref{lem:Lspec}, we determine the spectrum of $\bfL^{*}$.
\begin{lemma} \label{lem:L*spec}
    The essential spectrum of $\bfL^*$ is $(-i\infty, -i\gamma^{-1}] \cup [i\gamma^{-1},i\infty)$ without embedded eigenvalues and the discrete spectrum of $\bfL^*$ contains only $0$. The zero eigenvalue is of dimension $2$, i.e., $\ker ((\bfL^*)   ^n)= \ker ((\bfL^*) ^2)$ for all $n\geq2$. Moreover, $\ker \, ((\bfL^*) ^2) = \mathrm{span} \, \{\bfJ \bmY_{0}, \bfJ \bmY_{1}\}$.
\end{lemma}
One can also analogously define the adjoint Riesz projection onto the discrete spectrum of $\bfL^{*}$,
\begin{equation}\label{eq:rieszproj*}
  P_\mathrm{d}^*:= \frac{1}{2\pi i} \oint_\Gamma (\bfL^* - \lambda)^{-1} \, \ud \lambda,
\end{equation}
where $\Gamma$ is a simple closed curve within the resolvent set that encloses the zero eigenvalue. 
Analogously to the proof of Lemma~\ref{lem:ranPd}, we determine the range of the projection $P_{\mathrm{d}}^{*}$.

\begin{lemma} \label{lem:ranPd*}
We have $\mathrm{ran} (P_\mathrm{d}^*)=\calN^*:=\ker((\bfL^*)^2)=\mathrm{span}\{\bfJ\bmY_{0},\bfJ\bmY_{1}\}$.
\end{lemma}

Next, we derive a spectral decomposition for $\bm{L}^2(\bbR)$.

\begin{lemma}\label{lem:decom}
There is a direct sum decomposition
\begin{equation} \label{eq:decom}
      \bm{L}^2(\bbR)=\calN + (\calN^*)^{\perp},
\end{equation}
which is invariant under $\bfL$, whence, it is invariant under the flow $e^{t\bfL}$. Let $P_\mathrm{e}$ be the projection onto $(\calN^*)^{\perp}$, which is induced by the splitting \eqref{eq:decom}. One has $P_\mathrm{d}=1-P_\mathrm{e}$ and
\begin{equation}\label{eq:Pdexp}
    P_\mathrm{d} \bm{\phi} = \frac{\langle \bm{\phi}, \bfJ\bmY_{1}\rangle}{\langle \bmY_{0}, \bfJ\bmY_{1}\rangle} \bmY_{0} - \frac{\langle \bm{\phi}, \bfJ\bmY_{0}\rangle}{\langle \bmY_{0}, \bfJ\bmY_{1}\rangle} \bmY_{1}.
\end{equation}
Explicitly, we have for any $\bmphi \in \bm{L}^2(\bbR)$ the decomposition
\begin{equation}
    \bmphi = P_{\mathrm{e}} \bmphi + d_0 \bmY_0 + d_1 \bmY_1
\end{equation}
with 
\begin{equation}
    d_0 := \frac{\langle \bm{\phi}, \bfJ\bmY_{1}\rangle}{\langle \bmY_{0}, \bfJ\bmY_{1}\rangle}, \quad d_1 := - \frac{\langle \bm{\phi}, \bfJ\bmY_{0}\rangle}{\langle \bmY_{0}, \bfJ\bmY_{1}\rangle}.
\end{equation}
It follows that the projection $P_{\mathrm{e}}$ is bounded as an operator $H^k \times H^{k-1} \to H^k \times H^{k-1}$ for any integer~$k \geq 1$. 
\end{lemma}
\begin{proof}
By definition of the projection $\Pe$ in \eqref{eq:defPe}, we can write
\begin{equation}
    \bm{L}^2(\bbR)=\ker (P_\mathrm{e})+ \mathrm{ran} (P_\mathrm{e})=\ker (P_\mathrm{e})+ \ker(P^*_\mathrm{e})^{\perp}.
\end{equation}
Using that by Lemma \ref{lem:ranPd} and Lemma \ref{lem:ranPd*}, 
\begin{equation}
    \ker (P_\mathrm{e})= \mathrm{ran} (P_\mathrm{d})=\mathcal{N},\quad   \ker (P^*_\mathrm{e})= \mathrm{ran} (P_\mathrm{d}^*)=\mathcal{N}^*,
\end{equation}
the decomposition \eqref{eq:decom} follows. Clearly, one has $\bfL \calN\subset \calN$ and $\bfL (\calN^*)^{\perp}\subset (\calN^*)^{\perp} $. It follows that the splitting \eqref{eq:decom} is invariant under $\bfL$.
Next, for any $\bmphi\in \bm{L}^2(\bbR)$, by \eqref{eq:decom} and the argument above, we have
$\bmphi=P_\mathrm{e} \bmphi+P_\mathrm{d} \bmphi$
with $P_\mathrm{e} \bmphi \in(\calN^*)^{\perp}$. In view of Lemma~\ref{lem:Lspec} and Lemma~\ref{lem:ranPd}, $P_\mathrm{d} \bmphi$ can be further written as
\begin{equation*}
    P_\mathrm{d} \bmphi= d_0 \bmY_0 + d_1 \bmY_1.
\end{equation*}
Since $\calN^*=\mathrm{span}\{\bfJ\bmY_{0},\bfJ\bmY_{1}\}$ by Lemma~\ref{lem:ranPd*}, it follows that
\begin{equation*}
   \langle \bm{\phi}, \bfJ\bmY_{0}\rangle =  \langle P_\mathrm{e}\bm{\phi}, \bfJ\bmY_{0}\rangle + d_1 \langle \bmY_1, \bfJ\bmY_{0}\rangle. 
\end{equation*}
The first term on the right-hand side vanishes, whence
\begin{equation*}
    d_1 = \frac{ \langle \bm{\phi}, \bfJ\bmY_{0}\rangle}{\langle \bmY_1, \bfJ\bmY_{0}\rangle } = -\frac{ \langle \bm{\phi}, \bfJ\bmY_{0}\rangle}{\langle \bmY_0, \bfJ\bmY_{1}\rangle }.
\end{equation*}
Similarly, one has $d_0 = \frac{ \langle \bm{\phi}, \bfJ\bmY_{1}\rangle}{\langle \bmY_0, \bfJ\bmY_{1}\rangle }$. 
Thus, \eqref{eq:Pdexp} follows and the proof is finished.
\end{proof}

\begin{remark}\label{rem:Pdzero}
Since $\bmY_0$ and $\bmY_1$ are linearly independent, we have that $P_\mathrm{d}\bm\phi=0$  is equivalent to $\ \langle \bm{\phi}, \bfJ\bmY_{1}\rangle =  \langle \bm{\phi}, \bfJ\bmY_{0}\rangle =0.$
\end{remark}

\begin{remark}\label{rem:Pdadj}
Note that by \eqref{eq:Pdexp}, one can conclude that
\begin{equation}
    \bfJ \circ P_\mathrm{d} \circ \bfJ =P_\mathrm{d}^*.
\end{equation}
Moreover, we note that the dual space of $P_\mathrm{d} (\bm{L}^2(\bbR))$ with respect to the $\bm{L}^2(\bbR)$ inner product is $P_\mathrm{d}^*(\bm{L}^2(\bbR))$.  Also note that $\bfJ  (P_\mathrm{e} (\bm{L}^2(\bbR)))=P^*_\mathrm{e} (\bm{L}^2(\bbR)) $. In general, if $P_\mathrm{d} \bmf=0$, then $P_\mathrm{d}^*  (\bfJ \bmf)=0$.
\end{remark}

\begin{remark}
Recall that the the symplectic form associated with the wave equation is given by $\omega(\bm{\phi}, \bm{\psi} ) := \langle \bm{\phi}, \bfJ\bm{\psi}\rangle$. 
Thus, $\langle \bm{\phi}, \bfJ\bmY_{1}\rangle = \langle \bm{\phi}, \bfJ\bmY_{0} \rangle = 0$ means that $\bm{\phi}$ is symplectically orthogonal to $\bmY_{0}$ and $\bmY_{1}$. 
\end{remark}

\subsection{Resolvents and limiting absorption principle}

In this subsection we compute the resolvent of the matrix operator $\bfL$. Additionally, we establish the limiting absorption principle and further mapping properties for the resolvent of $\bfL$.
In general, these properties can be inferred from the corresponding properties of the free matrix resolvent using resolvent identities. We therefore first present the corresponding results for the free matrix resolvent.
Afterwards, we present a streamlined proof for the sine-Gordon case based on the explicit formulas for the integral kernel of the resolvent of the matrix linearized operator around the moving sine-Gordon kink.

In what follows, we use the notation 
\begin{equation} \label{eq:Sigmaell}
    \Sigma_\ell := \bigl( -\infty, -\gamma^{-1} \bigr] \cup \bigl[ \gamma^{-1}, \infty \bigr), \quad \ell \in (-1,1), \quad \gamma := \frac{1}{\sqrt{1-\ell^2}}.
\end{equation}

\subsubsection{Free resolvents}

We start with the computations for the free matrix operator $\bfL_0$. 
To this end we first recall the integral kernel of the resolvent of the one-dimensional Laplacian,
\begin{equation} \label{eq:1dimresolv}
    \frakR_0(\zeta)(x,z) := (-\partial^2 -\zeta^2)^{-1}(x,y) = - \frac{e^{i\zeta|x-y|}}{2i\zeta}, \quad \zeta \in \bbC_+.
\end{equation}
In the next lemma we determine the resolvent for the free matrix operator $\bfL_0$ and we compute its integral kernel.

\begin{lemma} \label{lem:freeresol}
    For all $\lambda\in\bbC\setminus \Sigma_\ell$,  one has 
    \begin{align} 
    \big( \bfL_0-i\lambda\mathbb{I}_{2\times2}\big)^{-1}&= \mathscr{R}_\ell(\lambda) \begin{bmatrix}
	\ell\partial -i\lambda & -1\\
	-\partial^2+1 & 
    \ell\partial -i\lambda
	\end{bmatrix}\\
    &=\begin{bmatrix}
	\mathscr{R}_\ell(\lambda) (\ell\partial -i\lambda) & -\mathscr{R}_\ell(\lambda) \\
	1- (\ell\partial -i\lambda)\mathscr{R}_\ell(\lambda) (\ell\partial -i\lambda)& 
    (\ell\partial -i\lambda) \mathscr{R}_\ell(\lambda) 
	\end{bmatrix} \label{eq:Rell}
\end{align}
where 
\begin{equation} \label{equ:mathscrRell_definition}
    \begin{aligned}
        \mathscr{R}_\ell(\lambda)(x,y) &:= \bigl( - (1-\ell^2) \partial^2 - 2i\ell\lambda\partial + 1 - \lambda^2 \bigr)^{-1}(x,y) \\
        &= e^{-i\gamma^2\ell\lambda x}\left(-(1-\ell^2)\partial^2+1-\gamma^2\lambda^2\right)^{-1}(x,y)e^{i\gamma^2\ell\lambda y} \\ 
        &= - \frac{e^{i(\gamma\zeta|x-y|+\ell\gamma^2\lambda(y-x))}}{2i\zeta}
    \end{aligned}
\end{equation}
with $\zeta^2 = \gamma^2 \lambda^2-1$ for $\Im \, \zeta>0$.  Finally, define $\mathscr{R}_\ell^{\pm}(\lambda) := \mathscr{R}_\ell(\lambda\pm i0)$. Then for $\lambda\gamma>1$ we have
\begin{equation} \label{eq:res_lim}
\mathscr{R}_\ell^\pm(\lambda)(x,y)= \mp \frac{e^{i\gamma(\pm \sqrt{\lambda^2\gamma^2-1}|x-y|+\ell\gamma\lambda(y-x))}}{2i\sqrt{\lambda^2\gamma^2-1}},
\end{equation}
were $x$ is the output variable, and $y$ is the input variable.
\end{lemma}
\begin{proof}
By the above spectral analysis over $\bbC$,   if $\lambda\in\bbC \setminus \Sigma_\ell$, then the resolvent $\big( \mathcal{L}-i\lambda\mathbb{I}_{2\times2}\big)^{-1}$ is bounded from $H^1(\bbR)\times L^2(\bbR)$ to $H^2(\bbR)\times H^1(\bbR)$.
   The symbol of the matrix operator
\begin{equation}
\bfL_0-i\lambda\mathbb{I}_{2\times2}=\begin{bmatrix}
	 \ell\partial-i\lambda & 1\\
	\partial^2-1 & 
 \ell\partial-i\lambda
	\end{bmatrix}
\end{equation}
equals
\begin{equation}
   \left(\bfL_0-i\lambda\mathbb{I}_{2\times2}\right)\widehat{\ }(\xi)=\begin{bmatrix}
	 i\xi \ell-i\lambda & 1\\
	-\xi^2-1 & 
 i\xi \ell-i\lambda
	\end{bmatrix}.
\end{equation}
Inverting the matrix yields 
\begin{equation} \label{eq:inversematrix}
   \begin{bmatrix}
	   i\xi \ell-i\lambda & 1 \\
	   -\xi^2-1 & i\xi \ell-i\lambda
	\end{bmatrix}^{-1}
    =
    \bigl( (1-\ell^2) \xi^2 + 2\ell\lambda\xi + 1 - \lambda^2 \bigr)^{-1} 
    \begin{bmatrix}
	   i\xi \ell-i\lambda  & -1 \\
	   \xi^2+1 & i\xi \ell-i\lambda
	\end{bmatrix}.
\end{equation}
Note that $(1-\ell^2) \xi^2 + 2\ell\lambda\xi + 1 - \lambda^2 = 0$ is equivalent to $\lambda =\xi\ell\pm\sqrt{1+\xi^2}$, which is exactly the same as~\eqref{eq:zspec} for $z=i\lambda$ and~$\xi = \sigma \in\bbR$. By our choice of $\lambda$, it follows that the matrix is invertible. 
We remove the first order derivative by conjugation by a complex phase. More precisely,  
we use that
\begin{align}
    e^{-i\gamma^2\ell\lambda x}\left(-(1-\ell^2)\partial^2+1-\gamma^2\lambda^2\right)e^{i\gamma^2\ell\lambda x} 
    &= - (1-\ell^2) \partial^2 - 2i \ell \lambda \partial + 1 - \lambda^2.
\end{align}
Define the multiplication operator $(M_a f)(x)=e^{ia x}f(x)$ and 
\begin{equation}
   \mathscr{R}_\ell(\lambda) :=  M_{\gamma^2\ell\lambda}^{-1}\left(-(1-\ell^2)\partial^2+1-\gamma^2\lambda^2\right)^{-1}M_{\gamma^2\ell\lambda} = \bigl( -(1-\ell^2)\partial^2 - 2i \ell\lambda\partial + 1-\lambda^2 \bigr)^{-1}.
\end{equation}
Taking the inverse Fourier transform in \eqref{eq:inversematrix} leads to~\eqref{eq:Rell}. To prove~\eqref{eq:res_lim}, we compute
\[
\lim_{\epsilon\to0+} \sqrt{(\lambda\pm i\epsilon)^2\gamma^2-1} = \pm \sqrt{\lambda^2\gamma^2-1}, 
\]
using that $\lambda\gamma>1$ and that the branch of the square root on the left-hand side takes its values in the upper half-plane. 
\end{proof}

For $\lambda\in\bbC\setminus\Sigma_\ell$, we define
\begin{equation}\label{eq:freesol}
   \bfR_0(\lambda):=(\bfL_0-i\lambda\mathbb{I}_{2\times2})^{-1},
\end{equation}
and we denote $\bfR_0^{\pm}(\lambda) := \bfR_0(\lambda\pm i0)$ with the limit being understood in the pointwise sense for the kernels, see Lemma~\ref{lem:freeresol}. 
Next, we study the limiting absorption principle for $\bfR_0$. For $\alpha\in\bbR$, we define the weighted spaces
\begin{equation}
    L^{2,\alpha} := \bigl\{ f \, \big|  \, \jx^\alpha f \in L^2 \bigr\},
\end{equation}
and for $k \in\mathbb{N}$ we set
\begin{equation}
     H^{k,\alpha} := \big\{ f \, \big| \, \jx^\alpha \px^j f \in L^2, \, 0 \leq j \leq k \big\}.
\end{equation}
For vector functions, we set
\begin{equation}
    \calH^\alpha:=H^{1,\alpha}\times L^{2,\alpha}.
\end{equation}
\begin{definition}  
We use $\mathcal{B}(X_1,X_2)$ to denote the space of bounded linear operators from a Banach space $X_1$ to another Banach space $X_2$.
\end{definition}

We start by recalling some standard limiting absorption mapping properties  for the free scalar resolvent \eqref{eq:1dimresolv} based on its explicit formula. See for example, \cite{Agmon, JenKa, KoKo1dkg, Murata}.

\begin{lemma} \label{lem:freeLAP}
For $\frakR_0(\zeta)$ from \eqref{eq:1dimresolv}, we have the following mapping properties:
\begin{itemize}
    \item For $\zeta>0$, $\frakR_0(\zeta\pm i\epsilon)\rightarrow \frakR_0(\zeta\pm i0)$  as $\epsilon\rightarrow0^+$ in $\mathcal{B}(H^{-1,\sigma},H^{1,-\sigma})$ for $\sigma>\frac{1}{2}$. Moreover, this convergence is uniform for $\zeta\geq r$ for any $r>0$.
    \item  For any $M\geq0$ and $\zeta\in\bbC\setminus[0,\infty)$, the following expansion holds
\begin{equation}\label{eq:freeresolexp}
    \frakR_0(\zeta)=\sum_{j=-1}^M A_j \zeta^{j}+\calO (\zeta^{M+1}), \quad \zeta\rightarrow 0
\end{equation}in the norm of  $\mathcal{B}(H^{-1,\sigma},H^{1,-\sigma})$ for $\sigma>\frac{3}{2}+M+1$. Moreover, one has that the operator $A_{-1}$ is multiplication by $-\frac{i}{2}$, $A_0$ is a convolution operator with kernel $\frac{1}{2}|x-y|$, and $A_j\in\mathcal{B}(H^{-1,\sigma},H^{1,-\sigma})$  for $\sigma>\frac{3}{2}+k$ and ever integer $j \geq -1$.
\item For any $r>0$, $\zeta\in\bbC\setminus[0,\infty)$, and $|\zeta|\geq r$, it holds that
\begin{equation}\label{eq:largezeta}
    \|\frakR_0(\zeta)\|_{\mathcal{B}(H^{m,\sigma},H^{m+n,-\sigma})}\lesssim |\zeta|^{-(1-n)}
\end{equation} 
for $m=0,1$ and $n=-1,0,1$ with $\sigma>\frac{1}{2}$.
\end{itemize}
\end{lemma}

Using the explicit formulas from Lemma \ref{lem:freeresol} and Lemma \ref{lem:freeLAP}, the mapping properties for the free matrix resolvent \eqref{eq:freesol} can be computed directly upong making the change of variables $\zeta^2=\gamma^2\lambda^2-1.$

\begin{corollary} \label{cor:freeMLAP}
 For $\bfR_0(\lambda)$ from \eqref{eq:freesol}, we have the following mapping properties.
\begin{itemize}
    \item For $\lambda\in \Sigma_\ell$, one has
\begin{equation}\label{eq:freeRLAP}
    \bfR_0(\lambda\pm i\epsilon)\rightarrow \bfR^{\pm}_0(\lambda) := \bfR_0(\lambda i \pm 0)
\end{equation}
   as $\epsilon\rightarrow 0^+$ in $\mathcal{B}(\calH^\sigma,\calH^{-\sigma})$ for $\sigma>\frac{1}{2}$.
   \item One has the following asymptotic expansion: for $\lambda\in\bbC\setminus \Sigma_\ell$ and $\lambda \rightarrow \pm\frac{1}{\gamma}$,
\begin{equation}\label{eq:freeRexp}
\bfR_0(\lambda)=B^{\pm}_{-1}\frac{1}{\sqrt{\eta}}+B^{\pm}_0+\mathcal{O}(|\eta|^{\frac{1}{2}}), \quad \eta =\lambda\pm \frac{1}{\gamma}
\end{equation}
in $\mathcal{B}(\calH^\sigma,\calH^{-\sigma})$ for $\sigma>\frac{5}{2}$, where $B^{\pm}_{-1}\in \mathcal{B}(\calH^\sigma,\calH^{-\sigma})$ for $\sigma>\frac{1}{2}$ and $B^{\pm}_{0}\in \mathcal{B}(\calH^\sigma,\calH^{-\sigma})$ for $\sigma>\frac{3}{2}$. 
    \item For any $r>0$, one has
\begin{equation}\label{eq:freeRlarge}
    \|\bfR_0(\lambda)\|_{\mathcal{B}(\calH^\sigma,\calH^{-\sigma})}<\infty
\end{equation}for $\lambda\in\bbC\setminus\Sigma_\ell$ and $|\lambda\pm\frac{1}{\gamma}|\geq r$ with $\sigma>\frac{1}{2}$.
\end{itemize}   
\end{corollary}
\begin{proof}
For the limiting absorption principle, we observe that scaling and complex phase multiplication do not change mapping properties, whence the convergence properties from Lemma~\ref{lem:freeLAP} hold for $\mathscr{R}_\ell(\lambda)$.
Since  scaling and complex phase multiplication do not change the order of the expansions, \eqref{eq:freeresolexp} holds for $\mathscr{R}_\ell(\lambda)$ as well.
Now we study the resolvent estimates corresponding to \eqref{eq:largezeta}.  Note that for $n=0,1$, the estimates for $\mathscr{R}_\ell(\lambda)$ follow from \eqref{eq:largezeta} and the explicit formula \eqref{equ:mathscrRell_definition}. So it suffices to consider $n=-1$. This does not follow from a direct computation since the additional differentiation in the input will hit the complex phase which will result in an additional growth in $\lambda$. To resolve this, we note that from \eqref{equ:mathscrRell_definition}, we have
\begin{equation}
    \bigl( - (1-\ell^2) \partial^2 - 2i\ell\lambda\partial + 1 - \lambda^2 \bigr) \mathscr{R}_\ell(\lambda)=1,
\end{equation}
which implies that 
\begin{equation}\label{eq:rewriteRell}
    \mathscr{R}_\ell(\lambda)=\frac{1}{1-\lambda^2}\Big(1+(1-\ell^2) \partial^2\mathscr{R}_\ell(\lambda) + 2i\ell\lambda\partial \mathscr{R}_\ell(\lambda) \Big).
\end{equation}
Applying the estimate with $n=1$ and $n=0$, it follows that
\begin{align}
\|\partial^2\mathscr{R}\|_{\mathcal{B}(H^{m,\sigma},H^{m-1,-\sigma})}&=   \|\mathscr{R}\|_{\mathcal{B}(H^{m,\sigma},H^{m+1,-\sigma})}\lesssim 1,\\
 \|\partial^2\mathscr{R}\|_{\mathcal{B}(H^{m,\sigma},H^{m-1,-\sigma})}&=   \|\mathscr{R}\|_{\mathcal{B}(H^{m,\sigma},H^{m,-\sigma})}\lesssim \frac{1}{|\lambda|}.
\end{align}
The two preceding bounds combined with \eqref{eq:rewriteRell} yield the desired estimate.

For the free matrix resolvent, we can apply the results for $\mathscr{R}$ in each entry. Then  \eqref{eq:freeRLAP}, \eqref{eq:freeRexp}, \eqref{eq:freeRlarge} follow directly.
\end{proof}

\subsubsection{Perturbed resolvents}

Using resolvent identities, we can extend the mapping properties for the free matrix resolvent to perturbed matrix resolvents. We start with general computations for a perturbed version of Lemma \ref{lem:freeresol}.

For $\lambda\in\bbC\setminus\Sigma_\ell$, we define the perturbed resolvent
\begin{equation}\label{eq:Presol}
   \bfR(\lambda):=(\bfL-i\lambda\mathbb{I}_{2\times2})^{-1}.
\end{equation}
By the resolvent identity $(A+B)^{-1} = (\mathrm{Id} + A^{-1} B)^{-1} A^{-1}$ with $A := \bfL_0 - i \lambda \bbI_{2\times 2}$ and $B := {\bf V}_\ell$, we can write
\begin{equation}\label{eq:resolindeM}
\bfR(\lambda)=\left(\bbI_{2 \times 2} + \bfR_0(\lambda){\bf V}_{\ell}\right)^{-1}\bfR_0(\lambda)
\end{equation}
where
\begin{equation}
    {\bf V}_{\ell} = \begin{bmatrix} 0 & 0 \\ -V_\ell & 0 \end{bmatrix}.
\end{equation}
In particular, in the sine-Gordon setting, $V_\ell$ is given by \eqref{eq:Lell}.
Recall from Lemma \ref{lem:freeresol} that
 \begin{equation} \label{eq:freesolM}
 \begin{aligned}
    \bfR_0(\lambda)= \big( \bfL_0-i\lambda\mathbb{I}_{2\times2}\big)^{-1} &= \mathscr{R}_\ell(\lambda) \begin{bmatrix}
	\ell\partial -i\lambda & -1 \\
	-\partial^2+1 & \ell\partial -i\lambda
	\end{bmatrix} \\
    &= \begin{bmatrix}
	\mathscr{R}_\ell(\lambda) (\ell\partial -i\lambda) & -\mathscr{R}_\ell(\lambda) \\
	1- (\ell\partial -i\lambda)\mathscr{R}_\ell(\lambda) (\ell\partial -i\lambda) & (\ell\partial -i\lambda) \mathscr{R}_\ell(\lambda) 
	\end{bmatrix} 
 \end{aligned}
 \end{equation}
with $\mathscr{R}_\ell(\lambda)$ defined in \eqref{equ:mathscrRell_definition}.
By direct computation we obtain
\begin{align}
    \bfR_0(\lambda){\bf V}_{\ell} &= \begin{bmatrix}
    \left(\ell\partial-i\lambda\right)\mathscr{R}_\ell\left(\lambda\right) & -\mathscr{R}_\ell\left(\lambda\right) \\
    1-\left(\ell\partial-i\lambda\right)\mathscr{R}_\ell\left(\lambda\right) \left(\ell\partial-i\lambda\right) & \left(\ell\partial-i\lambda\right)\mathscr{R}_\ell\left(\lambda\right)
    \end{bmatrix}
    \begin{bmatrix}
        0 & 0\\
        -V_{\ell} & 0
    \end{bmatrix} \\
    &= \begin{bmatrix}
        \mathscr{R}_\ell\left(\lambda\right)V_{\ell} & 0\\
        -\left(\ell\partial-i\lambda\right)\mathscr{R}_\ell\left(\lambda\right)V_{\ell} & 0
       \end{bmatrix},
\end{align}
whence 
\begin{equation} \label{eq:resolindMF}
    \begin{aligned}
    \bigl( \bbI_{2 \times 2} + \bfR_0(\lambda){\bf V}_{\ell} \bigr)^{-1} &= \begin{bmatrix}
        1+\mathscr{R}_\ell\left(\lambda\right)V_{\ell} & 0\\
        -\left(\ell\partial - i\lambda\right)\mathscr{R}_\ell\left(\lambda\right)V_{\ell} & 1
        \end{bmatrix} \\
    &= \begin{bmatrix}
        \left(1+\mathscr{R}_\ell\left(\lambda\right)V_{\ell}\right)^{-1} & 0\\
        \left(\ell\partial-i\lambda\right)\left(1-\left(1+\mathscr{R}_\ell\left(\lambda\right)V_{\ell}\right)^{-1}\right) & 1
      \end{bmatrix}.
    \end{aligned}
\end{equation}
Next, we define
\begin{equation}
  R_{\ell}\left(\lambda\right) := \bigl( - (1-\ell^2) \partial^2 + V_\ell - 2 i \ell\lambda \partial + 1 - \lambda^2 \bigr)^{-1}.
\end{equation}
By the resolvent identity, we can then write
\begin{equation}
    R_{\ell}\left(\lambda\right) = \bigl(1+\mathscr{R}_\ell\left(\lambda\right)V_{\ell} \bigr)^{-1} \mathscr{R}_\ell\left(\lambda\right).
\end{equation}
Then by \eqref{eq:resolindeM}, \eqref{eq:freesolM}, and \eqref{eq:resolindMF}, we have
{\footnotesize
\begin{equation} \label{eq:resolP}
\begin{aligned}
   \bfR\left(\lambda\right) &= \begin{bmatrix}
    \left(1+\mathscr{R}_\ell\left(\lambda\right)V_{\ell}\right)^{-1} & 0\\
    \left(\ell\partial-\lambda\right)\left(1-\left(1+\mathscr{R}_\ell\left(\lambda\right)V_{\ell}\right)^{-1}\right) & 1
    \end{bmatrix} 
    \begin{bmatrix}
        \left(\ell\partial-i\lambda\right)\mathscr{R}_\ell\left(\lambda\right) & -\mathscr{R}_\ell\left(\lambda\right)\\
        1-\left(\ell\partial-i\lambda\right)\mathscr{R}_\ell\left(\lambda\right) \left(\ell\partial-i\lambda\right)& \left(\ell\partial-i\lambda\right)\mathscr{R}_\ell\left(\lambda\right)
    \end{bmatrix} \\
    &= \begin{bmatrix}
        R_{\ell}\left(\lambda\right)\left(\ell\partial-i\lambda\right) & -R_{\ell}\left(\lambda\right)\\
        1-\left(\ell\partial-i\lambda\right)R_{\ell}\left(\lambda\right)\left(\ell\partial-i\lambda\right) & \left(\ell\partial-i\lambda\right)R_{\ell}\left(\lambda\right)
       \end{bmatrix}.
\end{aligned}
\end{equation}}

From the general identity \eqref{eq:resolP} for the matrix resolvent, we now infer the explicit integral kernel for the resolvent of the matrix linearized operator around the moving sine-Gordon kink. 

\begin{lemma} \label{lem:Rexpression}
    For all $\lambda \in \bbC \backslash \bigl( \{0\} \cup \Sigma_\ell \bigr)$ the resolvent of the matrix operator $\bfL$ defined in \eqref{eq:bfL} is given by
    \begin{equation} \label{equ:bfL_resolvent}
        \bfR(\lambda) := \bigl( \bfL - i \lambda \mathbb{I}_{2\times2} \bigr)^{-1} = \begin{bmatrix} R_\ell(\lambda) (\ell \partial - i \lambda) & -R_\ell(\lambda) \\ 1 - (\ell\partial - i \lambda) R_\ell(\lambda) (\ell \partial - i \lambda) & (\ell \partial - i \lambda) R_\ell(\lambda) \end{bmatrix}
    \end{equation}
    with
    \begin{equation} \label{equ:scalar_Rell_def}
        R_\ell(\lambda) := \bigl( - \gamma^{-2} \partial^2 - 2 \sech^2(\gamma \cdot) + 1 - \lambda^2 - 2 i \ell \lambda \partial \bigr)^{-1}.
    \end{equation}
    The kernel of $R_\ell(\lambda)$ is given by
    \begin{equation} \label{equ:scalar_Rell_kernel}
        R_\ell(\lambda)(x,y) := e^{-i\gamma^2 \ell \lambda x} \widetilde{R}_\ell(\lambda)(x,y) e^{i\gamma^2\ell \lambda y},
    \end{equation}
    where
    \begin{equation} \label{equ:wtilRell_Green_function}
     \begin{aligned}
        \widetilde{R}_\ell(\lambda)(x,y) &:= \bigl( - \gamma^{-2} \partial^2 - 2 \sech^2(\gamma \cdot) + 1 - \gamma^2 \lambda^2 \bigr)^{-1}(x,y) \\
        &= \gamma^2 \bigl( W\bigl( f_+(\gamma \cdot, \zeta), f_-(\gamma \cdot, \zeta) \bigr) \bigr)^{-1} \Bigl( f_+(\gamma x,\zeta) f_-(\gamma y,\zeta) \one_{[x>y]} + f_-(\gamma x,\zeta) f_+(\gamma y,\zeta) \one_{[x<y]} \Bigr)
     \end{aligned}
    \end{equation}
    for $\zeta^2 = \gamma^2 \lambda^2 - 1$ with $\zeta \in \bbC_+$, and
    \begin{equation}\label{eq:fpmzeta}
        f_\pm(\gamma x, \zeta) = c_\pm(\zeta) \bigl( \mp i \zeta + \tanh(\gamma x) \bigr) e^{\pm i \gamma x \zeta}, \quad c_\pm(\zeta) = \mp (i\zeta - 1)^{-1}.
    \end{equation}
    Here, $f_\pm(\cdot, \zeta)$ denote the usual Jost solutions associated with $-\partial^2 - 2 \sech^2(\cdot)$.
\end{lemma}
\begin{proof}
From the calculations above for \eqref{eq:resolP}, we have verified that for $\lambda \in \bbC \backslash \bigl( \{0\} \cup \Sigma_\ell \bigr)$ the resolvent $(\bfL - i\lambda)^{-1}$ is given by \eqref{equ:bfL_resolvent}.
 The identity \eqref{equ:scalar_Rell_kernel} for the kernel of $R_\ell(\lambda)$ is a consequence of the following computation
 \begin{equation*}
     \begin{aligned}
         e^{-i\gamma^2 \ell \lambda x} \bigl( - \gamma^{-2} \partial_x^2 - 2 \sech^2(\gamma x) + 1 - \gamma^2 \lambda^2 \bigr) e^{i\gamma^2 \ell \lambda x} = -\gamma^{-2} \partial_x^2 - 2 \sech^2(\gamma x) + 1 - \lambda^2 - 2i\ell \lambda \partial_x.
     \end{aligned}
 \end{equation*}
 In order to arrive at the expression for the Green's function \eqref{equ:wtilRell_Green_function}, we observe that the (rescaled) Schr\"odinger operator
 \begin{equation}\label{eq:Hell}
    H_\ell := - \gamma^{-2} \partial^2 - 2 \sech^2(\gamma \cdot)
 \end{equation}
 satisfies
 \begin{equation}
    H_\ell = T_\ell \circ H \circ T_\ell^{-1}, \quad (T_\ell f)(x) := f(\gamma x), \quad H := - \partial^2 - 2 \sech^2(\cdot).
 \end{equation}
 Correspondingly, it follows from \eqref{eq:JostH} that the functions 
 \begin{equation}
      g_\pm(x,\zeta) := f_\pm(\gamma x, \zeta) = T_\ell f_\pm(x,\zeta) = c_\pm(\zeta) \bigl( \mp i \zeta + \tanh(\gamma x) \bigr) e^{\pm i \gamma x \zeta}
 \end{equation}
 satisfy
 \begin{equation*}
    e^{\mp i \gamma x \zeta} g_{\pm}(x, \zeta) \to 1 \quad \text{as} \quad x \to \pm \infty,
 \end{equation*}
 and
 \begin{equation}
    H_\ell g_\pm(x,\zeta) = \zeta^2 g_\pm(x,\zeta).
 \end{equation}
 From the latter and \eqref{eq:Hresol}, we infer \eqref{equ:wtilRell_Green_function}.
\end{proof}

Using the explicit formulas for the integral kernel from Lemma~\ref{lem:Rexpression}, we now present a streamlined proof of the mapping properties for the resolvent of the matrix linearized operator around the moving sine-Gordon kink. 

\begin{corollary} \label{cor:LAP}
The matrix resolvent $\bfR(\lambda)$ defined in \eqref{equ:bfL_resolvent} satisfies the following mapping properties:
\begin{itemize}
    \item For $\lambda\in \Sigma_\ell$, one has
\begin{equation}\label{eq:RLAP}
    \bfR(\lambda\pm i\epsilon)\rightarrow \bfR^{\pm}(\lambda)=\bfR(\lambda\pm i0)
\end{equation}as $\epsilon\rightarrow 0^+$ in $\mathcal{B}(\calH^\sigma,\calH^{-\sigma})$ for $\sigma>\frac{1}{2}$.
   \item One has the following asymptotic expansion: for $\lambda\in\bbC\setminus \Sigma_\ell$ and $\lambda\rightarrow \pm \frac{1}{\gamma}$
\begin{equation}\label{eq:Rexp}
    \bfR(\lambda)={\bf B}^{\pm}_{-1}\frac{1}{\sqrt{\eta}}+{\bf B}^{\pm}_0+\mathcal{O}(|\eta|^\frac{1}{2}), \eta =\lambda\pm \frac{1}{\gamma}
\end{equation}
in $\mathcal{B}(\calH^\sigma,\calH^{-\sigma})$ for $\sigma>\frac{5}{2}$ where ${\bf B}^{\pm}_{-1}\in \mathcal{B}(\calH^\sigma,\calH^{-\sigma})$ for $\sigma>\frac{1}{2}$ and ${\bf B}^{\pm}_{0}\in \mathcal{B}(\calH^\sigma,\calH^{-\sigma})$ for $\sigma>\frac{3}{2}$. 
 \item As $|\lambda|\rightarrow\infty$, one has
\begin{equation}\label{eq:Rlarge}
    \|\bfR(\lambda)\|_{\mathcal{B}(\calH^\sigma,\calH^{-\sigma})}<\infty
\end{equation}
for $\lambda\in\bbC\setminus\Sigma_\ell$  with $\sigma>\frac{1}{2}$.    
\end{itemize}
\end{corollary}
\begin{proof}
All asserted mapping properties follow from the explicit formulas obtained in Lemma~\ref{lem:Rexpression} in a similar manner as in the proof of Corollary~\ref{cor:freeMLAP}. So we only provide a sketch of the proof here.
The first claim follows from the explicit formulas \eqref{equ:bfL_resolvent}, \eqref{equ:scalar_Rell_def}, \eqref{equ:scalar_Rell_kernel}, \eqref{equ:wtilRell_Green_function} upon passing to the limit.    
The second claim follows from \eqref{equ:wtilRell_Green_function} by a direct expansion using \eqref{eq:TW} and the observation that $W\bigl( f_+(\gamma \cdot, \zeta), f_-(\gamma \cdot, \zeta) \bigr)\sim \zeta= \sqrt{\gamma^2\lambda^2-1}=\sqrt{\gamma\lambda-1} \sqrt{\gamma\lambda+1}$ for $|\zeta|\rightarrow0$.
The last claim is a direct consequence of the explicit formulas \eqref{equ:bfL_resolvent}, \eqref{equ:scalar_Rell_def}, \eqref{equ:scalar_Rell_kernel}, \eqref{equ:wtilRell_Green_function}, \eqref{eq:fpmzeta}, as in the case of the proof of \eqref{eq:freeRlarge} for the free matrix resolvent.
\end{proof}

\begin{remark}
 Note that since $-1$ is an eigenvalue for $H_\ell$ defined in \eqref{eq:Hell}, the corresponding resolvent \eqref{equ:wtilRell_Green_function} blows up with order $\lambda^{-2}$ near $\lambda=0$, see Jensen-Kato \cite{JenKa}, and the same applies to \eqref{equ:scalar_Rell_kernel}. This property is also carried over to the matrix resolvent \eqref{equ:bfL_resolvent}. So we expect that the resolvent $\bfR(\lambda)$ has a pole of order two at $\lambda=0$.  This matches the calculations in Lemma \ref{lem:ranPd}.
\end{remark}

\begin{remark}
We emphasize that our streamlined proof of Corollary~\ref{cor:LAP} is based on the explicit formulas for the Jost functions associated with the scalar linearized operator around the (static) sine-Gordon kink. An alternative general way to show the asserted mapping properties is tu use the identity \eqref{eq:resolindeM} and the free estimates in Corollary \ref{cor:freeMLAP}. Note that the existence and mapping properties of $\left(\bbI_{2 \times 2} + \bfR_0(\lambda){\bf V}_{\ell}\right)^{-1}$ are related to the existence of resonances and eigenvalues. 
\end{remark}

\subsection{Representation of the evolution via the jump of the resolvent}

In this subsection, we develop a representation formula for the linear evolution $e^{t\bfL} \Pe$ in terms of the jump of the resolvent across the essential spectrum. Note that this type of representation follows from Stone's formula in the self-adjoint setting. 

\begin{lemma} \label{lem:jumpformula} 
It holds that
\begin{align}\label{eq:rep1}
    e^{t\bfL} = \frac{1}{2\pi} \int_{\Sigma_\ell} e^{i\lambda t} \bigl[ \bfR(\lambda+i0)-\bfR (\lambda-i0) \bigr] \, \ud \lambda + e^{t\bfL}P_\mathrm{d},
\end{align}
where $\Sigma_\ell$ is given by \eqref{eq:Sigmaell}. 
The representation above and the convergence hold in the following weak sense: for a Schwartz function $\bm{\phi}\in \mathcal{S}_x$ and a space-time Schwartz function $\bm{\psi}(t,x)\in \mathcal{S}_{t,x}$, the limit
\begin{align*}
\int \left\langle e^{t\bfL}\bm{\phi},\bm{\psi}\right\rangle_{\calH} \,\ud t &= \int \lim_{R\rightarrow\infty}\frac{1}{2\pi }\int_{R\geq|\lambda|,\,\lambda\in\Sigma_\ell}e^{i\lambda t}\left\langle  [\bfR(\lambda+i0)-\bfR (\lambda-i0)]\bm{\phi},\bm{\psi}(t)\right\rangle_{\calH} \, \ud \lambda \, \ud t \\
&\quad +\int \left\langle e^{t\bfL}P_\mathrm{d}\bm{\phi}, \bm{\psi}\right\rangle_{\calH}\,\ud t
\end{align*}
exists. 
\end{lemma}

\begin{remark}   
Note that \eqref{eq:freeRlarge} and \eqref{eq:Rlarge} are very different from the Schr\"odinger problem. The Schr\"odinger resolvent decays for large frequencies, while \eqref{eq:freeRlarge} and \eqref{eq:Rlarge} do not.
This is related to the fact that the  Klein-Gordon equation is a wave-type/relativistic equation. The $\frac{1}{2}$ decay for Schr\"odinger resolvents is related to the local smoothing effect. For wave-type equations, we have local energy decay. 
If we restrict the evolution onto the essential spectrum, i.e, $e^{t\bfL } P_\mathrm{e}$, then one can test against $\bm{\psi}\in L^2_t \calH^{\frac{1}{2}+}$ for $\bm{\phi}\in \calH^{\frac{1}{2}+}$. This is consistent  with the local energy decay for $e^{t\bfL} P_\mathrm{e}$.
\end{remark}
\begin{proof}
We first define the evolution $e^{t\bfL}$ by the Hille-Yosida theorem.  To verify the conditions of the Hille-Yosida theorem, we note that for  $a>0$ large enough, the resolvent of $\bfL-a$, $ \rho(\bfL-a)$, satisfies
\begin{equation}
  (0,\infty)\subset  \rho(\bfL-a)
\end{equation}
and one has the following decay estimate
\begin{equation}
    \|(\bfL-a-\lambda)^{-1}\|_\calH\leq \frac{1}{\lambda}
\end{equation} for all $\lambda>0$. Indeed, the inequality above follows from for $a$ large depending on ${\bf V}$ 
\begin{align}
    \|(\bfL-a-\lambda)^{-1}\|_\calH &= \Bigl\|(\bfL_0-a-\lambda)^{-1}\left(I+{\bf V}(\bfL_0-a-\lambda)^{-1}\right)^{-1}\Bigr\|_\calH \\
    &\leq \frac{1}{a+\lambda}\sum_{k=0}^\infty \Bigl(\frac{C}{a+\lambda}\Bigr)^k = \frac{1}{a-C+\lambda}\leq\frac{1}{\lambda},
\end{align}
where we used the decay for the free resolvent \eqref{eq:LowerR}. It follows that $e^{(\bf L-a)t}$ produces a contractive semigroup, whence
\begin{equation}\label{eq:energygrowth}
    \|e^{t\bfL}\|_{\calH}\leq e^{ta}, \quad\forall t\in \bbR.
\end{equation}
Let $\Re z >a$. Then we can use the Laplace transform to write 
\begin{equation}\label{eq:laptrans}
    (\bfL-z)^{-1}=-\int_0^\infty e^{-tz} e^{t\bfL}\,\ud t
\end{equation}where the convergence is in the norm sense.
Let $b>a$ and $t>0$. One can define the inverse Laplace transform as
\begin{equation}
    e^{t\bfL}=-\frac{1}{2\pi i} \int_{b-i\infty}^{b+i\infty}e^{tz}(\bfL-z)^{-1}\,\ud z.
\end{equation}
The convergence above can be justified in the standard manner as follows. For $\bm{\phi},\bm{\varsigma}\in\mathcal{S}$ one has for $t\geq0$
\begin{equation}\label{eq:inlapR}
    \langle e^{t\bfL} \bm{\phi}, \bm{\varsigma} \rangle_{\calH}=-\lim_{R\rightarrow\infty}\frac{1}{2\pi i}\int_{b-iR}^{b+iR} e^{tz} \langle (\bfL-z)^{-1} \bm{\phi}, \bm{\varsigma} \rangle_{\calH}\,\ud z.
\end{equation}
Indeed, for $t>0$, using \eqref{eq:laptrans}, we can write
\begin{align}
    -\frac{1}{2\pi i}\int_{b-iR}^{b+iR} e^{tz} \langle (\bfL-z)^{-1} \bm{\phi}, \bm{\varsigma} \rangle_{\calH}\,\ud z&=\frac{1}{2\pi i} \int_{b-iR}^{b+iR} e^{tz} \int_0^{\infty}e^{-sz}\langle  e^{s\bfL}\bm{\phi}, \bm{\varsigma} \rangle_{\calH}\,\ud s \, \ud z\\
    &=\frac{1}{\pi}\int_0^\infty e^{(t-s)b}\frac{\sin((t-s)R)}{t-s}\langle  e^{s\bfL}\bm{\phi}, \bm{\varsigma} \rangle_{\calH}\,\ud s\\
    &=e^{tb}\frac{1}{\pi}\int_0^\infty e^{-s b}\frac{\sin((t-s)R)}{t-s}\langle  e^{s\bfL}\bm{\phi}, \bm{\varsigma} \rangle_{\calH}\,\ud s.
\end{align}
Note that $e^{-s b}\langle e^{s\bfL}\bm{\phi}, \bm{\varsigma} \rangle_{\calH}$ is continuously differentiable with respect to the variable $s$ and decays exponentially in $s$ because of \eqref{eq:energygrowth} and the fact that $b>a$. Then \eqref{eq:inlapR} follows from the standard properties of the Dirichlet kernel. Also note that when $t<0$, the limit is zero. Combining with the integral corresponding to $t<0$, we conclude that
\begin{align}\label{eq:laplace1}
       \langle e^{t\bfL} \bm{\phi}, \bm{\varsigma} \rangle_{\calH}&=-\lim_{R\rightarrow\infty}\biggl(\frac{1}{2\pi i}\int_{b-iR}^{b+iR} e^{tz} \langle (\bfL-z)^{-1} \bm{\phi}, \bm{\varsigma} \rangle_{\calH}\,\ud z\\
       &\quad \quad \quad \quad \quad \quad -\frac{1}{2\pi i}\int_{-b-iR}^{-b+iR} e^{tz} \langle (\bfL-z)^{-1} \bm{\phi}, \bm{\varsigma} \rangle_{\calH}\,\ud z \biggr)\nonumber.
\end{align}

Next, we perform a contour deformation to push the integral onto the spectrum of $\bfL$.

\begin{figure}[ht]
\centering
\begin{tikzpicture}

    \def\R{3}      
    \def\b{5}      
    \def\u{1}    
    \def\d{0.3}    

    \draw[thick,->] (-0.5,0) -- (7,0) node[right] {$\Re$}; 
    \draw[thick,->] (0,-4) -- (0,4) node[above] {$\Im$}; 

    \draw[thick] (\b,-\R) -- (\b,\R) node[right] {$b + iR$};
    \node[right] at (\b,-\R) {$b - iR$};

    \draw[thick] (0.1,\R) -- (\b,\R); 
    \draw[thick] (0.1,-\R) -- (\b,-\R); 

    \draw[thick] (0.1, 0.7) arc[start angle=-90,end angle=90,radius=\d]; 
    \draw[thick] (0.1,-\u-\d) arc[start angle=-90,end angle=90,radius=\d]; 
    \draw[thick] (0.1,-0.3) arc[start angle=-90,end angle=90,radius=\d]; 

    \draw[thick] (0.1,\u+\d) -- (0.1,\R); 
    \draw[thick] (0.1,-\u-\d) -- (0.1,-\R); 
    \draw[thick] (0.1,-\d) -- (0.1,-\u+\d); 
    \draw[thick] (0.1,\d) -- (0.1,\u-\d); 

    \filldraw[blue] (0,0) circle (2pt) node[left]{}; 

    \filldraw[red] (0,\u) circle (2pt); 
    \filldraw[red] (0,-\u) circle (2pt); 
    \node[left] at (0,\u) {$i\gamma^{-1}$}; 
    \node[left] at (0,-\u) {$-i\gamma^{-1}$}; 

    \node at (-0.65,\R) {$\epsilon+iR$}; 
    \node at (-0.65,-\R) {$\epsilon-iR$}; 

\end{tikzpicture}
\caption{Contour: $\Gamma_{R,\delta,\epsilon}^+$.}
\label{fig:contour}
\end{figure}
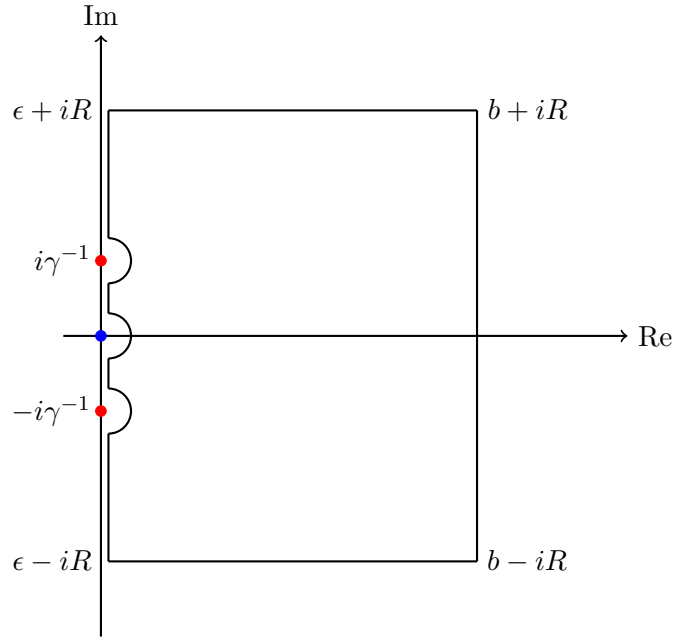

Since $(\bfL-z)^{-1}$ is not analytic at $z=0$ by Lemma \ref{lem:Lspec} and since it blows up near $\pm i \gamma^{-1}$ by \eqref{eq:Rexp} in Corollary \ref{cor:LAP}, for  $0<\epsilon\ll\delta\ll 1$ we consider the contour $\Gamma_{R,\delta,\epsilon}^+$, see Figure \ref{fig:contour},  which consists of the segment $b-iR$ to $b+iR$ as its right boundary, and the left boundary that contains semi-circular arcs on the right half of the complex plane of radius $\delta>0$ centered at the edges $\pm i \gamma^{-1}$ and $0$ and otherwise parallel segments  to the imaginary axis of distance $\epsilon$ to the imaginary axis, and finally the top boundary $\Im z =iR, \epsilon \leq \Re z\leq b, \epsilon \leq \Re z\leq b$, the bottom boundary $\Im z=-iR$. From the spectral information for $\bfL$, see Lemma \ref{lem:Lspec}, we conclude by Cauchy's theorem that for any $t\in\bbR$, $\bm{\phi}\in\mathcal{S}_x$, $\bm{\psi}\in\mathcal{S}_{t,x}$,
\begin{equation}
    \frac{1}{2\pi i} \oint_{\Gamma^+_{R,\delta,\epsilon}} e^{tz} \langle (\bfL-z)^{-1} \bm{\phi}, \bm{\psi}(t) \rangle_{\calH}\,\ud z=0.
 \end{equation}
 We can define a similar contour $\Gamma_{R,\delta}^{-}$, which is the reflection of $\Gamma^+_{R,\delta}$ about the imaginary axis $i\bbR$ on the left half of the complex plane. Again by Cauchy's theorem, we have
 \begin{equation}
    \frac{1}{2\pi i} \oint_{\Gamma^-_{R,\delta,\epsilon}} e^{tz} \langle (\bfL-z)^{-1} \bm{\phi}, \bm{\psi}(t) \rangle_{\calH}\,\ud z=0.
\end{equation}
Consider the integral over the top boundary. Integrating by parts in $t$, one can write
\begin{align}
&\int\frac{1}{2\pi i} \int_{\epsilon+iR,b+iR}  e^{tz} \langle (\bfL-z)^{-1} \bm{\phi}, \bm{\psi}(t) \rangle_{\calH}\,\ud z \, \ud t \\
&=  \frac{1}{2\pi i} \int_{\epsilon+iR,b+iR} \int e^{tz} \frac{1}{(z+\kappa)^2}  \langle (\bfL-z)^{-1} \bm{\phi}, (\partial_t+\kappa)^2 \bm{\psi}(t) \rangle_{\calH}\, \ud t \, \ud z 
\end{align}
for some $\kappa>0$. We note that using the boundedness of the resolvent \eqref{eq:Rlarge} and letting $R\rightarrow\infty$ in the formula above, it follows that
\begin{equation}\label{eq:laplace2}
   \lim_{R\rightarrow\infty} \int \frac{1}{2\pi i} \int_{\epsilon+iR,b+iR} e^{tz} \langle (\bfL-z)^{-1} \bm{\phi}, \bm{\psi}(t) \rangle_{\calH}\,\ud z \, \ud t=0.
\end{equation}
By the same reasons, one also has\begin{equation}\label{eq:laplace3}
   \lim_{R\rightarrow\infty} \int \frac{1}{2\pi i} \int_{\epsilon-iR,b-iR} e^{tz} \langle (\bfL-z)^{-1} \bm{\phi}, \bm{\psi}(t) \rangle_{\calH}\,\ud z \, \ud t=0.
\end{equation}
Replacing $\bm{\varsigma}$ in \eqref{eq:laplace1} by $\bm{\psi}$, integrating the resulting formula in $t$, and applying \eqref{eq:laplace2}, \eqref{eq:laplace3} to send $R\rightarrow\infty$, we obtain
\begin{align} \label{eq:Cauchy}
\int \left\langle e^{t\bfL}\bm{\phi},\bm{\psi}\right\rangle_{\calH} \,\ud t &= \int \frac{1}{2\pi i}\int_{c^+_{\delta, \epsilon}} e^{tz} \langle (\bfL-z)^{-1} \bm{\phi}, \bm{\psi}(t) \rangle_{\calH}\,\ud z \, \ud t \\
&\quad +\int \frac{1}{2\pi i}\int_{c^-_{\delta, \epsilon}} e^{tz} \langle (\bfL-z)^{-1} \bm{\phi}, \bm{\psi}(t) \rangle_{\calH}\,\ud z \, \ud t
\\
&\quad +\int \frac{1}{2\pi i}\oint_{\gamma}e^{tz} \langle (\bfL-z)^{-1} \bm{\phi}, \bm{\psi}(t) \rangle_{\calH}\,\ud z \, \ud t.
\end{align}
In the preceding integrals,  $c^+_{\delta, \epsilon}$ is a contour which consists of a straight line $is+\epsilon$, $\gamma^{-1}+\delta\leq s<\infty$, and a circular loop of radius $\delta$ centered at $i \gamma^{-1}$, and then a straight line $is-\epsilon$, $\gamma^{-1} + \delta \leq s < \infty$.
The contour $c^-_{\delta,\epsilon}$ is defined similarly, and finally $\gamma$ is  a circular loop of sufficiently small radius centered at $0$.

Now we study the integral around the origin. Near $0$, from our analysis of the generalized kernel of $\bfL$, see Lemma \ref{lem:ranPd} and Lemma \ref{lem:Lspec}, we have 
\begin{equation}\label{eq:orderpole}
    \|(\bfL-z)^{-1}\|_{\calH}\lesssim \frac{1}{|z|^2}.
\end{equation}
Now using the Taylor expansion, we write
\begin{equation}\label{eq:taylor}
    e^{tz}=1+tz+p_r(tz)
\end{equation}where $|p_r(x)|\lesssim |x|^2$. 
Then we write
\begin{align}
    \frac{1}{2\pi i}\oint_{\gamma}e^{tz} \langle (\bfL-z)^{-1} \bm{\phi}, \bm{\psi}(t) \rangle_{\calH}\,\ud z
    &= \frac{1}{2\pi i}\oint_{\gamma}(e^{tz}-(1+tz)) \langle (\bfL-z)^{-1} \bm{\phi}, \bm{\psi}(t) \rangle_{\calH}\,\ud z\label{eq:cont1}\\
    &\quad + \frac{1}{2\pi i}\oint_{\gamma}((1+tz)) \langle (\bfL-z)^{-1} \bm{\phi}, \bm{\psi}(t) \rangle_{\calH}\,\ud z\\
    &= \frac{1}{2\pi i}\oint_{\gamma} \langle (1+t\bfL)(\bfL-z)^{-1} \bm{\phi}, \bm{\psi}(t) \rangle_{\calH}\,\ud z\label{eq:cont2}\\
    &\quad +\frac{1}{2\pi i}\oint_{\gamma} \langle t(z-\bfL)(\bfL-z)^{-1} \bm{\phi}, \bm{\psi}(t) \rangle_{\calH}\,\ud z\label{eq:cont3}\\
    &=\frac{1}{2\pi i}\oint_{\gamma} \langle e^{t\bfL} (\bfL-z)^{-1} \bm{\phi}, \bm{\psi}(t) \rangle_{\calH}\,\ud z\label{eq:cont4}\\
    &=\langle e^{t\bfL}  P_\mathrm{d} \bm{\phi}, \bm{\psi}(t) \rangle_{\calH}.
\end{align}
The right-hand side of \eqref{eq:cont1} vanishes  since Cauchy's theorem can be applied after the cancellation between the order of the pole \eqref{eq:orderpole} and the order of the remainder from the Taylor expansion \eqref{eq:taylor}. By the same reasoning, \eqref{eq:cont3} is zero. To pass from \eqref{eq:cont2} to \eqref{eq:cont4}, we used the fact that $\ker(\bfL^n)=\ker(\bfL^2)$. The last identity is the definition of the Riesz projection \eqref{eq:rieszproj}.
 
Returning to \eqref{eq:Cauchy}, letting $\epsilon\rightarrow 0^+$, and then sending $\delta\rightarrow 0$, by the limiting absorption principle, \eqref{eq:RLAP}, and the expansion \eqref{eq:Rexp} respectively, we conclude after setting $z = i\lambda$ that
\begin{align}
\int \left\langle e^{t\bfL}\bm{\phi},\bm{\psi}\right\rangle_{\calH} \,\ud t &= \int \frac{1}{2\pi }\int_{\lambda\in\Sigma_\ell}e^{i\lambda t}\left\langle  [(\bfL-(i\lambda+0^+))^{-1}-(\bfL-(i\lambda-0^+))^{-1}]\bm{\phi},\bm{\psi}(t)\right\rangle_{\calH} \ud\lambda \, \ud t \\
&\quad +\int \left\langle e^{t\bfL}P_\mathrm{d}\bm{\phi}, \bm{\psi}\right\rangle_{\calH}\,\ud t.
\end{align}
Note that the last term above is finite since the range of $P_\mathrm{d}$ is the generalized kernel of finite algebraic multiplicity of $\bfL$ by Lemma \ref{lem:ranPd}. The desired result follows.
\end{proof}

\begin{remark}
 Note that the free operator $\bfL _0$ is actually selfadjoint with respect the inner product in $\calH$.
 Moreover, Lemma~\ref{lem:calLsa} gives for all $\lambda\in\bbR\setminus\{0\}$ that
 \begin{equation}
    \|(\bfL_0-\lambda)^{-1}\|_{\calH\to\calH}\le |\lambda|^{-1}.
 \end{equation}
 Then the preceding decay estimate combined with the closedness of $\Dom(\bfL_0)$ and the Hille-Yosida theorem, implies that on the real Hilbert space $\calH$, there exists a unitary group $e^{t\bfL_0}$ generated by $\bfL_0$ for all $t \in \bbR$.  Alternatively, in the complexified setting, we can use the unitary evolution $e^{it(-i\bfL_0)}=e^{t\bfL_0}$ provided by Stone's theorem, and note that this evolution preserves real-valued functions. 
\end{remark}

\subsection{Representation of the evolution via the distorted Fourier transform}

Combining the representation of the evolution $e^{\bfL} \Pe$ in terms of the jump of the resolvent across the essential spectrum provided by Lemma~\ref{lem:jumpformula} with the expressions for the integral kernel of the resolvent obtained in Lemma~\ref{lem:Rexpression}, we now derive an oscillatory integral representation for the evolution $e^{t\bfL} \Pe$ from which we can read off the definition of the distorted Fourier transform associated with the matrix operator $\bfL$.

In what follows, we write 
\begin{equation*}
    \Sigma_\ell := \Sigma_\ell^{-} \cup \Sigma_\ell^{+}, \quad \Sigma_\ell^{-} := \bigl( -\infty, -\gamma^{-1} \bigr], \quad \Sigma_\ell^{+} := \bigl[ \gamma^{-1}, \infty\bigr).
\end{equation*}
Restricting the evolution to the essential spectrum, we have by Lemma~\ref{lem:jumpformula} that
\begin{equation}
    \begin{aligned}
        e^{t\bfL} P_\mathrm{e} \bmf = \frac{1}{2\pi} \int_{\Sigma_\ell} e^{it\lambda} \Bigl( \bfR(\lambda+i0) - \bfR(\lambda-i0) \Bigr) \bmf \, \ud \lambda
    \end{aligned}
\end{equation}
with $\bfR(\lambda)$ defined in \eqref{equ:bfL_resolvent}.
Observe that for \emph{real-valued} $\bmf = (f_1, f_2)$ we have formally
\begin{equation}
    \begin{aligned}
        \frac{1}{2\pi} \int_{\Sigma^+_\ell} e^{it\lambda} \Bigl( \bfR(\lambda+i0)-\bfR (\lambda-i0)\Bigr) \bmf \, \ud \lambda = \overline{\frac{1}{2\pi} \int_{\Sigma^-_\ell} e^{it\lambda} \Bigl( \bfR(\lambda+i0)-\bfR (\lambda-i0)\Bigr) \bmf \, \ud \lambda},
    \end{aligned}
\end{equation}
whence
\begin{equation}\label{eq:ReJump}
 \begin{aligned}
  e^{t\bfL} P_\mathrm{e} \bmf = \Re \biggl( \frac{1}{\pi} \int_{\Sigma^+_\ell} e^{it\lambda} \Bigl( \bfR(\lambda+i0)-\bfR (\lambda-i0)\Bigr) \bmf \, \ud \lambda \biggr).
 \end{aligned}
\end{equation}
In view of Lemma~\ref{lem:Rexpression}, we have
\begin{equation}
    \begin{aligned}
        &\bfR(\lambda+i0)-\bfR (\lambda-i0) \\
        &= \begin{bmatrix}
            \bigl( R_\ell(\lambda+i0) - R_\ell(\lambda-i0) \bigr) (\ell \partial - i \lambda) & - \bigl( R_\ell(\lambda+i0) - R_\ell(\lambda-i0) \bigr) \\
            - (\ell \partial - i\lambda) \bigl( R_\ell(\lambda+i0) - R_\ell(\lambda-i0) \bigr) (\ell \partial - i \lambda) & (\ell \partial - i \lambda) \bigl( R_\ell(\lambda+i0) - R_\ell(\lambda-i0) \bigr)
        \end{bmatrix}.
    \end{aligned}
\end{equation}
On $\Sigma_\ell^+$, changing variables to $\sigma = \sqrt{\gamma^2 \lambda^2 - 1}$, whence $\lambda = \gamma^{-1} \jap{\sigma}$ and $\ud \lambda = \gamma^{-1} \jap{\sigma}^{-1} \sigma$, we have
\begin{equation}\label{eq:lambdatosigma}
    \begin{aligned}
        {\bf k} = (\bfk_1, \bfk_2) &:=\frac{1}{\pi} \int_{\Sigma_\ell^+} e^{it\lambda} \Bigl( \bfR(\lambda+i0) - \bfR(\lambda-i0) \Bigr) \bmf \, \ud \lambda \\
        &= \frac{1}{\pi} \int_0^\infty e^{it\gamma^{-1} \jap{\sigma}} \bigl( \gamma^{-1} \jap{\sigma}^{-1} \sigma \bigr) \Bigl( \bfR(\gamma^{-1} \jap{\sigma} + i0) - \bfR(\gamma^{-1} \jap{\sigma} - i0) \Bigr) \bmf \, \ud \sigma.
    \end{aligned}
\end{equation}
We  next compute the jump of the resolvent across $\Sigma_\ell^+$ in terms of the Jost functions. After applying the change of variables $\sigma = \sqrt{\gamma^2 \lambda^2 - 1}$ for $\lambda \in \Sigma_\ell^+$, \eqref{equ:wtilRell_Green_function} implies for $\lambda \in \Sigma_\ell^+$ that
\begin{equation*}
    \begin{aligned}
        &\widetilde{R}_\ell(\lambda+i0)(x,y)=\widetilde{R}_\ell(\gamma^{-1} \jap{\sigma}+i0)(x,y)  \\
        &\quad = \gamma^2 \Bigl( W\bigl( f_+(\gamma\cdot, \sigma), f_-(\gamma \cdot, \sigma) \bigr) \Bigr)^{-1} \Bigl( f_+(\gamma x, \sigma) f_-(\gamma y, \sigma) \one_{[x>y]} + f_-(\gamma x, \sigma) f_+(\gamma y, \sigma) \one_{[x<y]} \Bigr), \\
    \end{aligned}
\end{equation*}
as well as 
\begin{equation*}
    \begin{aligned}
        &\widetilde{R}_\ell(\lambda-i0)(x,y)=\widetilde{R}_\ell(\gamma^{-1} \jap{\sigma}-i0)(x,y)  \\
        &\quad = \gamma^2 \Bigl( W\bigl( f_+(\gamma \cdot, -\sigma), f_-(\gamma \cdot, -\sigma) \bigr) \Bigr)^{-1} \\
        &\qquad \qquad \qquad \qquad \qquad \qquad \times \Bigl( f_+(\gamma x,-\sigma) f_-(\gamma y,-\sigma) \one_{[x>y]} + f_-(\gamma x,-\sigma) f_+(\gamma y,-\sigma) \one_{[x<y]} \Bigr).
    \end{aligned}
\end{equation*}
By direct computation, using the explicit formula of the transmission coefficient \eqref{eq:Htrans} and its relation with the Wronskian \eqref{eq:TW}, one has
\begin{equation}
 W\bigl( f_+(\gamma \cdot, \sigma), f_-(\gamma \cdot, \sigma) \bigr) = \gamma W\bigl( f_+(\cdot, \sigma), f_-(\cdot, \sigma) \bigr) =  \gamma (-2i\sigma) \frac{\sigma^2+1}{\sigma^2+2i\sigma-1} = \gamma \frac{-2i\sigma}{T(\sigma)}.
\end{equation}
It follows that
\begin{equation*}
    \begin{aligned}
        \widetilde{R}_\ell(\gamma^{-1} \jap{\sigma}+i0)(x,y) &= - \gamma \frac{T(\sigma)}{2i\sigma} \Bigl( f_+(\gamma x,\sigma) f_-(\gamma y,\sigma) \one_{[x>y]} + f_-(\gamma x,\sigma) f_+(\gamma y,\sigma) \one_{[x<y]} \Bigr), \\
        \widetilde{R}_\ell(\gamma^{-1} \jap{\sigma}-i0)(x,y) &= -\gamma \frac{T(-\sigma)}{2i(-\sigma)} \Bigl( f_+(\gamma x,-\sigma) f_-(\gamma y,-\sigma) \one_{[x>y]} + f_-(\gamma x,-\sigma) f_+(\gamma y,-\sigma) \one_{[x<y]} \Bigr).
    \end{aligned}
\end{equation*}
Hence, using that $f_\pm(\gamma x,\sigma) = T(-\sigma) f_\mp(\gamma x,-\sigma)$ by \eqref{eq:scattrela}, we find
\begin{equation}
 \begin{aligned}
  &\wtilR_\ell(\gamma^{-1} \jap{\sigma}+i0)(x,y) - \wtilR_\ell(\gamma^{-1} \jap{\sigma}-i0)(x,y) \\
  &= - \gamma \frac{1}{2i\sigma} \Bigl( T(\sigma) f_+(\gamma x,\sigma) f_-(\gamma y,\sigma) + T(-\sigma) f_+(\gamma x,-\sigma) f_-(\gamma y,-\sigma) \Bigr) \one_{[x>y]} \\
  &\quad \, - \gamma \frac{1}{2i\sigma} \Bigl( T(\sigma) f_-(\gamma x, \sigma) f_+(\gamma y, \sigma) + T(-\sigma) f_-(\gamma x,-\sigma) f_+(\gamma y,-\sigma) \Bigr) \one_{[x<y]} \\
  &= - \gamma \frac{|T(\sigma)|^2}{2i\sigma} \Bigl( f_+(\gamma x,\sigma) f_+(\gamma y,-\sigma) + f_-(\gamma x,\sigma) f_-(\gamma y,-\sigma) \Bigr),
 \end{aligned}
\end{equation}
and therefore
\begin{equation} \label{equ:jump_scalar_Rell}
 \begin{aligned}
  &R_\ell(\gamma^{-1} \jap{\sigma}+i0)(x,y) - R_\ell(\gamma^{-1} \jap{\sigma}-i0)(x,y) \\
  &= e^{-i\ell\gamma \jap{\sigma} x} \bigl( \wtilR_\ell(\lambda+i0)(x,y) - \wtilR_\ell(\lambda-i0)(x,y) \bigr) e^{i\ell\gamma \jap{\sigma} y} \\
  &= - \gamma \frac{|T(\sigma)|^2}{2i\sigma} e^{-i\ell \gamma \jap{\sigma} (x-y)} \Bigl( f_+(\gamma x,\sigma) f_+(\gamma y,-\sigma) + f_-(\gamma x,\sigma) f_-(\gamma y, -\sigma) \Bigr).
 \end{aligned}
\end{equation}
Before inserting \eqref{equ:jump_scalar_Rell} into the second line of the expression \eqref{eq:lambdatosigma}, we pass to the new frequency variable
\begin{equation} \label{equ:new_freq_variable_xi}
    \xi := \gamma (\sigma - \ell \jap{\sigma}).
\end{equation}
Then we have the relations
\begin{equation}
    \sigma = \gamma (\xi + \ell \jxi), \quad \jap{\sigma} = \gamma (\jxi + \ell \xi), \quad \frac{\ud \sigma}{\ud \xi} = \gamma \bigl( 1 + \ell \jxi^{-1} \xi \bigr) = \jxi^{-1} \jap{\sigma}.
\end{equation}
Hence, the second line of \eqref{eq:lambdatosigma} reads with respect to the new frequency variables
\begin{equation}\label{eq:sigmatoxi}
    \begin{aligned}
        &{\bf k} = (\bfk_1, \bfk_2) = \frac{1}{\pi} \int_{-\gamma\ell}^\infty e^{it(\jxi + \ell \xi)} \Bigl( \frac{\jxi + \ell \xi}{\jxi}\Bigr) \Bigl( \bfR(\ell \xi + \jxi+ i0) - \bfR(\ell \xi + \jxi - i0) \Bigr) \bmf \, \ud \xi.
    \end{aligned}
\end{equation}
Recall the formulas for $f_\pm$ from \eqref{eq:fpmzeta}. 
With respect to the new variable $\xi$, we introduce
\begin{align} \label{equ:dFT_basis}
 e_\ell(x,\xi) := \frac{1}{\sqrt{2\pi}} e^{i x \xi}
  \begin{cases}
   m_\ell^+(\gamma x,\xi) & \xi \geq -\gamma \ell,
   \\
   \\
   m_\ell^-(\gamma x,\xi) & \xi < -\gamma \ell,
  \end{cases}
\end{align}
where 
\begin{equation} 
    \begin{aligned}
        m_\ell^+(\gamma x,\xi) &:= \frac{\gamma(\xi+\ell\jxi) + i \tanh(\gamma x)}{\gamma(\xi+\ell\jxi) - i}, \\
        m_\ell^-(\gamma x,\xi) &:= \frac{\gamma(\xi+\ell\jxi) + i \tanh(\gamma x)}{\gamma(\xi+\ell\jxi) + i}.
    \end{aligned}
\end{equation}
Inserting the jump formula \eqref{equ:jump_scalar_Rell} into \eqref{eq:sigmatoxi}, we obtain in terms of the new frequency variable $\xi$ that
\begin{equation} \label{eq:k1}
\begin{aligned}
    {\bf k}_1(t,x) &= \int_\bbR\int_\bbR  e_\ell(x,\xi) e^{it(\jxi + \ell\xi)} i\jxi^{-1} \overline{e_\ell(y,\xi)} \bigl(\ell\partial_y-i(\jxi+\ell \xi) \bigr) f_1(y)\,\ud y \, \ud \xi \\
    &\quad - \int_\bbR\int_\bbR  e_\ell(x,\xi) e^{it(\jxi + \ell\xi)} i\jxi^{-1} \overline{e_\ell(y,\xi)} f_2(y) \, \ud y \, \ud \xi
\end{aligned}
\end{equation}
and
\begin{equation} \label{eq:k2}
\begin{aligned}
    {\bf k}_2(t,x) &= \int_\bbR\int_\bbR \bigl(-\ell\px+i(\jxi + \ell\xi)\bigr) e_\ell(x,\xi) e^{it(\jxi+\ell\xi)} i\jxi^{-1} \\ 
    &\qquad \qquad \qquad \qquad \qquad \qquad \qquad \qquad \times \overline{e_\ell(y,\xi)} \bigl(\ell\partial_y-i(\jxi+\ell \xi)\bigr) f_1(y)\,\ud y \, \ud \xi \\
    &\quad - \int_\bbR\int_\bbR \bigl(-\ell\px+i(\jxi +\ell\xi)\bigr) e_\ell(x,\xi) e^{it(\jxi+\ell\xi)} i\jxi^{-1} \overline{e_\ell(y,\xi)} f_2(y)\,\ud y \, \ud \xi.
\end{aligned}
\end{equation}

Putting the preceding computations together, we can now read off the scalar distorted transforms from which the distorted Fourier transform associated with $\bfL$ is built.
Inspecting \eqref{eq:k1} and \eqref{eq:k2} motivates the following definitions.
We introduce the scalar distorted transform $\calF_\ell$ as
\begin{align}
    \calF_\ell[g](\xi) :=  \int_\bbR \overline{e_\ell(x,\xi)} g(x) \, \ud x,
\end{align}
whose adjoint with respect to the scalar inner product \eqref{eq:L2inner}
is given by
\begin{equation}
    \calF_\ell^{\ast}[h](x) = \int_{\bbR} e_\ell(x,\xi) h(\xi) \, \ud \xi.
\end{equation}
Moreover, we introduce the ``differentiated'' scalar distorted transform
\begin{align}
\begin{split}
    \calF_{\ell,D}[g](\xi) := \int_\bbR \overline{e_\ell(x,\xi)} D g(x) \, \ud x, \quad D := \ell \partial_x - i (\jxi + \ell \xi),
\end{split}
\end{align}
for which the adjoint is given by
\begin{equation}
    \calF_{\ell,D}^\ast[h](x) = \int_\bbR D^\ast e_\ell(x,\xi) h(\xi) \, \ud \xi, \quad D^\ast = -\ell \partial_x + i  (\jxi + \ell \xi).
\end{equation}
Combining \eqref{eq:ReJump}, \eqref{eq:sigmatoxi}, \eqref{eq:k1}, and \eqref{eq:k2}, we arrive at the following representation formula for the linear evolution in terms of the scalar distorted transforms $\calF_\ell$, $\calF_{\ell,D}$, and their adjoints.
\begin{proposition} \label{prop:rep_formula_propagator_dFT}
For any pair of real-valued Schwartz functions $\bm{f} = (f_1,f_2) \in \calS(\bbR) \times \calS(\bbR)$ and for any $t \in \bbR$, we have the following representation formula for the evolution $e^{t\bfL} P_\mathrm{e} \bmf$,
\begin{equation}
    \begin{aligned}
        \big( e^{t\bfL} P_\mathrm{e} \bm{f} \big)_1(x) &= \Re \, \biggl( \calF_\ell^\ast\Bigl[ e^{i t (\jxi + \ell \xi)} i\jxi^{-1} \Bigl( \calF_{\ell,D}[f_1](\xi) - \calF_\ell[f_2](\xi) \Bigr) \Bigr](x) \biggr), \\
        \big( e^{t\bfL} P_\mathrm{e} \bm{f} \big)_2(x) &= \Re \, \biggl( \calF_{\ell,D}^\ast \Bigl[ e^{i t (\jxi + \ell \xi)} i{\jxi}^{-1} \Bigl( \calF_{\ell,D}[f_1](\xi) - \calF_\ell[f_2](\xi) \Bigr) \Bigr](x) \biggr).
    \end{aligned}
\end{equation}
\end{proposition}

In the specific case of the sine-Gordon model, an inspection of \eqref{equ:dFT_basis} shows that there is a jump discontinuity 
\begin{equation}
    \lim_{\xi \downarrow -\gamma \ell} m_\ell^+(\gamma x,\xi) = - \lim_{\xi \uparrow -\gamma \ell} m_\ell^-(\gamma x,\xi) = - \tanh(\gamma x).
\end{equation}
This discontinuity is related to the threshold resonances of the linearized operator $\bfL$.
For the nonlinear analysis it is preferable to have better regularity properties of the distorted transforms on the frequency side.
Following \cite{CP22}, we pass to the modified distorted Fourier basis
\begin{align} \label{equ:modified_dFT_basis}
 e_\ell^\#(x,\xi) := \frac{1}{\sqrt{2\pi}} e^{i x \xi}
  \begin{cases}
   m_\ell^+(\gamma x,\xi) & \xi \geq -\gamma \ell,
   \\
   \\
   - m_\ell^-(\gamma x,\xi) & \xi < -\gamma \ell.
  \end{cases}
\end{align}
The latter can be written more succinctly as 
\begin{equation} \label{equ:modified_dFT_basis_succinct}
    e_\ell^\#(x,\xi) = \frac{1}{\sqrt{2\pi}} m_\ell^\#(\gamma x,\xi) e^{i x \xi}
\end{equation}
with
\begin{equation}
    m_\ell^\#(\gamma x,\xi) := \frac{\gamma(\xi+\ell\jxi) + i \tanh(\gamma x)}{|\gamma(\xi+\ell\jxi)| - i}.
\end{equation}
Then the associated modified distorted scalar transforms are 
\begin{equation}\label{eq:modF}
    \begin{aligned}
        \calF^\#_\ell[g](\xi) &:= \int_\bbR \overline{e^{\#}_\ell(x,\xi)} g(x) \, \ud x, \\
        \calF_\ell^{\#,\ast}[h](x) &= \int_\bbR e^{\#}_\ell(x,\xi) h(\xi) \, \ud \xi,
    \end{aligned}
\end{equation}
as well as
\begin{equation}\label{eq:modFD}
    \begin{aligned}
        \calF^{\#}_{\ell,D}[g](\xi) &:= \int_\bbR \overline{e^{\#}_\ell(x,\xi)} D g(x) \, \ud x, \\
        \calF^{\#, \ast}_{\ell,D}[h](x) &= \int_\bbR D^\ast e^{\#}_\ell(x,\xi) h(\xi) \, \ud \xi,
    \end{aligned}
\end{equation}
where we recall that
\begin{equation*}
    D = \ell \px - i (\jxi + \ell \xi), \quad D^\ast = -\ell \px + i (\jxi + \ell \xi).
\end{equation*}
Using the preceding definitions, we immediately obtain from Proposition~\ref{prop:rep_formula_propagator_dFT} the following representation formula for the linear evolution $e^{t\bfL} P_\mathrm{e} \bmf$ in terms of the modified scalar distorted transforms $\calF_{\ell}^{\#}$, $\calF_{\ell,D}^{\#}$, and their adjoints.
\begin{proposition} \label{prop:rep_formula_propagator_modified_dFT}
For any pair of real-valued Schwartz functions $\bm{f} = (f_1,f_2) \in \calS(\bbR) \times \calS(\bbR)$ and for any $t \in \bbR$, we have the following representation formula for the evolution $e^{t\bfL} P_\mathrm{e} \bmf$,
\begin{equation} \label{equ:rep_formula_propagator_modified_dFT}
    \begin{aligned}
        \big( e^{t\bfL} P_\mathrm{e} \bm{f} \big)_1(x) &= \Re \, \biggl( \calF_\ell^{\#,\ast}\Bigl[ e^{i t (\jxi + \ell \xi)} i\jxi^{-1} \Bigl( \calF_{\ell,D}^{\#}[f_1](\xi) - \calF_\ell^{\#}[f_2](\xi) \Bigr) \Bigr](x) \biggr), \\
        \big( e^{t\bfL} P_\mathrm{e} \bm{f} \big)_2(x) &= \Re \, \biggl( \calF_{\ell,D}^{\#,\ast} \Bigl[ e^{i t (\jxi + \ell \xi)} i{\jxi}^{-1} \Bigl( \calF_{\ell,D}^{\#}[f_1](\xi) - \calF_\ell^{\#}[f_2](\xi) \Bigr) \Bigr](x) \biggr).
    \end{aligned}
\end{equation}
\end{proposition}

\section{Distorted Fourier Theory} \label{sec:distorted_fourier_theory}

In the preceding section, we derived an oscillatory integral representation for the evolution $e^{t\bfL_\ell} \Pe$ in terms of the modified scalar distorted transforms introduced in \eqref{eq:modF} and \eqref{eq:modFD}.
In this section we study the resulting (vectorial) distorted Fourier transform associated with the matrix operator $\bfL_\ell$ more systematically. In particular, we establish a Plancherel theorem, Fourier inversion formulas, and Fourier duality between regularities and weights. 
Below, we will use the subscript $\ell$ again in $\bfL_\ell$ and $\bfH_\ell$ to indicate the dependence on $\ell$ explicitly. 
All arguments in this section are not specific to the sine-Gordon case and extend to general matrix operators obtained from linearizing around moving solitons in scalar field theories on the line.
Only in the last Subsection~\ref{subsec:mapping_properties} we give streamlined proofs of several mapping properties of the modified scalar distorted transforms using the explicit formulas for the distorted Fourier basis elements in the sine-Gordon case, but these proofs can be easily replaced by general arguments.

For the reader's orientation, we emphasize again that in the nonlinear analysis in the remainder of this paper, the results from this section will be applied in a moving frame with the spatial coordinate $y := x - q(t)$, as introduced in \eqref{equ:setting_up_moving_frame_coordinate}. 
However, since we hope that the results from this section are of independent interest, we use the more conventional notation $x$ for the spatial coordinate here.

\subsection{Vectorial distorted Fourier transform} \label{subsec:vectorial_dist_FT}

We define the vectorial distorted Fourier transform associated with $\bfL_\ell$ for a pair of real-valued functions $\bmf = (f_1, f_2) \in \calH$ by 
\begin{equation}\label{eq:Tell}
  \mathcal{T}_\ell [\bmf] (\xi):= \begin{bmatrix} \calF^{\#}_{\ell,D}\,\Re  & -\calF^{\#}_\ell \, \Re \end{bmatrix} \bm{f}= \calF^{\#}_{\ell,D}\bigl[\Re f_1\bigr](\xi) - \calF^{\#}_{\ell}\bigl[ \Re f_2\bigr](\xi),
\end{equation}
where the scalar distorted transforms $\calF_\ell^{\#}$ and $\calF_{\ell, D}^{\#}$ are defined in \eqref{eq:modF} and \eqref{eq:modFD}.
The formal $L^2$ adjoint of $\mathcal{T}_\ell$ with respect to the inner product \eqref{eq:L2inner} for $g \in L^2$ is given by
\begin{equation}\label{eq:adjointcalT}
   \calT_\ell^* [g] (x):= \begin{bmatrix} \Re\,\calF^{\#,*}_{\ell,D} \\ -\Re\,\calF^{\#,*}_\ell  \end{bmatrix} g = \begin{bmatrix} \Re\,\calF^{\#,*}_{\ell,D}[g](x) \\ -\Re\,\calF^{\#,*}_\ell[g](x) \end{bmatrix}.
\end{equation}
Correspondingly, we have 
\begin{equation*}
    \langle  \mathcal{T}_\ell [\bmf], g\rangle=  \langle \bmf,  \calT_\ell^* [g] \rangle,
\end{equation*}
where the inner product \eqref{eq:L2inner} is used on the left-hand side and the inner product \eqref{eq:L2L2inner} is used on the right-hand side.

\begin{remark}
Note that the use of $\Re$ in the definition of $\calT_\ell$ is redundant, but we keep it there to make notation more uniform.
Observe that $\calT_\ell$ maps a $2$-vector of real-valued functions to a complex-valued function, which gives consistent dimensions.  
\end{remark}

Using the vectorial Fourier transform $\calT_\ell$ and its adjoint $\calT_\ell^{\ast}$, the representation formula for the evolution from Proposition \ref{prop:rep_formula_propagator_modified_dFT} can be written succinctly as 
\begin{equation}\label{eq:evoFourier}
      e^{t\bfL_\ell}P_{\mathrm{e},\ell}\bmf   = \bfJ^* \calT_\ell^*\Bigl[ e^{i t (\jxi + \ell \xi)} i\jxi^{-1} \calT_\ell(\bmf) \Bigr].
\end{equation}
Setting $t=0$ in the preceding expression, we find 
\begin{equation}\label{eq:Pefiden}
     P_{\mathrm{e},\ell}\bmf  =- \bfJ \calT_\ell^*\Bigl[ i \jxi^{-1} \calT_\ell(\bmf) \Bigr] = \bfJ^* \calT_\ell^*\Bigl[ i \jxi^{-1} \calT_\ell(\bmf)\Bigr],
\end{equation}
where $P_{\mathrm{e},\ell}$ is the projection onto the essential spectrum with respect to $\bfL_\ell$, see \eqref{eq:defPe}. 
Moreover, taking the time derivative in the expansion from Proposition \ref{prop:rep_formula_propagator_modified_dFT} and evaluating the resulting identity at $t=0$, we obtain
\begin{equation}
    \bfL_\ell P_{\mathrm{e},\ell}\bmf = \bfJ \bfH_\ell P_{\mathrm{e},\ell}\bm f = -\bfJ^* \calT_\ell^*\Bigl[ (\jxi+\ell\xi) \jxi^{-1} \calT_\ell(\bmf) \Bigr],
\end{equation}
which implies
\begin{equation}\label{eq:Hfac}
    \bfH_\ell P_{\mathrm{e},\ell} \bmf = \calT_\ell^*\Bigl[ (\jxi+\ell\xi) \jxi^{-1} \calT_\ell(\bmf) \Bigr].
\end{equation}

Next, on the subspace $P_{\mathrm{e},\ell}\calH$ of $\calH$ given by
\begin{equation}
    P_{\mathrm{e},\ell}\calH := \bigl\{ P_{\mathrm{e},\ell}\bmf \, \big| \, \bmf\in\calH=H^1(\bbR)\times L^2(\bbR) \bigr\},
\end{equation}
we define the inner product 
\begin{equation} \label{eq:innerHell}
    \langle \bm{f},\bm{g}\rangle_{\bfH_\ell} := \langle\bfH_\ell \bm{f},\bm{g} \rangle.
\end{equation}
We note that the positive definiteness of \eqref{eq:innerHell} on $P_{\mathrm{e},\ell} \calH$ follows from Lemma~\ref{lem:coerH}.
Additionally, we introduce the short-hand notation
\begin{equation} \label{equ:innerL2ell_definition}
   \langle f,  g \rangle_{L^2_\ell} := \left\langle f, \jxi_\ell \jxi^{-1} g \right\rangle = \Re \int f(\xi) \overline{g(\xi)} \, \jxi_\ell \jxi^{-1} \, \ud \xi,
\end{equation}
where we set
\begin{equation}
    \jxi_\ell := \jxi+\ell \xi.
\end{equation}
Clearly, $\jxi_\ell\sim \jxi$. In terms of the inner products \eqref{eq:innerHell} and \eqref{equ:innerL2ell_definition}, we obtain the following Plancherel theorem for the vectorial Fourier transform $\calT_\ell$.

\begin{lemma}\label{lem:Plancherel}
For any $\bm{f},\bm{g}\in  \mathcal{H}$, it holds that
\begin{equation}\label{eq:planch}
      \langle P_{\mathrm{e},\ell} \bm{f},  \bm{g}\rangle_{\bfH_\ell} = \langle \bfH_\ell \bm{f}, \bm{g} \rangle = \langle \calT_\ell \bmf, \calT_\ell \bm{g} \rangle_{L^2_\ell }.
\end{equation}
\end{lemma}
\begin{proof}
    For any $\bmf\in P_{\mathrm{e},\ell} \calH $ and $\bmg\in\calH$, we have
    \begin{align}
      \langle \bm{f},\bm{g}\rangle_{\bfH_\ell} &= \langle\bfH_\ell \bm{f},\bm{g}\rangle= \Bigl\langle \calT_\ell^* \bigl[ \jxi_\ell \jxi^{-1} \calT_\ell \bmf \bigr], \bm{g} \Bigr\rangle = \langle \calT_\ell \bmf ,\calT_\ell \bm{g} \rangle_{L^2_\ell},
    \end{align}
    where in the second equality above we used \eqref{eq:Hfac}. The identity \eqref{eq:planch} follows.
\end{proof}

We record that by the Plancherel type identity \eqref{eq:planch}, $\calT_\ell$ is injective on $P_{\mathrm{e},\ell}\calH$, whence the range of $\calT_\ell^*$ is $P_\mathrm{e}^*\calH$. It is also straightforward to show that the domain of $\calT_\ell^*$ is $L^2$.

\subsection{Fourier inversion} \label{subsec:fourier_inversion} 
 
 In view of \eqref{eq:Pefiden}, the following inversion formula holds for $\bmf\in P_{\mathrm{e},\ell}\calH$ on the physical side
 \begin{equation} \label{eq:physicalidentity}
    \bmf = - \bfJ \calT_\ell^\ast \bigl[ i \jxi^{-1} \calT_\ell[\bmf] \bigr] = \bfJ^\ast \calT_\ell^*\bigl[ i\jxi^{-1}  \calT_\ell[\bmf] \bigr].
\end{equation}
In this subsection, we establish the inversion formula on the distorted frequency side.
Its derivation is based on the following lemma.

\begin{lemma} \label{lem:imageT}
    The image of $P_{\mathrm{e},\ell}\calH$ under the vectorial Fourier transform $\calT_\ell$ is $L^2$.
\end{lemma}

Before proving Lemma~\ref{lem:imageT}, we present several key results that follow directly from it.

\begin{lemma} \label{lem:transY}
It holds that
\begin{equation}\label{eq:transY}
  \calT_\ell  \bmY_{0}=\calT_\ell  \bmY_{1}=0.
\end{equation}    
\end{lemma}
\begin{proof}
Let $\bmg\in P_{\mathrm{d},\ell}\calH$. Then we have $\bfH_\ell \bmg =-\frac{1}{\langle \bmY_{0}, \bfJ\bmY_{1}\rangle}  \langle \bm{g}, \bfJ\bmY_{0}\rangle  \bfJ\bmY_{0}$. 
By \eqref{eq:planch} and Lemma~\ref{lem:decom} we conclude for all $\bmf \in P_{\mathrm{e},\ell} \calH $ that
    \begin{align}
      \langle \calT_\ell \bmf,  \calT_\ell \bm{g} \rangle_{L^2_\ell }=\langle \bm{f},\bm{g}\rangle_{\bfH_\ell} =\langle\bfH_\ell \bm{f},\bm{g}\rangle=\langle \bm{f},\bfH_\ell\bm{g}\rangle=-\frac{\langle \bm{g}, \bfJ\bmY_{0}\rangle }{\langle \bmY_{0}, \bfJ\bmY_{1}\rangle}\langle \bmf, \bfJ\bmY_{0}\rangle=0.
\end{align}
The range of $\calT_\ell$ is $L^2$ by Lemma \ref{lem:imageT}, whence $\calT_\ell \bmg=0$, and the assertion follows.
\end{proof}

An immediate consequence of Lemma~\ref{lem:transY} is that we can freely put $P_{\mathrm{e},\ell}$ in front of functions before taking the vectorial Fourier transform.

\begin{corollary}\label{cor:TellP}
For any $\bmf\in\calH$ it holds that
\begin{equation}
    \calT_\ell[\bmf](\xi)=    \calT_\ell[P_{\mathrm{e},\ell}\bmf](\xi).
\end{equation}
\end{corollary}
\begin{proof}
By Lemma \ref{lem:decom} we can write
\begin{align*}
    \bmf &= P_{\mathrm{e},\ell} \bmf + P_{\mathrm{d},\ell} \bmf = P_{\mathrm{e},\ell} \bmf+  \frac{1}{\langle \bmY_{0}, \bfJ\bmY_{1}\rangle} \langle \bmf, \bfJ\bmY_{1}\rangle \bmY_{0}-\frac{1}{\langle \bmY_{0}, \bfJ\bmY_{1}\rangle}  \langle \bmf, \bfJ\bmY_{0}\rangle  \bmY_{1}.
\end{align*}
Applying the vectorial Fourier transform and using Lemma~\ref{lem:transY}, we obtain
\begin{align*}
     \calT_\ell[\bmf](\xi)& = \calT_\ell[P_{\mathrm{e},\ell}\bmf](\xi)+ \frac{1}{\langle \bmY_{0}, \bfJ\bmY_{1}\rangle} \langle \bmf, \bfJ\bmY_{1}\rangle \calT_\ell[\bmY_{0}](\xi)-\frac{1}{\langle \bmY_{0}, \bfJ\bmY_{1}\rangle}  \langle \bmf, \bfJ\bmY_{0}\rangle  \calT_\ell[\bmY_{1}](\xi)=\calT_\ell[P_{\mathrm{e},\ell}\bmf](\xi),
\end{align*}
as desired.
\end{proof}

Next, we conclude the Fourier inversion formula on the distorted frequency side.

\begin{lemma} \label{lem:inversionXi} 
For all $h \in L^2$ it holds that
\begin{equation}\label{eq:xiinversion}
    \mathcal{T}_\ell \Bigl[ \bfJ \mathcal{T}^*_\ell \bigl[ i\jxi^{-1} h(\xi) \bigr] \Bigr] = h.
\end{equation}
\end{lemma}
\begin{proof}
 Applying the vectorial distorted Fourier transfrom $\calT_\ell$ to both sides of \eqref{eq:physicalidentity}, we obtain
 \begin{equation}
   \mathcal{T}_\ell\Bigl[ \bfJ^* \calT^*_\ell \bigl[ i \jxi^{-1} \mathcal{T}_\ell[\bmf] \bigr] \Bigr] =  \mathcal{T}_\ell[\bmf].
 \end{equation}
 By Lemma \ref{lem:imageT}, we can write any $h \in L^2$ as $h=\calT_\ell [\bmf]$ for some $\bmf\in P_{\mathrm{e},\ell}\calH$. It follows that 
 \begin{equation}
     \mathcal{T}_\ell\Bigl[ \bfJ^* \mathcal{T}^*_\ell \bigl[ i\jxi^{-1} h(\xi) \bigr] \Bigr] = h,
 \end{equation}
 as desired.
\end{proof}

Combing Proposition \ref{prop:rep_formula_propagator_modified_dFT} and the inversion formula \eqref{eq:xiinversion} on the distorted frequency side, we conclude the following identity for the evolution on the distorted frequency side.
\begin{lemma} \label{lem:propagator_on_dist_Fourier_side}
    Fix $\ell \in (-1,1)$. Denote by $P_{\mathrm{e},\ell}$ the projection to the essential spectrum with respect to the linearized operator $\bfL_\ell$. For any pair of real-valued Schwartz functions $\bmF = (F_1, F_2)$ and for any $t \in \bbR$, we have 
    \begin{equation}
        \calF_{\ell,D}^{\#}\Bigl[ \bigl( e^{\pm t\bfL_\ell} P_{\mathrm{e},\ell} \bmF \bigr)_1 \Bigr](\xi) - \calF_{\ell}^{\#}\Bigl[ \bigl( e^{\pm t\bfL_\ell} P_{\mathrm{e},\ell} \bmF \bigr)_2 \Bigr](\xi) = e^{\pm i t (\jxi + \ell \xi)} \Bigl( \calF_{\ell,D}^{\#}[F_1](\xi) - \calF_{\ell}^{\#}[F_2](\xi) \Bigr). 
    \end{equation}
\end{lemma}
\begin{proof}
 By Proposition~\ref{prop:rep_formula_propagator_modified_dFT} we have 
 \begin{align}
      e^{\pm t\bfL_\ell} P_{\mathrm{e},\ell} \bmF 
      &= \bfJ^* \calT_\ell^*\Bigl[ e^{i t (\jxi + \ell \xi)} i\jxi^{-1} \calT_\ell[\bmF] \Bigr] \\
      &= \bfJ^* \calT_\ell^*\Bigl[ e^{\pm i t (\jxi + \ell \xi)} i\jxi^{-1} \bigl( \calF_{\ell,D}^{\#}[F_1](\xi) - \calF_{\ell}^{\#}[F_2](\xi) \bigr) \Bigr].
 \end{align}
 Invoking Lemma~\ref{lem:inversionXi}, we conclude
 \begin{align}
    &\calF_{\ell,D}^{\#}\Bigl[ \bigl( e^{\pm t\bfL_\ell} P_{\mathrm{e},\ell} \bmF \bigr)_1 \Bigr](\xi) - \calF_{\ell}^{\#}\Bigl[ \bigl( e^{\pm t\bfL_\ell} P_{\mathrm{e},\ell} \bmF \bigr)_2 \Bigr](\xi) 
    = \calT_\ell\Bigl[ e^{\pm t\bfL_\ell} P_{\mathrm{e},\ell} \bmF \Bigr](\xi) \\
    &= \calT_\ell\biggl[ \bfJ^\ast \calT_\ell^\ast \Bigl[ i\jxi^{-1} e^{\pm i t(\jxi+\ulell\xi)} \bigl( \calFulellDsh[F_1](\xi) - \calFulellsh[F_2](\xi) \bigr) \Bigr] \biggr] \\ 
    &= e^{\pm i t (\jxi + \ell \xi)}\Bigl( \calF_{\ell,D}^{\#}\bigl[ F_1 \bigr](\xi) - \calF_{\ell}^{\#}\bigl[F_2\bigr](\xi) \Bigr).
 \end{align}
 The desired identity follows.
 \end{proof}
Now we prove Lemma \ref{lem:imageT}. It is established via a continuation argument. The idea is that as the velocity varies, the dimension of the kernel of $\calT_\ell^*$ or the co-kernel of $\calT_\ell$ is an integer that cannot be continuously varied.
\begin{proof}[Proof of Lemma \ref{lem:imageT}] 
 First of all, note that the desired result holds for $\ell=0$. This is the standard case passing from a solution to the Klein-Gordon equation to a complexified solution to the half-Klein-Gordon equation.
 
More precisely, when $\ell=0$, by a standard duality argument, it suffices to show that the adjoint $\calT_\ell^*$ defined in \eqref{eq:adjointcalT}, is injective, i.e., that
\begin{equation}
       \calT_0^*[g] := \begin{bmatrix} \Re\,\calF^{\#,*}_{0,D} \\ -\Re\,\calF^{\#,*}_0  \end{bmatrix} g = 0 \quad \Rightarrow \quad g=0. 
\end{equation}
Explicitly, it is reduced to show that
\begin{align}\label{eq:twozeros}
    \Re \int e^{\#}_0(x,\xi) g(\xi)\,\ud \xi=\Re \int  e^{\#}_0(x,\xi) i \jxi g(\xi)\,\ud \xi=0 
\end{align}implies $g=0$.
Note that
\begin{align}
    0=2\Re \int e^{\#}_0(x,\xi) g(\xi)\,\ud \xi= &\int e^{\#}_0(x,\xi) g(\xi)\,\ud \xi+ \int \overline{e^{\#}_0(x,\xi) g(\xi)}\,\ud \xi\\
    &=\int e^{\#}_0(x,\xi) g(\xi)\,\ud \xi+ \int e^{\#}_0(x,\xi)\overline{ g(-\xi)}\,\ud \xi\\
    &= \int e^{\#}_0(x,\xi)( g(\xi)+\overline{ g(-\xi)})\,\ud \xi
\end{align}
Since here  $e_0(x,\xi)=e^{\#}_0(x,\xi)\text{sgn}(\xi)$ is the distorted Fourier basis associated with a scalar self-adjoint operator $L$ defined in \eqref{eq:Hk}, it follows that
\begin{equation}
    g(\xi)+\overline{ g(-\xi)}=0.
\end{equation}By a similar argument, the second part of \eqref{eq:twozeros} implies
\begin{equation}
      \jxi g(\xi)-\jxi \overline{ g(-\xi)}=0.
\end{equation}
Therefore, we conclude $g=0$.

Next we study $\ell\neq 0$. Without loss of generality, we consider positive velocities.
We set
\begin{equation}\label{eq:defbeta}
    \beta := \sup \, \bigl\{ 0 \leq 0<1 \, \big| \, \text{ran}( \calT_\ell) = L^2 \bigr\}.
\end{equation}
To achieve the desired result, we claim that $\beta=1$. Assume $\beta<1$. Without loss of generality, one can further assume that  $\text{ran}( \calT_\beta)= L^2$. This can be justified by a limiting argument or  otherwise we can take $\tilde{\beta}<\beta$ such that  $\ell-\tilde{\beta}<\frac{3}{2}\epsilon$ to replace the role of $\beta$ in the later arguments.

In the next 6 steps, we will show that for $0<\epsilon\ll1$ small enough $\text{ran}( \calT_\ell)= L^2$ for $\ell=\beta+\epsilon<1$.

\smallskip

\noindent {\bf Step 1:} We claim there is  a bijection $\mathcal{M}$ such that $\mathcal{M} P_{\mathrm{e},\ell} \calH= P_{\mathrm{e},\beta}\calH$ if $|\ell-\beta|\leq\epsilon\ll1$. This  is verified by taking the difference of the generalized kernels in our settings, for $\bmf\in P_{\mathrm{e},\ell}\calH$, 
\begin{equation}
\mathcal{M}\bmf:=\bmf -P_{\mathrm{d},\beta}\bmf=\bmf -(P_{\mathrm{d},\beta}-P_{\mathrm{d},\ell})\bmf.
\end{equation}
The last term is small given $|\ell-\beta|\ll 1$. We can always invert $\mathcal{M}$.  The first claim is valid.

We can also abstractly define the adjoint of $\mathcal{M}$ from $P_{\mathrm{e},\beta}^* (H^{-1} \times L^2)\rightarrow P_{\mathrm{e},\ell}^* (H^{-1} \times L^2) $.

\smallskip

\noindent {\bf Step 2:} We claim $\bfL_\ell^{-1}$ is bounded  from    $P_{\mathrm{e},\ell} \calH$ to $P_{\mathrm{e},\ell}(H^2\times H^1)$  for all $|\ell|<1$.

Now suppose $\bm{g}=\binom{g_1}{g_2}\in P_{\mathrm{e},\ell}\calH$. We  first want to solve  
$\bfL_\ell\bm{\psi} = \bm{g}$
uniquely for $\bm{\psi}\in P_{\mathrm{e},\ell}(H^2\times H^1)$. To achieve our goal, explicitly, we need to solve the following system:
\begin{align}\label{eq:psisyst}
  \begin{cases}
 \ell\px \psi_1+\psi_2 =g_1
  \\
  \\
 -\bigl( - \px^2 +V_\ell + 1 \bigr)\psi_1 +\ell\px \psi_2=g_2.
\end{cases}
\end{align}
Plugging the first equation into the second equation and eliminating $g_2$, we get
\begin{equation}
   L_\gamma \psi_1:=-(1-\ell^2)\partial_x^2+V_\ell\psi_1+\psi_1=g_2-\ell \partial_x g_1.
\end{equation}
Note that $\langle \bm g, \bfJ\bmY_{0}\rangle=0$ implies $\langle g_2-\ell \partial_x g_1, \gamma K'(\gamma x)\rangle_{L^2}=0.$ Therefore, we can always invert the linear operator $L_\ell+\ell^2\px^2$ whose kernel is given by $ \gamma K'(\gamma x)$ to find the unique $\psi_1$ such that $\langle \psi_1, \gamma K'(\gamma x)\rangle_{L^2}=0$ via the standard elliptic theory, and the regularity also follows. After finding $\psi_1$, from the first equation of the system \eqref{eq:psisyst}, one can find $\psi_2$ and this gives $\bm{\psi}$. We still need to check $\bmpsi\in P_{\mathrm{e},\ell}(H^2\times H^1)$. To show this, we write $\bm{\psi}=P_{\mathrm{e},\ell}\bm{\psi}+P_{\mathrm{d},\ell}\bm{\psi}.$ For the discrete part, explicitly, one has
\begin{equation}
       P_{\mathrm{d},\ell} \bm{\psi}= \frac{1}{\langle \bmY_{0}, \bfJ\bmY_{1}\rangle} \langle \bm{\psi}, \bfJ\bmY_{1}\rangle \bmY_{0}-\frac{1}{\langle \bmY_{0}, \bfJ\bmY_{1}\rangle}  \langle \bm{\psi}, \bfJ\bmY_{0}\rangle  \bmY_{1}.
\end{equation}
Applying $\bfL$ to $ P_{\mathrm{d},\ell} \bm{\psi}$, we get
\begin{equation} 
       \bfL\left(P_{\mathrm{d},\ell} \bm{\psi}\right)= -\frac{1}{\langle \bmY_{0}, \bfJ\bmY_{1}\rangle}  \langle \bm{\psi}, \bfJ\bmY_{0}\rangle  \bmY_{0}.
\end{equation}But $  \bfL  \bm{\psi}=  \bm{g}\in P_{\mathrm{e},\ell}\calH$ which implies $  \langle \bm{\psi}, \bfJ\bmY_{0}\rangle =0$.  Since $\langle \psi_1, \gamma K'(\gamma x)\rangle=0$ by construction, it follows that $ \langle \bm{\psi}, \bfJ\bmY_{1}\rangle  =0$.  We conclude that $\bm{\psi}=P_{\mathrm{e},\ell}\bm{\psi}.$
 \smallskip
 
\noindent {\bf Step 3:} We also claim that $\bfL^{-1}_{\ell}$ boundedly maps  $P_{\mathrm{e},\ell} (L^2\times H^{-1})$ to $P_{\mathrm{e},\ell} \calH$.

This is achieved by a similar argument as the one above. Following the same steps above
\begin{equation}
    L_\gamma \psi_1= -(1-\ell^2)\partial_x^2-V_\ell\psi_1+\psi_1=g_2-\ell \partial_x g_1.
\end{equation}holds in the sense of $H^{-1}$. In other words, one has
\begin{equation}
    \langle   L_\gamma \psi_1, h\rangle = \langle g_2-\ell \partial_x g_1, h\rangle
\end{equation} for all $h\in \mathrm{P}_{c,\ell} H^1$ where $\mathrm{P}_{c,\ell}$ denotes the projection on the continuous spectrum of the self-adjoint operator $ L_\gamma$. Due to the coercivity on  the continuous spectrum of $ L_\gamma$, see Corollary \ref{cor:coerLell}, $\langle ( L_\gamma \cdot, \cdot\rangle$ defines an inner product over $\mathrm{P}_{c,\ell} H^1$. By construction, $\langle g_2-\ell \partial_x g_1, \cdot \rangle$ is a bounded linear functional over $\mathrm{P}_{c,\ell} H^1$, the desired result follows from the Riesz representation. $\bm{\psi}=P_{\mathrm{e},\ell}\bm{\psi}$ is obtained via the same argument as in {\bf Step 2}.

\smallskip

\noindent {\bf Step 4:} Now we consider the operator $\calQ=\calT_\ell+\eta \calT_\beta\mathcal{M}$ for $\ell=\beta+\epsilon$ for $\epsilon,\eta>0 $ small enough where $\calM$ is from  {\bf Step 1}.

We first note that $\left(\bfJ^*\calQ^{*}\frac{\jxi_\ell}{\jxi}\calQ\right)^{-1}$ is  bounded from  $P_{\mathrm{e},\ell} \calH$ to $P_{\mathrm{e},\ell}(H^2\times H^1)$  and from  $P_{\mathrm{e},\ell} (L^2\times H^{-1})$ to $P_{\mathrm{e},\ell} \calH$ for $\eta$ small enough since $\bfL_\ell^{-1}=(\bfJ^*\calT^*_\ell\frac{\jxi_\ell}{\jxi}\calT_\ell)^{-1}$ satisfies the desired properties by {\bf Step 2} and {\bf Step 3},  and other terms from $\bfJ^*\calQ^{*}\frac{\jxi_\ell}{\jxi}\calQ$ are small provided that $\eta$ is small enough. Using the smallness of $\eta$ and  the  injectivity of $\calT_\ell$ which follows \eqref{eq:planch}, we can again show that $\calQ$ is injective as well. 

Next, we claim that $\text{ran}\calQ$ is $L^2$ for $\eta$ small enough. Recall that by construction, \eqref{eq:defbeta},  the range of $\eta \calT_\beta$ is $L^2$. To show the desired result, we argue by contradiction. Suppose $\text{ran}\calQ\neq L^2$. There one can find $h\in L^2$, $h\neq0$ such that for all $\bmf\in P_{\mathrm{e},\ell}\calH$.
\begin{equation}
    \langle h, \calQ \bm f\rangle_{L^2_\ell}=0.
\end{equation}
If $h\in \text{ran} \calT_\ell$, clearly, the inner product above will not be zero given $\eta$ is small enough $\forall \bmf \in P_{\mathrm{e},\ell}\calH$ by the Plancherel identity in Lemma \ref{lem:Plancherel}. If $h\notin \text{ran} \calT_\ell$, without loss of generality, we can further assume that $h\in (\text{ran} \calT_\ell)^{\perp}$.  Then one has
\begin{equation}
    \langle h, \calQ \bm f\rangle_{L^2_\ell}= \eta \langle h, \calT_\beta \mathcal{M} \bm f\rangle_{L^2_\ell}
\end{equation}which will be nonzero for some $ \bmf \in P_{\mathrm{e},\ell}\calH=\mathcal{M}^{-1}P_{\mathrm{e},\beta}\calH$ provided $|\ell-\beta|\leq \epsilon\ll 1$ small enough since $\calT_\beta $ surjective by assumption.

 \noindent {\bf Step 5:} We claim that $(\bfJ^* \calQ^*\frac{\jxi_\ell}{\jxi})^{-1}$ is bounded from $P_{\mathrm{e},\ell}( L^2 \times H^{-1})$ to $L^2$. Denote $\bfJ^*\calQ^{*}\frac{\jxi_\ell}{\jxi}\calQ=:\mathcal{W}$. We claim that $$(\bfJ^* \calQ^*\frac{\jxi_\ell}{\jxi})^{-1}=\calQ\mathcal{W}^{-1}.$$Note that by the mapping properties of $\mathcal{W}$, the operator $\calQ\mathcal{W}^{-1}$ is well-defined from $P_{\mathrm{e},\ell}( L^2 \times H^{-1})$ to $L^2$.  To show our claim, we clearly note that $\bfJ^* \calQ^*\frac{\jxi_\ell}{\jxi} \circ(\calQ\mathcal{W}^{-1})=I$ in $P_{\mathrm{e},\ell}( L^2 \times H^{-1})$.  Now check the identity on the frequency side. For any $h\in L^2$,
\begin{align}
    \langle\calQ \bmf , \frac{\jxi_\ell}{\jxi} \calQ\mathcal{W}^{-1} (\bfJ^* \calQ^* )h\rangle&= \langle   \bmf , \calQ^* \frac{\jxi_\ell}{\jxi}  \calQ\mathcal{W}^{-1} (\bfJ^* \calQ^* )h\rangle\\
    &= \langle \bmf ,-\bfJ^* (\bfJ^* \calQ^* )h\rangle\\
    &=\langle   \bmf , \calQ^*h\rangle\\
    &= \langle  \calQ \bmf , h \rangle.
\end{align}Since $\calQ$ is surjective by {\bf Step 4}, one has $\frac{\jxi_\ell}{\jxi}\calQ\mathcal{W}^{-1} (\bfJ^* \calQ^* )=I $. Therefore, $(\bfJ^* \calQ^*\frac{\jxi_\ell}{\jxi})^{-1}=\calQ\mathcal{W}^{-1}$.

\smallskip

 \noindent {\bf Step 6: }
By construction, $\calQ\mathcal{W}^{-1}$ is   bounded from $P_{\mathrm{e},\ell}( L^2 \times H^{-1})$ to $L^2$ which shows that $\bfJ^* \calQ^* \frac{\jxi_\ell}{\jxi}$  has a bounded inverse. We note that
\begin{equation}
    \bfJ^* \calQ^*\frac{\jxi_\ell}{\jxi}= \bfJ^* \calT_\ell^*\frac{\jxi_\ell}{\jxi}+\eta \bfJ^* \mathcal{M}^* \calT^*_\beta \frac{\jxi_\ell}{\jxi}
\end{equation}Since $\bfJ^* \mathcal{M}^*\calT^*_\beta \frac{\jxi_\ell}{\jxi}$ is bounded, then for $\eta$ small enough (depending on $\beta<1$), from the mapping properties of $ \bfJ^* \calQ^* \frac{\jxi_\ell}{\jxi}$ and the Neumann series, it follows that  $(\bfJ^* \calT_\ell^*)^{-1}$ is bounded from $P_{\mathrm{e},\ell} (L^2 \times H^{-1})$ to $L^2$ since $|\frac{\jxi_\ell}{\jxi}|\sim 1$. In particular, $\calT_\ell^*$ is injective, which implies that the range of $\calT_\ell$ is $L^2$.  This is a contradiction with the definition of $\beta$. 

\bigskip

Combing {\bf Step 1} to {\bf Step 6}, we conclude that $\forall |\ell|<1$, the range of $\calT_\ell$ is $L^2$.
\end{proof}

\subsection{Higher order Sobolev estimates} \label{subsec:higher_order_sobolev}

Suppose $\bmf\in H^{k+1}\times H^{k}$ and $P_{\mathrm{e},\ell}\bmf=\bmf$. We would like to establish relations between the $ H^{k+1}\times H^{k}$ norm of $\bmf$ and the $L^2$ norm of $\jxi^k \calT_\ell[\bmf] (\xi)$ as in the standard Fourier theory.
\begin{lemma} \label{lem:higherorder}
Let $k \geq 0$ be an integer.
Suppose $\bmf\in H^{k+1}\times H^{k}$ and $P_{\mathrm{e},\ell}\bmf=\bmf$. Then it holds that
\begin{equation}\label{eq:highcomp}
    \|\bmf\|_{H^{k+1}\times H^{k}} \sim \|\jxi^k \calT_\ell[\bmf](\xi)\|_{L^2}.
\end{equation}
\end{lemma}
\begin{proof}
We start with  $k=0$. By Lemma \ref{lem:Plancherel}, for $P_{\mathrm{e},\ell}\bmf=\bmf\in P_{\mathrm{e},\ell}\calH$ by the coercivity of $\bfH_\ell$ over $P_{\mathrm{e},\ell}\calH$, see Lemma \ref{lem:coerH},
\begin{equation} 
      \| \bmf\|_{\calH}^2 \sim \langle \bm{f},\bm{f}\rangle_{\bfH_\ell} =\langle\bfH_\ell \bm{f},\bm{f} \rangle=\langle \calT_\ell \bmf ,\calT_\ell \bm{f}\rangle_{L^2_\ell}\sim \|  \calT_\ell (\bmf) (\xi)\|^2_{L^2}.
\end{equation}
Consider $k=1$. We observe that by taking the time derivative in \eqref{eq:evoFourier} three times and evaluating the resulting formula at $t=0$, it follows
\begin{equation}
    \bfL_\ell^3 \bmf= \bfJ^* \calT_\ell^*[ \jxi_\ell^3\jxi^{-1} \calT_\ell(\bmf)].
\end{equation}Recalling $\bfL_\ell=\bfJ\bfH_\ell$, we write
\begin{equation}
    \bfJ\bfL_\ell^3\bmf=\bfH_\ell\bfJ^*\bfH_\ell\bfL_\ell\bmf= (\bfL_\ell)^*\bfH_\ell\bfL_\ell\bmf=  \calT_\ell^*[ \jxi_\ell^3\jxi^{-1} \calT_\ell(\bmf)].
\end{equation}
Taking the inner product with $\bmf$ with the formula above, we get
\begin{align*}
   \langle(\bfL_\ell)^*\bfH_\ell\bfL_\ell\bmf  ,\bmf\rangle   &= \langle \bfL_\ell\bmf,\bfL_\ell\bmf\rangle_{\bfH_\ell}=\langle \calT_\ell^*[\jxi_\ell^3\jxi^{-1} \calT_\ell(\bmf)],\bm{f}\rangle= \langle\jxi_\ell\calT_\ell \bmf ,\jxi_\ell \calT_\ell \bm{f}\rangle_{L^2_\ell}.
\end{align*}
Note that $P_{\mathrm{e},\ell} \bfL_\ell\bmf=\bfL_\ell\bmf$ since the essential subspace is invariant under $\bfL_\ell$. By the coercivity of $\bfH_\ell$, Lemma \ref{lem:coerH}, applied to $\bfL_\ell \bmf$ and expressions above, as in the case $k=0$, and $\jxi_\ell\sim\jxi$,  one concludes
\begin{equation}\label{eq:k=1}
    \| \bfL_\ell \bmf \|_{\calH} \sim \|\jxi \calT_\ell \bmf\|_{L^2}.
\end{equation}
Write out $\bfL_\ell \bmf$ more explicitly, 
\begin{equation}
    \bm{h}:=\bfL_\ell \bmf=  \begin{bmatrix} \ell \px f_1 +f_2 \\ - \bigl( - \px^2 - 2 \sech^2(\gamma x) + 1 \bigr) f_1+ \ell \px f_2 \end{bmatrix}
\end{equation}
and set
\begin{equation}
    h_1:=  \ell \px f_1 +f_2,\,\,\,h_2:= - \bigl( - \px^2 - 2 \sech^2(\gamma x) + 1 \bigr) f_1+ \ell \px f_2.
\end{equation}
From \eqref{eq:k=1}, we know that 
\begin{equation}
    \|h_1\|_{H^1}+\|h_2\|_{L^2} \sim  \|\jxi \calT_\ell \bmf\|_{L^2}\lesssim \|\bmf\|_{H^2 \times H^1}.
\end{equation}
From expressions of $h_1$ and $h_2$, the equation for $f_1$ can be rewritten as
\begin{equation}
    - L_\gamma f_1= - \bigl( - \px^2 - 2 \sech^2(\gamma x) + 1 \bigr) f_1-\ell^2\px^2 f_1=-\ell\px h_1+h_2.
\end{equation}
Note that since $\bm{h}=P_{\mathrm{e},\ell} \bm{h}$, $(-\ell\px h_1+h_2)\perp \ker ( L_\gamma).$  By the standard elliptic theory, we can resolve $f_1$ in terms of $-\ell\px h_1+h_2\in L^2$ with
\begin{equation}
    \|f_1\|_{H^2} \sim \|(-\ell\px h_1+h_2)\|_{L^2}.
\end{equation}
Finally, from the expression for $h_2$ and $f_1\in H^2$, we conclude that
\begin{equation}
    \| \ell \px f_2 \|_{L^2} \lesssim \|h_2\|_{L^2}+ \|f_1\|_{H^2} \lesssim \|h_1\|_{H^1}+\|h_2\|_{L^2}.
\end{equation}
From the expression   for $h_1\in H^1$ and $f_1\in H^2$, it follows
\begin{equation}
    \|  f_2 \|_{H^1} \lesssim \|h_1\|_{H^1}+ \|\px f_1\|_{H^1} \lesssim \|h_1\|_{H^1}+\|h_2\|_{L^2}.
\end{equation}
Therefore, we conclude that
\begin{equation}
     \|f_1\|_{H^2}+   \|f_2\|_{H^1} \lesssim  \|h_1\|_{H^1}+\|h_2\|_{L^2}\lesssim \|\jxi \calT_\ell \bmf\|_{L^2}.
\end{equation}This finished the proof for $k=1$.

One can perform an induction argument to show similar arguments hold for $k\in\mathbb{N}$. For general $k$, we take $2k+1$ times $t$ derivative of $e^{t\bfL_\ell}\bmf$ and evaluate the resulting expression at $t=0$
\begin{equation}
     ((\bfL_\ell)^*)^k\bfH_\ell(\bfL_\ell)^k\bmf=  \calF_\ell^*[ \jxi_\ell^{2k+1} \jxi^{-1}\calT_\ell(\bmf)].
\end{equation}By the same argument above, we conclude that
\begin{equation}
    \| \bfL_\ell^k \bmf \|_{\calH} \sim \|\jxi^k \calT_\ell \bmf\|_{L^2}.
\end{equation}
We claim that
\begin{equation}
       \| \bfL_\ell^k \bmf \|_{\calH}\sim  \| \bfL_\ell^{k-1} \bmf \|_{H^2\times H^1}\sim \|\bmf\|_{H^{k+1}\times H^k}.
\end{equation}
The basic $k=0,1$ cases are shown above. Suppose the desired result holds for $k=m$. Consider $k=m+1$
\begin{equation}
       \| \bfL_\ell^{m+1} \bmf \|_{\calH}= \| \bfL_\ell^m (\bfL_\ell  \bmf )\|_{\calH}.
\end{equation}
By the induction hypothesis, we get
\begin{equation}
       \| \bfL_\ell^{m+1} \bmf \|_{\calH}= \| \bfL_\ell^m (\bfL_\ell  \bmf )\|_{\calH}\sim \|\bfL_\ell \bmf\|_{H^{m+1}\times H^m} .
\end{equation}
Setting $\bm{h}:=\bfL_\ell  \bmf$, then one has
\begin{equation}
    \ell \px f_1 +f_2=:  h_1 \in H^{m+1},\,\,\,- \bigl( - \px^2 - 2 \sech^2(\gamma x) + 1 \bigr) f_1+ \ell \px f_2=:h_2 \in H^m.
\end{equation}
From expressions of $h_1$ and $h_2$, the equation for $f_1$ can be rewritten as
\begin{equation}\label{eq:eqf1}
    -L_\ell f_1=- \bigl( - \px^2 - 2 \sech^2(\gamma x) + 1 \bigr) f_1-\ell^2\px^2 f_1=-\ell\px h_1+h_2.
\end{equation}
Note that since $\bm{h}=P_{\mathrm{e},\ell} \bm{h}$, $(-\ell\px h_1+h_2)\perp \ker L_\gamma$ and $(-\ell\px h_1+h_2)\in H^m$.  Note that from the standard elliptic theory, we can resolve $f_1$ in terms of $(-\ell\px h_1+h_2)\in H^m$ with
\begin{equation}\label{eq:esthigh}
    \|f_1\|_{H^{m+2}} \sim \|(-\ell\px h_1+h_2)\|_{H^m}.
\end{equation}
More precisely, to show the estimate \eqref{eq:esthigh} above, we  commute the equation for $f_1$, \eqref{eq:eqf1}, with  $\sqrt{\mathrm{P}_{c,\ell} L_\gamma}$ where  $\mathrm{P}_{c,\ell}$ denotes the projection on the continuous spectrum of the self-adjoint operator $ L_\gamma$. The higher order regularity estimates follow from the following estimate: from the $L^2$ mapping properties of Jost functions or even more generally by the boundedness of wave operators, see Weder \cite{Weder} and Galtbayar-Yajima \cite{GaYa}, one has for any $n\in\mathbb{N}$
\begin{equation}
    \|\mathrm{P}_{c,\ell} g\|_{H^n}\sim \sum_{j=0}^n\left\lVert \left(\sqrt{\mathrm{P}_{c,\ell}( L_\gamma})\right)^jg \right\rVert_{L^2}.
\end{equation}
From the equations of $h_1$ and $h_2$, we can obtain that
\begin{equation}
    \|f_2\|_{H^{m+1}} \lesssim \|f_1\|_{H^{m+2}} + \|h_1\|_{H^{m+1}}\lesssim   \|\bm{h}\|_{H^{m+1}\times H^m} =\|\bfL_\ell \bmf\|_{H^{m+1}\times H^m}.
\end{equation}
It follows that
\begin{equation}
      \| \bfL_\ell^{m+1} \bmf \|_{\calH}\sim \| \bmf\|_{H^{m+2}\times H^{m+1}}.
\end{equation}
Hence, we conclude for $k\in \mathbb{N}$ and $\bmf=P_{\mathrm{e},\ell}\bmf$ that
\begin{equation}
   \| \bmf\|_{H^{k+1}\times H^{k}} \sim \|\jxi^k \calT_\ell \bmf\|_{L^2},
\end{equation}
as desired.
\end{proof}

\begin{remark}
The generalized kernel of $\bfL_\ell$ is spanned by $\bmY_0$ and $\bmY_1$. The orthogonality to $\bmY_1$ ensures the coercivity of $\bfH_\ell$, see Lemma \ref{lem:coerH}, which makes $\langle\cdot, \cdot\rangle_{\bfH_\ell}$ a well-defined inner product. The orthogonality to $\bmY_0$ ensures the invertibility of $\bfL_\ell$.  $\bmY_0$ and $\bmY_1$ together produce an invariant subspace of $\bfL_\ell$.
\end{remark}

\subsection{Mapping properties of transforms} \label{subsec:mapping_properties}

In this subsection we establish mapping properties for the scalar transforms $\calF_\ell^{\#}$ and $\calF_{\ell,D}^{\#}$ as well as some related bounds, which will be needed for the nonlinear analysis in this paper.
Below we give streamlined proofs that exploit the explicit formulas for the distorted Fourier basis elements in the case of the sine-Gordon model, but these proofs can easily be replaced by robust arguments that apply to general moving soliton settings for many other scalar field theories on the line.

We recall from \eqref{equ:modified_dFT_basis_succinct} the explicit form of the distorted Fourier basis elements for the sine-Gordon model
\begin{equation} \label{eq:ellsplit}
    e_\ell^\#(x,\xi) = \frac{1}{\sqrt{2\pi}} m_\ell^\#(\gamma x,\xi) e^{i x \xi}
\end{equation}
with
\begin{equation} \label{eq:mImII}
\begin{aligned}
    m_\ell^\#(\gamma x,\xi) = \frac{\gamma(\xi+\ell\jxi) + i \tanh(\gamma x)}{|\gamma(\xi+\ell\jxi)| - i} 
    &=\frac{\gamma(\xi+\ell\jxi) }{|\gamma(\xi+\ell\jxi)| - i}+\tanh(\gamma x)\frac{ i }{|\gamma(\xi+\ell\jxi)| - i} \\
    &=: m_{\mathrm{I}}(\xi)+\tanh(\gamma x) m_{\mathrm{II}}(\xi).
\end{aligned}
\end{equation}
We also recall that
\begin{equation} \label{eq:Dsell}
    (D^\ast  e_\ell^\#)(x,\xi) = \frac{1}{\sqrt{2\pi}} \Bigl( \bigl(-\ell \px + i (\jxi + \ell \xi) \bigr) m_\ell^\#(\gamma x,\xi) \Bigr) e^{i x \xi} + \frac{1}{\sqrt{2\pi}} m_\ell^\#(\gamma x,\xi) (-i\ell \xi) e^{i x \xi}.
\end{equation}
Note that both $m_{\mathrm{I}}(\xi)$ and $m_{\mathrm{II}}(\xi)$ in \eqref{eq:mImII} are Lipschitz in $\xi$. By direct computation, we have
\begin{equation}\label{eq:mi}
    |m_{\mathrm{I}}(\xi)| \lesssim 1, \qquad |\partial_\xi m_{\mathrm{I}}(\xi)|\lesssim \frac{1}{\jxi}
\end{equation}
and
\begin{equation}\label{eq:mii}
    |m_{\mathrm{II}}(\xi)|\lesssim \frac{1}{\jxi}, \qquad|\partial_\xi m_{\mathrm{II}}(\xi)|\lesssim \frac{1}{\jxi^2}.
\end{equation}
By inspection of the preceding formulas \eqref{eq:ellsplit}, \eqref{eq:mImII}, \eqref{eq:Dsell} \eqref{eq:mi}, \eqref{eq:mii}, we infer the following bounds.
\begin{lemma}\label{lem:boundseDe}
For all integers $j \geq 0$ and for $k = 0, 1$, it holds that
\begin{align}
    \sup_{x,\xi \in \bbR} \, \bigl| \px^j \pxi^k \bigl[ e^{-ix\xi} e_\ell^{\#}(x,\xi) \bigr] \bigr| &\lesssim_{\ell,j, k}  1, \\ 
    \sup_{x,\xi \in \bbR} \, \bigl| \jxi^{-1} \px^j \pxi^k \bigl[ e^{-ix\xi} D^\ast e_\ell^{\#}(x,\xi) \bigr] \bigr| &\lesssim_{\ell, j, k} 1.
\end{align}
\end{lemma}

Next, we establish $L^2 \to L^2$ bounds for several pseudo-differential operators involving the distorted Fourier basis elements \eqref{eq:ellsplit} and \eqref{eq:Dsell}.
\begin{lemma} \label{lem:pseboundfreq}
The following bounds hold
\begin{align}
    \biggl\| \int_\bbR \jxi^{-k} \px^k e_\ell^\#(x,\xi) h(\xi) \, \ud \xi \biggr\|_{L^2_x} &\lesssim_{\ell,k}  \|h\|_{L^2_\xi}, \quad 0 \leq k \leq 3, \label{equ:mapping1_L2}  \\ 
    \biggl\| \int_\bbR \jx^{-1} \jxi^{-k} \pxi \px^k e_\ell^\#(x,\xi) h(\xi) \, \ud \xi \biggr\|_{L^2_x} &\lesssim_{\ell,k} \|h\|_{L^2_\xi}, \quad 0 \leq k \leq 3, \label{equ:mapping2_L2} \\ 
    \biggl\| \int_\bbR \jxi^{-k-1} \px^k \bigl( D^\ast e_\ell^\#(x,\xi) \bigr) h(\xi) \, \ud \xi \biggr\|_{L^2_x} &\lesssim_{\ell,k} \|h\|_{L^2_\xi}, \quad 0 \leq k \leq 2, \label{equ:mapping3_L2} \\  
    \biggl\| \int_\bbR \jx^{-1} \jxi^{-k-1} \pxi \px^k \bigl( D^\ast e_\ell^\#(x,\xi) \bigr) h(\xi) \, \ud \xi \biggr\|_{L^2_x} &\lesssim_{\ell,k} \|h\|_{L^2_\xi}, \quad 0 \leq k \leq 2. \label{equ:mapping4_L2}     
\end{align}
\end{lemma}
\begin{proof}
We begin with the proof of \eqref{equ:mapping1_L2} for $k=0$.
Inserting \eqref{eq:ellsplit}, \eqref{eq:mImII}, and using the standard Plancherel theorem, we immediately obtain the desired bound
\begin{align}
\biggl\| \int_\bbR e_\ell^\#(x,\xi) h(\xi) \, \ud \xi \biggr\|_{L^2_x}&\lesssim \biggl\| \int_\bbR e^{i x \xi}   m_{\mathrm{I}}(\xi) h(\xi) \, \ud \xi \biggr\|_{L^2_x}+\biggl\|  \tanh(\gamma x) \int_\bbR   e^{i x \xi}  m_{\mathrm{II}}(\xi) h(\xi) \, \ud \xi \biggr\|_{L^2_x}\\
&\lesssim \|m_{\mathrm{I}}(\xi) h(\xi)\|_{L^2_\xi}+ \|m_{\mathrm{II}}(\xi) h(\xi)\|_{L^2_\xi}\lesssim \|h\|_{L^2_\xi}.
\end{align}
For the proof of \eqref{equ:mapping1_L2} for $1 \leq k \leq 3$ and for the proofs of \eqref{equ:mapping2_L2}, \eqref{equ:mapping3_L2}, \eqref{equ:mapping4_L2}, we observe that for all integers $k \geq 0$ and for $j = 0, 1$,
\begin{equation} \label{eq:mbounds}
    \sup_{x,\xi\in\bbR} \, \bigl|\partial_\xi^j \px^j m^\#(\gamma x,\xi)\bigr| = \sup_{x,\xi\in\bbR} \, \bigl| \partial_\xi^j m_{\mathrm{I}}(\xi)+\partial_x^j\tanh(\gamma x)\partial_\xi^jm_{\mathrm{II}}(\xi) \bigr| \lesssim_{k,j} 1,
\end{equation}
and we compute
\begin{equation} \label{eq:pxiDe}
\begin{aligned}
    \partial_\xi D^\ast e_\ell^\#(x,\xi) &= \partial_\xi \bigl( -\ell \px + i (\jxi + \ell \xi) \bigr) e_\ell^{\#}(x,\xi) \\
    &= - \ell \partial_\xi \px e_\ell^\#(x,\xi) + i (\jxi + \ell \xi) \partial_\xi e_\ell^\#(x,\xi) + i ( \xi \jxi^{-1} + \ell) e_\ell^\#(x,\xi).
\end{aligned}
\end{equation}
By the product rule of differentiation, it follows by inspection that
\begin{equation}\label{eq:mbounds2}
    \jx^{-j} \jxi^{-k} \pxi^j \px^k e_\ell^\#(x,\xi) \qquad \text{and} \qquad \jx^{-j} \jxi^{-k-1} \pxi^j \px^k \bigl( D^\ast e_\ell^\#(x,\xi) \bigr)
\end{equation}
can be written as a linear combination of tensorized product terms $\fraka(x) \frakb(\xi) e^{ix\xi}$ for bounded functions $\fraka, \frakb \in W^{1,\infty}(\bbR)$. Then all the asserted bounds follow by the same argument as above from the standard Plancherel theorem.
\end{proof}

Next, we establish mapping properties for the scalar transforms $\calF_\ell^{\#}$ and $\calF_{\ell,D}^{\#}$.

\begin{lemma} \label{lem:mapping_properties_calF}
The scalar transform $\calF_\ell^{\#}$ defined in \eqref{eq:modF} satisfies
\begin{align}
    \bigl\| \calF_{\ell}^{\#}[f](\xi) \bigr\|_{L^2_\xi} &\lesssim \|f\|_{L^2_x}, \label{equ:mapping_property_calFulellsh_L2} \\
    \bigl\| \jxi \calF_{\ell}^{\#}[f](\xi) \bigr\|_{L^2_\xi} &\lesssim \|f\|_{H^1_x}, \label{equ:mapping_property_calFulellsh_jxi} \\    
    \bigl\| \jxi^2\calF_{\ell}^{\#}[f](\xi) \bigr\|_{L^2_\xi} &\lesssim \|f\|_{H^2_x}, \label{equ:mapping_property_calFulellsh_jxi2} \\
    \bigl\| \pxi \calF_{\ell}^{\#}[f](\xi) \bigr\|_{L^2_\xi} &\lesssim \|\jx f\|_{L^2_x}, \label{equ:mapping_property_calFulellsh_pxi} \\
    \bigl\| \jxi^2 \pxi \calF_{\ell}^{\#}[f](\xi) \bigr\|_{L^2_\xi} &\lesssim \|\jx f\|_{H^2_x}, \label{equ:mapping_property_calFulellsh_jxi2_pxi} \\
     \bigl\| \jxi^{\frac32} \calF_{\ell}^{\#}\bigl[f\bigr](\xi) \bigr\|_{L^\infty_\xi} &\lesssim \bigl\| \jDx^2 f \bigr\|_{L^1_x}, \label{equ:mapping_property_calFulellsh_jxi32_Linfty}
\end{align}
and the scalar transform $\calF_{\ell,D}^{\#}$ defined in \eqref{eq:modFD} satisfies
\begin{align}
    \bigl\| \calF_{\ell,D}^{\#}[f](\xi) \bigr\|_{L^2_\xi} &\lesssim \|f\|_{H^1_x}, \label{equ:mapping_property_calFulellDsh_L2} \\
    \bigl\| \jxi^2 \calF_{\ell,D}^{\#}[f](\xi) \bigr\|_{L^2_\xi} &\lesssim \|f\|_{H^3_x}, \label{equ:mapping_property_calFulellDsh_jxi2} \\
    \bigl\| \pxi \calF_{\ell,D}^{\#}[f](\xi) \bigr\|_{L^2_\xi} &\lesssim \|\jx f\|_{H^1_x}, \label{equ:mapping_property_calFulellDsh_pxi} \\
    \bigl\| \jxi^2 \pxi \calF_{\ell,D}^{\#}[f](\xi) \bigr\|_{L^2_\xi} &\lesssim \|\jx f\|_{H^3_x}, \label{equ:mapping_property_calFulellDsh_jxi2_pxi} \\
    \bigl\| \jxi^{\frac32} \calF_{\ell,D}^{\#}[f\bigr](\xi) \bigr\|_{L^\infty_\xi} &\lesssim \bigl\| \jDx^3 f \bigr\|_{L^1_x}. \label{equ:mapping_property_calFulellDsh_jxi32_Linfty} 
\end{align}
Moreover, the adjoint transform $\calF_{\ell}^{\#,\ast}$ satisfies the weighted Sobolev bound
\begin{equation} \label{equ:mapping_property_calFulellshast_jy_H2y}
    \bigl\| \jx \calF_{\ell}^{\#,\ast}[h] \bigr\|_{H^2_x} \lesssim \bigl\| \jxi^2 \pxi h(\xi) \bigr\|_{L^2_\xi} + \bigl\| \jxi^2 h(\xi) \bigr\|_{L^2_\xi}. 
\end{equation}   
\end{lemma}
\begin{proof}
The first five bounds \eqref{equ:mapping_property_calFulellsh_L2}, \eqref{equ:mapping_property_calFulellsh_jxi},  \eqref{equ:mapping_property_calFulellsh_jxi2}, \eqref{equ:mapping_property_calFulellsh_pxi}, \eqref{equ:mapping_property_calFulellsh_jxi2_pxi} follow from the tensorized product structure of $\pxi^j e_\ell^{\#}(x,\xi)$, $j = 0, 1$, the standard Plancherel theorem, and basic properties of the flat Fourier transform, as in the proof of Lemma~\ref{lem:pseboundfreq}.

Next, we discuss the $L^\infty_\xi$ bound \eqref{equ:mapping_property_calFulellsh_jxi32_Linfty}.
Clearly, we have 
\begin{equation}
     \bigl\| \jxi^{\frac32} \calF_{\ell}^{\#}\bigl[f\bigr](\xi) \bigr\|_{L^\infty_\xi} \leq    \bigl\| \jxi^{2} \calF_{\ell}^{\#}\bigl[f\bigr](\xi) \bigr\|_{L^\infty_\xi},
\end{equation}
so it suffices to bound the latter one.
Using \eqref{eq:ellsplit} and integrating by parts, we find
\begin{align}
\jxi^{2} \calF_{\ell}^{\#}\bigl[f\bigr](\xi) &= \frac{1}{\sqrt{2\pi}}\int (\xi^2+1) e^{-i\xi x}\Big (\overline{m}_{\mathrm{I}}(\xi)+\tanh(\gamma x) \overline{m}_{\mathrm{II}}(\xi)\Big) f(x)\,\ud x \\
&= \frac{1}{\sqrt{2\pi}} \int e^{-i\xi x}\Big (\overline{m}_{\mathrm{I}}(\xi)+\tanh(\gamma x) \overline{m}_{\mathrm{II}}(\xi)\Big) f(x)\,\ud x \\
&\quad - \frac{1}{\sqrt{2\pi}} \int  e^{-i\xi x}\Big (\overline{m}_{\mathrm{I}}(\xi)+\tanh(\gamma x) \overline{m}_{\mathrm{II}}(\xi)\Big) \partial_x^2 f(x)\,\ud x \\
&\quad - \frac{1}{\sqrt{2\pi}} \int  e^{-i\xi x}2\partial_x \tanh(\gamma x) \overline{m}_{\mathrm{II}}(\xi) \partial_x f(x)\,\ud x \\
&\quad - \frac{1}{\sqrt{2\pi}} \int  e^{-i\xi x}\partial_x^2 \tanh(\gamma x) \overline{m}_{\mathrm{II}}(\xi) f(x)\,\ud x.
\end{align}
Taking the $L^\infty_\xi$, we conclude the asserted bound \eqref{equ:mapping_property_calFulellDsh_jxi32_Linfty}.

The $L^2$-type bounds \eqref{equ:mapping_property_calFulellDsh_L2}, \eqref{equ:mapping_property_calFulellDsh_jxi2}, \eqref{equ:mapping_property_calFulellDsh_pxi}, \eqref{equ:mapping_property_calFulellDsh_jxi2_pxi} for $\calF_{\ell,D}^{\#}$ follow similarly from \eqref{eq:Dsell}, \eqref{eq:pxiDe}, and standard properties of the flat Fourier transform. The last $L^\infty_\xi$ bound \eqref{equ:mapping_property_calFulellDsh_jxi32_Linfty} follows from the explicit formula \eqref{eq:Dsell} and the same argument as for the $L^\infty_\xi$-type bound for $\calF_{\ell}^{\#}$. 

Finally the last bound \eqref{equ:mapping_property_calFulellshast_jy_H2y} again follows from the tensorized product structure of $e_\ell^{\#}(x,\xi)$ and basic properties of the flat Fourier transform. We skip the details.
This finishes the proof of the lemma.
\end{proof}

Occasionally, we will also need the following related $L^2$-type bound with a derivative gain.
\begin{lemma}
For any integer $k \geq 1$ it holds that 
\begin{align}
    \biggl\| \jxi^k \int_\bbR e^{-ix\xi} \overline{\px m_\ell^{\#}(\gamma x, \xi)} h(x) \, \ud x \biggr\|_{L^2_\xi} &\lesssim \|h\|_{H^{k-1}_x}, \\
    \biggl\| \jxi^k \int_\bbR e^{-ix\xi} \overline{\px^2 m_\ell^{\#}(\gamma x, \xi)} h(x) \, \ud x \biggr\|_{L^2_\xi} &\lesssim \|h\|_{H^{k-1}_x}.
\end{align}
\end{lemma}
\begin{proof}
In view of \eqref{eq:mImII}, we compute for $j=1,2$ that
\begin{align}
    \px^j m_\ell^\#(\gamma x,\xi) = \partial_x^j \bigl( \tanh(\gamma x) \bigr) m_{\mathrm{II}}(\xi).
\end{align}
Moreover, we recall that
\begin{equation}
  | m_{\mathrm{II}}(\xi)|= \left|\frac{ i }{|\gamma(\xi+\ell\jxi)| - i}\right|\lesssim\frac{1}{\jxi}.
\end{equation}
Then the asserted bounds follow from standard properties of the flat Fourier transform.
\end{proof}

Finally, we collect several mapping properties for the following transform that arises in Section~\ref{sec:setting_up}, specifically in \eqref{equ:setting_up_calL_definition}, from the preparation of a modulation term for the nonlinear analysis,
\begin{align} \label{eq:calelldef}
        \calL_\ell\bigl(\bmf\bigr)(\xi)  &:= \frac{1}{\sqrt{2\pi}} \int_\bbR e^{-ix\xi} \overline{\py m_\ell^{\#}(\gamma x, \xi)} \bigl( i \jxi f_1(t,x) + f_2(t,x) \bigr) \, \ud x \\ 
        &\quad + \frac{\ell}{\sqrt{2\pi}} \int_\bbR e^{-i x \xi} \overline{\py^2 m_\ell^{\#}(\gamma x, \xi)} f_1(t,x)  \, \ud x
    \end{align}
where $\bmf = (f_1, f_2)$ and 
\begin{equation} 
    \begin{aligned}
        \px m_\ell^{\#}(\gamma x,\xi) &= \frac{i \gamma}{\bigl| \gamma (\xi + \ell \jxi) \bigr| - i} \sech^2(\gamma x), \\
        \px^2 m_\ell^{\#}(\gamma x,\xi) &= \frac{-2i\gamma^2}{\bigl| \gamma (\xi + \ell \jxi) \bigr| - i} \sech^2(\gamma x) \tanh(\gamma x).
    \end{aligned}
\end{equation}

In the application of the following bounds in the nonlinear analysis it is important that the second component of the input vector is estimated in terms of one derivative less than the first component.

\begin{lemma} \label{lem:mappingcalLell}
For $\bmf = (f_1, f_2)$ the transform \eqref{eq:calelldef} satisfies 
    \begin{align}
            \bigl\| \jxi^{\frac32} \calL_\ell\bigl( \bmf \bigr)(\xi) \bigr\|_{L^\infty_\xi} &\lesssim \bigl\| \sech^2(\gamma \cdot) f_1 \bigr\|_{W^{2,1}_x} + \bigl\| \sech^2(\gamma \cdot) f_2 \bigr\|_{W^{1,1}_x},\label{eq:calLellLinfty} \\
            \bigl\| \jxi^2 \calL_\ell\bigl(\bmf \bigr)(\xi) \bigr\|_{L^2_\xi} &\lesssim \bigl\| \sech^2(\gamma \cdot) f_1 \bigr\|_{H^2_x} + \bigl\| \sech^2(\gamma \cdot) f_2 \bigr\|_{H^1_x}, \label{eq:calLellL2}\\ 
            \bigl\| \jxi^2 \pxi \calL_\ell\bigl(\bm f\bigr)(\xi) \bigr\|_{L^2_\xi} &\lesssim \bigl\| \jx \sech^2(\gamma \cdot) f_1 \bigr\|_{H^2_x} + \bigl\| \jx \sech^2(\gamma \cdot) f_2 \bigr\|_{H^1_x}\label{eq:pxicalLellL2}.
    \end{align}
\end{lemma}
\begin{proof}
All of these estimates can be proved in the same manner as the estimates in the preceding Lemma~\ref{lem:mapping_properties_calF}. Here the analysis is even simpler because both $\px m_\ell^{\#}(\gamma x,\xi)$ and $\px^2 m_\ell^{\#}(\gamma x,\xi)$ are  localized in space. Moreover, we observe that these two functions decay like $\jxi^{-1}$, which results in gaining one derivative on the right-hand sides of the estimates above. We sketch the proof of the first two estimates \eqref{eq:calLellLinfty}, \eqref{eq:calLellL2}.

For \eqref{eq:calLellLinfty}, as before, it is clear that 
\begin{equation}
    \bigl\| \jxi^{\frac32} \calL_\ell\bigl( \bmf \bigr)(\xi) \bigr\|_{L^\infty_\xi}\leq \bigl\| \jxi^2\calL_\ell\bigl( \bmf \bigr)(\xi) \bigr\|_{L^\infty_\xi}.
\end{equation}
Writing out the last term explicitly,
\begin{align}
    \jxi^2\calL_\ell\bigl( \bmf \bigr)(\xi) &= \frac{\xi^2+1}{\sqrt{2\pi}} \int_\bbR e^{-ix\xi} \overline{\px m_\ell^{\#}(\gamma x, \xi)} \bigl( i \jxi f_1(t,x) + f_2(t,x) \bigr) \, \ud x \\ 
        &\quad + \frac{\ell(1+\xi^2)}{\sqrt{2\pi}} \int_\bbR e^{-ix\xi} \overline{\px^2 m_\ell^{\#}(\gamma x, \xi)} f_1(t,x)  \, \ud x
\end{align}
and integrating by parts as in the analysis of the $L^\infty_\xi$ bounds in Lemma~\ref{lem:mapping_properties_calF}, we obtain
\begin{align}
    \left| \frac{\xi^2+1}{\sqrt{2\pi}} \int_\bbR e^{-ix\xi} \overline{\py m_\ell^{\#}(\gamma x, \xi)} \bigl( i \jxi f_1(t,x)  \bigr) \, \ud x \right|_{L^\infty_\xi} &\lesssim \left\| \sech^2(\gamma \cdot)f_1\right\|_{W^{2,1}_x}, \\
        \left| \frac{\xi^2+1}{\sqrt{2\pi}} \int_\bbR e^{-i x \xi} \overline{\px m_\ell^{\#}(\gamma x, \xi)} \bigl(  f_2(t,x)  \bigr) \, \ud x \right|_{L^\infty_\xi} &\lesssim \left\| \sech^2(\gamma \cdot)f_2\right\|_{W^{1,1}_x}, \\
   \left| \frac{\ell(1+\xi^2)}{\sqrt{2\pi}} \int_\bbR e^{-i x \xi} \overline{\py^2 m_\ell^{\#}(\gamma x, \xi)} f_1(t,x)  \, \ud x \right|_{L^\infty_\xi} &\lesssim \left\| \sech^2(\gamma \cdot)\tanh(\gamma \cdot)f_1\right\|_{W^{1,1}_x}.     
\end{align}
Putting all of these estimates together, the asserted bound \eqref{eq:calLellLinfty} follows.
Using the same formula for $\jxi^2\calL_\ell\bigl( \bmf \bigr)(\xi)$ and integrating by parts, Plancherel's theorem implies that
\begin{align}
    \left| \frac{\xi^2+1}{\sqrt{2\pi}} \int_\bbR e^{-ix\xi} \overline{\py m_\ell^{\#}(\gamma x, \xi)} \bigl( i \jxi f_1(t,x) \bigr) \, \ud x \right|_{L^2_\xi} &\lesssim \left\| \sech^2(\gamma \cdot)f_1\right\|_{H^2_x}, \\
    \left| \frac{\xi^2+1}{\sqrt{2\pi}} \int_\bbR e^{-ix\xi} \overline{\py m_\ell^{\#}(\gamma x, \xi)} \bigl(  f_2(t,x)  \bigr) \, \ud x \right|_{L^2_\xi} &\lesssim \left\| \sech^2(\gamma \cdot)f_2\right\|_{H^1_x}, \\
   \left| \frac{\ell(1+\xi^2)}{\sqrt{2\pi}} \int_\bbR e^{-ix\xi} \overline{\py^2 m_\ell^{\#}(\gamma x, \xi)} f_1(t,x)  \, \ud x \right|_{L^2_\xi} &\lesssim \left\| \sech^2(\gamma \cdot)\tanh(\gamma \cdot)f_1\right\|_{H^1_x},     
\end{align}
which imply \eqref{eq:calLellL2}. Finally \eqref{eq:pxicalLellL2} follows analogously.
\end{proof}

\section{Linear Decay Estimates} \label{sec:linear_decay_estimates}

In this section we establish linear decay estimates for the evolution generated by the matrix operator $\bfL_\ell$.
The analysis here is based on the oscillatory integral representation of the evolution obtained in Proposition~\ref{prop:rep_formula_propagator_modified_dFT} in terms of the modified distorted Fourier transform.

For the reader's orientation, we note once more that in the remainder of this paper, the results from this section will be applied in a moving frame with the spatial coordinate $y := x - q(t)$, as introduced in \eqref{equ:setting_up_moving_frame_coordinate}. 
However, since we hope that the results from this section are of independent interest, we use the more conventional notation $x$ for the spatial coordinate here.

\subsection{Pointwise decay} \label{subsec:pointwise_decay}

We start by establishing pointwise dispersive decay estimates and asymptotics for the evolution, which follow from the core technical lemma below.

\begin{lemma} \label{lem:core_linear_dispersive_decay}
Let $\ell \in (-1,1)$. Suppose $\fraka \colon \bbR^2 \to \bbC$ satisfies
\begin{equation} \label{equ:core_linear_lemma_coefficient_assumption}
    \sup_{x, \xi \in \bbR} \, \bigl( |\fraka(x,\xi)| + |\pxi \fraka(x,\xi)| \bigr) \lesssim 1.
\end{equation}
There exists $C \geq 1$ such that for any Schwartz function $h \in \calS(\bbR)$ we have for all $t \geq 1$,
\begin{equation} \label{equ:core_linear_asymptotics}
    \begin{aligned}
        &\biggl\| \int_\bbR e^{ix\xi} e^{it(\jxi+\ell\xi)} \fraka(x,\xi) h(\xi) \, \ud \xi - \sqrt{2\pi} \frac{e^{i\frac{\pi}{4}} e^{i\rho}}{t^{\frac12}} \fraka(x,\xi_\ast) \jap{\xi_\ast}^{\frac32} h(\xi_\ast) \one_{(-t,t)}(x+\ell t) \biggr\|_{L^\infty_x} \\
        &\qquad \qquad \qquad \qquad \qquad \qquad \qquad \qquad \qquad \qquad \leq \frac{C(\ell)}{t^{\frac23}} \Big( \bigl\|\jxi^2 \pxi h(\xi) \bigr\|_{L^2_\xi} + \bigl\| \jxi^2 h(\xi) \bigr\|_{L^2_\xi} \Bigr),
    \end{aligned}
\end{equation}
where for $-t < x+ \ell t < t$ we use the notation
\begin{equation*}
 \rho := \sqrt{t^2-(x+\ell t)^2}, \quad \xi_\ast = - \frac{x+\ell t}{\rho}.
\end{equation*}
In particular, for all $t \geq 1$ it holds that
\begin{equation}
    \begin{aligned}
        &\biggl\| \int_\bbR e^{ix\xi} e^{it(\jxi+\ell\xi)} \fraka(x,\xi) h(\xi) \, \ud \xi \biggr\|_{L^\infty_x} \\
        &\qquad \lesssim \frac{1}{t^{\frac12}} \bigl\| \jxi^{\frac32} h(\xi) \bigr\|_{L^\infty_\xi} + \frac{1}{t^{\frac23}} \Big( \bigl\|\jxi^2 \pxi h(\xi) \bigr\|_{L^2_\xi} + \bigl\| \jxi^2 h(\xi) \bigr\|_{L^2_\xi} \Bigr).
    \end{aligned}
\end{equation}
\end{lemma}

Before we turn to the proof of Lemma~\ref{lem:core_linear_dispersive_decay}, we record the following pointwise decay estimates and asymptotics for the evolution generated by the matrix operator $\bfL_\ell$.

\begin{corollary} \label{cor:linear_asymptotics}
Fix $\ell \in (-1,1)$ and denote by $\bfL_\ell$ the matrix operator defined in \eqref{eq:bfL}.
For a given pair of real-valued Schwartz functions $\bmf = (f_1, f_2) \in \calS(\bbR) \times \calS(\bbR)$, we define
\begin{equation*}
    g_\ell^{\#}(\xi) := \mathcal{T}_\ell [\bmf] (\xi)= \mathcal{F}_{\ell,D}^{\#}[f_1](\xi) - \mathcal{F}_{\ell}^{\#}[f_2](\xi).
\end{equation*}
Then there exists a constant $C(\ell) \geq 1$ such that for any pair of real-valued Schwartz functions $\bmf = (f_1, f_2) \in \calS(\bbR)\times\calS(\bbR)$, we have for all $t \geq 1$,
\begin{equation} \label{equ:linear_asymptotics_1st_comp}
    \begin{aligned}
        &\biggl\| \bigl( e^{t\bfL_\ell} P_\mathrm{e} \bm{f} \bigr)_1(x) + \Im \biggl( \frac{e^{i\frac{\pi}{4}} e^{i\rho}}{t^{\frac12}} m_\ell^\#(\gamma x, \xi_\ast) \jap{\xi_\ast}^{\frac12} g_\ell^{\#}(\xi_\ast) \one_{(-t,t)}(x+\ell t) \biggr) \biggr\|_{L^\infty_x} \\
        &\qquad \qquad \qquad \qquad \qquad \qquad \qquad \leq \frac{C(\ell)}{t^{\frac23}} \Big( \bigl\|\jxi \partial_\xi g_\ell^{\#}(\xi) \bigr\|_{L^2_\xi} + \bigl\|\jxi g_\ell^{\#}(\xi) \bigr\|_{L^2_\xi} \Big),
    \end{aligned}
\end{equation}
as well as 
\begin{equation} \label{equ:linear_asymptotics_2nd_comp_plus_derivative_of_1st_comp}
    \begin{aligned}
        &\biggl\| \px \Bigl( \bigl( e^{t\bfL_\ell} P_\mathrm{e} \bm{f} \bigr)_1(x) \Bigr) + \Im \biggl( \frac{e^{i\frac{\pi}{4}} e^{i\rho}}{t^{\frac12}} n_{\ell,1}^\#(\gamma x, \xi_\ast) \jap{\xi_\ast}^{\frac32} g_\ell^{\#}(\xi_\ast) \one_{(-t,t)}(x + \ell t) \biggr) \biggr\|_{L^\infty_x} \\
        &\, \, + \biggl\| \bigl( e^{t\bfL_\ell} P_\mathrm{e} \bm{f} \bigr)_2(x) + \Im \biggl( \frac{e^{i\frac{\pi}{4}} e^{i\rho}}{t^{\frac12}} n_{\ell,2}^\#(\gamma x, \xi_\ast) \jap{\xi_\ast}^{\frac32} g_\ell^{\#}(\xi_\ast) \one_{(-t,t)}(x + \ell t) \biggr) \biggr\|_{L^\infty_x} \\
        &\qquad \qquad \qquad \qquad \qquad \qquad \qquad \, \, \, \leq \frac{C(\ell)}{t^{\frac23}} \Big( \bigl\|\jxi^2 \partial_\xi g_\ell^{\#}(\xi) \bigr\|_{L^2_\xi} + \bigl\|\jxi^2 g_\ell^{\#}(\xi) \bigr\|_{L^2_\xi} \Big),
    \end{aligned}
\end{equation}
where 
\begin{equation*}
    \rho := \sqrt{t^2-(x+\ell t)^2}, \quad \xi_\ast = - \frac{x+\ell t}{\rho}, \quad -t < x+ \ell t < t,
\end{equation*}
and 
\begin{equation*}
    \begin{aligned}
        m_\ell^\#(\gamma x, \xi_\ast) &:= \frac{\gamma(\xi_\ast + \ell \jap{\xi_\ast}) + i \tanh(\gamma x)}{|\gamma(\xi_\ast + \ell \jap{\xi_\ast})| - i}, \\
        n_{\ell,2}^\#(\gamma x, \xi_\ast) &:= i m_\ell^\#(\gamma x, \xi_\ast) - \jap{\xi_\ast}^{-1} \ell \px  m_\ell^\#(\gamma x, \xi_\ast), \\ 
        n_{\ell,1}^\#(\gamma x, \xi_\ast) &:= i \jap{\xi_\ast}^{-1} \xi_\ast m_\ell^{\#}(\gamma x, \xi_\ast) + \jap{\xi_\ast}^{-1} \px m_\ell^{\#}(\gamma x, \xi_\ast).
    \end{aligned}
\end{equation*}
In particular, for all $t \geq 1$ it holds that
\begin{equation} \label{equ:linear_evolution_dispersive_decay}
    \begin{aligned}
        &\bigl\|\bigl( e^{t\bfL_\ell} P_\mathrm{e} \bm{f} \bigr)_1\bigr\|_{L^\infty_x} + \bigl\|\bigl( e^{t\bfL_\ell} P_\mathrm{e} \bm{f} \bigr)_2\bigr\|_{L^\infty_x} + \bigl\|\px \bigl( e^{t\bfL_\ell} P_\mathrm{e} \bm{f} \bigr)_1\bigr\|_{L^\infty_x} \\
        &\quad \lesssim \frac{1}{t^{\frac12}} \bigl\| \jxi^{\frac32} g_\ell^{\#}(\xi) \bigr\|_{L^\infty_\xi} + \frac{1}{t^{\frac23}} \Bigl( \bigl\|\jxi^2 \partial_\xi g_\ell^{\#}(\xi) \bigr\|_{L^2_\xi} + \bigl\|\jxi^2 g_\ell^{\#}(\xi) \bigr\|_{L^2_\xi} \Bigr).
    \end{aligned}
\end{equation}
\end{corollary}
\begin{proof}
    The asymptotics \eqref{equ:linear_asymptotics_1st_comp} for $( e^{t\bfL_\ell} P_\mathrm{e} \bm{f} )_1(x)$ follow from the representation formula \eqref{equ:rep_formula_propagator_modified_dFT}, the pointwise bounds for the distorted Fourier basis elements obtained in Lemma~\ref{lem:boundseDe}, and Lemma~\ref{lem:core_linear_dispersive_decay} with the choices $\fraka(x,\xi) := e^{-ix\xi} e_\ell^{\#}(x,\xi)$ and $h(\xi) := i \jxi^{-1} g_\ell^{\#}(\xi)$.
    Similarly, we infer the asymptotics \eqref{equ:linear_asymptotics_2nd_comp_plus_derivative_of_1st_comp} for $( e^{t\bfL_\ell} P_\mathrm{e} \bm{f} )_2(x)$ with the choices $\fraka(x,\xi) := \jxi^{-1} e^{-ix\xi} D^\ast e_\ell^{\#}(x,\xi)$ and $h(\xi) := i g_\ell^{\#}(\xi)$, while the asymptotics \eqref{equ:linear_asymptotics_2nd_comp_plus_derivative_of_1st_comp} for $\px ( e^{t\bfL_\ell} P_\mathrm{e} \bm{f} )_1(x)$ follow by taking $\fraka(x,\xi) := \jxi^{-1} e^{-ix\xi} \px e_\ell^{\#}(x,\xi)$ and $h(\xi) := i g_\ell^{\#}(\xi)$.
    Finally, the dispersive decay estimate \eqref{equ:linear_evolution_dispersive_decay} is a direct consequence of \eqref{equ:linear_asymptotics_1st_comp} and~\eqref{equ:linear_asymptotics_2nd_comp_plus_derivative_of_1st_comp}.
\end{proof}

Now we prove Lemma~\ref{lem:core_linear_dispersive_decay}.
\begin{proof}[Proof of Lemma~\ref{lem:core_linear_dispersive_decay}]
    We adapt the proof of \cite[Lemma 2.1]{LS1} to our setting.
    Introducing the phase function 
    \begin{equation*}
        \phi_\ell(\xi;v) := \jxi + (v+\ell) \xi, \quad v := \frac{x}{t},
    \end{equation*}    
    we write
    \begin{equation*}
        \int_\bbR e^{ix\xi} e^{it(\jxi+\ell\xi)} \fraka(x,\xi) h(\xi) \, \ud \xi = \int_\bbR e^{i t \phi_\ell(\xi; v)} \fraka(x,\xi) h(\xi) \, \ud \xi.
    \end{equation*} 
    Next, we compute 
    \begin{equation*}
        \pxi \phi_\ell(\xi;v) = \frac{\xi}{\jxi} + (v+\ell)
    \end{equation*}
    as well as
    \begin{equation} \label{equ:disp_decay_proof_pxi2_phi}
        \pxi^2 \phi_\ell(\xi;v) = \jxi^{-3}.
    \end{equation}
    It follows that the phase $\phi_\ell(\xi;v)$ is stationary $\pxi \phi_\ell(\xi_\ast;v) = 0$ at some point $\xi_\ast \in \bbR$ if and only if
    \begin{equation} \label{equ:disp_decay_proof_stat_point_cond}
        \frac{\xi_\ast}{\jap{\xi_\ast}} = - (v+\ell).
    \end{equation}
    Correspondingly, the phase $\phi_\ell(\xi;v)$ has a unique stationary point $\xi_\ast$ satisfying \eqref{equ:disp_decay_proof_stat_point_cond} if and only if $-1 < v+\ell < 1$, which is equivalent to $-t < x+\ell t < t$ for $t > 0$. By direct computation, we find for $t > 0$ that 
    \begin{equation} \label{equ:disp_decay_proof_xiast_explicit}
     \xi_\ast = - \frac{v+\ell}{\sqrt{1-(v+\ell)^2}} = - \frac{x + \ell t}{\sqrt{t^2 - (x+\ell t)^2}}
    \end{equation}
    as well as 
    \begin{equation*}
     \jap{\xi_\ast} = \frac{1}{\sqrt{1-(v+\ell)^2}} = \frac{t}{\sqrt{t^2-(x+\ell t)^2}}.
    \end{equation*}

    We first treat the case when the phase does not have a stationary point. If $|v+\ell| \geq 1$, then we have uniformly for all $\xi \in \bbR$ that
    \begin{equation} \label{equ:disp_decay_proof_pxi_phi_nonstationary}
        \bigl| \pxi \phi_\ell(\xi; v) \bigr| \geq |v+\ell| - \frac{|\xi|}{\jxi} \geq 1 - \frac{|\xi|}{\jxi} \geq \frac{1}{2\jxi^2}.
    \end{equation}
    We denote by $\chi_0$ a smooth cut-off to the unit interval $[-1,1]$, and set $\chi_1 := 1 - \chi_0$. Integrating by parts followed by Cauchy-Schwarz using \eqref{equ:core_linear_lemma_coefficient_assumption}, \eqref{equ:disp_decay_proof_pxi2_phi}, and \eqref{equ:disp_decay_proof_pxi_phi_nonstationary}, we obtain for $t \geq 1$,
    \begin{equation} \label{equ:disp_decay_proof_int_by_parts_outside_light_cone}
        \begin{aligned}
            \biggl| \int_\bbR e^{i t \phi_\ell(\xi; v)} \fraka(x,\xi) h(\xi) \, \ud \xi \biggr| &\lesssim \int_\bbR \chi_1 \bigl( \xi^2 t^{-1} \bigr) \bigl| \fraka(x,\xi) h(\xi) \bigr| \, \ud \xi \\
            &\quad \quad + t^{-1} \int_\bbR \frac{|\pxi^2 \phi_\ell(\xi;v)|}{|\pxi \phi_\ell(\xi;v)|^2} \chi_0\bigl( \xi^2 t^{-1} \bigr) \bigl| \fraka(x,\xi) h(\xi) \bigr| \, \ud \xi \\
            &\quad \quad + t^{-1}  \int_\bbR \frac{1}{|\pxi \phi_\ell(\xi;v)|} \Bigl| \pxi \Bigl( \chi_0\bigl( \xi^2 t^{-1} \bigr) \fraka(x,\xi) h(\xi) \Bigr) \Bigr| \, \ud \xi \\
            &\lesssim \frac{1}{t^{\frac34}} \bigl( \|\jxi^2 \pxi h\|_{L^2_\xi} + \|\jxi^2 h\|_{L^2_\xi} \bigr).
        \end{aligned}
    \end{equation}

    Now we turn to the case $-1 < v + \ell < 1$ where the phase has a unique stationary point.
    We claim that the bound \eqref{equ:disp_decay_proof_pxi_phi_nonstationary} continues to hold (up to multiplicative constants) for all $\xi \in \bbR \setminus I(\xi_\ast)$, where
    \[
        I(\xi_\ast) := \bigl[ \xi_\ast - \jap{\xi_\ast}/100, \xi_\ast + \jap{\xi_\ast}/100 \bigr].
    \]
    To see this we compute that
    \begin{equation} \label{equ:disp_decay_proof_pxi_phi_subtr_off}
        \begin{aligned}
            \pxi \phi_\ell(\xi;v) = \pxi \phi_\ell(\xi;v) - \pxi \phi_\ell(\xi_\ast;v) = \frac{\xi}{\jxi} - \frac{\xi_\ast}{\jap{\xi_\ast}} = \frac{\xi^2-\xi_\ast^2}{\jxi \jap{\xi_\ast} (\xi \jap{\xi_\ast} + \xi_\ast \jap{\xi})}.
        \end{aligned}
    \end{equation}
    We may assume that $\xi_\ast \geq 0$ without loss of generality. If $\xi \geq \xi_\ast + \jap{\xi_\ast}/100$, then \eqref{equ:disp_decay_proof_pxi_phi_subtr_off} implies
    \begin{equation*}
        |\pxi \phi_\ell(\xi;v)| \gtrsim \jap{\xi_\ast}^{-2} \gtrsim \jap{\xi}^{-2}.
    \end{equation*}
    Indeed, since $\xi\geq \max(1,101\xi_\ast)/100$, then for the absolute value of the denominator one has $\simeq \jxi^2 \jap{\xi_\ast}^2$, while the absolute value of the numerator exceeds $ \jap{\xi_\ast}^2$ up to a multiplicative constant.
    Next, if $\xi \leq \xi_\ast - \jap{\xi_\ast}/100$, then we assert that
    $|\partial_\xi \phi_\ell (\xi; v) | \gtrsim_{|\ell|} \jap{\xi}^{-2}$.
    Indeed, if $\xi \leq 0$, then we infer from the first line of \eqref{equ:disp_decay_proof_pxi_phi_subtr_off} that $|\partial_\xi \phi_\ell(\xi, v) | \ge \xi_\ast\jap{\xi_\ast}^{-1} \gtrsim 1$ unless $0 \leq \xi_\ast \ll 1$. In that latter case, one has  $\xi \leq \xi_\ast - \jap{\xi_\ast}/100$, whence $\xi \leq -1/200$, say. But now the first line of \eqref{equ:disp_decay_proof_pxi_phi_subtr_off} implies that $|\partial_\xi \phi_\ell(\xi, v)| \gtrsim 1$.
    On the other hand, consider $0 \leq \xi \leq \xi_\ast - \jap{\xi_\ast}/100$. Then $\xi_\ast\gtrsim 1$ and the numerator is bounded from below by $\jap{\xi_\ast}^2$, while the denominator is bounded from above by $(\jxi \jap{\xi_\ast})^2$. Thus, we have established the claim provided $\xi_\ast \geq 0$. Finally, changing  $\xi \mapsto -\xi$ reduces the case $\xi_\ast\le0$ to the previous one, and we have established the claim above in general.

    Now we define
    \begin{equation*}
     \chi_{\xi_\ast}(\xi) := \chi_0 \bigl(C_0(\xi-\xi_\ast)\jap{\xi_\ast}^{-1} \bigr)
    \end{equation*}
    for some large constant $C_0 \geq 100$, and $\chi_0$ as above.
    Integrating by parts as in~\eqref{equ:disp_decay_proof_int_by_parts_outside_light_cone} we therefore obtain
    \begin{equation} \label{equ:disp_decay_proof_main_contribution_isolated}
        \begin{aligned}
            \biggl| \int_\bbR e^{i t \phi_\ell(\xi; v)} \bigl( 1 - \chi_{\xi_\ast}(\xi) \bigr) \fraka(x,\xi) h(\xi) \, \ud \xi \biggr|
            \lesssim \frac{1}{t^{\frac34}} \bigl( \|\jxi^2 \pxi h\|_{L^2_\xi} + \|\jxi^2 h\|_{L^2_\xi} \bigr),
        \end{aligned}
    \end{equation}
    which holds uniformly in $t \geq 1$ and $x \in \bbR$. 
    To analyze the leading order contribution to the asymptotics \eqref{equ:core_linear_asymptotics} stemming from
    \begin{equation*}
        \int_\bbR e^{i t \phi_\ell(\xi; v)} \chi_{\xi_\ast}(\xi) \fraka(x,\xi) h(\xi) \, \ud \xi,
    \end{equation*}
    we consider the phase difference
    \begin{equation*}
        \begin{aligned}
            \phi_\ell(\xi;v) - \phi_\ell(\xi_\ast;v) = \jxi - \jap{\xi_\ast}  + (v+\ell) (\xi - \xi_\ast).
        \end{aligned}
    \end{equation*}
    Inserting \eqref{equ:disp_decay_proof_stat_point_cond}, we compute that
    \begin{equation*}
        \phi_\ell(\xi;v) - \phi_\ell(\xi_\ast;v) = \frac{(\xi-\xi_\ast)^2}{\jap{\xi_\ast} (1 + \xi \xi_\ast + \jap{\xi} \jap{\xi_\ast})} =: \zeta^2.
    \end{equation*}
    The change of variables $\xi\mapsto\zeta$ is a diffeomorphism on the support of $\chi_{\xi_\ast}(\xi)$ given by
    \begin{equation*}
        \begin{aligned}
          \zeta = \frac{\xi-\xi_\ast}{\sqrt{\jap{\xi_\ast} (1+\xi\xi_\ast+\jap{\xi}\jap{\xi_\ast})} }, \quad \frac{\ud \zeta}{\ud \xi} \simeq \jap{\xi}^{-\frac32}, \quad \frac{\ud \zeta}{\ud \xi}(\xi_\ast) =  \frac{1}{\sqrt{2}} \jap{\xi_\ast}^{-\frac32}.
        \end{aligned}
    \end{equation*}
    Moreover, note that by direct computation we have for $t > 0$,
    \begin{equation}
     \phi_\ell(\xi_\ast;v) = \sqrt{1-(v+\ell)^2}, \quad \rho := t \phi_\ell(\xi_\ast;v) = \sqrt{t^2 - (x + \ell t)^2}.
    \end{equation}    
    Then the change of variables reduces us to complex Gaussians with explicit Fourier transforms
    \begin{equation*}
        \begin{aligned}
            \int_\bbR e^{it\phi_\ell(\xi;v)} \chi_{\xi_\ast}(\xi) \fraka(x, \xi) h(\xi) \, \ud\xi &= e^{i\rho}  \int_\bbR e^{it\zeta^2}\; \overline{G}(\zeta;t,x) \, \ud \zeta = \frac{e^{i\rho}e^{i\frac{\pi}{4}}}{\sqrt{2t}} \int_\bbR e^{-\frac{iy^2}{4t}} \; \overline{\widehat{G}(y;t,x)} \, \ud  y,
        \end{aligned}
    \end{equation*}
    where 
    \begin{equation*}
        \overline{G}(\zeta;t,x) := \chi_{\xi_\ast}(\xi) \fraka(x,\xi) h(\xi) \frac{\ud \xi}{\ud \zeta}.
    \end{equation*}
    It follows from
    \begin{equation*}
        \begin{aligned}
            \int_\bbR \overline{\widehat{G}(y;t,x)} \, \ud  y = \sqrt{2\pi} \overline{G}(0;t,x) &= \sqrt{2\pi} \chi_{\xi_\ast}(\xi_\ast) \fraka(x,\xi_\ast) h(\xi_\ast) \frac{\ud \xi}{\ud \zeta}(\xi_\ast) 
            = \sqrt{4\pi} \jap{\xi_\ast}^{\frac32} \fraka(x,\xi_\ast) h(\xi_\ast)
        \end{aligned}
    \end{equation*}
    that
    \begin{equation} \label{equ:disp_decay_proof_leading_order_term_asymptotics}
        \begin{aligned}
            &\int_\bbR e^{it\phi_\ell(\xi;v)} \chi_{\xi_\ast}(\xi) \fraka(x,\xi) h(\xi) \, \ud\xi \\
            &= \sqrt{2\pi} \frac{e^{i\rho}e^{i\frac{\pi}{4}}}{\sqrt{t}} \jap{\xi_\ast}^{\frac32} \fraka(x,\xi_\ast)  h(\xi_\ast) + \calO_{L^\infty_x} \biggl( t^{-\frac12} \int_\bbR \Bigl| e^{-\frac{iy^2}{4t}} - 1 \Bigr| \, \bigl| {\widehat{G}(y;t,x)} \bigr| \, \ud y \biggr).
        \end{aligned}
    \end{equation}
    The error term in the last line is estimated just as in \cite[Lemma 2.1]{LS1}. 
    To be specific, using additionally \eqref{equ:core_linear_lemma_coefficient_assumption}, we find that
    \begin{equation} \label{equ:disp_decay_proof_leading_order_term_asymptotics_remainder_bound}
        \int_\bbR \Bigl| e^{-\frac{iy^2}{4t}} - 1 \Bigr| \, \bigl|{\widehat{G}(y;t,x)}  \bigr| \, \ud y \lesssim t^{-\frac16} \biggl( \int_\bbR \bigl( |h(\xi)|^2 + |\partial_\xi h(\xi)|^2 \bigr) \bigl|\fraka(x,\xi)\bigr|^2 \jap{\xi}^4 \, \ud \xi \bigg)^{\frac12}.
    \end{equation}
    Combining \eqref{equ:disp_decay_proof_int_by_parts_outside_light_cone}, 
    \eqref{equ:disp_decay_proof_main_contribution_isolated},
    \eqref{equ:disp_decay_proof_leading_order_term_asymptotics}, and \eqref{equ:disp_decay_proof_leading_order_term_asymptotics_remainder_bound} yields \eqref{equ:linear_asymptotics_1st_comp} uniformly for $t \geq 1$ and $x\in\bbR$.
\end{proof}

\subsection{Improved local decay} \label{subsec:improved_local_decay_estimates}

While the threshold resonances of the matrix operator $\bfL_\ell$ cause slow local decay of the evolution, one nevertheless gains improved local decay at high frequencies and upon substracting off the resonance.
For the reader's orientation we note that the choice of the cut-off function $\chi_{\{\leq 2\gamma|\ell|\}}(\xi)$ in the statement of the lemmas below is informed by the observation that the phase of $e^{i t (\jxi + \ell \xi)}$ has a unique stationary point at $\xi=-\gamma\ell$

The improved local decay for the evolution at high frequencies is based on the following technical lemma.

\begin{lemma} \label{lem:improved_local_decay}
Let $\ell \in (-1,1)$. Denote by $\chi_{\{\leq 2\gamma|\ell|\}}(\xi)$ a smooth even bump function satisfying $\chi_{\{\leq 2\gamma|\ell|\}}(\xi) = 1$ for $|\xi| \leq 2\gamma|\ell|$ and $\chi_{\{\leq 2\gamma|\ell|\}}(\xi) = 0$ for $|\xi| \geq 4\gamma|\ell|$.
Then there exists a constant $C(\ell) \geq 1$ such that for $0 \leq k \leq 3$ and for all $t \geq 0$,
\begin{equation} \label{equ:improved_local_decay_high_freq}
    \begin{aligned}
        \biggl\| \jx^{-1} \px^k \biggl( \int_\bbR e_\ell^{\#}(x,\xi) e^{i t (\jxi + \ell \xi)} &\bigl( 1 - \chi_{\{\leq 2\gamma|\ell|\}}(\xi) \bigr) \, i\jxi^{-1} g(\xi) \, \ud \xi \biggr) \biggr\|_{L^2_x} \\ 
        &\quad \quad \quad \quad \quad  \leq \frac{C(\ell)}{\jt} \Bigl( \bigl\|\jxi^{k-1} \pxi g(\xi) \bigr\|_{L^2_\xi} + \bigl\| \jxi^{k-1} g(\xi) \bigr\|_{L^2_\xi} \Bigr).
    \end{aligned}
\end{equation}
Similarly, it holds for $0 \leq k \leq 2$ and for all $t \geq 0$ that
\begin{equation} \label{equ:improved_local_decay_high_freq_for_2nd_comp}
    \begin{aligned}
        \biggl\| \jx^{-1} \px^k \biggl( \int_\bbR D^\ast e_\ell^{\#}(x,\xi) e^{i t (\jxi + \ell \xi)} &\bigl( 1 - \chi_{\{\leq 2\gamma|\ell|\}}(\xi) \bigr) \, i\jxi^{-1} g(\xi) \, \ud \xi \biggr) \biggr\|_{L^2_x} \\ 
        &\quad \quad \leq \frac{C(\ell)}{\jt} \Bigl( \bigl\|\jxi^k \pxi g(\xi) \bigr\|_{L^2_\xi} + \bigl\| \jxi^k g(\xi) \bigr\|_{L^2_\xi} \Bigr).
    \end{aligned}
\end{equation}
\end{lemma}
\begin{proof}
The bounds for short times $0 \leq t \leq 1$ just follow from straightforward energy estimates. So we assume that $t \geq 1$.
Integrating by parts in the frequency variable, we find 
\begin{equation*}
    \begin{aligned}
        &\jx^{-1} \px^k \biggl( \int_\bbR e_\ell^{\#}(x,\xi) e^{i t (\jxi + \ell \xi)} \bigl( 1 - \chi_{\{\leq 2\gamma|\ell|\}}(\xi) \bigr) \, i\jxi^{-1} g(\xi) \, \ud \xi \biggr) \\ 
        &= - \frac{1}{t} \int_\bbR \jx^{-1} \jxi^{-k} \pxi \px^k e_\ell^{\#}(x,\xi) e^{it(\jxi+\ell\xi)} \bigl( 1 - \chi_{\{\leq 2\gamma|\ell|\}}(\xi) \bigr) \bigl(\xi+\ell\jxi\bigr)^{-1} \jxi^{k-1} g(\xi) \, \ud \xi \\ 
        &\quad - \frac{1}{t} \jx^{-1} \int_\bbR \jxi^{-k} \px^k e_{\ell}^{\#}(x,\xi) e^{i t (\jxi + \ell \xi)} \jxi^k \pxi \biggl( \bigl( 1 - \chi_{\{\leq 2\gamma|\ell|\}}(\xi) \bigr) \bigl(\xi+\ell\jxi\bigr)^{-1} g(\xi) \biggr) \, \ud \xi.
    \end{aligned}
\end{equation*}
Observe that $\xi+\ell\jxi = 0$ if and only if $\xi=-\gamma\ell$, which is strictly away from the support of $1 - \chi_{\{\leq 2\gamma|\ell|\}}(\xi)$. 
Correspondingly, the asserted improved local decay estimate \eqref{equ:improved_local_decay_high_freq} follows from the mapping properties \eqref{equ:mapping1_L2}, \eqref{equ:mapping2_L2}, and the uniform bound
\begin{equation*}
    \sup_{|\xi| \geq 2 \gamma |\ell|} \, \bigl| (\xi+\ell \jxi)^{-1} \jxi \bigr| \lesssim_\ell 1.
\end{equation*}
The proof of \eqref{equ:improved_local_decay_high_freq_for_2nd_comp} proceeds analogously using the mapping properties \eqref{equ:mapping3_L2} and \eqref{equ:mapping4_L2}. 
\end{proof}

The improved local decay for the evolution upon subtracting off the threshold resonances is based on the following technical lemma.

\begin{lemma} \label{lem:local_decay_resonance_subtracted_off}
Let $\ell \in (-1,1)$. Denote by $\chi_{\{\leq 2\gamma|\ell|\}}(\xi)$ a smooth even bump function satisfying $\chi_{\{\leq 2\gamma|\ell|\}}(\xi) = 1$ for $|\xi| \leq 2\gamma|\ell|$ and $\chi_{\{\leq 2\gamma|\ell|\}}(\xi) = 0$ for $|\xi| \geq 4\gamma|\ell|$.
Then there exists a constant $C(\ell) \geq 1$ such that for $0 \leq k \leq 3$ and for all $t \geq 0$,
\begin{equation} \label{equ:local_decay_resonance_subtracted_off}
    \begin{aligned}
        \biggl\| \jx^{-3} \px^k \biggl( \int_\bbR \Bigl( e_\ell^{\#}(x,\xi) - e_\ell^{\#}(x,- \gamma \ell) \Bigr) e^{i t (\jxi + \ell \xi)} &\chi_{\{\leq 2\gamma|\ell|\}}(\xi) \, i\jxi^{-1} g(\xi) \, \ud \xi \biggr) \biggr\|_{L^\infty_x} \\ 
        &\quad \quad \, \, \, \, \, \leq \frac{C(\ell)}{\jt} \Bigl( \|\pxi g(\xi) \|_{L^2_\xi} + \|g(\xi)\|_{L^2_\xi} \Bigr).
    \end{aligned}
\end{equation}
Similarly, it holds for all $t \geq 0$ that
\begin{equation} \label{equ:local_decay_resonance_subtracted_off_for_2nd_comp}
    \begin{aligned}
        \biggl\| \jx^{-3} \biggl( \int_\bbR \Bigl( \jxi^{-1} (D^\ast e_\ell^{\#})(x,\xi) - \gamma^{-1} (D^\ast e_\ell^{\#})(x,- \gamma \ell) \Bigr) &e^{i t (\jxi + \ell \xi)} \chi_{\{\leq 2\gamma|\ell|\}}(\xi) \, i g(\xi) \, \ud \xi \biggr) \biggr\|_{L^\infty_x} \\ 
        &\leq \frac{C(\ell)}{\jt} \Bigl( \|\pxi g(\xi) \|_{L^2_\xi} + \|g(\xi)\|_{L^2_\xi} \Bigr).
    \end{aligned}
\end{equation}
\end{lemma}
\begin{proof}
    It suffices to consider times $t \geq 1$.
    We begin with the proof of \eqref{equ:local_decay_resonance_subtracted_off}.
    Integrating by parts in the frequency variable, we obtain for $0 \leq k \leq 3$,
    \begin{equation*}
        \begin{aligned}
            &\int_\bbR \px^k \Bigl( e_\ell^{\#}(x,\xi) - e_\ell^{\#}(x,-\gamma \ell) \Bigr) e^{i t (\jxi + \ell \xi)} \chi_{\{\leq 2\gamma|\ell|\}}(\xi) \, i\jxi^{-1} g(\xi) \, \ud \xi \\
            &= - \frac{1}{t} \int_\bbR \frac{ \px^k \bigl(e_\ell^{\#}(x,\xi) - e_\ell^{\#}(x,-\gamma \ell)\bigr)}{\xi+\ell\jxi} e^{i t (\jxi + \ell \xi)} \pxi \bigl( \chi_{\{\leq 2\gamma|\ell|\}}(\xi) \, g(\xi) \bigr) \, \ud \xi \\ 
            &\quad - \frac{1}{t} \int_\bbR \pxi \biggl( \frac{\px^k \bigl(e_\ell^{\#}(x,\xi) - e_\ell^{\#}(x,-\gamma \ell)\bigr)}{\xi+\ell\jxi} \biggr) e^{i t (\jxi + \ell \xi)} \chi_{\{\leq 2\gamma|\ell|\}}(\xi) \, g(\xi) \, \ud \xi \\ 
            &= I(t,x) + II(t,x).
        \end{aligned}
    \end{equation*}
    Next, we write
    \begin{equation*}
        \frac{\px^k \bigl( e_\ell^{\#}(x,\xi) - e_\ell^{\#}(x,-\gamma \ell) \bigr)}{\xi+\ell\jxi} = \biggl( \frac{1}{\xi+\gamma\ell} \int_{-\gamma\ell}^\xi \partial_\eta \px^k e_\ell^{\#}(x,\eta) \, \ud \eta \biggr) \frac{\xi+\gamma\ell}{\xi+\ell\jxi}    
    \end{equation*}
    as well as
    \begin{equation*}
        \begin{aligned}
        \pxi \biggl( \frac{ \px^k \bigl( e_\ell^{\#}(x,\xi) - e_\ell^{\#}(x,-\gamma \ell) \bigr)}{\xi+\ell\jxi} \biggr) 
        &= \biggl( \frac{1}{(\xi+\gamma\ell)^2} \int_{-\gamma \ell}^{\xi} \eta \, \bigl( \partial_\eta^2 \px^k e_\ell^{\#}(x,\eta) \bigr) \, \ud \eta \biggr) \frac{\xi+\gamma\ell}{\xi+\ell\jxi} \\ 
        &\quad + \biggl( \frac{\px^k \bigl( e_\ell^{\#}(x,\xi) - e_\ell^{\#}(x,-\gamma \ell) \bigr)}{\xi+\ell\jxi} \biggr) \pxi \biggl( \frac{\xi+\gamma\ell}{\xi+\ell\jxi} \biggr).
        \end{aligned}
    \end{equation*}    
    Invoking the bounds
    \begin{equation*}
        \sup_{x, \eta \in \bbR} \, \jx^{-1} \jap{\eta}^{-k} \bigl| \partial_\eta \px^k e_\ell^{\#}(x,\eta) \bigr| + \sup_{x, \eta \in \bbR} \, \jx^{-2} \jap{\eta}^{-k} \bigl| \partial_\eta^2 \px^k e_\ell^{\#}(x,\eta) \bigr| \lesssim_{\ell, k} 1,
    \end{equation*}    
    it follows that 
    \begin{equation*}
        \begin{aligned}
            \sup_{x,\xi \in \bbR} \, \Biggl( \biggl| \jx^{-1} \jxi^{-k} \frac{\px^k \bigl( e_\ell^{\#}(x,\xi) - e_\ell^{\#}(x,-\gamma \ell) \bigr)}{\xi+\ell\jxi} \biggr| + \biggl| \jx^{-2} \jxi^{-k} \pxi \biggl( \frac{\px^k \bigl( e_\ell^{\#}(x,\xi) - e_\ell^{\#}(x,-\gamma \ell) \bigr)}{\xi+\ell\jxi} \biggr) \biggr| \Biggr) \lesssim_\ell 1.
        \end{aligned}
    \end{equation*}
    Moreover, we have 
    \begin{equation*}
        \sup_{\xi \in \bbR} \, \biggl|\pxi \biggl( \frac{\xi+\gamma\ell}{\xi+\ell\jxi} \biggr)\biggr| \lesssim_{\ell, k} 1.
    \end{equation*}    
    Thus, using the Cauchy-Schwarz inequality in the frequency variable $\xi$ and exploiting the finite support of the cut-off $\chi_{\{\leq 2\gamma|\ell|\}}(\xi)$, we conclude 
    \begin{equation*}
        \begin{aligned}
            \bigl\|\jx^{-3} I(t,x)\bigr\|_{L^\infty_x} + \bigl\|\jx^{-3} II(t,x)\bigr\|_{L^2_x} &\lesssim \|\jx^{-1}\|_{L^2_x} \Bigl( \bigl\|\jx^{-2} I(t,x)\bigr\|_{L^\infty_x} + \bigl\|\jx^{-2} II(t,x)\bigr\|_{L^\infty_x} \Bigr) \\
            &\lesssim \frac{1}{t} \bigl( \|\pxi g(\xi)\|_{L^2_\xi} + \|g(\xi)\|_{L^2_\xi} \bigr),
        \end{aligned}
    \end{equation*}
    as desired. 

    The proof of \eqref{equ:local_decay_resonance_subtracted_off_for_2nd_comp} proceeds analogously, using the bounds 
    \begin{equation*}
        \sup_{x, \eta \in \bbR} \, \jx^{-1} \bigl| \partial_\eta \bigl( \jap{\eta}^{-1} (D^\ast e_\ell^{\#})(x,\eta) \bigr) \bigr| + \sup_{x, \eta \in \bbR} \, \jx^{-2} \bigl| \partial_\eta^2 \bigl( \jap{\eta}^{-1} (D^\ast e_\ell^{\#})(x,\eta) \bigr) \bigr| \lesssim_{\ell} 1.
    \end{equation*}  
    We note that due to the low frequency cut-off no frequency weights fall onto the input.
\end{proof}

\subsection{Integrated local energy decay with a moving center} \label{subsec:ILED}

The derivation of the weighted energy estimates in Subsection~\ref{subsec:weighted_energy_estimates} for the contributions of spatially localized terms with cubic-type time decay hinges on integrated local energy decay estimates with a moving center, which are stated on the frequency side in the next proposition.
We emphasize that incorporating the moving center here into these otherwise well-known estimates is at the expense of additional derivatives on the inputs.
While we hope the results from this subsection to be of independent interest, we point the reader to the proof of Lemma~\ref{lem:pxi_calR} to see how these integrated local energy decay estimates enter the nonlinear analysis in this paper.

\begin{proposition} \label{prop:ILED_moving_center_frequency_formulation}
    Fix $\ell \in (-1,1)$. Let $\fraka \in C^\infty(\bbR)$ with $\|\partial_x^k \fraka\|_{L^\infty_x} \lesssim_k 1$ for all integers $k \geq 0$ and let $\frakb \in L^{\infty}_\xi(\bbR)$.
    Suppose that $\theta \colon [0,\infty) \to \bbR$ satisfies $\theta(0) = 0$ and $|\theta'(s)| \lesssim \varepsilon \js^{-1+\delta}$ for all $s \geq 0$ for some $0 < \varepsilon, \delta \ll 1$.
    Then we have for all $t \geq 0$ that
    \begin{equation} \label{equ:ILED_moving_center_frequency_side}
        \begin{aligned}
            &\biggl\| \int_0^t e^{-i s (\jxi + \ell \xi)} (\xi+\ell \jxi) \jxi^{-1} e^{-i\xi\theta(s)} \biggl[ \int_\bbR e^{-ix\xi} \fraka(x) \frakb(\xi) F(s,x) \, \ud x \biggr] \, \ud s \biggr\|_{L^2_\xi} \\ 
            &\qquad \qquad \qquad \qquad \qquad \qquad \qquad \qquad \qquad \qquad \qquad \qquad \lesssim \bigl\| \jx^2 F(s,x) \bigr\|_{L^2_s([0,t]; H^1_x)},
        \end{aligned}
    \end{equation}
    as well as 
    \begin{equation} \label{equ:ILED_moving_center_frequency_side_with_jxi2}
        \begin{aligned}
            &\biggl\| \jxi^2 \int_0^t e^{-i s (\jxi + \ell \xi)} (\xi+\ell \jxi) \jxi^{-1} e^{-i\xi\theta(s)} \biggl[ \int_\bbR e^{-ix\xi} \fraka(x) \frakb(\xi) F(s,x) \, \ud x \biggr] \, \ud s \biggr\|_{L^2_\xi} \\ 
            &\qquad \qquad \qquad \qquad \qquad \qquad \qquad \qquad \qquad \qquad \qquad \qquad \lesssim \bigl\| \jx^2 F(s,x) \bigr\|_{L^2_s([0,t]; H^3_x)}.
        \end{aligned}
    \end{equation}
    Moreover, if $\frakc(\xi)$ satisfies $\|\jxi \frakc(\xi)\|_{L^\infty_\xi} \lesssim 1$, then the following variant holds for all $t \geq 0$,
    \begin{equation} \label{equ:ILED_moving_center_frequency_side_with_jxi2_plus_decaying_coefficient}
        \begin{aligned}
            &\biggl\| \jxi^2 \int_0^t e^{-i s (\jxi + \ell \xi)} (\xi+\ell \jxi) \jxi^{-1} e^{-i\xi\theta(s)} \biggl[ \int_\bbR e^{-ix\xi} \fraka(x) \frakc(\xi) F(s,x) \, \ud x \biggr] \, \ud s \biggr\|_{L^2_\xi} \\ 
            &\qquad \qquad \qquad \qquad \qquad \qquad \qquad \qquad \qquad \qquad \qquad \qquad \lesssim \bigl\| \jx^2 F(s,x) \bigr\|_{L^2_s([0,t]; H^2_x)},
        \end{aligned}
    \end{equation}
    while if $\frakd(\xi)$ satisfies $\|\jxi^{-1} \frakd(\xi)\|_{L^\infty_\xi} \lesssim 1$, then we have for all $t \geq 0$,
    \begin{equation} \label{equ:ILED_moving_center_frequency_side_with_jxi2_plus_growing_coefficient}
        \begin{aligned}
            &\biggl\| \jxi^2 \int_0^t e^{-i s (\jxi + \ell \xi)} (\xi+\ell \jxi) \jxi^{-1} e^{-i\xi\theta(s)} \biggl[ \int_\bbR e^{-ix\xi} \fraka(x) \frakd(\xi) F(s,x) \, \ud x \biggr] \, \ud s \biggr\|_{L^2_\xi} \\ 
            &\qquad \qquad \qquad \qquad \qquad \qquad \qquad \qquad \qquad \qquad \qquad \qquad \lesssim \bigl\| \jx^2 F(s,x) \bigr\|_{L^2_s([0,t]; H^4_x)},
        \end{aligned}
    \end{equation}    
\end{proposition}

The spatial weights and the Sobolev regularity on the right-hand sides of \eqref{equ:ILED_moving_center_frequency_side}, \eqref{equ:ILED_moving_center_frequency_side_with_jxi2}, \eqref{equ:ILED_moving_center_frequency_side_with_jxi2_plus_decaying_coefficient}, and \eqref{equ:ILED_moving_center_frequency_side_with_jxi2_plus_growing_coefficient} are not optimal, but they are convenient and sufficient for the nonlinear analysis.

The proof of Proposition~\ref{prop:ILED_moving_center_frequency_formulation} is a consequence of the following integrated local energy decay estimate with a moving center, which we hope to be of independent interest.
The proof is inspired by the proof of Lemma~9.6 in \cite{NakSchlag12}. See also Lemma~4.5 in \cite{CJ1} for a related computation.

\begin{lemma} \label{lem:ILED}
Fix $\ell \in (-1,1)$. Let $1>\nu>0$ and let $\sigma \geq \frac{3}{2\nu}$. 
Suppose that $\theta \colon [0,\infty) \to \bbR$ satisfies $\theta(0) = 0$ and $|\theta'(t)| \lesssim \varepsilon \jt^{-1+\delta}$ for all $t \geq 0$ for some $0 < \varepsilon, \delta \ll 1$. Define 
\begin{equation*}
    \bigl[\mathrm{T}_{\theta}(t,s) f\bigr](x) := f\bigl(x+(\theta(t)-\theta(s))\bigr).
\end{equation*}
Then we have the homogeneous integrated local energy decay estimate
\begin{equation} \label{equ:homogeneous_ILED} 
    \begin{aligned}
	   \biggl\| \jx^{-\sigma} e^{it(\jDx+\ell D_x)} (D_x +\ell \jap{D_x})\jDx^{-1} \bigl[ \mathrm{T}_{\theta}(t,0) f \bigr] \biggr\| _{L^2_t([0,\infty); L^2_x)} &\lesssim \|f\|_{H_{x}^{\frac{\nu}{2}}},
    \end{aligned}
\end{equation}
and for all $t \geq 0$ the dual homogeneous estimate
\begin{equation} \label{equ:dual_homogeneous_ILED}
    \begin{aligned}
        \biggl\| \int_0^t e^{-is(\jDx+\ell D_x)} (D_x+\ell \jDx) \jDx^{-1} \bigl[ \mathrm{T}_{\theta}(0,s) F(s) \bigr] \, \ud s \biggr\|_{L^2_x} \lesssim \bigl\| \jx^\sigma F \bigr\|_{L^2_s([0,t]; H^{\frac{\nu}{2}}_x)}.
    \end{aligned}
\end{equation}
\end{lemma}

We first present a succinct proof of Proposition~\ref{prop:ILED_moving_center_frequency_formulation} using Lemma~\ref{lem:ILED}. 

\begin{proof}[Proof of Proposition \ref{prop:ILED_moving_center_frequency_formulation}]
    By Plancherel's theorem and the dual homogeneous estimate \eqref{equ:dual_homogeneous_ILED} with the admissible parameter choices $\nu=\frac34$ and $\sigma = 2$, we obtain
    \begin{equation*}
        \begin{aligned}
            &\biggl\| \int_0^t e^{-i s (\jxi + \ell \xi)} (\xi+\ell \jxi) \jxi^{-1} e^{-i\xi\theta(s)} \biggl[ \int_\bbR e^{-ix\xi} \fraka(x) \frakb(\xi) F(s,x) \, \ud x \biggr] \, \ud s \biggr\|_{L^2_\xi} \\ 
            &\lesssim \|\frakb\|_{L^\infty_\xi} \biggl\| \int_0^t e^{-i s (\jxi + \ell \xi)} (\xi+\ell \jxi) \jxi^{-1} e^{-i\xi\theta(s)} \widehat{\calF}\bigl[ \fraka(\cdot) F(s,\cdot) \bigr](\xi) \, \ud s \biggr\|_{L^2_\xi} \\
            &\lesssim \biggl\| \int_0^t e^{-is(\jDx+\ell D_x)} (D_x+\ell\jDx)\jDx^{-1} \bigl[ T_\theta(0,s) \fraka(\cdot) F(s,\cdot) \bigr] \, \ud s \biggr\|_{L^2_x} \\ 
            &\lesssim \bigl\| \jx^2 \fraka(x) F(s,x) \bigr\|_{L^2_s([0,t]; H^{\frac38}_x)} \lesssim \bigl\| \jx^2 F(s,x) \bigr\|_{L^2_s([0,t];H^1_x)}.
        \end{aligned}
    \end{equation*}
    This proves \eqref{equ:ILED_moving_center_frequency_side}. Then writing 
    \begin{equation*}
        \jxi^2 \widehat{\calF}\bigl[ \fraka(\cdot) F(s,\cdot) \bigr](\xi) = \widehat{\calF}\Bigl[\jDx^2 \bigl(\fraka(\cdot) F(s,\cdot)\bigr) \Bigr](\xi),
    \end{equation*}
    and bounding 
    \begin{equation*}
        \bigl\| \jx^2 \jDx^2 \bigl(\fraka(\cdot) F(s,\cdot)\bigr) \bigr\|_{L^2_s([0,t];H^{\frac34}_x)} \lesssim \bigl\| \jx^2 F(s,x) \bigr\|_{L^2_s([0,t];H^3_x)},
    \end{equation*}
    the second asserted estimate \eqref{equ:ILED_moving_center_frequency_side_with_jxi2} follows analogously.
    The remaining asserted estimates \eqref{equ:ILED_moving_center_frequency_side_with_jxi2_plus_decaying_coefficient} and \eqref{equ:ILED_moving_center_frequency_side_with_jxi2_plus_growing_coefficient} also follow analogously by either absorbing one $\jxi$ weight using $\|\jxi \frakc(\xi)\|_{L^\infty_\xi} \lesssim 1$, or by picking up an additional $\jxi$ weight to access the bound $\|\jxi^{-1} \frakd(\xi)\|_{L^\infty_\xi} \lesssim 1$.
\end{proof}

The proof of Lemma~\ref{lem:ILED} can in turn be reduced to the following technical lemma.

\begin{lemma} \label{lem:ILED_proof_reduced_to}
    Fix $\ell \in (-1,1)$. Let $1>\nu>0$ and let $\sigma \geq \frac{3}{2\nu}$. Set $w(x) := \jx^{-\sigma}$.
    Suppose that $\theta \colon [0,\infty) \to \bbR$ satisfies $\theta(0) = 0$ and $|\theta'(t)| \lesssim \varepsilon \jt^{-1+\delta}$ for all $t \geq 0$ for some $0 < \varepsilon, \delta \ll 1$. Define 
    \begin{equation*}
        \bigl[\mathrm{T}_{\theta}(t,s) f\bigr](x) := f\bigl(x+(\theta(t)-\theta(s))\bigr).
    \end{equation*}    
    For any $F \in L^2_t([0,\infty);H^{\nu}_x)$ and any $G \in L^2_t([0,\infty); L^2_x)$ it holds that
    \begin{equation} \label{equ:ILED_proof_reduced_to1}
        \begin{aligned}
            &\biggl| \int_0^\infty \int_0^t \left\langle e^{i (t-s) (\jDx+\ell D_x)} (D_x+\ell \jDx)^2 \jDx^{-2} \bigl[ \mathrm{T}_{\theta}(t,s) w F(s) \bigr], w G(t) \right\rangle \, \ud s \, \ud t \biggr| \\
            &\qquad \qquad \qquad \qquad \qquad \qquad \qquad \qquad \qquad \qquad \lesssim \|F\|_{L^2_t([0,\infty);H^{\nu}_x)} \|G\|_{L^2_t([0,\infty); L^2_x)},
        \end{aligned}
    \end{equation}
    as well as 
    \begin{equation} \label{equ:ILED_proof_reduced_to2}
        \begin{aligned}
            &\biggl| \int_0^\infty \int_t^\infty \left\langle e^{i (t-s) (\jDx+\ell D_x)} (D_x+\ell \jDx)^2 \jDx^{-2} \bigl[ \mathrm{T}_{\theta}(t,s) w F(s) \bigr], w G(t) \right\rangle \, \ud s \, \ud t \biggr| \\
            &\qquad \qquad \qquad \qquad \qquad \qquad \qquad \qquad \qquad \qquad \lesssim \|F\|_{L^2_t([0,\infty);H^{\nu}_x)} \|G\|_{L^2_t([0,\infty); L^2_x)}.
        \end{aligned}
    \end{equation}    
\end{lemma}

Next, we present a proof of Lemma~\ref{lem:ILED} based on the statement of the preceding Lemma~\ref{lem:ILED_proof_reduced_to}. The proof of Lemma~\ref{lem:ILED_proof_reduced_to} will then occupy the remainder of this subsection.

\begin{proof}[Proof of Lemma~\ref{lem:ILED}]
We first observe that for any $t \geq 0$, the dual homogeneous estimate \eqref{equ:dual_homogeneous_ILED} follows by duality from the homogeneous estimate \eqref{equ:homogeneous_ILED} restricted to the time interval $[0,t]$.

For the proof of \eqref{equ:homogeneous_ILED} we introduce the operator
\begin{equation*}
    (\calK f)(t,x) := \jx^{-\sigma} \Bigl( e^{it(\jDx+\ell D_x)} (D_x+\ell \jDx)\jDx^{-1} \bigl[ \mathrm{T}_{\theta}(t,0) \jDx^{- \frac{\nu}{2}} f \bigr] \Bigr)(x).
\end{equation*}
Then the homogeneous integrated local energy decay estimate \eqref{equ:homogeneous_ILED} is equivalent to the boundedness of $\calK$ as an operator 
\begin{equation*}
    \calK \colon L^2_x \to L^2_t\bigl([0,\infty); L^2_x\bigr).
\end{equation*}
The adjoint of $\calK$ as an operator 
\begin{equation*}
    \calK^\ast \colon L^2_t\bigl([0,\infty); L^2_x\bigr) \to L^2_x
\end{equation*}
is given by
\begin{equation*}
    (\calK^\ast F)(x) := \biggl( \int_0^\infty e^{-is(\jDx+\ell D_x)} (D_x+\ell \jDx)\jDx^{-1} \bigl[ \mathrm{T}_{\theta}(0,s) \jDy^{-\frac{\nu}{2}} \langle \cdot \rangle^{-\sigma} F(s,\cdot) \bigr] \, \ud s \biggr)(x).
\end{equation*}
Thus, by standard duality arguments, proving the boundedness of $\calK \colon L^2_x \to L^2_t\bigl([0,\infty); L^2_x\bigr)$ is equivalent to proving the boundedness of 
\begin{equation} \label{equ:ILED_proof_calKcalKast_operator}
    \calK \calK^\ast \colon L^2_t\bigl([0,\infty); L^2_x\bigr) \to L^2_t\bigl([0,\infty); L^2_x\bigr)
\end{equation}
given by 
\begin{equation*} 
    \begin{aligned}
        &\bigl( \calK \calK^\ast F \bigr)(t,x) \\ 
        &\quad = \jx^{-\sigma} \biggl( \int_0^\infty e^{i(t-s)(\jDx+\ell D_x)} (D_x+\ell \jDx)^2 \jDx^{-2} \bigl[ \mathrm{T}_{\theta}(t,s) \jDy^{-\nu} \langle \cdot \rangle^{-\sigma} F(s,\cdot) \bigr] \, \ud s \biggr)(x).
    \end{aligned}
\end{equation*}
Distinguishing the time integration $0 \leq s \leq t$ and $t \leq s < \infty$, the boundedness \eqref{equ:ILED_proof_calKcalKast_operator} is in turn a consequence of the following two estimates
\begin{equation*} 
    \biggl\| \jx^{-\sigma} \int_0^t e^{i(t-s)(\jDx+\ell D_x)} (D_x+\ell \jDx)^2 \jDx^{-2} \bigl[ \mathrm{T}_{\theta}(t,s)  F(s,\cdot) \bigr] \, \ud s \biggr\|_{L^2_t([0,\infty); L^2_x)} \lesssim \bigl\| \jx^\sigma F \bigr\|_{L^2_t([0,\infty); H^\nu_x)}
\end{equation*}
and 
\begin{equation*} 
    \biggl\| \jx^{-\sigma} \int_t^\infty e^{i(t-s)(\jDx+\ell D_x)} (D_x+\ell \jDx)^2 \jDx^{-2} \bigl[ \mathrm{T}_{\theta}(t,s)  F(s,\cdot) \bigr] \, \ud s \biggr\|_{L^2_t([0,\infty); L^2_x)} \lesssim \bigl\| \jx^\sigma F \bigr\|_{L^2_t([0,\infty); H^\nu_x)}.
\end{equation*}    
By duality these two estimates follow directly from Lemma~\ref{lem:ILED_proof_reduced_to}.
\end{proof}

The proof of Lemma~\ref{lem:ILED_proof_reduced_to} requires some preparations.
In what follows, we will provide the details of the proof of the first asserted estimate \eqref{equ:ILED_proof_reduced_to1} in the statement of Lemma~\ref{lem:ILED_proof_reduced_to}, and we leave the analogous proof of the second asserted estimate \eqref{equ:ILED_proof_reduced_to2} to the reader. 

We assume throughout that $\theta \colon [0,\infty) \to \bbR$ satisfies $\theta(0) = 0$ and $|\theta'(t)| \lesssim \varepsilon \jt^{-1+\delta}$ for all $t \geq 0$ for some $0 < \varepsilon, \delta \ll 1$, and we consider spatial translations defined by
\begin{equation*}
    \bigl[\mathrm{T}_{\theta}(t,s) f\bigr](x) := f\bigl(x+(\theta(t)-\theta(s))\bigr).
\end{equation*}  

Let $\psi \in C_c^\infty(\bbR)$ be a smooth non-negative even bump function with $\psi(\xi) = 1$ for $|\xi| \leq 1$ and $\psi(\xi) = 0$ for $|\xi| \geq 2$. Set $\varphi(\xi) := \psi(\xi) - \psi(2\xi)$.
We denote by $\sum_{j=0}^{\infty} P_j = 1$ a standard dyadic Littlewood-Paley decomposition defined by 
\begin{align*}
    P_j u &:= \mathcal{F}^{-1}\bigl[ \varphi(2^{-j} \xi) \hatu(\xi) \bigr], \quad j \geq 1, \\ 
    P_0 u &:= \mathcal{F}^{-1}\bigl[ \psi(\xi) \hatu(\xi) \bigr].
\end{align*}

For any integer $j \geq 0$, any $M \geq 0$, and  any $F, G \in L^2_t L^2_x([0,\infty);L^2_x)$, we define
\begin{equation} \label{equ:ILED_proof_IjM_definition}
    \begin{aligned}
        &I_j^M(F, G) \\ 
        &:= \int_{0}^{\infty} \int_{0}^{\max\left\{ t-M, 0 \right\} } \Bigl\langle e^{i (t-s) (\jD+\ell D)}  (D+\ell \jD)^2 \jD^{-2} \bigl[ \mathrm{T}_{\theta}(t,s) P_j \bigl( w F(s) \bigr) \bigr], w G(t) \Bigr\rangle \, \ud s \, \ud t.
    \end{aligned}
\end{equation}
Later on we will have to sum over $j \geq 0$. To this end we will use the following almost orthogonality property
\begin{equation} \label{equ:ILED_proof_almost_orthogonality}
	I_{j}^{M}(F, G) = I_{j}^{M}\bigl( Q_j F, Q_j G \bigr)
\end{equation}
with
\begin{equation*}
 Q_j F := \sum_{\left|k-j\right| \leq 5} \Bigl( P_k F + w^{-1} [ P_k, w ] F \Bigr),
\end{equation*}
where the commutator term is small in the sense that
\begin{equation} \label{equ:ILED_proof_commutator_bound}
 \left\Vert w^{-1} [ P_k, w ]f \right \Vert_{L^2_x} \lesssim 2^{-k} \|f\|_{L^2_x}.
\end{equation}
For given $\ell \in (-1,1)$ we denote by $j_0$ the smallest integer such that
\begin{equation} \label{equ:ILED_proof_j0_definition}
    \frac{2^{j_0-1}}{\langle 2^{j_0-1} \rangle } \geq 1-\frac{1}{4}(1-|\ell|).
\end{equation}

The next lemma quantifies that the spatial shift caused by $T_\theta(t,s)$ is much smaller than $|t-s|$ when $|t-s| \gtrsim 1$.
\begin{lemma} \label{lem:thetadiff} 
Suppose that $\theta \colon [0,\infty) \to \bbR$ satisfies $\theta(0) = 0$ and $|\theta'(t)| \lesssim \varepsilon \jt^{-1+\delta}$ for all $t \geq 0$ for some $0 < \varepsilon, \delta \ll 1$.
Then we have for all $t \geq s \geq 0$,
\begin{equation} \label{eq:thetaTS}
    \left| \theta(t)-\theta(s) \right| \lesssim \varepsilon |t-s|^\delta.
\end{equation}
\end{lemma}
\begin{proof}
By the fundamental theorem of calculus we have 
\begin{align} \label{eq:diffestitheta}
    \left|\theta(t)-\theta(s)\right| \leq \int_{s}^t |\theta'(\tau)|\,\ud \tau \lesssim \int_{s}^t \varepsilon \langle \tau\rangle^{-1+\delta} \, \ud \tau \leq \int_{s}^t \varepsilon \tau^{-1+\delta} \, \ud \tau \lesssim \varepsilon (t^\delta-s^\delta).
\end{align}
To estimate $t^\delta-s^\delta$, we note that
\begin{equation*}
    t^\delta-s^\delta = |t-s|^\delta \left( \Bigl(1+ \frac{s}{t-s}\Bigr)^\delta - \Bigl( \frac{s}{t-s} \Bigr)^\delta \right),
\end{equation*}
so it suffices to bound the expression in the parenthesis. 
Setting $h(x) := (1+ x)^\delta - x^\delta$ for $x\geq0$, we have for all $x > 0$,
\begin{equation*}
    h'(x) = \delta \left( \bigl(1+ x\bigr)^{-1+\delta}-x^{-1+\delta}\right)<0.
\end{equation*}
It follows that $0 < h(x) \leq h(0)=1$ for all $x \geq 0$. Thus, we have for all $0 \leq s \leq t$,
\begin{equation*}
    t^\delta-s^\delta = |t-s|^\delta\left( \Bigl( 1+ \frac{s}{t-s} \Bigr)^\delta - \Bigl( \frac{s}{t-s} \Bigr)^\delta \right)\leq |t-s|^\delta,
\end{equation*}
which together with \eqref{eq:diffestitheta} implies the asserted estimate \eqref{eq:thetaTS}.
\end{proof}

Next, we establish pointwise decay estimates for the integral kernels of frequency localized Klein-Gordon evolutions.

\begin{lemma} \label{lem:ILED_proof_kernels_Kj_decay}
Fix $\ell \in (-1,1)$. 
Define the kernels
\begin{align}
        K_j(t,x) &:= \int_\bbR e^{ix\xi} e^{it(\jxi+\ell\xi)} (\xi + \ell \jxi)^2 \jxi^{-2} \varphi(2^{-j} \xi) \, \ud \xi, \quad j \geq 1, \label{equ:ILED_proof_Kj_kernel} \\
        K_{0}(t,x) &:= \int_\bbR e^{ix\xi} e^{it(\jxi+\ell\xi)} (\xi + \ell \jxi)^2 \jxi^{-2} \psi(\xi) \, \ud \xi. \label{equ:ILED_proof_K0_kernel}
\end{align}   
For $j \geq j_0$ with $j_0$ defined in \eqref{equ:ILED_proof_j0_definition}, we have for any integer $N \geq 1$ uniformly for all $t \geq 1$ the dispersive decay estimates
\begin{equation} \label{equ:ILED_proof_Kj_dispersive_decay}
    \bigl| K_j(t,x) \bigr| \lesssim_{\ell, N} \begin{cases}
	2^{j} (2^j t)^{-N},  &\bigl|\frac{x}{t}+\ell\bigr| \geq 1+\frac{3}{4} (1-|\ell|) \,\, \, \textrm{or} \, \, \, \bigl |\frac{x}{t}+\ell\bigr|\leq 1-\frac{3}{4} (1-|\ell|),\\
	2^{\frac32 j} t^{-\frac12},  &1-\frac{3}{4}(1-|\ell|) \leq \bigl |\frac{x}{t}+\ell\bigr|\leq 1+\frac{3}{4} (1-|\ell|).
	\end{cases}
\end{equation}
Moreover, for any $j \geq 0$ we have the improved local decay estimate
\begin{equation} \label{equ:ILED_proof_improved_local_decay}
	\bigl| \jx^{-1} K_j(t,x) \bigr| \lesssim 2^{\frac{3}{2}j} \jt^{-\frac{3}{2}}.
\end{equation}    
\end{lemma}
\begin{proof}
Assume $t \geq 1$.
We begin with the proof of the dispersive decay estimates \eqref{equ:ILED_proof_Kj_dispersive_decay}.
Write 
\begin{equation*}
    K_j(t,x) = \int_\bbR e^{i t \phi_\ell(\xi;v)} (\xi + \ell \jxi)^2 \jxi^{-2} \varphi(2^{-j} \xi) \, \ud \xi
\end{equation*}
with 
\begin{equation*}
    \phi_\ell(\xi; v) := \jxi + (v+\ell)\xi, \quad v := \frac{x}{t}.
\end{equation*}
\noindent {\bf Case 1: Away from the tilted light cone.} 
For $|v+\ell| \geq 1 + \frac34 (1-|\ell|)$ or for $|v+\ell| \leq 1 - \frac34 (1-|\ell|)$, we find by direct computation that for $j \geq j_0$ on the support of $\varphi(2^{-j}\xi)$ it holds that
\begin{equation*}
    \bigl|\pxi \phi_\ell(\xi;v)\bigr| \geq \frac34 (1-|\ell|).
\end{equation*}
Moreover, we have $\pxi^2 \phi_\ell(\xi;v) = \jxi^{-3}$.
Integrating by parts in $\xi$ repeatedly, we obtain for any integer $N \geq 1$ that
\begin{equation*}
    K_j(t,x) = \int_\bbR e^{it\phi_\ell(\xi;v)} L^N \Bigl( (\xi+\ell\jxi)^2 \jxi^{-2} \varphi(2^{-j} \xi) \Bigr) \, \ud \xi
\end{equation*}
with 
\begin{equation*}
    L := - \pxi \biggl( \frac{1}{i t \pxi \phi_\ell} \, \cdot \biggr).
\end{equation*}
Noting that $\pxi \bigl( (\pxi \phi_\ell)^{-1} \bigr) = - \jxi^{-3} (\pxi \phi_\ell)^{-2}$ and that $\pxi \bigl( (\xi+\ell\jxi) \jxi^{-1} \bigr) = \jxi^{-3}$, we conclude 
\begin{equation*}
    \bigl| K_j(t,x) \bigr| \lesssim_{\ell, N} (2^j t)^{-N} 2^j.
\end{equation*}

\noindent {\bf Case 2: Close to the tilted light cone.} 
For $1 - \frac34 (1-|\ell|) \leq |v+\ell| \leq 1 + \frac34 (1-|\ell|)$, the phase $\phi_\ell(\xi;v)$ possibly has a stationary point on the support of $\varphi(2^{-j} \xi)$. 
Here we proceed as in the proof of \cite[Lemma 2.2]{LS1}.
We have the trivial bound $|K_j(t,x)| \leq C 2^j$ for all $t > 0$ and $x \in \bbR$. 
For $t \geq 2^j$ we show the stronger bound $|K_j(t,x)| \leq C 2^{\frac32 j} t^{-\frac12}$. 
Now write 
\begin{equation} \label{equ:ILED_proof_dispersive_bounds_Kj_rescaled_time_s}
    K_j(t,x) = 2^j \int_\bbR e^{i t \phi_{\ell,j}(\xi;v)} m_\ell(2^j \xi) \varphi(\xi) \, \ud \xi
\end{equation}
with 
\begin{equation*}
    \phi_{\ell,j}(\xi;v) := 2^{2j} \bigl( 2^{-j} \jap{2^j \xi} + (v+\ell)\xi \bigr), \quad m_\ell(\zeta) := (\zeta + \ell\jap{\zeta}) \jap{\zeta}^{-1}. 
\end{equation*}
Then we have on the support of $\varphi$, which is contained in $[-2,-2^{-1}] \cup [2^{-1}, 2]$,
\begin{equation} \label{equ:ILED_proof_dispersive_bounds_near_lightcone_list}
    \begin{aligned}
        \pxi \phi_{\ell,j}(\xi;v) &= 2^{2j} \Bigl( \jap{2^j \xi}^{-1} (2^j \xi) + (v+\ell) \Bigr), \\
        \pxi^2 \phi_{\ell, j}(\xi;v) &= 2^{3j} \jap{2^j \xi}^3 \simeq 1, \\
        \bigl| \pxi^3 \phi_{\ell,j}(\xi;v) \bigr| &= 3 (2^{5j} |\xi|) \jap{2^j \xi}^{-5} \simeq 1.
    \end{aligned}
\end{equation}
Without loss of generality, by symmetry we can assume that the support of $\varphi$ is contained in $[2^{-1}, 2]$. We now distinguish the cases 
\begin{itemize}
 \item[(a)] $\mathrm{min} \, |\pxi \phi_{\ell,j}(\xi;v)| \gtrsim s^{-\frac12}$ on $[2^{-1}, 2]$,
 \item[(b)] $\mathrm{min} \, |\pxi \phi_{\ell,j}(\xi;v)| \ll s^{-\frac12}$ on $[2^{-1}, 2]$. 
\end{itemize}
In case (a) we conclude from the second derivative in \eqref{equ:ILED_proof_dispersive_bounds_near_lightcone_list} that 
\begin{equation*}
    \bigl| \pxi \phi_{\ell,j}(\xi;v) \bigr| \gtrsim s^{-\frac12} + \mathrm{min} \bigl\{ |\xi-2^{-1}|, |\xi-2| \bigr\} \quad \text{for all } \xi \in \bigl[2^{-1}, 2\bigr].
\end{equation*}
Integrating by parts once in \eqref{equ:ILED_proof_dispersive_bounds_Kj_rescaled_time_s} then gives the desired bound
\begin{equation*}
    \begin{aligned}
        \bigl| K_j(t,x) \bigr| \lesssim 2^j s^{-1} \int_{\frac12}^2 \biggl( \frac{|\pxi^2 \phi_{\ell,j}(\xi;v)|}{( \pxi \phi_{\ell,j}(\xi;v) )^2} + \frac{1}{|\pxi \phi_{\ell,j}(\xi;v)|} \biggr) \, \ud \xi \lesssim 2^j s^{-\frac12} \simeq 2^{\frac32 j} t^{-\frac12}. 
    \end{aligned}
\end{equation*}
In case (b) suppose the minimum of $|\pxi \phi_{\ell,j}(\xi;v)|$ is attained at the point $\xi_\ast \in [2^{-1}, 2]$. 
Then we infer from the second derivative bound in \eqref{equ:ILED_proof_dispersive_bounds_near_lightcone_list} that
\begin{equation*}
    |\pxi \phi_{\ell,j}(\xi;v)| \gtrsim |\xi-\xi_\ast| \quad \text{for } \xi \in \bigl[2^{-1}, 2\bigr] \text{ with } |\xi-\xi_\ast| \gtrsim s^{-\frac12}. 
\end{equation*}
With $\psi \in C_c^\infty(\bbR)$ the smooth even non-negative bump function defined previously, integrating by parts once in \eqref{equ:ILED_proof_dispersive_bounds_Kj_rescaled_time_s} when $|\xi-\xi_\ast| \gtrsim s^{-\frac12}$ yields
\begin{equation*}
    \begin{aligned}
        &\bigl| K_j(t,x) \bigr| \\
        &\lesssim 2^j \biggl| \int_0^\infty e^{is\phi_{\ell,j}(\xi;v)} \psi\bigl( s^{\frac12} (\xi-\xi_\ast) \bigr) m_\ell(2^j \xi) \varphi(\xi) \, \ud \xi \biggr| \\ 
        &\quad + 2^j s^{-1} \biggl| \int_0^\infty e^{is\phi_{\ell,j}(\xi;v)} \pxi \biggl( \bigl(\pxi \phi_{\ell,j}\bigr)^{-1} \Bigl( 1 - \psi\bigl( s^{\frac12} (\xi-\xi_\ast) \bigr) \Bigr) m_\ell(2^j \xi) \varphi(\xi) \biggr) \, \ud \xi \biggr| \\ 
        &\lesssim 2^j s^{-\frac12} + 2^j s^{-1} \int_0^\infty \biggl( \frac{|\pxi^2 \phi_{\ell,j}|}{|\pxi \phi_{\ell,j}|^2} \Bigl( 1 - \psi\bigl( s^{\frac12} (\xi-\xi_\ast) \bigr) \Bigr) + \frac{1}{|\pxi \phi_{\ell,j}|} \Bigl( s^{\frac12} \bigl| \psi'\bigl( s^{\frac12} (\xi-\xi_\ast) \bigr) \bigr| + |\varphi(\xi)| \Bigr) \biggr) \, \ud \xi \\ 
        &\lesssim 2^j s^{-\frac12} \simeq 2^{\frac32 j} t^{-\frac12}.
    \end{aligned}
\end{equation*}
This concludes the proof of the dispersive decay estimates \eqref{equ:ILED_proof_Kj_dispersive_decay}.

Finally, we establish the improved local decay estimate \eqref{equ:ILED_proof_improved_local_decay}. We treat the case $j \geq 1$, the case $j=0$ being analogous. 
Observing that $\pxi (\jxi + \ell \xi) = (\xi+\ell\jxi) \jxi^{-1}$, we first integrate by parts in $\xi$ in \eqref{equ:ILED_proof_Kj_kernel} to obtain
\begin{equation*}
    \begin{aligned}
        K_j(t,x) &= -\frac{1}{it} \int_\bbR e^{it(\jxi+\ell\xi)} \pxi \biggl( e^{ix\xi} (\xi + \ell \jxi) \jxi^{-1} \varphi(2^{-j} \xi) \biggr) \, \ud \xi \\ 
        &= -\frac{x}{it} \int_\bbR e^{it(\jxi+\ell\xi)} e^{ix\xi} (\xi + \ell \jxi) \jxi^{-1} \varphi(2^{-j} \xi) \, \ud \xi \\
        &\quad - \frac{1}{it} \int_\bbR e^{it(\jxi+\ell\xi)} e^{ix\xi} \pxi \bigl( (\xi + \ell \jxi) \jxi^{-1} \varphi(2^{-j} \xi) \bigr) \, \ud \xi.
    \end{aligned}
\end{equation*}
Then we run a stationary phase argument as in the preceding derivation of the dispersive decay estimates \eqref{equ:ILED_proof_Kj_dispersive_decay} to conclude the asserted improved local decay estimate \eqref{equ:ILED_proof_improved_local_decay}.
\end{proof}

Next, we derive key bounds on the interaction terms $I_j^M(F,G)$ defined in \eqref{equ:ILED_proof_IjM_definition}.

\begin{lemma} \label{lem:ILED_proof_Ij_bounds}
Let $\sigma>\frac{3}{2}$ and $M\geq1$.
Then for any $j\geq j_0 $ with $j_0$ defined in \eqref{equ:ILED_proof_j0_definition}, we have
	\begin{equation} \label{equ:ILED_proof_Ij_bounds_largej}
		\bigl| I_{j}^{M}(F, G)\bigr| \leq C(\sigma) 2^{\frac{3}{2}j} M^{1-\sigma} \|F\|_{L^2_t L^2_x} \|G\|_{L^2_t L^2_x}.
	\end{equation}
    Moreover, for $0 \leq j \leq j_0$ it holds that
	\begin{equation} \label{equ:ILED_proof_Ij_bounds_smallj}
	\bigl|I_{j}^{M}(F, G)\bigr| \leq C(\sigma) 2^{\frac{3}{2}j} M^{-\frac{1}{2}+\delta} \|F\|_{L^2_t L^2_x} \|G\|_{L^2_t L^2_x}.
\end{equation}
\end{lemma}
\begin{proof}
	We write
	\begin{equation} \label{equ:ILED_proof_Ij_as_sum_of_Jj}
	   I_{j}^{M}(F,G) = J_{j}^{M}(F,G) + \overline{J_{j}^{M}(G,F)}
	\end{equation}
	with 
    \begin{equation*} 
        \begin{aligned}
		      J_{j}^{M}(F,G) &:= \iint_{\{|y| > |x|\}} \int_{s=0}^\infty \int_{t=s+M}^\infty K_{j}(t-s,x-y-b) w(y) F(s,y) w(x) \overline{G}(t,x) \, \ud t \, \ud s \, \ud y \, \ud x
        \end{aligned}
    \end{equation*}
    and 
    \begin{equation*}
        b = b\left(t,s\right) := \int_s^t \theta'\left(\tau\right) \, \ud \tau.
    \end{equation*}
    We set 
    \begin{equation*}
        u := t-s, \quad z := x-y-b, \quad h := y+b.
    \end{equation*}
	Then we have 
	\begin{equation*}
	    h-b=y, \quad h+z=x.
	\end{equation*}
	In the new variables, the integral reads
    \begin{equation} \label{equ:ILED_proof_Jj_in_new_variables}
    \begin{aligned}
		J_{j}^{M}(F,G) &= \iint_{\{\left|h-b\right|>\left|h+z\right|\}} \int_{s=0}^\infty \int_{u=M}^\infty K_j(u,z) w(h-b) F(s,h-b) \\ 
        &\qquad \qquad \qquad \qquad \qquad \qquad \qquad \times w(h+z) \overline{G}(s+u,h+z) \, \ud u \, \ud s \, \ud h \, \ud z.
	\end{aligned}
    \end{equation}

We first consider the case $0 \leq j \leq j_0$. Using the improved local decay estimate \eqref{equ:ILED_proof_improved_local_decay} we bound \eqref{equ:ILED_proof_Jj_in_new_variables} by
\begin{equation*}
    \begin{aligned}
		\bigl| J_{j}^{M}(F,G) \bigr| &\lesssim \iint_{\{|h-b|>|h+z|\}} \int_{s=0}^\infty \int_{u=M}^\infty \langle u \rangle^{-\frac{3}{2}} \bigl| \langle z \rangle w(h-b) F(s,h-b) \\
        &\qquad \qquad \qquad \qquad \qquad \qquad \qquad \qquad \times w(h+z) \overline{G}(s+u,h+z) \bigr| \, \ud u \, \ud s \, \ud h \, \ud z.
    \end{aligned}
\end{equation*}
From the triangle inequality we obtain $\langle z\rangle \lesssim \langle h-b\rangle+\langle h+z\rangle+\langle b \rangle$, and thus
\begin{equation} \label{equ:ILED_proof_firstboundJ01}
	\begin{aligned}
		&\left|J_{j}^{M}(F,G)\right| \\ 
		&\quad \lesssim \iint_{\{\left|h-b\right|>\left|h+z\right|\}} \int_{s=0}^\infty \int_{u=M}^\infty \langle u\rangle^{-\frac{3}{2}} \bigl|\langle h-b\rangle w(h-b) F(s,h-b) \\
        &\qquad \qquad \qquad \qquad \qquad \qquad \qquad \qquad \times w(h+z) \overline{G}(s+u,h+z) \bigr| \, \ud u \, \ud s \, \ud h \, \ud z \\
		&\quad \quad + \iint_{\{\left|h-b\right|>\left|h+z\right|\}} \int_{s=0}^\infty \int_{u=M}^\infty \langle u\rangle^{-\frac{3}{2}} \bigl| w(h-b) F(s,h-b) \\
        &\qquad \qquad \qquad \qquad \qquad \qquad \qquad \qquad \times \langle h+z\rangle w(h+z) \overline{G}(s+u,h+z) \bigr| \, \ud u \, \ud s \, \ud h \, \ud z \\
		&\quad \quad + \iint_{\{\left|h-b\right|>\left|h+z\right|\}} \int_{s=0}^\infty \int_{u=M}^\infty \langle u\rangle^{-\frac{3}{2}} \langle b\rangle \bigl| w(h-b) F(s,h-b) \\
        &\qquad \qquad \qquad \qquad \qquad \qquad \qquad \qquad \times w(h+z) \overline{G}(s+u,h+z) \bigr| \, \ud u \, \ud s \, \ud h \, \ud z \\
		&\quad \lesssim \int_{u=M}^\infty \langle u\rangle^{-\frac{3}{2}+\delta} \|F\|_{L_{t}^{2}L_{x}^{2}} \|G\|_{L_{t}^{2}L_{x}^{2}} \|\langle x \rangle w(x)\|^2_{L^2_x} \,\ud u \\
		&\quad \lesssim M^{-\frac{1}{2}+\delta} \|F\|_{L^2_t L^2_x} \|G\|_{L^2_t L^2_x}
  \end{aligned}
\end{equation}
Note that in the second to last inequality, we applied the Cauchy-Schwarz inequality and used that $|b(t,s)| \lesssim \varepsilon |t-s|^\delta$ by Lemma~\ref{lem:thetadiff}. We also used that $x \mapsto \jx w(x)$ is $L^2_x$-integrable, which requires $-\sigma+1<-\frac{1}{2}$, whence the assumption $\sigma>\frac{3}{2}$.
The preceding bound \eqref{equ:ILED_proof_firstboundJ01} directly implies the asserted estimate \eqref{equ:ILED_proof_Ij_bounds_smallj} for $0 \leq j \leq j_0-1$.

Next we study the case $j \geq j_0$. 
We divide the range of integration in \eqref{equ:ILED_proof_Jj_in_new_variables} into two regions. Region $I$ corresponds to $\bigl|\frac{z}{u}+\ell\bigr| \geq 1+\frac{3}{4} (1-|\ell|) \, \textrm{or} \, \bigl|\frac{z}{u}+\ell\bigr| \leq 1-\frac{3}{4}(1-|\ell|)$, while region $II$ is defined by the restriction $1-\frac{3}{4} (1-|\ell|) \leq \bigl |\frac{z}{u}+\ell\bigr| \leq 1+\frac{3}{4} (1-|\ell|)$.
To distinguish the two regions we use (sharp) cut-off functions $\chi_I(u,z)$ and $\chi_{II}(u,z)$.

\noindent
{\bf Region I:}
Using the dispersive decay estimate in \eqref{equ:ILED_proof_Kj_dispersive_decay} for region $I$, 
we obtain for any $N \in \bbN$ that the integral \eqref{equ:ILED_proof_Jj_in_new_variables} restricted to region $I$ is bounded by
\begin{equation} \label{equ:ILED_proof_firstboundJ}
	\begin{aligned}
        &\iint_{\{\left|h-b\right|>\left|h+z\right|\}} \int_{s=0}^\infty \int_{u=M}^\infty \chi_I(u,z) \bigl|K_j(u,z) w(h-b) F(s,h-b) \\ 
        &\qquad \qquad \qquad \qquad \qquad \qquad \qquad \qquad \qquad \times w(h+z) \overline{G}(s+u,h+z) \bigr| \, \ud u \, \ud s \, \ud h \, \ud z \\
	    &\lesssim \iint_{\{|h-b| > |h+z|\}} \int_{s=0}^\infty \int_{u=M}^\infty 2^j (2^j u)^{-N-1} \bigl|w(h-b) F(s,h-b) \\
        &\qquad \qquad \qquad \qquad \qquad \qquad \qquad \qquad \qquad \times w(h+z) \overline{G}(s+u,h+z)\bigr| \, \ud u \, \ud s \, \ud h \, \ud z \\
		&\lesssim 2^{-jN} M^{-N} \|F\|_{L^2_t L^2_x} \|G\|_{L^2_t L^2_x}.
	\end{aligned}
\end{equation}
 
\noindent
{\bf Region II:}
If $1-\frac{3}{4}(1-|\ell|) \leq \bigl |\frac{z}{u}+\ell\bigr|\leq 1+\frac{3}{4}(1-|\ell|)$, we have
\begin{equation*}
 \frac{1}{4} (1-|\ell|) |u|\leq |z|\leq 2|u|.
\end{equation*}
It follows that
\[
 |b(t,s)| \lesssim \varepsilon |t-s|^{\delta} \ll u \simeq |z| \lesssim |h-b| + |h+z| + |b|.
\]
Since $|h+z| < |h-b|$ by construction and since $|b| \ll |u|$, the preceding estimate implies
\[
 u \lesssim |h-b|.
\]
Thus, by the dispersive decay estimate in \eqref{equ:ILED_proof_Kj_dispersive_decay} for region $II$, 
the integral \eqref{equ:ILED_proof_Jj_in_new_variables} restricted to region $II$ is bounded by
\begin{equation} \label{equ:ILED_proof_secondboundJ}
	\begin{aligned}
        &\iint_{\{\left|h-b\right|>\left|h+z\right|\}} \int_{s=0}^\infty \int_{u=M}^\infty \chi_{II}(u,z) \, \chi_{\{ u \, \lesssim \, |h-b|\}} \bigl|K_j(u,z) w(h-b) F(s,h-b) \\ 
        &\qquad \qquad \qquad \qquad \qquad \qquad \qquad \qquad \qquad \qquad \times w(h+z) \overline{G}(s+u,h+z) \bigr| \, \ud u \, \ud s \, \ud h \, \ud z \\
	    &\lesssim \iint_{\{\left|h-b\right|>\left|h+z\right|\}} \int_{s=0}^\infty \int_{u=M}^\infty \chi_{II}(u,z) \, \chi_{\{ u \, \lesssim \, |h-b|\}} \, 2^{\frac32 j} u^{-\frac12} \jap{h-b}^{-\sigma} |F(s, h-b)| \\
        &\qquad \qquad \qquad \qquad \qquad \qquad \qquad \qquad \qquad \qquad \times \jap{h+z}^{-\sigma} |G(s+u,h+z)| \, \ud u \, \ud s \, \ud h \, \ud z \\
        &\lesssim \int_{s=0}^\infty \int_{u=M}^\infty 2^{\frac32 j} u^{-\frac12} \bigl\| \jap{h-b}^{-\sigma} \jap{h+z}^{-\sigma} \bigr\|_{L^2_h(\{|h-b| \, \gtrsim \, u\}) \, L^2_z(\{ |z| \, \simeq \, u\})} \\
        &\qquad \qquad \qquad \qquad \qquad \qquad \qquad \qquad \qquad \qquad \times \|F(s,\cdot)\|_{L^2_x} \|G(s+u, \cdot)\|_{L^2_x} \, \ud s \, \ud u \\
        &\lesssim \int_{s=0}^\infty \int_{u=M}^\infty 2^{\frac32 j} u^{-\sigma} \|F(s,\cdot)\|_{L^2_x} \|G(s+u, \cdot)\|_{L^2_x} \, \ud s \, \ud u \\
        &\lesssim 2^{\frac32 j} M^{1-\sigma} \|F\|_{L^2_t L^2_x} \|G\|_{L^2_t L^2_x}.
    \end{aligned}
\end{equation}
Combining \eqref{equ:ILED_proof_firstboundJ} and \eqref{equ:ILED_proof_secondboundJ} yields the asserted estimate \eqref{equ:ILED_proof_Ij_bounds_largej} for $j \geq j_0$.
\end{proof}

Finally, we are in the position to establish the proof of Lemma~\ref{lem:ILED_proof_reduced_to}.

\begin{proof}[Proof of Lemma \ref{lem:ILED_proof_reduced_to}]
As indicated earlier, we prove the first asserted estimate \eqref{equ:ILED_proof_reduced_to1}, and we leave the analogous proof of the second asserted estimate \eqref{equ:ILED_proof_reduced_to2} to the reader.

Inserting the Littlewood-Paley decomposition $\sum_{j=0}^\infty P_j = 1$ and using the almost orthogonality \eqref{equ:ILED_proof_almost_orthogonality}, we have 
    \begin{equation} \label{equ:ILED_proof_reduced_to_chopped_up}
        \begin{aligned}
            &\biggl| \int_0^\infty \int_0^t \left\langle e^{i (t-s) (\jD+\ell D)} (D+\ell \jD)^2 \jD^{-2} \bigl[ \mathrm{T}_{\theta}(t,s) w F(s) \bigr], w G(t) \right\rangle \, \ud s \, \ud t \biggr| \\
            &\leq \sum_{j=0}^\infty \, \bigl| I_j^0(F, G) \bigr| = \sum_{j=0}^\infty \, \bigl| I_j^0(Q_j F, Q_j G) \bigr|.
        \end{aligned}
    \end{equation}
For each $j \geq 0$ we write for some $M \geq 1$ to be chosen sufficiently large below (depending on $j$),
\begin{equation} \label{equ:ILED_proof_Ij_splitting}
    I_j^0(Q_j F, Q_j G) = \Bigl( I_j^0(Q_j F, Q_j G) - I_j^M(Q_j F, Q_j G) \Bigr) + I_j^M(Q_j F, Q_j G).
\end{equation}
Observe that in the expression in the parenthesis on the right-hand side of \eqref{equ:ILED_proof_Ij_splitting} the interaction is strong in the sense that $|t-s| \leq M$. In this regime we just apply the trivial $L^2_x$ conservation for the Klein-Gordon propagator,
\begin{equation} \label{equ:ILED_proof_Ij_split1}
    \bigl| I_j^0(Q_j F, Q_j G) - I_j^M(Q_j F, Q_j G) \bigr| \lesssim M \| Q_j F \|_{L^2_t L^2_x} \| Q_j G \|_{L^2_t L^2_x}.
\end{equation}
Instead, the interactions in the second term on the right-hand side of \eqref{equ:ILED_proof_Ij_splitting} are milder in the sense that $|t-s| \geq M$. Here we distinguish the cases $0 \leq j \leq j_0 -1$ and $j \geq j_0$.
For $0 \leq j \leq j_0 -1$ we have by \eqref{equ:ILED_proof_Ij_bounds_smallj} with the trivial choice $M = 1$,
\begin{equation} \label{equ:ILED_proof_Ij_split2}
    \bigl| I_j^{M = 1}(Q_j F, Q_j G) \bigr| \lesssim_\sigma 2^{\frac32 j_0} \| F \|_{L^2_t L^2_x} \| G \|_{L^2_t L^2_x}. 
\end{equation} 
For $j \geq j_0$ we obtain from \eqref{equ:ILED_proof_Ij_bounds_smallj} that 
\begin{equation} \label{equ:ILED_proof_Ij_split3}
    \bigl| I_j^M(Q_j F, Q_j G) \bigr| \lesssim_\sigma 2^{\frac32 j} M^{1-\sigma} \| Q_j F \|_{L^2_t L^2_x} \| Q_j G \|_{L^2_t L^2_x}.
\end{equation}
For each fixed $j \geq j_0$ we choose $M := 2^{\nu j}$. Since by assumption $\sigma \geq \frac{3}{2\nu}$ for $0 < \nu < 1$, we have $2^{\frac32 j} M^{1-\sigma} \leq M$, whence 
\begin{equation} \label{equ:ILED_proof_Ij_split4}
    \bigl| I_j^{M = 2^{\nu j}}(Q_j F, Q_j G) \bigr| \lesssim_\sigma 2^{\nu j} \| Q_j F \|_{L^2_t L^2_x} \| Q_j G \|_{L^2_t L^2_x}.
\end{equation}
From \eqref{equ:ILED_proof_Ij_splitting}, the estimates \eqref{equ:ILED_proof_Ij_split1}, \eqref{equ:ILED_proof_Ij_split2}, \eqref{equ:ILED_proof_Ij_split3}, \eqref{equ:ILED_proof_Ij_split4}, the commutator bound \eqref{equ:ILED_proof_commutator_bound}, and the assumption $0 < \nu < 1$, we conclude that the right-hand side of \eqref{equ:ILED_proof_reduced_to_chopped_up} is bounded by
\begin{equation*}
    \begin{aligned}
        &\sum_{j=0}^\infty \, \bigl| I_j^0(Q_j F, Q_j G) \bigr| \\
        &\lesssim \sum_{0 \leq j \leq j_0 - 1} \, 2^{\frac32 j_0} \| F \|_{L^2_t L^2_x} \| G \|_{L^2_t L^2_x} + \sum_{j \geq j_0} \, 2^{\nu j} \| Q_j F \|_{L^2_t L^2_x} \| Q_j G \|_{L^2_t L^2_x} \\ 
        &\lesssim_{j_0} \| F \|_{L^2_t L^2_x} \| G \|_{L^2_t L^2_x} + \sum_{j \geq j_0} \sum_{|k_1-j| \leq 5} \sum_{|k_2-j| \leq 5} \, 2^{\nu j} \| P_{k_1} F \|_{L^2_t L^2_x} \| P_{k_2} G \|_{L^2_t L^2_x} \\ 
        &\quad \quad + \sum_{j \geq j_0} 2^{(\nu-1) j} \| F \|_{L^2_t L^2_x} \| G \|_{L^2_t L^2_x} \\  
        &\lesssim_{j_0} \|F\|_{L^2_t([0,\infty);H^{\nu}_x)} \|G\|_{L^2_t([0,\infty); L^2_x)},
    \end{aligned}
\end{equation*}
as desired.
\end{proof}

\section{Modified Nonlinear Spectral Distributions} \label{sec:nonlinear_spectral_distributions}

As a preparation for the nonlinear analysis in the remainder of this paper, in this section we determine the structure of cubic spectral distributions in the evolution equation for the effective profile and of quadratic spectral distributions in the modulation equations.

\subsection{Cubic spectral distributions for the evolution equation of the effective profile} \label{subsec:cubic_spectral_distributions}

The contributions of the constant coefficient cubic nonlinearities on the right-hand sides of the evolution equations \eqref{equ:setting_up_g_evol_equ3} and \eqref{equ:setting_up_g_evol_equ4} for the effective profile give rise to the following cubic spectral distributions for some $\ell \in (0,1)$, as detailed in Subsection~\ref{subsec:structure_cubic_nonlinearities},
\begin{equation} \label{equ:cubic_spectral_distributions_subsec}
\begin{aligned}
    \nu_{\ell, +++}(\xi, \xi_1, \xi_2, \xi_3) &:= \int_\bbR \overline{\elsharp(y,\xi)} \, \elsharp(y,\xi_1) \, \elsharp(y,\xi_2) \, \elsharp(y,\xi_3) \, \ud y, \\
    \nu_{\ell, +-+}(\xi, \xi_1, \xi_2, \xi_3) &:= \int_\bbR \overline{\elsharp(y,\xi)} \, \elsharp(y,\xi_1) \, \overline{\elsharp(y,\xi_2)} \, \elsharp(y,\xi_3) \, \ud y, \\
    \nu_{\ell, +--}(\xi, \xi_1, \xi_2, \xi_3) &:= \int_\bbR \overline{\elsharp(y,\xi)} \, \elsharp(y,\xi_1) \, \overline{\elsharp(y,\xi_2)} \, \overline{\elsharp(y,\xi_3)} \, \ud y, \\
    \nu_{\ell, ---}(\xi, \xi_1, \xi_2, \xi_3) &:= \int_\bbR \overline{\elsharp(y,\xi)} \, \overline{\elsharp(y,\xi_1)} \, \overline{\elsharp(y,\xi_2)} \, \overline{\elsharp(y,\xi_3)} \, \ud y.
\end{aligned}
\end{equation}
Below we use the shorthand notation
\begin{equation*}
    \eta := \gamma (\xi + \ell \jxi), \quad \eta_j := \gamma (\xi_j + \ell \jap{\xi_j}).
\end{equation*}
By direct computation we obtain in the sense of distributions the following decompositions of the cubic spectral distributions \eqref{equ:cubic_spectral_distributions_subsec}. We omit the details of these lengthy, but straightforward computations. 
Inserting the definition of the distorted Fourier basis elements \eqref{equ:modified_dFT_basis_succinct} into \eqref{equ:cubic_spectral_distributions_subsec} and expanding the products, the basic idea is just to use the elementary identity $\tanh^2(\cdot) = 1 - \sech^2(\cdot)$ to write the resulting list of terms as linear combinations of the functions $1$, $\tanh(\gamma y)$, and $\sech^m(\gamma y) \tanh^n(\gamma y)$ with $m \in \{2,4\}$ and $n \in \{0,1\}$ multiplied by the corresponding products of frequency terms. The Fourier transform identities \eqref{equ:prelim_FT_one} and \eqref{equ:prelim_FT_tanh} then lead to the following decompositions.

\medskip 

\noindent \underline{Type I: $+++$ interactions.}
\begin{align*}
    \nu_{\ell, +++}(\xi, \xi_1, \xi_2, \xi_3) &= \frakm_{\ell, +++}^{\delta_0}(\xi, \xi_1, \xi_2, \xi_3) \, \delta_0\bigl(-\xi+\xi_1+\xi_2+\xi_3\bigr) \\ 
    &\quad + \frakm_{\ell, +++}^{\pvdots}(\xi, \xi_1, \xi_2, \xi_3) \, \pvdots \cosech\Bigl( \frac{\pi}{2\gamma} \bigl(-\xi+\xi_1+\xi_2+\xi_3\bigr) \Bigr) \\
    &\quad + \nu_{\ell, +++}^{\mathrm{reg}}(\xi, \xi_1, \xi_2, \xi_3),
\end{align*}
where
\begin{align*}
    \frakm_{\ell, +++}^{\delta_0}(\xi, \xi_1, \xi_2, \xi_3) &=  \frac{1}{2\pi} \frac{1}{|\eta|+i} \frac{1}{|\eta_1|-i} \frac{1}{|\eta_2|-i} \frac{1}{|\eta_3|-i} \\
    &\qquad \qquad \times \Bigl( \eta \eta_1 \eta_2 \eta_3 -\eta \eta_1 -\eta \eta_2 + \eta_1 \eta_2 - \eta \eta_3 + \eta_1 \eta_3 + \eta_2 \eta_3 - 1 \Bigr), \\
    \frakm_{\ell, +++}^{\pvdots}(\xi, \xi_1, \xi_2, \xi_3) &= \frac{1}{4\pi\gamma} \frac{1}{|\eta|+i} \frac{1}{|\eta_1|-i} \frac{1}{|\eta_2|-i} \frac{1}{|\eta_3|-i} \\ 
    &\qquad \qquad \times \Bigl( - \eta + \eta_1 + \eta_2 + \eta_3 + \eta \eta_1 \eta_2 + \eta \eta_1 \eta_3 + \eta \eta_2 \eta_3 - \eta_1 \eta_2 \eta_3 \Bigr),
\end{align*}
and where $\nu_{\ell, +++}^{\mathrm{reg}}(\xi, \xi_1, \xi_2, \xi_3)$ is a linear combination of terms of the form
\begin{equation*}
    \gamma^{-1} \frakb(\xi) \frakb_1(\xi_1) \frakb_2(\xi_2) \frakb_3(\xi_3) \, \widehat{\calF}^{-1}\bigl[\varphi\bigr]\bigl( \gamma^{-1} (-\xi+\xi_1+\xi_2+\xi_3)\bigr), 
\end{equation*}
with $\varphi \in \calS(\bbR)$ some Schwartz function and
\begin{equation*}
    \frakb(\xi) = \frac{1}{|\gamma(\xi+\ell\jxi)|+i} \quad \quad \text{or} \quad \quad \frakb(\xi) = \frac{\gamma (\xi + \ell \jxi)}{|\gamma(\xi+\ell\jxi)|+i},
\end{equation*}
as well as
\begin{equation*}
    \frakb_j(\xi_j) = \frac{1}{|\gamma(\xi_j+\ell\jap{\xi_j})|-i} \quad \quad \text{or} \quad \quad \frakb_j(\xi_j) = \frac{\gamma (\xi_j+\ell \jap{\xi_j})}{|\gamma(\xi_j +\ell\jap{\xi_j})|-i}, \quad 1 \leq j \leq 3.
\end{equation*}

\medskip 

\noindent \underline{Type II: $+-+$ interactions.}
\begin{align*}
    \nu_{\ell, +-+}(\xi, \xi_1, \xi_2, \xi_3) &= \frakm_{\ell, +-+}^{\delta_0}(\xi, \xi_1, \xi_2, \xi_3) \, \delta_0\bigl(-\xi+\xi_1-\xi_2+\xi_3\bigr) \\ 
    &\quad + \frakm_{\ell, +-+}^{\pvdots}(\xi, \xi_1, \xi_2, \xi_3) \, \pvdots \cosech\Bigl( \frac{\pi}{2\gamma} \bigl(-\xi+\xi_1-\xi_2+\xi_3\bigr) \Bigr) \\
    &\quad + \nu_{\ell, +-+}^{\mathrm{reg}}(\xi, \xi_1, \xi_2, \xi_3),
\end{align*}
where
\begin{align*}
    \frakm_{\ell, +-+}^{\delta_0}(\xi, \xi_1, \xi_2, \xi_3) &=  \frac{1}{2\pi} \frac{1}{|\eta|+i} \frac{1}{|\eta_1|-i} \frac{1}{|\eta_2|+i} \frac{1}{|\eta_3|-i} \\
    &\qquad \qquad \times \Bigl( \eta \eta_1 \eta_2 \eta_3 + \eta \eta_1 - \eta \eta_2 + \eta_1 \eta_2 + \eta \eta_3 - \eta_2 \eta_4 + \eta_3 \eta_4 + 1 \Bigr), \\
    \frakm_{\ell, +-+}^{\pvdots}(\xi, \xi_1, \xi_2, \xi_3) &= \frac{1}{4\pi\gamma} \frac{1}{|\eta|+i} \frac{1}{|\eta_1|-i} \frac{1}{|\eta_2|+i} \frac{1}{|\eta_3|-i} \\ 
    &\qquad \qquad \times \Bigl( \eta - \eta_1 + \eta_2 - \eta_3 + \eta \eta_1 \eta_2 - \eta \eta_1 \eta_3 + \eta \eta_2 \eta_3 - \eta_1 \eta_2 \eta_3 \Bigr),
\end{align*}
and where $\nu_{\ell, +-+}^{\mathrm{reg}}(\xi, \xi_1, \xi_2, \xi_3)$ is a linear combination of terms of the form
\begin{equation*}
    \gamma^{-1} \frakb(\xi) \frakb_1(\xi_1) \frakb_2(\xi_2) \frakb_3(\xi_3) \, \widehat{\calF}^{-1}\bigl[\varphi\bigr]\bigl( \gamma^{-1} (-\xi+\xi_1-\xi_2+\xi_3)\bigr), 
\end{equation*}
with $\varphi \in \calS(\bbR)$ some Schwartz function and
\begin{equation*}
    \frakb(\xi) = \frac{1}{|\gamma(\xi+\ell\jxi)|+i} \quad \quad \text{or} \quad \quad \frakb(\xi) = \frac{\gamma (\xi + \ell \jxi)}{|\gamma(\xi+\ell\jxi)|+i},
\end{equation*}
as well as $\frakb_j(\xi_j)$, $1 \leq j \leq 3$, given by one of the following expressions or complex conjugates thereof
\begin{equation*}
    \frac{1}{|\gamma(\xi_j+\ell\jap{\xi_j})|\pm i} \quad \quad \text{or} \quad \quad \frac{\gamma (\xi_j+\ell \jap{\xi_j})}{|\gamma(\xi_j +\ell\jap{\xi_j})| \pm i}.
\end{equation*}

\medskip 

\noindent We observe that 
\begin{equation} \label{equ:cubic_spectral_distributions_pmp_delta_identity}
    \frakm_{\ell, +-+}^{\delta_0}(\xi, \xi, \xi, \xi) = \frac{1}{2\pi} \quad \text{for all } \xi \in \bbR,
\end{equation}
and that we have the cancellation property
\begin{equation} \label{equ:cubic_spectral_distributions_pv_vanishing}
    \frakm_{\ell, +-+}^{\pvdots}(\xi, \xi, \xi, \xi) = 0 \quad \text{for all } \xi \in \bbR.
\end{equation}

\medskip 

\noindent \underline{Type III: $+--$ interactions.}
\begin{align*}
    \nu_{\ell, +--}(\xi, \xi_1, \xi_2, \xi_3) &= \frakm_{\ell, +--}^{\delta_0}(\xi, \xi_1, \xi_2, \xi_3) \, \delta_0\bigl(-\xi+\xi_1-\xi_2-\xi_3\bigr) \\ 
    &\quad + \frakm_{\ell, +--}^{\pvdots}(\xi, \xi_1, \xi_2, \xi_3) \, \pvdots \cosech\Bigl( \frac{\pi}{2\gamma} \bigl(-\xi+\xi_1-\xi_2-\xi_3\bigr) \Bigr) \\
    &\quad + \nu_{\ell, +--}^{\mathrm{reg}}(\xi, \xi_1, \xi_2, \xi_3),
\end{align*}
where
\begin{align*}
    \frakm_{\ell, +--}^{\delta_0}(\xi, \xi_1, \xi_2, \xi_3) &=  \frac{1}{2\pi} \frac{1}{|\eta|+i} \frac{1}{|\eta_1|-i} \frac{1}{|\eta_2|+i} \frac{1}{|\eta_3|+i} \\
    &\qquad \qquad \times \Bigl( \eta \eta_1 \eta_2 \eta_3 - \eta \eta_1 + \eta \eta_2 - \eta_1 \eta_2 + \eta \eta_3 - \eta_1 \eta_3 + \eta_2 \eta_3 \Bigr), \\
    \frakm_{\ell, +--}^{\pvdots}(\xi, \xi_1, \xi_2, \xi_3) &= \frac{1}{4\pi\gamma} \frac{1}{|\eta|+i} \frac{1}{|\eta_1|-i} \frac{1}{|\eta_2|+i} \frac{1}{|\eta_3|+i} \\ 
    &\qquad \qquad \times \Bigl( -\eta + \eta_1 - \eta_2 - \eta_3 -\eta \eta_1 \eta_2 - \eta \eta_1 \eta_3 + \eta \eta_2 \eta_3 - \eta_1 \eta_2 \eta_3  \Bigr),
\end{align*}
and where $\nu_{\ell, +--}^{\mathrm{reg}}(\xi, \xi_1, \xi_2, \xi_3)$ is a linear combination of terms of the form
\begin{equation*}
    \gamma^{-1} \frakb(\xi) \frakb_1(\xi_1) \frakb_2(\xi_2) \frakb_3(\xi_3) \, \widehat{\calF}^{-1}\bigl[\varphi\bigr]\bigl( \gamma^{-1} (-\xi+\xi_1-\xi_2-\xi_3)\bigr), 
\end{equation*}
with $\varphi \in \calS(\bbR)$ some Schwartz function and
\begin{equation*}
    \frakb(\xi) = \frac{1}{|\gamma(\xi+\ell\jxi)|+i} \quad \quad \text{or} \quad \quad \frakb(\xi) = \frac{\gamma (\xi + \ell \jxi)}{|\gamma(\xi+\ell\jxi)|+i},
\end{equation*}
as well as $\frakb_j(\xi_j)$, $1 \leq j \leq 3$, given by one of the following expressions 
\begin{equation*}
    \frac{1}{|\gamma(\xi_j+\ell\jap{\xi_j})|\pm i} \quad \quad \text{or} \quad \quad \frac{\gamma (\xi_j+\ell \jap{\xi_j})}{|\gamma(\xi_j +\ell\jap{\xi_j})| \pm i}.
\end{equation*}

\medskip 

\noindent \underline{Type IV: $---$ interactions.}
\begin{align*}
    \nu_{\ell, ---}(\xi, \xi_1, \xi_2, \xi_3) &= \frakm_{\ell, ---}^{\delta_0}(\xi, \xi_1, \xi_2, \xi_3) \, \delta_0\bigl(-\xi-\xi_1-\xi_2-\xi_3\bigr) \\ 
    &\quad + \frakm_{\ell, ---}^{\pvdots}(\xi, \xi_1, \xi_2, \xi_3) \, \pvdots \cosech\Bigl( \frac{\pi}{2\gamma} \bigl(-\xi-\xi_1-\xi_2-\xi_3\bigr) \Bigr) \\
    &\quad + \nu_{\ell, ---}^{\mathrm{reg}}(\xi, \xi_1, \xi_2, \xi_3),
\end{align*}
where
\begin{align*}
    \frakm_{\ell, ---}^{\delta_0}(\xi, \xi_1, \xi_2, \xi_3) &=  \frac{1}{2\pi} \frac{1}{|\eta|+i} \frac{1}{|\eta_1|+i} \frac{1}{|\eta_2|+i} \frac{1}{|\eta_3|+i} \\
    &\qquad \qquad \times \Bigl( \eta \eta_1 \eta_2 \eta_3 - \eta \eta_1 - \eta \eta_2 - \eta_1 \eta_2 - \eta \eta_3 - \eta_1 \eta_3 - \eta_2 \eta_3 + 1 \Bigr), \\
    \frakm_{\ell, +--}^{\pvdots}(\xi, \xi_1, \xi_2, \xi_3) &= \frac{1}{4\pi\gamma} \frac{1}{|\eta|+i} \frac{1}{|\eta_1|+i} \frac{1}{|\eta_2|+i} \frac{1}{|\eta_3|+i} \\ 
    &\qquad \qquad \times \Bigl( \eta + \eta_ 1 + \eta_2 + \eta_3 - \eta \eta_1 \eta_2 - \eta \eta_1 \eta_3 - \eta \eta_2 \eta_3 - \eta_1 \eta_2 \eta_3 \Bigr),
\end{align*}
and where $\nu_{\ell, ---}^{\mathrm{reg}}(\xi, \xi_1, \xi_2, \xi_3)$ is a linear combination of terms of the form
\begin{equation*}
    \gamma^{-1} \frakb(\xi) \frakb_1(\xi_1) \frakb_2(\xi_2) \frakb_3(\xi_3) \, \widehat{\calF}^{-1}\bigl[\varphi\bigr]\bigl( \gamma^{-1} (-\xi-\xi_1-\xi_2-\xi_3)\bigr), 
\end{equation*}
with $\varphi \in \calS(\bbR)$ some Schwartz function and
\begin{equation*}
    \frakb(\xi) = \frac{1}{|\gamma(\xi+\ell\jxi)|+i} \quad \quad \text{or} \quad \quad \frakb(\xi) = \frac{\gamma (\xi + \ell \jxi)}{|\gamma(\xi+\ell\jxi)|+i},
\end{equation*}
as well as $\frakb_j(\xi_j)$, $1 \leq j \leq 3$, given by one of the following expressions 
\begin{equation*}
    \frac{1}{|\gamma(\xi_j+\ell\jap{\xi_j})| + i} \quad \quad \text{or} \quad \quad \frac{\gamma (\xi_j+\ell \jap{\xi_j})}{|\gamma(\xi_j +\ell\jap{\xi_j})| + i}.
\end{equation*}

\subsection{Quadratic spectral distributions for the modulation equations}  \label{subsec:quadratic_spectral_distributions}

The contribution of the leading order quadratic nonlinearity on the right-hand side of the modulation equation for $\dot{\ell}$ gives rise to the following quadratic spectral distributions, as detailed in the discussion of the term $I(s)$ in Section~\ref{sec:modulation_control},
\begin{equation}
    \begin{aligned}
        \mu_{\ell;+-}(\xi_1, \xi_2) &:= \frac12 \int_\bbR \gamma \sech^2(\gamma y) \tanh(\gamma y) e_\ell^{\#}(y,\xi_1) \overline{e_\ell^{\#}(y,\xi_2)} \, \ud y, \label{eq:muell+-} \\     
        \mu_{\ell;++}(\xi_1, \xi_2) &:= \frac14 \int_\bbR \gamma \sech^2(\gamma y) \tanh(\gamma y) e_\ell^{\#}(y,\xi_1) e_\ell^{\#}(y,\xi_2) \, \ud y, \\
        \mu_{\ell;--}(\xi_1, \xi_2) &:= \frac14 \int_\bbR \gamma \sech^2(\gamma y) \tanh(\gamma y) \overline{e_\ell^{\#}(y,\xi_1)} \, \overline{e_\ell^{\#}(y,\xi_2)} \, \ud y.
    \end{aligned}
\end{equation}
In this subsection we compute these three quadratic spectral distributions explicitly. We uncover remarkable null structures for all them. The proof of Theorem~\ref{thm:main} hinges on the null structure for the quadratic spectral distribution $\mu_{\ell;+-}(\xi_1, \xi_2)$ stemming from the resonant leading order quadratic nonlinearity in the modulation equation for $\dot{\ell}$. The favorable structures of the other two quadratic spectral distributions  $\mu_{\ell;++}(\xi_1, \xi_2)$ and  $\mu_{\ell;--}(\xi_1, \xi_2)$ simplify the treatment of the other (non-resonant) leading order quadratic nonlinearities, but they are not essential.
Note that $\overline{\mu_{\ell;--}(\xi_1, \xi_2)} = \mu_{\ell;++}(\xi_1, \xi_2)$, so it suffices to consider $\mu_{\ell;++}(\xi_1, \xi_2)$ and $\mu_{\ell;+-}(\xi_1, \xi_2)$.

\begin{lemma} \label{lem:null_structure2} 
    For any $\ell \in (-1,1)$ we have 
    \begin{align} 
        \mu_{\ell;+-}(\xi_1, \xi_2) &= \bigl( (\jxione + \ell \xi_1) - (\jxitwo + \ell \xi_2) \bigr) \kappa_{\ell;+-}(\xi_1, \xi_2), \label{equ:null_structure2} \\
        \mu_{\ell;++}(\xi_1, \xi_2) &= \bigl( (\jxione + \ell \xi_1) + (\jxitwo + \ell \xi_2) \bigr) \kappa_{\ell;++}(\xi_1, \xi_2),  \label{equ:null_structure2++}
    \end{align}
    with 
    \begin{equation*}
        \begin{aligned}
            \kappa_{\ell;+-}(\xi_1, \xi_2) &:= \frac{i}{192 \pi} \frac{1}{|\gamma (\xi_1+\ell\jxione)|-i} \frac{1}{|\gamma (\xi_2+\ell\jxitwo)|+i}  \gamma^{-1} (\xi_1-\xi_2) \cosech\Bigl( \frac{\pi}{2} \gamma^{-1} (\xi_1-\xi_2) \Bigr) \times \\
            &\quad \quad \times \biggl( 3\gamma^2 \bigl( \jxione + \ell \xi_1 + \jxitwo + \ell \xi_2 \bigr) \bigl( \gamma (\xi_1 + \ell \jxione) + \gamma (\xi_2+\ell\jxitwo) \bigr) \\ 
            &\qquad \qquad - 3 \gamma \ell \Bigl( \gamma (\xi_1 + \ell \jxione) + \gamma (\xi_2+\ell\jxitwo) \Bigr)^2 \\ 
            &\qquad \qquad + \ell (\xi_1-\xi_2) \Bigl( -6 \bigl( \gamma (\xi_1 + \ell \jxione) - \gamma (\xi_2+\ell\jxitwo) \bigr) + 3 \gamma \ell \bigl( \jxione + \ell \xi_1 - \jxitwo - \ell \xi_2 \bigr) \Bigr) \\
            &\qquad \qquad + 4 \gamma \ell \bigl( -2 + \gamma^{-2}(\xi_1-\xi_2)^2 \bigr) \biggr).
        \end{aligned}
    \end{equation*}
    and 
    \begin{equation*}
     \begin{aligned}
            \kappa_{\ell;++}(\xi_1, \xi_2) &:= \frac{i}{192 \pi} \frac{1}{|\gamma (\xi_1+\ell\jxione)|-i} \frac{1}{|\gamma (\xi_2+\ell\jxitwo)|-i}  \gamma^{-1} (\xi_1+\xi_2) \cosech\Bigl( \frac{\pi}{2} \gamma^{-1} (\xi_1+\xi_2) \Bigr) \times \\
            &\quad \quad \times \biggl(-3\gamma^2 \bigl( \jxione + \ell \xi_1 - \jxitwo - \ell \xi_2 \bigr) \bigl( \gamma (\xi_1 + \ell \jxione) - \gamma (\xi_2+\ell\jxitwo) \bigr) \\ 
            &\qquad \qquad + 3 \gamma \ell \Bigl( \gamma (\xi_1 - \ell \jxione) + \gamma (\xi_2+\ell\jxitwo) \Bigr)^2 \\ 
            &\qquad \qquad + \ell (\xi_1-\xi_2) \Bigl( 6 \bigl( \gamma (\xi_1 + \ell \jxione) - \gamma (\xi_2+\ell\jxitwo) \bigr) - 3 \gamma \ell \bigl( \jxione + \ell \xi_1 + \jxitwo + \ell \xi_2 \bigr) \Bigr) \\
            &\qquad \qquad + 4 \gamma \ell \bigl( -2 + \gamma^{-2}(\xi_1+\xi_2)^2 \bigr) \biggr).
        \end{aligned}
    \end{equation*}
\end{lemma}
\begin{proof}
We begin with the computations for $\mu_{\ell;+-}(\xi_1, \xi_2)$.
    Inserting the definition of the distorted Fourier basis elements \eqref{equ:modified_dFT_basis} into  \eqref{eq:muell+-} and making the change of variables $\gamma y \mapsto y$, we find 
    \begin{equation*}
        \begin{aligned}
            &\mu_{\ell;+-}(\xi_1, \xi_2) \\ 
            &= \frac{1}{4\pi} \frac{1}{|\gamma (\xi_1+\ell\jxione)|-i} \frac{1}{|\gamma (\xi_2+\ell\jxitwo)|+i} \times \\
            &\quad \quad \times \int_\bbR e^{iy\gamma^{-1}(\xi_1-\xi_2)} \sech^2(y) \tanh(y) \bigl( \gamma (\xi_1 + \ell \jxione) + i \tanh(y) \bigr) \bigl( \gamma(\xi_2+\ell\jxitwo) - i \tanh(y) \bigr) \, \ud y \\ 
            &= \frac{1}{4\pi} \frac{1}{|\gamma (\xi_1+\ell\jxione)|-i} \frac{1}{|\gamma (\xi_2+\ell\jxitwo)|+i} \times \\
            &\quad \quad \times \biggl( \gamma (\xi_1 + \ell \jxione) \, \gamma (\xi_2+\ell\jxitwo) \int_\bbR e^{iy\gamma^{-1}(\xi_1-\xi_2)} \sech^2(y) \tanh(y) \, \ud y \\
            &\qquad \quad \quad - i \bigl( \gamma (\xi_1 + \ell \jxione) - \gamma (\xi_2+\ell\jxitwo) \bigr) \int_\bbR e^{iy\gamma^{-1}(\xi_1-\xi_2)} \sech^2(y) \tanh^2(y) \, \ud y \\ 
            &\qquad \quad \quad + \int_\bbR e^{iy\gamma^{-1}(\xi_1-\xi_2)} \sech^2(y) \tanh^3(y) \, \ud y \biggr). 
        \end{aligned}
    \end{equation*}
    Next, using the Fourier transform identities \eqref{equ:appendix_FT_sech2_tanh1}, \eqref{equ:appendix_FT_sech2_tanh2}, \eqref{equ:appendix_FT_sech2_tanh3},   
    we obtain that 
    \begin{equation} \label{equ:null_structure2_derivation0}
        \begin{aligned}
            &\mu_{\ell;+-}(\xi_1, \xi_2) \\ 
            &= \frac{i}{192 \pi} \frac{1}{|\gamma (\xi_1+\ell\jxione)|-i} \frac{1}{|\gamma (\xi_2+\ell\jxitwo)|+i}  \gamma^{-1} (\xi_1-\xi_2) \cosech\Bigl( \frac{\pi}{2} \gamma^{-1} (\xi_1-\xi_2) \Bigr) \fraka_\ell(\xi_1,\xi_2)
        \end{aligned}
    \end{equation}
    with 
    \begin{equation} \label{equ:null_structure2_derivation1}
        \begin{aligned}
            \fraka_\ell(\xi_1,\xi_2) &:= 12 \gamma^{-1} (\xi_1-\xi_2) \gamma (\xi_1 + \ell \jxione) \, \gamma (\xi_2+\ell\jxitwo) \\ 
            &\quad \quad + 4 \bigl( -2 + \gamma^{-2}(\xi_1-\xi_2)^2 \bigr) \bigl( \gamma (\xi_1 + \ell \jxione) - \gamma (\xi_2+\ell\jxitwo) \bigr) \\ 
            &\quad \quad - \gamma^{-1} (\xi_1-\xi_2) \bigl( -8 + \gamma^{-2} (\xi_1-\xi_2)^2 \bigr).
        \end{aligned}
    \end{equation}
    By direct computation we have
    \begin{equation} \label{equ:null_structure2_derivation2}
        \gamma (\xi_1 + \ell \jxione) - \gamma (\xi_2+\ell\jxitwo) = \gamma^{-1} (\xi_1 - \xi_2) + \gamma \ell \bigl( \jxione + \ell \xi_1 - \jxitwo - \ell \xi_2 \bigr).
    \end{equation}
    Inserting this identity in the second term on the right-hand side of \eqref{equ:null_structure2_derivation1} and then factoring out $\gamma^{-1} (\xi_1-\xi_2)$, we find that
    \begin{equation} \label{equ:null_structure2_derivation3}
        \begin{aligned}
            \fraka_\ell(\xi_1,\xi_2) &= \gamma^{-1} (\xi_1-\xi_2) \Bigl( 12 \gamma (\xi_1 + \ell \jxione) \, \gamma (\xi_2+\ell\jxitwo) + 4 \bigl( -2 + \gamma^{-2}(\xi_1-\xi_2)^2 \bigr) \\ 
            &\qquad \qquad \qquad \qquad \qquad \qquad \qquad \qquad \qquad \qquad \qquad - \bigl( -8 + \gamma^{-2} (\xi_1-\xi_2)^2 \bigr) \Bigr) \\ 
            &\quad + 4 \bigl( -2 + \gamma^{-2}(\xi_1-\xi_2)^2 \bigr) \gamma \ell \bigl( \jxione + \ell \xi_1 - \jxitwo - \ell \xi_2 \bigr) \\ 
            &= \gamma^{-1} (\xi_1-\xi_2) \Bigl( 12 \gamma (\xi_1 + \ell \jxione) \, \gamma (\xi_2+\ell\jxitwo) + 3 \gamma^{-2}(\xi_1-\xi_2)^2 \Bigr) \\ 
            &\quad + 4 \bigl( -2 + \gamma^{-2}(\xi_1-\xi_2)^2 \bigr) \gamma \ell \bigl( \jxione + \ell \xi_1 - \jxitwo - \ell \xi_2 \bigr). 
        \end{aligned}
    \end{equation}
    The second term on the right-hand side of \eqref{equ:null_structure2_derivation3} already features the desired phase factor. 
    In the first term on the right-hand side of \eqref{equ:null_structure2_derivation3} we insert the identity for $\gamma^{-1} (\xi_1-\xi_2)$ from \eqref{equ:null_structure2_derivation2} to obtain
    \begin{equation} \label{equ:null_structure2_derivation4}
        \begin{aligned}
            &\gamma^{-1} (\xi_1-\xi_2) \Bigl( 12 \gamma (\xi_1 + \ell \jxione) \, \gamma (\xi_2+\ell\jxitwo) + 3 \gamma^{-2}(\xi_1-\xi_2)^2 \Bigr) \\ 
            &= \gamma^{-1} (\xi_1-\xi_2) \biggl( 12 \gamma (\xi_1 + \ell \jxione) \, \gamma (\xi_2+\ell\jxitwo) \\ 
            &\qquad \qquad \qquad \quad \quad + 3 \Bigl( \gamma (\xi_1 + \ell \jxione) - \gamma (\xi_2+\ell\jxitwo) - \gamma \ell \bigl( \jxione + \ell \xi_1 - \jxitwo - \ell \xi_2 \bigr) \Bigr)^2 \biggr) \\ 
            &= \gamma^{-1} (\xi_1-\xi_2) \biggl( 12 \gamma (\xi_1 + \ell \jxione) \, \gamma (\xi_2+\ell\jxitwo) 
            + 3 \Bigl( \gamma (\xi_1 + \ell \jxione) - \gamma (\xi_2+\ell\jxitwo) \Bigr)^2 \biggr) \\ 
            &\quad + \gamma \ell \bigl( \jxione + \ell \xi_1 - \jxitwo - \ell \xi_2 \bigr) \times \\ 
            &\quad \quad \quad \times \gamma^{-1} (\xi_1-\xi_2) \Bigl( -6 \bigl( \gamma (\xi_1 + \ell \jxione) - \gamma (\xi_2+\ell\jxitwo) \bigr) + 3 \gamma \ell \bigl( \jxione + \ell \xi_1 - \jxitwo - \ell \xi_2 \bigr) \Bigr).            
        \end{aligned}
    \end{equation}
    We record that the second term on the right-hand side of \eqref{equ:null_structure2_derivation4} features the desired phase factor. In the first term on the right-hand side of \eqref{equ:null_structure2_derivation4} we now complete the square and then insert the identity for $\gamma^{-1} (\xi_1-\xi_2)$ from \eqref{equ:null_structure2_derivation2} again to find that
    \begin{equation} \label{equ:null_structure2_derivation5}
        \begin{aligned}
            &\gamma^{-1} (\xi_1-\xi_2) \biggl( 12 \gamma (\xi_1 + \ell \jxione) \, \gamma (\xi_2+\ell\jxitwo) 
            + 3 \Bigl( \gamma (\xi_1 + \ell \jxione) - \gamma (\xi_2+\ell\jxitwo) \Bigr)^2 \biggr) \\             
            &= 3 \gamma^{-1} (\xi_1-\xi_2) \Bigl( \gamma (\xi_1 + \ell \jxione) + \gamma (\xi_2+\ell\jxitwo) \Bigr)^2 \\
            &= 3 \Bigl( \gamma (\xi_1 + \ell \jxione) - \gamma (\xi_2+\ell\jxitwo) - \gamma \ell \bigl( \jxione + \ell \xi_1 - \jxitwo - \ell \xi_2 \bigr) \Bigr) \times \\ 
            &\qquad \qquad \qquad \qquad \qquad \qquad \qquad \qquad \times \Bigl( \gamma (\xi_1 + \ell \jxione) + \gamma (\xi_2+\ell\jxitwo) \Bigr)^2 \\
            &= 3 \Bigl( \gamma^2 (\xi_1 + \ell \jxione)^2 - \gamma^2 (\xi_2+\ell\jxitwo)^2 \Bigr) \Bigl( \gamma (\xi_1 + \ell \jxione) + \gamma (\xi_2+\ell\jxitwo) \Bigr) \\ 
            &\quad - 3 \gamma \ell \bigl( \jxione + \ell \xi_1 - \jxitwo - \ell \xi_2 \bigr) \Bigl( \gamma (\xi_1 + \ell \jxione) + \gamma (\xi_2+\ell\jxitwo) \Bigr)^2.
        \end{aligned}
    \end{equation}
    The second term on the right-hand side of \eqref{equ:null_structure2_derivation5} again features the desired phase factor. Using the observation
    \begin{equation*}
        \begin{aligned}
            \gamma^2 (\xi_1 + \ell \jxione)^2 - \gamma^2 (\xi_2+\ell\jxitwo)^2 = \gamma^2 (\jxione + \ell \xi_1)^2 - \gamma^2 (\jxitwo + \ell \xi_2)^2,
        \end{aligned}
    \end{equation*}
    we can now also uncover the desired phase factor in the first term on the right-hand side of \eqref{equ:null_structure2_derivation5},
    \begin{equation} \label{equ:null_structure2_derivation6}
        \begin{aligned}
            &3 \Bigl( \gamma^2 (\xi_1 + \ell \jxione)^2 - \gamma^2 (\xi_2+\ell\jxitwo)^2 \Bigr) \Bigl( \gamma (\xi_1 + \ell \jxione) + \gamma (\xi_2+\ell\jxitwo) \Bigr) \\ 
            &= 3 \gamma^2 \bigl( \jxione + \ell \xi_1 - \jxitwo - \ell \xi_2 \bigr) \bigl( \jxione + \ell \xi_1 + \jxitwo + \ell \xi_2 \bigr) \bigl( \gamma (\xi_1 + \ell \jxione) + \gamma (\xi_2+\ell\jxitwo) \bigr).
        \end{aligned}
    \end{equation}
    Hence, combining \eqref{equ:null_structure2_derivation3}, \eqref{equ:null_structure2_derivation4}, \eqref{equ:null_structure2_derivation5}, and \eqref{equ:null_structure2_derivation6}, we find that 
    \begin{equation*}
        \begin{aligned}
            \fraka_\ell(\xi_1,\xi_2) &= 3 \gamma^2 \bigl( \jxione + \ell \xi_1 - \jxitwo - \ell \xi_2 \bigr) \bigl( \jxione + \ell \xi_1 + \jxitwo + \ell \xi_2 \bigr) \bigl( \gamma (\xi_1 + \ell \jxione) + \gamma (\xi_2+\ell\jxitwo) \bigr) \\ 
            &\quad - 3 \gamma \ell \bigl( \jxione + \ell \xi_1 - \jxitwo - \ell \xi_2 \bigr) \Bigl( \gamma (\xi_1 + \ell \jxione) + \gamma (\xi_2+\ell\jxitwo) \Bigr)^2 \\ 
            &\quad + \gamma \ell \bigl( \jxione + \ell \xi_1 - \jxitwo - \ell \xi_2 \bigr) \times \\ 
            &\quad \quad \quad \times \gamma^{-1} (\xi_1-\xi_2) \Bigl( -6 \bigl( \gamma (\xi_1 + \ell \jxione) - \gamma (\xi_2+\ell\jxitwo) \bigr) + 3 \gamma \ell \bigl( \jxione + \ell \xi_1 - \jxitwo - \ell \xi_2 \bigr) \Bigr) \\ 
            &\quad + 4 \bigl( -2 + \gamma^{-2}(\xi_1-\xi_2)^2 \bigr) \gamma \ell \bigl( \jxione + \ell \xi_1 - \jxitwo - \ell \xi_2 \bigr) \\ 
            &= \bigl( \jxione + \ell \xi_1 - \jxitwo - \ell \xi_2 \bigr) \times \\ 
            &\quad \quad \times \biggl(3\gamma^2 \bigl( \jxione + \ell \xi_1 + \jxitwo + \ell \xi_2 \bigr) \bigl( \gamma (\xi_1 + \ell \jxione) + \gamma (\xi_2+\ell\jxitwo) \bigr) \\ 
            &\qquad \qquad - 3 \gamma \ell \Bigl( \gamma (\xi_1 + \ell \jxione) + \gamma (\xi_2+\ell\jxitwo) \Bigr)^2 \\ 
            &\qquad \qquad + \ell (\xi_1-\xi_2) \Bigl( -6 \bigl( \gamma (\xi_1 + \ell \jxione) - \gamma (\xi_2+\ell\jxitwo) \bigr) + 3 \gamma \ell \bigl( \jxione + \ell \xi_1 - \jxitwo - \ell \xi_2 \bigr) \Bigr) \\
            &\qquad \qquad + 4 \gamma \ell \bigl( -2 + \gamma^{-2}(\xi_1-\xi_2)^2 \bigr) \biggr).
        \end{aligned}
    \end{equation*}
    Inserting the preceding identity back into \eqref{equ:null_structure2_derivation0} yields the asserted identity \eqref{equ:null_structure2}.

    Next, we determine the quadratic spectral distribution $\mu_{\ell;+-}(\xi_1, \xi_2)$. 
    Proceeding as above, we have
    \begin{equation*}
        \begin{aligned}
            &\mu_{\ell;++}(\xi_1, \xi_2) \\ 
            &= \frac{1}{4\pi} \frac{1}{|\gamma (\xi_1+\ell\jxione)|-i} \frac{1}{|\gamma (\xi_2+\ell\jxitwo)|-i} \times \\
            &\quad \quad \times \int_\bbR e^{iy\gamma^{-1}(\xi_1+\xi_2)} \sech^2(y) \tanh(y) \bigl( \gamma (\xi_1 + \ell \jxione) + i \tanh(y) \bigr) \bigl( \gamma(\xi_2+\ell\jxitwo) + i \tanh(y) \bigr) \, \ud y \\ 
            &= \frac{1}{4\pi} \frac{1}{|\gamma (\xi_1+\ell\jxione)|-i} \frac{1}{|\gamma (\xi_2+\ell\jxitwo)|-i} \times \\
            &\quad \quad \times \biggl( \gamma (\xi_1 + \ell \jxione) \, \gamma (\xi_2+\ell\jxitwo) \int_\bbR e^{iy\gamma^{-1}(\xi_1+\xi_2)} \sech^2(y) \tanh(y) \, \ud y \\
            &\qquad \quad \quad + i \bigl( \gamma (\xi_1 + \ell \jxione) + \gamma (\xi_2+\ell\jxitwo) \bigr) \int_\bbR e^{iy\gamma^{-1}(\xi_1+\xi_2)} \sech^2(y) \tanh^2(y) \, \ud y \\ 
            &\qquad \quad \quad - \int_\bbR e^{iy\gamma^{-1}(\xi_1+\xi_2)} \sech^2(y) \tanh^3(y) \, \ud y \biggr). 
        \end{aligned}
    \end{equation*}
Using the Fourier transform identities \eqref{equ:appendix_FT_sech2_tanh1}, \eqref{equ:appendix_FT_sech2_tanh2}, \eqref{equ:appendix_FT_sech2_tanh3}, 
it follows that
    \begin{equation} \label{equ:null_structure2++_derivation0}
        \begin{aligned}
            &\mu_{\ell;++}(\xi_1, \xi_2) \\ 
            &= \frac{i}{192 \pi} \frac{1}{|\gamma (\xi_1+\ell\jxione)|-i} \frac{1}{|\gamma (\xi_2+\ell\jxitwo)|-i}  \gamma^{-1} (\xi_1+\xi_2) \cosech\Bigl( \frac{\pi}{2} \gamma^{-1} (\xi_1+\xi_2) \Bigr) \fraka_{\ell,++}(\xi_1,\xi_2)
        \end{aligned}
    \end{equation}
    with 
    \begin{equation} \label{equ:null_structure2++_derivation1}
        \begin{aligned}
            \fraka_{\ell,++}(\xi_1,\xi_2) &:= 12 \gamma^{-1} (\xi_1+\xi_2) \gamma (\xi_1 + \ell \jxione) \, \gamma (\xi_2+\ell\jxitwo) \\ 
            &\quad \quad - 4 \bigl( -2 + \gamma^{-2}(\xi_1+\xi_2)^2 \bigr) \bigl( \gamma (\xi_1 + \ell \jxione) +\gamma (\xi_2+\ell\jxitwo) \bigr) \\ 
            &\quad \quad +\gamma^{-1} (\xi_1+\xi_2) \bigl( -8 + \gamma^{-2} (\xi_1+\xi_2)^2 \bigr).
        \end{aligned}
    \end{equation}
Again, by direct computation we have
    \begin{equation} \label{equ:null_structure2++_derivation2}
        \gamma (\xi_1 + \ell \jxione) + \gamma (\xi_2+\ell\jxitwo) = \gamma^{-1} (\xi_1 + \xi_2) + \gamma \ell \bigl( \jxione + \ell \xi_1 + \jxitwo + \ell \xi_2 \bigr).
    \end{equation}
Inserting the identity above into the second term on the right-hand side of \eqref{equ:null_structure2++_derivation1} and then factoring out $\gamma^{-1} (\xi_1+\xi_2)$, we get that
\begin{equation} \label{equ:null_structure2++_derivation3}
        \begin{aligned}
            \fraka_{\ell,++}(\xi_1,\xi_2) &= \gamma^{-1} (\xi_1+\xi_2) \Bigl( 12 \gamma (\xi_1 + \ell \jxione) \, \gamma (\xi_2+\ell\jxitwo) - 4 \bigl( -2 + \gamma^{-2}(\xi_1+\xi_2)^2 \bigr) \\ 
            &\qquad \qquad \qquad \qquad \qquad \qquad \qquad \qquad \qquad \qquad \qquad + \bigl( -8 + \gamma^{-2} (\xi_1+\xi_2)^2 \bigr) \Bigr) \\ 
            &\quad - 4 \bigl( -2 + \gamma^{-2}(\xi_1+\xi_2)^2 \bigr) \gamma \ell \bigl( \jxione + \ell \xi_1 +\jxitwo + \ell \xi_2 \bigr) \\ 
            &= \gamma^{-1} (\xi_1+\xi_2) \Bigl( 12 \gamma (\xi_1 + \ell \jxione) \, \gamma (\xi_2+\ell\jxitwo) - 3 \gamma^{-2}(\xi_1+\xi_2)^2 \Bigr) \\ 
            &\quad - 4 \bigl( -2 + \gamma^{-2}(\xi_1+\xi_2)^2 \bigr) \gamma \ell \bigl( \jxione + \ell \xi_1+ \jxitwo +\ell \xi_2 \bigr). 
        \end{aligned}
    \end{equation}
The second term on the right-hand side of \eqref{equ:null_structure2++_derivation3} already contains the desired phase factor. 
    In the first term on the right-hand side of \eqref{equ:null_structure2++_derivation3} we insert the identity for $\gamma^{-1} (\xi_1+\xi_2)$ from \eqref{equ:null_structure2++_derivation2} to obtain
    \begin{equation} \label{equ:null_structure2++_derivation4}
        \begin{aligned}
            &\gamma^{-1} (\xi_1+\xi_2) \Bigl( 12 \gamma (\xi_1 + \ell \jxione) \, \gamma (\xi_2+\ell\jxitwo) - 3 \gamma^{-2}(\xi_1+\xi_2)^2 \Bigr) \\ 
            &= \gamma^{-1} (\xi_1+\xi_2) \biggl( 12 \gamma (\xi_1 + \ell \jxione) \, \gamma (\xi_2+\ell\jxitwo) \\ 
            &\qquad \qquad \qquad \quad \quad - 3 \Bigl( \gamma (\xi_1 + \ell \jxione) + \gamma (\xi_2+\ell\jxitwo) - \gamma \ell \bigl( \jxione + \ell \xi_1 + \jxitwo + \ell \xi_2 \bigr) \Bigr)^2 \biggr) \\ 
            &= \gamma^{-1} (\xi_1+\xi_2) \biggl( 12 \gamma (\xi_1 + \ell \jxione) \, \gamma (\xi_2+\ell\jxitwo) 
            - 3 \Bigl( \gamma (\xi_1 + \ell \jxione) + \gamma (\xi_2+\ell\jxitwo) \Bigr)^2 \biggr) \\ 
            &\quad + \gamma \ell \bigl( \jxione + \ell \xi_1 + \jxitwo + \ell \xi_2 \bigr) \times \\ 
            &\quad \quad \quad \times \gamma^{-1} (\xi_1+\xi_2) \Bigl( 6 \bigl( \gamma (\xi_1 + \ell \jxione) - \gamma (\xi_2+\ell\jxitwo) \bigr) - 3 \gamma \ell \bigl( \jxione + \ell \xi_1 + \jxitwo +\ell \xi_2 \bigr) \Bigr).            
        \end{aligned}
    \end{equation}
    We record that the second term on the right-hand side of \eqref{equ:null_structure2++_derivation4} has the desired phase factor. In the first term on the right-hand side of \eqref{equ:null_structure2++_derivation4} we now complete the square and then insert the identity for $\gamma^{-1} (\xi_1+\xi_2)$ from \eqref{equ:null_structure2++_derivation2} again to find that
    \begin{equation} \label{equ:null_structure2++_derivation5}
        \begin{aligned}
            &\gamma^{-1} (\xi_1+\xi_2) \biggl( 12 \gamma (\xi_1 + \ell \jxione) \, \gamma (\xi_2+\ell\jxitwo) 
            - 3 \Bigl( \gamma (\xi_1 + \ell \jxione) + \gamma (\xi_2+\ell\jxitwo) \Bigr)^2 \biggr) \\             
            &= -3 \gamma^{-1} (\xi_1+\xi_2) \Bigl( \gamma (\xi_1 + \ell \jxione) - \gamma (\xi_2+\ell\jxitwo) \Bigr)^2 \\
            &= -3 \Bigl( \gamma (\xi_1 + \ell \jxione) + \gamma (\xi_2+\ell\jxitwo) - \gamma \ell \bigl( \jxione + \ell \xi_1 + \jxitwo + \ell \xi_2 \bigr) \Bigr) \times \\ 
            &\qquad \qquad \qquad \qquad \qquad \qquad \qquad \qquad \times \Bigl( \gamma (\xi_1 + \ell \jxione) - \gamma (\xi_2+\ell\jxitwo) \Bigr)^2 \\
            &= -3 \Bigl( \gamma^2 (\xi_1 + \ell \jxione)^2 - \gamma^2 (\xi_2+\ell\jxitwo)^2 \Bigr) \Bigl( \gamma (\xi_1 + \ell \jxione) - \gamma (\xi_2+\ell\jxitwo) \Bigr) \\ 
            &\quad + 3 \gamma \ell \bigl( \jxione + \ell \xi_1 + \jxitwo + \ell \xi_2 \bigr) \Bigl( \gamma (\xi_1 + \ell \jxione) - \gamma (\xi_2+\ell\jxitwo) \Bigr)^2.
        \end{aligned}
    \end{equation}
    The second term on the right-hand side of \eqref{equ:null_structure2++_derivation5} again features the desired phase factor. Using the observation
    \begin{equation*}
        \begin{aligned}
            \gamma^2 (\xi_1 + \ell \jxione)^2 - \gamma^2 (\xi_2+\ell\jxitwo)^2 = \gamma^2 (\jxione + \ell \xi_1)^2 - \gamma^2 (\jxitwo + \ell \xi_2)^2,
        \end{aligned}
    \end{equation*}
    we can now also uncover the desired phase factor in the first term on the right-hand side of \eqref{equ:null_structure2++_derivation5},
    \begin{equation} \label{equ:null_structure2++_derivation6}
        \begin{aligned}
            &-3  \Bigl( \gamma^2 (\xi_1 + \ell \jxione)^2 - \gamma^2 (\xi_2+\ell\jxitwo)^2 \Bigr) \Bigl( \gamma (\xi_1 + \ell \jxione) - \gamma (\xi_2+\ell\jxitwo) \Bigr) \\ 
            &= -3 \gamma^2\bigl( \jxione + \ell \xi_1 + \jxitwo + \ell \xi_2 \bigr)  \bigl( \jxione + \ell \xi_1 - \jxitwo - \ell \xi_2 \bigr) \bigl( \gamma (\xi_1 + \ell \jxione) - \gamma (\xi_2+\ell\jxitwo) \bigr).
        \end{aligned}
    \end{equation}
    Hence, combining \eqref{equ:null_structure2++_derivation3}, \eqref{equ:null_structure2++_derivation4}, \eqref{equ:null_structure2++_derivation5}, and \eqref{equ:null_structure2++_derivation6}, we conclude that 
    \begin{equation*}
        \begin{aligned}
            \fraka_{\ell,++}(\xi_1,\xi_2) &= -3 \gamma^2 \bigl( \jxione + \ell \xi_1 + \jxitwo + \ell \xi_2 \bigr)  \bigl( \jxione + \ell \xi_1 - \jxitwo - \ell \xi_2 \bigr) \bigl( \gamma (\xi_1 + \ell \jxione) - \gamma (\xi_2+\ell\jxitwo) \bigr) \\ 
            &\quad + 3 \gamma \ell \bigl( \jxione + \ell \xi_1 + \jxitwo + \ell \xi_2 \bigr) \Bigl( \gamma (\xi_1 + \ell \jxione) - \gamma (\xi_2+\ell\jxitwo) \Bigr)^2 \\ 
            &\quad + \gamma \ell \bigl( \jxione + \ell \xi_1 - \jxitwo - \ell \xi_2 \bigr) \times \\ 
            &\quad \quad \quad \times \gamma^{-1} (\xi_1+\xi_2) \Bigl( 6 \bigl( \gamma (\xi_1 + \ell \jxione) - \gamma (\xi_2+\ell\jxitwo) \bigr) - 3 \gamma \ell \bigl( \jxione + \ell \xi_1 + \jxitwo + \ell \xi_2 \bigr) \Bigr) \\ 
            &\quad +\bigl( \jxione + \ell \xi_1 +\jxitwo +\ell \xi_2 \bigr) 4  \gamma \ell\bigl( -2 + \gamma^{-2}(\xi_1+\xi_2)^2 \bigr)  \\ 
            &= \bigl( \jxione + \ell \xi_1 + \jxitwo + \ell \xi_2 \bigr) \\ 
            &\quad \quad \times \biggl(-3\gamma^2 \bigl( \jxione + \ell \xi_1 - \jxitwo - \ell \xi_2 \bigr) \bigl( \gamma (\xi_1 + \ell \jxione) - \gamma (\xi_2+\ell\jxitwo) \bigr) \\ 
            &\qquad \qquad +3 \gamma \ell \Bigl( \gamma (\xi_1 + \ell \jxione) - \gamma (\xi_2+\ell\jxitwo) \Bigr)^2 \\ 
            &\qquad \qquad + \ell (\xi_1-\xi_2) \Bigl( 6 \bigl( \gamma (\xi_1 + \ell \jxione) - \gamma (\xi_2+\ell\jxitwo) \bigr) - 3 \gamma \ell \bigl( \jxione + \ell \xi_1 + \jxitwo + \ell \xi_2 \bigr) \Bigr) \\
            &\qquad \qquad + 4 \gamma \ell \bigl( -2 + \gamma^{-2}(\xi_1+\xi_2)^2 \bigr) \biggr).
        \end{aligned}
    \end{equation*}
    Inserting the preceding identity back into \eqref{equ:null_structure2++_derivation0} yields the asserted identity \eqref{equ:null_structure2++}.
\end{proof}

\section{Setting up the Analysis} \label{sec:setting_up}

In this section we begin in earnest with the proof of the asymptotic stability of the family of moving kink solutions \eqref{equ:intro_vectorial_moving_kink_family} to the sine-Gordon equation~\eqref{equ:intro_sG_1st_order}. In Subsection~\ref{subsec:local_existence_in_neighborhood_of_kink} we record a local existence result for the sine-Gordon equation that suffices for our purposes, and in Subsection~\ref{subsec:modulation_and_orbital_stability} we establish a decomposition of a solution to the sine-Gordon equation~\eqref{equ:intro_sG_1st_order} in a neighborhood of a kink into a modulated kink and a radiation term.    
In Subsections~\ref{subsec:evolution_equation_profile}--\ref{subsec:structure_cubic_nonlinearities} we then prepare the evolution equation for the effective profile on the distorted Fourier side for the nonlinear analysis in the remainder of this paper.

\subsection{Local existence in a neighborhood of a kink} \label{subsec:local_existence_in_neighborhood_of_kink}
In what follows we will work with $H^3_x \times H^2_x$ solutions to the sine-Gordon equation in a neighborhood of a kink. For given $\ell_0 \in (-1,1)$, $x_0 \in \bbR$, and $\bmu_0 = (u_{0,1}, u_{0,2}) \in H^3_x(\bbR) \times H^2_x(\bbR)$, we construct a (local-in-time) solution $\bmphi(t) = \bigl(\phi_1(t), \phi_2(t)\bigr)$ to \eqref{equ:intro_sG_1st_order} with initial data 
\begin{equation}
    \bmphi(0,x) = \bmK_{\ell_0,x_0}(x) + \bmu_0(x-x_0),
\end{equation}
where we recall from \eqref{equ:statement_theorem_initial_data} that 
\begin{equation}
    \bmK_{\ell_0, x_0}(x) := \begin{bmatrix}
        K\bigl(\gamma_0 (x-x_0)\bigr) \\ -\gamma_0 \ell_0 K'\bigl(\gamma_0(x-x_0)\bigr)
    \end{bmatrix}.
\end{equation}
To this end we introduce the variables
\begin{equation}
    \begin{aligned}
        w_1(t,x) &:= \phi_1(t,x) - K\bigl( \gamma_0 (x - \ell_0 t - x_0) \bigr), \\
        w_2(t,x) &:= \phi_2(t,x) - \gamma_0 \ell_0 K'\bigl( \gamma_0 (x - \ell_0 t - x_0) \bigr).
    \end{aligned}
\end{equation}
These have to solve the system
\begin{equation} \label{equ:setting_up_system_for_v}
    \left\{ \begin{aligned}
        \pt w_1 &= w_2, \\
        \pt w_2 &= \px^2 w_1 - W'\Bigl( K\bigl(\gamma_0(\cdot-\ell_0 t - x_0)\bigr) + w_1\Bigr) + W'\Bigl( K\bigl(\gamma_0(\cdot-\ell_0 t - x_0)\bigr) \Bigr),
    \end{aligned} \right.
\end{equation}
where we remind the reader that the self-interaction potential for the sine-Gordon model is 
\begin{equation}
    W(\phi) = 1 - \cos(\phi).
\end{equation}
By a standard fixed-point argument we obtain the following local existence result.

\begin{lemma} \label{lem:setting_up_local_existence}
    For any $\bmu_0 = (u_{0,1}, u_{0,2}) \in H^3_x(\bbR) \times H^2_x(\bbR)$ there exists a unique solution $\bm{w} = (w_1, w_2) \in \calC\bigl([0,T_\ast); H^3_x \times H^2_x\bigr)$ to \eqref{equ:setting_up_system_for_v} with initial data $\bm{w}(0,x) = \bmu_0(x-x_0)$ on a maximal interval of existence $[0,T_\ast)$. Moreover, it holds that
    \begin{equation} \label{equ:setting_up_local_existence_continuation_criterion} 
        T_\ast < \infty \quad \Rightarrow \quad \limsup_{t \, \nearrow \, T_\ast} \, \|\bm{w}(t)\|_{H^3_x \times H^2_x} = \infty.
    \end{equation}
\end{lemma}

\subsection{Modulation and orbital stability} \label{subsec:modulation_and_orbital_stability}

The goal of this subsection is to establish a decomposition of a solution to the sine-Gordon equation in a neighborhood of a kink into a modulated kink and a radiation term that is orthogonal to directions related to the invariances of the sine-Gordon equation under spatial translations and under Lorentz transformations.
To this end we largely proceed along the lines of the modulation and orbital stability arguments in Subsection~2.8, Subsection~2.9, and Section~3 in \cite{KMMV20}. 

We recall from \eqref{equ:statement_theorem_modulated_kink} the expression for the modulated kink
\begin{equation} \label{eq:modkink}
    \bmK_{\ell,q}(x) := \begin{bmatrix}
                        K\bigl( \gamma (x-q) \bigr) \\
                        - \gamma \ell K'\bigl( \gamma (x-q) \bigr)
                     \end{bmatrix}, 
                     \quad \ell \in (-1,1), \quad q \in \bbR, \quad \gamma := \frac{1}{\sqrt{1-\ell^2}},
\end{equation}
and we remind the reader of the matrix $\bfJ$ defined in \eqref{equ:bfJ_definition}.
In what follows, $\ell = \ell(t)$ and $q = q(t)$ will usually be time-dependent continuously differentiable paths, but the time-dependence will often not be made explicit in the notation. Observe that if $\dot{\ell} = 0$ and if $\dot{q}-\ell = 0$, then \eqref{eq:modkink} is a moving kink solution to \eqref{equ:intro_sG_1st_order}.
Occasionally, we will also denote the scalar modulated kink by
\begin{equation*}
    K_{\ell,q}(x) := K\bigl( \gamma (x-q) \bigr).
\end{equation*}

The main result of this subsection is summarized in the following proposition.

\begin{proposition} \label{prop:setting_up_modulation_orbital}
    Fix $\ell_0 \in (-1,1)$. There exist constants $C_1 \geq 1$ and $0 < \varepsilon_1 \ll 1$ with the following properties:
    For any $x_0 \in \bbR$ and any $\bmu_0 \in H^3_x(\bbR) \times H^2_x(\bbR)$ with $\|\bmu_0\|_{H^1_x \times L^2_x} \leq \varepsilon_1$, Lemma~\ref{lem:setting_up_local_existence} furnishes a unique $H^3_x \times H^2_x$ solution to \eqref{equ:intro_sG_1st_order} with initial data $\bmphi(0,x) = \bmK_{\ell_0,x_0}(x) + \bmu_0(x-x_0)$ defined on a maximal interval of existence $[0, T_\ast)$ for some $T_\ast > 0$.
    Then there exist unique continuously differentiable paths $\ell \colon [0,T_\ast) \to (-1,1)$ and $q \colon [0,T_\ast) \to \bbR$ with 
    \begin{equation*}
     |\ell(0) - \ell_0| + |q(0)-x_0| \leq C_1 \|\bmu_0\|_{H^1_x \times L^2_x}
    \end{equation*}
    such that the following holds:
    \begin{itemize}[leftmargin=1.8em]
        \item[(1)] Decomposition of the solution into a modulated kink and a radiation term:
        \begin{equation}
            \bmphi(t,x) = \bmK_{\ell(t),q(t)}(x) + \bmu\bigl(t,x-q(t)\bigr), \quad 0 \leq t < T_\ast.
        \end{equation}        

        \item[(2)] Orthogonality:
        \begin{equation} \label{equ:orthogonality_radiation}
            \bigl\langle \bfJ \partial_q \bm{K}_{\ell,q}, \bmu(t,\cdot-q) \bigr\rangle = \bigl\langle \bfJ \partial_\ell \bm{K}_{\ell,q}, \bmu(t,\cdot-q) \bigr\rangle = 0, \quad 0 \leq t < T_\ast.
        \end{equation}        

        \item[(3)] Evolution equation for $\bmu(t,x-q(t))$:
            \begin{equation} \label{equ:setting_up_perturbation_equ}
                \partial_t \bmu 
                = \bfL_{\ell,q} \bmu + (\dot{q} - \ell) \px \bmu - (\dot q-\ell) \partial_q \bm{K}_{\ell,q} - \dot \ell \partial_\ell \bm{K}_{\ell,q} + \calN(\bmu)
            \end{equation} 
            with 
            \begin{equation} \label{equ:setting_up_L_ell_q_definition}
                \bfL_{\ell,q} =
                \begin{bmatrix}
	               \ell \px & 1 \\
	                 - L_{\ell,q} & \ell \px
	              \end{bmatrix}, 
                \quad
                L_{\ell,q} := - \px^2 + W''(K_{\ell,q}) = - \px^2 - 2\sech^2\bigl(\gamma(x-q)\bigr) + 1,
            \end{equation}
            and 
            \begin{equation} \label{equ:setting_up_definition_calN}
                \calN(\bmu) =   \begin{bmatrix}
                                 0 \\ - W'\bigl(K_{\ell,q}+u_1\bigr) + W'\bigl( K_{\ell,q} \bigr) + W''(K_{\ell,q}) u_1
                                \end{bmatrix}.
            \end{equation}

        \item[(4)] Modulation equations:
        \begin{equation} \label{equ:modulation_equ}
            \bbM_{\ell}[\bmu] \begin{bmatrix} \dot{\ell} \\ \dot{q}-\ell \end{bmatrix} =
        \begin{bmatrix}
            - \bigl\langle \bfJ \partial_q \bm{K}_{\ell,q}, \calN(\bm{u}) \bigr\rangle \\
            \bigl\langle \bfJ \partial_\ell \bm{K}_{\ell,q}, \calN(\bm{u}) \bigr\rangle 
        \end{bmatrix}
        \end{equation}
        with
        \begin{equation} \label{equ:modulation_equ_matrix_lhs}
        \bbM_{\ell}[\bmu] := 
        \begin{bmatrix} 
        \gamma^3 \| K' \|_{L^2}^2 + \bigl\langle \bfJ \partial_\ell \partial_q \bm{K}_{\ell,q}, \bm{u} \bigr\rangle &  \bigl\langle \bfJ \partial_q \partial_q \bm{K}_{\ell,q}, \bm{u} \bigr\rangle \\
        \bigl\langle \bfJ \partial_\ell \partial_\ell \bm{K}_{\ell,q}, \bm{u} \bigr\rangle & \gamma^3 \| K' \|_{L^2}^2 - \bigl\langle \bfJ \partial_q \partial_\ell \bm{K}_{\ell,q}, \bm{u} \bigr\rangle
        \end{bmatrix}.
        \end{equation}        

        \item[(5)] Stability:
        \begin{equation} \label{equ:smallness_orbital_proposition}
            \sup_{0 \leq t < T_\ast} \, \Bigl( \|\bmu(t)\|_{H^1_x \times L^2_x} + |\ell(t) - \ell_0| \Bigr) \leq C_1 \|\bmu_0\|_{H^1_x \times L^2_x}.
        \end{equation}

        \item[(6)] Comparison estimate: 
        \begin{equation} \label{equ:comparison_estimate}
            |\ell(t) - \ell_0| \leq \frac14 \gamma_0^{-2}, \quad \frac12 \gamma_0 \leq \gamma(t) \leq 2 \gamma_0, \quad 0 \leq t < T_\ast.
        \end{equation}
    \end{itemize}
\end{proposition}

\begin{remark}
    We emphasize that in \eqref{equ:orthogonality_radiation}, \eqref{equ:setting_up_perturbation_equ}, \eqref{equ:setting_up_L_ell_q_definition}, \eqref{equ:setting_up_definition_calN}, \eqref{equ:modulation_equ}, \eqref{equ:modulation_equ_matrix_lhs},  and in many parts of all subsequent sections, the time dependence of the modulation parameters $\ell = \ell(t)$ and $q=q(t)$ will not be made explicit in the notation.
\end{remark}

\begin{remark}
    Since we are working with $H^3_x(\bbR)\times H^2_x(\bbR)$ solutions to the sine-Gordon equation, one could actually show that the paths $\ell(t)$ and $q(t)$ are $\calC^3$, but we do not need to use this fact for the proof of Theorem~\ref{thm:main}, so we omit details.
\end{remark}

We begin with several preparations for the proof of Proposition~\ref{prop:setting_up_modulation_orbital}.
First, we note that the differential equation $-K'' + W'(K) = 0$ for the scalar kink implies that the modulated kink \eqref{eq:modkink} solves the system
\begin{equation} \label{eq:vecKeq}
 \begin{bmatrix} 
  0 & 1 \\ \px^2 & 0 
 \end{bmatrix} \bmK_{\ell, q}(x) 
 + 
 \begin{bmatrix}
     0 \\ - W'(K_{\ell,q}) 
 \end{bmatrix} 
 + \ell \px \bmK_{\ell, q}(x) = 0.  
\end{equation}
Viewing $\ell$ and $q$ as time-independent in \eqref{eq:vecKeq}, differentiating \eqref{eq:vecKeq} with respect to $q$ gives
\begin{equation}
   \bfL_{\ell,q}\left(\partial_q \bmK_{\ell, q}\right)=0,
\end{equation}
while differentiating \eqref{eq:vecKeq} with respect to $\ell$ yields
\begin{equation}
   \bfL_{\ell,q}\left(\partial_\ell \bmK_{\ell, q}\right) = \partial_q \bmK_{\ell, q}, \quad \text{whence} 
   \quad \bfL_{\ell,q}^2\left(\partial_\ell \bmK_{\ell, q}\right) = 0.
\end{equation}
By direct computation we infer from \eqref{eq:modkink} that 
\begin{equation} 
 \begin{aligned}
\partial_q \bmK_{\ell, q}(x) = \begin{bmatrix} -\gamma K'(\gamma( x-q)) \\ \gamma^2 \ell K''(\gamma (x-q)) \end{bmatrix}, \quad
 \partial_\ell \bmK_{\ell, q}(x) = \begin{bmatrix} \gamma^3 \ell (x-q) K'(\gamma (x-q)) \\ - \gamma^3 K'(\gamma (x-q)) - \gamma^4 \ell^2 (x-q) K''(\gamma (x-q)) \end{bmatrix}.
 \end{aligned}
\end{equation}
This computation gives the expressions for $\bmY_0$ and $\bmY_1$ introduced in \eqref{equ:Y0Y1_definition}.
By Lemma \ref{lem:Lspec} the generalized null space of $\bfL_{\ell,q}$ is spanned by $\partial_q \bmK_{\ell, q}$ and $\partial_\ell \bmK_{\ell, q}(x)$. 
The motivation for the choice of orthogonality \eqref{equ:orthogonality_radiation} for the radiation term is that then the projection of the radiation term $\bmu(t)$ to the discrete spectrum of the linearized operator $\bfL_{\ell(t), q(t)}$ is zero by Remark~\ref{rem:Pdzero}.

Moreover, we record that
\begin{equation} \label{equ:Jinnerproduct_generalized_kernel_elements}
    \langle \bfJ \partial_q \bmK_{\ell, q},\partial_\ell \bmK_{\ell, q}\rangle = -\gamma^4 \int_\bbR (K'(\gamma(x-q)))^2\,\ud x = -\gamma^3\|K'\|^2_{L^2}.
\end{equation}

Next, we establish a few comparison estimates.

\begin{lemma}[Basic comparison estimates] \label{lem:setting_up_basic_comparison}
    Let $\ell_0 \in (-1,1)$ and set $\gamma_0 := (1-\ell_0^2)^{-\frac12}$.
    \begin{itemize}[leftmargin=1.8em]
        \item[(a)] If $|\ell_1 - \ell_0| \leq \frac14 \gamma_0^{-2}$, then 
                    \begin{equation*}
                        1-|\ell_1| \geq \frac12 \bigl( 1-|\ell_0| \bigr),
                    \end{equation*}
                    whence $\ell_1 \in (-1,1)$. Moreover, $\gamma_1 := (1-\ell_1^2)^{-\frac12}$ satisfies
                    \begin{equation} \label{equ:setting_up_comparison_gammas}
                        \frac12 \gamma_0 \leq \gamma_1 \leq 2 \gamma_0, \quad |\gamma_1-\gamma_0| \leq 2 \gamma_0^3 |\ell_1-\ell_0|.
                    \end{equation}
        \item[(b)] For $\ell_0, \ell_1 \in (-1,1)$ it holds that
                    \begin{equation} \label{equ:setting_up_comparison_ellgamma_difference}
                        |\ell_1 - \ell_0| \leq 3 \gamma_0 \bigl| \ell_1 \gamma_1 - \ell_0 \gamma_0 \bigr|.
                    \end{equation}
        \item[(c)] There exists an absolute constant $C \geq 1$ such that for any $\ell, \ell_0 \in (-1,1)$ with $|\ell-\ell_0| \leq \frac14 \gamma_0^{-2}$ and for any $q, x_0 \in \bbR$,
                    \begin{equation} \label{equ:setting_up_comparison_kink_difference}
                        \bigl\| \bmK_{\ell,q} - \bmK_{\ell_0,x_0} \bigr\|_{H^1_x \times L^2_x} \leq C \Bigl( \gamma_0^{\frac52} |\ell-\ell_0| + \gamma_0^{\frac32} |q-x_0| \Bigr).
                    \end{equation}
    \end{itemize}
\end{lemma}
\begin{proof}
    Part (a) follows by direct computation. To prove part (b) we use the short-hand notation $d_j := \ell_j \gamma_j$. Then we have $\ell_j = d_j (1+d_j^2)^{-\frac12}$ and we compute that
    \begin{equation*}
        \begin{aligned}
            \ell_1 - \ell_0 
            &= \bigl( d_1 - d_0 \bigr) \biggl( \frac{1}{\sqrt{1+d_1^2}} - d_0 \frac{d_0 + d_1}{\sqrt{1+d_0^2} \sqrt{1+d_1^2} \bigl( \sqrt{1+d_0^2} + \sqrt{1+d_1^2} \bigr)} \biggr),            
        \end{aligned}
    \end{equation*}
    whence 
    \begin{equation*}
        |\ell_1 - \ell_0| \leq \bigl( 1 + 2 |d_0| \bigr) |d_1 - d_0| \leq 3 \gamma_0 |d_1 - d_0|,
    \end{equation*}    
    as claimed.
    Part (c) again follows by direct computation, for which we omit the details.
\end{proof}

While the completely integrable sine-Gordon equation admits infinitely many conservation laws, relevant to our analysis in this subsection is that solutions $\bmphi = (\phi_1, \phi_2)$ to \eqref{equ:intro_sG_1st_order} formally conserve the energy 
\begin{equation} \label{equ:def_energy_system}
    E[\bmphi] = \int_\bbR \biggl( \frac12 (\px \phi_1)^2 + \frac12 \phi_2^2 + W(\phi_1) \biggr) \, \ud x
\end{equation}
and the momentum 
\begin{equation} \label{equ:def_momentum_system}
    P[\bmphi] = \int_\bbR \phi_2 (\px \phi_1) \, \ud x.
\end{equation}

We also introduce the following conserved quantity, which corresponds to the magnitude of the energy-momentum $2$-vector and which is invariant under Lorentz transformations:
\begin{equation*}
    M[\bmphi] := \bigl( E[\bmphi] \bigr)^2 - \bigl( P[\bmphi] \bigr)^2.
\end{equation*}
Next, we study expansions of the conservation laws for \eqref{equ:intro_sG_1st_order} around the modulated kinks.
\begin{lemma} \label{lem:setting_up_expansion_conservation_laws} 
    Let $\ell \in (-1,1)$, $q \in \bbR$, and $\bmu = (u_1, u_2) \in H^1_x(\bbR) \times L^2_x(\bbR)$.
    \begin{itemize}[leftmargin=1.8em]
        \item[(a)] Conservation laws for modulated kinks:
        \begin{equation} \label{equ:setting_up_conservation_laws_kinks}
            E\bigl[ \bmK_{\ell,q} \bigr] = \gamma \|K'\|_{L^2_x}^2, \quad P\bigl[ \bmK_{\ell, q} \bigr] = - \gamma \ell \|K'\|_{L^2_x}^2, \quad M\bigl[\bmK_{\ell,q}\bigr] = \|K'\|_{L^2_x}^4. 
        \end{equation}

        \item[(b)] Expansion of energy around a modulated kink:
            \begin{equation} \label{equ:setting_uP_expansion_E_around_kink}
                \begin{aligned}
                    E\bigl[ \bmK_{\ell,q} + \bmu \bigr] - E\bigl[ \bmK_{\ell,q} \bigr]  
                    = \frac12 \langle L_{\ell,q} u_1, u_1 \rangle + \frac12 \|u_2\|_{L^2_x}^2 - \ell \bigl\langle \bfJ \partial_q \bmK_{\ell,q}, \bmu \bigr\rangle + \calR_1
                \end{aligned}
            \end{equation}
            with    
            \begin{equation*}
                L_{\ell,q} := -\px^2 + W''\bigl(K_{\ell,q}(x)\bigr)
            \end{equation*}
            and
            \begin{equation*}
                \calR_1 := \int_\bbR \biggl( W(K_{\ell,q} + u_1) - W(K_{\ell,q}) - W'(K_{\ell,q}) u_1 - \frac12 W''(K_{\ell,q}) u_1^2 \biggr) \, \ud x.
            \end{equation*}

        \item[(c)] Expansion of momentum around a modulated kink:
            \begin{equation} \label{equ:setting_up_expansion_P_around_kink}
                P\bigl[ \bmK_{\ell,q} + \bmu \bigr] - P\bigl[ \bmK_{\ell,q} \bigr] = \bigl\langle \bfJ \partial_q \bmK_{\ell,q}, \bmu \bigr\rangle + P[ \bmu ].
            \end{equation}        

        \item[(d)] Expansion of $M\bigl[\bmphi\bigr]$ around a modulated kink:
            If $\langle \bfJ \partial_q \bmK_{\ell,q}, \bmu \rangle = 0$, then we have 
            \begin{equation} \label{equ:setting_uP_expansion_M_around_kink}
                M\bigl[ \bmK_{\ell,q} + \bmu \bigr] - M\bigl[ \bmK_{\ell,q} \bigr] = \gamma \|K'\|_{L^2_x}^2 \bigl\langle \bfH_{\ell,q} \bmu, \bmu \rangle + \calR_2,
            \end{equation}
            where 
            \begin{equation*}
                \bfH_{\ell, q} :=   \begin{bmatrix}
                                L_{\ell,q} & -\ell \px \\ 
                                \ell \px & 1 
                            \end{bmatrix},
            \end{equation*}
            and
            \begin{equation*}
                \calR_2 := 2 \gamma \|K'\|_{L^2_x} \calR_1 + \frac14 \Bigl( \bigl\langle \bfH_{\ell,q} \bmu, \bmu \rangle - 2\ell P[\bmu] + 2 \calR_1 \Bigr)^2 - \bigl( P[\bmu] \bigr)^2.
            \end{equation*}
    \end{itemize}
\end{lemma}
\begin{proof}
    The asserted identities all follow by direct computation.
\end{proof}

Finally, we are in the position to prove Proposition~\ref{prop:setting_up_modulation_orbital}.

\begin{proof}[Proof of Proposition~\ref{prop:setting_up_modulation_orbital}]
Let $\ell_0 \in (-1,1)$, $x_0 \in \bbR$, and $\bmu_0 \in H^3_x(\bbR) \times H^2_x(\bbR)$ be given. 
Lemma~\ref{lem:setting_up_local_existence} furnishes a unique $H^3_x \times H^2_x$ solution $\bmphi(t)$ to \eqref{equ:intro_sG_1st_order} with initial data $\bmphi(0,x) = \bmK_{\ell_0,x_0}(x) + \bmu_0(x-x_0)$ defined on a maximal interval of existence $[0,T_\ast)$.

{\bf Step 1.} (Choice of parameters at initial time)
We show that there exist constants $\kappa_1 \equiv \kappa_1(\ell_0)$ with $0 < \kappa_1 \ll 1$ and $C_0 \equiv C_0(\ell_0) \geq 1$ such that if $\kappa := \|\bmu_0\|_{H^1_x \times L^2_x} \leq \kappa_1$, then there exist $\ell(0) \in (-1,1)$ and $q(0) \in \bbR$ with the properties that 
\begin{equation} \label{equ:setting_up_modulation_proof_smallness_from_IFT}
    |\ell(0)-\ell_0| + |q(0)-x_0| \leq C_0 \kappa,
\end{equation}
and 
\begin{equation} \label{equ:setting_up_modulation_proof_orthogonality_from_IFT}
    \bigl\langle \bfJ \partial_q \bmK_{\ell(0), q(0)}, \bmphi(0) - \bmK_{\ell(0), q(0)} \bigr\rangle = \bigl\langle \bfJ \partial_\ell \bmK_{\ell(0), q(0)}, \bmphi(0) - \bmK_{\ell(0), q(0)} \bigr\rangle = 0.
\end{equation}
To this end we consider the smooth functional
\begin{equation*}
    \bmG\bigl( \ell, q, \bmpsi \bigr) :=    \begin{bmatrix}
                                                 \bigl\langle \bfJ \partial_q \bmK_{\ell, q}, \bmpsi - \bmK_{\ell, q} \bigr\rangle \\ 
                                                 \bigl\langle \bfJ \partial_\ell \bmK_{\ell, q}, \bmpsi - \bmK_{\ell, q} \bigr\rangle
                                            \end{bmatrix}
\end{equation*}
defined on the domain 
\begin{equation*}
    \bigl( \ell, q, \bmpsi \bigr) \in (-1,1) \times \bbR \times \bigl( (K,0) + H^1_x \times L^2_x \bigr).
\end{equation*}
Clearly, we have $\bmG(\ell_0, x_0, \bmK_{\ell_0, x_0}) = 0$, and by direct computation we obtain that the Jacobian matrix of $\bfG$ with respect to the modulation parameters is invertible at $(\ell_0, x_0, \bmK_{\ell_0,x_0})$,
\begin{equation*}
    \begin{aligned}
     \frac{\partial \bmG}{\partial (\ell, q)} \bigg|_{(\ell,q, \bmpsi) = (\ell_0, x_0, \bmK_{\ell_0, x_0})} &= \begin{bmatrix} 
                                                            - \bigl\langle \bfJ \partial_q \bmK_{\ell, q}, \partial_\ell \bmK_{\ell, q} \bigr\rangle & 0 \\
                                                            0 & - \bigl\langle \bfJ \partial_\ell \bmK_{\ell, q}, \partial_q \bmK_{\ell, q} \bigr\rangle 
                                                        \end{bmatrix} \bigg|_{(\ell, q) = (\ell_0, x_0)} \\
                                                        &=
                                                        \begin{bmatrix}
                                                            -\gamma_0^3 \|K'\|_{L^2_x}^2 & 0 \\ 0 & \gamma_0^3 \|K'\|_{L^2_x}^2
                                                        \end{bmatrix}.
    \end{aligned}
\end{equation*}     
In fact, an analogous computation shows that the Jacobian matrix of $\bmG$ with respect to the modulation parameters is uniformly non-degenerate in a neighborhood of $(\ell_0, x_0, \bmK_{\ell_0, x_0})$.
The existence of $0 < \kappa_1 \ll 1$ and $C_0 \geq 1$ with the asserted properties now follows from a quantitative version of the implicit function theorem, see for instance \cite[Remark 3.2]{Jendrej15}.

{\bf Step 2.} (Derivation of modulation equations)
Next, we carry out formal computations. Given continuously differentiable paths $\ell(t)$, $q(t)$ defined for times $0 \leq t < T_\ast$ with the initial values $\ell(0)$, $q(0)$ obtained in the preceding step, we set
\begin{equation*}
    \bmu\bigl(t,x-q(t)\bigr) := \bmphi(t,x) - \bmK_{\ell(t), q(t)}(x).
\end{equation*}
Using the equation \eqref{equ:intro_sG_1st_order} for $\bmphi(t,x)$ and the equation \eqref{eq:vecKeq} for $\bmK_{\ell(t), q(t)}(x)$ we obtain by direction computation the evolution equation \eqref{equ:setting_up_perturbation_equ} for the radiation term $\bmu$.

Differentiating the orthogonality conditions
\begin{equation} \label{equ:setting_up_modulation_proof_orthogonality}
    \bigl\langle \bfJ \partial_q \bmK_{\ell(t), q(t)}, \bmu(t,\cdot-q(t)) \bigr\rangle = \bigl\langle \bfJ \partial_\ell \bmK_{\ell(t), q(t)}, \bmu(t,\cdot-q(t)) \bigr\rangle = 0,
\end{equation}
we arrive at the modulation equations \eqref{equ:modulation_equ}.
Indeed, differentiating the first equation in \eqref{equ:setting_up_modulation_proof_orthogonality} and inserting \eqref{equ:setting_up_perturbation_equ}, we get
\begin{equation}
 \begin{aligned}
   0 &= \frac{\ud}{\ud t} \, \bigl\langle \bfJ \partial_q \bmK_{\ell(t), q(t)}, \bmu(t,\cdot-q(t)) \bigr\rangle \\
   &= \dot{\ell} \bigl\langle \bfJ \partial_\ell \partial_q \bmK_{\ell(t), q(t)}, \bmu(t,\cdot-q(t)) \bigr\rangle + \bigl\langle \bfJ \partial_q \bmK_{\ell(t), q(t)}, (\partial_t \bmu)(t,\cdot-q(t)) \bigr\rangle \\ 
   &= \dot{\ell} \bigl\langle \bfJ \partial_\ell \partial_q \bmK_{\ell(t), q(t)}, \bmu(t,\cdot-q(t)) \bigr\rangle 
   + \bigl\langle \bfJ \partial_q \bmK_{\ell(t), q(t)}, \bfL_{\ell,q} \bmu (t,\cdot-q(t)) \bigr\rangle \\
   &\quad + (\dot{q} - \ell) \bigl\langle \bfJ \partial_q \bmK_{\ell(t), q(t)},  \px \bmu (t,\cdot-q(t))  \bigr\rangle 
   - (\dot q-\ell)\bigl\langle \bfJ \partial_q \bmK_{\ell(t), q(t)},    \partial_q \bmK_{\ell(t), q(t)}  \bigr\rangle \\
   &\quad - \dot \ell\bigl\langle \bfJ \partial_q \bmK_{\ell(t), q(t)},  \partial_\ell \bmK_{\ell(t), q(t)}\bigr\rangle + \bigl\langle \bfJ \partial_q \bmK_{\ell(t), q(t)}, \calN(\bmu) \bigr\rangle.
 \end{aligned}
\end{equation}
Next, we observe that $\bigl\langle \bfJ \partial_q \bmK_{\ell(t), q(t)}, \partial_q \bmK_{\ell(t), q(t)}  \bigr\rangle=0$, and that by Lemma~\ref{lem:L*spec},
\begin{equation}
 \bigl\langle \bfJ \partial_q \bmK_{\ell(t), q(t)}, \bfL_{\ell,q} \bmu (t,\cdot-q(t)) \bigr\rangle=\bigl\langle \bfL_{\ell,q}^*\bfJ \partial_q \bmK_{\ell(t), q(t)},  \bmu (t,\cdot-q(t))\bigr\rangle=0.
\end{equation}
Moreover, integration by parts gives
$$ \bigl\langle \bfJ \partial_q \bmK_{\ell(t), q(t)},  \px \bmu  (t,\cdot-q(t))\bigr\rangle= \bigl\langle \bfJ \partial_q \partial_q \bmK_{\ell(t), q(t)},   \bmu (t,\cdot-q(t)) \bigr\rangle. $$
Putting the preceding computations together and using \eqref{equ:Jinnerproduct_generalized_kernel_elements}, the first row of \eqref{equ:modulation_equ} follows.

The derivation of the second row of \eqref{equ:modulation_equ} is quite similar. 
Indeed, differentiating the second equation in \eqref{equ:setting_up_modulation_proof_orthogonality} and inserting \eqref{equ:setting_up_perturbation_equ} again, we find
\begin{equation}
\begin{aligned}
   0 &= \frac{\ud}{\ud t} \, \bigl\langle \bfJ \partial_\ell \bmK_{\ell(t), q(t)}, \bmu (t,\cdot-q(t) \bigr\rangle \\
   &= \dot{\ell} \bigl\langle \bfJ \partial_\ell \partial_\ell \bmK_{\ell(t), q(t)}, \bmu (t,\cdot-q(t) \bigr\rangle  + \bigl\langle \bfJ \partial_\ell \bmK_{\ell(t), q(t)}, (\partial_t \bmu)(t,\cdot-q(t) \bigr\rangle \\ 
   &= \dot{\ell} \bigl\langle \bfJ \partial_\ell \partial_\ell \bmK_{\ell(t), q(t)}, \bmu (t,\cdot-q(t) \bigr\rangle + \bigl\langle \bfJ \partial_\ell \bmK_{\ell(t), q(t)}, \bfL_{\ell,q} \bmu  (t,\cdot-q(t)\bigr\rangle \\
   &\quad +(\dot{q} - \ell) \bigl\langle \bfJ \partial_\ell \bmK_{\ell(t), q(t)},  \px \bmu  (t,\cdot-q(t) \bigr\rangle - (\dot q-\ell)\bigl\langle \bfJ \partial_\ell \bmK_{\ell(t), q(t)},    \partial_q \bmK_{\ell(t), q(t)} \bigr\rangle \\
   &\quad - \dot{\ell} \bigl\langle \bfJ \partial_\ell \bmK_{\ell(t), q(t)},  \partial_\ell \bmK_{\ell(t), q(t)}\bigr\rangle + \bigl\langle \bfJ \partial_\ell \bmK_{\ell(t), q(t)}, \calN(\bmu) \bigr\rangle.
\end{aligned}
\end{equation}
In view of Lemma~\ref{lem:L*spec} and the first equation in \eqref{equ:setting_up_modulation_proof_orthogonality}, we have 
\begin{equation}
\begin{aligned}
    \bigl\langle \bfJ \partial_\ell \bmK_{\ell(t), q(t)}, \bfL_{\ell,q} \bmu (t,\cdot-q(t) \bigr\rangle &= \bigl\langle \bfL_{\ell,q}^*\bfJ \partial_\ell \bmK_{\ell(t), q(t)},  \bmu (t,\cdot-q(t)\bigr\rangle \\
    &= \bigl\langle \bfJ \partial_q \bmK_{\ell(t), q(t)},  \bmu (t,\cdot-q(t)\bigr\rangle = 0.
\end{aligned}
\end{equation}
Moreover, we observe that $\bigl\langle \bfJ \partial_\ell \bmK_{\ell(t), q(t)}, \partial_\ell \bmK_{\ell(t), q(t)}  \bigr\rangle=0$, and that integration by parts gives
\begin{equation}
    \bigl\langle \bfJ \partial_\ell \bmK_{\ell(t), q(t)},  \px \bmu (t,\cdot-q(t) \bigr\rangle= \bigl\langle \bfJ \partial_q \partial_\ell \bmK_{\ell(t), q(t)},   \bmu (t,\cdot-q(t)  \bigr\rangle. 
\end{equation}
Combining the preceding computations then gives the second row of \eqref{equ:modulation_equ}.

{\bf Step 3.} (Solving the first-order ODEs)
Clearly, if $\|\bmu(t)\|_{H^1_x\times L^2_x} + |\ell(t)-\ell(0)| \ll 1$ is sufficiently small, the matrix $\bbM_{\ell}[\bmu]$ in \eqref{equ:modulation_equ} is invertible and we obtain the following non-autonomous system of differential equations for the modulation parameters $\ell(t)$ and $q(t)$,
\begin{equation} \label{equ:modulation_equ_inverted_M}
    \begin{bmatrix} \dot{\ell} \\ \dot{q} \end{bmatrix} = \begin{bmatrix} 0 \\ \ell \end{bmatrix} + \bbM_{\ell}[\bmu]^{-1} 
    \begin{bmatrix}
        - \bigl\langle \bfJ \partial_q \bm{K}_{\ell,q}, \calN(\bm{u}) \bigr\rangle \\
        \bigl\langle \bfJ \partial_\ell \bm{K}_{\ell,q}, \calN(\bm{u}) \bigr\rangle 
    \end{bmatrix}.
\end{equation}
Imposing the orthogonality conditions \eqref{equ:setting_up_modulation_proof_orthogonality} for times $0 \leq t < T_\ast$ is therefore equivalent to solving the system \eqref{equ:modulation_equ_inverted_M}.
The right-hand side of \eqref{equ:modulation_equ_inverted_M} depends Lipschitz-continuously on $\ell$ and $q$, and it depends continuously on time~$t$.
Thus, the Picard-Lindel\"of theorem furnishes a local solution to \eqref{equ:modulation_equ_inverted_M}. 
Moreover, there exists a small constant $0 < \kappa_2 \ll \kappa_1 \ll 1$ with $\kappa_2 \equiv \kappa_2(\ell_0)$ with the property that
as long as
\begin{equation*}
    \|\bmu(t)\|_{H^1_x \times L^2_x} + |\ell(t) - \ell(0)| \leq \kappa_2,
\end{equation*}
the matrix $\bbM_{\ell}[\bmu]$ is invertible and the solution to \eqref{equ:modulation_equ_inverted_M} may be continued. 

{\bf Step 4.} (Orbital stability argument)
In order to conclude that the solution to the system of differential equations \eqref{equ:modulation_equ_inverted_M} exists on the entire time interval $[0,T_\ast)$, we consider the exit time
\begin{equation*}
    T_1 := \sup \, \biggl\{ 0 < T < T_\ast \, \bigg| \, \sup_{0 \leq t \leq T} \, \Bigl( \|\bmu(t)\|_{H^1_x\times L^2_x} + |\ell(t) - \ell(0)| \Bigr) \leq \kappa_2 \biggr\}.
\end{equation*}
We now show that by choosing $0 < \kappa \ll \kappa_2$ sufficiently small, we must then have $T_1 = T_\ast$.
We may assume that $|\ell(t)-\ell_0| \leq\kappa_2 \leq \frac14 \gamma_0^{-2}$ for all $0 \leq t < T_1$, so that by \eqref{equ:setting_up_comparison_gammas} it holds that $\frac12 \gamma_0 \leq \gamma(t) \leq 2 \gamma_0$ for all $0 \leq t < T_1$. 

By the local theory, we must certainly have $T_1 > 0$. 
Suppose $T_1 < T_\ast$. 
Using the orthogonality \eqref{equ:setting_up_modulation_proof_orthogonality}, the coercivity bound \eqref{eq:H_lower_bound}, the upper bound \eqref{eq:H_upper_bound}, the fact that the remainder term $\calR_2$ in \eqref{equ:setting_uP_expansion_M_around_kink} satisfies $|\calR_2| \leq C(\ell_0) \|\bmu\|_{H^1_x \times L^2_x}^3$ for $0 < \|\bmu(t)\|_{H^1_x \times L^2_x} \leq 1$ and some $C(\ell_0) \geq 1$, and the conservation of energy and momentum, we find that for any $0 \leq t < T_1$,
\begin{equation*}
    \begin{aligned}
        \frac14 \mu \gamma_0^{-2} \|\bmu(t)\|_{H^1_x \times L^2_x}^2 - C(\ell_0) \|\bmu(t)\|_{H^1_x \times L^2_x}^3 &\leq M\bigl[ \bmK_{\ell(t), q(t)} + \bmu(t) \bigr] - M\bigl[ \bmK_{\ell(t),q(t)} \bigr] \\ 
        &= M\bigl[ \bmK_{\ell(0),q(0)} + \bmu(0) \bigr] - M\bigl[\bmK_{\ell(0),q(0)}\bigr] \\
        &\leq C \gamma_0 \|\bmu(0)\|_{H^1_x \times L^2_x}^2 + C(\ell_0) \|\bmu(0)\|_{H^1_x \times L^2_x}^3.
    \end{aligned}
\end{equation*}
Since by \eqref{equ:setting_up_comparison_kink_difference} and \eqref{equ:setting_up_modulation_proof_smallness_from_IFT},
\begin{equation*}
    \begin{aligned}
        \|\bmu(0)\|_{H^1_x \times L^2_x} &\leq \|\bmu_0\|_{H^1_x \times L^2_x} + \bigl\| \bmK_{\ell(0),q(0)} - \bmK_{\ell_0, x_0} \bigr\|_{H^1_x \times L^2_x} \\ 
        &\leq \|\bmu_0\|_{H^1_x \times L^2_x} + C\Bigl( \gamma_0^{\frac52} |\ell(0)-\ell_0| + \gamma_0^{\frac32} |q(0)-x_0| \Bigr) \\
        &\leq C(\ell_0) \kappa, 
    \end{aligned}
\end{equation*}
it follows by a continuity argument that 
\begin{equation} \label{equ:setting_up_modulation_proof_continuity_argument}
    \sup_{0 \leq t < T_1} \, \|\bmu(t)\|_{H^1_x \times L^2_x} \leq \widetilde{C}(\ell_0) \kappa
\end{equation}
for some constant $\widetilde{C}(\ell_0) \geq 1$.

Moreover, by \eqref{equ:setting_up_conservation_laws_kinks} and the conservation of momentum, we have for all $0 \leq t < T_1$ that
\begin{equation*}
    \begin{aligned}
        \bigl( \ell(t) \gamma(t) - \ell(0) \gamma(0) \bigr) \|K'\|_{L^2_x}^2 &= -P\bigl[\bmK_{\ell(t),q(t)}\bigr] + P\bigl[\bmK_{\ell(0),q(0)}\bigr] \\ 
        &= P\bigl[ \bmK_{\ell(t),q(t)} + \bmu(t) \bigr] - P\bigl[ \bmK_{\ell(t),q(t)} \bigr] \\ 
        &\quad - \Bigl( P\bigl[\bmK_{\ell(0),q(0)} + \bmu(0) \bigr] - P\bigl[ \bmK_{\ell(0),q(0)} \bigr] \Bigr).
    \end{aligned}
\end{equation*}
The expansion \eqref{equ:setting_up_expansion_P_around_kink} and the orthogonality \eqref{equ:setting_up_modulation_proof_orthogonality} together with \eqref{equ:setting_up_comparison_ellgamma_difference} now imply
\begin{equation*}
    \begin{aligned}
        \sup_{0 \leq t < T_1} \, |\ell(t) - \ell(0)| \lesssim \sup_{0 \leq t < T_1} \, \gamma_0 \bigl| \ell(t) \gamma(t) - \ell(0) \gamma(0) \bigr| &\lesssim \gamma_0 \sup_{0 \leq t < T_1} \, P[\bmu(t)] \\ 
        &\lesssim \gamma_0 \sup_{0 \leq t < T_1} \, \|\bmu(t)\|_{H^1_x \times L^2_x}^2.
    \end{aligned}
\end{equation*}
Hence, by \eqref{equ:setting_up_modulation_proof_continuity_argument} there exist constants $0 < \varepsilon_1 \ll 1$ with $\varepsilon_1 \equiv \varepsilon(\ell_0)$ and $C_1 \equiv C_1(\ell_0) \geq 1$ such that for all $0 < \kappa \leq \varepsilon_1$, we obtain that
\begin{equation*}
    \sup_{0 \leq t < T_1} \, \Bigl( \|\bmu(t)\|_{H^1_x \times L^2_x} + |\ell(t) - \ell(0)| \Bigr) \leq C_1\kappa \leq \kappa_2,
\end{equation*}
as desired. This finishes the proof of Proposition~\ref{prop:setting_up_modulation_orbital}.
\end{proof}

\subsection{Evolution equation for the profile on the distorted Fourier side} \label{subsec:evolution_equation_profile}

Our starting point for the study of the asymptotic stability of the family of moving kink solutions \eqref{equ:intro_vectorial_moving_kink_family} is the decomposition of the perturbed solution into a modulated kink and a radiation term furnished by Proposition~\ref{prop:setting_up_modulation_orbital}.
In order to infer decay and asymptotics of the radiation term, we will proceed in the spirit of the space-time resonances method and propagate pointwise as well as weighted energy bounds for the effective profile on the distorted Fourier side.
The goal of this subsection is to derive the evolution equation for the effective profile on the distorted Fourier side and to prepare it for the nonlinear analysis in all subsequent sections. 

From now on we work in the moving frame coordinate 
\begin{equation} \label{equ:setting_up_moving_frame_coordinate}
    y := x - q(t),
\end{equation}
and we view the radiation term $\bmu$ as a function of $t$ and $y$.
Then \eqref{equ:setting_up_perturbation_equ} implies that the evolution equation for the radiation term relative to the moving frame coordinate reads
\begin{equation} \label{equ:setting_up_perturbation_equ_moving_frame}
    \partial_t \bmu = \begin{bmatrix}
	\dotq \py & 1 \\
	- L_\ell & \dotq \py
	\end{bmatrix} \bmu + \calMod + \calN(\bmu), \quad L_\ell := -\py^2 - 2 \sech^2(\gamma y) + 1,  
\end{equation} 
where we introduce the short-hand notation
\begin{equation*}
    \calMod := - (\dotq - \ell) \bmY_{0,\ell} - \dotell \bmY_{1,\ell}
\end{equation*}
with 
\begin{equation*} 
 \begin{aligned}
  \bmY_{0,\ell}(y) &= \begin{bmatrix} -\gamma K'(\gamma y) \\ \gamma^2 \ell K''(\gamma y) \end{bmatrix}, \quad
  \bmY_{1,\ell}(y) = \begin{bmatrix} \ell \gamma^3 y K'(\gamma y) \\ - \gamma^3 K'(\gamma y) - \ell^2 \gamma^4 y K''(\gamma y) \end{bmatrix}.
 \end{aligned}
\end{equation*}
Moreover, in the moving frame coordinate the nonlinear term in \eqref{equ:setting_up_perturbation_equ} reads
\begin{equation*}
    \calN(\bmu) = \begin{bmatrix} 0 \\ -W'\bigl(K(\gamma y) + u_1\bigr) + W'\bigl(K(\gamma y)\bigr) + W''\bigl(K(\gamma y)\bigr) u_1 \end{bmatrix}.
\end{equation*}

In order to analyze the dispersive behavior of small solutions to \eqref{equ:setting_up_perturbation_equ_moving_frame}, we need to pass to a time-independent reference operator associated with a suitably chosen fixed Lorentz boost parameter. In what follows we will denote by $\ulell \in (-1,1)$ such a fixed Lorentz boost parameter. Then we write \eqref{equ:setting_up_perturbation_equ_moving_frame} as
\begin{equation} \label{equ:setting_up_perturbation_equ_with_frozen_operator}
    \begin{aligned}
        \partial_t \bmu = \bfL_{\underline{\ell}} \bmu + (\dotq - \ulell) \py \bmu + \calMod + \calN(\bmu) + \calE_1
    \end{aligned}
\end{equation}
with
\begin{equation*}
   \bfL_{\ulell} :=  \begin{bmatrix} \ulell \py & 1 \\ -L_{\ulell} & \ulell \py \end{bmatrix}, \quad L_{\ulell} := -\py^2 - 2 \sech^2(\underline{\gamma} y) + 1, \quad \underline{\gamma} := \frac{1}{\sqrt{1-\ulell^2}},
\end{equation*}
and 
\begin{equation*}
    \calE_1 := \begin{bmatrix} 0 \\ \bigl( 2 \sech^2(\gamma y) - 2 \sech^2(\underline{\gamma}y) \bigr) u_1 \end{bmatrix}.
\end{equation*}

Next, using Lemma~\ref{lem:decom} we decompose the radiation term $\bmu(t)$ into its projection to the essential spectrum and its discrete components relative to the reference operator $\bfL_\ulell$ for some $\ulell \in (-1,1)$, 
\begin{equation} \label{equ:setting_up_decomposition_radiation}
    \bmu(t) = \ulPe \bmu(t) + d_{0,\ulell}(t) \bmY_{0,\ulell} + d_{1,\ulell}(t) \bmY_{1,\ulell},
\end{equation}
where
\begin{equation*}
    \begin{aligned}
    d_{0,\ulell}(t) := \frac{\langle \bfJ \bmY_{1,\ulell}, \bmu(t) \rangle}{\langle \bfJ \bmY_{1,\ulell}, \bmY_{0,\ulell} \rangle}, 
    \quad
    d_{1,\ulell}(t) := - \frac{\langle \bfJ \bmY_{0,\ulell}, \bmu(t) \rangle}{\langle \bfJ \bmY_{1,\ulell}, \bmY_{0,\ulell} \rangle}.
    \end{aligned}
\end{equation*}
We define the profile of the radiation term $\bmu(t)$ relative to the reference operator $\bfL_\ulell$ by 
\begin{equation} \label{equ:setting_up_definition_profile}
    \bmf_{\ulell}(t) = \bigl( f_{\ulell,1}(t), f_{\ulell,2}(t) \bigr) := e^{-t\bfL_\ulell} \bigl( (\ulPe \bmu)(t) \bigr).
\end{equation}
It follows from \eqref{equ:setting_up_perturbation_equ_with_frozen_operator} that the evolution equation for the profile $\bmf_{\ulell}(t)$ is given by
\begin{equation} \label{equ:setting_uP_evolution_equation_profile}
    \pt \bmf_\ulell = e^{-t\bfL_\ulell} \ulPe \Bigl( (\dotq - \ulell) \py \bmu + \calMod + \calN(\bmu) + \calE_1 \Bigr).
\end{equation}

In view of the representation formula \eqref{equ:rep_formula_propagator_modified_dFT} for the evolution $e^{t\bfL_\ulell} \ulPe \bmf_\ulell$ in terms of the distorted Fourier transform and in view of the dispersive decay estimates from Corollary~\ref{cor:linear_asymptotics}, it turns out to be key to propagate pointwise as well as energy bounds on the following quantity on the distorted Fourier side. We refer to it as the effective profile of the radiation term
\begin{equation} \label{equ:setting_up_definition_g}
    \gulellsh(t,\xi) :=\calT_\ulell[ \bmf_\ulell (t)] (\xi)=\calFulellDsh\bigl[ f_{\ulell,1}(t) \bigr](\xi) - \calFulellsh\bigl[ f_{\ulell, 2}(t) \bigr](\xi).         
\end{equation}
It will be useful to introduce a short-hand notation for the two components of the vector $\ulPe \bmu$,
\begin{equation*}
    \ulPe \bmu = \begin{bmatrix} u_{\mathrm{e},1} \\ u_{\mathrm{e},2} \end{bmatrix}.
\end{equation*}
By Proposition~\ref{prop:rep_formula_propagator_modified_dFT} we arrive at the following representation formulas for the components of $\ulPe \bmu$ in terms of the effective profile
\begin{align}
        \usubeone(t,y) &:= \Re \, \Bigl( \calFulellshast \Bigl[ e^{i t (\jxi + \ulell \xi)} i\jxi^{-1} \gulellsh(t,\xi) \Bigr](y) \Bigr), \label{equ:setting_up_representation_formula_usubeone} \\
        \usubetwo(t,y) &:= \Re \, \Bigl( \calFulellDshast \Bigl[ e^{i t (\jxi + \ulell \xi)} i{\jxi}^{-1} \gulellsh(t,\xi) \Bigr](y) \Bigr). \label{equ:setting_up_representation_formula_usubetwo}
\end{align}

By Lemma~\ref{lem:propagator_on_dist_Fourier_side} the evolution equation for the effective profile $g_\ulell^{\#}(t,\xi)$ is given by
\begin{equation} \label{equ:setting_up_g_evol_equ1}
    \begin{aligned}
        \pt g_\ulell^{\#}(t,\xi) = e^{-i t (\jxi + \ulell \xi)} \Bigl( \calF_{\ulell,D}^{\#}\bigl[ F_1(t) \bigr](\xi) - \calF_{\ulell}^{\#}\bigl[ F_2(t) \bigr](\xi) \Bigr),
    \end{aligned}
\end{equation}
where 
\begin{equation*}
    \begin{aligned}
        \bigl( F_1(t), F_2(t) \bigr) := \ulPe \Bigl( (\dotq(t) - \ulell) \py \bmu(t) + \calMod(t) + \calN\bigl(\bmu(t)\bigr) + \calE_1(t) \Bigr).
    \end{aligned}
\end{equation*}
Thus, we find that 
\begin{equation} \label{equ:setting_up_g_evol_equ2}
    \begin{aligned}
        \pt g_\ulell^{\#}(t,\xi) &= e^{-i t (\jxi + \ulell \xi)} \biggl( (\dot{q}(t) - \ulell) \Bigl( \calFulellDsh\bigl[ \py u_1(t) \bigr](\xi) - \calFulellsh\bigl[ \py u_2(t) \bigr](\xi) \Bigr) \\ 
        &\qquad \qquad \qquad \quad + \calFulellDsh \bigl[ \bigl(\calMod(t)\bigr)_1 \bigr](\xi) - \calFulellsh\bigl[ \bigl(\calMod(t)\bigr)_2 \bigr](\xi) \\
        &\qquad \qquad \qquad \quad - \calFulellsh\bigl[ \bigl(\calN(\bmu(t))\bigr)_2 \bigr](\xi)  - \calFulellsh\bigl[ \bigl(\calE_1(t)\bigr)_2 \bigr](\xi) \biggr).
    \end{aligned}
\end{equation}

It remains to prepare some of the terms on the right-hand side of \eqref{equ:setting_up_g_evol_equ2} for the nonlinear analysis.
We begin with the first term on the right-hand side of \eqref{equ:setting_up_g_evol_equ2}.
Integrating by parts in the spatial variable, we find that
\begin{equation} \label{equ:setting_up_pyu_term1}
    \begin{aligned}
        \calF_{\ulell,D}^{\#}\bigl[ \py u_1(t) \bigr](\xi) - \calF_{\ulell}^{\#}\bigl[ \py u_2(t) \bigr](\xi) 
        = i \xi \Bigl( \calF_{\ulell,D}^{\#}\bigl[ u_1(t) \bigr](\xi) - \calF_{\ulell}^{\#}\bigl[ u_2(t) \bigr](\xi) \Bigr) + \calL_\ulell\bigl(\bmu(t)\bigr)(\xi) 
    \end{aligned}
\end{equation}
with 
\begin{equation} \label{equ:setting_up_calL_definition}
    \begin{aligned}
        \calL_\ulell\bigl(\bmu(t)\bigr)(\xi)  &:= \frac{1}{\sqrt{2\pi}} \int_\bbR e^{-iy\xi} \overline{\py m_\ulell^{\#}(\ulg y, \xi)} \bigl( i \jxi u_1(t,y) + u_2(t,y) \bigr) \, \ud y \\ 
        &\quad \quad + \frac{\ulell}{\sqrt{2\pi}} \int_\bbR e^{-iy\xi} \overline{\py^2 m_\ulell^{\#}(\ulg y, \xi)} u_1(t,y) \bigr) \, \ud y.
    \end{aligned}
\end{equation}
Here we emphasize the spatial localization of the coefficients
\begin{equation} \label{equ:setting_up_calL_coefficients_mulell}
    \begin{aligned}
        \py m_\ulell^{\#}(\ulg y,\xi) &= \frac{i \ulg}{\bigl| \ulg (\xi + \ulell \jxi) \bigr| - i} \sech^2(\ulg y), \\
        \py^2 m_\ulell^{\#}(\ulg y,\xi) &= \frac{-2i\ulg^2}{\bigl| \ulg (\xi + \ulell \jxi) \bigr| - i} \sech^2(\ulg y) \tanh(\ulg y).
    \end{aligned}
\end{equation}
Using Corollary~\ref{cor:TellP} and Lemma~\ref{lem:propagator_on_dist_Fourier_side}, we write the first term on the right-hand side of \eqref{equ:setting_up_pyu_term1} as
\begin{equation*}
    \begin{aligned}
        \calF_{\ulell,D}^{\#}\bigl[ u_1(t) \bigr](\xi) - \calF_{\ulell}^{\#}\bigl[ u_2(t) \bigr](\xi) &= \calF_{\ulell,D}^{\#}\bigl[ \bigl(\ulPe u_1(t)\bigr)_1 \bigr](\xi) - \calF_{\ulell}^{\#}\bigl[ \bigl( \ulPe u_2(t)\bigr)_2 \bigr](\xi) \\ 
        &= \calF_{\ulell,D}^{\#}\bigl[ \bigl( e^{t\bfL_\ulell} \ulPe \bmf(t) \bigr)_1 \bigr](\xi) - \calF_{\ulell}^{\#}\bigl[ \bigl( e^{t\bfL_\ulell} \ulPe \bmf(t) \bigr)_2 \bigr](\xi) \\ 
        &= e^{it(\jxi + \ulell \xi)} \gulellsh(t,\xi).
    \end{aligned}
\end{equation*}
Inserting the preceding identity back into \eqref{equ:setting_up_pyu_term1}, we arrive at
\begin{equation} \label{equ:setting_up_pyu_term2}
    \begin{aligned}
        \calF_{\ulell,D}^{\#}\bigl[ \py u_1(t) \bigr](\xi) - \calF_{\ulell}^{\#}\bigl[ \py u_2(t) \bigr](\xi) = i \xi e^{it(\jxi+\ulell \xi)} \gulellsh(t,\xi) + \calL_\ulell\bigl(\bmu(t)\bigr)(\xi). 
    \end{aligned}
\end{equation}

Next, we further decompose the nonlinearities $\calN\bigl(\bmu(t)\bigr)$.
First, we insert the decomposition \eqref{equ:setting_up_decomposition_radiation} of $\bmu(t)$ into its projection to the essential spectrum and its discrete components. 
Then the leading order contributions stem from those terms with all inputs given by the projection to the essential spectrum $(\ulPe \bmu)(t)$, because we can expect the discrete components to decay fast owing to the orthogonality property \eqref{equ:orthogonality_radiation} of the radiation term. Correspondingly, we write 
\begin{equation}
     \calN\bigl(\bmu(t)\bigr) = \calN\bigl(\ulPe \bmu(t)\bigr) + \calE_2(t).
\end{equation}
Then we further expand the leading order nonlinear terms $\calN\bigl(\ulPe \bmu(t)\bigr)$ as
\begin{equation} \label{equ:setting_up_decomposition_nonlinearity}
        \calN\bigl(\ulPe \bmu(t)\bigr) = \begin{bmatrix} 0 \\ \calQ_{\ell(t)}\bigl(\usubeone(t)\bigr) + \calC\bigl(\usubeone(t)\bigr) + \calR_1\bigl(\usubeone(t)\bigr) + \calR_2\bigl(\usubeone(t)\bigr) \end{bmatrix},
\end{equation}
where 
\begin{align*}
    \calQ_{\ell(t)}\bigl(\usubeone(t)\bigr) &:= -\frac12 W^{(3)}\bigl(K(\gamma(t) y)\bigr) \usubeone(t)^2 = -\sech(\gamma(t) y) \tanh(\gamma(t) y) \usubeone(t)^2, \\
    \calC\bigl(\usubeone(t)\bigr) &:= -\frac{1}{3!} W^{(4)}\bigl(K(\gamma(t) y)\bigr) \usubeone(t)^3 = \frac16 \usubeone(t)^3 - \frac13 \sech^2(\gamma(t) y) \usubeone(t)^3, \\
    \calR_1\bigl(\usubeone(t)\bigr) &:= -\frac{1}{4!} W^{(5)}\bigl(K(\gamma(t) y)\bigr) \usubeone(t)^4 = \frac{1}{12} \sech(\gamma(t) y) \tanh(\gamma(t) y) \usubeone(t)^4, \\
    \calR_2\bigl(\usubeone(t)\bigr) &:= -W'\bigl(K(\gamma(t) y) + \usubeone(t)\bigr) + \sum_{k=1}^4 W^{(k+1)}\bigl(K(\gamma(t) y)\bigr) \usubeone(t)^k \\ 
    &= - \frac{1}{4!} \biggl( \int_0^1 (1-r)^4 \cos\bigl(K(\gamma(t) y) + r \usubeone(t) \bigr) \, \ud r \biggr) \usubeone(t)^5.
\end{align*}
In what follows we will use the short-hand notation $\alpha(\gamma(t) y)$ for the variable coefficient of the quadratic nonlinearity, i.e., we write
\begin{equation*}
    \calQ_{\ell(t)}\bigl(\usubeone(t)\bigr) = \alpha(\gamma(t) y) \usubeone(t)^2, \quad \alpha(\gamma(t) y) := -\sech(\gamma(t) y) \tanh(\gamma(t) y).
\end{equation*}
For the subsequent analysis we also need to pass to a time-independent variable coefficient in the quadratic nonlinearity. Correspondingly, we write
\begin{equation*}
    \calQ_{\ell(t)}\bigl(\usubeone(t)\bigr) = \calQ_{\ulell}\bigl(\usubeone(t)\bigr) + \calE_3(t)
\end{equation*}
with
\begin{equation*}
    \begin{aligned}
    \calQ_{\ulell}\bigl(\usubeone(t)\bigr) := \alpha(\ulg y) \usubeone(t)^2, \quad 
    \calE_3(t) := \calQ_{\ell(t)}\bigl(\usubeone(t)\bigr) - \calQ_{\ulell}\bigl(\usubeone(t)\bigr).
    \end{aligned}
\end{equation*}
Moreover, it will be useful to split the cubic nonlinearity into a spatially localized and a spatially non-localized part 
\begin{equation*}
    \calC\bigl(\usubeone(t)\bigr) = \frac16 \usubeone(t)^3 + \calC_{\ell(t)}\bigl( \usubeone(t) \bigr), \quad \calC_{\ell(t)}\bigl( \usubeone(t) \bigr) := -\frac13 \sech^2(\gamma(t) y) \usubeone(t)^3.
\end{equation*}


Combining all of the preceding identities, we arrive at the following evolution equation for the effective profile 
\begin{equation} \label{equ:setting_up_g_evol_equ3}
    \begin{aligned}
        \pt \gulellsh(t,\xi) &= i\xi (\dot{q}(t) - \ulell) \gulellsh(t,\xi) - e^{-i t (\jxi + \ulell \xi)} \calFulellsh\bigl[ \calQ_{\ulell}\bigl(\usubeone(t)\bigr) \bigr](\xi) \\ 
        &\qquad \quad - e^{-i t (\jxi + \ulell \xi)} \calFulellsh\bigl[ {\textstyle \frac16} \usubeone(t)^3 \bigr](\xi) + e^{-i t (\jxi + \ulell \xi)} \widetilde{\calR}(t,\xi),
    \end{aligned}
\end{equation}
where 
\begin{equation} \label{equ:setting_up_definition_wtilR}
    \begin{aligned}
        \widetilde{\calR}(t,\xi) &:= (\dot{q}(t) - \ulell) \calL_\ulell\bigl(\bmu(t)\bigr)(\xi) + \calF_{\ulell,D}^{\#}\bigl[ \bigl(\calMod(t)\bigr)_1 \bigr](\xi) - \calF_{\ulell}^{\#}\bigl[ \bigl(\calMod(t)\bigr)_2 \bigr](\xi) \\ 
        &\quad \quad - \calF_{\ulell}^{\#}\bigl[ \calC_{\ell(t)}\bigl( \usubeone(t) \bigr) \bigr](\xi) - \calF_{\ulell}^{\#}\bigl[ \calR_1\bigl(\usubeone(t)\bigr) \bigr](\xi) - \calF_{\ulell}^{\#}\bigl[ \calR_2\bigl(\usubeone(t)\bigr) \bigr](\xi) \\ 
        &\quad \quad - \calFulellsh\bigl[ \bigl(\calE_1(t)\bigr)_2 \bigr](\xi) - \calFulellsh\bigl[ \bigl(\calE_2(t)\bigr)_2 \bigr](\xi) -\calFulellsh\bigl[ \calE_3(t) \bigr](\xi).
    \end{aligned}
\end{equation}
Heuristically, all terms in $\widetilde{\calR}(t,\xi)$ should be thought of as spatially localized terms with at least cubic-type time decay apart from the contributions of the higher order nonlinearities contained in the term $\calF_{\ulell}^{\#}\bigl[ \calR_2\bigl(u_1(t)\bigr) \bigr](\xi)$.

\subsection{Normal form transformation} \label{subsec:normal_form_transformation}

In this subsection we exploit a remarkable null structure in the quadratic nonlinearity $\calQ_\ulell(\usubeone(t))$ in the evolution equation for the effective profile to cast it into a better form via a normal form transformation. 

Due to the spatial localization of the variable coefficient $\alpha(\ulg y)$, the leading order behavior of the quadratic nonlinearity is determined by the slow local decay of $\usubeone(t,y)$, which is caused by the threshold resonances of the reference operator $\bfL_\ulell$.
Correspondingly, we decompose the radiation term into a term that captures the leading order local decay behavior and a remainder term that enjoys stronger local decay.
To this end we denote by $\chi_{\{\leq 2\ulg|\ulell|\}}(\xi)$ a smooth even bump function with $\chi_{\{\leq 2\ulg|\ulell|\}}(\xi) = 1$ for $|\xi| \leq 2\ulg|\ulell|$.
The choice of this cut-off is informed by the observation that the phase $e^{it(\jxi+\ulell\xi)}$ has a unique stationary point at $\xi=-\ulg\ulell$. 
Then we write 
\begin{equation}
    \usubeone(t,y) = \Re \Bigl( e_\ulell^{\#}\bigl(y, - \ulg \ulell\bigr) h_\ulell(t) \Bigr) + \Remusubeone(t,y), \quad 0 \leq t \leq T,
\end{equation}
with
\begin{equation} \label{equ:setting_up_hulell_definition}
    h_\ulell(t) := \int_\bbR e^{it(\jxi+\ulell \xi)} \chi_{\{\leq 2\ulg|\ulell|\}}(\xi) \, i \jxi^{-1} \gulellsh(t,\xi) \, \ud \xi.
\end{equation}
The remainder term $\Remusubeone(t,y)$ enjoys faster local decay in view of Lemma~\ref{lem:improved_local_decay} and Lemma~\ref{lem:local_decay_resonance_subtracted_off}.
As a result, we arrive at the following decomposition of the quadratic nonlinearity 
\begin{equation*}
    \begin{aligned}
        \alpha(\ulg y) \bigl( \usubeone(t,y) \bigr)^2 &= \alpha(\ulg y) \Bigl( \Re \, \bigl( e_\ulell^{\#}\bigl(y, - \ulg \ulell\bigr) h_\ulell(t) \bigr) \Bigr)^2 \\ 
        &\quad + \alpha(\ulg y) \biggl( \bigl( \usubeone(t,y) \bigr)^2 - \Bigl( \Re \, \bigl( e_\ulell^{\#}\bigl(y, - \ulg \ulell\bigr) h_\ulell(t) \bigr) \Bigr)^2 \biggr),
    \end{aligned}
\end{equation*}
where the leading order behavior stems from the first term on the right-hand side. 

It follows that 
\begin{equation*}
    \begin{aligned}
        \calFulellsh\bigl[ \calQ_{\ulell}\bigl(\usubeone(t)\bigr) \bigr](\xi) &= \frakq_{1,\ulell}(\xi) h_\ulell(t)^2 + \frakq_{2,\ulell}(\xi) h_\ulell(t) \overline{h_\ulell(t)} + \frakq_{3,\ulell}(\xi) \overline{h_\ulell(t)^2} + \calFulellsh\bigl[ \calQ_{\ulell,r}\bigl(\usubeone(t)\bigr) \bigr](\xi), 
    \end{aligned}
\end{equation*}
where we set 
\begin{equation} \label{equ:setting_up_coefficients_normal_form}
    \begin{aligned}
        \frakq_{1,\ulell}(\xi) &:= \frac14 \calFulellsh\Bigl[ \alpha(\ulg y) e_\ulell^{\#}\bigl(y, -\ulg \ulell\bigr)^2 \Bigr](\xi), \\ 
        \frakq_{2,\ulell}(\xi) &:= \frac12 \calFulellsh\Bigl[ \alpha(\ulg y) e_\ulell^{\#}\bigl(y, -\ulg \ulell\bigr) \overline{e_\ulell^{\#}\bigl(y, -\ulg \ulell\bigr)} \Bigr](\xi), \\ 
        \frakq_{3,\ulell}(\xi) &:= \frac14 \calFulellsh\Bigl[ \alpha(\ulg y) \overline{e_\ulell^{\#}\bigl(y, -\ulg \ulell\bigr)}^2 \Bigr](\xi),
    \end{aligned}
\end{equation}
as well as
\begin{equation} \label{equ:setting_up_definition_calQr}
    \begin{aligned}
        \calQ_{\ulell,r}\bigl(\usubeone(t)\bigr) &:= \alpha(\ulg y) \biggl( \bigl( \usubeone(t,y) \bigr)^2 - \Bigl( \Re \, \bigl( e_\ulell^{\#}\bigl(y, - \ulg \ulell\bigr) h_\ulell(t) \bigr) \Bigr)^2 \biggr).
    \end{aligned}
\end{equation}
By a formal stationary phase computation, see Lemma~\ref{lem:core_linear_dispersive_decay}, we expect $h_\ulell(t) \sim e^{it\ulg^{-1}} t^{-\frac12}$ to leading order for $t \gg 1$ so that $\pt \bigl( e^{-it\ulg^{-1}} h_\ulell(t) \bigr)$ should have more time decay. Correspondingly, we filter out the phase and write 
\begin{equation} \label{equ:setting_up_diff_by_parts1}
    \begin{aligned}
        &e^{-it(\jxi+\ulell \xi)} \Bigl( \frakq_{1,\ulell}(\xi) h_\ulell(t)^2 + \frakq_{2,\ulell}(\xi) h_\ulell(t) \overline{h_\ulell(t)} + \frakq_{3,\ulell}(\xi) \overline{h_\ulell(t)}^2 \Bigr) \\
        &= e^{-it(\jxi+\ulell \xi - 2\ulg^{-1})} \frakq_{1,\ulell}(\xi) \bigl( e^{-i t \ulg^{-1}} h_\ulell(t) \bigr)^2 \\ 
        &\quad + e^{-it(\jxi+\ulell \xi)} \frakq_{2,\ulell}(\xi) \bigl( e^{-i t \ulg^{-1}} h_\ulell(t) \bigr) \overline{\bigl( e^{-i t \ulg^{-1}} h_\ulell(t) \bigr)} \\
        &\quad + e^{-it(\jxi+\ulell \xi + 2\ulg^{-1})} \frakq_{3,\ulell}(\xi) \overline{\bigl( e^{-i t \ulg^{-1}} h_\ulell(t) \bigr)}^2.
    \end{aligned}
\end{equation}
Note that in the second and the third term on the right-hand side of \eqref{equ:setting_up_diff_by_parts1}, the phases do not have time resonances.
In fact, a straightforward computation shows that 
\begin{equation} \label{equ:setting_up_phases_lower_bound}
    \jxi+\ulell \xi + 2\ulg^{-1} \geq \jxi + \ulell \xi \geq \ulg^{-1} \quad \text{for all } \xi \in \bbR.
\end{equation}
However, we observe for the first term on the right-hand side of \eqref{equ:setting_up_diff_by_parts1} that
\begin{equation} \label{equ:setting_up_resonant_frequency}
    \jxi + \ulell \xi - 2 \ulg^{-1} = 0 \quad \Leftrightarrow \quad \xi = \ulg (-2\ulell \pm \sqrt{3}) \quad \Leftrightarrow \quad \ulg^{-1} \xi + 2 \ulell = \pm \sqrt{3}.
\end{equation}
At this point we exploit the following remarkable null structure of the quadratic nonlinearity.
It is a generalized version of the null structure observed in \cite[Remark 1.2]{LLSS} and \cite[Lemma 3.1]{LS1} for the special case of odd perturbations of the sine-Gordon kink.

\begin{lemma} \label{lem:null_structure1}
    For any $\ulell \in (-1,1)$ it holds that
    \begin{equation} \label{equ:null_structure1_FT_identity}
        \begin{aligned}
            \calFulellsh\Bigl[ \alpha(\ulg y) e_\ulell^{\#}(y, -\ulg \ulell)^2 \Bigr](\xi) &= \bigl( \jxi + \ulell \xi - 2 \ulg^{-1} \bigr) \psi_\ulell(\xi)
        \end{aligned}
    \end{equation}
    with 
    \begin{equation}
        \begin{aligned}
            \psi_\ulell(\xi) &:= \frac{1}{\sqrt{2\pi}} \frac{1}{48 i \ulg} \frac{1}{|\ulg (\xi + \ulell \jxi)| + i} \sech\Bigl( \frac{\pi}{2} \bigl( \ulg^{-1} \xi + 2 \ulell \bigr) \Bigr) \\ 
            &\quad \quad \times \biggl( 4 \ulell \ulg \Bigl( -5 + \bigl( \ulg^{-1} \xi + 2 \ulell \bigr)^2 \Bigr) + 3 \bigl( \jxi - \ulell \xi + 2 \ulg^{-1} \bigr) \Bigl( 1 + \bigl( \ulg^{-1} \xi + 2\ulell \bigr)^2 \Bigr) \biggr).
        \end{aligned}
    \end{equation}
    In particular, we have      
    \begin{equation*}
        \calFulellsh\Bigl[ \alpha(\ulg y) e_\ulell^{\#}(y, -\ulg \ulell)^2 \Bigr]\bigl( \ulg (-2\ulell \pm \sqrt{3} )\bigr) = 0.
    \end{equation*}
\end{lemma}
\begin{proof}
    By direct computation we obtain 
    \begin{equation*}
        \begin{aligned}
            &\calFulellsh\Bigl[ \alpha(\ulg y) e_\ulell^{\#}(y, -\ulg \ulell)^2 \Bigr](\xi) \\ 
            &= \frac{1}{(2\pi)^{\frac32}} \int_\bbR e^{-iy\xi} \frac{\ulg (\xi + \ulell \jxi) - i \tanh(\ulg y)}{|\ulg (\xi + \ulell \jxi)| + i} \bigl( - \sech(\ulg y) \tanh(\ulg y) \bigr) e^{-i2\ulell\ulg y} \tanh^2(\ulg y) \, \ud y \\ 
            &= \frac{1}{(2\pi)^{\frac32}} \frac{\ulg^{-1}}{|\ulg (\xi + \ulell \jxi)| + i} \biggl( - \ulg (\xi + \ulell \jxi) \int_\bbR e^{-iy ( \ulg^{-1} \xi + 2 \ulell)} \sech(y) \tanh^3(y) \, \ud y \\
            &\qquad \qquad \qquad \qquad \qquad \qquad \qquad \qquad \qquad + i \int_\bbR e^{-i y (\ulg^{-1} \xi + 2 \ulell)} \sech(y) \tanh^4(y) \, \ud y \biggr),
        \end{aligned}
    \end{equation*}
    where in the last step we also made a change of variables.
    Using the Fourier transform identities \eqref{equ:appendix_FT_sech_tanh3} and \eqref{equ:appendix_FT_sech_tanh4}, we find 
    \begin{equation*}
        \begin{aligned}
            &- \ulg (\xi + \ulell \jxi) \int_\bbR e^{-iy (\ulg^{-1}\xi + 2 \ulell)} \sech(y) \tanh^3(y) \, \ud y 
            + i \int_\bbR e^{-i y (\ulg^{-1}\xi+2\ulell)} \sech(y) \tanh^4(y) \, \ud y \\ 
            &=  - \ulg (\xi + \ulell \jxi) i \frac{\pi}{6} \bigl( \ulg^{-1} \xi + 2 \ulell \bigr) \Bigl( -5 + \bigl( \ulg^{-1} \xi + 2 \ulell \bigr)^2 \Bigr) \sech\Bigl( \frac{\pi}{2} \bigl( \ulg^{-1} \xi + 2 \ulell \bigr) \Bigr) \\ 
            &\quad + i \frac{\pi}{24} \Bigl( 9 - 14 \bigl( \ulg^{-1} \xi + 2 \ulell \bigr)^2 + \bigl( \ulg^{-1} \xi + 2 \ulell \bigr)^4 \Bigr) \sech\Bigl( \frac{\pi}{2} \bigl( \ulg^{-1} \xi + 2 \ulell \bigr) \Bigr) \\
            &= i \frac{\pi}{6} \Bigl( - \ulg (\xi+\ulell\jxi) + \bigl( \ulg^{-1} \xi + 2 \ulell \bigr) \Bigr) \Bigl( -5 + \bigl( \ulg^{-1} \xi + 2 \ulell \bigr)^2 \Bigr) \sech\Bigl( \frac{\pi}{2} \bigl( \ulg^{-1} \xi + 2 \ulell \bigr) \Bigr) \\
            &\quad + i \frac{\pi}{24} \biggl( -4 \bigl( \ulg^{-1} \xi + 2 \ulell \bigr)^2 \Bigl( -5 + \bigl( \ulg^{-1} \xi + 2 \ulell \bigr)^2 \Bigr) \\
            &\quad \quad \quad \quad \quad + 9 - 14 \bigl( \ulg^{-1} \xi + 2 \ulell \bigr)^2 + \bigl( \ulg^{-1} \xi + 2 \ulell \bigr)^4 \biggr) \sech\Bigl( \frac{\pi}{2} \ulg^{-1} (\xi + 2 \ulell \ulg) \Bigr).
        \end{aligned}
    \end{equation*}
    Now using the identities
    \begin{equation*}
        \begin{aligned}
        - \ulg (\xi+\ulell\jxi) + \ulg^{-1} (\xi + 2\ulell \ulg) &= - \ulell \ulg \bigl( \jxi + \ulell \xi - 2 \ulg^{-1} \bigr), \\
        3 - \bigl( \ulg^{-1} \xi + 2 \ulell \bigr)^2 &= - \bigl( \jxi + \ulell \xi - 2 \ulg^{-1} \bigr) \bigl( \jxi - \ulell \xi + 2 \ulg^{-1} \bigr),
        \end{aligned}
    \end{equation*}
    as well as 
    \begin{equation*}
        -4\eta^2 ( -5 + \eta^2 ) + 9 - 14 \eta^2 + \eta^4 = 3 (3-\eta^2) (1+\eta^2) \quad \text{with} \quad \eta := \ulg^{-1}(\xi+2\ulell\ulg),
    \end{equation*}
    we arrive at the asserted identity \eqref{equ:null_structure1_FT_identity}.
\end{proof}

\begin{remark} \label{rem:setting_up_rapid_decay_qjs}
    The preceding Lemma~\ref{lem:null_structure1} shows that the coefficient $\frakq_{1,\ulell}(\xi)$ as well as the coefficient $( \jxi + \ulell \xi - 2 \ulg^{-1} )^{-1} \frakq_{1,\ulell}(\xi)$ are rapidly decaying and smooth up to the factor $(|\ulg(\xi + \ulell\jxi)|+i)^{-1}$ stemming from the definition of the distorted Fourier basis element $e_\ulell^{\#}(x,\xi)$.

    One arrives at the same conclusions for the coefficients $\frakq_{2,\ulell}(\xi)$, $( \jxi + \ulell \xi )^{-1} \frakq_{2,\ulell}(\xi)$, as well as for the coefficients $\frakq_{3,\ulell}(\xi)$, $( \jxi + \ulell \xi + 2 \ulg^{-1} )^{-1} \frakq_{3,\ulell}(\xi)$.
\end{remark}

Thanks to Lemma~\ref{lem:null_structure1}, we may rewrite \eqref{equ:setting_up_diff_by_parts1} as
\begin{equation} \label{equ:setting_up_diff_by_parts2}
    \begin{aligned}
        &e^{-it(\jxi+\ulell \xi)} \Bigl( \frakq_{1,\ulell}(\xi) h_\ulell(t)^2 + \frakq_{2,\ulell}(\xi) h_\ulell(t) \overline{h_\ulell(t)} + \frakq_{3,\ulell}(\xi) \overline{h_\ulell(t)}^2 \Bigr) \\ 
        &\quad \quad = \pt \Bigl( B\bigl[ \gulellsh \bigr](t,\xi) \Bigr) + e^{-it(\jxi+\ulell \xi)} \calR_q(t,\xi)
    \end{aligned}
\end{equation}
with 
\begin{equation} \label{equ:setting_up_definition_B}
    \begin{aligned}
        B\bigl[ \gulellsh \bigr](t,\xi) &:= i e^{-it(\jxi+\ulell \xi - 2\ulg^{-1})} \bigl( \jxi + \ulell \xi - 2 \ulg^{-1} \bigr)^{-1} \frakq_{1,\ulell}(\xi) \bigl( e^{-i t \ulg^{-1}} h_\ulell(t) \bigr)^2 \\
        &\quad + i e^{-it(\jxi+\ulell \xi)} \bigl( \jxi + \ulell \xi \bigr)^{-1} \frakq_{2,\ulell}(\xi)  \bigl( e^{-i t \ulg^{-1}} h_\ulell(t) \bigr) \overline{\bigl( e^{-i t \ulg^{-1}} h_\ulell(t) \bigr)} \\
        &\quad + i e^{-it(\jxi+\ulell \xi + 2\ulg^{-1})} \bigl( \jxi + \ulell \xi + 2 \ulg^{-1} \bigr)^{-1} \frakq_{3,\ulell}(\xi) \overline{\bigl( e^{-i t \ulg^{-1}} h_\ulell(t) \bigr)}^2
    \end{aligned}
\end{equation}
and 
\begin{equation} \label{equ:setting_up_definition_calRq}
    \begin{aligned}
        \calR_q(t,\xi) := &- 2i e^{i t 2\ulg^{-1}} \bigl( \jxi + \ulell \xi - 2 \ulg^{-1} \bigr)^{-1} \frakq_{1,\ulell}(\xi) \, \pt \bigl( e^{-i t \ulg^{-1}} h_\ulell(t) \bigr) \bigl( e^{-i t \ulg^{-1}} h_\ulell(t) \bigr) \\ 
        &- 2i \bigl( \jxi + \ulell \xi \bigr)^{-1} \frakq_{2,\ulell}(\xi) \, \Re \Bigl( \pt \bigl( e^{-i t \ulg^{-1}} h_\ulell(t) \bigr) \, \overline{\bigl( e^{-i t \ulg^{-1}} h_\ulell(t) \bigr)} \Bigr) \\
        &- 2i e^{-it 2\ulg^{-1}} \bigl( \jxi + \ulell \xi + 2 \ulg^{-1} \bigr)^{-1} \frakq_{3,\ulell}(\xi) \, \overline{\pt \bigl( e^{-i t \ulg^{-1}} h_\ulell(t) \bigr)} \overline{\bigl( e^{-i t \ulg^{-1}} h_\ulell(t) \bigr)}. 
    \end{aligned}
\end{equation}
Inserting \eqref{equ:setting_up_diff_by_parts2} back into \eqref{equ:setting_up_g_evol_equ3}, we find that
\begin{equation*} 
    \begin{aligned}
        &\pt \Bigl( g_\ulell^{\#}(t,\xi) + B\bigl[ \gulellsh \bigr](t,\xi) \Bigr) \\
        &\quad =  i\xi (\dot{q}(t) - \ulell) g_\ulell^{\#}(t,\xi) - e^{-i t (\jxi + \ulell \xi)} \calFulellsh\bigl[ {\textstyle \frac16} \usubeone(t)^3 \bigr](\xi) + e^{-i t (\jxi + \ulell \xi)} \calR(t,\xi)
    \end{aligned}
\end{equation*}
with 
\begin{equation} \label{equ:setting_up_definition_calR}
    \begin{aligned}
        \calR(t,\xi) &:= \widetilde{\calR}(t,\xi) - \calR_q(t,\xi) - \calF_{\ulell}^{\#}\bigl[ \calQ_{\ulell,r}\bigl(\usubeone(t)\bigr) \bigr](\xi) \\
        &= (\dot{q}(t) - \ulell) \calL_\ulell\bigl(\bmu(t)\bigr)(\xi) + \calF_{\ulell,D}^{\#}\bigl[ \bigl(\calMod(t)\bigr)_1 \bigr](\xi) - \calF_{\ulell}^{\#}\bigl[ \bigl(\calMod(t)\bigr)_2 \bigr](\xi) \\ 
        &\quad \, \,- \calF_{\ulell}^{\#}\bigl[ \calC_{\ell(t)}\bigl( \usubeone(t) \bigr) \bigr](\xi)  - \calF_{\ulell}^{\#}\bigl[ \calR_1\bigl(\usubeone(t)\bigr) \bigr](\xi) - \calF_{\ulell}^{\#}\bigl[ \calR_2\bigl(\usubeone(t)\bigr) \bigr](\xi) \\
        &\quad \, \, - \calFulellsh\bigl[ \bigl(\calE_1(t)\bigr)_2 \bigr](\xi) - \calFulellsh\bigl[ \bigl(\calE_2(t)\bigr)_2 \bigr](\xi) -\calFulellsh\bigl[ \calE_3(t) \bigr](\xi) \\
        &\quad \, \, - \calR_q(t,\xi) - \calF_{\ulell}^{\#}\bigl[ \calQ_{\ulell,r}\bigl(\usubeone(t)\bigr) \bigr](\xi).
    \end{aligned}
\end{equation}
Finally, using an integrating factor we arrive at the following renormalized evolution equation for the effective profile 
\begin{equation} \label{equ:setting_up_g_evol_equ4}
    \begin{aligned}
        \pt \biggl( e^{-i\xi\theta(t)} \Bigl( g_\ulell^{\#}(t,\xi) + B\bigl[ \gulellsh \bigr](t,\xi) \Bigr) \biggr) 
        &= - e^{-i\xi\theta(t)} e^{-i t (\jxi + \ulell \xi)} \calFulellsh\bigl[ {\textstyle \frac16} \usubeone(t)^3 \bigr](\xi) \\
        &\quad \, + e^{-i\xi\theta(t)} e^{-i t (\jxi + \ulell \xi)} \calR(t,\xi) \\
        &\quad \, + e^{-i\xi\theta(t)} i\xi (\dot{q}(t) - \ulell) B\bigl[ \gulellsh \bigr](t,\xi) 
    \end{aligned}
\end{equation}
with
\begin{equation*}
    \theta(t) := \int_0^t \bigl( \dot{q}(s) - \ulell \bigr) \, \ud s.
\end{equation*}
On the right-hand side of \eqref{equ:setting_up_g_evol_equ4} we singled out the spatially non-localized cubic terms, while at least heuristically, all terms in $\calR(t,\xi)$ should be thought of as spatially localized terms with at least cubic-type time decay apart from the contributions of the higher order nonlinearities contained in the term $\calF_{\ulell}^{\#}\bigl[ \calR_2\bigl(u_1(t)\bigr) \bigr](\xi)$.

\subsection{Structure of the cubic nonlinearities} \label{subsec:structure_cubic_nonlinearities}

For the nonlinear analysis in the subsequent sections we need to uncover the fine structure of the non-localized cubic nonlinearities on the right-hand side of the evolution equations \eqref{equ:setting_up_g_evol_equ3} and \eqref{equ:setting_up_g_evol_equ4} for the effective profile.
To this end we write the representation formula \eqref{equ:setting_up_representation_formula_usubeone} for $\usubeone(t,y)$ as 
\begin{equation} \label{equ:setting_up_representation_formula_usubeone_expanded}
    \begin{aligned}
        \usubeone(t,y) &= \Re \, \Bigl( \calFulellshast \Bigl[ e^{i t (\jxi + \ulell \xi)} i\jxi^{-1} \gulellsh(t,\xi) \Bigr](y) \Bigr) \\ 
        &= \frac{i}{2} \int_\bbR \eulsharp(y,\xi) \, e^{it(\jxi+\ulell\xi)} \jxi^{-1} \gulellsh(t,\xi) \, \ud \xi \\
        &\quad - \frac{i}{2} \int_\bbR \overline{\eulsharp(y,\xi)} e^{-it(\jxi+\ulell\xi)} \jxi^{-1} \overline{\gulellsh(t,\xi)} \, \ud \xi.
    \end{aligned}
\end{equation}
Then inserting \eqref{equ:setting_up_representation_formula_usubeone_expanded} into the non-localized cubic nonlinearity, we find by direct computation 
\begin{equation} \label{equ:setting_up_dist_FT_of_cubic_expanded}
        e^{-i t (\jxi + \ulell \xi)} \calFulellsh\bigl[ \usubeone(t)^3 \bigr](\xi) 
        = - \frac{i}{8} \calI_1(t,\xi) + \frac{3i}{8} \calI_2(t,\xi) - \frac{3i}{8} \calI_3(t,\xi) + \frac{i}{8} \calI_4(t,\xi),
\end{equation}
where
\begin{equation} \label{equ:definition_calI}
    \begin{aligned}
        \calI_1(t,\xi) &:= \iiint e^{it \Phi_{1,\ulell}(\xi,\xi_1,\xi_2,\xi_3)} \gulellsh(t,\xi_1) \, \gulellsh(t,\xi_2) \, \gulellsh(t,\xi_3) \\ 
        &\qquad \qquad \qquad \qquad \qquad \qquad \times \jxione^{-1} \jxitwo^{-1} \jxithree^{-1} \, \nu_{\ulell, +++}(\xi, \xi_1, \xi_2, \xi_3) \, \ud \xi_1 \, \ud \xi_2 \, \ud \xi_3 \\ 
        \calI_2(t,\xi) &:= \iiint e^{it \Phi_{2,\ulell}(\xi,\xi_1,\xi_2,\xi_3)} \gulellsh(t,\xi_1) \, \overline{\gulellsh(t,\xi_2)} \, \gulellsh(t,\xi_3) \\
        &\qquad \qquad \qquad \qquad \qquad \qquad \times \jxione^{-1} \jxitwo^{-1} \jxithree^{-1} \nu_{\ulell, +-+}(\xi, \xi_1, \xi_2, \xi_3) \, \ud \xi_1 \, \ud \xi_2 \, \ud \xi_3 \\ 
        \calI_3(t,\xi) &:= \iiint e^{it\Phi_{3,\ulell}(\xi,\xi_1,\xi_2,\xi_3)} \gulellsh(t,\xi_1) \, \overline{\gulellsh(t,\xi_2)} \, \overline{\gulellsh(t,\xi_3)} \\
        &\qquad \qquad \qquad \qquad \qquad \qquad \times \jxione^{-1} \jxitwo^{-1} \jxithree^{-1} \, \nu_{\ulell, +--}(\xi, \xi_1, \xi_2, \xi_3) \, \ud \xi_1 \, \ud \xi_2 \, \ud \xi_3 \\ 
        \calI_4(t,\xi) &:= \iiint e^{it\Phi_{4,\ulell}(\xi,\xi_1,\xi_2,\xi_3)} \overline{\gulellsh(t,\xi_1)} \, \overline{\gulellsh(t,\xi_2)} \, \overline{\gulellsh(t,\xi_3)} \\
        &\qquad \qquad \qquad \qquad \qquad \qquad \times \jxione^{-1} \jxitwo^{-1} \jxithree^{-1} \, \nu_{\ulell, ---}(\xi, \xi_1, \xi_2, \xi_3) \, \ud \xi_1 \, \ud \xi_2 \, \ud \xi_3
    \end{aligned}
\end{equation}
with
\begin{equation}
\begin{aligned}
    \Phi_{1,\ulell}(\xi,\xi_1,\xi_2,\xi_3) &:= -(\jxi + \ulell \xi) + (\jxione + \ulell \xi_1) + (\jxitwo + \ulell \xi_2) + (\jxithree + \ulell \xi_3), \\
    \Phi_{2,\ulell}(\xi,\xi_1,\xi_2,\xi_3) &:= -(\jxi + \ulell \xi) + (\jxione + \ulell \xi_1) - (\jxitwo + \ulell \xi_2) + (\jxithree + \ulell \xi_3), \\
    \Phi_{3,\ulell}(\xi,\xi_1,\xi_2,\xi_3) &:= -(\jxi + \ulell \xi) + (\jxione + \ulell \xi_1) - (\jxitwo + \ulell \xi_2) - (\jxithree + \ulell \xi_3), \\
    \Phi_{4,\ulell}(\xi,\xi_1,\xi_2,\xi_3) &:= -(\jxi + \ulell \xi) - (\jxione + \ulell \xi_1) - (\jxitwo + \ulell \xi_2) - (\jxithree + \ulell \xi_3), 
\end{aligned}
\end{equation}
and
\begin{equation} \label{equ:definition_cubic_spectral_distributions}
\begin{aligned}
    \nu_{\ulell, +++}(\xi, \xi_1, \xi_2, \xi_3) &:= \int_\bbR \overline{\eulsharp(y,\xi)} \, \eulsharp(y,\xi_1) \, \eulsharp(y,\xi_2) \, \eulsharp(y,\xi_3) \, \ud y, \\
    \nu_{\ulell, +-+}(\xi, \xi_1, \xi_2, \xi_3) &:= \int_\bbR \overline{\eulsharp(y,\xi)} \, \eulsharp(y,\xi_1) \, \overline{\eulsharp(y,\xi_2)} \, \eulsharp(y,\xi_3) \, \ud y, \\
    \nu_{\ulell, +--}(\xi, \xi_1, \xi_2, \xi_3) &:= \int_\bbR \overline{\eulsharp(y,\xi)} \, \eulsharp(y,\xi_1) \, \overline{\eulsharp(y,\xi_2)} \, \overline{\eulsharp(y,\xi_3)} \, \ud y, \\
    \nu_{\ulell, ---}(\xi, \xi_1, \xi_2, \xi_3) &:= \int_\bbR \overline{\eulsharp(y,\xi)} \, \overline{\eulsharp(y,\xi_1)} \, \overline{\eulsharp(y,\xi_2)} \, \overline{\eulsharp(y,\xi_3)} \, \ud y.
\end{aligned}
\end{equation}

Next, we insert the decompositions of the cubic spectral distributions \eqref{equ:definition_cubic_spectral_distributions} determined in Subsection~\ref{subsec:cubic_spectral_distributions} into the cubic interaction terms $\calI_j(t,\xi)$, $1 \leq j \leq 4$, defined in \eqref{equ:definition_calI}. 
This leads to decompositions
\begin{equation} \label{equ:calIj_decomposition}
 \calI_j(t,\xi) = \calI_j^{\delta_0}(t,\xi) + \calI_j^{\pvdots}(t,\xi) + \calI_j^{\mathrm{reg}}(t,\xi), \quad 1 \leq j \leq 4.
\end{equation}
The cubic interaction terms with a Dirac kernel in \eqref{equ:calIj_decomposition} are given by 
\begin{equation} \label{equ:structure_cubic_nonlinearities_dirac_kernel}
    \begin{aligned}
        \calI_1^{\delta_0}(t,\xi) &= \iint e^{it \Psi_{1}(\xi,\xi_1,\xi_2)} \gulellsh(t,\xi_1) \, \gulellsh(t,\xi_2) \, \gulellsh(t,\xi_3) \, \jxione^{-1} \jxitwo^{-1} \jxithree^{-1} \\ 
        &\qquad \qquad \qquad \qquad \qquad \qquad \qquad \qquad \qquad \times \frakm_{\ulell, +++}^{\delta_0}(\xi,\xi_1,\xi_2,\xi_3) \, \ud \xi_1 \, \ud \xi_2 \\ 
        \calI_2^{\delta_0}(t,\xi) &= \iint e^{it \Psi_{2}(\xi,\xi_1,\xi_2)} \gulellsh(t,\xi_1) \, \overline{\gulellsh(t,\xi_2)} \, \gulellsh(t,\xi_3) \, \jxione^{-1} \jxitwo^{-1} \jxithree^{-1} \\ 
        &\qquad \qquad \qquad \qquad \qquad \qquad \qquad \qquad \qquad \times \frakm_{\ulell, +-+}^{\delta_0}(\xi,\xi_1,\xi_2,\xi_3) \, \ud \xi_1 \, \ud \xi_2 \\ 
        \calI_3^{\delta_0}(t,\xi) &= \iint e^{it \Psi_{3}(\xi,\xi_1,\xi_2)} \gulellsh(t,\xi_1) \, \overline{\gulellsh(t,\xi_2)} \, \overline{\gulellsh(t,\xi_3)} \, \jxione^{-1} \jxitwo^{-1} \jxithree^{-1} \\ 
        &\qquad \qquad \qquad \qquad \qquad \qquad \qquad \qquad \qquad \times \frakm_{\ulell, +--}^{\delta_0}(\xi,\xi_1,\xi_2,\xi_3) \, \ud \xi_1 \, \ud \xi_2 \\         
        \calI_4^{\delta_0}(t,\xi) &= \iint e^{it \Psi_{4}(\xi,\xi_1,\xi_2)} \overline{\gulellsh(t,\xi_1)} \, \overline{\gulellsh(t,\xi_2)} \, \overline{\gulellsh(t,\xi_3)} \, \jxione^{-1} \jxitwo^{-1} \jxithree^{-1} \\ 
        &\qquad \qquad \qquad \qquad \qquad \qquad \qquad \qquad \qquad \times \frakm_{\ulell, ---}^{\delta_0}(\xi,\xi_1,\xi_2,\xi_3) \, \ud \xi_1 \, \ud \xi_2,                 
    \end{aligned}
\end{equation}
where
\begin{equation*}
    \xi_3 := \left\{ \begin{aligned}
                        &\xi-\xi_1-\xi_2, \quad \quad \, j = 1, \\
                        &\xi-\xi_1+\xi_2, \quad \quad \, j = 2, \\
                        &-\xi+\xi_1-\xi_2, \quad j = 3, \\
                        &-\xi-\xi_1-\xi_2, \quad j = 4,
                     \end{aligned} \right.
\end{equation*}
and 
\begin{equation}
    \begin{aligned}
        \Psi_{1}(\xi,\xi_1,\xi_2) &:= \Phi_{1,\ulell}(\xi,\xi_1,\xi_2,\xi-\xi_1-\xi_2) = -\jxi + \jxione + \jxitwo + \jap{\xi-\xi_1-\xi_2}, \\
        \Psi_{2}(\xi,\xi_1,\xi_2) &:= \Phi_{2,\ulell}(\xi,\xi_1,\xi_2,\xi-\xi_1+\xi_2) = -\jxi + \jxione - \jxitwo + \jap{\xi-\xi_1+\xi_2}, \\
        \Psi_{3}(\xi,\xi_1,\xi_2) &:= \Phi_{3,\ulell}(\xi,\xi_1,\xi_2,-\xi+\xi_1-\xi_2) = -\jxi + \jxione - \jxitwo - \jap{\xi-\xi_1+\xi_2}, \\
        \Psi_{4}(\xi,\xi_1,\xi_2) &:= \Phi_{4,\ulell}(\xi,\xi_1,\xi_2,-\xi-\xi_1-\xi_2) = -\jxi - \jxione - \jxitwo - \jap{\xi+\xi_1+\xi_2}.
    \end{aligned}
\end{equation}
The cubic interaction terms with a Hilbert-type kernel in \eqref{equ:calIj_decomposition} can be written as
\begin{equation} \label{equ:structure_cubic_nonlinearities_PV}
    \begin{aligned}
        \calI_1^{\pvdots}(t,\xi) &= \iiint e^{it \Omega_{1,\ulell}(\xi,\xi_1,\xi_2,\xi_4)} \gulellsh(t,\xi_1) \, \gulellsh(t,\xi_2) \, \gulellsh(t,\xi_3) \, \jxione^{-1} \jxitwo^{-1} \jxithree^{-1} \\ 
        &\qquad \qquad \qquad \qquad \times \frakm_{\ulell, +++}^{\pvdots}(\xi, \xi_1, \xi_2, \xi_3) \, \pvdots \cosech\Bigl( \frac{\pi}{2\ulg} \xi_4 \Bigr) \, \ud \xi_1 \, \ud \xi_2 \, \ud \xi_4 \\
        \calI_2^{\pvdots}(t,\xi) &= \iiint e^{it \Omega_{2,\ulell}(\xi,\xi_1,\xi_2,\xi_4)} \gulellsh(t,\xi_1) \, \overline{\gulellsh(t,\xi_2)} \, \gulellsh(t,\xi_3) \, \jxione^{-1} \jxitwo^{-1} \jxithree^{-1} \\ 
        &\qquad \qquad \qquad \qquad \times \frakm_{\ulell, +-+}^{\pvdots}(\xi, \xi_1, \xi_2, \xi_3) \, \pvdots \cosech\Bigl( \frac{\pi}{2\ulg} \xi_4 \Bigr) \, \ud \xi_1 \, \ud \xi_2 \, \ud \xi_4 \\        
        \calI_3^{\pvdots}(t,\xi) &= \iiint e^{it \Omega_{3,\ulell}(\xi,\xi_1,\xi_2,\xi_4)} \gulellsh(t,\xi_1) \, \overline{\gulellsh(t,\xi_2)} \, \overline{\gulellsh(t,\xi_3)} \, \jxione^{-1} \jxitwo^{-1} \jxithree^{-1} \\ 
        &\qquad \qquad \qquad \qquad \times \frakm_{\ulell, +--}^{\pvdots}(\xi, \xi_1, \xi_2, \xi_3) \, \pvdots \cosech\Bigl( \frac{\pi}{2\ulg} \xi_4 \Bigr) \, \ud \xi_1 \, \ud \xi_2 \, \ud \xi_4 \\                
        \calI_4^{\pvdots}(t,\xi) &= \iiint e^{it \Omega_{4,\ulell}(\xi,\xi_1,\xi_2,\xi_4)} \overline{\gulellsh(t,\xi_1)} \, \overline{\gulellsh(t,\xi_2)} \, \overline{\gulellsh(t,\xi_3)} \, \jxione^{-1} \jxitwo^{-1} \jxithree^{-1} \\ 
        &\qquad \qquad \qquad \qquad \times \frakm_{\ulell, ---}^{\pvdots}(\xi, \xi_1, \xi_2, \xi_3) \, \pvdots \cosech\Bigl( \frac{\pi}{2\ulg} \xi_4 \Bigr) \, \ud \xi_1 \, \ud \xi_2 \, \ud \xi_4,                       
    \end{aligned}
\end{equation}
where
\begin{equation*}
    \xi_3 := \left\{ \begin{aligned}
                        &\xi-\xi_1-\xi_2+\xi_4, \quad \quad \, j = 1, \\
                        &\xi-\xi_1+\xi_2+\xi_4, \quad \quad \, j = 2, \\
                        &-\xi+\xi_1-\xi_2-\xi_4, \quad j = 3, \\
                        &-\xi-\xi_1-\xi_2-\xi_4, \quad j = 4,
                     \end{aligned} \right.
\end{equation*}
and
\begin{equation}
    \begin{aligned}
        \Omega_{1,\ulell}(\xi,\xi_1,\xi_2,\xi_4) &:= \Phi_{1,\ulell}(\xi,\xi_1,\xi_2,\xi-\xi_1-\xi_2+\xi_4) 
        = -\jxi + \jxione + \jxitwo + \jxithree + \ulell \xi_4, \\
        \Omega_{2,\ulell}(\xi,\xi_1,\xi_2,\xi_4) &:= \Phi_{2,\ulell}(\xi,\xi_1,\xi_2,\xi-\xi_1+\xi_2+\xi_4) 
        = -\jxi + \jxione - \jxitwo + \jxithree + \ulell \xi_4, \\
        \Omega_{3,\ulell}(\xi,\xi_1,\xi_2,\xi_4) &:= \Phi_{3,\ulell}(\xi,\xi_1,\xi_2,-\xi+\xi_1-\xi_2-\xi_4)  
        = -\jxi + \jxione - \jxitwo - \jxithree + \ulell \xi_4, \\
        \Omega_{4,\ulell}(\xi,\xi_1,\xi_2,\xi_4) &:= \Phi_{4,\ulell}(\xi,\xi_1,\xi_2,-\xi-\xi_1-\xi_2-\xi_4) 
        = -\jxi - \jxione - \jxitwo - \jxithree + \ulell \xi_4.
    \end{aligned}
\end{equation}

Finally, the regular cubic interaction terms in \eqref{equ:calIj_decomposition} are given by
\begin{equation} \label{equ:structure_cubic_nonlinearities_regular}
    \begin{aligned}
        \calI_1^{\mathrm{reg}}(t,\xi) &:= \iiint e^{it \Phi_{1,\ulell}(\xi,\xi_1,\xi_2,\xi_3)} \gulellsh(t,\xi_1) \, \gulellsh(t,\xi_2) \, \gulellsh(t,\xi_3) \\ 
        &\qquad \qquad \qquad \qquad \qquad \qquad \times \jxione^{-1} \jxitwo^{-1} \jxithree^{-1} \, \nu_{\ulell, +++}^{\mathrm{reg}}(\xi, \xi_1, \xi_2, \xi_3) \, \ud \xi_1 \, \ud \xi_2 \, \ud \xi_3 \\ 
        \calI_2^{\mathrm{reg}}(t,\xi) &:= \iiint e^{it \Phi_{2,\ulell}(\xi,\xi_1,\xi_2,\xi_3)} \gulellsh(t,\xi_1) \, \overline{\gulellsh(t,\xi_2)} \, \gulellsh(t,\xi_3) \\
        &\qquad \qquad \qquad \qquad \qquad \qquad \times \jxione^{-1} \jxitwo^{-1} \jxithree^{-1} \nu_{\ulell, +-+}^{\mathrm{reg}}(\xi, \xi_1, \xi_2, \xi_3) \, \ud \xi_1 \, \ud \xi_2 \, \ud \xi_3 \\ 
        \calI_3^{\mathrm{reg}}(t,\xi) &:= \iiint e^{it\Phi_{3,\ulell}(\xi,\xi_1,\xi_2,\xi_3)} \gulellsh(t,\xi_1) \, \overline{\gulellsh(t,\xi_2)} \, \overline{\gulellsh(t,\xi_3)} \\
        &\qquad \qquad \qquad \qquad \qquad \qquad \times \jxione^{-1} \jxitwo^{-1} \jxithree^{-1} \, \nu_{\ulell, +--}^{\mathrm{reg}}(\xi, \xi_1, \xi_2, \xi_3) \, \ud \xi_1 \, \ud \xi_2 \, \ud \xi_3 \\ 
        \calI_4^{\mathrm{reg}}(t,\xi) &:= \iiint e^{it\Phi_{4,\ulell}(\xi,\xi_1,\xi_2,\xi_3)} \overline{\gulellsh(t,\xi_1)} \, \overline{\gulellsh(t,\xi_2)} \, \overline{\gulellsh(t,\xi_3)} \\
        &\qquad \qquad \qquad \qquad \qquad \qquad \times \jxione^{-1} \jxitwo^{-1} \jxithree^{-1} \, \nu_{\ulell, ---}^{\mathrm{reg}}(\xi, \xi_1, \xi_2, \xi_3) \, \ud \xi_1 \, \ud \xi_2 \, \ud \xi_3.
    \end{aligned}
\end{equation}

\section{Bootstrap Setup and Proof of Theorem~\ref{thm:main}} \label{sec:bootstrap_setup_proof_thm}

The strategy of the proof of Theorem~\ref{thm:main} is an intertwined bootstrap argument that is structured around two central bootstrap propositions.
In Subsection~\ref{subsec:bootstrap_propositions}, we present the statements of the two propositions.
Moreover, as consequences of the bootstrap assumptions we derive several decay estimates and related bounds that will be utilized repeatedly in all subsequent sections. 
Subsection~\ref{subsec:proof_of_main_theorem} provides a concise proof of Theorem 1.1 based on the conclusions of the two bootstrap propositions.
The proofs of the two propositions then comprise the remainder of the paper.

\subsection{The main bootstrap propositions} \label{subsec:bootstrap_propositions}

The starting point for the proof of Theorem~\ref{thm:main} are the local existence result from Lemma~\ref{lem:setting_up_local_existence} and the modulational setup established in Proposition~\ref{prop:setting_up_modulation_orbital}.
For given $\ell_0 \in (-1,1)$, $x_0 \in \bbR$, and inital data $\bmu_0 = (u_{0,1}, u_{0,2}) \in H^3_x \times H^2_x$,
Lemma~\ref{lem:setting_up_local_existence} furnishes a unique $H^3_x \times H^2_x$-solution $\bmphi(t)$ to the sine-Gordon equation~\eqref{equ:intro_sG_1st_order} with initial data
\begin{equation}
    \bmphi(0,x) = \bmK_{\ell_0,x_0}(x) + \bmu_0(x-x_0)
\end{equation}
defined on a maximal interval of existence $[0,T_\ast)$ for some $0 < T_\ast \leq \infty$. 
Under the smallness assumption $\|\bmu_0\|_{H^1_x\times L^2_x} \leq \varepsilon_1$, by Proposition~\ref{prop:setting_up_modulation_orbital} there exist unique continuously differentiable paths $\ell \colon [0,T_\ast) \to (-1,1)$ and $q\colon [0,T_\ast) \to \bbR$ satisfying $(1)$--$(6)$ in the statement of Proposition~\ref{prop:setting_up_modulation_orbital}.
In particular, the solution $\bmphi(t)$ can be decomposed into a modulated kink and a radiation term
\begin{equation} 
    \bmphi(t,x) = \bmK_{\ell(t), q(t)}(x) + \bmu\bigl(t, x-q(t)\bigr), \quad 0 \leq t < T_\ast.
\end{equation}
Here, $0 < \varepsilon_1 \ll 1$ is the small absolute constant from the statement of Proposition~\ref{prop:setting_up_modulation_orbital}.

Our task is now to simultaneously establish decay as well as asymptotics for the radiation term~$\bmu(t)$ and to determine the asymptotic behavior of the modulation parameters $\ell(t)$ and $q(t)$. Along the way we of course also need to conclude that $T_\ast = \infty$.
In order to analyze the long-time behavior of the radiation term, we need to pass to a time-independent reference operator $\bfL_{\ulell}$ associated with a suitably chosen fixed Lorentz boost parameter $\ulell \in (-1,1)$ with $|\ulell-\ell_0| \leq \frac14 \gamma_0^{-2}$. 
We recall from Subsection~\ref{subsec:evolution_equation_profile} that the profile of the radiation term relative to the reference operator $\bfL_\ulell$ is defined by
\begin{equation}
    \bmf_\ulell(t) = \bigl( f_{\ulell,1}(t), f_{\ulell,2}(t) \bigr) := e^{-t\bfL_\ulell} \bigl( \ulPe \bmu(t) \bigr).
\end{equation}
The associated effective profile on the distorted Fourier side is defined as
\begin{equation}
    g_\ulell^\#(t,\xi) := \calF_{\ulell, D}^{\#}\bigl[ f_{\ulell,1}(t) \bigr](\xi) - \calF_\ulell^\#\bigl[f_{\ulell,2}(t)](\xi).
\end{equation}
In order to infer decay and ultimately asymptotics of the radiation term, we proceed in the spirit of the space-time resonances method and seek to control the following pointwise as well as energy bounds for the effective profile on a given time interval $[0,T]$ for some $0 < T < T_\ast$,
\begin{equation} \label{equ:bootstrap_norm_XT_definition}
    \begin{aligned}
        \bigl\| g_\ulell^\# \bigr\|_{X(T)} &:= \sup_{0 \leq t \leq T} \, \biggl( \bigl\| \jxi^{\frac32} g_\ulell^\#(t,\xi) \bigr\|_{L^\infty_\xi} + \jt^{-\delta} \bigl\| \jxi^2 g_\ulell^\#(t,\xi) \bigr\|_{L^2_\xi} \\
        &\qquad \qquad \quad + \jt^{-\delta} \bigl\| \pxi g_\ulell^\#(t,\xi) \bigr\|_{L^2_\xi} + \jt^{-2\delta} \bigl\| \jxi^2 \pxi g_\ulell^\#(t,\xi) \bigr\|_{L^2_\xi} \biggr).
    \end{aligned}
\end{equation}
Here $0 < \delta \ll 1$ is a small absolute constant.
The entire analysis in the remainder of this paper is based on the evolution equation for the effective profile $\gulellsh(t,\xi)$, which has been set up in Subsections~\ref{subsec:evolution_equation_profile}--\ref{subsec:structure_cubic_nonlinearities}.
We refer the reader to Subsection~\ref{subsec:overview_bootstrap_setup} for a discussion of the origin of the two differing slow growth rates of the weighted energy norms in the second line of \eqref{equ:bootstrap_norm_XT_definition}.

In the following first core bootstrap proposition we deduce a convergence rate of the Lorentz boost modulation parameter $\ell(t)$ to its final value $\ell(T)$ on a given time interval $[0,T]$ with $0 < T < T_\ast$.
We defer the proof of the proposition to Section~\ref{sec:modulation_control}.

\begin{proposition}[Control of the modulation parameters] \label{prop:modulation_parameter_control}
    Let $\ell_0 \in (-1,1)$ and let $0 < \varepsilon_1 \ll 1$ be the small constant from the statement of Proposition~\ref{prop:setting_up_modulation_orbital}.
    There exist constants $C_0 \geq 1$ and $0 < \varepsilon_0 \ll \varepsilon_1 \ll 1$ with the following properties:
    Let $x_0 \in \bbR$ and let $\bmu_0 = (u_{0,1}, u_{0,2}) \in H^3_x(\bbR) \times H^2_x(\bbR)$ with $\varepsilon := \|\jx \bmu_0\|_{H^3_x \times H^2_x} \leq \varepsilon_0$. 
    Denote by $\bmphi(t)$ the $H^3_x \times H^2_x$ solution to \eqref{equ:intro_sG_1st_order} with initial data $\bmphi(0,x) = \bmK_{\ell_0,x_0}(x) + \bmu_0(x-x_0)$ on its maximal interval of existence $[0,T_\ast)$ furnished by Lemma~\ref{lem:setting_up_local_existence}. 
    Let $\ell \colon [0,T_\ast) \to (-1,1)$ and $q \colon [0,T_\ast) \to \bbR$ be the unique continuously differentiable paths satisfying (1)--(6) in the statement of Proposition~\ref{prop:setting_up_modulation_orbital}.
    Fix arbitrary $0 < T < T_\ast$ and $\ulell \in (-1,1)$ with $|\ulell - \ell_0| \leq \frac14 \gamma_0^{-2}$.
    Suppose 
    \begin{align}
        \sup_{0 \leq t \leq T} \, \jt^{1-\delta} |\ell(t) - \ulell| &\leq 2 C_0 \varepsilon, \label{equ:prop_modulation_parameters_assumption1} \\ 
        \bigl\| g_\ulell^\# \bigr\|_{X(T)} &\leq 2 C_0 \varepsilon. \label{equ:prop_modulation_parameters_assumption2}
    \end{align}
    Then it follows that  
    \begin{equation} \label{equ:prop_modulation_parameters_conclusion}
        \sup_{0 \leq t \leq T} \, \jt^{1-\delta} |\ell(t) - \ell(T)| \leq C_0 \varepsilon. 
    \end{equation}
\end{proposition}

In the following second core bootstrap proposition we obtain control of the norm~\eqref{equ:bootstrap_norm_XT_definition} for the effective profile. We defer the proof to Section~\ref{sec:energy_estimates} and Section~\ref{sec:pointwise_profile_bounds}. Specifically, it is an immediate consequence of Proposition~\ref{prop:sobolev_profile_bounds}, Proposition~\ref{prop:pxi_profile_bounds}, Proposition~\ref{prop:pointwise_profile_bounds}, and the local existence theory.

\begin{proposition}[Profile bounds] \label{prop:profile_bounds}
    Let $\ell_0 \in (-1,1)$ and let $0 < \varepsilon_1 \ll 1$ be the small constant from the statement of Proposition~\ref{prop:setting_up_modulation_orbital}.
    There exist constants $C_0 \geq 1$ and $0 < \varepsilon_0 \ll \varepsilon_1 \ll 1$ with the following properties:
    Let $x_0 \in \bbR$ and let $\bmu_0 = (u_{0,1}, u_{0,2}) \in H^3_x(\bbR) \times H^2_x(\bbR)$ with $\varepsilon := \|\jx \bmu_0\|_{H^3_x \times H^2_x} \leq \varepsilon_0$. 
    Denote by $\bmphi(t)$ the $H^3_x \times H^2_x$ solution to \eqref{equ:intro_sG_1st_order} with initial data $\bmphi(0,x) = \bmK_{\ell_0,x_0}(x) + \bmu_0(x-x_0)$ on its maximal interval of existence $[0,T_\ast)$ furnished by Lemma~\ref{lem:setting_up_local_existence}. 
    Let $\ell \colon [0,T_\ast) \to (-1,1)$ and $q \colon [0,T_\ast) \to \bbR$ be the unique  continuously differentiable paths satisfying (1)--(6) in the statement of Proposition~\ref{prop:setting_up_modulation_orbital}.
    Fix arbitrary $0 < T < T_\ast$ and $\ulell \in (-1,1)$ with $|\ulell - \ell_0| \leq \frac14 \gamma_0^{-2}$. 
    Suppose 
    \begin{align} 
        \sup_{0 \leq t \leq T} \, \jt^{1-\delta} |\ell(t) - \ulell| &\leq 2 C_0 \varepsilon, \label{equ:prop_profile_bounds_assumption1} \\
        \bigl\| g_\ulell^\# \bigr\|_{X(T)} &\leq 2 C_0 \varepsilon. \label{equ:prop_profile_bounds_assumption2}
    \end{align}
    Then we have 
    \begin{equation} \label{equ:prop_profile_bounds_conclusion}
        \bigl\| g_\ulell^\# \bigr\|_{X(T)} \leq C_0 \varepsilon.
    \end{equation}
\end{proposition}

In the next corollary we collect decay estimates and related bounds that are immediate consequences of the bootstrap assumptions in the statements of Proposition~\ref{prop:modulation_parameter_control} and Proposition~\ref{prop:profile_bounds}.
These estimates will be used over and over in the remainder of this paper. 

\begin{corollary} \label{cor:consequences_bootstrap_assumptions}
    Let $\ell_0 \in (-1,1)$, and let $C_0 \geq 1$ and $0 < \varepsilon_0 \ll \varepsilon_1 \ll 1$ be the constants from the statements of Proposition~\ref{prop:setting_up_modulation_orbital}, Proposition~\ref{prop:modulation_parameter_control}, and Proposition~\ref{prop:profile_bounds}. Let $x_0 \in \bbR$ and let $\bmu_0 = (u_{0,1}, u_{0,2}) \in H^3_x(\bbR) \times H^2_x(\bbR)$ with $\varepsilon := \|\jx \bmu_0\|_{H^3_x \times H^2_x} \leq \varepsilon_0$. 
    Denote by $\bmphi(t)$ the $H^3_x \times H^2_x$ solution to \eqref{equ:intro_sG_1st_order} with initial data $\bmphi(0,x) = \bmK_{\ell_0,x_0}(x) + \bmu_0(x-x_0)$ on its maximal interval of existence $[0,T_\ast)$ furnished by Lemma~\ref{lem:setting_up_local_existence}. 
    Let $\ell \colon [0,T_\ast) \to (-1,1)$ and $q \colon [0,T_\ast) \to \bbR$ be the unique continuously differentiable paths satisfying (1)--(6) in the statement of Proposition~\ref{prop:setting_up_modulation_orbital}.
    Fix arbitrary $0 < T < T_\ast$ and $\ulell \in (-1,1)$ with $|\ulell - \ell_0| \leq \frac14 \gamma_0^{-2}$.
    Suppose  
    \begin{align} 
        \sup_{0 \leq t \leq T} \, \jt^{1-\delta} |\ell(t) - \ulell| &\leq 2 C_0 \varepsilon, \label{equ:consequences_cor_assumption1} \\ 
        \sup_{0 \leq t \leq T} \, \bigl\| g_\ulell^\# \bigr\|_{X(T)} &\leq 2 C_0 \varepsilon. \label{equ:consequences_cor_assumption2}
    \end{align}
    Then the following estimates hold:
    \begin{itemize}[leftmargin=1.8em]
        
        \item[(1)] Decomposition of the radiation term:
        \begin{equation} \label{equ:consequences_decomposition_radiation}
            \bmu(t) = \ulPe \bmu(t) + d_{0,\ulell}(t) \bmY_{0,\ulell} + d_{1,\ulell}(t) \bmY_{1,\ulell}, \quad \ulPe \bmu(t) = \begin{bmatrix} \usubeone(t) \\ \usubetwo(t) \end{bmatrix}, \quad 0 \leq t \leq T,
        \end{equation}
        with 
        \begin{align}
            \sup_{0 \leq t \leq T} \, \jt^{\frac12} \Bigl( \| \usubeone(t) \|_{L^\infty_y} + \| \usubetwo(t) \|_{L^\infty_y} + \| \py \usubeone(t) \|_{L^\infty_y} \Bigr) &\lesssim \varepsilon, \label{equ:consequences_dispersive_decay_ulPe_radiation} \\
            \sup_{0 \leq t \leq T} \, \jt^{\frac32 - \delta} \bigl( |d_{0,\ulell}(t)| + |d_{1,\ulell}(t)| \bigr) &\lesssim \varepsilon. \label{equ:consequences_discrete_comp_decay}
        \end{align}    
        
        \item[(2)] Dispersive decay:
        \begin{align}
            \sup_{0 \leq t \leq T} \, \jt^{\frac12} \Bigl( \|u_1(t)\|_{L^\infty_y} + \|u_2(t)\|_{L^\infty_y} + \|\py u_1(t)\|_{L^\infty_y} \Bigr) &\lesssim \varepsilon, \label{equ:consequences_dispersive_decay_radiation} 
        \end{align}

        \item[(3)] Sobolev bounds for the radiation term:
        \begin{align}
            \sup_{0 \leq t \leq T} \, \Bigl( \|\ulPe \bmu(t)\|_{H^1_y \times L^2_y} + \| \bmu(t) \|_{H^1_y \times L^2_y} \Bigr) &\lesssim \varepsilon. \label{equ:consequences_energy_bound_radiation} \\
            \sup_{0 \leq t \leq T} \, \jt^{-\delta} \Bigl( \|\ulPe \bmu(t)\|_{H^3_y \times H^2_y} + \| \bmu(t) \|_{H^3_y \times H^2_y} \Bigr) &\lesssim \varepsilon. \label{equ:consequences_sobolev_bound_radiation}
        \end{align}

        \item[(4)] Uniform-in-time $L^2_\xi$-bound for the effective profile:
            \begin{equation} \label{equ:consequences_uniform_in_time_L2_bound} 
                \sup_{0 \leq t \leq T} \, \bigl\| \gulellsh(t,\xi) \bigr\|_{L^2_\xi} \lesssim \varepsilon.
            \end{equation}               

        \item[(5)] Auxiliary decay estimates for the modulation parameters:
        \begin{align}
             \sup_{0 \leq t \leq T} \, \jt \bigl( |\dot{\ell}(t)| + |\dot{q}(t)-\ell(t)| \bigr) &\lesssim  \varepsilon, \label{equ:consequences_crude_decay_modulation_parameters} \\ 
             \sup_{0 \leq t \leq T} \, \jt^{1-\delta} |\dot{q}(t) - \ulell| &\lesssim \varepsilon. \label{equ:consequences_qdot_minus_ulell_decay}
        \end{align}

        \item[(6)] Growth bound for the phase:
            \begin{equation} \label{equ:consequences_theta_growth_bound}
                \sup_{0 \leq t \leq T} \, \jt^{-\delta} |\theta(t)| \lesssim \varepsilon.
            \end{equation}

        \item[(7)] Local decay:
            \begin{align}
                \sup_{0 \leq t \leq T} \, \jt^{\frac12} \Bigl( \bigl\| \jy^{-1} \usubeone(t) \bigr\|_{H^3_y} + \bigl\| \jy^{-1} \usubetwo(t) \bigr\|_{H^2_y} \Bigr) &\lesssim \varepsilon, \label{equ:consequences_local_decay_usubeonetwo} \\ 
                \sup_{0 \leq t \leq T} \, \jt^{\frac12} \Bigl( \bigl\| \jy^{-1} u_1(t) \bigr\|_{H^3_y} + \bigl\| \jy^{-1} u_2(t) \bigr\|_{H^2_y} \Bigr) &\lesssim \varepsilon. \label{equ:consequences_local_decay_uonetwo}
            \end{align}
        
        \item[(8)] Auxiliary growth bound for a weighted Sobolev norm of the radiation term:
        \begin{equation} \label{equ:consequences_weighted_sobolev_strong_growth}
            \sup_{0 \leq t \leq T} \, \jt^{-1-\delta} \bigl\| \jy \usubeone(t) \bigr\|_{H^2_y} \lesssim \varepsilon.
        \end{equation}

        \item[(9)] Auxiliary bounds for the spatially localized cubic-type nonlinearities $\widetilde{\calR}(t,\xi)$:
            \begin{align}
                \sup_{0 \leq t \leq T} \, \jt^{\frac32-\delta} \bigl\| \jxi^2 \widetilde{\calR}(t,\xi) \bigr\|_{L^2_\xi} &\lesssim \varepsilon^2, \label{equ:consequences_wtilcalR_L2xi_bound} \\
                \sup_{0 \leq t \leq T} \, \jt^{1-\delta} \bigl\| \jxi^2 \pxi  \widetilde{\calR}(t,\xi) \bigr\|_{L^2_\xi} &\lesssim \varepsilon^2. \label{equ:consequences_wtilcalR_pxi_L2xi_bound}
            \end{align}            

        \item[(10)] Leading order local decay:
            \begin{align}
                \usubeone(t,y) &= \Re \Bigl( e_\ulell^{\#}\bigl(y, - \ulg \ulell\bigr) h_\ulell(t) \Bigr) + \Remusubeone(t,y), \quad \quad \quad \quad \, \, \, 0 \leq t \leq T, \label{equ:consequences_decomposition_leading_order_local_decay} \\
                \usubetwo(t,y) &= \Re \Bigl( \ulg^{-1} (D^\ast e_\ulell^{\#})\bigl(y, - \ulg \ulell\bigr) p_\ulell(t) \Bigr) + \Remusubetwo(t,y), \quad 0 \leq t \leq T, \label{equ:consequences_decomposition_leading_order_local_decay_usubetwo}
            \end{align}
        with
        \begin{align}
                \sup_{0 \leq t \leq T} \, \jt^{\frac12} \bigl| h_\ulell(t) \bigr| &\lesssim \varepsilon, \label{equ:consequences_hulell_decay} \\
                \sup_{0 \leq t \leq T} \, \jt^{1-\delta} \bigl| \pt \bigl( e^{-it\ulg^{-1}} h_\ulell(t) \bigr) \bigr| &\lesssim \varepsilon, \label{equ:consequences_hulell_phase_filtered_decay} \\
                \sup_{0 \leq t \leq T} \, \jt^{1-\delta} \bigl\| \jy^{-3} \Remusubeone(t) \bigr\|_{H^1_y} &\lesssim \varepsilon, \label{equ:consequences_Remusubeone_H1y_local_decay} \\ 
                \sup_{0 \leq t \leq T} \, \jt^{1-2\delta} \bigl\| \jy^{-3} \Remusubeone(t) \bigr\|_{H^3_y} &\lesssim \varepsilon, \label{equ:consequences_Remusubeone_H3y_local_decay}
        \end{align}  
        as well as 
        \begin{align}
            \sup_{0 \leq t \leq T} \, \jt^{\frac12} |p_\ulell(t)| &\lesssim \varepsilon, \label{equ:consequences_pulell_decay} \\ 
            \sup_{0 \leq t \leq T} \, \jt^{1-\delta} \bigl| \pt \bigl( e^{-it\ulg^{-1}} p_\ulell(t) \bigr) \bigr| &\lesssim \varepsilon, \label{equ:consequences_pulell_phase_filtered_decay} \\ 
            \sup_{0 \leq t \leq T} \, \jt^{1-\delta} \bigl\| \jy^{-3} \Remusubetwo(t,y) \bigr\|_{L^2_y} &\lesssim \varepsilon. \label{equ:consequences_Remusubetwo_L2y_local_decay}
        \end{align}

        \item[(11)] Auxiliary linear Klein-Gordon evolutions: 

        \noindent For $\frakb \in W^{1,\infty}_\xi(\bbR)$ define 
        \begin{equation} \label{equ:consequences_aux_KG_evolutions_def}
            v(t,y) := \frac{1}{\sqrt{2\pi}} \int_\bbR e^{i y \xi} e^{it(\jxi+\ulell\xi)} \frakb(\xi) \jxi^{-1} \gulellsh(t,\xi) \, \ud \xi, \quad (t,y) \in [0,T] \times \bbR.
        \end{equation}
        Then it holds that
        \begin{align}
            \sup_{0 \leq t \leq T} \, \jt^{\frac12} \bigl( \|v(t)\|_{L^\infty_y} + \|\py v(t)\|_{L^\infty_y} \bigr) &\lesssim \varepsilon, \label{equ:consequences_aux_KG_disp_decay} \\ 
            \sup_{0 \leq t \leq T} \, \bigl( \|v(t)\|_{H^1_y} + \jt^{-\delta} \|v(t)\|_{H^3_y} \bigr) &\lesssim \varepsilon, \label{equ:consequences_aux_KG_energy_bounds} \\
            \sup_{0 \leq t \leq T} \, \jt^{\frac12} \bigl\| \jy^{-1} v(t) \bigr\|_{H^3_y} &\lesssim \varepsilon. \label{equ:consequences_aux_KG_local_H3y_decay}
        \end{align}
        Moreover, we have the leading order local decay decomposition
        \begin{equation} \label{equ:consequences_aux_KG_leading_order_local_decay_decomp}
            v(t,y) = e^{-i\ulg\ulell y} \tilde{h}_\ulell(t) + R_{v}(t,y), \quad 0 \leq t \leq T,
        \end{equation}
        with 
        \begin{align}
                \sup_{0 \leq t \leq T} \, \jt^{\frac12} | \tilde{h}_\ulell(t) | &\lesssim \varepsilon, \label{equ:consequences_aux_KG_tildeh_ulell_decay} \\
                \sup_{0 \leq t \leq T} \, \jt^{1-\delta} \bigl| \pt \bigl( e^{-it\ulg^{-1}} \tilde{h}_\ulell(t) \bigr) \bigr| &\lesssim \varepsilon, \label{equ:consequences_aux_KG_tildeh_ulell_phase_filtered_decay} \\
                \sup_{0 \leq t \leq T} \, \jt^{1-\delta} \bigl\| \jy^{-3} R_{v}(t) \bigr\|_{H^1_y} &\lesssim \varepsilon. \label{equ:consequences_aux_KG_Remv_improved_H1y_local_decay}
        \end{align}  
        In the special case
        \begin{equation*}
            \frakb(\xi) = \frac{\ulg(\xi+\ulell\jxi)}{|\ulg(\xi+\ulell\jxi)| \pm i},
        \end{equation*}
        the following improved local decay estimate holds
        \begin{align}
            \sup_{0 \leq t \leq T} \, \jt^{1-\delta} \bigl\| \jy^{-1} v(t) \bigr\|_{H^1_y} &\lesssim \varepsilon.\label{equ:consequences_aux_KG_improved_local_H1y_decay} 
        \end{align}
    \end{itemize}
\end{corollary}

\begin{proof}[Proof of Corollary \ref{cor:consequences_bootstrap_assumptions}]
We establish the asserted bounds item by item.
Throughout we consider times $0 \leq t \leq T$.

\noindent \underline{Proof of (1).} 
The decomposition \eqref{equ:consequences_decomposition_radiation} of the radiation term follows from Lemma~\ref{lem:decom}.
We begin with the proof of \eqref{equ:consequences_dispersive_decay_ulPe_radiation}.
For short times $0 \leq t \leq 1$ we have by the Sobolev embedding $H^1_y(\bbR) \hookrightarrow L^\infty_y(\bbR)$, by Lemma~\ref{lem:higherorder}, and by the bootstrap assumption \eqref{equ:consequences_cor_assumption2},
\begin{equation*}
    \begin{aligned}
        \|\usubeone(t)\|_{L^\infty_y} + \|\usubetwo(t)\|_{L^\infty_y} + \|\py \usubeone(t)\|_{L^\infty_y}  
        \lesssim \|\ulPe \bmu(t)\|_{H^2_y \times H^1_y} \lesssim \bigl\| \jxi \gulellsh(t,\xi) \bigr\|_{L^2_\xi} \lesssim \varepsilon.
    \end{aligned}
\end{equation*}
For times $1 \leq t \leq T$, using the bootstrap assumption \eqref{equ:consequences_cor_assumption2}, we obtain from the dispersive decay estimate \eqref{equ:linear_evolution_dispersive_decay} that 
\begin{equation*}
    \begin{aligned}
        &\|\usubeone(t)\|_{L^\infty_y} + \|\usubetwo(t)\|_{L^\infty_y} + \|\py \usubeone(t)\|_{L^\infty_y}  \\
        &= \bigl\| \bigl( e^{t\bfL_\ulell} \ulPe \bmf_\ulell(t) \bigr)_1 \bigr\|_{L^\infty_y} + \bigl\| \bigl( e^{t\bfL_\ulell} \ulPe \bmf_\ulell(t) \bigr)_2 \bigr\|_{L^\infty_y} + \bigl\| \py \bigl( e^{t\bfL_\ulell} \ulPe \bmf_\ulell(t) \bigr)_1 \bigr\|_{L^\infty_y} \\
        &\lesssim t^{-\frac12} \bigl\| \jxi^{\frac32} \gulellsh(t,\xi) \bigr\|_{L^\infty_\xi} + t^{-\frac23} \Bigl( \bigl\| \jxi^2 \pxi \gulellsh(t,\xi) \bigr\|_{L^2_\xi} + \bigl\| \jxi^2 \gulellsh(t,\xi) \bigr\|_{L^2_\xi} \Bigr) \\ 
        &\lesssim t^{-\frac12} \cdot \varepsilon + t^{-\frac23} \cdot \varepsilon \jt^{2\delta} \lesssim \varepsilon t^{-\frac12}.
    \end{aligned}
\end{equation*}
Combining the preceding estimates establishes \eqref{equ:consequences_dispersive_decay_ulPe_radiation}.

For the proof of \eqref{equ:consequences_discrete_comp_decay} we test the decomposition \eqref{equ:consequences_decomposition_radiation} against $\bfJ \bmY_{0,\ell(t)}$, respectively against $\bfJ \bmY_{1,\ell(t)}$, and use the orthogonality properties \eqref{equ:orthogonality_radiation} of the radiation term to find that
\begin{equation*}
    \begin{aligned}
    \begin{bmatrix}
        \langle \bfJ \bmY_{0,\ell(t)}, \bmY_{0,\ulell} \rangle & \langle \bfJ \bmY_{0,\ell(t)}, \bmY_{1,\ulell} \rangle \\
        \langle \bfJ \bmY_{1,\ell(t)}, \bmY_{0,\ulell} \rangle & \langle \bfJ \bmY_{1,\ell(t)}, \bmY_{1,\ulell} \rangle
    \end{bmatrix}
    \begin{bmatrix}
        d_{0,\ulell}(t) \\
        d_{1,\ulell(t)} 
    \end{bmatrix}
    &= 
    \begin{bmatrix}
        -\langle \bfJ \bmY_{0,\ell(t)}, \ulPe \bmu(t) \rangle \\
        - \langle \bfJ \bmY_{1,\ell(t)}, \ulPe \bmu(t) \rangle 
    \end{bmatrix} 
    = 
    \begin{bmatrix}
        -\langle \bfJ (\bmY_{0,\ell(t)} - \bmY_{0,\ulell}), \ulPe \bmu(t) \rangle \\
        - \langle \bfJ (\bmY_{1,\ell(t)} - \bmY_{0,\ulell}), \ulPe \bmu(t) \rangle 
    \end{bmatrix}. 
    \end{aligned}
\end{equation*}
In view of the orthogonality properties~\eqref{equ:orthogonality_radiation} and using \eqref{equ:Jinnerproduct_generalized_kernel_elements}, we can write the matrix on the left-hand side of the preceding equation as 
\begin{equation*}
    \begin{aligned}
        &\begin{bmatrix}
            \langle \bfJ \bmY_{0,\ell(t)}, \bmY_{0,\ulell} \rangle & \langle \bfJ \bmY_{0,\ell(t)}, \bmY_{1,\ulell} \rangle \\
            \langle \bfJ \bmY_{1,\ell(t)}, \bmY_{0,\ulell} \rangle & \langle \bfJ \bmY_{1,\ell(t)}, \bmY_{1,\ulell} \rangle
        \end{bmatrix} \\ 
        &=         
        \begin{bmatrix}
            0 & \langle \bfJ \bmY_{0,\ulell}, \bmY_{1,\ulell} \rangle \\
            \langle \bfJ \bmY_{1,\ulell}, \bmY_{0,\ulell} \rangle & 0
        \end{bmatrix} 
        +
        \begin{bmatrix}
            \langle \bfJ (\bmY_{0,\ell(t)} - \bmY_{0,\ulell}), \bmY_{0,\ulell} \rangle & \langle \bfJ (\bmY_{0,\ell(t)} - \bmY_{0,\ulell}), \bmY_{1,\ulell} \rangle \\
            \langle \bfJ (\bmY_{1,\ell(t)} - \bmY_{1,\ulell}), \bmY_{0,\ulell} \rangle & \langle \bfJ (\bmY_{1,\ell(t)} - \bmY_{1,\ulell}), \bmY_{1,\ulell} \rangle
        \end{bmatrix} \\         
        &=
        \ulg^3 \|K'\|_{L^2_y}^2 \begin{bmatrix} 0 & 1 \\ -1 & 0 \end{bmatrix} + \calO_{\ell_0} \bigl( |\ell(t) - \ulell| \bigr) \begin{bmatrix} 1 & 1 \\ 1 & 1 \end{bmatrix}.
    \end{aligned}
\end{equation*}
Since $|\ell(t)-\ulell| \lesssim \varepsilon$ by \eqref{equ:consequences_cor_assumption1} and since the assumption $|\ulell-\ell_0| \leq \frac14 \gamma_0^{-2}$ implies $\frac12 \gamma_0 \leq \ulg \leq 2\gamma_0$ by \eqref{equ:setting_up_comparison_gammas}, it follows that the matrix on the left-hand side is invertible and that there is a uniform upper bound on the operator norm of its inverse whose size depends only on $\ell_0$. 
Hence, using assumption \eqref{equ:consequences_cor_assumption1} and the dispersive decay estimates \eqref{equ:consequences_dispersive_decay_ulPe_radiation}, we obtain for $0 \leq t \leq T$ the asserted decay estimates \eqref{equ:consequences_discrete_comp_decay} for the discrete components,
\begin{equation*}
    \begin{aligned}
        |d_{0,\ulell}(t)| + |d_{1,\ulell}(t)|  
        &\lesssim_{\ell_0} \bigl|\langle \bfJ (\bmY_{0,\ell(t)} - \bmY_{0,\ulell}), \ulPe \bmu(t) \rangle\bigr| + \bigl|\langle \bfJ (\bmY_{1,\ell(t)} - \bmY_{0,\ulell}), \ulPe \bmu(t) \rangle\bigr| \\ 
        &\lesssim_{\ell_0} \Bigl( \bigl\|\bfJ (\bmY_{0,\ell(t)} - \bmY_{0,\ulell}) \bigr\|_{L^1_y} + \bigl\| \bfJ (\bmY_{1,\ell(t)} - \bmY_{1,\ulell}) \bigr\|_{L^1_y} \Bigr) \bigl\| \ulPe \bmu(t) \bigr\|_{L^\infty_y} \\
        &\lesssim_{\ell_0} |\ell(t)-\ulell| \bigl\| \ulPe \bmu(t) \bigr\|_{L^\infty_y} \\
        &\lesssim_{\ell_0} \varepsilon \jt^{-1+\delta} \cdot \varepsilon \jt^{-\frac12} \lesssim_{\ell_0} \varepsilon^2 \, \jt^{-\frac32 + \delta}.
    \end{aligned}
\end{equation*}

\noindent \underline{Proof of (2).} 
The dispersive decay estimates \eqref{equ:consequences_dispersive_decay_radiation} follow immediately from the decomposition \eqref{equ:consequences_decomposition_radiation} of the radiation term and the decay estimates \eqref{equ:consequences_dispersive_decay_ulPe_radiation}, \eqref{equ:consequences_discrete_comp_decay} established in the preceding step.

\noindent \underline{Proof of (3).} 
The uniform-in-time energy bound \eqref{equ:consequences_energy_bound_radiation} for the radiation term follows from the stability bound~\eqref{equ:smallness_orbital_proposition} and the boundedness of the projection $\ulPe$ in $H^1_y \times L^2_y$ by Lemma~\ref{lem:decom}.
Next, using \eqref{equ:consequences_cor_assumption2} we conclude from Lemma~\ref{lem:higherorder} that
\begin{equation*}
    \|\ulPe \bmu(t)\|_{H^3_y \times H^2_y} \lesssim \bigl\| \jxi^2 \gulellsh(t,\xi) \bigr\|_{L^2_\xi} \lesssim \varepsilon \jt^\delta.
\end{equation*}
Then the decomposition \eqref{equ:consequences_decomposition_radiation} and \eqref{equ:consequences_discrete_comp_decay} imply the asserted estimate
\begin{equation*}
    \begin{aligned}
    \|\bmu(t)\|_{H^3_y \times H^2_y} \lesssim \bigl\|\ulPe \bmu(t)\bigr\|_{H^3_y \times H^2_y} + |d_{0,\ulell}(t)| \bigl\|\bmY_{0,\ulell}\bigr\|_{H^3_y \times H^2_y} + |d_{1,\ulell}(t)| \bigl\|\bmY_{1,\ulell}\bigr\|_{H^3_y \times H^2_y} 
    &\lesssim \varepsilon \jt^\delta + \varepsilon \lesssim \varepsilon \jt^\delta.
    \end{aligned}
\end{equation*}

\noindent \underline{Proof of (4).} 
The uniform-in-time $L^2_\xi$-bound \eqref{equ:consequences_uniform_in_time_L2_bound} follows from 
Lemma~\ref{lem:higherorder}, the boundedness of the projection $\ulPe$ in $H^1_y \times L^2_y$ by Lemma~\ref{lem:decom}, and the stability bound \eqref{equ:smallness_orbital_proposition},
\begin{equation*}
    \bigl\| \gulellsh(t,\xi) \bigr\|_{L^2_\xi} \lesssim \|\ulPe \bmu(t)\|_{H^1_y \times L^2_y} \lesssim \|\bmu(t)\|_{H^1_y \times L^2_y} \lesssim \varepsilon.
\end{equation*}

\noindent \underline{Proof of (5).}
It follows from \eqref{equ:smallness_orbital_proposition} that the matrix \eqref{equ:modulation_equ_matrix_lhs} on the left-hand side of the modulation equations \eqref{equ:modulation_equ} is invertible with a uniform bound on its operator norm that depends only on~$\ell_0$. Hence, the decay estimates \eqref{equ:consequences_crude_decay_modulation_parameters} for the modulation parameters follow directly from the modulation equations \eqref{equ:modulation_equ} and the dispersive decay estimate~\eqref{equ:consequences_dispersive_decay_radiation} for the radiation,
\begin{equation*}
    \begin{aligned}
        |\dot{\ell}(t)| + |\dot{q}(t)-\ell(t)| &\lesssim \bigl\| \bfM_{\ell(t)}[\bmu(t)]^{-1} \bigr\| \Bigl( \bigl| \bigl\langle \bfJ \partial_q \bmK_{\ell,q}, \calN\bigl( \bmu(t) \bigr) \bigr\rangle \bigr| + \bigl| \bigl\langle \bfJ \partial_\ell \bmK_{\ell,q}, \calN\bigl( \bmu(t) \bigr) \bigr\rangle \bigr| \Bigr) \\ 
        &\lesssim_{\ell_0} \|u_1(t)\|_{L^\infty_y}^2 \lesssim \varepsilon^2 \jt^{-1}.
    \end{aligned}
\end{equation*}
Then the estimate \eqref{equ:consequences_qdot_minus_ulell_decay} is an immediate consequence of the assumption \eqref{equ:consequences_cor_assumption1} and \eqref{equ:consequences_crude_decay_modulation_parameters},
\begin{equation*}
    \begin{aligned}
        |\dot{q}(t)-\ulell| \leq |\dot{q}(t)-\ell(t)| + |\ell(t)-\ulell| \lesssim \varepsilon \jt^{-1} +  \varepsilon \jt^{-1+\delta} \lesssim \varepsilon \jt^{-1+\delta}.
    \end{aligned}
\end{equation*}

\noindent \underline{Proof of (6).} 
The decay estimate \eqref{equ:consequences_qdot_minus_ulell_decay} implies straightaway
\begin{equation*}
    \begin{aligned}
        |\theta(t)| \leq  \int_0^t |\dot{q}(s) - \ulell| \, \ud s \lesssim \int_0^t \varepsilon \js^{-1+\delta} \, \ud s \lesssim \varepsilon \jt^{\delta}.
    \end{aligned}
\end{equation*}

\noindent \underline{Proof of (7).} 
It suffices to establish \eqref{equ:consequences_local_decay_usubeonetwo}. Then \eqref{equ:consequences_local_decay_uonetwo} follows from \eqref{equ:consequences_local_decay_usubeonetwo}, the decomposition of the radiation term \eqref{equ:consequences_decomposition_radiation}, and the faster decay of the discrete components \eqref{equ:consequences_discrete_comp_decay}.
In order to prove \eqref{equ:consequences_local_decay_usubeonetwo}, we decompose $\usubeone(t,y)$ and $\usubetwo(t,y)$ into a low frequency and a high frequency component. 
Denoting by $\chi_{\{\leq 2\ulg|\ulell|\}}(\xi)$ a smooth even bump function satisfying $\chi_{\{\leq 2\ulg|\ulell|\}}(\xi) = 1$ for $|\xi| \leq 2\ulg|\ulell|$ and $\chi_{\{\leq 2\ulg|\ulell|\}}(\xi) = 0$ for $|\xi| \geq 4\ulg|\ulell|$, we write
\begin{equation*}
    \begin{aligned}
        \usubeone(t,y) &= \usubeonemusFlat(t,y) + \usubeonemusSharp(t,y), \\
        \usubetwo(t,y) &= \usubetwomusFlat(t,y) + \usubetwomusSharp(t,y), 
    \end{aligned}
\end{equation*}
with 
\begin{align*}
    \usubeonemusFlat(t,y) &:= \Re \, \Bigl( \calFulellshast \Bigl[ e^{i t (\jxi + \ulell \xi)} \chi_{\{\leq 2\ulg|\ulell|\}}(\xi) i\jxi^{-1} \gulellsh(t,\xi) \Bigr](y) \Bigr), \\
    \usubeonemusSharp(t,y) &:= \Re \, \Bigl( \calFulellshast \Bigl[ e^{i t (\jxi + \ulell \xi)} \bigl( 1 - \chi_{\{\leq 2\ulg|\ulell|\}}(\xi) \bigr) i\jxi^{-1} \gulellsh(t,\xi) \Bigr](y) \Bigr), \\
    \usubetwomusFlat(t,y) &:= \Re \, \Bigl( \calFulellDshast \Bigl[ e^{i t (\jxi + \ulell \xi)} \chi_{\{\leq 2\ulg|\ulell|\}}(\xi) i\jxi^{-1} \gulellsh(t,\xi) \Bigr](y) \Bigr), \\
    \usubetwomusSharp(t,y) &:= \Re \, \Bigl( \calFulellDshast \Bigl[ e^{i t (\jxi + \ulell \xi)} \bigl( 1 - \chi_{\{\leq 2\ulg|\ulell|\}}(\xi) \bigr) i\jxi^{-1} \gulellsh(t,\xi) \Bigr](y) \Bigr).    
\end{align*}
In view of \eqref{equ:consequences_sobolev_bound_radiation}, it suffices to consider times $t \geq 1$. 
Thanks to the bounds on the distorted Fourier basis elements from Lemma~\ref{lem:boundseDe}, we can invoke the dispersive decay estimate from Lemma~\ref{lem:core_linear_dispersive_decay}, and then conclude from the finite support of the frequency cut-off $\chi_{\{\leq 2\ulg|\ulell|\}}(\xi)$ and the bootstrap assumptions \eqref{equ:consequences_cor_assumption2} that 
\begin{equation*}
    \begin{aligned}
        \bigl\| \usubeonemusFlat(t) \bigr\|_{W^{3,\infty}_y} &\lesssim t^{-\frac12} \bigl\| \chi_{\{\leq 2\ulg|\ulell|\}}(\xi) \jxi^{\frac72} \gulellsh(t,\xi) \bigr\|_{L^\infty_\xi} \\
        &\quad + t^{-\frac23} \Bigl( \bigl\| \jxi^4 \pxi \bigl( \chi_{\{\leq 2\ulg|\ulell|\}}(\xi) \gulellsh(t,\xi) \bigr) \bigr\|_{L^2_\xi}  + \bigl\| \jxi^4 \chi_{\{\leq 2\ulg|\ulell|\}}(\xi) \gulellsh(t,\xi) \bigr\|_{L^2_\xi} \Bigr)  \\
        &\lesssim_{\ell_0} t^{-\frac12} \bigl\| \gulellsh(t,\xi) \bigr\|_{L^\infty_\xi} + t^{-\frac23} \Bigl( \bigl\|\gulellsh(t,\xi) \bigr) \bigr\|_{L^2_\xi}  + \bigl\| \gulellsh(t,\xi) \bigr\|_{L^2_\xi} \Bigr)  \\
        &\lesssim t^{-\frac12} \cdot \varepsilon + t^{-\frac23} \cdot \varepsilon \jt^\delta \lesssim \varepsilon t^{-\frac12}.
    \end{aligned}
\end{equation*}
Moreover, the improved local decay estimate \eqref{equ:improved_local_decay_high_freq} with $0 \leq k \leq 3$ and the bootstrap assumptions \eqref{equ:consequences_cor_assumption2} imply 
\begin{equation*}
    \bigl\| \jy^{-1} \usubeonemusSharp(t) \bigr\|_{H^3_y} \lesssim \jt^{-1} \Bigl( \bigl\| \jxi^2 \pxi \gulellsh(t,\xi) \bigr\|_{L^2_\xi} + \bigl\| \jxi^2 \gulellsh(t,\xi) \bigr\|_{L^2_\xi} \Bigr) \lesssim \jt^{-1} \cdot \varepsilon \jt^{2\delta} \lesssim \varepsilon t^{-1+2\delta},
\end{equation*}
whence
\begin{equation*}
    \begin{aligned}
        \bigl\| \jy^{-1} \usubeone(t) \bigr\|_{H^3_y} &\lesssim \bigl\| \usubeonemusFlat(t) \bigr\|_{W^{3,\infty}_y} + \bigl\| \jy^{-1} \usubeonemusSharp(t) \bigr\|_{H^3_y} \lesssim \varepsilon t^{-\frac12}.
    \end{aligned}
\end{equation*}
This gives \eqref{equ:consequences_local_decay_usubeonetwo} for $\usubeone(t,y)$.

Proceeding analogously for $\usubetwo(t,y)$, we infer from Lemma~\ref{lem:core_linear_dispersive_decay} using the bounds from Lemma~\ref{lem:boundseDe}, the finite support of the frequency cut-off $\chi_{\{\leq 2\ulg\ulell\}}(\xi)$, the improved local decay estimate \eqref{equ:improved_local_decay_high_freq_for_2nd_comp} with $0 \leq k \leq 2$, and the bootstrap assumptions \eqref{equ:consequences_cor_assumption2} that
\begin{equation*}
    \begin{aligned}
        \bigl\| \jy^{-1} \usubetwo(t) \bigr\|_{H^2_y} &\lesssim \bigl\| \usubetwomusFlat(t) \bigr\|_{W^{2,\infty}_y} + \bigl\| \jy^{-1} \usubetwomusSharp(t) \bigr\|_{H^2_y} \lesssim \varepsilon t^{-\frac12}.
    \end{aligned}
\end{equation*}
This proves \eqref{equ:consequences_local_decay_usubeonetwo} for $\usubetwo(t,y)$.

\noindent \underline{Proof of (8).} 
From the representation formula \eqref{equ:setting_up_representation_formula_usubeone} for $\usubeone(t,y)$,
the mapping property \eqref{equ:mapping_property_calFulellshast_jy_H2y} for $\calFulellshast$ from Lemma~\ref{lem:mapping_properties_calF},
and the bootstrap assumptions \eqref{equ:consequences_cor_assumption2}, we infer the asserted growth bound
\begin{equation*}
    \begin{aligned}
        \bigl\| \jy \usubeone(t,y) \bigr\|_{H^2_y} &\lesssim \Bigl\| \jxi^2 \pxi \Bigl( e^{it(\jxi+\ulell\xi)} i\jxi^{-1} \gulellsh(t,\xi) \Bigr) \Bigr\|_{L^2_\xi} + \Bigl\| \jxi^2 \Bigl( e^{it(\jxi+\ulell\xi)} i\jxi^{-1} \gulellsh(t,\xi) \Bigr) \Bigr\|_{L^2_\xi} \\ 
        &\lesssim \jt \bigl\| \jxi \gulellsh(t,\xi) \bigr\|_{L^2_\xi} + \bigl\| \jxi \pxi \gulellsh(t,\xi) \bigr\|_{L^2_\xi} \\ 
        &\lesssim \jt \cdot \varepsilon \jt^\delta + \varepsilon \jt^{2\delta} \lesssim \varepsilon \jt^{1+\delta}.
    \end{aligned}
\end{equation*}

\noindent \underline{Proof of (9).}
    We separately consider the contributions of each term in the definition \eqref{equ:setting_up_definition_wtilR} of $\widetilde{\calR}(t,\xi)$ to the bounds \eqref{equ:consequences_wtilcalR_L2xi_bound} and \eqref{equ:consequences_wtilcalR_pxi_L2xi_bound}. 
    Apart from the quintic and higher order nonlinearities in $\calR_2\bigl(\usubeone(t)\bigr)$, the corresponding estimates are largely identical for all other terms in $\widetilde{\calR}(t,\xi)$ due to their spatial localization and their cubic-type decay.

    \noindent {\it Contribution of $(\dot{q}(t) - \ulell) \calL_\ulell\bigl(\bmu(t)\bigr)(\xi)$}:
    Using \eqref{equ:consequences_qdot_minus_ulell_decay}, \eqref{equ:consequences_local_decay_uonetwo}, 
    and the mapping property \eqref{eq:calLellL2}, we obtain
    \begin{equation*}
        \begin{aligned}
            \bigl\| \jxi^2 (\dot{q}(t) - \ulell) \calL_\ulell\bigl(\bmu(t)\bigr)(\xi) \bigr\|_{L^2_\xi} 
            &\lesssim |\dot{q}(t)-\ulell| \Bigl( \bigl\| \sech^2(\ulg y) u_1(t,y) \bigr\|_{H^2_y} + \bigl\| \sech^2(\ulg y) u_2(t,y) \bigr\|_{H^1_y} \Bigr) \\
            &\lesssim |\dot{q}(t)-\ulell| \Bigl( \bigl\| \jy^{-1} u_1(t) \bigr\|_{H^2_y} + \bigl\| \jy^{-1} u_2(t) \bigr\|_{H^1_y} \Bigr) \\ 
            &\lesssim \varepsilon \jt^{-1+\delta} \cdot \varepsilon \jt^{-\frac12} \lesssim \varepsilon^2 \jt^{-\frac32+\delta}.
        \end{aligned}
    \end{equation*}
    Analogously, we obtain by the mapping property \eqref{eq:pxicalLellL2} that
    \begin{equation*}
        \begin{aligned}
            \bigl\| \jxi^2 \pxi \bigl( (\dot{q}(t) - \ulell) \calL_\ulell\bigl(\bmu(t)\bigr)(\xi) \bigr) \bigr\|_{L^2_\xi} 
            &\lesssim |\dot{q}(t)-\ulell| \Bigl( \bigl\| \jy \sech^2(\ulg y) u_1(t,y) \bigr\|_{H^2_y} + \bigl\| \jy \sech^2(\ulg y) u_2(t,y) \bigr\|_{H^1_y} \Bigr) \\
            &\lesssim |\dot{q}(t)-\ulell| \Bigl( \bigl\| \jy^{-1} u_1(t) \bigr\|_{H^2_y} + \bigl\| \jy^{-1} u_2(t) \bigr\|_{H^1_y} \Bigr) \\ 
            &\lesssim \varepsilon \jt^{-1+\delta} \cdot \varepsilon \jt^{-\frac12} \lesssim \varepsilon^2 \jt^{-\frac32+\delta}.
        \end{aligned}
    \end{equation*}

    \noindent {\it Contribution of $\calMod(t)$}: 
    By Corollary~\ref{cor:TellP} we have
    \begin{equation*}
        \begin{aligned}
            \calF_{\ulell,D}^{\#}\bigl[ \bigl(\calMod\bigr)_1(t) \bigr](\xi) - \calF_{\ulell}^{\#}\bigl[ \bigl(\calMod\bigr)_2(t) \bigr](\xi) = \calF_{\ulell,D}^{\#}\bigl[ \bigl(\ulPe \calMod\bigr)_1(t) \bigr](\xi) - \calF_{\ulell}^{\#}\bigl[ \bigl(\ulPe \calMod\bigr)_2(t) \bigr](\xi).
        \end{aligned}
    \end{equation*}
    Moreover, since $\ulPe \bmY_{j,\ulell} = 0$ for $j=0,1$, we can write
    \begin{equation*}
        \begin{aligned}
            \ulPe \calMod(t,y) = - \bigl(\dot{q}(t)-\ell(t)\bigr) \ulPe \bigl( \bmY_{0,\ell(t)}(y) - \bmY_{0,\ulell}(y) \bigr) - \dot{\ell}(t) \ulPe \bigl( \bmY_{1,\ell(t)}(y) - \bmY_{1,\ulell}(y) \bigr).
        \end{aligned}
    \end{equation*}
    By the preceding identities, the mapping properties \eqref{equ:mapping_property_calFulellsh_jxi2}, \eqref{equ:mapping_property_calFulellsh_jxi2_pxi}, \eqref{equ:mapping_property_calFulellDsh_jxi2}, \eqref{equ:mapping_property_calFulellDsh_jxi2_pxi} from Lemma~\ref{lem:mapping_properties_calF},
    and the decay estimates for the modulation parameters \eqref{equ:consequences_cor_assumption1}, \eqref{equ:consequences_crude_decay_modulation_parameters}, we obtain that
    \begin{equation*}
        \begin{aligned}
            &\Bigl\| \jxi^2 \Bigl( \calF_{\ulell,D}^{\#}\bigl[ \bigl(\calMod\bigr)_1(t) \bigr](\xi) - \calF_{\ulell}^{\#}\bigl[ \bigl(\calMod\bigr)_2(t) \bigr](\xi) \Bigr) \Bigr\|_{L^2_\xi} \\
            &\quad + \Bigl\| \jxi^2 \pxi \Bigl( \calF_{\ulell,D}^{\#}\bigl[ \bigl(\calMod\bigr)_1(t) \bigr](\xi) - \calF_{\ulell}^{\#}\bigl[ \bigl(\calMod\bigr)_2(t) \bigr](\xi) \Bigr) \Bigr\|_{L^2_\xi} \\
            &\lesssim |\dot{q}(t) - \ell(t)| \bigl\| \jy \ulPe \bigl( \bmY_{0,\ell(t)} - \bmY_{0,\ulell} \bigr) \bigr\|_{H^3_y} + |\dot{\ell}(t)| \bigl\| \jy \ulPe \bigl( \bmY_{1,\ell(t)} - \bmY_{1,\ulell} \bigr) \bigr\|_{H^3_y} \\ 
            &\lesssim \bigl( |\dot{q}(t) - \ell(t)| + |\dot{\ell}(t)| \bigr) |\ell(t)-\ulell| \\ 
            &\lesssim \varepsilon \jt^{-1} \cdot \varepsilon \jt^{-1+\delta} \lesssim \varepsilon^2 \jt^{-2+\delta}.
        \end{aligned}
    \end{equation*}

    \noindent {\it Contribution of $\calC_{\ell(t)}\bigl( \usubeone(t) \bigr)$}:    
    By the mapping properties \eqref{equ:mapping_property_calFulellsh_jxi2}, \eqref{equ:mapping_property_calFulellsh_jxi2_pxi} and the bounds  \eqref{equ:consequences_dispersive_decay_ulPe_radiation}, \eqref{equ:consequences_local_decay_usubeonetwo}, we obtain
    \begin{equation*}
        \begin{aligned}
            &\bigl\| \jxi^2 \calF_{\ulell}^{\#}\bigl[ \calC_{\ell(t)}\bigl( \usubeone(t) \bigr) \bigr](\xi) \bigr\|_{L^2_\xi} + \bigl\| \jxi^2 \pxi \calF_{\ulell}^{\#}\bigl[ \calC_{\ell(t)}\bigl( \usubeone(t) \bigr) \bigr](\xi) \bigr\|_{L^2_\xi} \\
            &\lesssim \bigl\| \jy \sech^2\bigl(\gamma(t) y\bigr) \usubeone(t)^3 \bigr\|_{H^2_y} \\ 
            &\lesssim_{\ell_0} \bigl\| \jy^{-1} \usubeone(t) \bigr\|_{H^2_y} \|\usubeone(t)\|_{W^{1,\infty}_y}^2 \\ 
            &\lesssim_{\ell_0} \varepsilon \jt^{-\frac12} \cdot \varepsilon^2 \jt^{-1} \lesssim_{\ell_0} \varepsilon^3 \jt^{-\frac32}.
        \end{aligned}
    \end{equation*}

    \noindent {\it Contribution of $\calR_1\bigl(\usubeone(t)\bigr)$}:
    Similarly to the preceding treatment of the contributions of the localized cubic terms $\calC_{\ell(t)}\bigl( \usubeone(t) \bigr)$, using \eqref{equ:mapping_property_calFulellsh_jxi2}, \eqref{equ:mapping_property_calFulellsh_jxi2_pxi}, \eqref{equ:consequences_dispersive_decay_ulPe_radiation}, \eqref{equ:consequences_local_decay_usubeonetwo}, we find 
    \begin{equation*}
        \begin{aligned}
            &\bigl\| \jxi^2 \calF_{\ulell}^{\#}\bigl[ \calR_1\bigl(\usubeone(t)\bigr) \bigr](\xi) \bigr\|_{L^2_\xi} + \bigl\| \jxi^2 \pxi \calF_{\ulell}^{\#}\bigl[ \calR_1\bigl(\usubeone(t)\bigr) \bigr](\xi) \bigr\|_{L^2_\xi} \\
            &\lesssim \bigl\| \jy \sech^2\bigl(\gamma(t) y\bigr) \tanh\bigl(\gamma(t) y\bigr) \usubeone(t)^4 \bigr\|_{H^2_y} \\
            &\lesssim_{\ell_0} \bigl\| \jy^{-1} \usubeone(t) \bigr\|_{H^2_y} \|\usubeone(t)\|_{W^{1,\infty}_y}^3 \\ 
            &\lesssim_{\ell_0} \varepsilon \jt^{-\frac12} \cdot \varepsilon^3 \jt^{-\frac32} \lesssim_{\ell_0} \varepsilon^4 \jt^{-2}.
        \end{aligned}
    \end{equation*}

    \noindent {\it Contribution of $\calR_2\bigl(\usubeone(t)\bigr)$}:    
    Since the quintic and higher nonlinearities in $\calR_2\bigl(\usubeone(t)\bigr)$ are not spatially localized, the proofs of \eqref{equ:energy_preparations_jxi2_L2_calR} and \eqref{equ:energy_preparations_jxi2_pxi_L2_calR} are different from the preceding ones. 
    First, using the mapping property \eqref{equ:mapping_property_calFulellsh_jxi2} and the bounds \eqref{equ:consequences_dispersive_decay_ulPe_radiation}, \eqref{equ:consequences_sobolev_bound_radiation}, we obtain
    \begin{equation*}
        \begin{aligned}
            \bigl\| \jxi^2 \calF_{\ulell}^{\#}\bigl[ \calR_2\bigl(\usubeone(t)\bigr) \bigr](\xi) \bigr\|_{L^2_\xi} 
            &\lesssim \biggl\| \biggl( \int_0^1 (1-r)^4 \cos\bigl(K(\gamma(t)y) + r \usubeone(t)\bigr) \, \ud r \biggr) \usubeone(t)^5 \biggr\|_{H^2_y} \\
            &\lesssim_{\ell_0} \|\usubeone(t)\|_{H^2_y} \|\usubeone(t)\|_{W^{1,\infty}_y}^4 \Bigl( 1 + \|\usubeone(t)\|_{W^{1,\infty}_y}^2 \Bigr)  \\ 
            &\lesssim_{\ell_0} \varepsilon \jt^\delta \cdot \varepsilon^4 \jt^{-2} \lesssim_{\ell_0} \varepsilon^5 \jt^{-2 + \delta}.
        \end{aligned}
    \end{equation*}
    Then, using the mapping property \eqref{equ:mapping_property_calFulellsh_jxi2_pxi} along with \eqref{equ:consequences_dispersive_decay_ulPe_radiation}, \eqref{equ:consequences_weighted_sobolev_strong_growth}, we find 
    \begin{equation*}
        \begin{aligned}
            \bigl\| \jxi^2 \pxi \calF_{\ulell}^{\#}\bigl[ \calR_2\bigl(\usubeone(t)\bigr) \bigr](\xi) \bigr\|_{L^2_\xi} 
            &\lesssim \biggl\| \jy \biggl( \int_0^1 (1-r)^4 \cos\bigl(K(\gamma(t)y) + r \usubeone(t)\bigr) \, \ud r \biggr) \usubeone(t)^5 \biggr\|_{H^2_y} \\
            &\lesssim_{\ell_0} \|\jy \usubeone(t)\|_{H^2_y} \|\usubeone(t)\|_{W^{1,\infty}_y}^4 \Bigl( 1 + \|\usubeone(t)\|_{W^{1,\infty}_y}^2 \Bigr)  \\ 
            &\lesssim_{\ell_0} \varepsilon \jt^{1+\delta} \cdot \varepsilon^4 \jt^{-2} \lesssim_{\ell_0} \varepsilon^5 \jt^{-1 + \delta}.
        \end{aligned}
    \end{equation*}

    \noindent {\it Contribution of $\calE_j(t)$, $1 \leq j \leq 3$}:
    For the contribution of $\calE_1(t)$ we have by the mapping properties \eqref{equ:mapping_property_calFulellsh_jxi2}, \eqref{equ:mapping_property_calFulellsh_jxi2_pxi} and by \eqref{equ:consequences_cor_assumption1}, \eqref{equ:consequences_local_decay_uonetwo} that
    \begin{equation*}
        \begin{aligned}
            &\bigl\| \jxi^2 \calFulellsh\bigl[ \bigl( \calE_1(t) \bigr)_2 \bigr](\xi) \bigr\|_{L^2_\xi} + \bigl\| \jxi^2 \pxi \calFulellsh\bigl[ \bigl( \calE_1(t) \bigr)_2 \bigr](\xi) \bigr\|_{L^2_\xi} \\
            &\lesssim \bigl\| \jy \bigl( \sech^2(\gamma(t) y) - \sech^2(\ulg y) \bigr) u_1(t,y) \bigr\|_{H^2_y} \\ 
            &\lesssim |\ell(t) - \ulell| \bigl\| \jy^{-1} u_1(t,y) \bigr\|_{H^2_y} \lesssim \varepsilon \jt^{-1+\delta} \cdot \varepsilon \jt^{-\frac12} \lesssim \varepsilon^2 \jt^{-\frac32+\delta}.
        \end{aligned}
    \end{equation*}
    Next, we observe that $\calE_2(t) = \calN(\bmu(t)) - \calN( \ulPe u(t) )$ is a linear combination of product terms, where at least one input is given by a discrete component $d_{k,\ulell}(t) \bmY_{k,\ulell}$, $0 \leq k \leq 1$, which decays faster and provides spatial localization. Correspondingly, we have by the mapping properties \eqref{equ:mapping_property_calFulellsh_jxi2}, \eqref{equ:mapping_property_calFulellsh_jxi2_pxi} and by  \eqref{equ:consequences_discrete_comp_decay}, \eqref{equ:consequences_local_decay_uonetwo} that
    \begin{equation*}
    \begin{aligned}
        &\bigl\| \jxi^2 \calFulellsh\bigl[ \bigl( \calE_2(t) \bigr)_2 \bigr](\xi) \bigr\|_{L^2_\xi} + \bigl\| \jxi^2 \pxi \calFulellsh\bigl[ \bigl( \calE_2(t) \bigr)_2 \bigr](\xi) \bigr\|_{L^2_\xi} \\
        &\lesssim \bigl( |d_{0,\ulell}(t)| + |d_{1,\ulell}(t)| \bigr) \Bigl( |d_{0,\ulell}(t)| + |d_{1,\ulell}(t)| + \bigl\|\jy^{-1} u_1(t)\bigr\|_{H^2_y} \Bigr) \\ 
        &\lesssim \varepsilon \jt^{-\frac32+\delta} \cdot \varepsilon \jt^{-\frac12} \lesssim \varepsilon^2 \jt^{-2+\delta}.
    \end{aligned}
    \end{equation*}
    Similarly, we estimate the contribution of $\calE_3(t)$ using \eqref{equ:consequences_cor_assumption1}, \eqref{equ:consequences_local_decay_usubeonetwo}, \eqref{equ:consequences_dispersive_decay_ulPe_radiation} by
    \begin{equation*}
    \begin{aligned}
        &\bigl\| \jxi^2 \calFulellsh\bigl[ \calE_3(t) \bigr](\xi) \bigr\|_{L^2_\xi} + \bigl\| \jxi^2 \pxi \calFulellsh\bigl[ \calE_3(t) \bigr](\xi) \bigr\|_{L^2_\xi} \\
        &\lesssim \Bigl\| \jy \bigl( \alpha\bigl(\gamma(t) y\bigr) - \alpha(\ulg y) \bigr) \bigl( \usubeone(t,y) \bigr)^2 \Bigr\|_{H^2_y} \\ 
        &\lesssim |\ell(t) - \ulell| \bigl\| \jy^{-1} \usubeone(t) \bigr\|_{H^2_y} \|\usubeone(t)\|_{W^{1,\infty}_y} \\
        &\lesssim \varepsilon \jt^{-1+\delta} \cdot \varepsilon \jt^{-\frac12} \cdot \varepsilon \jt^{-\frac12} \lesssim \varepsilon^3 \jt^{-2+\delta}.
    \end{aligned}
    \end{equation*} 
Combining all of the preceding estimates proves the asserted bounds \eqref{equ:consequences_wtilcalR_L2xi_bound} and \eqref{equ:consequences_wtilcalR_pxi_L2xi_bound}.

\noindent \underline{Proof of (10).} 
In what follows we freely assume $t \geq 1$, and we leave the easier case of short times $0 \leq t \leq 1$ to the reader.
First, we establish the improved local decay estimates \eqref{equ:consequences_Remusubeone_H1y_local_decay}, \eqref{equ:consequences_Remusubeone_H3y_local_decay} for the remainder term $\Remusubeone(t,y)$.
In view of the representation formula \eqref{equ:setting_up_representation_formula_usubeone} for $\usubeone(t,y)$ and the definition \eqref{equ:setting_up_hulell_definition} of $h_\ulell(t)$, we find that $\Remusubeone(t,y)$ is given by 
\begin{equation} \label{equ:consequences_proof_Remusubeone_expression}
    \begin{aligned}
        \Remusubeone(t,y) &= \Re \, \Bigl( \calFulellshast \Bigl[ e^{i t (\jxi + \ulell \xi)} \bigl( 1 - \chi_{\{\leq 2\ulg\ulell\}}(\xi) \bigr) i\jxi^{-1} \gulellsh(t,\xi) \Bigr](y) \Bigr) \\ 
        &\quad + \Re \, \biggl( \int_\bbR \Bigl( \eulsharp(y,\xi) - \eulsharp(y,-\ulg \ulell) \Bigr) e^{i t (\jxi + \ulell \xi)} \chi_{\{\leq 2\ulg\ulell\}}(\xi) i\jxi^{-1} \gulellsh(t,\xi) \, \ud \xi \biggr).
    \end{aligned}
\end{equation}
Thus, using the bootstrap assumptions \eqref{equ:consequences_cor_assumption2}, we deduce the improved local decay estimate \eqref{equ:consequences_Remusubeone_H1y_local_decay} from \eqref{equ:improved_local_decay_high_freq} with $0 \leq k \leq 1$ and from \eqref{equ:local_decay_resonance_subtracted_off}.
Analogously, we conclude the improved local decay estimate at higher regularity \eqref{equ:consequences_Remusubeone_H3y_local_decay}  from \eqref{equ:improved_local_decay_high_freq} with $0 \leq k \leq 3$ and from \eqref{equ:local_decay_resonance_subtracted_off}.

Similar to the definition of $h_\ulell(t)$ in \eqref{equ:setting_up_hulell_definition}, we set 
    \begin{equation*}
        p_\ulell(t) := \int_\bbR e^{it(\jxi+\ulell \xi)} \chi_{\{\leq 2\ulg|\ulell|\}}(\xi) \, i \gulellsh(t,\xi) \, \ud \xi. 
    \end{equation*}
Then the remainder term $\Remusubetwo(t,y)$ in the decomposition \eqref{equ:consequences_decomposition_leading_order_local_decay_usubetwo} of $\usubetwo(t,y)$ is given by
        \begin{equation*}
            \begin{aligned}
                \Remusubetwo(t,y) 
                &=  \Re \, \biggl( \int_\bbR \jxi^{-1} (D^\ast \eulsharp)(y,\xi) e^{i t (\jxi + \ulell \xi)} \bigl( 1 - \chi_{\{\leq 2\ulg\ulell\}}(\xi) \bigr) i \gulellsh(t,\xi) \Bigr](y) \Bigr) \\ 
                &\quad + \Re \, \biggl( \int_\bbR \Bigl( \jxi^{-1} (D^\ast \eulsharp)(y,\xi) - \ulg^{-1} (D^\ast \eulsharp)(y,-\ulg \ulell) \Bigr) e^{i t (\jxi + \ulell \xi)} \chi_{\{\leq 2\ulg\ulell\}}(\xi) i \gulellsh(t,\xi) \, \ud \xi \biggr).
            \end{aligned}
        \end{equation*}
Hence, using the bootstrap assumptions \eqref{equ:consequences_cor_assumption2}, we infer the improved local decay estimate \eqref{equ:consequences_Remusubetwo_L2y_local_decay} from \eqref{equ:improved_local_decay_high_freq_for_2nd_comp} with $k=0$ and from \eqref{equ:local_decay_resonance_subtracted_off_for_2nd_comp}.

The decay estimate \eqref{equ:consequences_hulell_decay} for $h_\ulell(t)$ and the decay estimate \eqref{equ:consequences_pulell_decay} for $p_\ulell(t)$ follow from Lemma~\ref{lem:core_linear_dispersive_decay} and the bootstrap assumptions \eqref{equ:consequences_cor_assumption2}.

Finally, we establish the improved decay estimates \eqref{equ:consequences_hulell_phase_filtered_decay} and \eqref{equ:consequences_pulell_phase_filtered_decay} for the time derivatives of the phase-filtered terms $\pt \bigl( e^{-it\ulg^{-1}} h_\ulell(t) \bigr)$, respectively $\pt \bigl( e^{-it\ulg^{-1}} p_\ulell(t) \bigr)$. 
Due to the presence of the low frequency cut-off in the definitions of $h_\ulell(t)$ and $p_\ulell(t)$, the proofs of \eqref{equ:consequences_hulell_phase_filtered_decay} and \eqref{equ:consequences_pulell_phase_filtered_decay} are essentially identical. 
We provide the details for \eqref{equ:consequences_hulell_phase_filtered_decay}.
By direct computation we find
\begin{equation*}
    \begin{aligned}
        \pt \bigl( e^{-it\ulg^{-1}} h_\ulell(t) \bigr) &= e^{-it\ulg^{-1}} \int_\bbR (\jxi+\ulell\xi-\ulg^{-1}) e^{it(\jxi+\ulell\xi)} \chi_{\{\leq 2\ulg\ulell\}}(\xi) \, \jxi^{-1} \gulellsh(t,\xi) \, \ud \xi \\ 
        &\quad + e^{-it\ulg^{-1}}  \int_\bbR e^{it(\jxi+\ulell\xi)} \chi_{\{\leq 2\ulg\ulell\}}(\xi) \, i \jxi^{-1} \pt \gulellsh(t,\xi) \, \ud \xi \\
        &=: I(t) + II(t).
    \end{aligned}
\end{equation*}
To bound the first term $I(t)$ we exploit that the factor $(\jxi+\ulell\xi-\ulg^{-1})$ vanishes at $\xi = -\ulg\ulell$, where the phase $e^{it(\jxi+\ulell\xi)}$ is stationary. Correspondingly, we can integrate by parts in $\xi$, and then conclude by the Cauchy-Schwarz inequality, the finite support of $\chi_{\{\leq 2\ulg|\ulell|\}}(\xi)$, and the assumption \eqref{equ:consequences_cor_assumption2} that
\begin{equation*}
    \begin{aligned}
        |I(t)| &\lesssim \frac{1}{t} \Bigl\| \pxi \Bigl( \chi_{\{\leq 2\ulg\ulell\}}(\xi) \bigl(\jxi+\ulell\xi-\ulg^{-1}\bigr) (\xi+\ulell\jxi)^{-1} \Bigr) \Bigr\|_{L^2_\xi} \bigl\|\gulellsh(t,\xi)\bigr\|_{L^2_\xi} \\ 
        &\quad + \frac{1}{t} \Bigl\| \chi_{\{\leq 2\ulg\ulell\}}(\xi) \bigl(\jxi+\ulell\xi-\ulg^{-1}\bigr) (\xi+\ulell\jxi)^{-1} \Bigr\|_{L^2_\xi} \bigl\| \pxi \gulellsh(t,\xi) \bigr\|_{L^2_\xi} \\ 
        &\lesssim_{\ell_0} \varepsilon t^{-1+\delta}.
    \end{aligned}
\end{equation*}
Next, we consider the second term $II(t)$. Inserting the evolution equation \eqref{equ:setting_up_g_evol_equ3} for the profile $\gulellsh(t,\xi)$, we find that
\begin{equation*}
    \begin{aligned}
        II(t) &= - e^{-it\ulg^{-1}} (\dotq(t) - \ulell) \int_\bbR e^{it(\jxi+\ulell\xi)} \chi_{\{\leq 2\ulg\ulell\}}(\xi) \jxi^{-1} \xi \,  \gulellsh(t,\xi) \, \ud \xi \\ 
        &\quad - e^{-it\ulg^{-1}} \int_\bbR \chi_{\{\leq 2\ulg\ulell\}}(\xi) \, i \jxi^{-1} \calFulellsh\bigl[ \calQ_{\ulell}\bigl(\usubeone(t)\bigr) \bigr](\xi) \, \ud \xi \\ 
        &\quad - e^{-it\ulg^{-1}} \int_\bbR \chi_{\{\leq 2\ulg\ulell\}}(\xi) \, i \jxi^{-1} \calFulellsh\bigl[ {\textstyle \frac16} \usubeone(t)^3 \bigr](\xi) \, \ud \xi \\ 
        &\quad + e^{-it\ulg^{-1}} \int_\bbR \chi_{\{\leq 2\ulg\ulell\}}(\xi) \, i \jxi^{-1} \widetilde{\calR}(t,\xi) \, \ud \xi \\ 
        &=: II_1(t) + II_2(t) + II_3(t) + II_4(t).
    \end{aligned}
\end{equation*}
For the first term $II_1(t)$ we conclude from \eqref{equ:consequences_qdot_minus_ulell_decay}, \eqref{equ:consequences_cor_assumption2}, and the dispersive decay estimate from Lemma~\ref{lem:core_linear_dispersive_decay} that
\begin{equation*}
    \begin{aligned}
        |II_1(t)| \lesssim |\dotq(t)-\ulell| \biggl| \int_\bbR e^{it(\jxi+\ulell\xi)} \chi_{\{\leq 2\ulg\ulell\}}(\xi) \, \jxi^{-1} \xi \gulellsh(t,\xi) \, \ud \xi \biggr| \lesssim_{\ell_0} \varepsilon t^{-1+\delta} \cdot \varepsilon t^{-\frac12} \lesssim_{\ell_0} \varepsilon t^{-\frac32+\delta}.
    \end{aligned}
\end{equation*}
Next, we consider the second term $II_2(t)$. Using the Cauchy-Schwarz inequality, the mapping property \eqref{equ:mapping_property_calFulellsh_L2}, and \eqref{equ:consequences_dispersive_decay_ulPe_radiation}, we obtain that
\begin{equation*}
    \begin{aligned}
        |II_2(t)| &\lesssim \bigl\| \calFulellsh\bigl[ \calQ_{\ulell}\bigl(\usubeone(t)\bigr) \bigr](\xi) \bigr\|_{L^2_\xi} \lesssim \| \alpha(\ulg y) \|_{L^2_y} \|\usubeone(t)\|_{L^\infty_y}^2 \lesssim \varepsilon^2 \jt^{-1}.
    \end{aligned}
\end{equation*}
Similarly, by the mapping property \eqref{equ:mapping_property_calFulellsh_L2}, \eqref{equ:smallness_orbital_proposition}, and \eqref{equ:consequences_dispersive_decay_ulPe_radiation}, we conclude
\begin{equation*}
    \begin{aligned}
        |II_3(t)| &\lesssim \bigl\| \calFulellsh\bigl[ {\textstyle \frac16} \usubeone(t)^3 \bigr](\xi) \bigr\|_{L^2_\xi} \lesssim \|\usubeone(t)\|_{L^2_y} \|\usubeone(t)\|_{L^\infty_y}^2 \lesssim \varepsilon \cdot \varepsilon^2 \jt^{-1} \lesssim \varepsilon^3 \jt^{-1}.
    \end{aligned}
\end{equation*}
Finally, we turn to the fourth term $II_4(t)$. 
Using the Cauchy-Schwarz inequality and the bound \eqref{equ:consequences_wtilcalR_L2xi_bound}, we obtain
\begin{equation*}
    |II_4(t)| \lesssim \bigl\| \widetilde{\calR}(t,\xi) \bigr\|_{L^2_\xi} \lesssim \varepsilon^2 \jt^{-\frac32+\delta}.
\end{equation*}
Combining the preceding estimates proves \eqref{equ:consequences_hulell_phase_filtered_decay}.

\noindent \underline{Proof of (11).} 
We discuss the proofs of the decay estimates and Sobolev bounds for the auxiliary Klein-Gordon evolutions~\eqref{equ:consequences_aux_KG_evolutions_def} item by item.

The dispersive decay estimates~\eqref{equ:consequences_aux_KG_disp_decay} follow from Lemma~\ref{lem:core_linear_dispersive_decay} using assumption \eqref{equ:consequences_cor_assumption2} and the fact that $\frakb \in W^{1,\infty}_\xi$.
Next, the Sobolev bounds~\eqref{equ:consequences_aux_KG_energy_bounds} are a consequence of an application of the standard Plancherel theorem combined with the uniform-in-time $L^2_\xi$ bound \eqref{equ:consequences_uniform_in_time_L2_bound} for the effective profile, respectively the slow growth bound for the frequency weighted norm $\|\jxi^2 \gulellsh(t,\xi)\|_{L^2_\xi}$ encoded in the assumed bounds \eqref{equ:consequences_cor_assumption2} on the effective profile.
Finally, the local decay estimate~\eqref{equ:consequences_aux_KG_local_H3y_decay} can be proved along the lines of the proof of the local decay estimate~\eqref{equ:consequences_local_decay_usubeonetwo} in item (7) above, working with the standard flat Fourier transform instead of the distorted Fourier transform $\calFulellsh$. We omit the details of the minor variants of the arguments.

To define the leading order local decay decomposition~\eqref{equ:consequences_aux_KG_leading_order_local_decay_decomp} for the auxiliary Klein-Gordon evolutions, we denote again by $\chi_{\{\leq 2\ulg|\ulell|\}}(\xi)$ a smooth even bump function satisfying $\chi_{\{\leq 2\ulg|\ulell|\}}(\xi) = 1$ for $|\xi| \leq 2\ulg|\ulell|$ and $\chi_{\{\leq 2\ulg|\ulell|\}}(\xi) = 0$ for $|\xi| \geq 4\ulg|\ulell|$.
We define the decomposition~\eqref{equ:consequences_aux_KG_leading_order_local_decay_decomp} implicitly by setting
\begin{equation}
    \tilde{h}_{\ulell}(t) := \frac{1}{\sqrt{2\pi}} \int_\bbR e^{it(\jxi+\ulell\xi)} \chi_{\{\leq 2 \ulg|\ulell|\}}(\xi) \jxi^{-1} \gulellsh(t,\xi) \, \ud \xi.
\end{equation}
Since the definition of $\tilde{h}_\ulell(t)$ is identical to that of $h_\ulell(t)$ in \eqref{equ:setting_up_hulell_definition} up to a constant factor, the decay estimates \eqref{equ:consequences_aux_KG_tildeh_ulell_decay}, \eqref{equ:consequences_aux_KG_tildeh_ulell_phase_filtered_decay} for $\tilde{h}_\ulell(t)$ just follow from the corresponding estimates \eqref{equ:consequences_hulell_decay}, \eqref{equ:consequences_hulell_phase_filtered_decay} for~$h_\ulell(t)$.
We turn to the proof of the improved local decay estimate \eqref{equ:consequences_aux_KG_Remv_improved_H1y_local_decay} for the remainder term $R_v(t,y)$, which is explicitly given by
\begin{equation}
    \begin{aligned}
        R_v(t,y) &= \frac{1}{\sqrt{2\pi}} \int_\bbR e^{i y \xi} e^{it(\jxi+\ulell\xi)} \bigl( 1 - \chi_{\{\leq 2\ulg|\ulell|\}}(\xi) \bigr) \frakb(\xi) \jxi^{-1} \gulellsh(t,\xi) \, \ud \xi \\ 
        &\quad + \frac{1}{\sqrt{2\pi}} \int_\bbR \bigl( e^{i y \xi} - e^{-iy\ulg\ulell} \bigr) e^{it(\jxi+\ulell\xi)} \chi_{\{\leq 2\ulg|\ulell|\}}(\xi) \frakb(\xi) \jxi^{-1} \gulellsh(t,\xi) \, \ud \xi \\ 
        &=: R_v^{(a)}(t,y) + R_v^{(b)}(t,y).
    \end{aligned}
\end{equation}
The proof of the asserted improved local decay estimate~\eqref{equ:consequences_aux_KG_Remv_improved_H1y_local_decay} for the contribution of $R_v^{(a)}(t,y)$ follows along the lines of the proof of \eqref{equ:improved_local_decay_high_freq} in Lemma~\ref{lem:improved_local_decay} with $e_\ell^{\#}(y,\xi)$ replaced by $e^{iy\xi}$, while the proof for the contribution of $R_v^{(b)}(t,y)$ follows along the lines of the proof of \eqref{equ:local_decay_resonance_subtracted_off} in Lemma~\ref{lem:local_decay_resonance_subtracted_off} again with $e_\ell^{\#}(y,\xi)$ replaced by $e^{iy\xi}$. We omit a further discussion of the minor variants of the arguments. 

Finally, we discuss the improved local decay estimate~\eqref{equ:consequences_aux_KG_improved_local_H1y_decay} in the special case where the coefficient $\frakb(\xi)$ is given by $\ulg(\xi+\ulell\jxi) \bigl( |\ulg(\xi+\ulell\jxi)|\pm i \bigr)^{-1}$, which vanishes at the problematic frequency $\xi = -\ulg \ulell$.
Integrating by parts in the frequency variable, we find that
\begin{equation}
    \begin{aligned}
        v(t,y) &= - \frac{1}{\sqrt{2\pi}} \frac{1}{it} \int_\bbR e^{it(\jxi+\ulell\xi)} \pxi \biggl( \frac{\jxi}{\xi+\ulell\jxi} e^{iy\xi} \frakb(\xi) \jxi^{-1} \gulellsh(t,\xi) \biggr) \, \ud \xi.
    \end{aligned}
\end{equation}
Since the $(\xi+\ulell\jxi)^{-1}$ singularity gets cancelled by the coefficient $\frakb(\xi)$, the asserted improved local decay \eqref{equ:consequences_aux_KG_improved_local_H1y_decay} then follows from Plancherel's theorem and the assumption \eqref{equ:consequences_cor_assumption2}.    
\end{proof}

\subsection{Proof of Theorem \ref{thm:main}} \label{subsec:proof_of_main_theorem}

Now we are in the position to establish the proof of Theorem~\ref{thm:main} based on the conclusions of Proposition~\ref{prop:modulation_parameter_control}, Proposition~\ref{prop:profile_bounds}, and  Corollary~\ref{cor:consequences_bootstrap_assumptions}.

\begin{proof}[Proof of Theorem \ref{thm:main}] 
Fix $\ell_0 \in (-1,1)$. Let $0 < \varepsilon_0 \ll 1$ and $C_0 \geq 1$ be the constants from the statements of Proposition~\ref{prop:modulation_parameter_control} and Proposition~\ref{prop:profile_bounds}. 
Let $x_0 \in \bbR$ and assume that $\bmu_0 = (u_{0,1}, u_{0,2}) \in H^3_x(\bbR) \times H^2_x(\bbR)$ satisfies the smallness condition in the weighted Sobolev norm \eqref{equ:statement_theorem_smallness_data}.
By Lemma~\ref{lem:setting_up_local_existence} the sine-Gordon equation \eqref{equ:intro_sG_1st_order} has a unique $H^3_x \times H^2_x$-solution $\bmphi(t)$ with initial data
\begin{equation*}
    \bmphi(0,x) = \bmK_{\ell_0,x_0}(x) + \bmu_0(x-x_0)
\end{equation*}
defined on a maximal interval of existance $[0,T_\ast)$ for some $0 < T_\ast \leq \infty$. 
Moreover, the continuation criterion \eqref{equ:setting_up_local_existence_continuation_criterion} holds.
By Proposition~\ref{prop:setting_up_modulation_orbital} there exist unique continuously differentiable paths $\ell \colon [0,T_\ast) \to (-1,1)$ and $q \colon [0,T_\ast) \to \bbR$ so that  $(1)$--$(6)$ in the statement of Proposition~\ref{prop:setting_up_modulation_orbital} hold. 
In particular, the solution $\bmphi(t)$ admits a decomposition into a modulated kink and a radiation term
\begin{equation} \label{equ:proof_of_theorem_decomposition_modulated_kink_plus_radiation}
    \bmphi(t,x) = \bmK_{\ell(t), q(t)}(x) + \bmu\bigl(t, x-q(t)\bigr), \quad 0 \leq t < T_\ast.
\end{equation}
We now consider the exit time 
\begin{equation*}
    \begin{aligned}
        T_0 := \sup \, \biggl\{ 0 \leq T < T_\ast \, \bigg| \, \bigl\| g_{\ell(T)}^\# \bigr\|_{X(T)} \leq C_0 \varepsilon, \quad &\sup_{0 \leq t \leq T} \, \jt^{1-\delta} \bigl| \ell(t) - \ell(T) \bigr| \leq C_0 \varepsilon \biggr\},
    \end{aligned}
\end{equation*}
where $\|\cdot\|_{X(T)}$ is the norm defined in \eqref{equ:bootstrap_norm_XT_definition} and where $g_{\ell(T)}^{\#}(t,\xi)$ denotes the effective profile relative to the linearized operator~$\bfL_{\ell(T)}$ defined in \eqref{equ:setting_up_definition_g}.
By the local existence theory we have $T_0 > 0$ as long as we choose $C_0 \gg 1$ sufficiently large.

We first show that $T_0 = T_\ast$. To this end we argue by contradiction. Suppose that $0 < T_0 < T_\ast$.
Then we pick a strictly monotone increasing sequence $\{T_n\}_{n \in \bbN} \subset [0,T_0)$ with $T_n \nearrow T_0$ as $n \to \infty$. Since the paths $\ell \colon [0,T_\ast) \to (-1,1)$ and $q \colon [0,T_\ast) \to \bbR$ are continuously differentiable, it follows that
\begin{equation}
    \sup_{0 \leq t \leq T_0} \, \jt^{1-\delta} |\ell(t) - \ell(T_0)| \leq C_0 \varepsilon.
\end{equation}
Since we have by the mapping properties \eqref{equ:mapping_property_calFulellsh_jxi2}, \eqref{equ:mapping_property_calFulellsh_jxi2_pxi}, \eqref{equ:mapping_property_calFulellsh_jxi32_Linfty}, \eqref{equ:mapping_property_calFulellDsh_jxi2}, \eqref{equ:mapping_property_calFulellDsh_jxi2_pxi}, \eqref{equ:mapping_property_calFulellDsh_jxi32_Linfty} that
\begin{equation} \label{equ:proof_of_theorem_initial_profile_norms_bound}
    \begin{aligned}
        \bigl\| \jxi^{\frac32} g_\ulell^\#(0,\xi) \bigr\|_{L^\infty_\xi} + \bigl\| \jxi^2 g_\ulell^\#(0,\xi) \bigr\|_{L^2_\xi} + \bigl\| \jxi^2 \pxi g_\ulell^\#(0,\xi) \bigr\|_{L^2_\xi} \lesssim \| \jx \bmu_0 \|_{H^3_x \times H^2_x} \lesssim \varepsilon,
    \end{aligned}
\end{equation}
we can use Proposition~\ref{prop:profile_bounds} with $\ulell = \ell(T_0)$ fixed to conclude via a continuity argument that in fact
\begin{equation}
    \bigl\| g_{\ell(T_0)}^\# \bigr\|_{X(T_0)} \leq C_0 \varepsilon. 
\end{equation}
Keeping $\ell(T_0)$ fixed, it follows by continuity that there exists $T_0 < \widetilde{T} < T_\ast$ such that 
\begin{align}
    \sup_{0 \leq t \leq \wtilT} \, \jt^{1-\delta} |\ell(t) - \ell(T_0)| &\leq 2C_0\varepsilon, \label{equ:proof_of_theorem_bootstrap1} \\
    \bigl\| g_{\ell(T_0)}^\# \bigr\|_{X(\wtilT)} &\leq 2C_0\varepsilon. \label{equ:proof_of_theorem_bootstrap2}
\end{align}
Correspondingly, thanks to \eqref{equ:proof_of_theorem_bootstrap1}, \eqref{equ:proof_of_theorem_bootstrap2}, the assumptions of Proposition~\ref{prop:modulation_parameter_control} are satisfied on $[0,\wtilT]$ with $\ulell = \ell(T_0)$ fixed, and we obtain that
\begin{equation} \label{equ:proof_of_theorem_wtilT_bound1}
    \begin{aligned} 
        \sup_{0 \leq t \leq \wtilT} \, \jt^{1-\delta} |\ell(t) - \ell(\wtilT)| &\leq C_0 \varepsilon.
    \end{aligned}
\end{equation}
Having established the bound \eqref{equ:proof_of_theorem_wtilT_bound1}, we can now run a continuity argument using Proposition~\ref{prop:profile_bounds} with $\ulell = \ell(\wtilT)$ fixed, to conclude that 
\begin{equation} \label{equ:proof_of_theorem_wtilT_bound2}
    \bigl\| g_{\ell(\wtilT)}^\# \bigr\|_{X(\wtilT)} \leq C_0 \varepsilon.
\end{equation}
But then \eqref{equ:proof_of_theorem_wtilT_bound1} and \eqref{equ:proof_of_theorem_wtilT_bound2} are a contradiction to $\wtilT > T_0$. Hence, we must have that $T_0 = T_\ast$. 

Next, we conclude that $T_\ast = \infty$. Since $\|g_{\ell(T)}^{\#}\|_{X(T)} \leq C_0 \varepsilon$ for all $0 < T < T_\ast$, the Sobolev bound \eqref{equ:consequences_sobolev_bound_radiation} implies for all $0 \leq t < T_\ast$ that $\|\bmu(t)\|_{H^3_x \times H^2_x} \lesssim \varepsilon \jt^{\delta} \lesssim \varepsilon \jap{T_\ast}^{\delta}$.
Using the decomposition \eqref{equ:proof_of_theorem_decomposition_modulated_kink_plus_radiation}, we can then infer from the continuation criterion~\eqref{equ:setting_up_local_existence_continuation_criterion}  that $T_\ast = \infty$.

Now we determine the Lorentz boost parameter of the final soliton.
To this end we pick a strictly increasing sequence $\{T_n\}_{n\in\bbN} \subset [0,\infty)$ with $T_n \nearrow T_\ast = \infty$ as $n\to\infty$. Using that for each $n \in \bbN$ we have 
\begin{equation}
    \sup_{0 \leq t \leq T_n} \, \jt^{1-\delta} |\ell(t) - \ell(T_n)| \leq C_0 \varepsilon,
\end{equation}
we conclude that $\{ \ell(T_n) \}_n \subset (-1,1)$ is a Cauchy sequence with $|\ell(T_n) - \ell_0| \leq \frac14 \gamma_0^{-2}$. Hence, there exists $\ell_\infty \in (-1,1)$ with $|\ell_\infty - \ell_0| \leq \frac14 \gamma_0^{-2}$ such that $\ell(T_n) \to \ell_\infty$ as $n \to \infty$. Since the path $\ell \colon [0,\infty) \to (-1,1)$ is continuous, the value of the asymptotic Lorentz boost parameter $\ell_\infty$ is independent of the choice of sequence $\{T_n\}_n$. Moreover, it follows that
\begin{equation} \label{equ:proof_of_theorem_ell_infty_bound1}
    \sup_{0 \leq t < \infty} \, \jt^{1-\delta} |\ell(t) - \ell_\infty| \leq C_0 \varepsilon.
\end{equation}
Finally, owing to the preceding bound \eqref{equ:proof_of_theorem_ell_infty_bound1} as well as \eqref{equ:proof_of_theorem_initial_profile_norms_bound}, we can run one more continuity argument using Proposition~\ref{prop:profile_bounds} with $\ulell = \ell_\infty$ fixed to obtain uniformly for all $0 < T < \infty$ that
\begin{equation} \label{equ:proof_of_theorem_ell_infty_bound2}
    \bigl\| g_{\ell_\infty}^\# \bigr\|_{X(T)} \leq C_0 \varepsilon.
\end{equation}

Thanks to \eqref{equ:proof_of_theorem_ell_infty_bound1} and \eqref{equ:proof_of_theorem_ell_infty_bound2}, the assumptions \eqref{equ:consequences_cor_assumption1} and \eqref{equ:consequences_cor_assumption2} of Corollary~\ref{cor:consequences_bootstrap_assumptions} are satisfied with the choice $\ulell = \ell_\infty$. 
Then the asserted dispersive decay estimates~\eqref{equ:statement_theorem_dispersive_decay_estimate} in the statement of Theorem~\ref{thm:main} are a consequence of the dispersive decay estimates~\eqref{equ:consequences_dispersive_decay_radiation} from Corollary~\ref{cor:consequences_bootstrap_assumptions}.
Moreover, the asserted asymptotics of the modulation parameters \eqref{equ:statement_theorem_modulation_parameters_asymptotics} in the statement of Theorem~\ref{thm:main} are a consequence of \eqref{equ:proof_of_theorem_ell_infty_bound1} and of \eqref{equ:consequences_qdot_minus_ulell_decay} from Corollary~\ref{cor:consequences_bootstrap_assumptions}.

It remains to deduce the asserted asymptotics~\eqref{equ:statement_theorem_radiation_asymptotics} of the radiation term.
To this end we invoke results that will be proved later in Section~\ref{sec:pointwise_profile_bounds}.
Using the difference bound \eqref{equ:pointwise_proposition_asserted_difference_bound} in the statement of Proposition~\ref{prop:pointwise_profile_bounds} with the choice $\ulell = \ell_\infty$, we infer that there exists an asymptotic effective profile $\widetilde{g}_\infty \in L^\infty_\xi(\bbR)$ such that
\begin{equation} \label{equ:proof_of_theorem_limiting_profile1}
    \Bigl\| e^{-i\widetilde{\Lambda}_\infty(t,\xi)} e^{-i\xi\theta_\infty(t)} \jxi^{\frac32} g_{\ell_\infty}^{\#}(t,\xi) - \widetilde{g}_\infty(\xi) \Bigr\|_{L^\infty_\xi} \lesssim \varepsilon^2 t^{-2\delta}, \quad t \geq 1,
\end{equation}
where 
\begin{equation*}
    \widetilde{\Lambda}_\infty(t,\xi) := \frac{1}{16} \jxi^{-3} \int_1^t \frac{1}{s} \bigl| \jxi^{\frac32} g_{\ell_\infty}^{\#}(s,\xi) \bigr|^2 \, \ud s
\end{equation*}
and 
\begin{equation*}
    \theta_\infty(t) := \int_0^t \bigl( \dot{q}(s) - \ell_\infty \bigr) \, \ud s.
\end{equation*}
Next, we multiply the differential equation \eqref{equ:pointwise_evol_equ_g_renorm_leading_order} for the effective profile $g_{\ell_\infty}^{\#}(t,\xi)$ obtained at the beginning of the proof of Proposition~\ref{prop:pointwise_profile_bounds} by the integrating factor $e^{-i\Gamma_\infty(\xi) \log(t)}$ with 
\begin{equation}
    \Gamma_\infty(\xi) := \frac{1}{16} \jxi^{-3} \bigl| \widetilde{g}_\infty(\xi) \bigr|^2.
\end{equation}
Then we repeat the arguments in the proof of Proposition~\ref{prop:pointwise_profile_bounds} and invoke \eqref{equ:proof_of_theorem_limiting_profile1} to infer that there exists an asymptotic effective profile $g_\infty \in L^\infty_\xi(\bbR)$ with $|g_\infty(\xi)| = |\tilde{g}_\infty(\xi)|$ such that
\begin{equation} \label{equ:proof_of_theorem_limiting_profile2}
    \Bigl\| e^{-i\Gamma_\infty(\xi)\log(t)} e^{-i\xi\theta_\infty(t)} \jxi^{\frac32} g_{\ell_\infty}^{\#}(t,\xi) - g_\infty(\xi) \Bigr\|_{L^\infty_\xi} \lesssim \varepsilon^2 t^{-2\delta}, \quad t \geq 1.
\end{equation}
In the moving frame coordinate $y := x-q(t)$ we have by \eqref{equ:consequences_decomposition_radiation} from Corollary~\ref{cor:consequences_bootstrap_assumptions} the decomposition
\begin{equation}
    \bmu(t,y) = \bigl( e^{t\bfL_{\ell_\infty}} \ulPe \bmf_\infty(t) \bigr)(y) + d_{0,\ell_\infty}(t) \bmY_{0,\ell_\infty}(y) + d_{1,\ell_\infty}(t) \bmY_{1,\ulell_\infty}(y),
\end{equation}
where the discrete components enjoy the fast decay $|d_{0,\ell_\infty}(t)| + |d_{1,\ell_\infty}(t)| \lesssim \varepsilon \jt^{-\frac32+\delta}$.
The asserted asymptotics for the radiation term \eqref{equ:statement_theorem_radiation_asymptotics} in the statement of Theorem~\ref{thm:main} now follow from Corollary~\ref{cor:linear_asymptotics} using that 
\begin{equation*}
    y + \ell_\infty t = x - q(t) + \ell_\infty t = x - \theta_\infty(t) - q(0).
\end{equation*}
This finishes the proof of Theorem~\ref{thm:main}.
\end{proof}

\section{Control of the Modulation Parameters} \label{sec:modulation_control}

In this section we establish the proof of Proposition~\ref{prop:modulation_parameter_control}. 

\begin{proof}[Proof of Proposition~\ref{prop:modulation_parameter_control}]  
We start off with an elementary application of the fundamental theorem of calculus to write for $0 \leq t \leq T$,
\begin{equation} \label{equ:modulation_proof_ellt_minus_ellT}
    \ell(t) - \ell(T) = - \int_t^T \dot{\ell}(s) \, \ud s.
\end{equation}
Now the idea is to insert the equation for $\dot{\ell}(s)$ stemming from the modulation equations \eqref{equ:modulation_equ}, and to exploit the oscillations in $\dot{\ell}(s)$ to infer the decay estimate 
\begin{equation} \label{equ:modulation_proof_goal}
    |\ell(t)-\ell(T)| \lesssim \varepsilon^2 \jt^{-1+\delta}.
\end{equation}
The latter then implies the asserted (improved) bound \eqref{equ:prop_modulation_parameters_conclusion} in the statement of Proposition~\ref{prop:modulation_parameter_control}.
However, the equation for $\dot{\ell}(s)$ features problematic resonant terms in the leading order quadratic nonlinearity. Their worst effects are suppressed by a remarkable null structure specific to the sine-Gordon model. 

From the modulation equations \eqref{equ:modulation_equ} we obtain
       \begin{equation}  \label{equ:modulation_proof_recall_modulation_equations}
            \begin{bmatrix} \dot{\ell} \\ \dot{q}-\ell \end{bmatrix} =  \bigl[ \bbM_{\ell}[\bmu] \bigr]^{-1}
        \begin{bmatrix}
            - \bigl\langle \bfJ \bmY_{0,\ell}, \calN(\bm{u}) \bigr\rangle \\
            \bigl\langle \bfJ \bmY_{1,\ell}, \calN(\bm{u}) \bigr\rangle 
        \end{bmatrix}
        \end{equation}
        with
    \begin{equation} \label{equ:modulation_proof_bbM_matrix_recalled}
        \begin{aligned}
        \bbM_{\ell}[\bmu] &= 
        \begin{bmatrix}
            \gamma^3 \|K'\|_{L^2}^2 & 0 \\ 0 & \gamma^3 \|K'\|_{L^2}^2
        \end{bmatrix}
        + 
        \begin{bmatrix} 
        \bigl\langle \bfJ \bmZ_{3,\ell}, \bm{u} \bigr\rangle &  \bigl\langle \bfJ \bmZ_{1,\ell}, \bm{u} \bigr\rangle \\
        \bigl\langle \bfJ \bmZ_{2,\ell}, \bm{u} \bigr\rangle & - \bigl\langle \bfJ \bmZ_{3,\ell}, \bm{u} \bigr\rangle
        \end{bmatrix},
        \end{aligned}
    \end{equation}
    where we recall the notation
\begin{equation} \label{equ:modulation_proof_recall_bmY_def}
    \bmY_{0,\ell}(y) := \partial_q \bm{K}_{\ell,q}(y), \quad \bmY_{1,\ell}(y) := \partial_\ell \bm{K}_{\ell,q}(y), 
\end{equation}
and where we introduce the short-hand notation
\begin{equation} \label{equ:modulation_proof_bmZ_def}
    \bmZ_{1,\ell}(y) := \partial_q \partial_q \bm{K}_{\ell,q}(y), \quad 
    \bmZ_{2,\ell}(y) := \partial_\ell \partial_\ell \bm{K}_{\ell,q}(y), \quad
    \bmZ_{3,\ell}(y) := \partial_q \partial_\ell \bm{K}_{\ell,q}(y).
\end{equation}
Observe that the modulation parameter $\ell \equiv \ell(s)$ is time-dependent in \eqref{equ:modulation_proof_recall_modulation_equations}, \eqref{equ:modulation_proof_bbM_matrix_recalled}, \eqref{equ:modulation_proof_recall_bmY_def}, \eqref{equ:modulation_proof_bmZ_def}.
The matrix $\bbM_{\ell}[\bmu]$ is invertible for all times $0 \leq s \leq T$ and by \eqref{equ:smallness_orbital_proposition} we have the uniform-in-time operator norm bound
\begin{equation} \label{equ:modulation_proof_op_norm_inv_bbM_bmu}
    \sup_{0 \leq s \leq T} \, \bigl\| \bigl[ \bbM_{\ell(s)}[\bmu(s)] \bigr]^{-1} \bigr\| \lesssim_{\ell_0} 1.
\end{equation}
In order to accurately capture the oscillations of the leading order quadratic and cubic nonlinearities on the right-hand side of \eqref{equ:modulation_proof_recall_modulation_equations}, we need to decompose the inverse $\bigl[ \bbM_{\ell(s)}[\bmu(s)] \bigr]^{-1}$ into a time-independent leading order term, a small slowly decaying term, and a remainder term with faster decay.
To this end, we write 
        \begin{equation} \label{equ:modulation_proof_decomposition_bbM}
            \bbM_{\ell}[\bmu(s)] = \bbM_\ulell + \bbA_\ulell(s) + \bbB(s),
        \end{equation}  
        where
        \begin{equation*}
            \begin{aligned}
                \bbM_\ulell &:= \begin{bmatrix} \ulg^3 \|K'\|_{L^2}^2 & 0 \\ 0 & \ulg^3 \|K'\|_{L^2}^2 \end{bmatrix}, \quad
                \bbA_\ulell(s) 
                    &:= \begin{bmatrix} 
                            \bigl\langle \bfJ \bmZ_{3,\ulell}, (\ulPe \bmu)(s) \bigr\rangle &  \bigl\langle \bfJ \bmZ_{1,\ulell}, (\ulPe \bmu)(s) \bigr\rangle \\
                            \bigl\langle \bfJ \bmZ_{2,\ulell}, (\ulPe \bmu)(s) \bigr\rangle & - \bigl\langle \bfJ \bmZ_{3,\ulell}, (\ulPe \bmu)(s) \bigr\rangle
                        \end{bmatrix}, 
            \end{aligned}
        \end{equation*}
        and 
        \begin{equation*}
            \bbB(s) := \bbM_{\ell}[\bmu(s)] - \bbM_\ulell - \bbA_\ulell(s).
        \end{equation*}
        Clearly, we have 
        \begin{equation} \label{equ:modulation_proof_bbM_ulell_inverse_op_norm}
            \bigl\| \bbM_\ulell^{-1} \bigr\| \lesssim_{\ell_0} 1.
        \end{equation}
        Moreover, the dispersive decay estimates \eqref{equ:consequences_dispersive_decay_ulPe_radiation} for the radiation term imply the following slowly decaying operator norm bound
        \begin{equation} \label{equ:modulation_proof_op_norm_bbA}
            \sup_{0 \leq s \leq T} \, \js^{\frac12} \bigl\| \bbA_\ulell(s) \bigr\| \lesssim \varepsilon.
        \end{equation}
        Upon closer inspection, one sees that every entry of the matrix $\bbB(s)$ features a discrete component from the decomposition~\eqref{equ:consequences_decomposition_radiation} of $\bmu(s)$ or it can be bounded by $|\ell(s)-\ulell|$.
        Thus, \eqref{equ:prop_modulation_parameters_assumption1} and \eqref{equ:consequences_discrete_comp_decay} imply the following operator norm bound with faster decay
        \begin{equation} \label{equ:modulation_proof_op_norm_bbB}
            \sup_{0\leq s \leq T} \, \js^{1-\delta} \bigl\| \bbB(s) \bigr\| \lesssim \varepsilon.
        \end{equation}

        From \eqref{equ:modulation_proof_decomposition_bbM} we obtain the following expansion 
        \begin{equation} \label{equ:modulation_proof_bbM_inverse_expansion}
            \bigl[ \bbM_{\ell}[\bmu(s)] \bigr]^{-1} = \bbM_\ulell^{-1} - \bbM_\ulell^{-1} \bbA_\ulell(s) \bbM_\ulell^{-1} + \bbD(s)
        \end{equation}
        with 
        \begin{equation*}
            \begin{aligned}
                \bbD(s) := \bigl[ \bbM_{\ell}[\bmu(s)] \bigr]^{-1} \Bigl( -\bbB(s) \bbM_\ulell^{-1} + \bbA_\ulell(s) \bbM_\ulell^{-1} \bbA_\ulell(s) \bbM_\ulell^{-1} + \bbB(s) \bbM_\ulell^{-1} \bbA_\ulell(s) \bbM_\ulell^{-1} \Bigr).
            \end{aligned}
        \end{equation*}
        By \eqref{equ:modulation_proof_bbM_ulell_inverse_op_norm}, \eqref{equ:modulation_proof_op_norm_bbA}, and \eqref{equ:modulation_proof_op_norm_bbB}, the remainder $\bbD(s)$ satisfies the faster decaying operator norm bound 
        \begin{equation} \label{equ:modulation_proof_op_norm_bbD}
            \sup_{0 \leq s \leq T} \, \js^{1-\delta} \bigl\| \bbD(s) \bigr\| \lesssim \varepsilon.
        \end{equation}

        Next, we decompose the nonlinearities on the right-hand side of the modulation equations \eqref{equ:modulation_proof_recall_modulation_equations} into leading order quadratic as well as cubic terms, and a remainder term.
        Using the notation and definitions introduced in Section~\ref{sec:setting_up}, we write
        \begin{equation} \label{equ:modulation_proof_nonlinearities_expansion}
            \begin{aligned}
                \begin{bmatrix}
                - \bigl\langle \bfJ \bmY_{0, \ell}, \calN(\bm{u}) \bigr\rangle \\
                \bigl\langle \bfJ \bmY_{1, \ell}, \calN(\bm{u}) \bigr\rangle 
                \end{bmatrix}
                &=  \begin{bmatrix} 
                        \bigl\langle \bigl( \bmY_{0, \ulell} \bigr)_1, \calQ_\ulell\bigl( \usubeone(s) \bigr) \bigr\rangle \\
                        - \bigl\langle \bigl( \bmY_{1, \ulell} \bigr)_1, \calQ_\ulell\bigl( \usubeone(s) \bigr) \bigr\rangle
                    \end{bmatrix} 
                    +
                    \begin{bmatrix}
                        \bigl\langle \bigl( \bmY_{0, \ulell} \bigr)_1, \calC\bigl( \usubeone(s) \bigr) \bigr\rangle \\ 
                        - \bigl\langle \bigl( \bmY_{1, \ulell} \bigr)_1, \calC\bigl( \usubeone(s) \bigr) \bigr\rangle
                    \end{bmatrix}
                    +
                    \begin{bmatrix}
                        - E_1(s) \\
                        E_2(s)
                    \end{bmatrix},
            \end{aligned}
        \end{equation}
        where 
        \begin{equation*}
            \begin{aligned}
                E_1(s) &:= - \bigl\langle \bigl( \bmY_{0, \ulell} \bigr)_1, \bigl( \calE_3 + \calR_1 + \calR_2 \bigr) \bigr\rangle + \bigl\langle \bfJ \bmY_{0,\ell}, \calE_2 \bigr\rangle + \bigl\langle \bfJ (\bmY_{0,\ell} - \bmY_{0,\ulell}), \calN\bigl( (\ulPe \bmu)(s) \bigr) \bigr\rangle, \\
                E_2(s) &:= - \bigl\langle \bigl( \bmY_{1, \ulell} \bigr)_1, \bigl( \calE_3 + \calR_1 + \calR_2 \bigr) \bigr\rangle + \bigl\langle \bfJ \bmY_{1,\ell}, \calE_2 \bigr\rangle + \bigl\langle \bfJ (\bmY_{1,\ell} - \bmY_{1,\ulell}), \calN\bigl( (\ulPe \bmu)(s) \bigr) \bigr\rangle.
            \end{aligned}
        \end{equation*}
        Using \eqref{equ:prop_modulation_parameters_assumption1}, \eqref{equ:consequences_dispersive_decay_ulPe_radiation}, \eqref{equ:consequences_discrete_comp_decay}, it is straightforward to see that 
        \begin{equation} \label{equ:modulation_proof_higher_order_nonlin_decay} 
            |E_1(s)| + |E_2(s)| \lesssim \varepsilon^2 \js^{-2+\delta}.
        \end{equation}

        Inserting the decompositions \eqref{equ:modulation_proof_bbM_inverse_expansion} and \eqref{equ:modulation_proof_nonlinearities_expansion} on the right-hand side of \eqref{equ:modulation_proof_recall_modulation_equations}, we find
        \begin{equation}
            \begin{aligned}
            &\begin{bmatrix} \dot{\ell}(s) \\ \dot{q}(s) - \ell(s) \end{bmatrix} 
            =
            \bigl[ \bbM_\ell[\bmu(s)] \bigr]^{-1} 
            \begin{bmatrix}
            - \bigl\langle \bfJ \bmY_{0,\ell}, \calN(\bm{u}) \bigr\rangle \\
            \bigl\langle \bfJ \bmY_{1,\ell}, \calN(\bm{u}) \bigr\rangle 
            \end{bmatrix} \\
            &= \bbM_\ulell^{-1} \begin{bmatrix}
                                    \bigl\langle \bigl( \bmY_{0, \ulell} \bigr)_1, \calQ_\ulell\bigl( \usubeone(s) \bigr) \bigr\rangle \\ 
                                    -\bigl\langle \bigl( \bmY_{1, \ulell} \bigr)_1, \calQ_\ulell\bigl( \usubeone(s) \bigr) \bigr\rangle 
                                \end{bmatrix}
            +
            \bbM_\ulell^{-1} \begin{bmatrix}
                                    \bigl\langle \bigl( \bmY_{0, \ulell} \bigr)_1, \calC\bigl( \usubeone(s) \bigr) \bigr\rangle \\ 
                                    -\bigl\langle \bigl( \bmY_{1, \ulell} \bigr)_1, \calC\bigl( \usubeone(s) \bigr) \bigr\rangle 
                                \end{bmatrix} \\ 
            &\quad -
            \bbM_\ulell^{-1} \bbA_\ulell(s) \bbM_\ulell^{-1}    \begin{bmatrix}
                                                                \bigl\langle \bigl( \bmY_{0, \ulell} \bigr)_1, \calQ_\ulell\bigl( \usubeone(s) \bigr) \bigr\rangle \\ 
                                                                -\bigl\langle \bigl( \bmY_{1, \ulell} \bigr)_1, \calQ_\ulell\bigl( \usubeone(s) \bigr) \bigr\rangle 
                                                                \end{bmatrix}                                
            -
            \bbM_\ulell^{-1} \bbA_\ulell(s) \bbM_\ulell^{-1}    \begin{bmatrix}
                                                                \bigl\langle \bigl( \bmY_{0, \ulell} \bigr)_1, \calC\bigl( \usubeone(s) \bigr) \bigr\rangle \\ 
                                                                -\bigl\langle \bigl( \bmY_{1, \ulell} \bigr)_1, \calC\bigl( \usubeone(s) \bigr) \bigr\rangle 
                                                                \end{bmatrix} \\
            &\quad + 
            \bbD(s) \begin{bmatrix}
                    \bigl\langle \bigl( \bmY_{0, \ulell} \bigr)_1, \calQ_\ulell\bigl( \usubeone(s) \bigr) \bigr\rangle \\ 
                    -\bigl\langle \bigl( \bmY_{1, \ulell} \bigr)_1, \calQ_\ulell\bigl( \usubeone(s) \bigr) \bigr\rangle 
                    \end{bmatrix}
            +
            \bbD(s) \begin{bmatrix}
                    \bigl\langle \bigl( \bmY_{0, \ulell} \bigr)_1, \calC\bigl( \usubeone(s) \bigr) \bigr\rangle \\ 
                    -\bigl\langle \bigl( \bmY_{1, \ulell} \bigr)_1, \calC\bigl( \usubeone(s) \bigr) \bigr\rangle 
                    \end{bmatrix} \\
            &\quad + \bigl[ \bbM_\ell[\bmu(s)] \bigr]^{-1} 
                     \begin{bmatrix}
                     -E_1(s) \\
                     E_2(s) 
                     \end{bmatrix} \\ 
            &=: I(s) + II(s) + III(s) + IV(s) + V(s) + VI(s) + VII(s).
            \end{aligned}
        \end{equation}
        The goal is now to show that the contributions of all terms $I(s), \ldots, VII(s)$ to the integral \eqref{equ:modulation_proof_ellt_minus_ellT} decay at least at the rate $\varepsilon^2 \jt^{-1+\delta}$ asserted in \eqref{equ:modulation_proof_goal}.
        In what follows, most of the work will go into deriving sufficient decay for the contribution of the first term $I(s)$ to the integral \eqref{equ:modulation_proof_ellt_minus_ellT}. This part of our analysis hinges on a remarkable quadratic null structure. To deal with the contributions of the cubic terms $II(s)$ and $III(s)$, we still need to exploit their (non-resonant) oscillatory behavior, while the remaining terms $IV(s)$, $V(s)$, $VI(s)$, $VII(s)$ can be bounded quite crudely.

    \noindent \underline{Contribution of the term $I(s)$.}
    The contribution of (the first component of) the term $I(s)$ to the integral \eqref{equ:modulation_proof_ellt_minus_ellT} is given by
    \begin{equation}
        - \ulg^{-3} \|K'\|_{L^2}^{-2} \int_t^T \bigl\langle \bigl( \bmY_{0, \ulell} \bigr)_1, \calQ_\ulell\bigl( \usubeone(s) \bigr) \bigr\rangle \, \ud s.
    \end{equation}
    Inserting the representation formula \eqref{equ:setting_up_representation_formula_usubeone} for $\usubeone(s,y)$, we compute 
    \begin{equation*}
    \begin{aligned}
        &\int_t^T \bigl\langle \bigl( \bmY_{0, \ulell} \bigr)_1, \calQ_\ulell\bigl( \usubeone(s) \bigr) \bigr\rangle \, \ud s
        = \int_t^T \int_\bbR \bigl(\bmY_{0,\ulell}\bigr)_1(y) \calQ_\ulell\bigl( \usubeone(s, y) \bigr) \, \ud y \, \ud s \\ 
        &= \int_t^T \int_\bbR -\ulg K'(\ulg y) \alpha(\ulg y) \usubeone(s,y)^2 \, \ud y \, \ud s \\
        &= \int_t^T \int_\bbR \ulg \sech^2(\ulg y) \tanh(\ulg y) \Biggl( \Re \biggl( \calFulellshast\Bigl[ e^{is(\jxi+\ulell\xi)} i \jxi^{-1} \gulellsh(s,\xi) \Bigr](y) \biggr) \Biggr)^2 \, \ud y \, \ud s \\
        &= \int_t^T \iint e^{is(\jxione + \ulell \xi_1)} e^{-is(\jxitwo + \ulell \xi_2)} i \jxione^{-1} (-i) \jxitwo^{-1} \gulellsh(s,\xi_1) \overline{\gulellsh(s,\xi_2)} \mu_{\ulell;+-}(\xi_1, \xi_2) \, \ud \xi_1 \, \ud \xi_2 \, \ud s \\ 
        &\quad + \int_t^T \iint e^{is(\jxione + \ulell \xi_1)} e^{is(\jxitwo + \ulell \xi_2)} i \jxione^{-1} i \jxitwo^{-1} \gulellsh(s,\xi_1) \gulellsh(s,\xi_2) \mu_{\ulell;++}(\xi_1, \xi_2) \, \ud \xi_1 \, \ud \xi_2 \, \ud s \\ 
        &\quad + \int_t^T \iint e^{-is(\jxione + \ulell \xi_1)} e^{-is(\jxitwo + \ulell \xi_2)} (-i) \jxione^{-1} (-i) \jxitwo^{-1} \overline{\gulellsh(s,\xi_1)} \, \overline{\gulellsh(s,\xi_2)} \mu_{\ulell;--}(\xi_1, \xi_2) \, \ud \xi_1 \, \ud \xi_2 \, \ud s \\ 
        &=: I_1(t) + I_2(t) + I_3(t),
    \end{aligned}
\end{equation*}
where 
\begin{equation}
    \begin{aligned}
        \mu_{\ulell;+-}(\xi_1, \xi_2) &:= \frac12 \int_\bbR \ulg \sech^2(\ulg y) \tanh(\ulg y) e_\ulell^{\#}(y,\xi_1) \overline{e_\ulell^{\#}(y,\xi_2)} \, \ud y, \\ 
        \mu_{\ulell;++}(\xi_1, \xi_2) &:= \frac14 \int_\bbR \ulg \sech^2(\ulg y) \tanh(\ulg y) e_\ulell^{\#}(y,\xi_1) e_\ulell^{\#}(y,\xi_2) \, \ud y, \\
        \mu_{\ulell;--}(\xi_1, \xi_2) &:= \frac14 \int_\bbR \ulg \sech^2(\ulg y) \tanh(\ulg y) \overline{e_\ulell^{\#}(y,\xi_1)} \, \overline{e_\ulell^{\#}(y,\xi_2)} \, \ud y.
    \end{aligned}
\end{equation}
The fine structure of these quadratic spectral distributions has already been determined in Subsection~\ref{subsec:quadratic_spectral_distributions}.
We observe that the phase of $e^{is(\jxione+\ulell\xi_1-\jxitwo-\ulell\xi_2)}$ in the integrand of the first term $I_1(t)$ has a time resonance, while the phases in the terms $I_2(t)$ and $I_3(t)$ do not.

\medskip 

\noindent \underline{Contribution of the resonant term $I_1(t)$.}
We begin with the analysis of the most delicate resonant term $I_1(t)$, where we exploit a remarkable null structure in the quadratic spectral distribution $\mu_{\ulell;+-}(\xi_1,\xi_2)$. Lemma~\ref{lem:null_structure2} uncovers that $\mu_{\ulell;+-}(\xi_1,\xi_2)$ is a multiple of $(\jxione+\ulell\xi_1-\jxitwo-\ulell\xi_2)$. This suppresses all resonant frequency interactions and makes the term $I_1(t)$ amenable to integration by parts in time.

In fact, a closer inspection of the expression for $\mu_{\ulell;+-}(\xi_1,\xi_2)$ reveals that it can be written as a linear combination of terms of the form
\begin{equation*}
    \bigl( (\jxione + \ulell \xi_1) - (\jxitwo + \ulell \xi_2) \bigr) \fraka_{\ulell}(\xi_1) \frakb_{\ulell}(\xi_2) \check{\varphi}(\xi_1-\xi_2), \quad \fraka_\ulell, \frakb_\ulell \in W^{1,\infty}(\bbR), \quad \varphi \in \calS(\bbR).
\end{equation*}
In order to estimate the contribution of the resonant term $I_1(t)$, it therefore suffices to consider 
\begin{equation*}
    \begin{aligned}
    &\int_t^T \iint e^{is(\jxione + \ulell \xi_1 - \jxitwo - \ulell \xi_2)} \bigl( (\jxione + \ulell \xi_1) - (\jxitwo + \ulell \xi_2) \bigr) \bigl( \jxione^{-1} \fraka_\ulell(\xi_1) \gulellsh(s,\xi_1) \bigr) \\ 
    &\qquad \qquad \qquad \qquad \qquad \qquad \qquad \qquad \qquad \qquad \times \overline{\bigl( \jxitwo^{-1} \frakb_\ulell(\xi_2) \gulellsh(s,\xi_2) \bigr)} \, \check{\varphi}(\xi_1-\xi_2) \, \ud \xi_1 \, \ud \xi_2 \, \ud s.
    \end{aligned}
\end{equation*}
Integrating by parts in time, we obtain 
\begin{equation*}
    \begin{aligned}
        &\int_t^T \iint e^{is(\jxione + \ulell \xi_1 - \jxitwo - \ulell \xi_2)} \bigl( (\jxione + \ulell \xi_1) - (\jxitwo + \ulell \xi_2) \bigr) \bigl( \jxione^{-1} \fraka_\ulell(\xi_1) \gulellsh(s,\xi_1) \bigr) \\ 
        &\qquad \qquad \qquad \qquad \qquad \qquad \qquad \qquad \qquad \qquad \times \overline{\bigl( \jxitwo^{-1} \frakb_\ulell(\xi_2) \gulellsh(s,\xi_2) \bigr)} \, \check{\varphi}(\xi_1-\xi_2) \, \ud \xi_1 \, \ud \xi_2 \, \ud s \\ 
        &= -i \biggl[ \iint e^{is(\jxione + \ulell \xi_1 - \jxitwo - \ulell \xi_2)} \bigl( \jxione^{-1} \fraka_\ulell(\xi_1) \gulellsh(s,\xi_1) \bigr)  \overline{\bigl( \jxitwo^{-1} \frakb_\ulell(\xi_2) \gulellsh(s,\xi_2) \bigr)} \, \check{\varphi}(\xi_1-\xi_2) \, \ud \xi_1 \, \ud \xi_2 \biggr]_{s=t}^{s=T} \\
        &\quad + i \int_t^T \iint e^{is(\jxione + \ulell \xi_1 - \jxitwo - \ulell \xi_2)} \bigl( \jxione^{-1} \fraka_\ulell(\xi_1) \ps \gulellsh(s,\xi_1) \bigr)  \overline{\bigl( \jxitwo^{-1} \frakb_\ulell(\xi_2) \gulellsh(s,\xi_2) \bigr)} \\ 
        &\qquad \qquad \qquad \qquad \qquad \qquad \qquad \qquad \qquad \qquad \qquad \qquad \qquad \qquad \qquad \times \check{\varphi}(\xi_1-\xi_2) \, \ud \xi_1 \, \ud \xi_2 \, \ud s \\ 
        &\quad + i \int_t^T \iint e^{is(\jxione + \ulell \xi_1 - \jxitwo - \ulell \xi_2)} \bigl( \jxione^{-1} \fraka_\ulell(\xi_1) \gulellsh(s,\xi_1) \bigr)  \overline{\bigl( \jxitwo^{-1} \frakb_\ulell(\xi_2) \ps \gulellsh(s,\xi_2) \bigr)} \\
        &\qquad \qquad \qquad \qquad \qquad \qquad \qquad \qquad \qquad \qquad \qquad \qquad \qquad \qquad \qquad \times \check{\varphi}(\xi_1-\xi_2) \, \ud \xi_1 \, \ud \xi_2 \, \ud s \\ 
        &= I_{1,1}(t) + I_{1,2}(t) + I_{1,3}(t).
    \end{aligned}
\end{equation*}
In order to faciliate the analysis of the contributions of the terms $I_{1,k}(t)$, $1 \leq k \leq 3$, it is convenient to introduce the following auxiliary linear Klein-Gordon evolutions 
\begin{equation*}
\begin{aligned}
    v_1(s,y) &:= \frac{1}{\sqrt{2\pi}} \int_\bbR e^{iy\xi_1} e^{is(\jxione + \ulell \xi_1)} \jxione^{-1} \fraka_\ulell(\xi_1) \gulellsh(s,\xi_1) \, \ud \xi_1, \\
    v_2(s,y) &:= \frac{1}{\sqrt{2\pi}} \int_\bbR e^{iy\xi_2} e^{is(\jxitwo + \ulell \xi_2)} \jxitwo^{-1} \frakb_{\ulell}(\xi_2) \gulellsh(s,\xi_2) \, \ud \xi_2.
\end{aligned}
\end{equation*}
They satisfy all the decay estimates and bounds listed in item (11) in the statement of Corollary~\ref{cor:consequences_bootstrap_assumptions}.
The first term $I_{1,1}(t)$ is straightforward to estimate.
Rewriting $I_{1,1}(t)$ as 
\begin{equation*}
    \begin{aligned}
        I_{1,1}(t) = -i \biggl[ \sqrt{2\pi} \int_\bbR v_1(s,y) \overline{v_2(s,y)} \varphi(y) \, \ud y \biggr]_{s=t}^{s=T}, 
    \end{aligned}
\end{equation*}
we obtain from \eqref{equ:consequences_aux_KG_disp_decay} the sufficient decay estimate
\begin{equation*}
    \bigl| I_{1,1}(t) \bigr| \lesssim \sup_{t \leq s \leq T} \, \|v_1(s)\|_{L^\infty_y} \|v_2(s)\|_{L^\infty_y} \|\varphi\|_{L^1_y} \lesssim \sup_{t \leq s \leq T} \, \varepsilon^2 \js^{-1} \lesssim \varepsilon^2 \jt^{-1}.
\end{equation*}
We now turn to the term $I_{1,2}(t)$. The term $I_{1,3}(t)$ can be handled analogously, and we omit the details. Inserting the evolution equation \eqref{equ:setting_up_g_evol_equ3} for the effective profile $\gulellsh(s,\xi_1)$, we find 
\begin{equation}
    \begin{aligned}
        I_{1,2}(t) &= i \int_t^T (\dot{q}(s) - \ulell) \iint e^{is(\jxione + \ulell \xi_1 - \jxitwo - \ulell \xi_2)} i\xi_1 \bigl( \jxione^{-1} \fraka_\ulell(\xi_1) \gulellsh(s,\xi_1) \bigr) \\
        &\qquad \qquad \qquad \qquad \qquad \qquad \times \overline{\bigl( \jxitwo^{-1} \frakb_\ulell(\xi_2) \gulellsh(s,\xi_2) \bigr)} \, \check{\varphi}(\xi_1-\xi_2) \, \ud \xi_1 \, \ud \xi_2 \, \ud s \\
        &\quad - i \int_t^T \iint e^{-is(\jxitwo + \ulell \xi_2)} \jxione^{-1} \fraka_\ulell(\xi_1) \calFulellsh\bigl[ \calQ_{\ulell}\bigl(\usubeone(s)\bigr) \bigr](\xi_1) \\ 
        &\qquad \qquad \qquad \qquad \qquad \qquad \times \overline{\bigl( \jxitwo^{-1} \frakb_\ulell(\xi_2) \gulellsh(s,\xi_2) \bigr)} \, \check{\varphi}(\xi_1-\xi_2) \, \ud \xi_1 \, \ud \xi_2 \, \ud s \\
        &\quad - \frac{i}{6} \int_t^T \iint e^{-is(\jxitwo + \ulell \xi_2)} \jxione^{-1} \fraka_\ulell(\xi_1) \calFulellsh\bigl[ \usubeone(s)^3 \bigr](\xi_1) \\ 
        &\qquad \qquad \qquad \qquad \qquad \qquad \times \overline{\bigl( \jxitwo^{-1} \frakb_\ulell(\xi_2) \gulellsh(s,\xi_2) \bigr)} \, \check{\varphi}(\xi_1-\xi_2) \, \ud \xi_1 \, \ud \xi_2 \, \ud s \\ 
        &\quad + i \int_t^T \iint e^{-is(\jxitwo + \ulell \xi_2)} \jxione^{-1} \fraka_\ulell(\xi_1) \widetilde{\calR}(s,\xi_1) \\ 
        &\qquad \qquad \qquad \qquad \qquad \qquad \times \overline{\bigl( \jxitwo^{-1} \frakb_\ulell(\xi_2) \gulellsh(s,\xi_2) \bigr)} \, \check{\varphi}(\xi_1-\xi_2) \, \ud \xi_1 \, \ud \xi_2 \, \ud s \\   
        &=: I_{1,2}^{(a)}(t) + I_{1,2}^{(b)}(t) + I_{1,2}^{(c)}(t) + I_{1,2}^{(d)}(t).
    \end{aligned}
\end{equation}

We begin with the estimate for the term $I_{2,2}^{(a)}(t)$, which we can write as
\begin{equation*}
    \begin{aligned}
        I_{1,2}^{(a)}(t) = i \sqrt{2\pi} \int_t^T (\dot{q}(s) - \ulell) \int_\bbR (\py v_1)(s,y) \overline{v_2(s,y)} \varphi(y) \, \ud y \, \ud s.
    \end{aligned}
\end{equation*}
Thus, using \eqref{equ:consequences_qdot_minus_ulell_decay} and \eqref{equ:consequences_aux_KG_disp_decay}, we obtain the acceptable bound
\begin{equation*}
    \begin{aligned}
        \bigl| I_{1,2}^{(a)}(t) \bigr| &\lesssim \int_t^T |\dot{q}(s)-\ulell| \| (\py v_1)(s) \|_{L^\infty_y} \|v_2(s)\|_{L^\infty_y} \|\varphi\|_{L^1_y} \, \ud s \\ 
        &\lesssim \int_t^T \varepsilon \js^{-1+\delta} \cdot \varepsilon \js^{-\frac12} \cdot \varepsilon \js^{-\frac12} \, \ud s \lesssim \varepsilon^3 \jt^{-1+\delta}.
    \end{aligned}
\end{equation*}

We continue with the term $I_{1,2}^{(b)}(t)$. It can be written as
\begin{equation*}
    \begin{aligned}
        I_{1,2}^{(b)}(t) &= -i \int_t^T \int_\bbR \biggl( \int_\bbR e^{iy\xi_1} \jxione^{-1} \fraka_\ulell(\xi_1) \calFulellsh\bigl[ \calQ_\ulell(s)\bigr](\xi_1) \, \ud \xi_1 \biggr) \overline{v_2(s,y)} \varphi(y) \, \ud y \, \ud s.
    \end{aligned}
\end{equation*}
Inserting the leading order local decay decomposition \eqref{equ:consequences_decomposition_leading_order_local_decay} for $\usubeone(s)$ given in expanded form by
\begin{equation*}
    \usubeone(s,y) = \frac12 e_\ulell^{\#}(y, -\ulg \ulell) h_\ulell(s) + \frac12 \overline{e_\ulell^{\#}(y, -\ulg \ulell)} \overline{h_\ulell(s)} + \Remusubeone(s,y),
\end{equation*}
and inserting the leading order local decay decomposition \eqref{equ:consequences_aux_KG_leading_order_local_decay_decomp} for $v_2(s,y)$ given by
\begin{equation*}
    v_2(s,y) = e^{-i\ulg\ulell y} \tilde{h}_\ulell(s) + R_{v_2}(s,y),
\end{equation*}
we find that $I_{1,2}^{(b)}(t)$ is a linear combination of two types of terms. 

The first type of term is of the form 
\begin{equation} \label{equ:modulation_proof_resonant_term_I1_type1}
    C_{\iota_1 \iota_2 -} \cdot \int_t^T H_{\iota_1 \iota_2 -}(s) \, \ud s
\end{equation}
for $\iota_1, \iota_2 \in \{\pm\}$ with 
\begin{equation*}  
    H_{\iota_1 \iota_2 -}(s) := h_\ulell^{\iota_1}(s) h_\ulell^{\iota_2}(s) \overline{\tilde{h}_\ulell(s)} 
\end{equation*}
and 
\begin{equation*}
    \begin{aligned}
        C_{\iota_1 \iota_2 -} &:= \int_\bbR \biggl( \int_\bbR e^{iy\xi_1} \jxione^{-1} \fraka_\ulell(\xi_1) \calFulellsh\Bigl[ \alpha(\ulg \cdot) e_\ulell^{\#,\iota_1}(\cdot, -\ulg\ulell) e_\ulell^{\#,\iota_2}(\cdot, -\ulg\ulell) \Bigr](\xi_1) \, \ud \xi_1 \biggr) e^{i \ulg \ulell y} \varphi(y) \, \ud y,
    \end{aligned}
\end{equation*}
where we use the short-hand notations
    \begin{equation}
        \begin{aligned}
        &h^{+}_\ulell(s) := h_\ulell(s), \quad h^{-}_\ulell(s) := \overline{h_\ulell(s)}, \\
        &e_\ulell^{\#, +}(\cdot, -\ulg \ulell) := e_\ulell^{\#}(\cdot, -\ulg \ulell), \quad e_\ulell^{\#, -}(\cdot, -\ulg \ulell) := \overline{e_\ulell^{\#}(\cdot, -\ulg \ulell)}.
        \end{aligned}
    \end{equation}
One checks that $|C_{\iota_1 \iota_2 -}| \lesssim_{\ell_0} 1$ in view of the rapid decay of $\alpha(\ulg \cdot)$ and $\varphi(\cdot)$.
Now the key observation is that by a formal stationary phase analysis we have to leading order that $h_\ulell^{\pm}(s) \sim s^{-\frac12} e^{\pm i s \ulg^{-1}}$ as well as $\overline{\tilde{h}(s)} \sim s^{-\frac12} e^{-is\ulg^{-1}}$. Hence, to leading order any of the terms $H_{\iota_1 \iota_2 -}(s)$, $\iota_1, \iota_2 \in \{\pm\}$, is of the form
\begin{equation} \label{equ:modulation_proof_resonant_term_I1_cubic_oscillations}
    H_{\iota_1 \iota_2 -}(s) \sim s^{-\frac32} e^{ims\ulg^{-1}}, \quad m \in \{-3,-1,1\}.
\end{equation}
The decay rate $s^{-\frac32}$ is not enough to obtain upon integration the asserted decay estimate \eqref{equ:modulation_proof_goal}. Fortunately, the phase does not vanish in the leading order behavior \eqref{equ:modulation_proof_resonant_term_I1_cubic_oscillations} of any of the terms $H_{\iota_1 \iota_2 -}(s)$. We can therefore integrate by parts in time in \eqref{equ:modulation_proof_resonant_term_I1_type1}, which then suffices to deduce the desired decay estimate \eqref{equ:modulation_proof_goal}. We present the details for the concrete case $\iota_1 = \iota_2 = +$. All other cases can be handled analogously.
Integrating by parts in time, we find 
\begin{equation*}
    \begin{aligned}
        \int_t^T H_{++-}(s) \, \ud s &= \int_t^T e^{is\ulg^{-1}} \bigl( e^{-is\ulg^{-1}} h_\ulell(s) \bigr)^2 \overline{\bigl( e^{-is\ulg^{-1}} \tilde{h}_\ulell(s) \bigr)} \, \ud s \\ 
        &= \biggl[ -i\ulg e^{is\ulg^{-1}} \bigl( e^{-is\ulg^{-1}} h_\ulell(s) \bigr)^2 \overline{\bigl( e^{-is\ulg^{-1}} \tilde{h}_\ulell(s) \bigr)} \biggr]_{s=t}^{s=T} \\ 
        &\quad + 2 i \ulg \int_t^T e^{is\ulg^{-1}} \bigl( e^{-is\ulg^{-1}} h_\ulell(s) \bigr) \cdot \ps \bigl( e^{-is\ulg^{-1}} h_\ulell(s) \bigr) \cdot \overline{\bigl( e^{-is\ulg^{-1}} \tilde{h}_\ulell(s) \bigr)} \, \ud s + \bigl\{ \text{similar} \bigr\}.
        &\quad 
    \end{aligned}
\end{equation*}
Thus, invoking the decay estimates \eqref{equ:consequences_hulell_decay}, \eqref{equ:consequences_hulell_phase_filtered_decay}, \eqref{equ:consequences_aux_KG_tildeh_ulell_decay}, \eqref{equ:consequences_aux_KG_tildeh_ulell_phase_filtered_decay}, we arrive at the acceptable bound
\begin{equation*}
    \begin{aligned}
        \biggl| \int_t^T H_{++-}(s) \, \ud s \biggr| &\lesssim \sup_{t \leq s \leq T} \, |h_\ulell(s)|^2 |\tilde{h}_\ulell(s)| + \int_t^T |h_\ulell(s)| \cdot \bigl| \ps \bigl( e^{-is\ulg^{-1}} h_\ulell(s) \bigr) \bigr| \cdot |\tilde{h}_\ulell(s)| \, \ud s + \ldots \\ 
        &\lesssim \sup_{t \leq s \leq T} \, \varepsilon^3 \js^{-\frac32} + \int_t^T \varepsilon \js^{-\frac12} \cdot \varepsilon \js^{-1+\delta} \cdot \varepsilon \js^{-\frac12} \, \ud s 
        \lesssim \varepsilon \jt^{-1+\delta}.
    \end{aligned}
\end{equation*}

The second type of term is of the form 
\begin{equation*}
    \begin{aligned}
        \int_t^T \int_\bbR \biggl( \int_\bbR e^{iy\xi_1} \jxione^{-1} \fraka_\ulell(\xi_1) \calFulellsh\Bigl[ \alpha(\ulg \cdot) U_1(s,\cdot) U_2(s,\cdot) \Bigr](\xi_1) \, \ud \xi_1 \biggr) \overline{V_2(s,y)} \varphi(y) \, \ud y \, \ud s,
    \end{aligned}
\end{equation*}
where $U_1(s,y), U_2(s,y)$ are given by $\eulsharp(y,-\ulg\ulell) h_\ulell(s)$ or $\Remusubeone(s,y)$ or complex conjugates thereof, while $V_2(s,y)$ is given by $e^{-i\ulg\ulell y} \tilde{h}_\ulell(s)$ or $R_{v_2}(s,y)$. Importantly, either at least one of the inputs $U_1(s,y), U_2(s,y)$ must be $\Remusubeone(s,y)$ (or its complex conjugate) or $V_2(s,y)$ must be given by $R_{v_2}(s,y)$. 
Exploiting the improved local decay \eqref{equ:consequences_Remusubeone_H1y_local_decay} of the remainder term $\Remusubeone(s,y)$, respectively the improved local decay \eqref{equ:consequences_aux_KG_Remv_improved_H1y_local_decay} of the remainder term $R_{v_2}(s,y)$, then gives an acceptable bound upon integrating in time. We provide the details for the concrete case $U_1(s,y) = U_2(s,y) = \eulsharp(y,-\ulg\ulell) h_\ulell(s)$ and $V_2(s,y) = R_{v_2}(s,y)$. All other cases can be dealt with analogously.
By H\"older's inequality, Plancherel's theorem, the mapping property~\eqref{equ:mapping_property_calFulellsh_L2}, and the decay estimates \eqref{equ:consequences_hulell_decay}, \eqref{equ:consequences_aux_KG_Remv_improved_H1y_local_decay}, we find that
\begin{equation*}
    \begin{aligned}
        &\biggl| \int_t^T \int_\bbR \biggl( \int_\bbR e^{iy\xi_1} \jxione^{-1} \fraka_\ulell(\xi_1) \calFulellsh\Bigl[ \alpha(\ulg \cdot) \eulsharp(\cdot,-\ulg\ulell) h_\ulell(s) \eulsharp(\cdot,-\ulg\ulell) h_\ulell(s) \Bigr](\xi_1) \, \ud \xi_1 \biggr) \overline{R_{v_2}(s,y)} \varphi(y) \, \ud y \, \ud s \biggr| \\ 
        &\lesssim \int_t^T \biggl\| \int_\bbR e^{iy\xi_1} \jxione^{-1} \fraka_\ulell(\xi_1) \calFulellsh\Bigl[ \alpha(\ulg \cdot) \eulsharp(\cdot,-\ulg\ulell) h_\ulell(s) \eulsharp(\cdot,-\ulg\ulell) h_\ulell(s) \Bigr](\xi_1) \, \ud \xi_1 \biggr\|_{L^2_y} \\ 
        &\qquad \qquad \qquad \qquad \qquad \qquad \qquad \qquad \qquad \qquad \qquad \qquad \qquad \times \bigl\| \jy^{-3} R_{v_2}(s) \bigr\|_{L^2_y} \bigl\| \jy^3 \varphi \bigr\|_{L^\infty_y} \, \ud s \\
        &\lesssim \int_t^T \Bigl\| \alpha(\ulg \cdot) \eulsharp(\cdot,-\ulg\ulell) \eulsharp(\cdot,-\ulg\ulell) \Bigr\|_{L^2_y} |h_\ulell(s)|^2 \bigl\| \jy^{-3} R_{v_2}(s) \bigr\|_{L^2_y} \bigl\| \jy^3 \varphi \bigr\|_{L^\infty_y} \, \ud s \\
        &\lesssim \int_t^T \varepsilon^2 \js^{-1} \cdot \varepsilon \js^{-1+\delta} \, \ud s \lesssim \varepsilon^3 \jt^{-1+\delta},
    \end{aligned}
\end{equation*}
as desired.

Next, we turn to the term $I_{1,2}^{(c)}(t)$, which we write as
\begin{equation} \label{equ:modulation_proof_I12c_rewritten}
    \begin{aligned}
        I_{1,2}^{(c)}(t) &= - \frac{i}{6} \int_t^T \int_\bbR \biggl( \int_\bbR e^{iy\xi_1} \jxione^{-1} \fraka_\ulell(\xi_1) \calFulellsh\bigl[ \usubeone(s)^3 \bigr](\xi_1) \, \ud \xi_1 \biggr) \overline{v_2(s,y)} \varphi(y) \, \ud y \, \ud s.
    \end{aligned}
\end{equation}
Here we need to insert the fine structure of the distorted Fourier transform of the cubic interactions $\calFulellsh\bigl[ \usubeone(s)^3 \bigr](\xi_1)$ determined in Subsection~\ref{subsec:cubic_spectral_distributions} and Subsection~\ref{subsec:structure_cubic_nonlinearities}. By \eqref{equ:setting_up_dist_FT_of_cubic_expanded} we have 
\begin{equation*}
        \calFulellsh\bigl[ \usubeone(s)^3 \bigr](\xi_1) 
        = e^{i s (\jxione + \ulell \xi_1)} \biggl( - \frac{i}{8} \calI_1(s,\xi_1) + \frac{3i}{8} \calI_2(s,\xi_1) - \frac{3i}{8} \calI_3(s,\xi_1) + \frac{i}{8} \calI_4(s,\xi_1) \biggr).
\end{equation*}
We only present the details for the contributions of $\calI_2(s,\xi_1)$, which concerns the $+-+$ cubic interactions. The contributions of all other cubic interactions can be dealt with analogously. By \eqref{equ:calIj_decomposition} we have the following refined decomposition into cubic interactions with a Dirac kernel, cubic interactions with a Hilbert-type kernel, and regular cubic interactions,
\begin{equation*}
    e^{i s (\jxione + \ulell \xi_1)} \calI_2(s,\xi_1) = e^{i s (\jxione + \ulell \xi_1)} \Bigl( \calI_2^{\delta_0}(s,\xi_1) + \calI_2^{\pvdots}(s,\xi_1) + \calI_2^{\mathrm{reg}}(s,\xi_1) \Bigr).
\end{equation*}
We consider these terms separately. In view of the fine structure of $\frakm_{\ulell,+-+}^{\delta_0}$ determined in Subsection~\ref{subsec:cubic_spectral_distributions}, the contributions of the cubic interactions with a Dirac kernel $e^{i s (\jxione + \ulell \xi_1)} \calI_2^{\delta_0}(s,\xi_1)$ to \eqref{equ:modulation_proof_I12c_rewritten} are linear combinations of terms of the form
\begin{equation*}
    \begin{aligned}
        I_{1,2}^{(c),+-+,\delta_0}(t) := \int_t^T \biggl( \int_\bbR e^{iy\xi_1} \frakc(\xi_1) \widehat{\calF}\Bigl[ \tilde{v}_1(s) \overline{\tilde{v}_2(s)} \tilde{v}_3(s) \Bigr](\xi_1) \, \ud \xi_1 \biggr) \overline{v_2(s,y)} \varphi(y) \, \ud y
    \end{aligned}
\end{equation*}
with $\frakc(\xi_1) := \jxione^{-1} \fraka_\ulell(\xi_1) \frakb(\xi_1)$ and
\begin{equation} \label{equ:modulation_proof_tildev_definition}
    \tilde{v}_j(s,y) := \frac{1}{\sqrt{2\pi}} \int_\bbR e^{iy\eta} e^{is(\jap{\eta}+\ulell\eta)} \jeta^{-1} \frakb_j(\eta) \gulellsh(s,\eta) \, \ud \eta, \quad 1 \leq j \leq 3,
\end{equation}
for some $\frakb, \frakb_1, \frakb_2, \frakb_3 \in W^{1,\infty}(\bbR)$. 
We note that the linear Klein-Gordon evolutions $\tilde{v}_j(s,y)$, $1 \leq j \leq 3$, satisfy all decay estimates and bounds listed in item (11) of Corollary~\ref{cor:consequences_bootstrap_assumptions}.
By the regularity assumptions, both $\frakc(\xi_1)$ and $\partial_{\xi_1} \frakc(\xi_1)$ are $L^2_{\xi_1}$-integrable, whence $\widehat{\calF}^{-1}[\frakc](y)$ is $L^1_y$ integrable. Hence, by H\"older's inequality, Young's inequality, and \eqref{equ:consequences_aux_KG_disp_decay}, we obtain the sufficient bound
\begin{equation*}
    \begin{aligned}
        \bigl| I_{1,2}^{(c),+-+,\delta_0}(t) \bigr| &\lesssim \int_t^T \biggl\| \int_\bbR e^{iy\xi_1} \frakc(\xi_1) \widehat{\calF}\Bigl[ \tilde{v}_1(s) \overline{\tilde{v}_2(s)} \tilde{v}_3(s) \Bigr](\xi_1) \, \ud \xi_1 \biggr\|_{L^\infty_y} \|v_2(s)\|_{L^\infty_y} \|\varphi\|_{L^1_y} \, \ud s \\ 
        &\lesssim \int_t^T \bigl\| \widehat{\calF}^{-1}[\frakc] \bigr\|_{L^1_y} \|\tilde{v}_1(s)\|_{L^\infty_y} \|\tilde{v}_2(s)\|_{L^\infty_y} \|\tilde{v}_3(s)\|_{L^\infty_y} \|v_2(s)\|_{L^\infty_y} \, \ud s \\
        &\lesssim \int_t^T \varepsilon^4 \js^{-2} \, \ud s \lesssim \varepsilon^4 \jt^{-1}.
    \end{aligned}
\end{equation*}

Similarly, in view of the fine structure of $\frakm_{\ulell,+-+}^{\pvdots}$ determined in Subsection~\ref{subsec:cubic_spectral_distributions}, the contributions of the cubic interactions with a Hilbert-type kernel $e^{i s (\jxione + \ulell \xi_1)} \calI_2^{\pvdots}(s,\xi_1)$ to \eqref{equ:modulation_proof_I12c_rewritten} are linear combinations of terms of the form
\begin{equation*}
    \begin{aligned}
        I_{1,2}^{(c),+-+,\pvdots}(t) := \int_t^T \biggl( \int_\bbR e^{iy\xi_1} \frakc(\xi_1) \widehat{\calF}\Bigl[ \tilde{v}_1(s) \overline{\tilde{v}_2(s)} \tilde{v}_3(s) \tanh(\ulg \cdot) \Bigr](\xi_1) \, \ud \xi_1 \biggr) \overline{v_2(s,y)} \varphi(y) \, \ud y
    \end{aligned}
\end{equation*}
with the same definitions for $\frakc(\xi_1)$ and $\tilde{v}_j(s,y)$, $1 \leq j \leq 3$ as before. 
Thus, we obtain the same sufficient decay estimate for the term $I_{1,2}^{(c),+-+,\pvdots}(t)$ as for the contributions of the cubic interactions with a Dirac kernel. We omit further details.

Finally, in view of the structure of $\nu_{\ulell,+-+}^{\reg}$ the contributions coming from the regular cubic interactions are of the form 
\begin{equation*}
    \begin{aligned}
        I_{1,2}^{(c),+-+,\reg}(t) := \int_t^T \biggl( \int_\bbR e^{iy\xi_1} \jxione^{-1} \fraka_\ulell(\xi_1) \frakb(\xi_1) \widehat{\calF}\Bigl[ \tilde{v}_1(s) \overline{\tilde{v}_2(s)} \tilde{v}_3(s) \widetilde{\varphi}(\cdot) \Bigr](\xi_1) \, \ud \xi_1 \biggr) \overline{v_2(s,y)} \varphi(y) \, \ud y,
    \end{aligned}
\end{equation*}
where $\widetilde{\varphi} \in \calS(\bbR)$ is some Schwartz function, $\frakb \in W^{1,\infty}(\bbR)$, and the linear Klein-Gordon evolutions $\tilde{v}_j(s)$ are defined as in \eqref{equ:modulation_proof_tildev_definition}.
Thus, by H\"older's inequality and repeated application of Plancherel's theorem we arrive at the acceptable bound
\begin{equation*}
    \begin{aligned}
        &\bigl| I_{1,2}^{(c),+-+,\mathrm{reg}}(t) \bigr| \\
        &\lesssim \int_t^T \biggl\| \biggl( \int_\bbR e^{iy\xi_1} \jxione^{-1} \fraka_\ulell(\xi_1) \frakb(\xi_1) \widehat{\calF}\Bigl[ \tilde{v}_1(s) \overline{\tilde{v}_2(s)} \tilde{v}_3(s) \widetilde{\varphi}(\cdot) \Bigr](\xi_1) \, \ud \xi_1 \biggr) \biggr\|_{L^2_y} \|v_2(s)\|_{L^\infty_y} \|\varphi\|_{L^2_y} \, \ud y \\
        &\lesssim \int_t^T \bigl\| \tilde{v}_1(s) \overline{\tilde{v}_2(s)} \tilde{v}_3(s) \widetilde{\varphi}(\cdot) \bigr\|_{L^2_y} \|v_2(s)\|_{L^\infty_y} \|\varphi\|_{L^2_y} \, \ud s \\
        &\lesssim \int_t^T \|\tilde{v}_1(s)\|_{L^\infty_y} \|\tilde{v}_2(s)\|_{L^\infty_y} \|\tilde{v}_3(s)\|_{L^\infty_y} \|\widetilde{\varphi}\|_{L^2_y} \|v_2(s)\|_{L^\infty_y} \|\varphi\|_{L^2_y} \, \ud s \lesssim \int_t^T \varepsilon^4 \js^{-2} \, \ud s \lesssim \varepsilon^4 \jt^{-1}.
    \end{aligned}
\end{equation*}
This concludes the treatment of the term $I_{1,2}^{(c)}(t)$. It remains to estimate the term $I_{1,2}^{(d)}(t)$. We write it as
\begin{equation*}
    \begin{aligned}
        I_{1,2}^{(d)}(t) &= i \int_t^T \biggl( \int_\bbR e^{iy\xi_1} \jxione^{-1} \fraka_\ulell(\xi_1) \widetilde{\calR}(s,\xi_1) \, \ud \xi_1 \biggr) \overline{v_2(s,y)} \varphi(y) \, \ud y \, \ud s.
    \end{aligned}
\end{equation*}
Using H\"older's inequality and Plancherel's theorem combined with the bounds \eqref{equ:consequences_wtilcalR_L2xi_bound} and \eqref{equ:consequences_aux_KG_disp_decay}, we arrive at the acceptable bound 
\begin{equation*}
\begin{aligned}
    |I_{1,2}^{(d)}(t)| &\lesssim \int_t^T \biggl\| \int_\bbR e^{iy\xi_1} \jxione^{-1} \fraka_\ulell(\xi_1) \widetilde{\calR}(s,\xi_1) \, \ud \xi_1 \biggr\|_{L^2_y} \|v_2(s)\|_{L^\infty_y} \|\varphi\|_{L^2_y} \, \ud s \\ 
    &\lesssim \int_t^T \bigl\| \widetilde{\calR}(s,\xi_1) \bigr\|_{L^2_{\xi_1}} \|v_2(s)\|_{L^\infty_y} \, \ud s  
    \lesssim \int_t^T \varepsilon^2 \js^{-\frac32+\delta} \cdot \varepsilon \js^{-\frac12} \, \ud s \lesssim \varepsilon^3 \jt^{-1+\delta}.
\end{aligned}    
\end{equation*}
Combining all of the preceding estimates yields an acceptable decay estimate for the contribution of the resonant term $I_1(t)$.

\medskip 

\noindent \underline{Contributions of the non-resonant terms $I_2(t)$ and $I_3(t)$.}
In Lemma~\ref{lem:null_structure2} we determined that the quadratic spectral distributions $\mu_{\ulell;++}(\xi_1,\xi_2)$ and $\mu_{\ulell;--}(\xi_1,\xi_2)$ are multiples of $(\jxione +\ulell\xi_1 + \jxitwo +\ulell\xi_2)$, which are the phases in $I_2(t)$ and $I_3(t)$. While these structures are not necessary for the analysis since these phases do not have time resonances, nonetheless it means that the terms $I_2(t)$, $I_3(t)$ can be treated exactly like the term $I_1(t)$ and we can omit any further details.

\medskip 

    \noindent \underline{Contribution of the term $II(s)$.}
    The contribution of the (first component of the) term $II(s)$ to the integral~\eqref{equ:modulation_proof_ellt_minus_ellT} is
    \begin{equation*}
        - \int_t^t \ulg^{-3} \|K'\|_{L^2}^{-2} \bigl\langle \bigl( \bmY_{0, \ulell} \bigr)_1, \calC\bigl( \usubeone(s) \bigr) \bigr\rangle \, \ud s = \frac16 \ulg^{-2} \|K'\|_{L^2}^{-2} \int_t^T \bigl\langle K'(\ulg \cdot) W^{(4)}\bigl(K(\ulg \cdot)\bigr), \usubeone(s)^3 \bigr\rangle \, \ud s.
    \end{equation*}
    We proceed similarly to the preceding treatment of the term $I_{1,2}^{(b)}(t)$. 
    Inserting the leading order local decay decomposition \eqref{equ:consequences_decomposition_leading_order_local_decay} for $\usubeone(s)$ given in expanded form by
    \begin{equation*}
        \usubeone(s,y) = \frac12 e_\ulell^{\#}(y, -\ulg \ulell) h_\ulell(s) + \frac12 \overline{e_\ulell^{\#}(y, -\ulg \ulell)} \overline{h_\ulell(s)} + \Remusubeone(s,y),
    \end{equation*}
    we obtain a linear combination of two types of terms.
    The first type is of the form
    \begin{equation} \label{equ:modulation_proof_contribution_II_type1}
        C_{\iota_1 \iota_2 \iota_3} \int_t^T H_{\iota_1 \iota_2 \iota_3}(s) \, \ud s
    \end{equation}
    for $\iota_1, \iota_2, \iota_3 \in \{\pm\}$ with 
    \begin{equation*}
        H_{\iota_1 \iota_2 \iota_3}(s) := h_\ulell^{\iota_1}(s) h_{\ulell}^{\iota_2}(s) h_{\ulell}^{\iota_3}(s)
    \end{equation*}
    and 
    \begin{equation*}
        C_{\iota_1 \iota_2 \iota_3} := \Bigl\langle K'(\ulg \cdot) W^{(4)}\bigl(K(\ulg \cdot)\bigr), e_\ulell^{\#, \iota_1}(\cdot, -\ulg \ulell) e_\ulell^{\#, \iota_2}(\cdot, -\ulg \ulell) e_\ulell^{\#, \iota_3}(\cdot, -\ulg \ulell) \Bigr\rangle,
    \end{equation*}
    where we use the short-hand notations 
    \begin{equation} \label{equ:modulation_proof_term_two_shorthand_def}
        \begin{aligned}
        &h^{+}_\ulell(s) := h_\ulell(s), \quad h^{-}_\ulell(s) := \overline{h_\ulell(s)}, \\
        &e_\ulell^{\#, +}(\cdot, -\ulg \ulell) := e_\ulell^{\#}(\cdot, -\ulg \ulell), \quad e_\ulell^{\#, -}(\cdot, -\ulg \ulell) := \overline{e_\ulell^{\#}(\cdot, -\ulg \ulell)}.
        \end{aligned}
    \end{equation}
    By a formal stationary phase analysis we find that to leading order $h_\ulell^{\pm}(s) \sim s^{-\frac12} e^{\pm i s \ulg^{-1}}$, whence to leading order any of the terms $H_{\iota_1 \iota_2 \iota_3}(s)$, $\iota_1, \iota_2, \iota_3 \in \{\pm\}$, is of the form
    \begin{equation} \label{equ:modulation_proof_Hiotas_leading_order_behavior}
        H_{\iota_1 \iota_2 \iota_3}(s) \sim s^{-\frac32} e^{i m s \ulg^{-1}}, \quad m \in \{-3, -1, 1, 3\}.
    \end{equation}
    The decay rate $s^{-\frac32}$ is not enough to recover upon integration the asserted decay estimate \eqref{equ:modulation_proof_goal}.
    But since the phase does not vanish in the leading order behavior \eqref{equ:modulation_proof_Hiotas_leading_order_behavior} of any of the terms $H_{\iota_1 \iota_2 \iota_3}(s)$, we can integrate by parts in time in \eqref{equ:modulation_proof_contribution_II_type1}, which then suffices to deduce the desired decay estimate \eqref{equ:modulation_proof_goal}.
    We present the details for the concrete case $\iota_1 = \iota_2 = \iota_3 = +$. All other cases can be handled analogously.
    Integrating by parts in time, we find 
    \begin{equation*}
        \begin{aligned}
            \int_t^T H_{+++}(s) \, \ud s &= \int_t^T e^{3is\ulg^{-1}} \bigl( e^{-is\ulg^{-1}} h_\ulell(s) \bigr)^3 \, \ud s \\ 
            &= -\frac{i\ulg}{3} \Bigl[ h_\ulell(s)^3 \Bigr]_{s=t}^{s=T} + i\ulg \int_t^T e^{3is\ulg^{-1}} \cdot \ps \bigl( e^{-is\ulg^{-1}} h_\ulell(s) \bigr) \cdot \bigl( e^{-is\ulg^{-1}} h_\ulell(s) \bigr)^2 \, \ud s. 
        \end{aligned}
    \end{equation*}
    Using \eqref{equ:consequences_hulell_decay} and \eqref{equ:consequences_hulell_phase_filtered_decay}, we arrive at the sufficient decay estimate
    \begin{equation*}
        \begin{aligned}
            \biggl| \int_t^T H_{+++}(s) \, \ud s \biggr| &\lesssim \sup_{t \leq s \leq T} \, \bigl| h_\ulell(s) \bigr|^3 + \int_t^T \bigl| \ps \bigl( e^{-is\ulg^{-1}} h_\ulell(s) \bigr) \bigr| \cdot |h_\ulell(s)|^2 \, \ud s \\ 
            &\lesssim \sup_{t \leq s \leq T} \, \varepsilon^3 \js^{-\frac32} + \int_t^T \varepsilon \js^{-1+\delta} \cdot \varepsilon^2 \js^{-1} \, \ud s \lesssim \varepsilon^3 \jt^{-1+\delta}.
        \end{aligned}
    \end{equation*}

    The second type is of the schematic form     
    \begin{equation} \label{equ:modulation_proof_contributionII_type2}
        \int_t^T \bigl\langle \varphi(\cdot), U_1(s,\cdot) U_2(s,\cdot) U_3(s,\cdot) \bigr\rangle \, \ud s,
    \end{equation}
    where $\varphi(y)$ is a rapidly decaying function and where $U_j(s,y)$, $1 \leq j \leq 3$, are given by $e_\ulell^{\#}(y, -\ulg \ulell) h_\ulell(s)$, $\overline{e_\ulell^{\#}(y, -\ulg \ulell) h_\ulell(s)}$, or $\Remusubeone(s,y)$, with at least one of the inputs given by $\Remusubeone(s,y)$. 
    For this second type we can deduce sufficient decay in time more crudely by exploiting the spatial localization of $\varphi(y)$, the improved local decay estimate~\eqref{equ:consequences_Remusubeone_H1y_local_decay} for $\Remusubeone(s,y)$ along with the bound \eqref{equ:consequences_hulell_decay}. 
    Specifically, using also Sobolev embedding if necessary, we obtain
    \begin{equation*}
        \begin{aligned}
             \biggl| \int_t^T \bigl\langle \varphi(\cdot), U_1(s,\cdot) U_2(s,\cdot) U_3(s,\cdot) \bigr\rangle \, \ud s \biggr| 
             &\lesssim \int_t^T \Bigl( |h_\ulell(s)| + \bigl\| \jy^{-3} \Remusubeone(s) \bigr\|_{H^1_y} \Bigr)^2 \bigl\| \jy^{-3} \Remusubeone(s) \bigr\|_{L^2_y} \, \ud s \\ 
             &\lesssim \int_t^T \varepsilon^2 \js^{-1} \cdot \varepsilon \js^{-1+\delta} \, \ud s \lesssim \varepsilon^3 \jt^{-1+\delta},
        \end{aligned}
    \end{equation*}
    as desired.

    \medskip 
    
    \noindent \underline{Contribution of the term $III(s)$.}
        The third term $III(s)$ can be treated similarly to the preceding term $II(s)$ as well as to the term $I_{1,2}^{(b)}(t)$.
        Its contribution to the integral \eqref{equ:modulation_proof_ellt_minus_ellT} is the first component of the vector quantity
        \begin{equation} \label{equ:modulation_proof_integral_of_third_term}
         \int_t^T \bbM_\ulell^{-1} \bbA_\ulell(s) \bbM_\ulell^{-1}    \begin{bmatrix}
                                                                \bigl\langle \bigl( \bmY_{0, \ulell} \bigr)_1, \calQ_\ulell\bigl( \usubeone(s) \bigr) \bigr\rangle \\ 
                                                                -\bigl\langle \bigl( \bmY_{1, \ulell} \bigr)_1, \calQ_\ulell\bigl( \usubeone(s) \bigr) \bigr\rangle 
                                                                \end{bmatrix} 
         \ud s, \\
       \end{equation}
       where we recall that
       \begin{equation*}
            \bbA_\ulell(s) 
                    := \begin{bmatrix} 
                            \bigl\langle \bfJ \bmZ_{3,\ulell}, (\ulPe \bmu)(s) \bigr\rangle &  \bigl\langle \bfJ \bmZ_{1,\ulell}, (\ulPe \bmu)(s) \bigr\rangle \\
                            \bigl\langle \bfJ \bmZ_{2,\ulell}, (\ulPe \bmu)(s) \bigr\rangle & - \bigl\langle \bfJ \bmZ_{3,\ulell}, (\ulPe \bmu)(s) \bigr\rangle
                        \end{bmatrix}.
       \end{equation*}
        In view of the rapid decay of the functions $\bmZ_{j,\ulell}$, $1 \leq j \leq 3$, the entries of the matrix $\bbA_\ulell(s)$ are a linear combination of terms $\langle \varphi_1(\cdot), u_{\mathrm{e},k}(s,\cdot) \rangle$ for $k \in \{1,2\}$ and for some Schwartz function $\varphi_1$, while the entries $\bigl\langle \bigl( \bmY_{k, \ulell} \bigr)_1(\cdot), \calQ_\ulell\bigl( \usubeone(s) \bigr) \bigr\rangle$, $k \in \{0,1\}$, of the $2$-vector in the integrand in \eqref{equ:modulation_proof_integral_of_third_term} are of the form $\langle \varphi_2(\cdot), \usubeone(s)^2 \rangle$ for some Schwartz function $\varphi_2$. 
        
        Inserting the leading order local decay decompositions \eqref{equ:consequences_decomposition_leading_order_local_decay} and \eqref{equ:consequences_decomposition_leading_order_local_decay_usubetwo} for $\usubeone(s)$, respectively for $\usubetwo(s)$, we obtain a linear combination of two types of terms. The first type of term is of the form
        \begin{equation} \label{equ:modulation_proof_third_term_first_type}
            D_{\iota_1 \iota_2 \iota_3} \int_t^T E_{\iota_1 \iota_2 \iota_3}(s) \, \ud s \quad \quad \text{or} \quad \quad F_{\iota_1 \iota_2 \iota_3} \int_t^T G_{\iota_1 \iota_2 \iota_3}(s) \, \ud s
        \end{equation}
        for $\iota_1, \iota_2, \iota_3 \in \{\pm\}$ with
        \begin{equation*}
            \begin{aligned}
            D_{\iota_1 \iota_2 \iota_3} &:= \bigl\langle \varphi_1(\cdot), e_\ulell^{\#,\iota_1}(\cdot, -\ulg\ulell) \bigr\rangle \, \bigl\langle \varphi_2(\cdot), e_\ulell^{\#,\iota_2}(\cdot, -\ulg\ulell) e_\ulell^{\#,\iota_3}(\cdot, -\ulg\ulell) \bigr\rangle, \\  
            E_{\iota_1 \iota_2 \iota_3}(s) &:= h_\ulell^{\iota_1}(s) h_\ulell^{\iota_2}(s) h_\ulell^{\iota_3}(s),
            \end{aligned}
        \end{equation*}
        and 
        \begin{equation*}
            \begin{aligned}            
            F_{\iota_1 \iota_2 \iota_3} &:= \bigl\langle \varphi_1(\cdot), \ulg^{-1} \bigl( D^\ast e_\ulell^{\#} \bigr)^{\iota_1}(\cdot,-\ulg\ulell) \bigr\rangle \, \bigl\langle \varphi_2(\cdot), e_\ulell^{\#,\iota_2}(\cdot, -\ulg\ulell) e_\ulell^{\#,\iota_3}(\cdot, -\ulg\ulell) \bigr\rangle, \\
            G_{\iota_1 \iota_2 \iota_3}(s) &:= p_\ulell^{\iota_1}(s) h_\ulell^{\iota_2}(s) h_\ulell^{\iota_3}(s),
            \end{aligned}
        \end{equation*}
        where we use the short-hand notations \eqref{equ:modulation_proof_term_two_shorthand_def} as well as 
        \begin{equation*}
            \begin{aligned}
                &p_\ulell^+(s) := p_\ulell(s), \quad p_\ulell^{-}(s) := \overline{p_\ulell(s)}, \\
                &\bigl( D^\ast e_\ulell^{\#} \bigr)^{+}(y,-\ulg\ulell) := \bigl( D^\ast e_\ulell^{\#} \bigr)(y,-\ulg\ulell), \quad \bigl( D^\ast e_\ulell^{\#} \bigr)^{-}(y,-\ulg\ulell) := \overline{\bigl( D^\ast e_\ulell^{\#} \bigr)(y,-\ulg\ulell)}.
            \end{aligned}
        \end{equation*}
        In view of the oscillatory behavior of $h_\ulell(s)$ and $p_\ulell(s)$ captured in \eqref{equ:consequences_hulell_phase_filtered_decay}, respectively in \eqref{equ:consequences_pulell_phase_filtered_decay}, together with the decay estimates \eqref{equ:consequences_hulell_decay}, \eqref{equ:consequences_pulell_decay}, it is clear that the contributions of the first type of terms \eqref{equ:modulation_proof_third_term_first_type} can be estimated analogously to the first type of terms \eqref{equ:modulation_proof_contribution_II_type1} among the contributions of the term $II(s)$. 
        We omit the details.

        The second type of term is of the form
        \begin{equation} \label{equ:modulation_proof_third_term_second_type}
            \begin{aligned}
                \int_t^T \langle \varphi_1(\cdot), V_1(s, \cdot) \rangle \, \langle \varphi_2(\cdot), U_2(s, \cdot) U_3(s, \cdot) \rangle \, \ud s,
            \end{aligned}
        \end{equation}
        where $\varphi_1(y), \varphi_2(y)$ are Schwartz functions, where $V_1(s,y)$ is given by $e_\ulell^{\#}(y, -\ulg \ulell) h_\ulell(s)$, $\overline{e_\ulell^{\#}(y, -\ulg \ulell) h_\ulell(s)}$, $\Remusubeone(s,y)$, $\ulg^{-1} \bigl(D^\ast e_\ulell^{\#}\bigr)(y,-\ulg\ulell) p_\ulell(s)$, $\overline{\ulg^{-1} \bigl(D^\ast e_\ulell^{\#}\bigr)(y,-\ulg\ulell) p_\ulell(s)}$, or $\Remusubetwo(s,y)$, and where $U_k(s,y)$, $2 \leq k \leq 3$, is given by $e_\ulell^{\#}(y, -\ulg \ulell) h_\ulell(s)$, $\overline{e_\ulell^{\#}(y, -\ulg \ulell) h_\ulell(s)}$, or $\Remusubeone(s,y)$,
        with at least one of the inputs $V_1(s,y)$, $U_1(s,y)$, $U_2(s,y)$ given by $\Remusubeone(s,y)$ or $\Remusubetwo(s,y)$. In view of the improved local decay estimates \eqref{equ:consequences_Remusubeone_H1y_local_decay} for $\Remusubeone(s,y)$ and \eqref{equ:consequences_Remusubetwo_L2y_local_decay} for $\Remusubetwo(s,y)$, it is clear that the contributions of the second type of terms \eqref{equ:modulation_proof_third_term_second_type} can also be estimated analogously to the corresponding second type of terms \eqref{equ:modulation_proof_contributionII_type2} among the contributions of the term $II(s)$. Again we omit further details.       

        We note that the treatment of the contributions of the term $III(s)$ is the only instance in this work, where a leading order local decay decomposition of the second component $\usubetwo(s)$ of the radiation term is needed.

\medskip 

    \noindent \underline{Contributions of the terms $IV(s)$, $V(s)$, $VI(s)$, and $VII(s)$.}
        Using \eqref{equ:consequences_dispersive_decay_ulPe_radiation}, \eqref{equ:modulation_proof_op_norm_inv_bbM_bmu}, \eqref{equ:modulation_proof_op_norm_bbD}, \eqref{equ:modulation_proof_higher_order_nonlin_decay} it is straightforward to conclude that
        \begin{equation*}
            |IV(s)| + |V(s)| + |VI(s)| + |VII(s)| \lesssim \varepsilon^2 \js^{-2+\delta}.
        \end{equation*}
        This pointwise-in-time bound suffices to deduce the asserted decay estimate \eqref{equ:modulation_proof_goal} upon integration in time.
\end{proof}

\section{Energy Estimates for the Profile} \label{sec:energy_estimates}

In this section we establish the Sobolev bounds as well as the weighted energy estimates for the effective profile. 

\subsection{Preparations}

We begin with some preparations.
In the next lemma we obtain Sobolev and weighted energy bounds for the quadratic normal form $B[\gulellsh](t,\xi)$ defined in \eqref{equ:setting_up_definition_B}.

\begin{lemma}
    Suppose the assumptions in the statement of Proposition~\ref{prop:profile_bounds} are in place.
    Then we have for all $0 \leq t \leq T$ that
    \begin{align}
        \bigl\| \jxi^3 B\bigl[\gulellsh\bigr](t,\xi) \bigr\|_{L^2_\xi} &\lesssim \varepsilon^2 \jt^{-1}, \label{equ:energy_preparations_jxi3_B_bound} \\
        \bigl\| \jxi^3 \pxi B\bigl[\gulellsh\bigr](t,\xi) \bigr\|_{L^2_\xi} &\lesssim \varepsilon^2. \label{equ:energy_preparations_jxi3_pxi_B_bound} 
    \end{align}
\end{lemma}
\begin{proof}
    We begin with the proof of \eqref{equ:energy_preparations_jxi3_B_bound}. 
    By Remark~\ref{rem:setting_up_rapid_decay_qjs}, \eqref{equ:setting_up_phases_lower_bound}, and \eqref{equ:consequences_hulell_decay}, we infer from \eqref{equ:setting_up_definition_B} uniformly for $0 \leq t \leq T$ that
    \begin{equation*}
        \begin{aligned}
            &\bigl\| \jxi^3 B\bigl[\gulellsh\bigr](t,\xi) \bigr\|_{L^2_\xi} \\ 
            &\lesssim |h_\ulell(t)|^2 \Bigl( \bigl\| \jxi^3 \bigl( \jxi + \ulell \xi - 2 \ulg^{-1} \bigr)^{-1} \frakq_{1,\ulell}(\xi) \bigr\|_{L^2_\xi} + \bigl\| \jxi^3 \frakq_{2,\ulell}(\xi) \bigr\|_{L^2_\xi} + \bigl\| \jxi^3 \frakq_{3,\ulell}(\xi) \bigr\|_{L^2_\xi} \Bigr) \\ 
            &\lesssim \varepsilon^2 \jt^{-1}.
        \end{aligned}
    \end{equation*}
    This proves \eqref{equ:energy_preparations_jxi3_B_bound}.
    Next, we deduce the bound \eqref{equ:energy_preparations_jxi3_pxi_B_bound}. 
    Again invoking Remark~\ref{rem:setting_up_rapid_decay_qjs}, \eqref{equ:setting_up_phases_lower_bound}, and \eqref{equ:consequences_hulell_decay}, we find for $0 \leq t \leq T$,
    \begin{equation*}
        \begin{aligned}
            &\bigl\| \jxi^3 \pxi B\bigl[\gulellsh\bigr](t,\xi) \bigr\|_{L^2_\xi} \\ 
            &\lesssim \jt |h_\ulell(t)|^2 \sum_{k=0,1} \, \biggl( \Bigl\| \jxi^3 \pxi^k \Bigl( \bigl( \jxi + \ulell \xi - 2 \ulg^{-1} \bigr)^{-1} \frakq_{1,\ulell}(\xi) \Bigr) \Bigr\|_{L^2_\xi} \\
            &\qquad \qquad \qquad \qquad \qquad \qquad \qquad \qquad + \bigl\| \jxi^3 \pxi^k \frakq_{2,\ulell}(\xi) \bigr\|_{L^2_\xi} + \bigl\| \jxi^3 \pxi^k \frakq_{3,\ulell}(\xi) \bigr\|_{L^2_\xi} \biggr) \\ 
            &\lesssim \jt \cdot \varepsilon^2 \jt^{-1} \lesssim_{\ell_0} \varepsilon^2,
        \end{aligned}
    \end{equation*}
    as desired.
\end{proof}

Next, we infer Sobolev bounds and weighted energy bounds for $\calR(t,\xi)$ defined in \eqref{equ:setting_up_definition_calR}. 
We recall that the collection of nonlinear terms in $\calR(t,\xi)$ should all be thought of as spatially localized with cubic-type $\jt^{-\frac32+\delta}$ decay apart from $\calR_2(\usubeone(t))$.

\begin{lemma} \label{lem:energy_estimates_calR_bounds}
    Suppose the assumptions in the statement of Proposition~\ref{prop:profile_bounds} are in place.
    Then we have for all $0 \leq t \leq T$ that
    \begin{align}
        \bigl\| \jxi^2 \calR(t,\xi) \bigr\|_{L^2_\xi} &\lesssim \varepsilon^2 \jt^{-\frac32+2\delta}, \label{equ:energy_preparations_jxi2_L2_calR} \\
        \bigl\| \jxi^2 \pxi \calR(t,\xi) \bigr\|_{L^2_\xi} &\lesssim \varepsilon^2 \jt^{-1+\delta}. \label{equ:energy_preparations_jxi2_pxi_L2_calR}
    \end{align}
\end{lemma}
\begin{proof}
    Throughout we consider times $0 \leq t \leq T$.
    Recall from \eqref{equ:setting_up_definition_calR} that 
    $$\calR(t,\xi) := \widetilde{\calR}(t,\xi) - \calR_q(t,\xi) - \calFulellsh\bigl[ \calQ_{\ulell,r}\bigl(\usubeone(t)\bigr) \bigr](\xi).$$
    In view of the previously established bounds \eqref{equ:consequences_wtilcalR_L2xi_bound}, \eqref{equ:consequences_wtilcalR_pxi_L2xi_bound} for $\widetilde{\calR}(t,\xi)$, it suffices to prove \eqref{equ:energy_preparations_jxi2_L2_calR} and \eqref{equ:energy_preparations_jxi2_pxi_L2_calR} for the terms $\calR_q(t,\xi)$ and $\calFulellsh\bigl[ \calQ_{\ulell,r}\bigl(\usubeone(t)\bigr)\bigr](\xi)$.
    By Remark~\ref{rem:setting_up_rapid_decay_qjs}, \eqref{equ:consequences_hulell_decay}, and \eqref{equ:consequences_hulell_phase_filtered_decay}, we obtain from the definition \eqref{equ:setting_up_definition_calRq} of $\calR_q(t,\xi)$ that
    \begin{equation*}
        \begin{aligned}
            &\bigl\| \jxi^2 \calR_q(t,\xi) \bigr\|_{L^2_\xi} + \bigl\| \jxi^2 \pxi \calR_q(t,\xi) \bigr\|_{L^2_\xi} \\
            &\lesssim |h_\ulell(t)| \bigl| \pt \bigl( e^{-it\ulg^{-1}} h_\ulell(t) \bigr) \bigr| \\ 
            &\qquad \times \sum_{k=0,1} \, \biggl( \, \Bigl\| \jxi^2 \pxi^k \Bigl( \bigl( \jxi + \ulell \xi - 2 \ulg^{-1} \bigr)^{-1} \frakq_{1,\ulell}(\xi) \Bigr) \Bigr\|_{L^2_\xi} + \bigl\| \jxi^2 \pxi^k \frakq_{2,\ulell}(\xi) \bigr\|_{L^2_\xi} + \bigl\| \jxi^2 \pxi^k \frakq_{3,\ulell}(\xi) \bigr\|_{L^2_\xi} \biggr) \\ 
            &\lesssim \varepsilon \jt^{-\frac12} \cdot \varepsilon \jt^{-1+\delta} \lesssim \varepsilon^2 \jt^{-\frac32+\delta}.
        \end{aligned}
    \end{equation*}
    Moreover, using the mapping properties \eqref{equ:mapping_property_calFulellsh_jxi2}, \eqref{equ:mapping_property_calFulellsh_jxi2_pxi}, we infer from the definition \eqref{equ:setting_up_definition_calQr} of $\calQ_{\ulell,r}\bigl(\usubeone(t)\bigr)$ and from the leading order local decay decomposition \eqref{equ:consequences_decomposition_leading_order_local_decay} of $\usubeone(t,y)$ that
    \begin{equation} \label{equ:prop_sobolev_profile_bounds_calQr_contributions1}
        \begin{aligned}
            &\bigl\| \jxi^2 \calF_{\ulell}^{\#}\bigl[ \calQ_{\ulell,r}\bigl(\usubeone(t)\bigr)\bigr](\xi) \bigr\|_{L^2_\xi} + \bigl\| \jxi^2 \pxi \calF_{\ulell}^{\#}\bigl[ \calQ_{\ulell,r}\bigl(\usubeone(t)\bigr)\bigr](\xi) \bigr\|_{L^2_\xi} \\
            &\lesssim \biggl\| \jy \alpha(\ulg y) \biggl( \bigl( \usubeone(t,y) \bigr)^2 - \Bigl( \Re \, \bigl( e_\ulell^{\#}\bigl(y, - \ulg \ulell\bigr) h_\ulell(t) \bigr) \Bigr)^2 \biggr) \biggr\|_{H^2_y} \\ 
            &\lesssim \bigl\| \jy^{-3} \Remusubeone(t,y) \bigr\|_{H^2_y} \Bigl( \bigl\| \jy^{-1} \usubeone(t,y) \bigr\|_{H^2_y} + |h_\ulell(t)| \Bigr) \\
            &\lesssim \varepsilon \jt^{-1+2\delta} \cdot \varepsilon \jt^{-\frac12} \lesssim \varepsilon^2 \jt^{-\frac32+2\delta}.
        \end{aligned}
    \end{equation}
    Here we made use of the spatial localization of the coefficient $\alpha(\ulg y)$, the Sobolev embedding $H^1_y(\bbR) \hookrightarrow L^\infty_y(\bbR)$, and the local decay estimates \eqref{equ:consequences_Remusubeone_H3y_local_decay}, \eqref{equ:consequences_local_decay_usubeonetwo}, \eqref{equ:consequences_hulell_decay}.  
    This finishes the proof of the lemma.
\end{proof}

\subsection{Sobolev estimates}

The derivation of the Sobolev bounds for the effective profile in the following proposition is straightforward.

\begin{proposition} \label{prop:sobolev_profile_bounds}
    Suppose the assumptions in the statement of Proposition~\ref{prop:profile_bounds} are in place.
    Then we have for all $0 \leq t \leq T$ that
    \begin{equation} \label{equ:prop_sobolev_profile_bounds}
        \bigl\|\jxi^2 \gulellsh(t,\xi)\bigr\|_{L^2_\xi} \lesssim \varepsilon + \varepsilon^2 \jt^{\delta}.
    \end{equation}
\end{proposition}
\begin{proof}
    Upon integrating \eqref{equ:setting_up_g_evol_equ4} in time, we obtain for $0 \leq t \leq T$ that
    \begin{equation} \label{equ:sobolev_profile_bounds_list}
        \begin{aligned}
            &\bigl\| \jxi^2 \gulellsh(t,\xi) \bigr\|_{L^2_\xi} \\ 
            &\lesssim \bigl\| \jxi^2 \gulellsh(0,\xi) \bigr\|_{L^2_\xi} + \sup_{0 \leq s \leq t} \, \bigl\| \jxi^2 B\bigl[\gulellsh\bigr](s,\xi) \bigr\|_{L^2_\xi} + \int_0^t |\dot{q}(s) - \ulell| \bigl\| \jxi^3 B\bigl[ \gulellsh \bigr](s,\xi) \bigr\|_{L^2_\xi} \, \ud s \\ 
            &\quad + \int_0^t \bigl\| \jxi^2 \calF_{\ulell}^{\#}\bigl[ \usubeone(s)^3 \bigr](\xi) \bigr\|_{L^2_\xi} \, \ud s + \int_0^t \bigl\| \jxi^2 \calR(s,\xi) \bigr\|_{L^2_\xi} \, \ud s.
        \end{aligned}
    \end{equation}
    We now estimate each term on the right-hand side of \eqref{equ:sobolev_profile_bounds_list}.

    By Lemma~\ref{lem:higherorder} we have 
    \begin{equation*}
        \bigl\| \jxi^2 \gulellsh(0,\xi) \bigr\|_{L^2_\xi} \lesssim \|\bmu_0\|_{H^3_y \times H^2_y} \lesssim \varepsilon.
    \end{equation*}
    Invoking \eqref{equ:energy_preparations_jxi3_B_bound}, the second term on the right-hand side of \eqref{equ:sobolev_profile_bounds_list} is trivially bounded by
    \begin{equation*}
        \sup_{0 \leq s \leq t} \, \bigl\| \jxi^2 B\bigl[\gulellsh\bigr](s,\xi) \bigr\|_{L^2_\xi} \lesssim \varepsilon^2.
    \end{equation*}
    Similarly, using \eqref{equ:energy_preparations_jxi3_B_bound} and \eqref{equ:consequences_qdot_minus_ulell_decay} we obtain the following sufficient bound for the third term on the right-hand side of \eqref{equ:sobolev_profile_bounds_list} 
    \begin{equation*}
            \int_0^t |\dot{q}(s) - \ulell| \bigl\| \jxi^3 B\bigl[ \gulellsh \bigr](s,\xi) \bigr\|_{L^2_\xi} \, \ud s \lesssim_{\ell_0} \int_0^t \varepsilon \js^{-1+\delta} \cdot \varepsilon^2 \js^{-1} \, \ud s \lesssim \varepsilon^3.
    \end{equation*}
    Using the mapping property \eqref{equ:mapping_property_calFulellsh_jxi2}, the product rule for differentiation, \eqref{equ:consequences_dispersive_decay_ulPe_radiation}, and \eqref{equ:consequences_sobolev_bound_radiation}, we can estimate the fourth term on the right-hand side of \eqref{equ:sobolev_profile_bounds_list} by
    \begin{equation*}
        \begin{aligned}
            \int_0^t \bigl\| \jxi^2 \calF_{\ulell}^{\#}\bigl[ \usubeone(s)^3 \bigr](\xi) \bigr\|_{L^2_\xi} \, \ud s \lesssim \int_0^t \bigl\| \usubeone(s)^3 \bigr\|_{H^2_y} \, \ud s 
            &\lesssim \int_0^t \|\usubeone(s)\|_{H^2_y} \|\usubeone(s)\|_{W^{1,\infty}_y}^2 \, \ud s \\ 
            &\lesssim \int_0^t \varepsilon \js^{\delta} \cdot \varepsilon^2 \js^{-1} \, \ud s \lesssim \varepsilon^3 \jt^\delta.
        \end{aligned}
    \end{equation*}
    Finally, we use \eqref{equ:energy_preparations_jxi2_L2_calR} to estimate the fifth term on the right-hand side of \eqref{equ:sobolev_profile_bounds_list} by
    \begin{equation*}
        \begin{aligned}
            \int_0^t \bigl\| \jxi^2 \calR(s,\xi) \bigr\|_{L^2_\xi} \, \ud s \lesssim \int_0^t \varepsilon^2 \js^{-\frac32+2\delta} \, \ud s \lesssim \varepsilon^2.
        \end{aligned}
    \end{equation*}

    Combining all of the preceding estimates yields the asserted bound \eqref{equ:prop_sobolev_profile_bounds}.    
\end{proof}

\subsection{Weighted energy estimates} \label{subsec:weighted_energy_estimates}

The main work in this section now goes into the derivation of the weighted energy estimates for the effective profile.
Due to the differing slow growth rates for the weighted energy bounds without additional Sobolev regularity and with additional Sobolev regularity this step requires great care. We therefore provide full details.

\begin{proposition} \label{prop:pxi_profile_bounds}
    Suppose the assumptions in the statement of Proposition~\ref{prop:profile_bounds} are in place.
    Then we have for all $0 \leq t \leq T$ that
    \begin{equation} \label{equ:pxi_profile_bounds}
        \bigl\|\pxi \gulellsh(t,\xi)\bigr\|_{L^2_\xi} \lesssim \varepsilon + \varepsilon^2 \jt^{\delta},
    \end{equation}
    and
    \begin{equation} \label{equ:pxi_with_jxi_profile_bounds}
        \bigl\|\jxi^2 \pxi \gulellsh(t,\xi)\bigr\|_{L^2_\xi} \lesssim \varepsilon + \varepsilon^2 \jt^{2\delta},
    \end{equation}    
\end{proposition}
\begin{proof}
We begin with the proof of the weighted energy bound without additional Sobolev regularity \eqref{equ:pxi_profile_bounds} by writing
\begin{equation*}
    \pxi \gulellsh(t,\xi) = \pxi \Bigl( e^{i\xi\theta(t)} e^{-i\xi\theta(t)} \gulellsh(t,\xi) \Bigr) = i \theta(t) \gulellsh(t,\xi) + e^{i\xi\theta(t)} \pxi \Bigl( e^{-i\xi\theta(t)} \gulellsh(t,\xi) \Bigr).
\end{equation*}
Using the growth bound \eqref{equ:consequences_theta_growth_bound} for the phase $\theta(t)$ along with the uniform-in-time bound \eqref{equ:consequences_uniform_in_time_L2_bound}, it follows that
\begin{equation} \label{equ:pxi_profile_differing_slow_growth_rates1}
    \begin{aligned}
        \|\pxi \gulellsh(t,\xi)\|_{L^2_\xi} &\lesssim |\theta(t)| \|\gulellsh(t,\xi)\|_{L^2_\xi} + \bigl\| \pxi \bigl( e^{-i\xi\theta(t)} \gulellsh(t,\xi) \bigr) \bigr\|_{L^2_\xi} \\ 
        &\lesssim \varepsilon \jt^{\delta} \cdot \varepsilon + \bigl\| \pxi \bigl( e^{-i\xi\theta(t)} \gulellsh(t,\xi) \bigr) \bigr\|_{L^2_\xi}.
    \end{aligned}
\end{equation}
Correspondingly, all the work now goes into deducing a sufficient bound on the second term on the right-hand side of the preceding estimate. Integrating \eqref{equ:setting_up_g_evol_equ4} in time, we obtain for all $0 \leq t \leq T$
\begin{equation} \label{equ:pxi_profile_bound_list}
    \begin{aligned}
        \bigl\| \pxi \bigl( e^{-i\xi\theta(t)} \gulellsh(t,\xi) \bigr) \bigr\|_{L^2_\xi} &\lesssim \bigl\| \pxi \gulellsh(0,\xi) \bigr\|_{L^2_\xi} 
        + \sup_{0 \leq s \leq t} \, \Bigl\| \pxi \Bigl( e^{-i\xi\theta(s)} B\bigl[ \gulellsh \bigr](s,\xi) \Bigr) \Bigr\|_{L^2_\xi} \\ 
        &\quad + \biggl\| \int_0^t \pxi \biggl( e^{-i\xi\theta(s)} ( \dot{q}(s) - \ulell) i\xi B\bigl[\gulellsh\bigr](s,\xi) \biggr) \, \ud s \biggr\|_{L^2_\xi} \\        
        &\quad + \biggl\| \int_0^t \pxi \biggl( e^{-i\xi\theta(s)} e^{-is(\jxi+\ulell\xi)} \calFulellsh\bigl[ \usubeone(s)^3 \bigr](\xi) \biggr) \, \ud s \biggr\|_{L^2_\xi} \\ 
        &\quad + \biggl\| \int_0^t \pxi \biggl( e^{-i\xi\theta(s)} e^{-is(\jxi+\ulell\xi)} \calR(s,\xi) \biggr) \, \ud s \biggr\|_{L^2_\xi}.
    \end{aligned}
\end{equation}
By the mapping properties \eqref{equ:mapping_property_calFulellsh_pxi}, \eqref{equ:mapping_property_calFulellDsh_pxi} the first term on the right-hand side of \eqref{equ:pxi_profile_bound_list} is bounded by
\begin{equation*}
    \bigl\| \pxi \gulellsh(0,\xi) \bigr\|_{L^2_\xi} \lesssim \|\jy \bmu_0\|_{H^1_y \times L^2_y} \lesssim \varepsilon.
\end{equation*}
Invoking the estimates \eqref{equ:energy_preparations_jxi3_B_bound} and \eqref{equ:energy_preparations_jxi3_pxi_B_bound} along with the growth bound \eqref{equ:consequences_theta_growth_bound} for the phase $\theta(s)$, we obtain for $0 \leq s \leq t$ that
\begin{equation*}
    \begin{aligned}
        \Bigl\| \pxi \Bigl( e^{-i\xi\theta(s)} B\bigl[ \gulellsh \bigr](s,\xi) \Bigr) \Bigr\|_{L^2_\xi} 
        &\lesssim |\theta(s)| \bigl\| B\bigl[ \gulellsh \bigr](s,\xi) \bigr\|_{L^2_\xi} + \bigl\| \pxi B\bigl[ \gulellsh \bigr](s,\xi) \bigr\|_{L^2_\xi} \\ 
        &\lesssim \varepsilon \js^\delta \cdot \varepsilon^2 \js^{-1} + \varepsilon^2 \lesssim \varepsilon^2.
    \end{aligned}
\end{equation*}
This gives a sufficient bound for the second term on the right-hand side of \eqref{equ:pxi_profile_bound_list}.
Similarly, using \eqref{equ:energy_preparations_jxi3_B_bound}, \eqref{equ:energy_preparations_jxi3_pxi_B_bound}, \eqref{equ:consequences_theta_growth_bound}, and \eqref{equ:consequences_qdot_minus_ulell_decay}, we can estimate the third term on the right-hand side of \eqref{equ:pxi_profile_bound_list} by 
\begin{equation*}
    \begin{aligned}
        &\biggl\| \int_0^t \pxi \biggl( e^{-i\xi\theta(s)} ( \dot{q}(s) - \ulell) i\xi B\bigl[\gulellsh\bigr](s,\xi) \biggr) \, \ud s \biggr\|_{L^2_\xi} \\
        &\lesssim \int_0^t \bigl( 1 + |\theta(s)| \bigr) \, |\dot{q}(s)-\ulell| \, \bigl\| \jxi B\bigl[\gulellsh\bigr](s,\xi) \bigr\|_{L^2_\xi} \, \ud s + \int_0^t |\dot{q}(s)-\ulell| \, \bigl\| \jxi \pxi B\bigl[\gulellsh\bigr](s,\xi) \bigr\|_{L^2_\xi} \, \ud s \\ 
        &\lesssim \int_0^t \bigl( 1 + \varepsilon \js^\delta \bigr) \cdot \varepsilon \js^{-1+\delta} \cdot \varepsilon^2 \js^{-1} \, \ud s + \int_0^t \varepsilon \js^{-1+\delta} \cdot \varepsilon^2 \, \ud s 
        \lesssim \varepsilon^3 \jt^\delta.
    \end{aligned}
\end{equation*}
The estimates for the fourth term and for the fifth term on the right-hand side of \eqref{equ:pxi_profile_bound_list} are much more involved. The corresponding details are deferred to Lemma~\ref{lem:pxi_cubic} and Lemma~\ref{lem:pxi_calR} below. 
Combining all of the preceding estimates with the bound \eqref{equ:pxi_cubic_claim} from Lemma~\ref{lem:pxi_cubic} below and with the bound \eqref{equ:pxi_calR_duhamel_bound} from Lemma~\ref{lem:pxi_calR} below yields the asserted estimate \eqref{equ:pxi_profile_bounds}.

We proceed similarly for the proof of the second weighted energy estimate with additional Sobolev regularity \eqref{equ:pxi_with_jxi_profile_bounds}.
First, we write
\begin{equation*}
    \jxi^2 \pxi \gulellsh(t,\xi) = \jxi^2 \pxi \Bigl( e^{i\xi\theta(t)} e^{-i\xi\theta(t)} \gulellsh(t,\xi) \Bigr) = i \theta(t) \jxi^2 \gulellsh(t,\xi) + e^{i\xi\theta(t)} \jxi^2 \pxi \Bigl( e^{-i\xi\theta(t)} \gulellsh(t,\xi) \Bigr).
\end{equation*}
Using the growth bound \eqref{equ:consequences_theta_growth_bound} for the phase $\theta(t)$ along with \eqref{equ:prop_profile_bounds_assumption2}, it follows that
\begin{equation} \label{equ:pxi_profile_differing_slow_growth_rates2}
    \begin{aligned}
        \bigl\|\jxi^2 \pxi \gulellsh(t,\xi) \bigr\|_{L^2_\xi} &\lesssim |\theta(t)| \bigl\|\jxi^2 \gulellsh(t,\xi)\bigr\|_{L^2_\xi} + \bigl\| \jxi^2 \pxi \bigl( e^{-i\xi\theta(t)} \gulellsh(t,\xi) \bigr) \bigr\|_{L^2_\xi} \\ 
        &\lesssim \varepsilon \jt^{\delta} \cdot \varepsilon \jt^{\delta} + \bigl\| \jxi^2 \pxi \bigl( e^{-i\xi\theta(t)} \gulellsh(t,\xi) \bigr) \bigr\|_{L^2_\xi} \\ 
        &\lesssim \varepsilon^2 \jt^{2\delta} + \bigl\| \jxi^2 \pxi \bigl( e^{-i\xi\theta(t)} \gulellsh(t,\xi) \bigr) \bigr\|_{L^2_\xi}.
    \end{aligned}
\end{equation}
Now all the work again goes into deducing a sufficient bound on the second term on the right-hand side of the preceding estimate. Integrating \eqref{equ:setting_up_g_evol_equ4} in time, we obtain for all $0 \leq t \leq T$
\begin{equation} \label{equ:pxi_with_jxi_profile_bound_list}
    \begin{aligned}
        &\bigl\| \jxi^2 \pxi \bigl( e^{-i\xi\theta(t)} \gulellsh(t,\xi) \bigr) \bigr\|_{L^2_\xi} \\
        &\quad \lesssim \bigl\| \jxi^2 \pxi \gulellsh(0,\xi) \bigr\|_{L^2_\xi} 
        + \sup_{0 \leq s \leq t} \, \Bigl\| \jxi^2 \pxi \Bigl( e^{-i\xi\theta(s)} B\bigl[ \gulellsh \bigr](s,\xi) \Bigr) \Bigr\|_{L^2_\xi} \\ 
        &\quad \quad + \biggl\| \jxi^2 \int_0^t \pxi \biggl( e^{-i\xi\theta(s)} ( \dot{q}(s) - \ulell) i\xi B\bigl[\gulellsh\bigr](s,\xi) \biggr) \, \ud s \biggr\|_{L^2_\xi} \\        
        &\quad \quad + \biggl\| \jxi^2 \int_0^t \pxi \biggl( e^{-i\xi\theta(s)} e^{-is(\jxi+\ulell\xi)} \calFulellsh\bigl[ \bigl( \usubeone(s) \bigr)^3 \bigr](\xi) \biggr) \, \ud s \biggr\|_{L^2_\xi} \\ 
        &\quad \quad + \biggl\| \jxi^2 \int_0^t \pxi \biggl( e^{-i\xi\theta(s)} e^{-is(\jxi+\ulell\xi)} \calR(s,\xi) \biggr) \, \ud s \biggr\|_{L^2_\xi}.
    \end{aligned}
\end{equation}
By \eqref{equ:mapping_property_calFulellsh_jxi2_pxi}, \eqref{equ:mapping_property_calFulellDsh_jxi2_pxi} the first term on the right-hand side of \eqref{equ:pxi_with_jxi_profile_bound_list} is bounded by
\begin{equation*}
    \bigl\| \jxi^2 \pxi \gulellsh(0,\xi) \bigr\|_{L^2_\xi} \lesssim \|\jy \bmu_0\|_{H^3_y \times H^2_y} \lesssim \varepsilon.
\end{equation*}
Invoking the estimates \eqref{equ:energy_preparations_jxi3_B_bound} and \eqref{equ:energy_preparations_jxi3_pxi_B_bound} along with the growth bound \eqref{equ:consequences_theta_growth_bound} for the phase $\theta(s)$, we obtain for $0 \leq s \leq t$ that
\begin{equation*}
    \begin{aligned}
        \Bigl\| \jxi^2 \pxi \Bigl( e^{-i\xi\theta(s)} B\bigl[ \gulellsh \bigr](s,\xi) \Bigr) \Bigr\|_{L^2_\xi} 
        &\lesssim |\theta(s)| \bigl\| \jxi^2 B\bigl[ \gulellsh \bigr](s,\xi) \bigr\|_{L^2_\xi} + \bigl\| \jxi^2 \pxi B\bigl[ \gulellsh \bigr](s,\xi) \bigr\|_{L^2_\xi} \\ 
        &\lesssim \varepsilon \js^\delta \cdot \varepsilon^2 \js^{-1} + \varepsilon^2 \lesssim \varepsilon^2.
    \end{aligned}
\end{equation*}
This gives a sufficient bound for the second term on the right-hand side of \eqref{equ:pxi_with_jxi_profile_bound_list}.
Similarly, using \eqref{equ:energy_preparations_jxi3_B_bound}, \eqref{equ:energy_preparations_jxi3_pxi_B_bound}, \eqref{equ:consequences_theta_growth_bound}, and \eqref{equ:consequences_qdot_minus_ulell_decay}, we can estimate the third term on the right-hand side of \eqref{equ:pxi_profile_bound_list} by 
\begin{equation*}
    \begin{aligned}
        &\biggl\| \jxi^2 \int_0^t \pxi \biggl( e^{-i\xi\theta(s)} ( \dot{q}(s) - \ulell) i\xi B\bigl[\gulellsh\bigr](s,\xi) \biggr) \, \ud s \biggr\|_{L^2_\xi} \\
        &\lesssim \int_0^t (1+ |\theta(s)|) \, |\dot{q}(s)-\ulell| \, \bigl\| \jxi^3 B\bigl[\gulellsh\bigr](s,\xi) \bigr\|_{L^2_\xi} \, \ud s + \int_0^t |\dot{q}(s)-\ulell| \, \bigl\| \jxi^3 \pxi B\bigl[\gulellsh\bigr](s,\xi) \bigr\|_{L^2_\xi} \, \ud s \\ 
        &\lesssim \int_0^t \bigl( 1 + \varepsilon \js^\delta \bigr) \cdot \varepsilon \js^{-1+\delta} \cdot \varepsilon^2 \js^{-1} \, \ud s + \int_0^t \varepsilon \js^{-1+\delta} \cdot \varepsilon^2 \, \ud s 
        \lesssim \varepsilon^3 \jt^\delta.
    \end{aligned}
\end{equation*}
The estimates for the fourth term and for the fifth term on the right-hand side of \eqref{equ:pxi_profile_bound_list} are again much more involved. The corresponding details are deferred to Lemma~\ref{lem:pxi_cubic} and Lemma~\ref{lem:pxi_calR} below. 
Combining all of the preceding estimates with the bound \eqref{equ:pxi_cubic_with_jxi_claim} from Lemma~\ref{lem:pxi_cubic} below and with the bound \eqref{equ:pxi_with_jxi_calR_duhamel_bound} from Lemma~\ref{lem:pxi_calR} below establishes the second asserted bound \eqref{equ:pxi_with_jxi_profile_bounds}.
\end{proof}

In the next lemma we deduce the weighted energy bounds for the contributions of the (non-localized) cubic nonlinearities on the right-hand side of the renormalized evolution equation \eqref{equ:setting_up_g_evol_equ4} for the effective profile.

\begin{lemma} \label{lem:pxi_cubic}
    Suppose the assumptions in the statement of Proposition~\ref{prop:profile_bounds} are in place.
    Then we have for all $0 \leq t \leq T$ that
    \begin{equation} \label{equ:pxi_cubic_claim}
        \biggl\| \int_0^t \pxi \biggl( e^{-i\xi\theta(s)} e^{-is(\jxi+\ulell\xi)} \calFulellsh\bigl[ \usubeone(s)^3 \bigr](\xi) \biggr) \, \ud s \biggr\|_{L^2_\xi} \lesssim \varepsilon^3 \jt^{\delta},
    \end{equation}
    and 
    \begin{equation} \label{equ:pxi_cubic_with_jxi_claim}
        \biggl\| \jxi^2 \int_0^t \pxi \biggl( e^{-i\xi\theta(s)} e^{-is(\jxi+\ulell\xi)} \calFulellsh\bigl[ \usubeone(s)^3 \bigr](\xi) \biggr) \, \ud s \biggr\|_{L^2_\xi} \lesssim \varepsilon^3 \jt^{2\delta}.
    \end{equation}    
\end{lemma}
\begin{proof}
Throughout we consider times $0 \leq t \leq T$.
We begin by computing
\begin{equation} \label{equ:pxi_cubic_first_computation}
    \begin{aligned}
        &\int_0^t \pxi \biggl( e^{-i\xi\theta(s)} e^{-is(\jxi+\ulell\xi)} \calFulellsh\bigl[ \usubeone(s)^3 \bigr](\xi) \biggr) \, \ud s \\  
        &= \int_0^t (-i) \theta(s) e^{-i\xi\theta(s)} e^{-is(\jxi+\ulell\xi)} \calFulellsh\bigl[ \usubeone(s)^3 \bigr](\xi) \, \ud s \\ 
        &\quad + \int_0^t e^{-i\xi\theta(s)} \pxi \Bigl( e^{-is(\jxi+\ulell\xi)} \calFulellsh\bigl[ \usubeone(s)^3 \bigr](\xi) \Bigr) \, \ud s. 
    \end{aligned}
\end{equation}
It is straightforward to see that the first term on the right-hand side of \eqref{equ:pxi_cubic_first_computation} enjoys acceptable bounds.
Indeed, using \eqref{equ:mapping_property_calFulellsh_L2}, \eqref{equ:consequences_dispersive_decay_ulPe_radiation}, \eqref{equ:consequences_energy_bound_radiation}, and \eqref{equ:consequences_theta_growth_bound}, we find that
\begin{equation*}
    \begin{aligned}
        &\biggl\| \int_0^t (-i) \theta(s) e^{-i\xi\theta(s)} e^{-is(\jxi+\ulell\xi)} \calFulellsh\bigl[ \usubeone(s)^3 \bigr](\xi) \, \ud s \biggr\|_{L^2_\xi} \\
        &\lesssim \int_0^t |\theta(s)| \bigl\| \calFulellsh\bigl[ \usubeone(s)^3 \bigr](\xi) \bigr\|_{L^2_\xi} \, \ud s  
        \lesssim \int_0^t |\theta(s)| \bigl\| \usubeone(s)^3 \bigr\|_{L^2_y} \, \ud s \\
        &\lesssim \int_0^t |\theta(s)| \|\usubeone(s)\|_{L^2_y} \|\usubeone(s)\|_{L^\infty_y}^2 \, \ud s 
        \lesssim \int_0^t \varepsilon \js^{\delta} \cdot \varepsilon \cdot \varepsilon^2 \js^{-1} \, \ud s \lesssim \varepsilon^4 \jt^\delta.
    \end{aligned}
\end{equation*}
Similarly, using \eqref{equ:mapping_property_calFulellsh_jxi2}, \eqref{equ:consequences_dispersive_decay_ulPe_radiation}, \eqref{equ:consequences_sobolev_bound_radiation}, and \eqref{equ:consequences_theta_growth_bound}, we obtain that
\begin{equation*}
    \begin{aligned}
        &\biggl\| \jxi^2 \int_0^t (-i) \theta(s) e^{-i\xi\theta(s)} e^{-is(\jxi+\ulell\xi)} \calFulellsh\bigl[ \usubeone(s)^3 \bigr](\xi) \, \ud s \biggr\|_{L^2_\xi} \\
        &\lesssim \int_0^t |\theta(s)| \bigl\| \jxi^2 \calFulellsh\bigl[ \usubeone(s)^3 \bigr](\xi) \bigr\|_{L^2_\xi} \, \ud s 
        \lesssim \int_0^t |\theta(s)| \bigl\| \usubeone(s)^3 \bigr\|_{H^2_y} \, \ud s \\
        &\lesssim \int_0^t |\theta(s)| \|\usubeone(s)\|_{H^2_y} \|\usubeone(s)\|_{W^{1,\infty}_y}^2 \, \ud s 
        \lesssim \int_0^t \varepsilon \js^\delta \cdot \varepsilon \js^\delta \cdot \varepsilon^2 \js^{-1} \, \ud s \\
        &\lesssim \int_0^t \varepsilon^4 \js^{-1+2\delta} \, \ud s \lesssim \varepsilon^4 \jt^{2\delta}.
    \end{aligned}
\end{equation*}
Thus, to complete the proofs of \eqref{equ:pxi_cubic_claim} and of \eqref{equ:pxi_cubic_with_jxi_claim}, the main work now goes into establishing the following two slowly growing weighted energy bounds for $0 \leq t \leq T$,
    \begin{align}
        \biggl\| \int_0^t e^{-i\xi\theta(s)} \pxi \biggl( e^{-is(\jxi+\ulell\xi)} \calFulellsh\bigl[ \usubeone(s)^3 \bigr](\xi) \biggr) \, \ud s \biggr\|_{L^2_\xi} &\lesssim \varepsilon^3 \jt^\delta, \label{equ:pxi_cubic_reduced_claim} \\
        \biggl\| \jxi^2 \int_0^t e^{-i\xi\theta(s)} \pxi \biggl( e^{-is(\jxi+\ulell\xi)} \calFulellsh\bigl[ \usubeone(s)^3 \bigr](\xi) \biggr) \, \ud s \biggr\|_{L^2_\xi} &\lesssim \varepsilon^3 \jt^{2\delta}. \label{equ:pxi_cubic_with_jxi_reduced_claim}
    \end{align}
To this end we recall from \eqref{equ:setting_up_dist_FT_of_cubic_expanded} the expansion
\begin{equation} 
        e^{-i s (\jxi + \ulell \xi)} \calFulellsh\bigl[ \usubeone(s)^3 \bigr](\xi) 
        = - \frac{i}{8} \calI_1(s,\xi) + \frac{3i}{8} \calI_2(s,\xi) - \frac{3i}{8} \calI_3(s,\xi) + \frac{i}{8} \calI_4(s,\xi),
\end{equation}
where the cubic interaction terms $\calI_j(s,\xi)$, $1 \leq j \leq 4$, are defined in \eqref{equ:definition_calI}.
Moreover, we recall from \eqref{equ:calIj_decomposition} their decompositions
\begin{equation*}
    \calI_j(s,\xi) = \calI_j^{\delta_0}(s,\xi) + \calI_j^{\pvdots}(s,\xi) + \calI_j^{\mathrm{reg}}(s,\xi), \quad 1 \leq j \leq 4,
\end{equation*}
where the three terms on the right-hand side are trilinear expressions in terms of the effective profile $\gulellsh(s,\xi)$ involving a Dirac kernel, a Hilbert-type kernel, and a regular kernel, respectively.

In what follows, we present the details of the derivation of the bounds \eqref{equ:pxi_cubic_reduced_claim} and \eqref{equ:pxi_cubic_with_jxi_reduced_claim} for the contribution of $\calI_2(s,\xi)$. The contributions of the other terms $\calI_j(s,\xi)$, $j = 1, 3, 4$, can be handled analogously, and we leave the details to the reader.

For the subsequent analysis we introduce some short-hand notation for auxiliary linear Klein-Gordon evolutions. 
For coefficients $\frakb_j \in W^{1,\infty}(\bbR)$, $1 \leq j \leq 3$, that will be further specified below, we define 
\begin{align}
        v_j(s,y) &:= \widehat{\calF}^{-1}\Bigl[ e^{is(\jap{\xi_j}+\ulell\xi_j)} \jap{\xi_j}^{-1} \frakb_j(\xi_j) \gulellsh(s,\xi_j) \Bigr](y), \quad (s,y) \in [0,T] \times \bbR, \label{equ:weighted_energies_proof_vj_def} \\
        v_j^{\musSharp}(s,y) &:= \widehat{\calF}^{-1}\Bigl[ e^{is(\jap{\xi_j}+\ulell\xi_j)} \frakb_j(\xi_j) \gulellsh(s,\xi_j) \Bigr](y), \quad \quad \quad \, \, \, \,  (s,y) \in [0,T] \times \bbR. \label{equ:weighted_energies_proof_vj_musSharp_def}
\end{align}
Here the musical notation $v_j^{\musSharp}(s,y)$ shall indicate that the regularity of $v_j^{\musSharp}(s,y)$ is raised with respect to the original evolution $v_j(s,y)$ due to the missing inverse Sobolev weight $\jap{\xi_j}^{-1}$ in its definition \eqref{equ:weighted_energies_proof_vj_musSharp_def}.
The free Klein-Gordon evolutions $v_j(s,y)$, $1 \leq j \leq 3$, enjoy all of the decay estimates in item (11) in the statement of Corollary~\ref{cor:consequences_bootstrap_assumptions}.
By Lemma~\ref{lem:core_linear_dispersive_decay} and by the bootstrap assumptions~\eqref{equ:prop_profile_bounds_assumption2},
the evolutions $v_j^{\musSharp}(s,y)$, $ 1 \leq j \leq 3$, satisfy for $0 \leq s \leq T$ the dispersive decay estimate
\begin{equation} \label{equ:weighted_energies_proof_vjmusSharp_decay}
    \| v_j^{\musSharp}(s) \|_{L^\infty_y} \lesssim \varepsilon \js^{-\frac12}.
\end{equation}
But due to their higher Sobolev regularity, under our bootstrap assumptions we do not have access to a dispersive decay estimate for their spatial derivatives $\py v_j^{\musSharp}(s,y)$.
We denote by $\tilde{v}_j(s,y)$ minor variants of the linear evolutions $v_j(s,y)$ that enjoy the same decay estimates, such as for instance
\begin{equation}
    \tilde{v}_j(s,y) := \widehat{\calF}^{-1}\Bigl[ e^{is(\jap{\xi_j}+\ulell\xi_j)} (\jap{\xi_j}^{-1}  \xi_j)\jap{\xi_j}^{-1} \frakb_j(\xi_j) \gulellsh(s,\xi_j) \Bigr](y), \quad (s,y) \in [0,T] \times \bbR. \label{equ:weighted_energies_proof_tildevj_def}
\end{equation}
Moreover, for coefficients $\frakb_j \in W^{1,\infty}(\bbR)$, $1 \leq j \leq 3$, we define the linear Klein-Gordon evolutions
\begin{align}
        w_j(s,y) &:= \widehat{\calF}^{-1}\Bigl[ e^{is(\jap{\xi_j}+\ulell\xi_j)} \partial_{\xi_j} \bigl( \jap{\xi_j}^{-1} \frakb_j(\xi_j) \gulellsh(s,\xi_j) \bigr) \Bigr](y), \quad \quad \, \, \, \, (s,y) \in [0,T] \times \bbR, \label{equ:weighted_energies_proof_wj_def} \\ 
        w_j^{\musSharp}(s,y) &:= \widehat{\calF}^{-1}\Bigl[ e^{is(\jap{\xi_j}+\ulell\xi_j)} \jap{\xi_j} \partial_{\xi_j} \bigl( \jap{\xi_j}^{-1} \frakb_j(\xi_j) \gulellsh(s,\xi_j) \bigr) \Bigr](y), \quad (s,y) \in [0,T] \times \bbR. \label{equ:weighted_energies_proof_wj_musSharp_def}
\end{align}
By Plancherel's theorem and by the bootstrap assumptions \eqref{equ:prop_profile_bounds_assumption1}, they satisfy for $0 \leq s \leq T$ the Sobolev bounds
\begin{align}
    \|w_j(s)\|_{H^1_y} &\lesssim \bigl\| \gulellsh(s,\xi_j) \bigr\|_{L^2_{\xi_j}} + \bigl\| \partial_{\xi_j} \gulellsh(s,\xi_j) \bigr\|_{L^2_{\xi_j}} \lesssim \varepsilon \js^{\delta}, \label{equ:weighted_energies_proof_wj_H1y_growth} \\
    \|w_j(s)\|_{H^3_y} &\lesssim \bigl\| \jxi^2 \gulellsh(s,\xi_j) \bigr\|_{L^2_{\xi_j}} + \bigl\| \jxi^2 \partial_{\xi_j} \gulellsh(s,\xi_j) \bigr\|_{L^2_{\xi_j}} \lesssim \varepsilon \js^{2\delta}, \label{equ:weighted_energies_proof_wj_H3y_growth}
\end{align}
while the evolutions $w_j^{\musSharp}(s,y)$, $1 \leq j \leq 3$, satisfy for $0 \leq s \leq T$,
\begin{align}
    \|w_j^{\musSharp}(s)\|_{L^2_y} &\lesssim \bigl\| \gulellsh(s,\xi_j) \bigr\|_{L^2_{\xi_j}} + \bigl\| \partial_{\xi_j} \gulellsh(s,\xi_j) \bigr\|_{L^2_{\xi_j}} \lesssim \varepsilon \js^{\delta}, \label{equ:weighted_energies_proof_wjmusSharp_L2y_growth} \\
    \|w_j^{\musSharp}(s)\|_{H^2_y} &\lesssim \bigl\| \jxi^2 \gulellsh(s,\xi_j) \bigr\|_{L^2_{\xi_j}} + \bigl\| \jxi^2 \partial_{\xi_j} \gulellsh(s,\xi_j) \bigr\|_{L^2_{\xi_j}} \lesssim \varepsilon \js^{2\delta}. \label{equ:weighted_energies_proof_wjmusSharp_H2y_growth}
\end{align}

\medskip 
\noindent \underline{Case 1: Cubic interactions with a Dirac kernel.}
We first consider the contributions of the cubic interactions with a Dirac kernel $\calI_2^{\delta_0}(s,\xi)$.
In view of the fine structure of the cubic spectral distributions determined in Subsection~\ref{subsec:cubic_spectral_distributions}, 
the term $\calI_2^{\delta_0}(s,\xi)$ is a linear combination of terms of the following schematic form 
\begin{equation}
    \begin{aligned}
        \calI_{2, \mathrm{schem}}^{\delta_0}(s,\xi) := \frakb(\xi) \iint e^{i s \Psi_{2}(\xi,\xi_1,\xi_2)} h_1(s,\xi_1) \, \overline{h_2(s,\xi_2)} h_3(s,\xi_3) \, \ud \xi_1 \, \ud \xi_2,
    \end{aligned}
\end{equation}
with the phase
\begin{equation*}
    \Psi_2(\xi,\xi_1,\xi_2) := -\jxi + \jxione - \jxitwo + \jxithree, \quad \quad \xi_3 := \xi - \xi_1 + \xi_2,
\end{equation*}
and for some coefficients $\frakb, \frakb_1, \frakb_2, \frakb_3 \in W^{1,\infty}(\bbR)$ such that the inputs are given by
\begin{equation*}
    h_j(s,\xi_j) := \jap{\xi_j}^{-1} \frakb_j(\xi_j) \gulellsh(s,\xi_j), \quad 1 \leq j \leq 3.
\end{equation*}
We begin by computing
\begin{equation} \label{equ:weighted_energy_proof_decomposition_jxi_pxi_calI2delta_schem}
    \begin{aligned}
        \jxi \pxi \calI_{2, \mathrm{schem}}^{\delta_0}(s,\xi) &= \frakb(\xi) \iint i s \cdot \jxi (\pxi \Psi_2) \cdot e^{i s \Psi_2} h_1(s,\xi_1) \, \overline{h_2(s,\xi_2)} h_3(s,\xi_3) \, \ud \xi_1 \, \ud \xi_2 \\
        &\quad + \jxi \frakb(\xi) \iint e^{i s \Psi_2} h_1(s,\xi_1) \, \overline{h_2(s,\xi_2)} (\pxithree h_3)(s,\xi_3) \, \ud \xi_1 \, \ud \xi_2 \\
        &\quad + \jxi (\pxi \frakb)(\xi) \iint e^{i s \Psi_2} h_1(s,\xi_1) \, \overline{h_2(s,\xi_2)} h_3(s,\xi_3) \, \ud \xi_1 \, \ud \xi_2 \\ 
        &=: \calJ_{2,1}^{\delta_0}(s,\xi) + \calJ_{2,2}^{\delta_0}(s,\xi) + \calJ_{2,3}^{\delta_0}(s,\xi).
    \end{aligned}
\end{equation}
In the first term $\calJ_{2,1}^{\delta_0}(s,\xi)$ on the right-hand side of \eqref{equ:weighted_energy_proof_decomposition_jxi_pxi_calI2delta_schem}, 
we insert the following identity
\begin{equation} \label{equ:weighted_energy_proof_phase_derivative_identity_Dirac}
    \jxi \pxi \Psi_2 = - \jxione \pxione \Psi_2 - \jxitwo \pxitwo \Psi_2 - \Psi_2 \frac{\xi_3}{\jxithree},
\end{equation}
which can be verified by direct computation.
Integrating by parts in the frequency variables $\xi_1$ and $\xi_2$ then gives
\begin{equation} \label{equ:weighted_energy_proof_calJ21delta_ibp}
    \begin{aligned}
        \calJ_{2,1}^{\delta_0}(s,\xi) &= \frakb(\xi) \iint e^{i s \Psi_{2}} \pxione \Bigl( \jxione h_1(s,\xi_1) \, \overline{h_2(s,\xi_2)} \, h_3(s,\xi_3) \Bigr) \, \ud \xi_1 \, \ud \xi_2 \\
        &\quad + \frakb(\xi) \iint e^{i s \Psi_{2}} \pxitwo \Bigl( h_1(s,\xi_1) \, \overline{\jxitwo h_2(s,\xi_2)} \, h_3(s,\xi_3) \Bigr) \, \ud \xi_1 \, \ud \xi_2 \\
        &\quad - \frakb(\xi) \iint i s \, \Psi_2 \, e^{i s \Psi_{2}} h_1(s,\xi_1) \, \overline{h_2(s,\xi_2)} \, (\jxithree^{-1} \xi_3) h_3(s,\xi_3) \, \ud \xi_1 \, \ud \xi_2 \\
        &=: \calJ_{2,1,1}^{\delta_0}(s,\xi) + \calJ_{2,1,2}^{\delta_0}(s,\xi) + \calJ_{2,1,3}^{\delta_0}(s,\xi).
    \end{aligned}
\end{equation}
The first term $\calJ_{2,1,1}^{\delta_0}(s,\xi)$ on the right-hand side of \eqref{equ:weighted_energy_proof_calJ21delta_ibp} reads in expanded form 
\begin{equation*}
    \begin{aligned}
    \calJ_{2,1,1}^{\delta_0}(s,\xi) &= \frakb(\xi) \iint e^{i s \Psi_{2}} \jxione (\pxione h_1)(s,\xi_1) \, \overline{h_2(s,\xi_2)} \, h_3(s,\xi_3) \, \ud \xi_1 \, \ud \xi_2 \\
    &\quad + \frakb(\xi) \iint e^{i s \Psi_{2}} (\jxione^{-1} \xi_1) h_1(s,\xi_1) \, \overline{h_2(s,\xi_2)} \, h_3(s,\xi_3) \, \ud \xi_1 \, \ud \xi_2 \\
    &\quad + \frakb(\xi) \iint e^{i s \Psi_{2}} \jxione h_1(s,\xi_1) \, \overline{h_2(s,\xi_2)} \, (-1) (\pxithree h_3)(s,\xi_3) \, \ud \xi_1 \, \ud \xi_2 \\
    &=: \calJ_{2,1,1,1}^{\delta_0}(s,\xi) + \calJ_{2,1,1,2}^{\delta_0}(s,\xi) + \calJ_{2,1,1,3}^{\delta_0}(s,\xi).
    \end{aligned}
\end{equation*}
Analogously, the second term $\calJ_{2,1,2}^{\delta_0}(s,\xi)$ on the right-hand side of \eqref{equ:weighted_energy_proof_calJ21delta_ibp} is given in expanded form by
\begin{equation*}
    \begin{aligned}
    \calJ_{2,1,2}^{\delta_0}(s,\xi) &= \frakb(\xi) \iint e^{i s \Psi_{2}} h_1(s,\xi_1) \, \overline{\jxitwo (\pxitwo h_2)(s,\xi_2)} \, h_3(s,\xi_3) \, \ud \xi_1 \, \ud \xi_2 \\
    &\quad + \frakb(\xi) \iint e^{i s \Psi_{2}} h_1(s,\xi_1) \, \overline{(\jxitwo^{-1} \xi_2) h_2(s,\xi_2)} \, h_3(s,\xi_3) \, \ud \xi_1 \, \ud \xi_2 \\
    &\quad + \frakb(\xi) \iint e^{i s \Psi_{2}} h_1(s,\xi_1) \, \overline{\jxitwo h_2(s,\xi_2)} \, (+1) (\pxithree h_3)(s,\xi_3) \, \ud \xi_1 \, \ud \xi_2 \\
    &=: \calJ_{2,1,2,1}^{\delta_0}(s,\xi) + \calJ_{2,1,2,2}^{\delta_0}(s,\xi) + \calJ_{2,1,2,3}^{\delta_0}(s,\xi).
    \end{aligned}
\end{equation*}
In order to obtain acceptable bounds for the third term $\calJ_{2,1,3}^{\delta_0}(s,\xi)$ on the right-hand side of \eqref{equ:weighted_energy_proof_calJ21delta_ibp}, we will have to exploit the presence of the phase $\Psi_2$ in the integrand and integrate by parts in time. To this end we compute 
\begin{equation*}
    \begin{aligned}
        &\int_0^t e^{-i\xi\theta(s)} \calJ_{2,1,3}^{\delta_0}(s,\xi) \, \ud s \\
        &= - \biggl[ s \cdot e^{-i\xi\theta(s)} \frakb(\xi) \iint e^{i s \Psi_{2}} h_1(s,\xi_1) \, \overline{h_2(s,\xi_2)} \, (\jxithree^{-1} \xi_3) h_3(s,\xi_3) \, \ud \xi_1 \, \ud \xi_2 \biggr]_{s=0}^{s=t} \\
        &\quad + \int_0^t 1 \cdot e^{-i\xi\theta(s)} \frakb(\xi) \iint e^{i s \Psi_{2}} h_1(s,\xi_1) \, \overline{h_2(s,\xi_2)} \, (\jxithree^{-1} \xi_3) h_3(s,\xi_3) \, \ud \xi_1 \, \ud \xi_2 \, \ud s \\
        &\quad + \int_0^t s \cdot (-i) \xi \cdot (\dot{q}(s)-\ulell) \cdot e^{-i\xi\theta(s)} \frakb(\xi) \iint e^{i s \Psi_{2}} h_1(s,\xi_1) \, \overline{h_2(s,\xi_2)} \, (\jxithree^{-1} \xi_3) h_3(s,\xi_3) \, \ud \xi_1 \, \ud \xi_2 \, \ud s \\        
        &\quad + \int_0^t s \cdot e^{-i\xi\theta(s)} \frakb(\xi) \iint e^{i s \Psi_{2}} (\ps h_1)(s,\xi_1) \, \overline{h_2(s,\xi_2)} \, (\jxithree^{-1} \xi_3) h_3(s,\xi_3) \, \ud \xi_1 \, \ud \xi_2 \, \ud s \\
        &\quad + \bigl\{ \text{similar terms} \bigr\} \\ 
        &=: \calK_{2,1,3,1}^{\delta_0}(t,\xi) + \calK_{2,1,3,2}^{\delta_0}(t,\xi) + \calK_{2,1,3,3}^{\delta_0}(t,\xi) + \calK_{2,1,3,4}^{\delta_0}(t,\xi) + \{\text{similar terms}\}.
    \end{aligned}
\end{equation*}
In the last term $\calK_{2,1,3,4}^{\delta_0}(t,\xi)$, we have to insert the evolution equation \eqref{equ:setting_up_g_evol_equ3} for the effective profile, which yields  
\begin{equation*}
    \begin{aligned}
        \calK_{2,1,3,4}^{\delta_0}(t,\xi) 
        &= \int_0^t s \cdot (\dot{q}(s)-\ulell) \cdot e^{-i\xi\theta(s)} \frakb(\xi) \iint e^{i s \Psi_{2}} \jxione^{-1} \frakb_1(\xi_1) i\xi_1 \gulellsh(s,\xi_1) \\
        &\qquad \qquad \qquad \qquad \qquad \qquad \qquad \qquad \times \overline{h_2(s,\xi_2)} \, (\jxithree^{-1} \xi_3) h_3(s,\xi_3) \, \ud \xi_1 \, \ud \xi_2 \, \ud s \\
        &\quad - \int_0^t s \cdot e^{-i\xi\theta(s)} \frakb(\xi) \iint e^{i s \Psi_{2}} \jxione^{-1} \frakb_1(\xi_1) e^{-is(\jxione+\ulell\xi_1)} \calFulellsh\bigl[ \calQ_\ulell\bigl( \usubeone(s) \bigr)\bigr](\xi_1) \\
        &\qquad \qquad \qquad \qquad \qquad \qquad \qquad \qquad \times \overline{h_2(s,\xi_2)} \, (\jxithree^{-1} \xi_3) h_3(s,\xi_3) \, \ud \xi_1 \, \ud \xi_2 \, \ud s \\
        &\quad - \int_0^t s \cdot e^{-i\xi\theta(s)} \frakb(\xi) \iint e^{i s \Psi_{2}} \jxione^{-1} \frakb_1(\xi_1) e^{-is(\jxione+\ulell\xi_1)} \calFulellsh\bigl[ {\textstyle \frac16} \usubeone(s)^3 \bigr](\xi_1) \\
        &\qquad \qquad \qquad \qquad \qquad \qquad \qquad \qquad \times \overline{h_2(s,\xi_2)} \, (\jxithree^{-1} \xi_3) h_3(s,\xi_3) \, \ud \xi_1 \, \ud \xi_2 \, \ud s \\        
        &\quad - \int_0^t s \cdot e^{-i\xi\theta(s)} \frakb(\xi) \iint e^{i s \Psi_{2}} \jxione^{-1} \frakb_1(\xi_1) e^{-is(\jxione+\ulell\xi_1)} \widetilde{\calR}(s,\xi_1) \\
        &\qquad \qquad \qquad \qquad \qquad \qquad \qquad \qquad \times \overline{h_2(s,\xi_2)} \, (\jxithree^{-1} \xi_3) h_3(s,\xi_3) \, \ud \xi_1 \, \ud \xi_2 \, \ud s \\
        &=: \calK_{2,1,3,4,1}^{\delta_0}(t,\xi) + \calK_{2,1,3,4,2}^{\delta_0}(t,\xi) + \calK_{2,1,3,4,3}^{\delta_0}(t,\xi) + \calK_{2,1,3,4,4}^{\delta_0}(t,\xi).
    \end{aligned}
\end{equation*}

Finally, we observe that the two terms $\calJ_{2,1,1,3}^{\delta_0}(s,\xi)$ and $\calJ_{2,1,2,3}^{\delta_0}(s,\xi)$ can be combined with the term $\calJ_{2,2}^{\delta_0}(s,\xi)$ from \eqref{equ:weighted_energy_proof_decomposition_jxi_pxi_calI2delta_schem}, to give
\begin{equation} \label{equ:weighted_energy_proof_Dirac_combined_terms}
    \begin{aligned}
        &\calJ_{2,1,1,3}^{\delta_0}(s,\xi) + \calJ_{2,1,2,3}^{\delta_0}(s,\xi) + \calJ_{2,2}^{\delta_0}(s,\xi) \\ 
        &= \frakb(\xi) \iint e^{i s \Psi_{2}} \bigl( \jxi - \jxione + \jxitwo \bigr) h_1(s,\xi_1) \, \overline{h_2(s,\xi_2)} \, (\pxithree h_3)(s,\xi_3) \, \ud \xi_1 \, \ud \xi_2 \\
        &= \frakb(\xi) \iint e^{i s \Psi_{2}} h_1(s,\xi_1) \, \overline{h_2(s,\xi_2)} \, \jxithree (\pxithree h_3)(s,\xi_3) \, \ud \xi_1 \, \ud \xi_2 \\
        &\quad + \frakb(\xi) \iint e^{i s \Psi_{2}} (-\Psi_2) h_1(s,\xi_1) \, \overline{h_2(s,\xi_2)} \, (\pxithree h_3)(s,\xi_3) \, \ud \xi_1 \, \ud \xi_2 \\ 
        &=: \calJ_{2,\mathrm{comb}, 1}^{\delta_0}(s,\xi) + \calJ_{2,\mathrm{comb}, 2}^{\delta_0}(s,\xi).
    \end{aligned}
\end{equation}
Importantly, in the term $\calJ_{2,\mathrm{comb}, 1}^{\delta_0}(s,\xi)$ the Sobolev weight $\jxithree$ is paired with the frequency derivative $\pxithree$ on the third input. 
Moreover, thanks to the phase $\Psi_2$ in the integrand, the other term $\calJ_{2,\mathrm{comb}, 2}^{\delta_0}(s,\xi)$ is amenable to integration by parts in time. This will be necessary to obtain acceptable bounds for the contribution to \eqref{equ:pxi_cubic_with_jxi_reduced_claim}, i.e., for the contributions to the weighted energy estimates with additional Sobolev weights.
We therefore record that integrating by parts in time for the term $\calJ_{2,\mathrm{comb}, 2}^{\delta_0}(s,\xi)$ leads to
\begin{equation*}
    \begin{aligned}
        &\int_0^t e^{-i\xi\theta(s)} \calJ_{2,\mathrm{comb}, 2}^{\delta_0}(s,\xi) \, \ud s  \\
        &= \biggl[ e^{-i\xi\theta(s)} i \frakb(\xi) \iint e^{is\Psi_2} h_1(s,\xi_1) \, \overline{h_2(s,\xi_2)} \, (\pxithree h_3)(s,\xi_3) \, \ud \xi_1 \, \ud \xi_2 \biggr]_{s=0}^{s=t} \\
        &\quad - \int_0^t (\dot{q}(s)-\ulell) e^{-i\xi\theta(s)} \frakb(\xi) \cdot \xi \iint e^{is\Psi_2} h_1(s,\xi_1) \, \overline{h_2(s,\xi_2)} \, (\pxithree h_3)(s,\xi_3) \, \ud \xi_1 \, \ud \xi_2 \, \ud s \\
        &\quad - \int_0^t e^{-i\xi\theta(s)} i \frakb(\xi) \iint e^{is\Psi_2} (\ps h_1)(s,\xi_1) \, \overline{h_2(s,\xi_2)} \, (\pxithree h_3)(s,\xi_3) \, \ud \xi_1 \, \ud \xi_2 \, \ud s \\ 
        &\quad - \int_0^t e^{-i\xi\theta(s)} i \frakb(\xi) \iint e^{is\Psi_2} h_1(s,\xi_1) \, \overline{(\ps h_2)(s,\xi_2)} \, (\pxithree h_3)(s,\xi_3) \, \ud \xi_1 \, \ud \xi_2 \, \ud s \\
        &\quad - \int_0^t e^{-i\xi\theta(s)} i \frakb(\xi) \iint e^{is\Psi_2} h_1(s,\xi_1) \, \overline{h_2(s,\xi_2)} \, (\pxithree \ps h_3)(s,\xi_3) \, \ud \xi_1 \, \ud \xi_2 \, \ud s \\ 
        &=: \calK_{2,\mathrm{comb},2,1}^{\delta_0}(t,\xi) + \calK_{2,\mathrm{comb},2,2}^{\delta_0}(t,\xi) + \calK_{2,\mathrm{comb},2,3}^{\delta_0}(t,\xi) + \calK_{2,\mathrm{comb},2,4}^{\delta_0}(t,\xi) + \calK_{2,\mathrm{comb},2,5}^{\delta_0}(t,\xi).
    \end{aligned}
\end{equation*}
Note that the terms $\calK_{2,\mathrm{comb},2,3}^{\delta_0}(t,\xi)$ and $\calK_{2,\mathrm{comb},2,4}^{\delta_0}(t,\xi)$ can be dealt with analogously.
In what follows, we only consider $\calK_{2,\mathrm{comb},2,3}^{\delta_0}(t,\xi)$. 
Inserting the evolution equation \eqref{equ:setting_up_g_evol_equ3} for the effective profile in the term $\calK_{2,\mathrm{comb},2,3}^{\delta_0}(t,\xi)$ gives 
\begin{equation*}
    \begin{aligned}
        &\calK_{2,\mathrm{comb},2,3}^{\delta_0}(t,\xi) \\
        &= - \int_0^t e^{-i\xi\theta(s)} i \frakb(\xi) \iint e^{is\Psi_2} (\ps h_1)(s,\xi_1) \, \overline{h_2(s,\xi_2)} \, (\pxithree h_3)(s,\xi_3) \, \ud \xi_1 \, \ud \xi_2 \, \ud s \\   
        &= - \int_0^t  e^{-i\xi\theta(s)} (\dot{q}(s)-\ulell) i \frakb(\xi) \iint e^{is\Psi_2} \jxione^{-1}  \frakb_1(\xi_1) (i\xi_1) \gulellsh(s,\xi_1) \, \overline{h_2(s,\xi_2)} \, (\pxithree h_3)(s,\xi_3) \, \ud \xi_1 \, \ud \xi_2 \, \ud s \\ 
        &\quad - \int_0^t  e^{-i\xi\theta(s)} i \frakb(\xi) \iint e^{is\Psi_2} \jxione^{-1} \frakb_1(\xi_1) e^{-is(\jxione+\ulell\xi_1)} \calFulellsh\bigl[ \calQ_\ulell\bigl( \usubeone(s) \bigr)\bigr](\xi_1) \, \overline{h_2(s,\xi_2)} \, (\pxithree h_3)(s,\xi_3) \, \ud \xi_1 \, \ud \xi_2 \, \ud s \\ 
        &\quad - \int_0^t  e^{-i\xi\theta(s)} i \frakb(\xi) \iint e^{is\Psi_2} \jxione^{-1} \frakb_1(\xi_1) e^{-is(\jxione+\ulell\xi_1)} \calFulellsh\bigl[ {\textstyle \frac16} \usubeone(s)^3 \bigr](\xi_1) \, \overline{h_2(s,\xi_2)} \, (\pxithree h_3)(s,\xi_3) \, \ud \xi_1 \, \ud \xi_2 \, \ud s \\ 
        &\quad - \int_0^t  e^{-i\xi\theta(s)} i \frakb(\xi) \iint e^{is\Psi_2} \jxione^{-1} \frakb_1(\xi_1) e^{-is(\jxione+\ulell\xi_1)} \widetilde{\calR}(s,\xi_1) \, \overline{h_2(s,\xi_2)} \, (\pxithree h_3)(s,\xi_3) \, \ud \xi_1 \, \ud \xi_2 \, \ud s \\
        &=: \calK_{2,\mathrm{comb},2,3,1}^{\delta_0}(t,\xi) + \calK_{2,\mathrm{comb},2,3,2}^{\delta_0}(t,\xi) + \calK_{2,\mathrm{comb},2,3,3}^{\delta_0}(t,\xi) + \calK_{2,\mathrm{comb},2,3,4}^{\delta_0}(t,\xi).
    \end{aligned}
\end{equation*}
Similarly, inserting the evolution equation \eqref{equ:setting_up_g_evol_equ3} for the effective profile in the term $\calK_{2,\mathrm{comb},2,5}^{\delta_0}(t,\xi)$ leads to 
\begin{equation*}
    \begin{aligned}
        &\calK_{2,\mathrm{comb},2,5}^{\delta_0}(t,\xi) \\
        &= - \int_0^t e^{-i\xi\theta(s)} i \frakb(\xi) \iint e^{is\Psi_2} h_1(s,\xi_1) \, \overline{h_2(s,\xi_2)} \, (\pxithree \ps h_3)(s,\xi_3) \, \ud \xi_1 \, \ud \xi_2 \, \ud s \\ 
        &= - \int_0^t e^{-i\xi\theta(s)} (\dot{q}(s) - \ulell) i \frakb(\xi) \iint e^{is\Psi_2} h_1(s,\xi_1) \, \overline{h_2(s,\xi_2)} \, \pxithree \Bigl( \jxithree^{-1} \frakb_3(\xi_3) (i\xi_3) \gulellsh(s,\xi_3) \Bigr) \, \ud \xi_1 \, \ud \xi_2 \, \ud s \\ 
        &\quad + \int_0^t e^{-i\xi\theta(s)} i \frakb(\xi) \iint e^{is\Psi_2} h_1(s,\xi_1) \, \overline{h_2(s,\xi_2)} \, \pxithree \Bigl( \jxithree^{-1} \frakb_3(\xi_3) e^{-is(\jxithree+\ulell\xi_3)} \calFulellsh\bigl[ \calQ_\ulell\bigl( \usubeone(s) \bigr)\bigr](\xi_3) \Bigr) \, \ud \xi_1 \, \ud \xi_2 \, \ud s \\ 
        &\quad + \int_0^t e^{-i\xi\theta(s)} i \frakb(\xi) \iint e^{is\Psi_2} h_1(s,\xi_1) \, \overline{h_2(s,\xi_2)} \, \pxithree \Bigl( \jxithree^{-1} \frakb_3(\xi_3) e^{-is(\jxithree+\ulell\xi_3)} \calFulellsh\bigl[ {\textstyle \frac16} \usubeone(s)^3 \bigr](\xi_3) \Bigr) \, \ud \xi_1 \, \ud \xi_2 \, \ud s \\ 
        &\quad - \int_0^t e^{-i\xi\theta(s)} i \frakb(\xi) \iint e^{is\Psi_2} h_1(s,\xi_1) \, \overline{h_2(s,\xi_2)} \, \pxithree \Bigl( \jxithree^{-1} \frakb_3(\xi_3) e^{-is(\jxithree+\ulell\xi_3)} \widetilde{\calR}(s,\xi_3) \Bigr) \, \ud \xi_1 \, \ud \xi_2 \, \ud s \\ 
        &=: \calK_{2,\mathrm{comb},2,5,1}^{\delta_0}(t,\xi) + \calK_{2,\mathrm{comb},2,5,2}^{\delta_0}(t,\xi) + \calK_{2,\mathrm{comb},2,5,3}^{\delta_0}(t,\xi) + \calK_{2,\mathrm{comb},2,5,4}^{\delta_0}(t,\xi).
    \end{aligned}
\end{equation*}

\begin{figure} \label{figure:decomposition_weighted_energy_terms}

\centering 

\tikzset{every picture/.style={line width=0.75pt}} 

\begin{tikzpicture}[x=0.75pt,y=0.75pt,yscale=-1,xscale=1]

\draw    (111.33,241.67) -- (160.33,203.36) ;
\draw    (111.33,241.67) -- (160.33,282.69) ;
\draw    (111.33,241.67) -- (163,242.69) ;
\draw    (201.67,202.02) -- (251.67,114.02) ;
\draw    (201.67,202.02) -- (251.67,202.02) ;
\draw    (201.67,202.02) -- (251.67,292.69) ;
\draw    (301.33,112.67) -- (341,82.69) ;
\draw    (301.33,112.67) -- (341,143.36) ;
\draw    (301.33,112.67) -- (341.67,112.36) ;
\draw    (301.33,203) -- (341.67,203.36) ;
\draw    (301.33,203) -- (341,233.36) ;
\draw    (301.33,203) -- (341,172.69) ;
\draw    (335.67,292.02) -- (348.34,291.45) ;
\draw [shift={(350.33,291.36)}, rotate = 177.4] [color={rgb, 255:red, 0; green, 0; blue, 0 }  ][line width=0.75]    (8.74,-2.63) .. controls (5.56,-1.12) and (2.65,-0.24) .. (0,0) .. controls (2.65,0.24) and (5.56,1.12) .. (8.74,2.63)   ;
\draw    (301,293.36) .. controls (331.67,277.36) and (307.67,306.02) .. (335.67,292.02) ;
\draw    (356,293.33) -- (401,254.02) ;
\draw    (356,293.33) -- (400.33,278.69) ;
\draw    (356,293.33) -- (401,307.36) ;
\draw    (356,293.33) -- (401,334.02) ;
\draw    (503.67,340.02) -- (517,340.02) ;
\draw [shift={(519,340.02)}, rotate = 180] [color={rgb, 255:red, 0; green, 0; blue, 0 }  ][line width=0.75]    (10.93,-3.29) .. controls (6.95,-1.4) and (3.31,-0.3) .. (0,0) .. controls (3.31,0.3) and (6.95,1.4) .. (10.93,3.29)   ;
\draw    (459.67,339.36) .. controls (498.33,328.02) and (469.67,349.36) .. (503.67,340.02) ;
\draw    (523.67,340.69) -- (561,299.36) ;
\draw    (523.67,340.69) -- (560.33,324.69) ;
\draw    (523.67,340.69) -- (561.67,354.69) ;
\draw    (523.67,340.69) -- (561,380.69) ;
\draw    (250.83,432.86) -- (262.5,432.57) ;
\draw [shift={(264.5,432.52)}, rotate = 178.6] [color={rgb, 255:red, 0; green, 0; blue, 0 }  ][line width=0.75]    (10.93,-3.29) .. controls (6.95,-1.4) and (3.31,-0.3) .. (0,0) .. controls (3.31,0.3) and (6.95,1.4) .. (10.93,3.29)   ;
\draw    (211.83,433.19) .. controls (254.5,418.52) and (207.83,447.19) .. (250.83,432.86) ;
\draw    (269.5,433.52) -- (309.83,371.86) ;
\draw    (269.5,433.52) -- (310.1,402.79) ;
\draw    (269.5,433.52) -- (310.5,433.19) ;
\draw    (269.5,433.52) -- (310.1,463.19) ;
\draw    (269.5,433.52) -- (310.5,493.19) ;
\draw    (449.5,432.79) -- (479.7,412.79) ;
\draw    (449.5,432.79) -- (479.7,425.99) ;
\draw    (449.5,432.79) -- (480.1,438.79) ;
\draw    (449.5,432.79) -- (479.7,452.79) ;
\draw    (449.7,491.99) -- (479.3,472.39) ;
\draw    (449.7,491.99) -- (480.1,486.79) ;
\draw    (449.7,491.99) -- (480.1,497.99) ;
\draw    (449.7,491.99) -- (480.1,512.79) ;
\draw    (430.5,433.19) -- (442.9,433.19) ;
\draw [shift={(444.9,433.19)}, rotate = 180] [color={rgb, 255:red, 0; green, 0; blue, 0 }  ][line width=0.75]    (10.93,-3.29) .. controls (6.95,-1.4) and (3.31,-0.3) .. (0,0) .. controls (3.31,0.3) and (6.95,1.4) .. (10.93,3.29)   ;
\draw    (391.1,432.59) .. controls (429.7,420.39) and (390.1,445.99) .. (430.5,433.19) ;
\draw    (432.5,492.79) -- (443.3,492.79) ;
\draw [shift={(445.3,492.79)}, rotate = 180] [color={rgb, 255:red, 0; green, 0; blue, 0 }  ][line width=0.75]    (10.93,-3.29) .. controls (6.95,-1.4) and (3.31,-0.3) .. (0,0) .. controls (3.31,0.3) and (6.95,1.4) .. (10.93,3.29)   ;
\draw    (391.1,492.99) .. controls (430.5,479.99) and (390.9,505.59) .. (432.5,492.79) ;

\draw (25.67,227.73) node [anchor=north west][inner sep=0.75pt]  {$\langle \xi \rangle \partial_{\xi}\mathcal{I}_{2,\mathrm{schem}}^{\delta_{0}}$};
\draw (169.33,187.73) node [anchor=north west][inner sep=0.75pt]    {$\mathcal{J}_{2,1}^{\delta _{0}}$};
\draw (170,229.4) node [anchor=north west][inner sep=0.75pt]  [color={rgb, 255:red, 208; green, 2; blue, 27 }  ,opacity=1 ]  {$\mathcal{J}_{2,2}^{\delta _{0}}$};
\draw (168.67,269.73) node [anchor=north west][inner sep=0.75pt]    {$\mathcal{J}_{2,3}^{\delta _{0}}$};
\draw (258.67,97.4) node [anchor=north west][inner sep=0.75pt]    {$\mathcal{J}_{2,1,1}^{\delta _{0}}$};
\draw (257.33,186.07) node [anchor=north west][inner sep=0.75pt]    {$\mathcal{J}_{2,1,2}^{\delta _{0}}$};
\draw (257.33,277.4) node [anchor=north west][inner sep=0.75pt]    {$\mathcal{J}_{2,1,3}^{\delta _{0}}$};
\draw (346,68.07) node [anchor=north west][inner sep=0.75pt]    {$\mathcal{J}_{2,1,1,1}^{\delta _{0}}$};
\draw (345.33,98.07) node [anchor=north west][inner sep=0.75pt]    {$\mathcal{J}_{2,1,1,2}^{\delta _{0}}$};
\draw (345.33,127.73) node [anchor=north west][inner sep=0.75pt]  [color={rgb, 255:red, 208; green, 2; blue, 27 }  ,opacity=1 ]  {$\mathcal{J}_{2,1,1,3}^{\delta _{0}}$};
\draw (346.67,156.73) node [anchor=north west][inner sep=0.75pt]    {$\mathcal{J}_{2,1,2,1}^{\delta _{0}}$};
\draw (347.33,188.07) node [anchor=north west][inner sep=0.75pt]    {$\mathcal{J}_{2,1,2,2}^{\delta _{0}}$};
\draw (347.33,218.07) node [anchor=north west][inner sep=0.75pt]  [color={rgb, 255:red, 208; green, 2; blue, 27 }  ,opacity=1 ]  {$\mathcal{J}_{2,1,2,3}^{\delta _{0}}$};
\draw (403.33,237.07) node [anchor=north west][inner sep=0.75pt]    {$\mathcal{K}_{2,1,3,1}^{\delta _{0}}$};
\draw (402.67,264.73) node [anchor=north west][inner sep=0.75pt]    {$\mathcal{K}_{2,1,3,2}^{\delta _{0}}$};
\draw (404,293.73) node [anchor=north west][inner sep=0.75pt]    {$\mathcal{K}_{2,1,3,3}^{\delta _{0}}$};
\draw (404.67,321.73) node [anchor=north west][inner sep=0.75pt]    {$\mathcal{K}_{2,1,3,4}^{\delta _{0}}$};
\draw (293.00,298.00) node [anchor=north west][inner sep=0.75pt]   [align=left] {{\tiny i.b.p.~in time}};
\draw (563.33,283.07) node [anchor=north west][inner sep=0.75pt]    {$\mathcal{K}_{2,1,3,4,1}^{\delta _{0}}$};
\draw (562.67,311.07) node [anchor=north west][inner sep=0.75pt]    {$\mathcal{K}_{2,1,3,4,2}^{\delta _{0}}$};
\draw (562.67,339.73) node [anchor=north west][inner sep=0.75pt]    {$\mathcal{K}_{2,1,3,4,3}^{\delta _{0}}$};
\draw (562.67,367.4) node [anchor=north west][inner sep=0.75pt]    {$\mathcal{K}_{2,1,3,4,4}^{\delta _{0}}$};
\draw (14.17,398.23) node [anchor=north west][inner sep=0.75pt]    {$ \begin{array}{l}
\textcolor[rgb]{0.82,0.01,0.11}{\mathcal{J}_{2,1,1,3}^{\delta _{0}}}\textcolor[rgb]{0.82,0.01,0.11}{+}\textcolor[rgb]{0.82,0.01,0.11}{\mathcal{J}_{2,1,2,3}^{\delta _{0}}}\textcolor[rgb]{0.82,0.01,0.11}{+}\textcolor[rgb]{0.82,0.01,0.11}{\mathcal{J}_{2,2}^{\delta _{0}}}\\
\ \ \ \ \ \ \ \ =\mathcal{J}_{2,\mathrm{comb} ,1}^{\delta _{0}} +\mathcal{J}_{2,\mathrm{comb} ,2}^{\delta _{0}}
\end{array}$};
\draw (460.00,346.33) node [anchor=north west][inner sep=0.75pt]   [align=left] {\begin{minipage}[lt]{45.9pt}\setlength\topsep{0pt}
{\tiny profile equ.}
\end{minipage}};
\draw (309.5,356.39) node [anchor=north west][inner sep=0.75pt]    {$\mathcal{K}_{2,\mathrm{comb} ,2,1}^{\delta _{0}}$};
\draw (310.7,386.59) node [anchor=north west][inner sep=0.75pt]    {$\mathcal{K}_{2,\mathrm{comb} ,2,2}^{\delta _{0}}$};
\draw (311.1,417.59) node [anchor=north west][inner sep=0.75pt]    {$\mathcal{K}_{2,\mathrm{comb} ,2,3}^{\delta _{0}}$};
\draw (312.7,446.79) node [anchor=north west][inner sep=0.75pt]    {$\mathcal{K}_{2,\mathrm{comb} ,2,4}^{\delta _{0}}$};
\draw (312.3,476.59) node [anchor=north west][inner sep=0.75pt]    {$\mathcal{K}_{2,\mathrm{comb} ,2,5}^{\delta _{0}}$};
\draw (203.9,439.19) node [anchor=north west][inner sep=0.75pt]   [align=left] {{\tiny i.b.p.~in time}};
\draw (388.0,441.59) node [anchor=north west][inner sep=0.75pt]   [align=left] {\begin{minipage}[lt]{40.49pt}\setlength\topsep{0pt}
\begin{center}
{\tiny profile equ.}
\end{center}

\end{minipage}};
\draw (388.0,501.59) node [anchor=north west][inner sep=0.75pt]   [align=left] {\begin{minipage}[lt]{40.49pt}\setlength\topsep{0pt}
\begin{center}
{\tiny profile equ.}
\end{center}

\end{minipage}};
\draw (482.7,405.79) node [anchor=north west][inner sep=0.75pt]  [font=\tiny]  {$\mathcal{K}_{2,\mathrm{comb} ,2,3,1}^{\delta _{0}}$};
\draw (482.7,419.79) node [anchor=north west][inner sep=0.75pt]  [font=\tiny]  {$\mathcal{K}_{2,\mathrm{comb} ,2,3,2}^{\delta _{0}}$};
\draw (482.7,433.59) node [anchor=north west][inner sep=0.75pt]  [font=\tiny]  {$\mathcal{K}_{2,\mathrm{comb} ,2,3,3}^{\delta _{0}}$};
\draw (482.3,448.59) node [anchor=north west][inner sep=0.75pt]  [font=\tiny]  {$\mathcal{K}_{2,\mathrm{comb} ,2,3,4}^{\delta _{0}}$};
\draw (482.3,464.79) node [anchor=north west][inner sep=0.75pt]  [font=\tiny]  {$\mathcal{K}_{2,\mathrm{comb} ,2,5,1}^{\delta _{0}}$};
\draw (482.3,478.59) node [anchor=north west][inner sep=0.75pt]  [font=\tiny]  {$\mathcal{K}_{2,\mathrm{comb} ,2,5,2}^{\delta _{0}}$};
\draw (482.3,493.39) node [anchor=north west][inner sep=0.75pt]  [font=\tiny]  {$\mathcal{K}_{2,\mathrm{comb} ,2,5,3}^{\delta _{0}}$};
\draw (482.7,507.19) node [anchor=north west][inner sep=0.75pt]  [font=\tiny]  {$\mathcal{K}_{2,\mathrm{comb} ,2,5,4}^{\delta _{0}}$};

\end{tikzpicture}
\caption{Summary of the decomposition of $\jxi \pxi \calI_{2,\mathrm{schem}}^{\delta_0}(s,\xi)$.}
\end{figure}

The decomposition of $\jxi \pxi \calI_{2,\mathrm{schem}}^{\delta_0}(s,\xi)$ into the preceding terms is illustrated in Figure~\ref{figure:decomposition_weighted_energy_terms}.
Next, we separately estimate the contributions of all of these terms to the two asserted slowly growing weighted energy bounds \eqref{equ:pxi_cubic_reduced_claim} and \eqref{equ:pxi_cubic_with_jxi_reduced_claim}.

\medskip 
\noindent \underline{Case 1.1: Contribution to \eqref{equ:pxi_cubic_reduced_claim}.}
Here our goal is to prove for $0 \leq t \leq T$ that
\begin{equation*}
    \biggl\| \int_0^t e^{-i\xi\theta(s)} \pxi \calI_{2, \mathrm{schem}}^{\delta_0}(s,\xi) \, \ud s \biggr\|_{L^2_\xi} \lesssim \varepsilon^3 \jt^\delta.
\end{equation*}
For technical convenience, in this case we first freely ``spend one derivative'' and start from the bound
\begin{equation}
    \begin{aligned}
     \biggl\| \int_0^t e^{-i\xi\theta(s)} \pxi \calI_{2, \mathrm{schem}}^{\delta_0}(s,\xi) \, \ud s \biggr\|_{L^2_\xi} &= \biggl\| \jxi^{-1} \int_0^t e^{-i\xi\theta(s)} \jxi \pxi \calI_{2, \mathrm{schem}}^{\delta_0}(s,\xi) \, \ud s \biggr\|_{L^2_\xi} \\
     &\lesssim  \biggl\| \int_0^t e^{-i\xi\theta(s)} \jxi \pxi \calI_{2, \mathrm{schem}}^{\delta_0}(s,\xi) \, \ud s \biggr\|_{L^2_\xi}. 
     \end{aligned}
\end{equation}
Now we are in the position to invoke the decomposition of $\jxi \pxi \calI_{2, \mathrm{schem}}^{\delta_0}(s,\xi)$ into the terms illustrated in Figure~\ref{figure:decomposition_weighted_energy_terms}. We proceed term by term.
To facilitate the derivation of acceptable bounds, it will be useful to rewrite these expressions in terms of the auxiliary linear Klein-Gordon evolutions $v_j(s,y)$, $\tilde{v}_j(s,y)$, $v_j^{\musSharp}(s,y)$, $w_j(s,y)$, and $w_j^{\musSharp}(s,y)$ defined in \eqref{equ:weighted_energies_proof_vj_def}, \eqref{equ:weighted_energies_proof_vj_musSharp_def}, \eqref{equ:weighted_energies_proof_tildevj_def}, \eqref{equ:weighted_energies_proof_wj_def}, and \eqref{equ:weighted_energies_proof_wj_musSharp_def}.

\medskip 

\noindent {\it Bounds for $\calJ_{2,1,k,1}^{\delta_0}(s,\xi)$, $1 \leq k \leq 2$:}
We provide the details for $k=1$. The case $k=2$ can be treated analogously. 
By Plancherel's theorem, \eqref{equ:consequences_aux_KG_disp_decay}, and \eqref{equ:weighted_energies_proof_wjmusSharp_L2y_growth}, we have 
\begin{equation*}
    \begin{aligned}
        &\biggl\| \int_0^t e^{-i\xi\theta(s)} \calJ_{2,1,1,1}(s,\xi) \, \ud s \biggr\|_{L^2_\xi} \\ 
        &= \biggl\| \int_0^t (2\pi)^{\frac52} \frakb(\xi) e^{-i\xi\theta(s)} e^{-is(\jxi+\ulell\xi)} \whatcalF\bigl[ w_1^{\musSharp}(s,\cdot) \overline{v_2(s,\cdot)} v_3(s,\cdot) \bigr](\xi) \, \ud s \biggr\|_{L^2_\xi} \\ 
        &\lesssim \int_0^t \Bigl\| \whatcalF\bigl[ w_1^{\musSharp}(s,\cdot) \overline{v_2(s,\cdot)} v_3(s,\cdot) \bigr](\xi) \Bigr\|_{L^2_\xi} \, \ud s 
        \lesssim \int_0^t \bigl\| w_1^{\musSharp}(s,y) \overline{v_2(s,y)} v_3(s,y) \bigr\|_{L^2_y} \, \ud s \\ 
        &\lesssim \int_0^t \|w_1^{\musSharp}(s)\|_{L^2_y} \|v_2(s)\|_{L^\infty_y} \|v_3(s)\|_{L^\infty_y} \, \ud s 
        \lesssim \int_0^t \varepsilon \js^{\delta} \cdot \varepsilon^2 \js^{-1} \, \ud s
        \lesssim \varepsilon^3 \jt^\delta.
    \end{aligned}
\end{equation*}

\noindent {\it Bounds for $\calJ_{2,1,k,2}^{\delta_0}(s,\xi)$, $1 \leq k \leq 2$:}
Again we provide the details for the case $k=1$, while the analogous case $k=2$ is left to the reader. By Plancherel's theorem, \eqref{equ:consequences_aux_KG_disp_decay}, \eqref{equ:consequences_aux_KG_energy_bounds} we obtain
\begin{equation*}
    \begin{aligned}
        &\biggl\| \int_0^t e^{-i\xi\theta(s)} \calJ_{2,1,1,2}(s,\xi) \, \ud s \biggr\|_{L^2_\xi} = \biggl\| \int_0^t (2\pi)^{\frac52} \frakb(\xi) e^{-i\xi\theta(s)} e^{-is(\jxi+\ulell\xi)} \whatcalF\bigl[ \tilde{v}_1(s,\cdot) \overline{v_2(s,\cdot)} v_3(s,\cdot) \bigr](\xi) \, \ud s \biggr\|_{L^2_\xi} \\ 
        &\lesssim \int_0^t \Bigl\| \whatcalF\bigl[ \tilde{v}_1(s,\cdot) \overline{v_2(s,\cdot)} v_3(s,\cdot) \bigr](\xi) \Bigr\|_{L^2_\xi} \, \ud s 
        \lesssim \int_0^t \bigl\| \tilde{v}_1(s,y) \overline{v_2(s,y)} v_3(s,y) \bigr\|_{L^2_y} \, \ud s \\ 
        &\lesssim \int_0^t \|\tilde{v}_1(s)\|_{L^2_y} \|v_2(s)\|_{L^\infty_y} \|v_3(s)\|_{L^\infty_y} \, \ud s 
        \lesssim \int_0^t \varepsilon \cdot \varepsilon^2 \js^{-1} \, \ud s
        \lesssim \varepsilon^3 \log(\jt).
    \end{aligned}
\end{equation*}

\noindent {\it Bound for $\calJ_{2,3}^{\delta_0}(s,\xi)$:}
In the same manner as for the preceding term, we infer by \eqref{equ:consequences_aux_KG_disp_decay}, \eqref{equ:consequences_aux_KG_energy_bounds} that
\begin{equation*}
    \begin{aligned}
        &\biggl\| \int_0^t e^{-i\xi\theta(s)} \calJ_{2,3}^{\delta_0}(s,\xi) \, \ud s \biggr\|_{L^2_\xi} \\ 
        &= \biggl\| \int_0^t e^{-i\xi\theta(s)} (2\pi)^{\frac52} (\pxi \frakb)(\xi) \, e^{-is(\jxi+\ulell\xi)} \, \jxi \whatcalF\bigl[ v_1(s,\cdot) \overline{v_2(s,\cdot)} v_3(s,\cdot) \bigr](\xi) \, \ud s \biggr\|_{L^2_\xi} \\ 
        &\lesssim \int_0^t \bigl\| \jxi \whatcalF\bigl[ v_1(s,\cdot) \overline{v_2(s,\cdot)} v_3(s,\cdot) \bigr](\xi) \bigr\|_{L^2_\xi} \, \ud s  
        \lesssim \int_0^t \bigl\| v_1(s,y) \overline{v_2(s,y)} v_3(s,y) \bigr\|_{H^1_y} \, \ud s \\
        &\lesssim \int_0^t \|v_1(s)\|_{W^{1,\infty}_y} \|v_2(s)\|_{W^{1,\infty}_y} \|v_3(s)\|_{H^1_y} \, \ud s \lesssim \int_0^t \varepsilon^2 \js^{-1} \cdot \varepsilon \, \ud s \lesssim \varepsilon^3 \log(\jt).
    \end{aligned}
\end{equation*}

\noindent {\it Bound for $\calJ_{2,\mathrm{comb}, 1}^{\delta_0}(s,\xi)$:}
Using \eqref{equ:consequences_aux_KG_disp_decay} and \eqref{equ:weighted_energies_proof_wjmusSharp_L2y_growth}, we obtain the acceptable bound
\begin{equation*}
    \begin{aligned}
        &\biggl\| \int_0^t e^{-i\xi\theta(s)} \calJ_{2,\mathrm{comb}, 1}^{\delta_0}(s,\xi) \, \ud s \biggr\|_{L^2_\xi} \\
        &= \biggl\| \int_0^t e^{-i\xi\theta(s)} (2\pi)^{\frac52} \frakb(\xi) e^{-is(\jxi+\ulell\xi)} \whatcalF\bigl[ v_1(s,\cdot) \overline{v_2(s,\cdot)} w_3^{\musSharp}(s,\cdot) \bigr](\xi) \biggr\|_{L^2_\xi} \\ 
        &\lesssim \int_0^t \bigl\| v_1(s,y) \overline{v_2(s,y)} w_3^{\musSharp}(s,y) \bigr\|_{L^2_y} \, \ud s 
        \lesssim \int_0^t \|v_1(s)\|_{L^\infty_y} \|v_2(s)\|_{L^\infty_y} \|w_3^{\musSharp}(s)\|_{L^2_y} \, \ud s \\ 
        &\lesssim \int_0^t \varepsilon^2 \js^{-1} \cdot \varepsilon \js^\delta \, \ud s \lesssim \varepsilon^3 \jt^\delta.
    \end{aligned}
\end{equation*}

\noindent {\it Bound for $\calJ_{2,\mathrm{comb}, 2}^{\delta_0}(s,\xi)$:}
By Plancherel's theorem and \eqref{equ:consequences_aux_KG_disp_decay}, \eqref{equ:weighted_energies_proof_vjmusSharp_decay}, \eqref{equ:weighted_energies_proof_wj_H1y_growth}, \eqref{equ:weighted_energies_proof_wjmusSharp_L2y_growth}, we arrive at the acceptable bound
\begin{equation*}
    \begin{aligned}
        &\biggl\| \int_0^t e^{-i\xi\theta(s)} \calJ_{2,\mathrm{comb}, 2}^{\delta_0}(s,\xi) \, \ud s \biggr\|_{L^2_\xi} \\
        &\lesssim \int_0^t \bigl\| \jxi \whatcalF\bigl[ v_1(s,\cdot) \overline{v_2(s,\cdot)} w_3(s,\cdot) \bigr](\xi) \bigr\|_{L^2_\xi} \, \ud s 
        + \int_0^t \bigl\| \whatcalF\bigl[ v_1^{\musSharp}(s,\cdot) \overline{v_2(s,\cdot)} w_3(s,\cdot) \bigr](\xi) \bigr\|_{L^2_\xi} \, \ud s \\
        &\quad + \int_0^t \bigl\| \whatcalF\bigl[ v_1(s,\cdot) \overline{v_2^{\musSharp}(s,\cdot)} w_3(s,\cdot) \bigr](\xi) \bigr\|_{L^2_\xi} \, \ud s 
        + \int_0^t \bigl\| \whatcalF\bigl[ v_1(s,\cdot) \overline{v_2(s,\cdot)} w_3^{\musSharp}(s,\cdot) \bigr](\xi) \bigr\|_{L^2_\xi} \, \ud s \\ 
        &\lesssim \int_0^t \|v_1(s)\|_{W^{1,\infty}_y} \|v_2(s)\|_{W^{1,\infty}_y} \|w_3(s)\|_{H^1_y} \, \ud s  
        + \int_0^t \|v_1^{\musSharp}(s)\|_{L^\infty_y} \|v_2(s)\|_{L^\infty_y} \|w_3(s)\|_{L^2_y} \, \ud s \\ 
        &\quad + \int_0^t \|v_1(s)\|_{L^\infty_y} \|v_2^{\musSharp}(s)\|_{L^\infty_y} \|w_3(s)\|_{L^2_y} \, \ud s 
        + \int_0^t \|v_1(s)\|_{L^\infty_y} \|v_2(s)\|_{L^\infty_y} \|w_3^{\musSharp}(s)\|_{L^2_y} \, \ud s \\ 
        &\lesssim \int_0^t \varepsilon^2 \js^{-1} \cdot \varepsilon \js^\delta \, \ud s \lesssim \varepsilon^3 \jt^\delta.
    \end{aligned}
\end{equation*}

\medskip 
\noindent {\it Bound for $\calK_{2,1,3,1}^{\delta_0}(t,\xi)$:}
An application of Plancherel's theorem together with \eqref{equ:consequences_aux_KG_disp_decay}, \eqref{equ:consequences_aux_KG_energy_bounds} gives the sufficient bound
\begin{equation*}
    \begin{aligned}
        &\bigl\| \calK_{2,1,3,1}^{\delta_0}(t,\xi) \bigr\|_{L^2_\xi} = \Bigl\| t \cdot e^{-i\xi\theta(t)} (2\pi)^{\frac52} \frakb(\xi) \, e^{-it(\jxi+\ulell\xi)}  \whatcalF\bigl[ v_1(t,\cdot) \overline{v_2(t,\cdot)} \tilde{v}_3(t,\cdot) \bigr](\xi) \Bigr\|_{L^2_\xi} \\ 
        &\lesssim t \cdot \bigl\|  \whatcalF\bigl[ v_1(t,\cdot) \overline{v_2(t,\cdot)} \tilde{v}_3(t,\cdot) \bigr](\xi) \bigr\|_{L^2_\xi} \lesssim t \cdot \bigl\| v_1(t,\cdot) \overline{v_2(t,\cdot)} \tilde{v}_3(t,\cdot) \bigr\|_{L^2_y} \\ 
        &\lesssim t \cdot \|v_1(t)\|_{L^2_y} \|v_2(t)\|_{L^\infty_y} \|\tilde{v}_3(t)\|_{L^\infty_y} \lesssim t \cdot \varepsilon \cdot \varepsilon^2 \jt^{-1} \lesssim \varepsilon^3.
    \end{aligned}
\end{equation*}

\noindent {\it Bound for $\calK_{2,1,3,2}^{\delta_0}(t,\xi)$:}
As for the preceding term, we use \eqref{equ:consequences_aux_KG_disp_decay} and \eqref{equ:consequences_aux_KG_energy_bounds} to infer the sufficient bound
\begin{equation*}
    \begin{aligned}
        &\bigl\| \calK_{2,1,3,2}^{\delta_0}(t,\xi) \bigr\|_{L^2_\xi} = \biggl\| (2\pi)^{\frac52} \frakb(\xi) \int_0^t e^{-i\xi\theta(s)} e^{-is(\jxi+\ulell\xi)} \whatcalF\bigl[ v_1(s,\cdot) \overline{v_2(s,\cdot)} \tilde{v}_3(s,\cdot) \bigr](\xi) \, \ud s \biggr\|_{L^2_\xi} \\ 
        &\lesssim \int_0^t \Bigl\| \whatcalF\bigl[ v_1(s,\cdot) \overline{v_2(s,\cdot)} \tilde{v}_3(s,\cdot) \bigr](\xi) \Bigr\|_{L^2_\xi} \, \ud s \lesssim \int_0^t \bigl\| v_1(s,y) \overline{v_2(s,y)} \tilde{v}_3(s,y) \bigr\|_{L^2_y} \, \ud s \\ 
        &\lesssim \int_0^t \|v_1(s)\|_{L^2_y} \|v_2(s)\|_{L^\infty_y} \|\tilde{v}_3(s)\|_{L^\infty_y} \, \ud s \lesssim \int_0^t \varepsilon \cdot \varepsilon^2 \js^{-1} \, \ud s \lesssim \varepsilon^3 \log(\jt).
    \end{aligned}
\end{equation*}

\noindent {\it Bound for $\calK_{2,1,3,3}^{\delta_0}(t,\xi)$:}
Combining \eqref{equ:consequences_qdot_minus_ulell_decay}, \eqref{equ:consequences_aux_KG_disp_decay}, \eqref{equ:consequences_aux_KG_energy_bounds} leads to the acceptable bound
\begin{equation*}
    \begin{aligned}
        &\bigl\| \calK_{2,1,3,3}^{\delta_0}(t,\xi) \bigr\|_{L^2_\xi} \\
        &= \biggl\| (2\pi)^{\frac52} \frakb(\xi) \int_0^t s \cdot (\dot{q}(s)-\ulell) \cdot e^{-i\xi\theta(s)} e^{-is(\jxi+\ulell\xi)} \cdot \xi \cdot \whatcalF\bigl[ v_1(s,\cdot) \overline{v_2(s,\cdot)} \tilde{v}_3(s,\cdot) \bigr](\xi) \, \ud s \biggr\|_{L^2_\xi} \\ 
        &\lesssim \int_0^t s \cdot |\dot{q}(s)-\ulell| \cdot \Bigl\| \xi \cdot \whatcalF\bigl[ v_1(s,\cdot) \overline{v_2(s,\cdot)} \tilde{v}_3(s,\cdot) \bigr](\xi) \Bigr\|_{L^2_\xi} \, \ud s \\ 
        &\lesssim \int_0^t s \cdot |\dot{q}(s)-\ulell| \cdot \bigl\| v_1(s,y) \overline{v_2(s,y)} \tilde{v}_3(s,y) \bigr\|_{H^1_y} \, \ud s \\
        &\lesssim \int_0^t s \cdot |\dot{q}(s)-\ulell| \cdot \|v_1(s)\|_{H^1_y} \|v_2(s)\|_{W^{1,\infty}_y} \|\tilde{v}_3(s)\|_{W^{1,\infty}_y} \, \ud s \\ &\lesssim \int_0^t s \cdot \varepsilon \js^{-1+\delta} \cdot \varepsilon \cdot \varepsilon^2 \js^{-1} \, \ud s \lesssim \varepsilon^4 \jt^\delta.
    \end{aligned}
\end{equation*}

\noindent {\it Bound for $\calK_{2,1,3,4,1}^{\delta_0}(t,\xi)$:}
Analogously to the preceding term, the combination of \eqref{equ:consequences_qdot_minus_ulell_decay}, \eqref{equ:consequences_aux_KG_disp_decay}, \eqref{equ:consequences_aux_KG_energy_bounds} gives the acceptable bound
\begin{equation*}
    \begin{aligned}
        &\bigl\| \calK_{2,1,3,4,1}^{\delta_0}(t,\xi) \bigr\|_{L^2_\xi} \\
        &= \biggl\| (2\pi)^{\frac52} \frakb(\xi) \int_0^t s \cdot (\dot{q}(s)-\ulell) \cdot e^{-i\xi\theta(s)} e^{-is(\jxi+\ulell\xi)} \whatcalF\bigl[ (\py v_1)(s,\cdot) \overline{v_2(s,\cdot)} \tilde{v}_3(s,\cdot) \bigr](\xi) \, \ud s \biggr\|_{L^2_\xi} \\
        &\lesssim \int_0^t s \cdot |\dot{q}(s)-\ulell| \cdot \Bigl\| \whatcalF\bigl[ (\py v_1)(s,\cdot) \overline{v_2(s,\cdot)} \tilde{v}_3(s,\cdot) \bigr](\xi) \Bigr\|_{L^2_\xi} \, \ud s \\
        &\lesssim \int_0^t s \cdot |\dot{q}(s)-\ulell| \cdot \bigl\| (\py v_1)(s,y) \overline{v_2(s,y)} \tilde{v}_3(s,y) \bigr\|_{L^2_y} \, \ud s \\ 
        &\lesssim \int_0^t s \cdot |\dot{q}(s)-\ulell| \cdot \|(\py v_1)(s)\|_{L^2_y} \|v_2(s)\|_{L^\infty_y} \|\tilde{v}_3(s)\|_{L^\infty_y} \, \ud s \\ 
        &\lesssim \int_0^t s \cdot \varepsilon \js^{-1+\delta} \cdot \varepsilon \cdot \varepsilon^2 \js^{-1} \, \ud s \lesssim \varepsilon^4 \jt^\delta.
    \end{aligned}
\end{equation*}

\noindent {\it Bound for $\calK_{2,1,3,4,2}^{\delta_0}(t,\xi)$:}
Using Plancherel's theorem and the mapping property~\eqref{equ:mapping_property_calFulellsh_L2} along with the decay estimates \eqref{equ:consequences_dispersive_decay_ulPe_radiation},  \eqref{equ:consequences_aux_KG_disp_decay} we deduce the sufficient bound 
\begin{equation*}
    \begin{aligned}
        &\bigl\| \calK_{2,1,3,4,2}^{\delta_0}(t,\xi) \bigr\|_{L^2_\xi} \\
        &= \biggl\| (2\pi)^{\frac52} \frakb(\xi) \int_0^t s \cdot e^{-i\xi\theta(s)} e^{-is(\jxi+\ulell\xi)} \\
        &\qquad \qquad \qquad \times \whatcalF\biggl[ \whatcalF^{-1}\Bigl[ \jxione^{-1} \frakb_1(\xi_1) \calFulellsh\bigl[ \calQ_\ulell\bigl( \usubeone(s) \bigr) \bigr](\xi_1) \Bigr](\cdot) \, \overline{v_2(s,\cdot)} \, \tilde{v}_3(s,\cdot) \biggr](\xi) \, \ud s \biggr\|_{L^2_\xi} \\ 
        &\lesssim \int_0^t s \cdot \biggl\| \whatcalF\biggl[ \whatcalF^{-1}\Bigl[ \jxione^{-1} \frakb_1(\xi_1) \calFulellsh\bigl[ \calQ_\ulell\bigl( \usubeone(s) \bigr)\bigr](\xi_1) \Bigr](\cdot) \, \overline{v_2(s,\cdot)} \, \tilde{v}_3(s,\cdot) \biggr](\xi) \biggr\|_{L^2_\xi} \, \ud s \\ 
        &\lesssim \int_0^t s \cdot \Bigl\| \whatcalF^{-1}\Bigl[ \jxione^{-1} \frakb_1(\xi_1) \calFulellsh\bigl[ \calQ_\ulell\bigl( \usubeone(s) \bigr)\bigr](\xi_1) \Bigr](y) \, \overline{v_2(s,y)} \, \tilde{v}_3(s,y) \Bigr\|_{L^2_y} \, \ud s \\ 
        &\lesssim \int_0^t s \cdot \bigl\| \calQ_\ulell\bigl( \usubeone(s) \bigr) \bigr\|_{L^2_y} \|v_2(s)\|_{L^\infty_y} \|\tilde{v}_3(s)\|_{L^\infty_y} \, \ud s \\ 
        &\lesssim \int_0^t s \cdot \|\alpha(\ulg y)\|_{L^2_y} \|\usubeone(s)\|_{L^\infty_y}^2 \|v_2(s)\|_{L^\infty_y} \|\tilde{v}_3(s)\|_{L^\infty_y} \, \ud s \\ 
        &\lesssim \int_0^t s \cdot \varepsilon^4 \js^{-2} \, \ud s \lesssim \varepsilon^4 \log(\jt).
    \end{aligned}
\end{equation*}

\noindent {\it Bound for $\calK_{2,1,3,4,3}^{\delta_0}(t,\xi)$:}
Similarly to the preceding term, we use Plancherel's theorem and the mapping property~\eqref{equ:mapping_property_calFulellsh_L2} along with the decay estimates \eqref{equ:consequences_dispersive_decay_ulPe_radiation}, \eqref{equ:consequences_aux_KG_disp_decay} and the uniform-in-time energy bound \eqref{equ:consequences_energy_bound_radiation}, to conclude the sufficient bound
\begin{equation*}
    \begin{aligned}
        &\bigl\| \calK_{2,1,3,4,3}^{\delta_0}(t,\xi) \bigr\|_{L^2_\xi} \\
        &= \biggl\| (2\pi)^{\frac52} \frakb(\xi) \int_0^t s \cdot e^{-i\xi\theta(s)} e^{-is(\jxi+\ulell\xi)} \\
        &\qquad \qquad \qquad \times \whatcalF\biggl[ \whatcalF^{-1}\Bigl[ \jxione^{-1} \frakb_1(\xi_1) \calFulellsh\bigl[ {\textstyle \frac16} \usubeone(s)^3 \bigr](\xi_1)\Bigr](\cdot) \, \overline{v_2(s,\cdot)} \, \tilde{v}_3(s,\cdot) \biggr](\xi) \, \ud s \biggr\|_{L^2_\xi} \\ 
        &\lesssim \int_0^t s \cdot \biggl\| \whatcalF\biggl[ \whatcalF^{-1}\Bigl[ \jxione^{-1} \frakb_1(\xi_1) \calFulellsh\bigl[ {\textstyle \frac16} \usubeone(s)^3 \bigr](\xi_1)\Bigr](\cdot) \, \overline{v_2(s,\cdot)} \, \tilde{v}_3(s,\cdot) \biggr](\xi) \biggr\|_{L^2_\xi} \, \ud s \\ 
        &\lesssim \int_0^t s \cdot \Bigl\| \whatcalF^{-1}\Bigl[ \jxione^{-1} \frakb_1(\xi_1) \calFulellsh\bigl[ {\textstyle \frac16} \usubeone(s)^3 \bigr](\xi_1)\Bigr](y) \, \overline{v_2(s,y)} \, \tilde{v}_3(s,y) \Bigr\|_{L^2_y} \, \ud s \\ 
        &\lesssim \int_0^t s \cdot \bigl\|\usubeone(s)^3 \bigr\|_{L^2_y} \|v_2(s)\|_{L^\infty_y} \|\tilde{v}_3(s)\|_{L^\infty_y} \, \ud s \\
        &\lesssim \int_0^t s \cdot \|\usubeone(s)\|_{L^2_y} \|\usubeone(s)\|_{L^\infty_y}^2 \|v_2(s)\|_{L^\infty_y} \|\tilde{v}_3(s)\|_{L^\infty_y} \, \ud s \\
        &\lesssim \int_0^t s \cdot \varepsilon \cdot \varepsilon^4 \js^{-2} \, \ud s \lesssim \varepsilon^5 \log(\jt).
    \end{aligned}
\end{equation*}

\noindent {\it Bound for $\calK_{2,1,3,4,4}^{\delta_0}(t,\xi)$:}
Here we use Plancherel's theorem, the decay estimate \eqref{equ:consequences_aux_KG_disp_decay} and the $L^2_{\xi_1}$-decay bound \eqref{equ:consequences_wtilcalR_L2xi_bound} for $\widetilde{\calR}(s,\xi_1)$ to infer the sufficient estimate
\begin{equation*}
    \begin{aligned}
        \bigl\| \calK_{2,1,3,4,4}^{\delta_0}(t,\xi) \bigr\|_{L^2_\xi} 
        &= \biggl\| (2\pi)^{\frac52} \frakb(\xi) \int_0^t s \cdot e^{-i\xi\theta(s)} e^{-is(\jxi+\ulell\xi)} \\
        &\qquad \qquad \qquad \times \whatcalF\biggl[ \whatcalF^{-1}\Bigl[ \jxione^{-1} \frakb_1(\xi_1) \widetilde{\calR}(s,\xi_1) \Bigr](\cdot) \, \overline{v_2(s,\cdot)} \, \tilde{v}_3(s,\cdot) \biggr](\xi) \, \ud s \biggr\|_{L^2_\xi} \\ 
        &\lesssim \int_0^t s \cdot \biggl\| \whatcalF\biggl[ \whatcalF^{-1}\Bigl[ \jxione^{-1} \frakb_1(\xi_1) \widetilde{\calR}(s,\xi_1) \Bigr](\cdot) \, \overline{v_2(s,\cdot)} \, \tilde{v}_3(s,\cdot) \biggr](\xi) \biggr\|_{L^2_\xi} \, \ud s \\ 
        &\lesssim \int_0^t s \cdot \Bigl\| \whatcalF^{-1}\Bigl[ \jxione^{-1} \frakb_1(\xi_1) \widetilde{\calR}(s,\xi_1) \Bigr](y) \, \overline{v_2(s,y)} \, \tilde{v}_3(s,y) \Bigr\|_{L^2_y} \, \ud s \\ 
        &\lesssim \int_0^t s \cdot \Bigl\| \whatcalF^{-1}\Bigl[ \jxione^{-1} \frakb_1(\xi_1) \widetilde{\calR}(s,\xi_1) \Bigr](y) \Bigr\|_{L^2_y} \|v_2(s)\|_{L^\infty_y} \|\tilde{v}_3(s)\|_{L^\infty_y} \, \ud s \\
        &\lesssim \int_0^t s \cdot \bigl\| \widetilde{\calR}(s,\xi_1) \bigr\|_{L^2_{\xi_1}} \|v_2(s)\|_{L^\infty_y} \|\tilde{v}_3(s)\|_{L^\infty_y} \, \ud s \\
        &\lesssim \int_0^t s \cdot \varepsilon^2 \js^{-\frac32+\delta} \cdot \varepsilon^2 \js^{-1} \, \ud s \lesssim \varepsilon^4.
    \end{aligned}
\end{equation*}

\medskip 
\noindent \underline{Case 1.2: Contribution to \eqref{equ:pxi_cubic_with_jxi_reduced_claim}.}
Now our goal is to prove the following (still slowly growing) weighted energy bound with additional Sobolev regularity for $0 \leq t \leq T$,
\begin{equation*}
    \biggl\| \jxi^2 \int_0^t e^{-i\xi\theta(s)} \pxi \calI_{2, \mathrm{schem}}^{\delta_0}(s,\xi) \, \ud s \biggr\|_{L^2_\xi} \lesssim \varepsilon^3 \jt^{2\delta}.
\end{equation*}
We write 
\begin{equation}
    \begin{aligned}
     \biggl\| \jxi^2 \int_0^t e^{-i\xi\theta(s)} \pxi \calI_{2, \mathrm{schem}}^{\delta_0}(s,\xi) \, \ud s \biggr\|_{L^2_\xi} &= \biggl\| \jxi \int_0^t e^{-i\xi\theta(s)} \jxi \pxi \calI_{2, \mathrm{schem}}^{\delta_0}(s,\xi) \, \ud s \biggr\|_{L^2_\xi},
     \end{aligned}
\end{equation}
so that we can work with the previously established decomposition of $\jxi \pxi \calI_{2, \mathrm{schem}}^{\delta_0}(s,\xi)$, see Figure~\ref{figure:decomposition_weighted_energy_terms}.
Again, it is useful to rewrite these expressions in terms of the auxiliary linear Klein-Gordon evolutions $v_j(s,y)$, $\tilde{v}_j(s,y)$, $v_j^{\musSharp}(s,y)$, $w_j(s,y)$, and $w_j^{\musSharp}(s,y)$ defined in \eqref{equ:weighted_energies_proof_vj_def}, \eqref{equ:weighted_energies_proof_vj_musSharp_def}, \eqref{equ:weighted_energies_proof_tildevj_def}, \eqref{equ:weighted_energies_proof_wj_def}, and \eqref{equ:weighted_energies_proof_wj_musSharp_def}.

\medskip 

\noindent {\it Bounds for $\jxi \calJ_{2,1,k,1}^{\delta_0}(s,\xi)$, $1 \leq k \leq 2$:}
We provide the details for $k=1$. The case $k=2$ can be treated analogously. By Plancherel's theorem, \eqref{equ:consequences_aux_KG_disp_decay}, and \eqref{equ:weighted_energies_proof_wjmusSharp_H2y_growth} we obtain the acceptable bound
\begin{equation*}
    \begin{aligned}
        &\biggl\| \jxi \int_0^t e^{-i\xi\theta(s)} \calJ_{2,1,1,1}(s,\xi) \, \ud s \biggr\|_{L^2_\xi} \\
        &= \biggl\| \int_0^t (2\pi)^{\frac52} \frakb(\xi) e^{-i\xi\theta(s)} e^{-is(\jxi+\ulell\xi)} \jxi \whatcalF\bigl[ w_1^{\musSharp}(s,\cdot) \overline{v_2(s,\cdot)} v_3(s,\cdot) \bigr](\xi) \, \ud s \biggr\|_{L^2_\xi} \\ 
        &\lesssim \int_0^t \Bigl\| \jxi \whatcalF\bigl[ w_1^{\musSharp}(s,\cdot) \overline{v_2(s,\cdot)} v_3(s,\cdot) \bigr](\xi) \Bigr\|_{L^2_\xi} \, \ud s 
        \lesssim \int_0^t \bigl\| w_1^{\musSharp}(s,y) \overline{v_2(s,y)} v_3(s,y) \bigr\|_{H^1_y} \, \ud s \\ 
        &\lesssim \int_0^t \|w_1^{\musSharp}(s)\|_{H^1_y} \|v_2(s)\|_{W^{1,\infty}_y} \|v_3(s)\|_{W^{1,\infty}_y} \, \ud s 
        \lesssim \int_0^t \varepsilon \js^{2\delta} \cdot \varepsilon^2 \js^{-1} \, \ud s
        \lesssim \varepsilon^3 \jt^{2\delta}.
    \end{aligned}
\end{equation*}

\noindent {\it Bounds for $\jxi \calJ_{2,1,k,2}^{\delta_0}(s,\xi)$, $1 \leq k \leq 2$:}
Again we provide the details for the case $k=1$, while the analogous case $k=2$ is left to the reader. Using Plancherel's theorem, \eqref{equ:consequences_aux_KG_disp_decay}, and \eqref{equ:consequences_aux_KG_energy_bounds}, we conclude the sufficient bound
\begin{equation*}
    \begin{aligned}
        &\biggl\| \jxi \int_0^t e^{-i\xi\theta(s)} \calJ_{2,1,1,2}(s,\xi) \, \ud s \biggr\|_{L^2_\xi} \\ 
        &= \biggl\| \int_0^t (2\pi)^{\frac52} \frakb(\xi) e^{-i\xi\theta(s)} e^{-is(\jxi+\ulell\xi)} \jxi \whatcalF\bigl[ \tilde{v}_1(s,\cdot) \overline{v_2(s,\cdot)} v_3(s,\cdot) \bigr](\xi) \, \ud s \biggr\|_{L^2_\xi} \\ 
        &\lesssim \int_0^t \Bigl\| \jxi \whatcalF\bigl[ \tilde{v}_1(s,\cdot) \overline{v_2(s,\cdot)} v_3(s,\cdot) \bigr](\xi) \Bigr\|_{L^2_\xi} \, \ud s 
        \lesssim \int_0^t \bigl\| \tilde{v}_1(s,y) \overline{v_2(s,y)} v_3(s,y) \bigr\|_{H^1_y} \, \ud s \\ 
        &\lesssim \int_0^t \|\tilde{v}_1(s)\|_{H^1_y} \|v_2(s)\|_{W^{1,\infty}_y} \|v_3(s)\|_{W^{1,\infty}_y} \, \ud s 
        \lesssim \int_0^t \varepsilon \cdot \varepsilon^2 \js^{-1} \, \ud s
        \lesssim \varepsilon^3 \log(\jt).
    \end{aligned}
\end{equation*}

\noindent {\it Bound for $\jxi \calJ_{2,3}^{\delta_0}(s,\xi)$:}
Here we obtain a sufficient estimate using Plancherel's theorem and the bounds \eqref{equ:consequences_aux_KG_disp_decay}, \eqref{equ:consequences_aux_KG_energy_bounds}, 
\begin{equation*}
    \begin{aligned}
        &\biggl\| \jxi \int_0^t e^{-i\xi\theta(s)} \calJ_{2,3}^{\delta_0}(s,\xi) \, \ud s \biggr\|_{L^2_\xi} \\ 
        &= \biggl\| \jxi \int_0^t e^{-i\xi\theta(s)} (2\pi)^{\frac52} (\pxi \frakb)(\xi) \, e^{-is(\jxi+\ulell\xi)} \, \jxi \whatcalF\bigl[ v_1(s,\cdot) \overline{v_2(s,\cdot)} v_3(s,\cdot) \bigr](\xi) \, \ud s \biggr\|_{L^2_\xi} \\ 
        &\lesssim \int_0^t \bigl\| \jxi^2 \whatcalF\bigl[ v_1(s,\cdot) \overline{v_2(s,\cdot)} v_3(s,\cdot) \bigr](\xi) \bigr\|_{L^2_\xi} \, \ud s 
        \lesssim \int_0^t \bigl\| v_1(s,y) \overline{v_2(s,y)} v_3(s,y) \bigr\|_{H^2_y} \, \ud s \\
        &\lesssim \int_0^t \|v_1(s)\|_{H^2_y} \|v_2(s)\|_{W^{1,\infty}_y} \|v_3(s)\|_{W^{1,\infty}_y} \, \ud s + \bigl\{\text{similar terms}\bigr\} \\
        &\lesssim \int_0^t \varepsilon \js^\delta \cdot \varepsilon^2 \js^{-1} \, \ud s \lesssim \varepsilon^3 \jt^\delta.
    \end{aligned}
\end{equation*}

\noindent {\it Bound for $\jxi \calJ_{2,\mathrm{comb}, 1}^{\delta_0}(s,\xi)$:}
For this term we obtain by Plancherel's theorem and the bounds \eqref{equ:consequences_aux_KG_disp_decay}, \eqref{equ:weighted_energies_proof_wjmusSharp_H2y_growth}, the acceptable estimate
\begin{equation*}
    \begin{aligned}
        &\biggl\| \jxi \int_0^t e^{-i\xi\theta(s)} \calJ_{2,\mathrm{comb}, 1}^{\delta_0}(s,\xi) \, \ud s \biggr\|_{L^2_\xi} \\ 
        &= \biggl\| \jxi \int_0^t e^{-i\xi\theta(s)} (2\pi)^{\frac52} \frakb(\xi) e^{-is(\jxi+\ulell\xi)} \whatcalF\bigl[ v_1(s,\cdot) \overline{v_2(s,\cdot)} w_3^{\musSharp}(s,\cdot) \bigr](\xi) \biggr\|_{L^2_\xi} \\ 
        &\lesssim \int_0^t \bigl\| v_1(s,y) \overline{v_2(s,y)} w_3^{\musSharp}(s,y) \bigr\|_{H^1_y} \, \ud s 
        \lesssim \int_0^t \|v_1(s)\|_{W^{1,\infty}_y} \|v_2(s)\|_{W^{1,\infty}_y} \|w_3^{\musSharp}(s)\|_{H^1_y} \, \ud s \\ 
        &\lesssim \int_0^t \varepsilon^2 \js^{-1} \cdot \varepsilon \js^{2\delta} \, \ud s \lesssim \varepsilon^3 \jt^{2\delta}.
    \end{aligned}
\end{equation*}

\medskip 
\noindent {\it Bound for $\jxi \calK_{2,1,3,1}^{\delta_0}(t,\xi)$:}
Here an application of Plancherel's theorem combined with \eqref{equ:consequences_aux_KG_disp_decay}, \eqref{equ:consequences_aux_KG_energy_bounds} gives the sufficient bound
\begin{equation*}
    \begin{aligned}
        &\bigl\| \jxi \calK_{2,1,3,1}^{\delta_0}(t,\xi) \bigr\|_{L^2_\xi} = \Bigl\| t \cdot e^{-i\xi\theta(t)} (2\pi)^{\frac52} \frakb(\xi) \, e^{-it(\jxi+\ulell\xi)} \jxi  \whatcalF\bigl[ v_1(t,\cdot) \overline{v_2(t,\cdot)} \tilde{v}_3(t,\cdot) \bigr](\xi) \Bigr\|_{L^2_\xi} \\ 
        &\lesssim t \cdot \bigl\| \jxi \whatcalF\bigl[ v_1(t,\cdot) \overline{v_2(t,\cdot)} \tilde{v}_3(t,\cdot) \bigr](\xi) \bigr\|_{L^2_\xi} \lesssim t \cdot \bigl\| v_1(t,\cdot) \overline{v_2(t,\cdot)} \tilde{v}_3(t,\cdot) \bigr\|_{H^1_y} \\ 
        &\lesssim t \cdot \|v_1(t)\|_{H^1_y} \|v_2(t)\|_{W^{1,\infty}_y} \|\tilde{v}_3(t)\|_{W^{1,\infty}_y} \lesssim t \cdot \varepsilon \cdot \varepsilon^2 \jt^{-1} \lesssim \varepsilon^3.
    \end{aligned}
\end{equation*}

\noindent {\it Bound for $\jxi \calK_{2,1,3,2}^{\delta_0}(t,\xi)$:}
Similarly to the preceding term, using Plancherel's theorem and  \eqref{equ:consequences_aux_KG_disp_decay}, \eqref{equ:consequences_aux_KG_energy_bounds}, we obtain the sufficient bound
\begin{equation*}
    \begin{aligned}
        &\bigl\| \jxi \calK_{2,1,3,2}^{\delta_0}(t,\xi) \bigr\|_{L^2_\xi} = \biggl\| (2\pi)^{\frac52} \frakb(\xi) \int_0^t e^{-i\xi\theta(s)} e^{-is(\jxi+\ulell\xi)} \jxi \whatcalF\bigl[ v_1(s,\cdot) \overline{v_2(s,\cdot)} \tilde{v}_3(s,\cdot) \bigr](\xi) \, \ud s \biggr\|_{L^2_\xi} \\ 
        &\lesssim \int_0^t \Bigl\| \jxi \whatcalF\bigl[ v_1(s,\cdot) \overline{v_2(s,\cdot)} \tilde{v}_3(s,\cdot) \bigr](\xi) \Bigr\|_{L^2_\xi} \, \ud s \lesssim \int_0^t \bigl\| v_1(s,y) \overline{v_2(s,y)} \tilde{v}_3(s,y) \bigr\|_{H^1_y} \, \ud s \\ 
        &\lesssim \int_0^t \|v_1(s)\|_{H^1_y} \|v_2(s)\|_{W^{1,\infty}_y} \|\tilde{v}_3(s)\|_{W^{1,\infty}_y} \, \ud s \lesssim \int_0^t \varepsilon \cdot \varepsilon^2 \js^{-1} \, \ud s \lesssim \varepsilon^3 \log(\jt).
    \end{aligned}
\end{equation*}

\noindent {\it Bound for $\jxi \calK_{2,1,3,3}^{\delta_0}(t,\xi)$:}
The bound for this term is again tight. 
Using \eqref{equ:consequences_qdot_minus_ulell_decay}, \eqref{equ:consequences_aux_KG_disp_decay}, and \eqref{equ:consequences_aux_KG_energy_bounds}, we obtain via Plancherel's theorem the acceptable bound
\begin{equation*}
    \begin{aligned}
        &\bigl\| \jxi \calK_{2,1,3,3}^{\delta_0}(t,\xi) \bigr\|_{L^2_\xi} \\
        &= \biggl\| (2\pi)^{\frac52} \frakb(\xi) \int_0^t s \cdot (\dot{q}(s)-\ulell) \cdot e^{-i\xi\theta(s)} e^{-is(\jxi+\ulell\xi)} \cdot \jxi \xi \cdot \whatcalF\bigl[ v_1(s,\cdot) \overline{v_2(s,\cdot)} \tilde{v}_3(s,\cdot) \bigr](\xi) \, \ud s \biggr\|_{L^2_\xi} \\ 
        &\lesssim \int_0^t s \cdot |\dot{q}(s)-\ulell| \cdot \Bigl\| \jxi \xi \cdot \whatcalF\bigl[ v_1(s,\cdot) \overline{v_2(s,\cdot)} \tilde{v}_3(s,\cdot) \bigr](\xi) \Bigr\|_{L^2_\xi} \, \ud s \\ 
        &\lesssim \int_0^t s \cdot |\dot{q}(s)-\ulell| \cdot \bigl\| v_1(s,y) \overline{v_2(s,y)} \tilde{v}_3(s,y) \bigr\|_{H^2_y} \, \ud s \\
        &\lesssim \int_0^t s \cdot |\dot{q}(s)-\ulell| \cdot \|v_1(s)\|_{H^2_y} \|v_2(s)\|_{W^{1,\infty}_y} \|\tilde{v}_3(s)\|_{W^{1,\infty}_y} \, \ud s + \bigl\{\text{similar terms}\bigr\} \\
        &\lesssim \int_0^t s \cdot \varepsilon \js^{-1+\delta} \cdot \varepsilon \js^\delta \cdot \varepsilon^2 \js^{-1} \, \ud s \lesssim \varepsilon^4 \jt^{2\delta}.
    \end{aligned}
\end{equation*}

\noindent {\it Bound for $\jxi \calK_{2,1,3,4,1}^{\delta_0}(t,\xi)$:}
Similarly to the preceding term, we conclude by Plancherel's theorem and \eqref{equ:consequences_qdot_minus_ulell_decay}, \eqref{equ:consequences_aux_KG_disp_decay}, \eqref{equ:consequences_aux_KG_energy_bounds}, the acceptable bound
\begin{equation*}
    \begin{aligned}
        &\bigl\| \jxi \calK_{2,1,3,4,1}^{\delta_0}(t,\xi) \bigr\|_{L^2_\xi} \\
        &= \biggl\| (2\pi)^{\frac52} \frakb(\xi) \int_0^t s \cdot (\dot{q}(s)-\ulell) \cdot e^{-i\xi\theta(s)} e^{-is(\jxi+\ulell\xi)} \jxi \whatcalF\bigl[ (\py v_1)(s,\cdot) \overline{v_2(s,\cdot)} \tilde{v}_3(s,\cdot) \bigr](\xi) \, \ud s \biggr\|_{L^2_\xi} \\
        &\lesssim \int_0^t s \cdot |\dot{q}(s)-\ulell| \cdot \Bigl\| \jxi \whatcalF\bigl[ (\py v_1)(s,\cdot) \overline{v_2(s,\cdot)} \tilde{v}_3(s,\cdot) \bigr](\xi) \Bigr\|_{L^2_\xi} \, \ud s \\
        &\lesssim \int_0^t s \cdot |\dot{q}(s)-\ulell| \cdot \bigl\| (\py v_1)(s,y) \overline{v_2(s,y)} \tilde{v}_3(s,y) \bigr\|_{H^1_y} \, \ud s \\ 
        &\lesssim \int_0^t s \cdot |\dot{q}(s)-\ulell| \cdot \| v_1(s)\|_{H^2_y} \|v_2(s)\|_{W^{1,\infty}_y} \|\tilde{v}_3(s)\|_{W^{1,\infty}_y} \, \ud s \\ 
        &\lesssim \int_0^t s \cdot \varepsilon \js^{-1+\delta} \cdot \varepsilon \js^\delta \cdot \varepsilon^2 \js^{-1} \, \ud s \lesssim \varepsilon^4 \jt^{2\delta}.
    \end{aligned}
\end{equation*}

\noindent {\it Bound for $\jxi \calK_{2,1,3,4,2}^{\delta_0}(t,\xi)$:}
Invoking Plancherel's theorem and the mapping property~\eqref{equ:mapping_property_calFulellsh_L2} along with the decay estimates \eqref{equ:consequences_dispersive_decay_ulPe_radiation}, \eqref{equ:consequences_aux_KG_disp_decay}, we infer the sufficient bound
\begin{equation*}
    \begin{aligned}
        &\bigl\| \jxi \calK_{2,1,3,4,2}^{\delta_0}(t,\xi) \bigr\|_{L^2_\xi} \\
        &= \biggl\| (2\pi)^{\frac52} \frakb(\xi) \int_0^t s \cdot e^{-i\xi\theta(s)} e^{-is(\jxi+\ulell\xi)} \\
        &\qquad \qquad \qquad \times \jxi \whatcalF\biggl[ \whatcalF^{-1}\Bigl[ \jxione^{-1} \frakb_1(\xi_1) \calFulellsh\bigl[ \calQ_\ulell\bigl( \usubeone(s) \bigr) \bigr](\xi_1) \Bigr](\cdot) \, \overline{v_2(s,\cdot)} \, \tilde{v}_3(s,\cdot) \biggr](\xi) \, \ud s \biggr\|_{L^2_\xi} \\ 
        &\lesssim \int_0^t s \cdot \biggl\| \jxi \whatcalF\biggl[ \whatcalF^{-1}\Bigl[ \jxione^{-1} \frakb_1(\xi_1) \calFulellsh\bigl[ \calQ_\ulell\bigl( \usubeone(s) \bigr)\bigr](\xi_1) \Bigr](\cdot) \, \overline{v_2(s,\cdot)} \, \tilde{v}_3(s,\cdot) \biggr](\xi) \biggr\|_{L^2_\xi} \, \ud s \\ 
        &\lesssim \int_0^t s \cdot \Bigl\| \whatcalF^{-1}\Bigl[ \jxione^{-1} \frakb_1(\xi_1) \calFulellsh\bigl[ \calQ_\ulell\bigl( \usubeone(s) \bigr)\bigr](\xi_1) \Bigr](y) \, \overline{v_2(s,y)} \, \tilde{v}_3(s,y) \Bigr\|_{H^1_y} \, \ud s \\ 
        &\lesssim \int_0^t s \cdot \Bigl\| \whatcalF^{-1}\Bigl[ \jxione^{-1} \frakb_1(\xi_1) \calFulellsh\bigl[ \calQ_\ulell\bigl( \usubeone(s) \bigr)\bigr](\xi_1) \Bigr](y) \Bigr\|_{H^1_y} \, \|v_2(s)\|_{W^{1,\infty}_y} \|\tilde{v}_3(s)\|_{W^{1,\infty}_y} \, \ud s \\ 
        &\lesssim \int_0^t s \cdot \bigl\| \calQ_\ulell\bigl( \usubeone(s) \bigr) \bigr\|_{L^2_y} \|v_2(s)\|_{W^{1,\infty}_y} \|\tilde{v}_3(s)\|_{W^{1,\infty}_y} \, \ud s \\ 
        &\lesssim \int_0^t s \cdot \|\alpha(\ulg y)\|_{L^2_y} \|\usubeone(s)\|_{L^{\infty}_y}^2 \|v_2(s)\|_{W^{1,\infty}_y} \|\tilde{v}_3(s)\|_{W^{1,\infty}_y} \, \ud s \\ 
        &\lesssim \int_0^t s \cdot \varepsilon^4 \js^{-2} \, \ud s \lesssim \varepsilon^4 \log(\jt).
    \end{aligned}
\end{equation*}

\noindent {\it Bound for $\jxi \calK_{2,1,3,4,3}^{\delta_0}(t,\xi)$:}
Similarly to the preceding term, we use Plancherel's theorem and the mapping property~\eqref{equ:mapping_property_calFulellsh_L2} together with the decay estimates \eqref{equ:consequences_dispersive_decay_ulPe_radiation}, \eqref{equ:consequences_aux_KG_disp_decay} as well as the uniform-in-time energy bound \eqref{equ:consequences_energy_bound_radiation} to arrive at the sufficient bound
\begin{equation*}
    \begin{aligned}
        &\bigl\| \jxi \calK_{2,1,3,4,3}^{\delta_0}(t,\xi) \bigr\|_{L^2_\xi} \\
        &= \biggl\| (2\pi)^{\frac52} \frakb(\xi) \int_0^t s \cdot e^{-i\xi\theta(s)} e^{-is(\jxi+\ulell\xi)} \\
        &\qquad \qquad \qquad \times \jxi \whatcalF\biggl[ \whatcalF^{-1}\Bigl[ \jxione^{-1} \frakb_1(\xi_1) \calFulellsh\bigl[ {\textstyle \frac16} \usubeone(s)^3 \bigr](\xi_1)\Bigr](\cdot) \, \overline{v_2(s,\cdot)} \, \tilde{v}_3(s,\cdot) \biggr](\xi) \, \ud s \biggr\|_{L^2_\xi} \\ 
        &\lesssim \int_0^t s \cdot \biggl\| \jxi \whatcalF\biggl[ \whatcalF^{-1}\Bigl[ \jxione^{-1} \frakb_1(\xi_1) \calFulellsh\bigl[ {\textstyle \frac16} \usubeone(s)^3 \bigr](\xi_1)\Bigr](\cdot) \, \overline{v_2(s,\cdot)} \, \tilde{v}_3(s,\cdot) \biggr](\xi) \biggr\|_{L^2_\xi} \, \ud s \\ 
        &\lesssim \int_0^t s \cdot \Bigl\| \whatcalF^{-1}\Bigl[ \jxione^{-1} \frakb_1(\xi_1) \calFulellsh\bigl[ {\textstyle \frac16} \usubeone(s)^3 \bigr](\xi_1)\Bigr](y) \, \overline{v_2(s,y)} \, \tilde{v}_3(s,y) \Bigr\|_{H^1_y} \, \ud s \\ 
        &\lesssim \int_0^t s \cdot \Bigl\| \whatcalF^{-1}\Bigl[ \jxione^{-1} \frakb_1(\xi_1) \calFulellsh\bigl[ {\textstyle \frac16} \usubeone(s)^3 \bigr](\xi_1)\Bigr](y) \Bigr\|_{H^1_y} \|v_2(s)\|_{W^{1,\infty}_y} \|\tilde{v}_3(s)\|_{W^{1,\infty}_y} \, \ud s \\
        &\lesssim \int_0^t s \cdot \bigl\|\usubeone(s)^3 \bigr\|_{L^2_y} \|v_2(s)\|_{W^{1,\infty}_y} \|\tilde{v}_3(s)\|_{W^{1,\infty}_y} \, \ud s \\
        &\lesssim \int_0^t s \cdot \|\usubeone(s)\|_{L^2_y} \|\usubeone(s)\|_{L^\infty_y}^2 \|v_2(s)\|_{W^{1,\infty}_y} \|\tilde{v}_3(s)\|_{W^{1,\infty}_y} \, \ud s \\
        &\lesssim \int_0^t s \cdot \varepsilon \cdot \varepsilon^4 \js^{-2} \, \ud s \lesssim \varepsilon^5 \log(\jt).
    \end{aligned}
\end{equation*}

\noindent {\it Bound for $\jxi \calK_{2,1,3,4,4}^{\delta_0}(t,\xi)$:}
Here we use Plancherel's theorem, the decay estimate \eqref{equ:consequences_aux_KG_disp_decay}, and the $L^2_{\xi_1}$-decay estimate \eqref{equ:consequences_wtilcalR_L2xi_bound} for $\widetilde{\calR}(s,\xi_1)$ to infer the sufficient bound
\begin{equation*}
    \begin{aligned}
        &\bigl\| \jxi \calK_{2,1,3,4,4}^{\delta_0}(t,\xi) \bigr\|_{L^2_\xi} \\
        &= \biggl\| (2\pi)^{\frac52} \frakb(\xi) \int_0^t s \cdot e^{-i\xi\theta(s)} e^{-is(\jxi+\ulell\xi)} \\
        &\qquad \qquad \qquad \times \jxi \whatcalF\biggl[ \whatcalF^{-1}\Bigl[ \jxione^{-1} \frakb_1(\xi_1) \widetilde{\calR}(s,\xi_1) \Bigr](\cdot) \, \overline{v_2(s,\cdot)} \, \tilde{v}_3(s,\cdot) \biggr](\xi) \, \ud s \biggr\|_{L^2_\xi} \\ 
        &\lesssim \int_0^t s \cdot \biggl\| \jxi \whatcalF\biggl[ \whatcalF^{-1}\Bigl[ \jxione^{-1} \frakb_1(\xi_1) \widetilde{\calR}(s,\xi_1) \Bigr](\cdot) \, \overline{v_2(s,\cdot)} \, \tilde{v}_3(s,\cdot) \biggr](\xi) \biggr\|_{L^2_\xi} \, \ud s \\ 
        &\lesssim \int_0^t s \cdot \Bigl\| \whatcalF^{-1}\Bigl[ \jxione^{-1} \frakb_1(\xi_1) \widetilde{\calR}(s,\xi_1) \Bigr](y) \, \overline{v_2(s,y)} \, \tilde{v}_3(s,y) \Bigr\|_{H^1_y} \, \ud s \\ 
        &\lesssim \int_0^t s \cdot \Bigl\| \whatcalF^{-1}\Bigl[ \jxione^{-1} \frakb_1(\xi_1) \widetilde{\calR}(s,\xi_1) \Bigr](y) \Bigr\|_{H^1_y} \|v_2(s)\|_{W^{1,\infty}_y} \|v_3(s)\|_{W^{1,\infty}_y} \, \ud s \\
        &\lesssim \int_0^t s \cdot \bigl\| \widetilde{\calR}(s,\xi_1) \bigr\|_{L^2_{\xi_1}} \|v_2(s)\|_{W^{1,\infty}_y} \|v_3(s)\|_{W^{1,\infty}_y} \, \ud s \\
        &\lesssim \int_0^t s \cdot \varepsilon^2 \js^{-\frac32+\delta} \cdot \varepsilon^2 \js^{-1} \, \ud s \lesssim \varepsilon^4.
    \end{aligned}
\end{equation*}

\noindent {\it Bound for $\jxi \calK_{2,\mathrm{comb},2,1}^{\delta_0}(t,\xi)$:}
Using \eqref{equ:consequences_aux_KG_disp_decay} and \eqref{equ:weighted_energies_proof_wj_H1y_growth}, we find that this mild term is uniformly bounded in time
\begin{equation*}
    \begin{aligned}
        &\bigl\| \jxi \calK_{2,\mathrm{comb},2,1}^{\delta_0}(t,\xi) \bigr\|_{L^2_\xi} 
        = \biggl\| \jxi \biggl[ e^{-i\xi\theta(s)} (2\pi)^{\frac52} i \frakb(\xi) e^{-is(\jxi+\ulell\xi)} \whatcalF\bigl[ v_1(s,\cdot) \overline{v_2(s,\cdot)} w_3(s,\cdot) \bigr](\xi) \biggr]_{s=0}^{s=t} \biggr\|_{L^2_\xi} \\ 
        &\lesssim \sup_{0 \leq s \leq t} \, \Bigl\| \jxi \whatcalF\bigl[ v_1(s,\cdot) \overline{v_2(s,\cdot)} w_3(s,\cdot) \bigr](\xi) \Bigr\|_{L^2_\xi} 
        \lesssim \sup_{0 \leq s \leq t} \, \bigl\| v_1(s,y) \overline{v_2(s,y)} w_3(s,y) \bigr\|_{H^1_y} \\ 
        &\lesssim \sup_{0 \leq s \leq t} \, \|v_1(s)\|_{W^{1,\infty}_y} \|v_2(s)\|_{W^{1,\infty}_y} \|w_3(s)\|_{H^1_y} \lesssim \sup_{0 \leq s \leq t} \, \varepsilon^2 \js^{-1} \cdot \varepsilon \js^\delta \lesssim \varepsilon^3.
    \end{aligned}
\end{equation*}

\noindent {\it Bound for $\jxi \calK_{2,\mathrm{comb},2,2}^{\delta_0}(t,\xi)$:}
By Plancherel's theorem, Sobolev embedding, and \eqref{equ:consequences_qdot_minus_ulell_decay}, \eqref{equ:consequences_aux_KG_disp_decay}, \eqref{equ:consequences_aux_KG_energy_bounds}, \eqref{equ:weighted_energies_proof_wj_H1y_growth}, \eqref{equ:weighted_energies_proof_wj_H3y_growth} we arrive at the sufficient bound 
\begin{equation*}
    \begin{aligned}
        &\bigl\| \jxi \calK_{2,\mathrm{comb},2,2}^{\delta_0}(t,\xi) \bigr\|_{L^2_\xi} \\
        &= \biggl\| \jxi \int_0^t (\dot{q}(s)-\ulell) e^{-i\xi\theta(s)} (2\pi)^{\frac52} \frakb(\xi) e^{-is(\jxi+\ulell\xi)} \cdot \xi \cdot \whatcalF\bigl[ v_1(s,\cdot) \overline{v_2(s,\cdot)} w_3(s,\cdot) \bigr](\xi) \, \ud s \biggr\|_{L^2_\xi} \\
        &\lesssim \int_0^t |\dot{q}(s)-\ulell| \Bigl\| \jxi^2 \whatcalF\bigl[ v_1(s,\cdot) \overline{v_2(s,\cdot)} w_3(s,\cdot) \bigr](\xi) \Bigr\|_{L^2_\xi} \, \ud s 
        \lesssim \int_0^t |\dot{q}(s)-\ulell| \bigl\| v_1(s,y) \overline{v_2(s,y)} w_3(s,y) \bigr\|_{H^2_y} \, \ud s \\ 
        &\lesssim \int_0^t |\dot{q}(s)-\ulell| \Bigl( \|v_1(s)\|_{H^2_y} \|v_2(s)\|_{L^\infty_y} \| w_3(s,y)\|_{H^1_y} +  \|v_1(s)\|_{L^\infty_y} \|v_2(s)\|_{H^2_y} \| w_3(s,y)\|_{H^1_y} \\ 
        &\qquad \qquad \qquad \qquad \qquad \qquad \qquad \qquad \qquad \qquad \qquad \qquad + \|v_1(s)\|_{W^{1,\infty}_y} \|v_2(s)\|_{W^{1,\infty}_y} \| w_3(s,y)\|_{H^2_y} \Bigr) \\ 
        &\lesssim \int_0^t \varepsilon \js^{-1+\delta} \Bigl( \varepsilon \js^\delta \cdot \varepsilon \js^{-\frac12} \cdot \varepsilon \js^\delta + \varepsilon \js^{-\frac12} \cdot \varepsilon \js^\delta \cdot \varepsilon \js^\delta + \varepsilon^2 \js^{-1} \cdot \varepsilon \js^{2\delta} \Bigr) \, \ud s 
        \lesssim \varepsilon^4.
    \end{aligned}
\end{equation*}

\noindent {\it Bound for $\jxi \calK_{2,\mathrm{comb},2,3,1}^{\delta_0}(t,\xi)$:}
Similarly to the preceding term, we use Plancherel's theorem and \eqref{equ:consequences_qdot_minus_ulell_decay}, \eqref{equ:consequences_aux_KG_disp_decay}, \eqref{equ:consequences_aux_KG_energy_bounds}, \eqref{equ:weighted_energies_proof_wj_H1y_growth}, to deduce the sufficient bound
\begin{equation*}
    \begin{aligned}
        &\bigl\| \jxi \calK_{2,\mathrm{comb},2,3,1}^{\delta_0}(t,\xi) \bigr\|_{L^2_\xi} \\
        &= \biggl\| \jxi i (2\pi)^{\frac52} \frakb(\xi) \int_0^t e^{-i\xi\theta(s)} (\dot{q}(s)-\ulell)  e^{-is(\jxi+\ulell\xi)} \whatcalF\bigl[ (\py v_1)(s,\cdot) \overline{v_2(s,\cdot)} w_3(s,\cdot) \bigr](\xi) \, \ud s \biggr\|_{L^2_\xi} \\ 
        &\lesssim \int_0^t \Bigl\| \jxi \whatcalF\bigl[ (\py v_1)(s,\cdot) \overline{v_2(s,\cdot)} w_3(s,\cdot) \bigr](\xi) \Bigr\|_{L^2_\xi} \, \ud s \\
        &\lesssim \int_0^t |\dot{q}(s)-\ulell| \bigl\| (\py v_1)(s,y) \overline{v_2(s,y)} w_3(s,y) \bigr\|_{H^1_y} \, \ud s \\ 
        &\lesssim \int_0^t |\dot{q}(s)-\ulell| \Bigl( \|v_1(s)\|_{H^2_y} \|v_2(s)\|_{L^\infty_y} \|w_3(s,y)\|_{H^1_y} + \|v_1(s)\|_{W^{1,\infty}_y} \|v_2(s)\|_{W^{1,\infty}_y} \|w_3(s,y)\|_{H^1_y} \Bigr) \\ 
        &\lesssim \int_0^t \varepsilon \js^{-1+\delta} \Bigl( \varepsilon \js^\delta \cdot \varepsilon \js^{-\frac12} \cdot \varepsilon \js^\delta + \varepsilon^2 \js^{-1} \cdot \varepsilon \js^{\delta} \Bigr) \, \ud s 
        \lesssim \varepsilon^4.
    \end{aligned}
\end{equation*}

\noindent {\it Bound for $\jxi \calK_{2,\mathrm{comb},2,3,2}^{\delta_0}(t,\xi)$:}
Here we use Plancherel's theorem, the mapping property \eqref{equ:mapping_property_calFulellsh_jxi}, Sobolev embedding combined with \eqref{equ:consequences_dispersive_decay_ulPe_radiation}, \eqref{equ:consequences_aux_KG_disp_decay}, \eqref{equ:weighted_energies_proof_wj_H1y_growth} to arrive at the sufficient bound
\begin{equation*}
    \begin{aligned}
        &\bigl\| \jxi \calK_{2,\mathrm{comb},2,3,2}^{\delta_0}(t,\xi) \bigr\|_{L^2_\xi} \\
        &= \biggl\| \jxi i (2\pi)^{\frac52} \frakb(\xi) \int_0^t e^{-i\xi\theta(s)} e^{-is(\jxi+\ulell\xi)} \\ 
        &\qquad \qquad \qquad \qquad \times \whatcalF\biggl[ \whatcalF^{-1}\Bigl[ \jxione^{-1} \frakb_1(\xi_1) \calFulellsh\bigl[ \calQ_\ulell\bigl( \usubeone(s) \bigr) \bigr](\xi_1) \Bigr](\cdot) \overline{v_2(s,\cdot)} w_3(s,\cdot) \biggr](\xi) \, \ud s \biggr\|_{L^2_\xi} \\ 
        &\lesssim \int_0^t \, \biggl\| \jxi \whatcalF\biggl[ \whatcalF^{-1}\Bigl[ \jxione^{-1} \frakb_1(\xi_1) \calFulellsh\bigl[ \calQ_\ulell\bigl( \usubeone(s) \bigr) \bigr](\xi_1) \Bigr](\cdot) \overline{v_2(s,\cdot)} w_3(s,\cdot) \biggr](\xi) \biggr\|_{L^2_\xi} \, \ud s \\ 
        &\lesssim \int_0^t \, \Bigl\| \whatcalF^{-1}\Bigl[ \jxione^{-1} \frakb_1(\xi_1) \calFulellsh\bigl[ \calQ_\ulell\bigl( \usubeone(s) \bigr) \bigr](\xi_1) \Bigr](y) \overline{v_2(s,y)} w_3(s,y) \Bigr\|_{H^1_y} \, \ud s \\ 
        &\lesssim \int_0^t \Bigl\| \whatcalF^{-1}\Bigl[ \jxione^{-1} \frakb_1(\xi_1) \calFulellsh\bigl[ \calQ_\ulell\bigl( \usubeone(s) \bigr) \bigr](\xi_1) \Bigr](y) \Bigr\|_{W^{1,\infty}_y} \|v_2(s)\|_{W^{1,\infty}_y} \|w_3(s)\|_{H^1_y} \, \ud s \\ 
        &\lesssim \int_0^t \bigl\| \calQ_\ulell\bigl( \usubeone(s) \bigr) \bigr\|_{H^1_y} \|v_2(s)\|_{W^{1,\infty}_y} \|w_3(s)\|_{H^1_y} \, \ud s  
        \lesssim \int_0^t \|\usubeone(s)\|_{W^{1,\infty}_y}^2 \|v_2(s)\|_{W^{1,\infty}_y} \|w_3(s)\|_{H^1_y} \, \ud s \\
        &\lesssim \int_0^t \varepsilon^3 \js^{-\frac32} \cdot \varepsilon \js^\delta \, \ud s \lesssim \varepsilon^4.
    \end{aligned}
\end{equation*}

\noindent {\it Bound for $\jxi \calK_{2,\mathrm{comb},2,3,3}^{\delta_0}(t,\xi)$:}
Similarly to the preceding term, combining Plancherel's theorem, the mapping property \eqref{equ:mapping_property_calFulellsh_jxi}, Sobolev embedding with the bounds \eqref{equ:consequences_dispersive_decay_ulPe_radiation}, \eqref{equ:consequences_energy_bound_radiation}, \eqref{equ:consequences_aux_KG_disp_decay}, \eqref{equ:weighted_energies_proof_wj_H1y_growth} leads to the sufficient estimate
\begin{equation*}
    \begin{aligned}
        &\bigl\| \jxi \calK_{2,\mathrm{comb},2,3,3}^{\delta_0}(t,\xi) \bigr\|_{L^2_\xi} \\
        &= \biggl\| \jxi (2\pi)^{\frac52} \frakb(\xi) \int_0^t e^{-i\xi\theta(s)} e^{-is(\jxi+\ulell\xi)} \\ 
        &\qquad \qquad \qquad \qquad \times \whatcalF\biggl[ \whatcalF^{-1}\Bigl[ \jxione^{-1} \frakb_1(\xi_1) \calFulellsh\bigl[ {\textstyle \frac16} \usubeone(s)^3 \bigr](\xi_1) \Bigr](\cdot) \overline{v_2(s,\cdot)} w_3(s,\cdot) \biggr](\xi) \, \ud s \biggr\|_{L^2_\xi} \\
        &\lesssim \int_0^t \Bigl\| \whatcalF^{-1}\Bigl[ \jxione^{-1} \frakb_1(\xi_1) \calFulellsh\bigl[ {\textstyle \frac16} \usubeone(s)^3 \bigr](\xi_1) \Bigr](y) \overline{v_2(s,y)} w_3(s,y) \Bigr\|_{H^1_y} \, \ud s \\ 
        &\lesssim \int_0^t \Bigl\| \whatcalF^{-1}\Bigl[ \jxione^{-1} \frakb_1(\xi_1) \calFulellsh\bigl[ {\textstyle \frac16} \usubeone(s)^3 \bigr](\xi_1) \Bigr](y) \Bigr\|_{W^{1,\infty}_y} \|v_2(s)\|_{W^{1,\infty}_y} \|w_3(s)\|_{H^1_y} \, \ud s \\ 
        &\lesssim \int_0^t \bigl\| \usubeone(s)^3 \bigr\|_{H^1_y} \|v_2(s)\|_{W^{1,\infty}_y} \|w_3(s)\|_{H^1_y} \, \ud s \\ 
        &\lesssim \int_0^t \|\usubeone(s)\|_{H^1_y} \|\usubeone(s)\|_{L^\infty_y}^2 \|v_2(s)\|_{W^{1,\infty}_y} \|w_3(s)\|_{H^1_y} \, \ud s \\ 
        &\lesssim \int_0^t \varepsilon \cdot \varepsilon^2 \js^{-1} \cdot \varepsilon \js^{-\frac12} \cdot \varepsilon \js^{\delta} \, \ud s \lesssim \varepsilon^5.
    \end{aligned}
\end{equation*}

\noindent {\it Bound for $\jxi \calK_{2,\mathrm{comb},2,3,4}^{\delta_0}(t,\xi)$:}
Here we use Plancherel's theorem, Sobolev embedding, H\"older's inequality in the frequency variable $\xi_1$, and the bounds \eqref{equ:consequences_wtilcalR_L2xi_bound}, \eqref{equ:consequences_aux_KG_disp_decay}, \eqref{equ:weighted_energies_proof_wj_H1y_growth} to get the sufficient bound
\begin{equation*}
    \begin{aligned}
        &\bigl\| \jxi \calK_{2,\mathrm{comb},2,3,4}^{\delta_0}(t,\xi) \bigr\|_{L^2_\xi} \\ 
        &= \biggl\| \jxi (2\pi)^{\frac52} \frakb(\xi) \int_0^t e^{-i\xi\theta(s)} e^{-is(\jxi+\ulell\xi)} \\ 
        &\qquad \qquad \qquad \qquad \times \whatcalF\biggl[ \whatcalF^{-1}\Bigl[ \jxione^{-1} \frakb_1(\xi_1) \widetilde{\calR}(s,\xi_1) \Bigr](\cdot) \overline{v_2(s,\cdot)} w_3(s,\cdot) \biggr](\xi) \, \ud s \biggr\|_{L^2_\xi} \\ 
        &\lesssim \int_0^t \Bigl\| \whatcalF^{-1}\Bigl[ \jxione^{-1} \frakb_1(\xi_1) \widetilde{\calR}(s,\xi_1) \Bigr](y) \overline{v_2(s,y)} w_3(s,y) \Bigr\|_{H^1_y} \, \ud s \\ 
        &\lesssim \int_0^t \Bigl\| \whatcalF^{-1}\Bigl[ \jxione^{-1} \frakb_1(\xi_1) \widetilde{\calR}(s,\xi_1) \Bigr](y) \Bigr\|_{W^{1,\infty}_y} \|v_2(s)\|_{W^{1,\infty}_y} \|w_3(s)\|_{H^1_y} \, \ud s \\ 
        &\lesssim \int_0^t \bigl\| \jxione \widetilde{\calR}(s,\xi_1) \bigr\|_{L^1_{\xi_1}} \|v_2(s)\|_{W^{1,\infty}_y} \|w_3(s)\|_{H^1_y} \, \ud s \\ 
        &\lesssim \int_0^t \varepsilon^2 \js^{-\frac32+\delta} \cdot \varepsilon \js^{-\frac12} \cdot \varepsilon \js^\delta \, \ud s \lesssim \varepsilon^4.
    \end{aligned}
\end{equation*}

\noindent {\it Bound for $\jxi \calK_{2,\mathrm{comb},2,5,1}^{\delta_0}(t,\xi)$:}
For this term we obtain a sufficient bound from \eqref{equ:consequences_qdot_minus_ulell_decay} and \eqref{equ:consequences_aux_KG_disp_decay},
\begin{equation*}
    \begin{aligned}
        &\bigl\| \jxi \calK_{2,\mathrm{comb},2,5,1}^{\delta_0}(t,\xi) \bigr\|_{L^2_\xi} \\
        &= \biggl\| \jxi i (2\pi)^{\frac52} \frakb(\xi) \int_0^t (\dot{q}(s)-\ulell) e^{-i\xi\theta(s)} e^{-is(\jxi+\ulell\xi)} \\
        &\qquad \times \whatcalF\biggl[ v_1(s,\cdot) \overline{v_2(s,\cdot)} \whatcalF^{-1}\Bigl[ e^{is(\jxithree+\ulell\xi_3)} \pxithree \Bigl( \jxithree^{-1} \frakb_3(\xi_3) (i\xi_3) \gulellsh(s,\xi_3) \Bigr) \Bigr](\cdot) \biggr](\xi) \, \ud s \biggr\|_{L^2_\xi} \\ 
        &\lesssim \int_0^t |\dot{q}(s)-\ulell|  \Bigl\| v_1(s,y) \overline{v_2(s,y)} \whatcalF^{-1}\Bigl[ e^{is(\jxithree+\ulell\xi_3)} \pxithree \Bigl( \jxithree^{-1} \frakb_3(\xi_3) (i\xi_3) \gulellsh(s,\xi_3) \Bigr) \Bigr](y) \Bigr\|_{H^1_y} \, \ud s \\ 
        &\lesssim \int_0^t |\dot{q}(s)-\ulell| \|v_1(s)\|_{W^{1,\infty}_y} \|v_2(s)\|_{W^{1,\infty}_y} \Bigl\| \whatcalF^{-1}\Bigl[ e^{is(\jxithree+\ulell\xi_3)} \pxithree \Bigl( \jxithree^{-1} \frakb_3(\xi_3) (i\xi_3) \gulellsh(s,\xi_3) \Bigr) \Bigr](y) \Bigr\|_{H^1_y} \, \ud s \\ 
        &\lesssim \int_0^t \varepsilon \js^{-1+\delta} \cdot \varepsilon^2 \js^{-1} \cdot \varepsilon \js^{2\delta} \, \ud s \lesssim \varepsilon^4,
    \end{aligned}
\end{equation*}
where we used that by Plancherel's theorem and the bootstrap assumptions \eqref{equ:prop_profile_bounds_assumption2}, it holds that 
\begin{equation*}
    \begin{aligned}
        &\Bigl\| \whatcalF^{-1}\Bigl[ e^{is(\jxithree+\ulell\xi_3)} \pxithree \Bigl( \jxithree^{-1} \frakb_3(\xi_3) (i\xi_3) \gulellsh(s,\xi_3) \Bigr) \Bigr](y) \Bigr\|_{H^1_y} \\ 
        &\quad \lesssim \bigl\| \jxithree \gulellsh(s,\xi_3) \bigr\|_{L^2_{\xi_3}} + \bigl\| \jxithree \pxithree \gulellsh(s,\xi_3) \bigr\|_{L^2_{\xi_3}} \lesssim \varepsilon \js^{2\delta}.
    \end{aligned}
\end{equation*}

\noindent {\it Bound for $\jxi \calK_{2,\mathrm{comb},2,5,2}^{\delta_0}(t,\xi)$:}
Similarly, using the decay estimate \eqref{equ:consequences_aux_KG_disp_decay} we infer the sufficient bound
\begin{equation*}
    \begin{aligned}
        &\bigl\| \jxi \calK_{2,\mathrm{comb},2,5,2}^{\delta_0}(t,\xi) \bigr\|_{L^2_\xi} \\
        &\lesssim \biggl\| \jxi (2\pi)^{\frac52} \frakb(\xi) \int_0^t e^{-i\xi\theta(s)} e^{-is(\jxi+\ulell\xi)} \\
        &\qquad \times \whatcalF\biggl[ v_1(s,\cdot) \overline{v_2(s,\cdot)} \whatcalF^{-1}\Bigl[ \bigl( \pxithree -i s (\jxithree^{-1} \xi_3 + \ulell) \bigr) \jxithree^{-1} \frakb_3(\xi_3) \calFulellsh\bigl[ \calQ_\ulell\bigl( \usubeone(s) \bigr)\bigr](\xi_3) \Bigr) \Bigr](\cdot) \biggr](\xi) \, \ud s \biggr\|_{L^2_\xi} \\
        &\lesssim \int_0^t \Bigl\| v_1(s,y) \overline{v_2(s,y)} \whatcalF^{-1}\Bigl[ \bigl( \pxithree -i s (\jxithree^{-1} \xi_3 + \ulell) \bigr) \jxithree^{-1} \frakb_3(\xi_3) \calFulellsh\bigl[ \calQ_\ulell\bigl( \usubeone(s) \bigr)\bigr](\xi_3) \Bigr) \Bigr](y) \Bigr\|_{H^1_y} \, \ud s \\ 
        &\lesssim \int_0^t \|v_1(s)\|_{W^{1,\infty}_y} \|v_2(s)\|_{W^{1,\infty}_y} \Bigl\| \whatcalF^{-1}\Bigl[ \bigl( \pxithree -i s (\jxithree^{-1} \xi_3 + \ulell) \bigr) \jxithree^{-1} \frakb_3(\xi_3) \calFulellsh\bigl[ \calQ_\ulell\bigl( \usubeone(s) \bigr)\bigr](\xi_3) \Bigr) \Bigr](y) \Bigr\|_{H^1_y} \, \ud s \\ 
        &\lesssim \int_0^t \varepsilon^2 \js^{-1} \cdot \varepsilon^2 \, \ud s \lesssim \varepsilon^4 \log(\jt),
    \end{aligned}
\end{equation*}
where we used that by Plancherel's theorem, the mapping properties \eqref{equ:mapping_property_calFulellsh_L2}, \eqref{equ:mapping_property_calFulellsh_pxi}, and \eqref{equ:consequences_dispersive_decay_ulPe_radiation} we have
\begin{equation*}
    \begin{aligned}
        &\Bigl\| \whatcalF^{-1}\Bigl[ \bigl( \pxithree -i s (\jxithree^{-1} \xi_3 + \ulell) \bigr) \jxithree^{-1} \frakb_3(\xi_3) \calFulellsh\bigl[ \calQ_\ulell\bigl( \usubeone(s) \bigr)\bigr](\xi_3) \Bigr) \Bigr](y) \Bigr\|_{H^1_y} \\ 
        &\lesssim \bigl\| \jy \calQ_\ulell\bigl( \usubeone(s) \bigr) \bigr\|_{L^2_y} + \js \bigl\| \calQ_\ulell\bigl( \usubeone(s) \bigr) \bigr\|_{L^2_y} \\ 
        &\lesssim \bigl\|\jy \alpha(\ulg y)\bigr\|_{L^2_y} \|\usubeone(s)\|_{L^\infty_y}^2 + \js \bigl\|\alpha(\ulg y)\bigr\|_{L^2_y} \|\usubeone(s)\|_{L^\infty_y}^2 \lesssim \varepsilon^2.
    \end{aligned}
\end{equation*}

\noindent {\it Bound for $\jxi \calK_{2,\mathrm{comb},2,5,3}^{\delta_0}(t,\xi)$:}
In the same manner, we obtain a sufficient bound for this term using the decay estimate~\eqref{equ:consequences_aux_KG_disp_decay}
\begin{equation*}
    \begin{aligned}
        &\bigl\| \jxi \calK_{2,\mathrm{comb},2,5,3}^{\delta_0}(t,\xi) \bigr\|_{L^2_\xi} \\
        &= \biggl\| \jxi (2\pi)^{\frac52} \frakb(\xi) \int_0^t e^{-i\xi\theta(s)} e^{-is(\jxi+\ulell\xi)} \\
        &\qquad \times \whatcalF\biggl[ v_1(s,\cdot) \overline{v_2(s,\cdot)} \whatcalF^{-1}\Bigl[ \bigl( \pxithree -i s (\jxithree^{-1} \xi_3 + \ulell) \bigr) \jxithree^{-1} \frakb_3(\xi_3) \calFulellsh\bigl[ {\textstyle \frac16} \usubeone(s)^3 \bigr](\xi_3) \Bigr](\cdot) \biggr](\xi) \, \ud s \biggr\|_{L^2_\xi} \\ 
        &\lesssim \int_0^t \Bigl\| v_1(s,y) \overline{v_2(s,y)} \whatcalF^{-1}\Bigl[ \bigl( \pxithree -i s (\jxithree^{-1} \xi_3 + \ulell) \bigr) \jxithree^{-1} \frakb_3(\xi_3) \calFulellsh\bigl[ {\textstyle \frac16} \usubeone(s)^3 \bigr](\xi_3) \Bigr](y) \Bigr\|_{H^1_y} \, \ud s \\ 
        &\lesssim \int_0^t \|v_1(s)\|_{W^{1,\infty}_y} \|v_2(s)\|_{W^{1,\infty}_y} \Bigl\| \whatcalF^{-1}\Bigl[ \bigl( \pxithree -i s (\jxithree^{-1} \xi_3 + \ulell) \bigr) \jxithree^{-1} \frakb_3(\xi_3) \calFulellsh\bigl[ {\textstyle \frac16} \usubeone(s)^3 \bigr](\xi_3) \Bigr](y) \Bigr\|_{H^1_y} \, \ud s \\ 
        &\lesssim \int_0^t \varepsilon^2 \js^{-1} \cdot \varepsilon^3 \js^\delta \, \ud s \lesssim \varepsilon^5 \jt^\delta,
    \end{aligned}
\end{equation*}
where we used that by Plancherel's theorem, the mapping properties \eqref{equ:mapping_property_calFulellsh_L2}, \eqref{equ:mapping_property_calFulellsh_pxi}, and \eqref{equ:consequences_dispersive_decay_ulPe_radiation}, \eqref{equ:consequences_energy_bound_radiation} we have
\begin{equation*}
    \begin{aligned}
        &\Bigl\| \whatcalF^{-1}\Bigl[ \bigl( \pxithree -i s (\jxithree^{-1} \xi_3 + \ulell) \bigr) \jxithree^{-1} \frakb_3(\xi_3) \calFulellsh\bigl[ {\textstyle \frac16} \usubeone(s)^3 \bigr](\xi_3) \Bigr](y) \Bigr\|_{H^1_y} \\
        &\lesssim \bigl\| \jy \usubeone(s)^3 \bigr\|_{L^2_y} + \js \bigl\| \usubeone(s)^3 \bigr\|_{L^2_y}  
        \lesssim \|\jy \usubeone(s)\|_{L^2_y} \|\usubeone(s)\|_{L^\infty_y}^2 + \js \|\usubeone(s)\|_{L^2_y} \|\usubeone(s)\|_{L^\infty_y}^2 \\ 
        &\lesssim \varepsilon \js^{1+\delta} \cdot \varepsilon^2 \js^{-1} + \js \cdot \varepsilon \cdot \varepsilon^2 \js^{-1} \lesssim \varepsilon^3 \js^\delta.
    \end{aligned}
\end{equation*}

\noindent {\it Bound for $\jxi \calK_{2,\mathrm{comb},2,5,4}^{\delta_0}(t,\xi)$:}
Finally, we also obtain a sufficient bound for this term from the decay estimate \eqref{equ:consequences_aux_KG_disp_decay},
\begin{equation*}
    \begin{aligned}
        &\bigl\| \jxi \calK_{2,\mathrm{comb},2,5,4}^{\delta_0}(t,\xi)\bigr\|_{L^2_\xi} \\
        &= \biggl\| \jxi (2\pi)^{\frac52} \frakb(\xi) \int_0^t e^{-i\xi\theta(s)} e^{-is(\jxi+\ulell\xi)} \\
        &\qquad \times \whatcalF\biggl[ v_1(s,\cdot) \overline{v_2(s,\cdot)} \whatcalF^{-1}\Bigl[ \bigl( \pxithree -i s (\jxithree^{-1} \xi_3 + \ulell) \bigr) \jxithree^{-1} \frakb_3(\xi_3) \widetilde{\calR}(s,\xi_3) \Bigr) \Bigr](\cdot) \biggr](\xi) \, \ud s \biggr\|_{L^2_\xi} \\
        &\lesssim \int_0^t \|v_1(s)\|_{W^{1,\infty}_y} \|v_2(s)\|_{W^{1,\infty}_y} \Bigl\|\whatcalF^{-1}\Bigl[ \bigl( \pxithree -i s (\jxithree^{-1} \xi_3 + \ulell) \bigr) \jxithree^{-1} \frakb_3(\xi_3) \widetilde{\calR}(s,\xi_3) \Bigr) \Bigr](y) \Bigr\|_{H^1_y} \, \ud s \\ 
        &\lesssim \int_0^t \varepsilon^2 \js^{-1} \cdot \varepsilon^2 \js^{-\frac12+\delta} \, \ud s \lesssim \varepsilon^4,
    \end{aligned}
\end{equation*}
where we used that by Plancherel's theorem, and \eqref{equ:consequences_wtilcalR_L2xi_bound}, \eqref{equ:consequences_wtilcalR_pxi_L2xi_bound} we have
\begin{equation*}
    \begin{aligned}
        &\Bigl\|\whatcalF^{-1}\Bigl[ \bigl( \pxithree -i s (\jxithree^{-1} \xi_3 + \ulell) \bigr) \jxithree^{-1} \frakb_3(\xi_3) \widetilde{\calR}(s,\xi_3) \Bigr) \Bigr](y) \Bigr\|_{H^1_y} \\
        &\lesssim \bigl\| \pxithree \widetilde{\calR}(s,\xi_3) \bigr\|_{L^2_{\xi_3}} + \js \bigl\| \widetilde{\calR}(s,\xi_3) \bigr\|_{L^2_{\xi_3}}  
        \lesssim \varepsilon^2 \js^{-1+\delta} + \js \cdot \varepsilon^2 \js^{-\frac32+\delta} \lesssim \varepsilon^2 \js^{-\frac12+\delta}.
    \end{aligned}
\end{equation*}

\medskip
\noindent \underline{Case 2: Cubic interactions with a Hilbert-type kernel.}
Next, we consider the cubic interactions with a Hilbert-type kernel. In view of the fine structure of the cubic spectral distributions analyzed in Subsection~\ref{subsec:cubic_spectral_distributions}, the term $\calI_2^{\pvdots}(s,\xi)$ is a linear combination of terms of the schematic form
\begin{equation*}
    \begin{aligned}
        \calI_{2, \mathrm{schem}}^{\pvdots}(s,\xi) := \frakb(\xi) \iiint e^{i s \Omega_{2,\ulell}(\xi,\xi_1,\xi_2,\xi_4)} h_1(s,\xi_1) \overline{h_2(s,\xi_2)} h_3(s,\xi_3) \, \pvdots \frac{\varphi(\xi_4)}{\xi_4} \, \ud \xi_1 \, \ud \xi_2 \, \ud \xi_4,
    \end{aligned}
\end{equation*}
with the phase
\begin{equation*}
    \Omega_{2,\ulell}(\xi,\xi_1,\xi_2) := -\jxi + \jxione - \jxitwo + \jxithree + \ulell \xi_4, \quad \quad \xi_3 := \xi - \xi_1 + \xi_2 + \xi_4,
\end{equation*}
with a Hilbert-type kernel involving the even Schwartz function
\begin{equation*}
    \varphi(\xi_4) := \xi_4 \, \cosech\Bigl( \frac{\pi}{2\ulg} \xi_4 \Bigr),   
\end{equation*}
and for some coefficients $\frakb, \frakb_1, \frakb_2, \frakb_3 \in W^{1,\infty}(\bbR)$ such that the inputs are given by
\begin{equation*}
    h_j(s,\xi_j) := \jap{\xi_j}^{-1} \frakb_j(\xi_j) \gulellsh(s,\xi_j), \quad 1 \leq j \leq 3.
\end{equation*}
Moreover, by the structure of $\frakm_{\ulell,+-+}^{\pvdots}$ uncovered in Subsection~\ref{subsec:cubic_spectral_distributions}, at least one of the coefficients must have an improved low frequency behavior. Specifically, it must either hold that
\begin{equation*}
    \frakb(\xi) = \frac{\ulg(\xi+\ulell\jxi)}{|\ulg(\xi+\ulell\jxi)|+i},
\end{equation*}
or at least one of the coefficients $\frakb_j(\xi_j)$, $1 \leq j \leq 3$, is given by
\begin{equation*}
    \frac{\ulg(\xi_j+\ulell\jap{\xi_j})}{|\ulg(\xi_j+\ulell\jap{\xi_j})|\pm i}.
\end{equation*}
We now proceed along the lines of the treatment of the cubic interactions with a Dirac kernel in the previous case, and compute 
\begin{equation*}
    \begin{aligned}
        \jxi \pxi \calI_{2, \mathrm{schem}}^{\pvdots}(s,\xi) &= \frakb(\xi) \iiint i s \cdot \jxi (\pxi \Omega_{2,\ulell}) e^{i s \Omega_{2,\ulell}} h_1(s,\xi_1) \overline{h_2(s,\xi_2)} h_3(s,\xi_3) \, \pvdots \frac{\varphi(\xi_4)}{\xi_4} \, \ud \xi_1 \, \ud \xi_2 \, \ud \xi_4 \\ 
        &\quad + \jxi \frakb(\xi) \iiint e^{i s \Omega_{2,\ulell}} h_1(s,\xi_1) \overline{h_2(s,\xi_2)} (\pxithree h_3)(s,\xi_3) \, \pvdots \frac{\varphi(\xi_4)}{\xi_4} \, \ud \xi_1 \, \ud \xi_2 \, \ud \xi_4 \\
        &\quad + \jxi (\pxi \frakb)(\xi) \iiint e^{i s \Omega_{2,\ulell}} h_1(s,\xi_1) \overline{h_2(s,\xi_2)} (\pxithree h_3)(s,\xi_3) \, \pvdots \frac{\varphi(\xi_4)}{\xi_4} \, \ud \xi_1 \, \ud \xi_2 \, \ud \xi_4 \\
        &=: \calJ_{2,1}^{\pvdots}(s,\xi) + \calJ_{2,2}^{\pvdots}(s,\xi) + \calJ_{2,3}^{\pvdots}(s,\xi).
    \end{aligned}
\end{equation*}
In the first term $\calJ_{2,1}^{\pvdots}(s,\xi)$ we insert the following identity
\begin{equation} \label{equ:weighted_energy_proof_phase_derivative_identity_Hilbert}
    \jxi \pxi \Omega_{2,\ulell} = - \jxione \pxione \Omega_{2,\ulell} - \jxitwo \pxitwo \Omega_{2,\ulell} - \Omega_{2,\ulell} \frac{\xi_3}{\jxithree} + \xi_4 \frac{\jxithree + \ulell \xi_3}{\jxithree},
\end{equation}
which can be verified by direct computation. 
In comparison with the corresponding identity \eqref{equ:weighted_energy_proof_phase_derivative_identity_Dirac} for the cubic interactions with a Dirac kernel, we observe the additional fourth term on the right-hand side of \eqref{equ:weighted_energy_proof_phase_derivative_identity_Hilbert}. 
Integrating by parts in the frequency variables $\xi_1$ and $\xi_2$ then gives 
\begin{equation} \label{equ:weighted_energy_proof_calJ21pv_phase_identity_inserted}
    \begin{aligned}
        &\calJ_{2,1}^{\pvdots}(s,\xi) \\ 
        &= \frakb(\xi) \iiint e^{i s \Omega_{2,\ulell}} \pxione \Bigl( \jxione h_1(s,\xi_1)  \overline{h_2(s,\xi_2)} h_3(s,\xi_3) \Bigr) \, \pvdots \frac{\varphi(\xi_4)}{\xi_4} \, \ud \xi_1 \, \ud \xi_2 \, \ud \xi_4 \\ 
        &\quad + \frakb(\xi) \iiint e^{i s \Omega_{2,\ulell}} \pxitwo \Bigl( h_1(s,\xi_1)  \overline{\jxitwo h_2(s,\xi_2)} h_3(s,\xi_3) \Bigr) \, \pvdots \frac{\varphi(\xi_4)}{\xi_4} \, \ud \xi_1 \, \ud \xi_2 \, \ud \xi_4 \\ 
        &\quad - \frakb(\xi) \iiint i s \, \Omega_{2,\ulell} \, e^{i s \Omega_{2,\ulell}} h_1(s,\xi_1) \overline{h_2(s,\xi_2)} (\jxithree^{-1} \xi_3) h_3(s,\xi_3) \, \pvdots \frac{\varphi(\xi_4)}{\xi_4} \, \ud \xi_1 \, \ud \xi_2 \, \ud \xi_4 \\ 
        &\quad + \frakb(\xi) \iiint i s \, e^{i s \Omega_{2,\ulell}} h_1(s,\xi_1) \overline{h_2(s,\xi_2)} \jxithree^{-1} (\jxithree + \ulell \xi_3) h_3(s,\xi_3) \, \varphi(\xi_4) \, \ud \xi_1 \, \ud \xi_2 \, \ud \xi_4 \\ 
        &=: \calJ_{2,1,1}^{\pvdots}(s,\xi) + \calJ_{2,1,2}^{\pvdots}(s,\xi) + \calJ_{2,1,3}^{\pvdots}(s,\xi) + \calJ_{2,1,4}^{\pvdots}(s,\xi).
    \end{aligned}
\end{equation}
The terms $\calJ_{2,1,1}^{\pvdots}(s,\xi)$ and $\calJ_{2,1,2}^{\pvdots}(s,\xi)$ read in expanded form
\begin{equation*}
    \begin{aligned}
        \calJ_{2,1,1}^{\pvdots}(s,\xi) &= \frakb(\xi) \iiint e^{i s \Omega_{2,\ulell}}  \jxione (\pxione h_1)(s,\xi_1) \overline{h_2(s,\xi_2)} h_3(s,\xi_3) \, \pvdots \frac{\varphi(\xi_4)}{\xi_4} \, \ud \xi_1 \, \ud \xi_2 \, \ud \xi_4 \\ 
        &\quad + \frakb(\xi) \iiint e^{i s \Omega_{2,\ulell}} (\jxione^{-1} \xi_1) h_1(s,\xi_1)  \overline{h_2(s,\xi_2)} h_3(s,\xi_3) \, \pvdots \frac{\varphi(\xi_4)}{\xi_4} \, \ud \xi_1 \, \ud \xi_2 \, \ud \xi_4 \\ 
        &\quad + \frakb(\xi) \iiint e^{i s \Omega_{2,\ulell}} \jxione h_1(s,\xi_1)  \overline{h_2(s,\xi_2)} (-1) (\pxithree h_3)(s,\xi_3) \, \pvdots \frac{\varphi(\xi_4)}{\xi_4} \, \ud \xi_1 \, \ud \xi_2 \, \ud \xi_4 \\ 
        &=: \calJ_{2,1,1,1}^{\pvdots}(s,\xi) + \calJ_{2,1,1,2}^{\pvdots}(s,\xi) + \calJ_{2,1,1,3}^{\pvdots}(s,\xi),
    \end{aligned}
\end{equation*}
and 
\begin{equation*}
    \begin{aligned}
        \calJ_{2,1,2}^{\pvdots}(s,\xi) &= \frakb(\xi) \iiint e^{i s \Omega_{2,\ulell}} h_1(s,\xi_1)  \overline{\jxitwo \pxitwo h_2(s,\xi_2)} h_3(s,\xi_3) \, \pvdots \frac{\varphi(\xi_4)}{\xi_4} \, \ud \xi_1 \, \ud \xi_2 \, \ud \xi_4 \\
        &\quad + \frakb(\xi) \iiint e^{i s \Omega_{2,\ulell}} h_1(s,\xi_1) \overline{(\jxitwo^{-1} \xi_2) h_2(s,\xi_2)} h_3(s,\xi_3) \, \pvdots \frac{\varphi(\xi_4)}{\xi_4} \, \ud \xi_1 \, \ud \xi_2 \, \ud \xi_4 \\ 
        &\quad + \frakb(\xi) \iiint e^{i s \Omega_{2,\ulell}} h_1(s,\xi_1) \overline{\jxitwo h_2(s,\xi_2)} (\pxithree h_3)(s,\xi_3) \, \pvdots \frac{\varphi(\xi_4)}{\xi_4} \, \ud \xi_1 \, \ud \xi_2 \, \ud \xi_4 \\
        &=: \calJ_{2,1,2,1}^{\pvdots}(s,\xi) + \calJ_{2,1,2,2}^{\pvdots}(s,\xi) + \calJ_{2,1,2,3}^{\pvdots}(s,\xi).
    \end{aligned}
\end{equation*}
Similarly to the case of the cubic interactions with a Dirac kernel, the terms $\calJ_{2,1,1,3}^{\pvdots}(s,\xi)$, $\calJ_{2,1,2,3}^{\pvdots}(s,\xi)$, and $\calJ_{2,3}^{\pvdots}(s,\xi)$ combine to give 
\begin{equation*} 
    \begin{aligned}
        &\calJ_{2,1,1,3}^{\pvdots}(s,\xi) + \calJ_{2,1,2,3}^{\pvdots}(s,\xi) + \calJ_{2,3}^{\pvdots}(s,\xi) \\
        &= \frakb(\xi) \iiint e^{i s \Omega_{2,\ulell}} \bigl( \jxi - \jxione + \jxitwo \bigr) h_1(s,\xi_1) \overline{h_2(s,\xi_2)} (\pxithree h_3)(s,\xi_3) \, \pvdots \frac{\varphi(\xi_4)}{\xi_4} \, \ud \xi_1 \, \ud \xi_2 \, \ud \xi_4 \\
        &= \frakb(\xi) \iiint e^{i s \Omega_{2,\ulell}} h_1(s,\xi_1) \overline{h_2(s,\xi_2)} \jxithree (\pxithree h_3)(s,\xi_3) \, \pvdots \frac{\varphi(\xi_4)}{\xi_4} \, \ud \xi_1 \, \ud \xi_2 \, \ud \xi_4 \\
        &\quad + \frakb(\xi) \iiint e^{i s \Omega_{2,\ulell}} (-\Omega_{2,\ulell}) h_1(s,\xi_1) \overline{h_2(s,\xi_2)} (\pxithree h_3)(s,\xi_3) \, \pvdots \frac{\varphi(\xi_4)}{\xi_4} \, \ud \xi_1 \, \ud \xi_2 \, \ud \xi_4 \\
        &\quad + \ulell \, \frakb(\xi) \iiint e^{i s \Omega_{2,\ulell}} h_1(s,\xi_1) \overline{h_2(s,\xi_2)}(\pxithree h_3)(s,\xi_3) \, \varphi(\xi_4) \, \ud \xi_1 \, \ud \xi_2 \, \ud \xi_4 \\  &=: \calJ_{2,\mathrm{comb},1}^{\pvdots}(s,\xi) + \calJ_{2,\mathrm{comb},2}^{\pvdots}(s,\xi) + \calJ_{2,\mathrm{comb},3}^{\pvdots}(s,\xi).
    \end{aligned}
\end{equation*}
In comparison with the corresponding identity \eqref{equ:weighted_energy_proof_Dirac_combined_terms} for the cubic interactions with a Dirac kernel, the right-hand side of the predecing equation features the additional mild term $\calJ_{2,\mathrm{comb},3}^{\pvdots}(s,\xi)$.

Note that the Hilbert-type kernel $\pvdots \cosech(\frac{\pi}{2\ulg}\cdot)$ just corresponds to multiplication by $\tanh(\ulg \cdot)$ on the physical side up to constants. 
Apart from the terms $\calJ_{2,1,4}^{\pvdots}(s,\xi)$ and $\calJ_{2,\mathrm{comb},3}^{\pvdots}(s,\xi)$, which have no counter-parts in the previous case of the cubic interactions with a Dirac kernel, we can therefore estimate the contributions of all other terms by proceeding essentially in the same manner as for the corresponding terms in the previous case.
We emphasize that to obtain acceptable bounds for the contributions of the terms $\calJ_{2,1,3}^{\pvdots}(s,\xi)$ and $\calJ_{2,\mathrm{comb},2}^{\pvdots}(s,\xi)$, integrating by parts in time is necessary again.
In what follows we only discuss the contributions of the new terms $\calJ_{2,1,4}^{\pvdots}(s,\xi)$ and $\calJ_{2,\mathrm{comb},3}^{\pvdots}(s,\xi)$, and we omit the lengthy, but largely identical details for all other terms.

We write $\calJ_{2,1,4}^{\pvdots}(s,\xi)$ in terms of the auxiliary linear Klein-Gordon evolutions as
\begin{equation*}
    \begin{aligned}
    \calJ_{2,1,4}^{\pvdots}(s,\xi) 
    &= i s \cdot (2\pi)^{\frac32} \frakb(\xi) e^{-is(\jxi+\ulell\xi)} \whatcalF\bigl[ v_1(s,\cdot) \overline{v_2(s,\cdot)} \tilde{v}_3(s,\cdot) \widehat{\varphi}(\cdot) \bigr](\xi) 
    \end{aligned}
\end{equation*}
with
\begin{equation*}
    \tilde{v}_3(s,y) := \widehat{\calF}^{-1}\bigl[ e^{is(\jxithree+\ulell\xi_3)} \jxithree^{-1} (\jxithree + \ulell \xi_3) h_3(s,\xi_3) \bigr](y). 
\end{equation*}
Clearly, the linear Klein-Gordon evolution $\tilde{v}_3(s,y)$ enjoys the same bounds as the evolution $v_3(s,y)$.
We also express $\calJ_{2,\mathrm{comb},3}^{\pvdots}(s,\xi)$ in terms of the auxiliary linear Klein-Gordon evolutions as
\begin{equation*}
    \begin{aligned}
        \calJ_{2,\mathrm{comb},3}^{\pvdots}(s,\xi) 
        &= \ulell \, (2\pi)^{\frac32} \frakb(\xi) e^{-is(\jxi+\ulell\xi)} \whatcalF\bigl[ v_1(s,\cdot) \overline{v_2(s,\cdot)} w_3(s,\cdot) \widehat{\varphi}(\cdot) \bigr](\xi).
    \end{aligned}
\end{equation*}

\medskip 
\noindent \underline{Case 2.1: Contribution to \eqref{equ:pxi_cubic_reduced_claim}.}
Here our goal is to prove for $0 \leq t \leq T$,
\begin{equation} \label{equ:weighted_energy_proof_goal_case21}
    \biggl\| \int_0^t e^{-i\xi\theta(s)} \pxi \calI_{2, \mathrm{schem}}^{\pvdots}(s,\xi) \, \ud s \biggr\|_{L^2_\xi} \lesssim \varepsilon^3 \jt^\delta.
\end{equation}
We can freely spend one Sobolev weight $\jxi$ as in Case~1.1 so that we can work with the decomposition of $\jxi \pxi \calI_{2, \mathrm{schem}}^{\pvdots}(s,\xi)$ determined above. 
As explained before, we only discuss the contributions of the terms $\calJ_{2,1,4}^{\pvdots}(s,\xi)$ and $\calJ_{2,\mathrm{comb},3}^{\pvdots}(s,\xi)$.

\medskip 
\noindent {\it Bound for $\calJ_{2,1,4}^{\pvdots}(s,\xi)$:}
The improved low frequency behavior of one of the coefficients $\frakb(\xi)$, $\frakb_1(\xi_1)$, $\frakb_2(\xi_2)$, $\frakb_3(\xi_3)$ is key to obtain acceptable bounds for the contribution of $\calJ_{2,1,4}^{\pvdots}(s,\xi)$.

First, suppose that (at least) one of the coefficients $\frakb_1(\xi_1)$, $\frakb_2(\xi_2)$, $\frakb_3(\xi_3)$ exhibits the low frequency improvement.
For concreteness, say $\frakb_1(\xi_1)$ exhibits the low frequency improvement
\begin{equation*}
    \frakb_1(\xi_1) = \frac{\ulg(\xi_1+\ulell\jap{\xi_1})}{|\ulg(\xi_1+\ulell\jap{\xi_1})|\pm i}.
\end{equation*}
Then the linear evolution $v_1(s,y)$ enjoys the improved local decay \eqref{equ:consequences_aux_KG_improved_local_H1y_decay}.
Using also \eqref{equ:consequences_aux_KG_disp_decay}, we then obtain the acceptable bound
\begin{equation*}
    \begin{aligned}
        \biggl\| \int_0^t e^{-i\xi\theta(s)} \calJ_{2,1,4}^{\pvdots}(s,\xi) \, \ud s \biggr\|_{L^2_\xi} &\lesssim \int_0^t \js \cdot \bigl\| \whatcalF\bigl[ v_1(s,\cdot) \overline{v_2(s,\cdot)} \tilde{v}_3(s,\cdot) \widehat{\varphi}(\cdot) \bigr](\xi) \bigr\|_{L^2_\xi} \, \ud s \\ 
        &\lesssim \int_0^t \js \cdot \bigl\| v_1(s,y) \overline{v_2(s,y)} \tilde{v}_3(s,y) \widehat{\varphi}(y) \bigr\|_{L^2_y} \, \ud s \\ 
        &\lesssim \int_0^t \js \cdot \bigl\|\jy^{-1} v_1(s)\bigr\|_{L^2_y} \|v_2(s)\|_{L^\infty_y} \|\tilde{v}_3(s)\|_{L^\infty_y} \|\jy \widehat{\varphi}(y)\|_{L^\infty_y} \, \ud s \\ 
        &\lesssim \int_0^t \js \cdot \varepsilon \js^{-1+\delta} \cdot \varepsilon \js^{-\frac12} \cdot \varepsilon \js^{-\frac12} \, \ud s \lesssim \varepsilon^3 \jt^\delta.
    \end{aligned}
\end{equation*}
Second, suppose instead that the coefficient $\frakb(\xi)$ of the output frequency has the improved low frequency behavior.
Then we can invoke the dual integrated local energy decay estimate \eqref{equ:ILED_moving_center_frequency_side} together with the decay estimate \eqref{equ:consequences_aux_KG_disp_decay} to conclude the acceptable bound
\begin{equation*}
    \begin{aligned}
        &\biggl\| \int_0^t e^{-i\xi\theta(s)} \calJ_{2,1,4}^{\pvdots}(s,\xi) \, \ud s \biggr\|_{L^2_\xi} \\
        &=\biggl\| \int_0^t e^{-i\xi\theta(s)} \cdot i s \cdot (2\pi)^{\frac32} \frakb(\xi) e^{-is(\jxi+\ulell\xi)} \whatcalF\bigl[ v_1(s,\cdot) \overline{v_2(s,\cdot)} \tilde{v}_3(s,\cdot) \widehat{\varphi}(\cdot) \bigr](\xi) \, \ud s \biggr\|_{L^2_\xi} \\
        &\lesssim \Bigl\| \jy^2 \cdot s \cdot v_1(s,y) \overline{v_2(s,y)} \tilde{v}_3(s,y) \widehat{\varphi}(y) \Bigr\|_{L^2_s([0,t];H^1_y)} \\ 
        &\lesssim \bigl\| \jy^2 \widehat{\varphi}(y) \bigr\|_{H^1_y} \Bigl\| s \cdot \|v_1(s)\|_{W^{1,\infty}_y} \|v_2(s)\|_{W^{1,\infty}_y} \|\tilde{v}_3(s)\|_{W^{1,\infty}_y} \Bigr\|_{L^2_s([0,t])} \\ 
        &\lesssim \Bigl\| s \cdot \varepsilon^3 \js^{-\frac32} \Bigr\|_{L^2_s([0,t])} \lesssim \varepsilon^3 \log(\jt) \lesssim \varepsilon^3 \jt^\delta.     
    \end{aligned}
\end{equation*}

\medskip 
\noindent {\it Bound for $\calJ_{2,\mathrm{comb},3}^{\pvdots}(s,\xi)$:}
Using \eqref{equ:consequences_aux_KG_disp_decay} and \eqref{equ:weighted_energies_proof_wj_H1y_growth} we obtain via Plancherel's theorem the acceptable bound
\begin{equation*}
    \begin{aligned}
        &\biggl\| \int_0^t e^{-i\xi\theta(s)} \calJ_{2,\mathrm{comb},3}^{\pvdots}(s,\xi) \, \ud s \biggr\|_{L^2_\xi}  
        \lesssim \int_0^t \bigl\| \whatcalF\bigl[ v_1(s,\cdot) \overline{v_2(s,\cdot)} w_3(s,\cdot) \widehat{\varphi}(\cdot) \bigr](\xi) \bigr\|_{L^2_\xi} \, \ud s \\
        &\lesssim \int_0^t \bigl\| v_1(s,y) \overline{v_2(s,y)} w_3(s,y) \widehat{\varphi}(y) \bigr\|_{L^2_y} \, \ud s 
        \lesssim \int_0^t \|v_1(s)\|_{L^\infty_y} \|v_2(s)\|_{L^\infty_y} \|w_3(s)\|_{L^2_y} \|\widehat{\varphi}\|_{L^\infty_y} \, \ud s \\
        &\lesssim \int_0^t \varepsilon^2 \js^{-1} \cdot \varepsilon \js^\delta \, \ud s \lesssim \varepsilon^3 \jt^\delta.
    \end{aligned}
\end{equation*}

\medskip 
\noindent \underline{Case 2.2: Contribution to \eqref{equ:pxi_cubic_with_jxi_reduced_claim}.}
Now our goal is to prove the following (still slowly growing) weighted energy bounds with additional Sobolev regularity for $0 \leq t \leq T$,
\begin{equation} \label{equ:weighted_energy_proof_goal_case22}
    \biggl\| \jxi^2 \int_0^t e^{-i\xi\theta(s)} \pxi \calI_{2, \mathrm{schem}}^{\pvdots}(s,\xi) \, \ud s \biggr\|_{L^2_\xi} \lesssim \varepsilon^3 \jt^{2\delta}.
\end{equation}
As indicated before, we only discuss the contributions of the terms $\calJ_{2,1,4}^{\pvdots}(s,\xi)$ and $\calJ_{2,\mathrm{comb},3}^{\pvdots}(s,\xi)$.

\medskip 
\noindent {\it Bound for $\jxi \calJ_{2,1,4}^{\pvdots}(s,\xi)$:}
Again, the improved low frequency behavior of one of the coefficients $\frakb(\xi)$, $\frakb_1(\xi_1)$, $\frakb_2(\xi_2)$, $\frakb_3(\xi_3)$ is key to obtain acceptable bounds for the contribution of $\jxi \calJ_{2,1,4}^{\pvdots}(s,\xi)$.

First, suppose that (at least) one of the coefficients $\frakb_1(\xi_1)$, $\frakb_2(\xi_2)$, $\frakb_3(\xi_3)$ exhibits the low frequency improvement.
For concreteness, say $\frakb_1(\xi_1)$ has the low frequency improvement.
Using again the improved local decay estimate \eqref{equ:consequences_aux_KG_improved_local_H1y_decay} along with \eqref{equ:consequences_aux_KG_disp_decay}, we arrive at the sufficient bound
\begin{equation*}
    \begin{aligned}
        \biggl\| \jxi \int_0^t e^{-i\xi\theta(s)} \calJ_{2,1,4}^{\pvdots}(s,\xi) \, \ud s \biggr\|_{L^2_\xi} &\lesssim \int_0^t \js \cdot \bigl\| \jxi \whatcalF\bigl[ v_1(s,\cdot) \overline{v_2(s,\cdot)} \tilde{v}_3(s,\cdot) \widehat{\varphi}(\cdot) \bigr](\xi) \bigr\|_{L^2_\xi} \, \ud s \\ 
        &\lesssim \int_0^t \js \cdot \bigl\| v_1(s,y) \overline{v_2(s,y)} \tilde{v}_3(s,y) \widehat{\varphi}(y) \bigr\|_{H^1_y} \, \ud s \\ 
        &\lesssim \int_0^t \js \cdot \bigl\|\jy^{-1} v_1(s)\bigr\|_{H^1_y} \|v_2(s)\|_{W^{1,\infty}_y} \|\tilde{v}_3(s)\|_{W^{1,\infty}_y} \|\jy \widehat{\varphi}(y)\|_{W^{1,\infty}_y} \, \ud s \\ 
        &\lesssim \int_0^t \js \cdot \varepsilon \js^{-1+\delta} \cdot \varepsilon \js^{-\frac12} \cdot \varepsilon \js^{-\frac12} \, \ud s \lesssim \varepsilon^3 \jt^\delta.
    \end{aligned}
\end{equation*}
Second, suppose instead that the coefficient $\frakb(\xi)$ of the output frequency has the improved low frequency behavior.
In this case we can invoke the dual integrated local energy decay estimate \eqref{equ:ILED_moving_center_frequency_side_with_jxi2} with additional Sobolev regularity. Using H\"older's inequality, Sobolev embedding, and the local decay estimate \eqref{equ:consequences_aux_KG_local_H3y_decay}, we then conclude the sufficient bound
\begin{equation*}
    \begin{aligned}
        &\biggl\| \jxi \int_0^t e^{-i\xi\theta(s)} \calJ_{2,1,4}^{\pvdots}(s,\xi) \, \ud s \biggr\|_{L^2_\xi} \\
        &=\biggl\| \jxi \int_0^t e^{-i\xi\theta(s)} \cdot i s \cdot (2\pi)^{\frac32} \frakb(\xi) e^{-is(\jxi+\ulell\xi)} \whatcalF\bigl[ v_1(s,\cdot) \overline{v_2(s,\cdot)} \tilde{v}_3(s,\cdot) \widehat{\varphi}(\cdot) \bigr](\xi) \, \ud s \biggr\|_{L^2_\xi} \\
        &\lesssim \Bigl\| \jy^2 \cdot s \cdot v_1(s,y) \overline{v_2(s,y)} \tilde{v}_3(s,y) \widehat{\varphi}(y) \Bigr\|_{L^2_s([0,t];H^3_y)} \\ 
        &\lesssim \bigl\| \jy^5 \widehat{\varphi}(y) \bigr\|_{H^3_y} \Bigl\| s \cdot \|\jy^{-1} v_1(s)\|_{H^3_y} \|\jy^{-1} v_2(s)\|_{H^3_y} \|\jy^{-1} \tilde{v}_3(s)\|_{H^3_y} \Bigr\|_{L^2_s([0,t])} \\ 
        &\lesssim \Bigl\| s \cdot \varepsilon^3 \js^{-\frac32} \Bigr\|_{L^2_s([0,t])} \lesssim \varepsilon^3 \log(\jt) \lesssim \varepsilon^3 \jt^\delta.     
    \end{aligned}
\end{equation*}

\medskip 
\noindent {\it Bound for $\jxi \calJ_{2,\mathrm{comb},3}^{\pvdots}(s,\xi)$:}
Using \eqref{equ:consequences_aux_KG_disp_decay} and \eqref{equ:weighted_energies_proof_wj_H1y_growth} we obtain via Plancherel's theorem the sufficient bound
\begin{equation*}
    \begin{aligned}
        &\biggl\| \jxi \int_0^t e^{-i\xi\theta(s)} \calJ_{2,\mathrm{comb},3}^{\pvdots}(s,\xi) \, \ud s \biggr\|_{L^2_\xi}  
        \lesssim \int_0^t \bigl\| \jxi \whatcalF\bigl[ v_1(s,\cdot) \overline{v_2(s,\cdot)} w_3(s,\cdot) \widehat{\varphi}(\cdot) \bigr](\xi) \bigr\|_{L^2_\xi} \, \ud s \\
        &\lesssim \int_0^t \bigl\| v_1(s,y) \overline{v_2(s,y)} w_3(s,y) \widehat{\varphi}(y) \bigr\|_{H^1_y} \, \ud s 
        \lesssim \int_0^t \|v_1(s)\|_{W^{1,\infty}_y} \|v_2(s)\|_{W^{1,\infty}_y} \|w_3(s)\|_{H^1_y} \|\widehat{\varphi}\|_{W^{1,\infty}_y} \, \ud s \\
        &\lesssim \int_0^t \varepsilon^2 \js^{-1} \cdot \varepsilon \js^\delta \, \ud s \lesssim \varepsilon^3 \jt^\delta.
    \end{aligned}
\end{equation*}

This finishes the discussion of the weighted energy estimates for the cubic interactions with a Hilbert-type kernel.

\medskip
\noindent \underline{Case 3: Regular cubic interactions.}
Finally, we consider the regular cubic interactions. In view of the fine structure of the cubic spectral distributions analyzed in Subsection~\ref{subsec:cubic_spectral_distributions}, the term $\calI_2^{\mathrm{reg}}(s,\xi)$ is a linear combination of terms of the schematic form
\begin{equation*}
    \begin{aligned}
        \calI_{2, \mathrm{schem}}^{\mathrm{reg}}(s,\xi) &:= \frakb(\xi) \iiint e^{is \Phi_{2,\ulell}(\xi,\xi_1,\xi_2,\xi_3)} h_1(s,\xi_1) \, \overline{h_2(s,\xi_2)} h_3(s,\xi_3) \\ 
        &\quad \quad \quad \quad \quad \quad \quad \quad \times \ulg^{-1} \widehat{\calF}^{-1}\bigl[\varphi\bigr]\bigl( \ulg^{-1} (-\xi+\xi_1-\xi_2+\xi_3) \bigr) \, \ud \xi_1 \, \ud \xi_2 \, \ud \xi_3,
    \end{aligned}
\end{equation*}
with the phase
\begin{equation*}
    \Phi_{2,\ulell}(\xi,\xi_1,\xi_2,\xi_3) := -(\jxi + \ulell \xi) + (\jxione + \ulell \xi_1) - (\jxitwo + \ulell \xi_2) + (\jxithree + \ulell \xi_3),
\end{equation*}
and for some coefficients $\frakb, \frakb_1, \frakb_2, \frakb_3 \in W^{1,\infty}(\bbR)$ such that the inputs are given by 
\begin{equation*}
    h_j(s,\xi_j) := \jap{\xi_j}^{-1} \frakb_j(\xi_j) \gulellsh(s,\xi_j), \quad 1 \leq j \leq 3,
\end{equation*}
and where $\varphi \in \calS(\bbR)$ is some Schwartz function.
In terms of the linear Klein-Gordon evolutions defined in \eqref{equ:weighted_energies_proof_vj_def}, we have
\begin{equation*}
    \begin{aligned}
        \calI_{2, \mathrm{schem}}^{\mathrm{reg}}(s,\xi) &= (2\pi)^{\frac32} e^{-is(\jxi+\ulell\xi)} \frakb(\xi) \widehat{\calF}\bigl[ v_1(s,\cdot) \overline{v_2(s,\cdot)} v_3(s,\cdot) \varphi(\ulg \cdot) \bigr](\xi). 
    \end{aligned}
\end{equation*}
By direct computation, we find that
\begin{equation*}
    \begin{aligned}
        \pxi \calI_{2, \mathrm{schem}}^{\mathrm{reg}}(s,\xi) = \sum_{k=1}^3 \calJ_{2,k}^{\mathrm{reg}}(s,\xi),
    \end{aligned}
\end{equation*}
where 
\begin{align*}
    \calJ_{2,1}^{\mathrm{reg}}(s,\xi) &:= (2\pi)^{\frac32} (-i) (\xi+\ulell\jxi) \jxi^{-1} e^{-is(\jxi+\ulell\xi)} \cdot s \cdot \frakb(\xi) \widehat{\calF}\bigl[ v_1(s,\cdot) \overline{v_2(s,\cdot)} v_3(s,\cdot) \varphi(\ulg \cdot) \bigr](\xi), \\
    \calJ_{2,2}^{\mathrm{reg}}(s,\xi) &:= (2\pi)^{\frac32} e^{-is(\jxi+\ulell\xi)} (\pxi \frakb)(\xi) \widehat{\calF}\bigl[ v_1(s,\cdot) \overline{v_2(s,\cdot)} v_3(s,\cdot) \varphi(\ulg \cdot) \bigr](\xi), \\           
    \calJ_{2,3}^{\mathrm{reg}}(s,\xi) &:= (2\pi)^{\frac32} e^{-is(\jxi+\ulell\xi)} \frakb(\xi) \pxi \widehat{\calF}\bigl[ v_1(s,\cdot) \overline{v_2(s,\cdot)} v_3(s,\cdot) \varphi(\ulg \cdot) \bigr](\xi).
\end{align*}

\medskip 
\noindent \underline{Case 3.1: Contribution to \eqref{equ:pxi_cubic_reduced_claim}.}
Here our goal is to prove for $0 \leq t \leq T$,
\begin{equation} \label{equ:weighted_energy_proof_goal_case31}
    \biggl\| \int_0^t e^{-i\xi\theta(s)} \pxi \calI_{2, \mathrm{schem}}^{\mathrm{reg}}(s,\xi) \, \ud s \biggr\|_{L^2_\xi} \lesssim \varepsilon^3 \jt^\delta.
\end{equation}
We estimate the contributions of $\calJ^{\mathrm{reg}}_{2,k}(s,\xi)$, $1 \leq k \leq 3$, to \eqref{equ:weighted_energy_proof_goal_case31} separately.
For the contribution of the more delicate term $\calJ^{\mathrm{reg}}_{2,1}(s,\xi)$ we use the dual integrated local energy decay estimate \eqref{equ:ILED_moving_center_frequency_side} from Proposition~\ref{prop:ILED_moving_center_frequency_formulation} along with \eqref{equ:consequences_aux_KG_disp_decay} to obtain the acceptable bound 
\begin{equation*}
    \begin{aligned}
        \biggl\| \int_0^t e^{-i\xi\theta(s)} \calJ_{2,1}^{\reg}(s,\xi) \, \ud s \biggr\|_{L^2_\xi} &\lesssim \Bigl\| s \cdot \jy^2 v_1(s,y) \overline{v_2(s,y)} v_3(s,y) \varphi(\ulg y) \Bigr\|_{L^2_s([0,t];H^1_y)} \\ 
        &\lesssim \bigl\| \jy^2 \varphi(\ulg \cdot) \bigr\|_{H^1_y} \Bigl\| s \cdot \|v_1(s)\|_{W^{1,\infty}_y} \|v_2(s)\|_{W^{1,\infty}_y} \|v_3(s)\|_{W^{1,\infty}_y} \Bigr\|_{L^2_s([0,t])} \\ 
        &\lesssim \bigl\| s \cdot \varepsilon^3 \js^{-\frac32} \bigr\|_{L^2_s([0,t])} \lesssim \varepsilon^3 \log(\jt) \lesssim \varepsilon^3 \jt^{\delta}.
    \end{aligned}
\end{equation*}
Next, it is straightforward to bound the contribution of the second term $\calJ^{\mathrm{reg}}_{2,2}(s,\xi)$ using \eqref{equ:consequences_aux_KG_disp_decay} and Plancherel's theorem. We find 
\begin{equation*}
    \begin{aligned}
        \biggl\| \int_0^t e^{-i\xi\theta(s)} \calJ_{2,2}^{\reg}(s,\xi) \, \ud s \biggr\|_{L^2_\xi} &\lesssim \int_0^t \bigl\| \widehat{\calF}\bigl[ v_1(s,\cdot) \overline{v_2(s,\cdot)} v_3(s,\cdot) \varphi(\ulg \cdot) \bigr](\xi) \bigr\|_{L^2_\xi} \, \ud s \\
        &\lesssim \int_0^t \bigl\| \varphi(\ulg y) \bigr\|_{L^2_y} \|v_1(s)\|_{L^\infty_y} \|v_2(s)\|_{L^\infty_y} \|v_3(s)\|_{L^\infty_y} \, \ud s \\ 
        &\lesssim \int_0^t \varepsilon^3 \js^{-\frac32} \, \ud s \lesssim \varepsilon^3.
    \end{aligned}
\end{equation*}
The third term $\calJ^{\mathrm{reg}}_{2,3}(s,\xi)$ can be treated in the same manner.

\medskip 
\noindent \underline{Case 3.2: Contribution to \eqref{equ:pxi_cubic_with_jxi_reduced_claim}.}
Now our goal is to prove the (still slowly growing) weighted energy bound with additional Sobolev regularity for $0 \leq t \leq T$,
\begin{equation*}
    \biggl\| \jxi^2 \int_0^t e^{-i\xi\theta(s)} \pxi \calI_{2, \mathrm{schem}}^{\mathrm{reg}}(s,\xi) \, \ud s \biggr\|_{L^2_\xi} \lesssim \varepsilon^3 \jt^{2\delta}.
\end{equation*}
For the contribution of the more delicate term $\calJ^{\mathrm{reg}}_{2,1}(s,\xi)$ we now use the dual integrated local energy decay estimate \eqref{equ:ILED_moving_center_frequency_side_with_jxi2} with additional Sobolev regularity from Proposition~\ref{prop:ILED_moving_center_frequency_formulation}. Using the local decay estimate \eqref{equ:consequences_aux_KG_local_H3y_decay} along with Sobolev embedding yields the sufficient bound 
\begin{equation*}
    \begin{aligned}
        &\biggl\| \jxi^2 \int_0^t e^{-i\xi\theta(s)} \calJ_{2,1}^{\reg}(s,\xi) \, \ud s \biggr\|_{L^2_\xi} \lesssim \Bigl\| s \cdot \jy^2  v_1(s,y) \overline{v_2(s,y)} v_3(s,y) \varphi(\ulg y) \Bigr\|_{L^2_s([0,t];H^3_y)} \\ 
        &\lesssim \bigl\| \jy^5 \varphi(\ulg y) \bigr\|_{H^3_y} \Bigl\| s \cdot \bigl\|\jy^{-1} v_1(s)\bigr\|_{H^3_y} \bigl\|\jy^{-1} v_2(s)\bigr\|_{H^3_y}\bigl\|\jy^{-1} v_3(s)\bigr\|_{H^3_y} \Bigr\|_{L^2_s([0,t])} \\ 
        &\lesssim \bigl\| s \cdot \varepsilon^3 \js^{-\frac32} \bigr\|_{L^2_s([0,t])} \lesssim \varepsilon^3 \log(\jt) \lesssim \varepsilon^3 \jt^{\delta}.
    \end{aligned}
\end{equation*}
Then it is again straightforward to bound the contribution of the second term $\calJ^{\mathrm{reg}}_{2,2}(s,\xi)$ via Plancherel's theorem, Sobolev embedding, and the local decay estimate~\eqref{equ:consequences_aux_KG_local_H3y_decay},
\begin{equation*}
    \begin{aligned}
        &\biggl\| \jxi^2 \int_0^t e^{-i\xi\theta(s)} \calJ_{2,2}^{\reg}(s,\xi) \, \ud s \biggr\|_{L^2_\xi} \lesssim \int_0^t \bigl\| \jxi^2 \widehat{\calF}\bigl[ v_1(s,\cdot) \overline{v_2(s,\cdot)} v_3(s,\cdot) \varphi(\ulg \cdot) \bigr](\xi) \bigr\|_{L^2_\xi} \, \ud s \\
        &\lesssim \int_0^t \bigl\| v_1(s,\cdot) \overline{v_2(s,\cdot)} v_3(s,\cdot) \varphi(\ulg \cdot) \bigr\|_{H^2_y} \, \ud s \\
        &\lesssim \int_0^t \bigl\| \jy^3 \varphi(\ulg \cdot) \bigr\|_{H^2_y} \bigl\| \jy^{-1} v_1(s) \bigr\|_{H^2_y} \bigl\| \jy^{-1} v_2(s) \bigr\|_{H^2_y} \bigl\| \jy^{-1} v_3(s) \bigr\|_{H^2_y}\, \ud s \\ 
        &\lesssim \int_0^t \varepsilon^3 \js^{-\frac32} \, \ud s \lesssim \varepsilon^3.
    \end{aligned}
\end{equation*}
The third term $\calJ^{\mathrm{reg}}_{2,3}(s,\xi)$ can again be treated in the same manner.

This finishes the discussion of the weighted energy estimates for the regular cubic interaction terms, and therefore concludes the proof of the lemma.  
\end{proof}

In the final lemma of this section we conclude the weighted energy bounds for the contributions of all nonlinear terms in $\calR(t,\xi)$ defined in \eqref{equ:setting_up_definition_calR}. We recall that these should be thought of as spatially localized with cubic-type $\jt^{-\frac32+\delta}$ decay apart from $\calR_2(\usubeone(t))$. 
The integrated local energy decay estimates with a moving center from Proposition~\ref{prop:ILED_moving_center_frequency_formulation} formulated on the frequency side are key to obtain acceptable bounds for the contributions of these spatially localized terms.
Due to the differing slow growth rates for the weighted energy estimates without and with additional Sobolev regularity this requires some care, in particular for the contributions of the renormalized quadratic nonlinearity $\calQ_{\ulell,r}\bigl(\usubeone(t)\bigr)$ defined in \eqref{equ:setting_up_definition_calQr}.

\begin{lemma} \label{lem:pxi_calR}
    Suppose the assumptions in the statement of Proposition~\ref{prop:profile_bounds} are in place.
    Then we have for all $0 \leq t \leq T$ that
    \begin{equation} \label{equ:pxi_calR_duhamel_bound}
        \biggl\| \int_0^t \pxi \biggl( e^{-i\xi\theta(s)} e^{-is(\jxi+\ulell\xi)} \calR(s,\xi) \biggr) \, \ud s \biggr\|_{L^2_\xi} \lesssim \varepsilon^2 \jt^\delta,
    \end{equation}
    and
    \begin{equation} \label{equ:pxi_with_jxi_calR_duhamel_bound}
        \biggl\| \jxi^2 \int_0^t \pxi \biggl( e^{-i\xi\theta(s)} e^{-is(\jxi+\ulell\xi)} \calR(s,\xi) \biggr) \, \ud s \biggr\|_{L^2_\xi} \lesssim \varepsilon^2 \jt^{2\delta}.
    \end{equation}    
\end{lemma}
\begin{proof}
Throughout we consider times $0 \leq t \leq T$. 
By the product rule for differentiation we have 
\begin{equation} \label{equ:pxi_without_reg_calR_list}
    \begin{aligned}
        &\biggl\| \int_0^t \pxi \biggl( e^{-i\xi\theta(s)} e^{-is(\jxi+\ulell\xi)} \calR(s,\xi) \biggr) \, \ud s \biggr\|_{L^2_\xi} \\
        &\lesssim \int_0^t |\theta(s)| \bigl\|\calR(s,\xi)\bigr\|_{L^2_\xi} \, \ud s + \int_0^t \bigl\| \pxi \calR(s,\xi) \bigr\|_{L^2_\xi} \, \ud s \\ 
        &\quad + \biggl\| \int_0^t e^{-i\xi\theta(s)} e^{-is(\jxi+\ulell\xi)} \cdot s \cdot \bigl(\xi+\ulell\jxi\bigr) \jxi^{-1} \calR(s,\xi) \, \ud s \biggr\|_{L^2_\xi},
    \end{aligned}
\end{equation}
and analogously
\begin{equation} \label{equ:pxi_with_jxi_calR_list}
    \begin{aligned}
        &\biggl\| \jxi^2 \int_0^t \pxi \biggl( e^{-i\xi\theta(s)} e^{-is(\jxi+\ulell\xi)} \calR(s,\xi) \biggr) \, \ud s \biggr\|_{L^2_\xi} \\
        &\lesssim \int_0^t |\theta(s)| \bigl\|\jxi^2 \calR(s,\xi)\bigr\|_{L^2_\xi} \, \ud s + \int_0^t \bigl\| \jxi^2 \pxi \calR(s,\xi) \bigr\|_{L^2_\xi} \, \ud s \\ 
        &\quad + \biggl\| \jxi^2 \int_0^t  e^{-i\xi\theta(s)} e^{-is(\jxi+\ulell\xi)} \cdot s \cdot \bigl(\xi+\ulell\jxi\bigr) \jxi^{-1} \calR(s,\xi) \, \ud s \biggr\|_{L^2_\xi}.
    \end{aligned}
\end{equation}
Using \eqref{equ:consequences_theta_growth_bound}, \eqref{equ:energy_preparations_jxi2_L2_calR}, and \eqref{equ:energy_preparations_jxi2_pxi_L2_calR}, the first two terms on the right-hand side of \eqref{equ:pxi_without_reg_calR_list} as well as the first two terms on the right-hand side of \eqref{equ:pxi_with_jxi_calR_list} are easily bounded by 
\begin{equation*}
    \begin{aligned}
        &\int_0^t |\theta(s)| \bigl\|\calR(s,\xi)\bigr\|_{L^2_\xi} \, \ud s + \int_0^t \bigl\| \pxi \calR(s,\xi) \bigr\|_{L^2_\xi} \, \ud s \\
        &\quad + \int_0^t |\theta(s)| \bigl\|\jxi^2 \calR(s,\xi)\bigr\|_{L^2_\xi} \, \ud s + \int_0^t \bigl\| \jxi^2 \pxi \calR(s,\xi) \bigr\|_{L^2_\xi} \, \ud s \\
        &\lesssim \int_0^t \varepsilon \js^\delta \cdot \varepsilon^2 \js^{-\frac32+2\delta} \, \ud s + \int_0^t \varepsilon^2 \js^{-1+\delta} \, \ud s \lesssim \varepsilon^2 \jt^\delta.
    \end{aligned}
\end{equation*}
This bound suffices for both asserted estimates \eqref{equ:pxi_calR_duhamel_bound} and \eqref{equ:pxi_with_jxi_calR_duhamel_bound}.
Now we turn to the main third terms on the right-hand side of \eqref{equ:pxi_without_reg_calR_list} and on the right-hand side of \eqref{equ:pxi_with_jxi_calR_list}.

First, we recall that all terms in $\calR(s,\xi)$ apart from the contribution of $\calR_2\bigl(\usubeone(s)\bigr)$ are spatially localized with at least cubic-type $\js^{-\frac32+\delta}$ decay. 
We will therefore use the dual integrated local energy decay estimates from Proposition~\ref{prop:ILED_moving_center_frequency_formulation} to bound their contributions, where we exploit that the distorted Fourier basis element $\overline{\eulsharp(y,\xi)}$ can be written as a sum of terms of the form $e^{-iy\xi} \fraka(y) \frakb(\xi)$ with $\fraka(y)$ smooth with bounded derivatives and $\frakb \in W^{1,\infty}$.

In what follows we consider step by step each term in the definition \eqref{equ:setting_up_definition_calR} of $\calR(s,\xi)$. If needed, we separately deduce sufficient bounds for the contribution of each term to the first asserted bound \eqref{equ:pxi_calR_duhamel_bound} and to the second asserted bound \eqref{equ:pxi_with_jxi_calR_duhamel_bound}.

\medskip
\noindent \underline{Contribution of $(\dot{q}(s) - \ulell) \calL_\ulell\bigl(\bmu(s)\bigr)(\xi)$.}
Here we can simultaneously deal with the contributions to \eqref{equ:pxi_calR_duhamel_bound} and to \eqref{equ:pxi_with_jxi_calR_duhamel_bound}.
We observe that the coefficients \eqref{equ:setting_up_calL_coefficients_mulell} in the definition \eqref{equ:setting_up_calL_definition} of $\calL_\ulell(\bmu(s))(\xi)$ satisfy
\begin{equation*}
    \sup_{y, \xi \in \bbR} \, \Bigl( \bigl| \jxi \py m_\ulell^{\#}(\ulg y, \xi) \bigr| + \bigl| \jxi \py^2 m_\ulell^{\#}(\ulg y,\xi) \bigr| \Bigr) \lesssim 1.
\end{equation*}
Thus, using the dual integrated local energy decay estimate \eqref{equ:ILED_moving_center_frequency_side_with_jxi2} from Proposition~\ref{prop:ILED_moving_center_frequency_formulation} for the contribution of $u_1(s,y)$ in $\calL_\ulell(\bmu(s))(\xi)$, while using the dual integrated local energy decay estimate variant \eqref{equ:ILED_moving_center_frequency_side_with_jxi2_plus_decaying_coefficient} for the contributions of $u_2(s,y)$ in $\calL_\ulell(\bmu(s))(\xi)$, we infer from \eqref{equ:consequences_qdot_minus_ulell_decay} and the local decay estimates \eqref{equ:consequences_local_decay_uonetwo} that
\begin{equation*}
    \begin{aligned}
        &\biggl\| \int_0^t e^{-i\xi\theta(s)} e^{-is(\jxi+\ulell\xi)} \cdot s \cdot \bigl(\xi+\ulell\jxi\bigr) \jxi^{-1}  \calL_\ulell\bigl(\bmu(s)\bigr)(\xi) \, \ud s \biggr\|_{L^2_\xi} \\ 
        &\quad + \biggl\| \jxi^2 \int_0^t e^{-i\xi\theta(s)} e^{-is(\jxi+\ulell\xi)} \cdot s \cdot \bigl(\xi+\ulell\jxi\bigr) \jxi^{-1}  \calL_\ulell\bigl(\bmu(s)\bigr)(\xi) \, \ud s \biggr\|_{L^2_\xi} \\ 
        &\lesssim \Bigl\| \jy^2 \cdot s \cdot (\dot{q}(s)-\ulell) \cdot \sech^2(\ulg y) u_1(s,y) \Bigr\|_{L^2_s([0,t]; H^3_y)} \\ 
        &\qquad + \Bigl\| \jy^2 \cdot s \cdot (\dot{q}(s)-\ulell) \cdot \sech^2(\ulg y) u_2(s,y) \Bigr\|_{L^2_s([0,t]; H^2_y)} \\
        &\lesssim \Bigl\| s \cdot |\dot{q}(s)-\ulell| \cdot \bigl\| \jy^{-1} u_1(s) \bigr\|_{H^3_y} \Bigr\|_{L^2_s([0,t])} + \Bigl\| s \cdot |\dot{q}(s)-\ulell| \cdot \bigl\| \jy^{-1} u_2(s) \bigr\|_{H^2_y} \Bigr\|_{L^2_s([0,t])} \\
        &\lesssim \Bigl\| s \cdot \varepsilon \js^{-1+\delta} \cdot \varepsilon \js^{-\frac12} \Bigr\|_{L^2_s([0,t])} \\
        &\lesssim \varepsilon^2 \jt^\delta.
    \end{aligned}
\end{equation*}
This bound suffices for the contributions to both \eqref{equ:pxi_calR_duhamel_bound} and \eqref{equ:pxi_with_jxi_calR_duhamel_bound}.

\medskip 
\noindent \underline{Contribution of $\calMod(s)$.}
    By Corollary~\ref{cor:TellP} we have
    \begin{equation*}
        \begin{aligned}
            \calF_{\ulell,D}^{\#}\bigl[ \bigl(\calMod\bigr)_1(t) \bigr](\xi) - \calF_{\ulell}^{\#}\bigl[ \bigl(\calMod\bigr)_2(t) \bigr](\xi) = \calF_{\ulell,D}^{\#}\bigl[ \bigl(\ulPe \calMod\bigr)_1(t) \bigr](\xi) - \calF_{\ulell}^{\#}\bigl[ \bigl(\ulPe \calMod\bigr)_2(t) \bigr](\xi).
        \end{aligned}
    \end{equation*}
    Moreover, since $\ulPe \bmY_{j,\ulell} = 0$ for $j=0,1$, we can write
    \begin{equation*}
        \begin{aligned}
            \ulPe \calMod(t,y) = - \bigl(\dot{q}(t)-\ell(t)\bigr) \ulPe \bigl( \bmY_{0,\ell(t)}(y) - \bmY_{0,\ulell}(y) \bigr) - \dot{\ell}(t) \ulPe \bigl( \bmY_{1,\ell(t)}(y) - \bmY_{1,\ulell}(y) \bigr).
        \end{aligned}
    \end{equation*}
    Then we can simultaneously deal with the contributions to \eqref{equ:pxi_calR_duhamel_bound} and to \eqref{equ:pxi_with_jxi_calR_duhamel_bound}.
    For the contribution of $(\ulPe \calMod)_1(s)$ we first observe that the distorted Fourier basis element in $\calFulellDsh$ satisfies
    \begin{equation*}
        \sup_{y,\xi \in \bbR} \, \bigl| \jxi^{-1} D^\ast e_\ulell^{\#}(y,\xi) \bigr| \lesssim 1.
    \end{equation*}
    Hence, using the dual integrated local energy decay estimate variant \eqref{equ:ILED_moving_center_frequency_side_with_jxi2_plus_growing_coefficient} along with the bounds \eqref{equ:prop_profile_bounds_assumption1}, \eqref{equ:consequences_crude_decay_modulation_parameters} for the modulation parameters, we find
    \begin{equation*}
        \begin{aligned}
            &\biggl\| \int_0^t e^{-i\xi\theta(s)} e^{-is(\jxi+\ulell\xi)} \cdot s \cdot \bigl(\xi+\ulell\jxi\bigr) \jxi^{-1}  \calFulellDsh\bigl[ \bigl(\ulPe \calMod\bigr)_1(s) \bigr](\xi) \, \ud s \biggr\|_{L^2_\xi} \\
            &\quad + \biggl\| \jxi^2 \int_0^t e^{-i\xi\theta(s)} e^{-is(\jxi+\ulell\xi)} \cdot s \cdot \bigl(\xi+\ulell\jxi\bigr) \jxi^{-1}  \calFulellDsh\bigl[ \bigl(\ulPe \calMod\bigr)_1(s) \bigr](\xi) \, \ud s \biggr\|_{L^2_\xi} \\
            &\lesssim \Bigl\| \jy^2 \cdot s \cdot \bigl(\ulPe \calMod\bigr)_1(s) \Bigr\|_{L^2_s([0,t]; H^4_y)} + \bigl\| s \cdot \bigl( |(\dot{q}-\ell)(s)| + |\dot{\ell}(s)| \bigr) \cdot |\ell(s)-\ulell| \bigr\|_{L^2_s([0,t])} \\ 
            &\lesssim \bigl\| s \cdot \varepsilon \js^{-1} \cdot \varepsilon \js^{-1+\delta} \bigr\|_{L^2_s([0,t])} \lesssim \varepsilon^2.
        \end{aligned}
    \end{equation*}
    The contribution of $(\ulPe \calMod)_2(s)$ can be treated analogously using the dual integrated local energy decay estimate \eqref{equ:ILED_moving_center_frequency_side_with_jxi2}.

\medskip 
\noindent \underline{Contribution of $\calC_{\ell(s)}\bigl( \usubeone(s) \bigr)$.}
Using \eqref{equ:ILED_moving_center_frequency_side_with_jxi2}, the Sobolev embedding $H^1_y(\bbR) \hookrightarrow L^\infty_y(\bbR)$, and the local decay estimates \eqref{equ:consequences_local_decay_usubeonetwo}, we find that
\begin{equation*}
    \begin{aligned}
        &\biggl\| \int_0^t e^{-i\xi\theta(s)} e^{-is(\jxi+\ulell\xi)} \cdot s \cdot \bigl(\xi+\ulell\jxi\bigr) \jxi^{-1}  \calFulellsh\bigl[\calC_{\ell(t)}\bigl( \usubeone(s) \bigr)\bigr](\xi) \, \ud s \biggr\|_{L^2_\xi} \\
        &\quad + \biggl\| \jxi^2 \int_0^t e^{-i\xi\theta(s)} e^{-is(\jxi+\ulell\xi)} \cdot s \cdot \bigl(\xi+\ulell\jxi\bigr) \jxi^{-1}  \calFulellsh\bigl[\calC_{\ell(t)}\bigl( \usubeone(s) \bigr)\bigr](\xi) \, \ud s \biggr\|_{L^2_\xi} \\
        &\lesssim \Bigl\| \jy^2 \cdot s \cdot \sech^2\bigl(\gamma(s) y\bigr) \usubeone(s,y)^3 \Bigr\|_{L^2_s([0,t]; H^3_y)} 
        \lesssim \Bigl\| s \cdot \bigl\|\jy^{-1} \usubeone(s)\bigr\|_{H^3_y}^3 \Bigr\|_{L^2_s([0,t])} \\ 
        &\lesssim \Bigl\| s \cdot \varepsilon^3 \js^{-\frac32} \Bigr\|_{L^2_s([0,t])} \lesssim \varepsilon^3 \log(2+\jt).
    \end{aligned}
\end{equation*}
This suffices for the contributions to both \eqref{equ:pxi_calR_duhamel_bound} and \eqref{equ:pxi_with_jxi_calR_duhamel_bound}.

\medskip 
\noindent \underline{Contribution of $\calR_1\bigl(\usubeone(s)\bigr)$.}
The quartic nonlinearity $\calR_1\bigl(\usubeone(s)\bigr)$ is also spatially localized. Its contribution can therefore be estimated as in the preceding treatment of the spatially localized cubic nonlinearity $\calC_{\ell(s)}\bigl( \usubeone(s) \bigr)$.

\medskip 
\noindent \underline{Contribution of $\calR_2\bigl(\usubeone(s)\bigr)$.}
We estimate the contributions of the quintic and higher nonlinearities directly using the mapping property \eqref{equ:mapping_property_calFulellsh_jxi2} and the bounds \eqref{equ:consequences_dispersive_decay_ulPe_radiation}, \eqref{equ:consequences_sobolev_bound_radiation},
\begin{equation*}
    \begin{aligned}
        &\biggl\| \int_0^t e^{-i\xi\theta(s)} e^{-is(\jxi+\ulell\xi)} \cdot s \cdot \bigl(\xi+\ulell\jxi\bigr) \jxi^{-1}  \calFulellsh\bigl[\calR_2\bigl(\usubeone(s)\bigr)\bigr](\xi) \, \ud s \biggr\|_{L^2_\xi} \\
        &\quad + \biggl\| \jxi^2 \int_0^t e^{-i\xi\theta(s)} e^{-is(\jxi+\ulell\xi)} \cdot s \cdot \bigl(\xi+\ulell\jxi\bigr) \jxi^{-1}  \calFulellsh\bigl[\calR_2\bigl(\usubeone(s)\bigr)\bigr](\xi) \, \ud s \biggr\|_{L^2_\xi} \\
        &\lesssim \int_0^t s \cdot \bigl\| \calR_2\bigl(\usubeone(s)\bigr) \bigr\|_{H^2_y} \, \ud s 
        \lesssim \int_0^t s \cdot \|\usubeone(s)\|_{H^2_y} \|\usubeone(s)\|_{W^{1,\infty}_y}^4 \Bigl( 1 + \|\usubeone(s)\|_{W^{1,\infty}_y}^2 \Bigr) \, \ud s \\
        &\lesssim \int_0^t s \cdot \varepsilon \js^\delta \cdot \varepsilon^4 \js^{-2} \, \ud s \lesssim \varepsilon \jt^\delta.
    \end{aligned}
\end{equation*}
This estimate suffices for the contributions to both \eqref{equ:pxi_calR_duhamel_bound} and \eqref{equ:pxi_with_jxi_calR_duhamel_bound}.

\medskip 
\noindent \underline{Contributions of $\calE_j(s)$, $1 \leq j \leq 3$.}
Here we can again simultaneously deal with the contributions to \eqref{equ:pxi_calR_duhamel_bound} and to \eqref{equ:pxi_with_jxi_calR_duhamel_bound}.
We begin with the term $\calE_1(s)$.
Using \eqref{equ:ILED_moving_center_frequency_side_with_jxi2} and the bounds \eqref{equ:prop_profile_bounds_assumption1}, \eqref{equ:consequences_local_decay_usubeonetwo}, we conclude 
\begin{equation*}
    \begin{aligned}
        &\biggl\| \int_0^t e^{-i\xi\theta(s)} e^{-is(\jxi+\ulell\xi)} \cdot s \cdot \bigl(\xi+\ulell\jxi\bigr) \jxi^{-1}  \calFulellsh\bigl[ \bigl( \calE_1(s) \bigr)_2 \bigr](\xi) \, \ud s \biggr\|_{L^2_\xi} \\
        &\quad + \biggl\| \jxi^2 \int_0^t e^{-i\xi\theta(s)} e^{-is(\jxi+\ulell\xi)} \cdot s \cdot \bigl(\xi+\ulell\jxi\bigr) \jxi^{-1} \calFulellsh\bigl[ \bigl( \calE_1(s) \bigr)_2 \bigr](\xi) \, \ud s \biggr\|_{L^2_\xi} \\
        &\lesssim \Bigl\| \jy^2 \cdot s \cdot \bigl( \sech^2(\gamma(s) y) - \sech^2(\ulg y) \bigr) \usubeone(s,y) \Bigr\|_{L^2_s([0,t]; H^3_y)} \\ 
        &\lesssim \Bigl\| s \cdot |\ell(s) - \ulell| \cdot \bigl\| \jy^{-1} \usubeone(s) \bigr\|_{H^3_y} \Bigr\|_{L^2_s([0,t])} 
        \lesssim \Bigl\| s \cdot \varepsilon \js^{-1+\delta} \cdot \varepsilon \js^{-\frac12} \Bigr\|_{L^2_s([0,t])} \lesssim \varepsilon^2 \jt^\delta.
    \end{aligned}
\end{equation*}

Next, we observe that $\calE_2(s) = \calN(\bmu(s)) - \calN( \ulPe u(s) )$ is a linear combination of nonlinear terms, where at least one input is given by a discrete component $d_{k,\ulell}(s) \bmY_{k,\ulell}$, $0 \leq k \leq 1$, which decays faster and is spatially localized. Correspondingly, we have by  \eqref{equ:consequences_discrete_comp_decay}, \eqref{equ:consequences_local_decay_usubeonetwo} that
\begin{equation*}
    \begin{aligned}
        &\bigl\| \jy^2 \cdot s \cdot \bigl( \calE_2(s) \bigr)_2 \bigr\|_{L^2_s([0,t];H^3_y)} \\
        &\lesssim \Bigl\| s \cdot \bigl( |d_{0,\ulell}(s)| + |d_{1,\ulell}(s)| \bigr) \Bigl( |d_{0,\ulell}(s)| + |d_{1,\ulell}(s)| + \bigl\|\jy^{-1} \usubeone(s)\|_{H^3_y} \Bigr) \Bigr\|_{L^2_s([0,t])} \\
        &\lesssim \Bigl\| s \cdot \varepsilon \js^{-\frac32+\delta} \cdot \varepsilon \js^{-\frac12} \Bigr\|_{L^2_s([0,t])} \lesssim \varepsilon^2.
    \end{aligned}
\end{equation*}
Similarly, we obtain for $\calE_3(s)$ by \eqref{equ:prop_profile_bounds_assumption1}, \eqref{equ:consequences_local_decay_usubeonetwo} that
\begin{equation*}
    \begin{aligned}
        \bigl\| \jy^2 \cdot s \cdot \calE_3(s) \bigr\|_{L^2_s([0,t];H^3_y)} &\lesssim \Bigl\| \jy^2 \cdot s \cdot \bigl( \alpha(\gamma(s)) - \alpha(\ulg y) \bigr) \bigl( \usubeone(s,y) \bigr)^2 \Bigr\|_{L^2_s([0,t];H^3_y)} \\ 
        &\lesssim \bigl\| s \cdot |\ell(s) - \ulell| \cdot \bigl\| \jy^{-1} \usubeone(s) \|_{H^3_y}^2 \bigr\|_{L^2_s([0,t])} \\ 
        &\lesssim \bigl\| s \cdot \varepsilon \js^{-1+\delta} \cdot \varepsilon^2 \js^{-1} \bigr\|_{L^2_s([0,t])} \lesssim \varepsilon^3.
    \end{aligned}
\end{equation*}
Thus, the contributions of $\calE_2(s)$ and $\calE_3(s)$ to both \eqref{equ:pxi_calR_duhamel_bound} and \eqref{equ:pxi_with_jxi_calR_duhamel_bound} can be estimated analogously to the contribution of $\calE_1(s)$ using the dual integrated local energy decay estimate \eqref{equ:ILED_moving_center_frequency_side_with_jxi2}.

\medskip 
\noindent \underline{Contribution of $\calR_q(s,\xi)$.}
In view of Remark~\ref{rem:setting_up_rapid_decay_qjs} the coefficients $( \jxi + \ulell \xi - 2 \ulg^{-1} )^{-1} \frakq_{1,\ulell}(\xi)$, $( \jxi + \ulell \xi )^{-1} \frakq_{2,\ulell}(\xi)$, and $( \jxi + \ulell \xi + 2 \ulg^{-1} )^{-1} \frakq_{3,\ulell}(\xi)$ in the definition~\eqref{equ:setting_up_definition_calRq} of $\calR_q(s,\xi)$ can all be written as $(|\ulg(\xi+\ulell\jxi)|+i)^{-1} \widehat{\varphi}_k(\xi)$, $1 \leq k \leq 3$, where $\widehat{\varphi}_k$ is the (flat) Fourier transform of some Schwartz function $\varphi_k \in \calS(\bbR)$.
Using \eqref{equ:ILED_moving_center_frequency_side_with_jxi2} and the decay estimates \eqref{equ:consequences_hulell_decay}, \eqref{equ:consequences_hulell_phase_filtered_decay}, we thus conclude 
\begin{equation*}
    \begin{aligned}
        &\biggl\| \int_0^t e^{-i\xi\theta(s)} e^{-is(\jxi+\ulell\xi)} \cdot s \cdot \bigl(\xi+\ulell\jxi\bigr) \jxi^{-1}  \calR_q(s,\xi) \, \ud s \biggr\|_{L^2_\xi} \\
        &\lesssim \biggl\| \jxi^2 \int_0^t e^{-i\xi\theta(s)} e^{-is(\jxi+\ulell\xi)} \cdot s \cdot \bigl(\xi+\ulell\jxi\bigr) \jxi^{-1} \calR_q(s,\xi) \, \ud s \biggr\|_{L^2_\xi} \\
        &\lesssim \sum_{1\leq k \leq 3} \, \bigl\| (|\ulg(\xi+\ulell\jxi)|+i)^{-1} \bigr\|_{L^\infty_\xi} \Bigl\| \jy^2 \varphi_k(y) \ps \bigl( e^{-i s \ulg^{-1}} h_\ulell(s) \bigr) \bigl( e^{-i s \ulg^{-1}} h_\ulell(s) \bigr) \Bigr\|_{L^2_s([0,t]; H^3_y)} \\ 
        &\lesssim \sum_{1\leq k \leq 3} \, \bigl\| \jy^2 \varphi_k(y) \bigr\|_{H^3_y} \bigl\| \varepsilon \js^{-1+\delta} \cdot \varepsilon \js^{-\frac12} \bigr\|_{L^2_s([0,t])} \lesssim \varepsilon^2 \jt^{\delta}.
    \end{aligned}
\end{equation*}
This bound suffices for the contributions to both \eqref{equ:pxi_calR_duhamel_bound} and \eqref{equ:pxi_with_jxi_calR_duhamel_bound}.

\medskip 
\noindent \underline{Contribution of $\calQ_{\ulell,r}\bigl(\usubeone(s)\bigr)$.}
Here we have to take some care care to separately obtain sufficient bounds for the contributions of $\calQ_{\ulell,r}\bigl(\usubeone(s)\bigr)$ defined in \eqref{equ:setting_up_definition_calQr} to \eqref{equ:pxi_calR_duhamel_bound} and to \eqref{equ:pxi_with_jxi_calR_duhamel_bound}.
Using the bounds \eqref{equ:consequences_local_decay_usubeonetwo}, \eqref{equ:consequences_hulell_decay}, \eqref{equ:consequences_Remusubeone_H1y_local_decay}, we obtain for $0 \leq s \leq t$,
\begin{equation*}
    \begin{aligned}
        \bigl\| \jy^2 \calQ_{\ulell,r}\bigl(\usubeone(s)\bigr) \bigr\|_{H^1_y}  
        &\lesssim \biggl\| \jy^2 \alpha(\ulg y) \biggl( \bigl( \usubeone(s,y) \bigr)^2 - \Bigl( \Re \, \bigl( e_\ulell^{\#}\bigl(y, - \ulg \ulell\bigr) h_\ulell(s) \bigr) \Bigr)^2 \biggr) \biggr\|_{H^1_y} \\
        &\lesssim \bigl\| \jy^{-3} \Remusubeone(s,y) \bigr\|_{H^1_y} \Bigl( \bigl\| \jy^{-1} \usubeone(s,y) \bigr\|_{H^1_y} + |h_\ulell(s)| \Bigr) \\
        &\lesssim \varepsilon \js^{-1+\delta} \cdot \varepsilon \js^{-\frac12} \lesssim \varepsilon^2 \js^{-\frac32+\delta}.        
    \end{aligned}
\end{equation*}
Thus, invoking the dual integrated local energy decay estimate \eqref{equ:ILED_moving_center_frequency_side} we obtain 
\begin{equation*}
    \begin{aligned}
        &\biggl\| \int_0^t e^{-i\xi\theta(s)} e^{-is(\jxi+\ulell\xi)} \cdot s \cdot \bigl(\xi+\ulell\jxi\bigr) \jxi^{-1}  \calFulellsh\bigl[ \calQ_{\ulell,r}\bigl(\usubeone(s)\bigr) \bigr](\xi) \, \ud s \biggr\|_{L^2_\xi} \\
        &\lesssim \Bigl\| \jy^2 \calQ_{\ulell,r}\bigl(\usubeone(s)\bigr) \Bigr\|_{L^2_s([0,t];H^1_y)} \lesssim \Bigl\| s \cdot \varepsilon^2 \js^{-\frac32+\delta} \Bigr\|_{L^2_s([0,t])} \lesssim \varepsilon^2 \jt^\delta,
    \end{aligned}
\end{equation*}
as desired for the contribution to \eqref{equ:pxi_calR_duhamel_bound}.
   
Instead at higher regularity, using \eqref{equ:consequences_local_decay_usubeonetwo}, \eqref{equ:consequences_hulell_decay}, \eqref{equ:consequences_Remusubeone_H3y_local_decay}, we find for $0 \leq s \leq t$ that 
\begin{equation*}
    \begin{aligned}
        \bigl\| \jy^2 \calQ_{\ulell,r}\bigl(\usubeone(s)\bigr) \bigr\|_{H^3_y}  
        &\lesssim \biggl\| \jy^2 \alpha(\ulg y) \biggl( \bigl( \usubeone(s,y) \bigr)^2 - \Bigl( \Re \, \bigl( e_\ulell^{\#}\bigl(y, - \ulg \ulell\bigr) h_\ulell(s) \bigr) \Bigr)^2 \biggr) \biggr\|_{H^3_y} \\
        &\lesssim \bigl\| \jy^{-3} \Remusubeone(s,y) \bigr\|_{H^3_y} \Bigl( \bigl\| \jy^{-1} \usubeone(s,y) \bigr\|_{H^3_y} + |h_\ulell(s)| \Bigr) \\
        &\lesssim \varepsilon \js^{-1+2\delta} \cdot \varepsilon \js^{-\frac12} \lesssim \varepsilon^2 \js^{-\frac32+2\delta}.        
    \end{aligned}
\end{equation*}
By the dual integrated local energy decay estimate \eqref{equ:ILED_moving_center_frequency_side_with_jxi2} we then deduce 
\begin{equation*}
    \begin{aligned}
        &\biggl\| \jxi^2 \int_0^t e^{-i\xi\theta(s)} e^{-is(\jxi+\ulell\xi)} \cdot s \cdot \bigl(\xi+\ulell\jxi\bigr) \jxi^{-1}  \calFulellsh\bigl[ \calQ_{\ulell,r}\bigl(\usubeone(s)\bigr) \bigr](\xi) \, \ud s \biggr\|_{L^2_\xi} \\
        &\lesssim \Bigl\| \jy^2 \calQ_{\ulell,r}\bigl(\usubeone(s)\bigr) \Bigr\|_{L^2_s([0,t];H^3_y)} \lesssim \Bigl\| s \cdot \varepsilon^2 \js^{-\frac32+2\delta} \Bigr\|_{L^2_s([0,t])} \lesssim \varepsilon^2 \jt^{2\delta}.
    \end{aligned}
\end{equation*}
This bound is sufficient for the contribution to \eqref{equ:pxi_with_jxi_calR_duhamel_bound}.
\end{proof}

\section{Pointwise Estimates for the Profile} \label{sec:pointwise_profile_bounds}

In this section we establish uniform-in-time pointwise bounds for the effective profile.
The main outcome is summarized in the following proposition.

\begin{proposition} \label{prop:pointwise_profile_bounds}
    Suppose the assumptions in the statement of Proposition~\ref{prop:profile_bounds} are in place.
    Assume $T \geq 1$. Then we have for all $1 \leq t \leq T$ that
    \begin{equation} \label{equ:pointwise_proposition_asserted_bound}
        \bigl\|\jxi^{\frac32} \gulellsh(t,\xi)\bigr\|_{L^\infty_\xi} \lesssim \bigl\|\jxi^{\frac32} \gulellsh(1,\xi)\bigr\|_{L^\infty_\xi} + \varepsilon^2.
    \end{equation}
    Moreover, we have for all $1 \leq t_1 \leq t_2 \leq T$ that 
    \begin{equation} \label{equ:pointwise_proposition_asserted_difference_bound}      
        \begin{aligned}
            \Bigl\| e^{-i\Lambda(t_2,\xi)} e^{-i\xi\theta(t_2)} \jxi^{\frac32} \gulellsh(t_2,\xi) - e^{-i\Lambda(t_1,\xi)} e^{-i\xi\theta(t_1)} \jxi^{\frac32} \gulellsh(t_1,\xi) \Bigr\|_{L^\infty_\xi} \lesssim \varepsilon^2 t_1^{-2\delta},
        \end{aligned}
    \end{equation}
    where 
    \begin{equation*}
        \Lambda(t,\xi) := \frac{1}{16} \jxi^{-3} \int_1^t \frac{1}{s} \bigl| \jxi^{\frac32} \gulellsh(s,\xi) \bigr|^2 \, \ud s.
    \end{equation*}
\end{proposition}

In the proof of Proposition~\ref{prop:pointwise_profile_bounds} we view the renormalized evolution equation for the effective profile \eqref{equ:setting_up_g_evol_equ4} as an ordinary differential equation in time for every fixed frequency. 
Multiplying \eqref{equ:setting_up_g_evol_equ4} by the frequency weight $\jxi^{\frac32}$, we obtain 
\begin{equation} \label{equ:pointwise_evol_equ_g_multiplied_by_jxi}
    \begin{aligned}
        &\pt \biggl( e^{-i\xi\theta(t)} \Bigl( \jxi^{\frac32} g_\ulell^{\#}(t,\xi) + \jxi^{\frac32} B\bigl[ \gulellsh \bigr](t,\xi) \Bigr) \biggr) \\
        &\quad = - \frac16 e^{-i\xi\theta(t)} e^{-i t (\jxi + \ulell \xi)} \jxi^{\frac32} \calFulellsh\bigl[ \usubeone(t)^3 \bigr](\xi) \\
         &\quad \quad + e^{-i\xi\theta(t)} e^{-i t (\jxi + \ulell \xi)} \jxi^{\frac32} \calR(t,\xi) \\
        &\quad \quad + e^{-i\xi\theta(t)} i\xi (\dot{q}(t) - \ulell) \jxi^{\frac32} B\bigl[ \gulellsh \bigr](t,\xi),
    \end{aligned}
\end{equation}
where we recall that $\calR(t,\xi)$ and $B[\gulellsh](t,\xi)$ were defined in \eqref{equ:setting_up_definition_calR}, respectively in \eqref{equ:setting_up_definition_B}.
Most of the work now goes into extracting the critically decaying leading order behavior of the cubic interactions in the first term on the right-hand side of \eqref{equ:pointwise_evol_equ_g_multiplied_by_jxi} via a stationary phase analysis. The resonant part of these critically decaying contributions leads to logarithmic phase corrections in the asymptotics of the radiation term.
The second and the third term on the right-hand side of \eqref{equ:pointwise_evol_equ_g_multiplied_by_jxi} turn out to have integrable time decay. 

In the next Subsection~\ref{subsec:stationary_phase_cubic_terms} we carry out the stationary phase analysis of the contributions of the (non-localized) cubic interactions. 
Subsequently, in Subsection~\ref{subsec:auxiliary_pointwise_bounds} we derive auxiliary pointwise bounds for the second and the third term on the right-hand side of \eqref{equ:pointwise_evol_equ_g_multiplied_by_jxi}. 
Finally, in Subsection~\ref{subsec:proof_pointwise_profile_bounds} we present the proof of Proposition~\ref{prop:pointwise_profile_bounds}.

\subsection{Stationary phase analysis of the non-localized cubic interactions} \label{subsec:stationary_phase_cubic_terms}

Building on the structure of the cubic interactions determined in Subsection~\ref{subsec:structure_cubic_nonlinearities}, we now extract the critically decaying leading order behavior of the first term on the right-hand side of \eqref{equ:pointwise_evol_equ_g_multiplied_by_jxi} via a stationary phase analysis. 

\begin{lemma} \label{lem:pointwise_cubic_terms_asymptotics}
    Suppose the assumptions in the statement of Proposition~\ref{prop:profile_bounds} are in place.
    Assume $T \geq 1$. Then we have for all $1 \leq t \leq T$ that
    \begin{equation}
        \begin{aligned}
            &\biggl\| e^{-i\xi\theta(t)} e^{-i t (\jxi + \ulell \xi)} \jxi^{\frac32} \calFulellsh\bigl[ \usubeone(t)^3 \bigr](\xi) \\ 
            &\qquad \quad - \frac{3i}{8} \frac{1}{t}  \jxi^{\frac32} \bigl|\gulellsh(t,\xi)\bigr|^2 e^{-i\xi\theta(t)} \gulellsh(t,\xi) - e^{-i\xi\theta(t)} \calT_{osc}(t,\xi) \biggr\|_{L^\infty_\xi} \lesssim \varepsilon^3 t^{-1-2\delta}
        \end{aligned}
    \end{equation}
    with 
    \begin{equation} \label{equ:pointwise_Tosc_definition}
        \begin{aligned}
            \calT_{osc}(t,\xi) &= -\frac{i}{8} \Bigl( \calI_{1,asympt}^{\delta_0}(t,\xi) + \calI_{1,asympt}^{\pvdots}(t,\xi) \Bigr) -\frac{3i}{8} \Bigl( \calI_{3,asympt}^{\delta_0}(t,\xi) + \calI_{3,asympt}^{\pvdots}(t,\xi) \Bigr) \\
            &\qquad +\frac{i}{8} \Bigl( \calI_{4,asympt}^{\delta_0}(t,\xi) + \calI_{4,asympt}^{\pvdots}(t,\xi) \Bigr),
        \end{aligned}
    \end{equation}
    where 
    \begin{align*}
        \calI_{1,asympt}^{\delta_0}(t,\xi) &:= \frac{2\pi}{t} \frac{1}{\sqrt{3}} e^{i\frac{\pi}{2}} e^{it(-\jxi+3\jap{\frac{\xi}{3}})}  \frakm_{\ulell,+++}^{\delta_0}\bigl(\xi, \tstyfrakxithree, \tstyfrakxithree, \tstyfrakxithree \bigr) \jxi^{\frac32} \Bigl( \gulellsh\bigl(t, \tstyfrakxithree \bigr) \Bigr)^3, \\
        \calI_{3,asympt}^{\delta_0}(t,\xi) &:= \frac{2\pi}{t} e^{it(-2\jxi)} \frakm_{\ulell,+--}^{\delta_0}\bigl(\xi,-\xi,-\xi,-\xi\bigr) \jxi^{\frac32} \bigl| \gulellsh(t,-\xi) \bigr|^2 \overline{\gulellsh(t,-\xi)}, \\
        \calI_{4,asympt}^{\delta_0}(t,\xi) &:= \frac{2\pi}{t} \frac{1}{\sqrt{3}} e^{-i\frac{\pi}{2}} e^{it(-\jxi-3\jap{\frac{\xi}{3}})} \frakm_{\ulell,---}^{\delta_0}\bigl(\xi, -\tstyfrakxithree, -\tstyfrakxithree, -\tstyfrakxithree \bigr) \jxi^{\frac32} \Bigl( \overline{\gulellsh\bigl(t, -\tstyfrakxithree \bigr)} \Bigr)^3,
    \end{align*}
    and
    \begin{align*}
            \calI^{\pvdots}_{1, asympt}(t,\xi) &:= - \frac{4\pi i \ulg}{t} \frac{1}{\sqrt{3}} e^{i\frac{\pi}{2}} e^{it(-\jxi+3\jap{\frac{\xi}{3}})} \frakm_{\ulell,+++}^{\pvdots}\bigl(\xi, \tstyfrakxithree, \tstyfrakxithree, \tstyfrakxithree \bigr) \jxi^{\frac32} \Bigl( \gulellsh\bigl(t, \tstyfrakxithree \bigr) \Bigr)^3 \qquad \qquad \\ 
            &\qquad \qquad \qquad \qquad \qquad \qquad \qquad \times \tanh\bigl( \ulg (\partial_\xi \phi_1)(\xi) t \bigr) \chi_1\bigl( (\tstyfrakxithree + \ulg\ulell) t^{5\delta} \bigr), \\
            \calI^{\pvdots}_{3, asympt}(t,\xi) &:= - \frac{4\pi i \ulg}{t} e^{it(-2\jxi))} \frakm_{\ulell,+--}^{\pvdots}\bigl(\xi, -\xi, -\xi, -\xi \bigr) \jxi^{\frac32} \bigl| \gulellsh(t,-\xi) \bigr|^2 \overline{\gulellsh(t,-\xi)} \\ 
            &\qquad \qquad \qquad \qquad \qquad \qquad \qquad \times \tanh\bigl( \ulg (\partial_\xi \phi_3)(\xi) t \bigr) \chi_1\bigl( (-\xi+ \ulg\ulell) t^{5\delta} \bigr), \\      
            \calI^{\pvdots}_{4, asympt}(t,\xi) &:= - \frac{4\pi i \ulg}{t} \frac{1}{\sqrt{3}} e^{-i\frac{\pi}{2}} e^{it(-\jxi-3\jap{\frac{\xi}{3}})} \frakm_{\ulell,---}^{\pvdots}\bigl(\xi, -\tstyfrakxithree, -\tstyfrakxithree, -\tstyfrakxithree \bigr) \jxi^{\frac32} \Bigl( \overline{\gulellsh\bigl(t, -\tstyfrakxithree \bigr)} \Bigr)^3 \\ 
            &\qquad \qquad \qquad \qquad \qquad \qquad \qquad \times \tanh\bigl( \ulg (\partial_\xi \phi_4)(\xi) t \bigr) \chi_1\bigl( (-\tstyfrakxithree + \ulg\ulell) t^{5\delta} \bigr).
    \end{align*}
    In the preceding list, we set $\chi_1(\cdot) := 1 - \chi_0(\cdot)$, where $\chi_0(\cdot)$ is a smooth even cut-off supported on $[-1,1]$. Moreover, the phases $\phi_j(\xi)$, $j = 1, 3, 4$, are given by
    \begin{equation*}
        \begin{aligned}
            \phi_1(\xi) := 3 \jap{ \tstyfracthird \xi } + \ulell \xi, \quad  
            \phi_3(\xi) := -\jap{\xi} + \ulell \xi, \quad
            \phi_4(\xi) := -3 \jap{ \tstyfracthird \xi } + \ulell \xi.
        \end{aligned}
    \end{equation*}
\end{lemma}

Our starting point for the proof of Lemma~\ref{lem:pointwise_cubic_terms_asymptotics} is the expansion~\eqref{equ:setting_up_dist_FT_of_cubic_expanded} of the distorted Fourier transform of the non-localized cubic interactions $\usubeone(t)^3$. It gives
\begin{equation}
    \begin{aligned}
        e^{-i t (\jxi + \ulell \xi)} \jxi^{\frac32} \calFulellsh\bigl[ \usubeone(t)^3 \bigr](\xi) &= - \frac{i}{8} \jxi^{\frac32} \calI_1(t,\xi) + \frac{3i}{8} \jxi^{\frac32} \calI_2(t,\xi) - \frac{3i}{8} \jxi^{\frac32} \calI_3(t,\xi) + \frac{i}{8} \jxi^{\frac32} \calI_4(t,\xi),
    \end{aligned}
\end{equation}
where the cubic interaction terms $\calI_j(t,\xi)$, $1 \leq j \leq 4$, are defined in \eqref{equ:definition_calI}.
Moreover, we recall from \eqref{equ:calIj_decomposition} their decompositions
\begin{equation*}
    \calI_j(t,\xi) = \calI_j^{\delta_0}(t,\xi) + \calI_j^{\pvdots}(t,\xi) + \calI_j^{\reg}(t,\xi), \quad 1 \leq j \leq 4,
\end{equation*}
where the three terms on the right-hand side are trilinear expressions in terms of the effective profile $\gulellsh(t,\xi)$ involving a Dirac kernel, a Hilbert-type kernel, and a regular kernel, respectively.
We analyze their contributions separately in the next three subsections.

\subsubsection{Contributions of cubic interactions with a Dirac kernel} \label{subsubsec:stationary_phase_Dirac}

In the next lemma we determine the critically decaying leading order behavior of all cubic interactions with a Dirac kernel.

\begin{lemma} \label{lem:cubic_dirac_stationary_phase}
    Suppose the assumptions in the statement of Proposition~\ref{prop:profile_bounds} are in place.
    Assume $T \geq 1$.
    Then we have for $1 \leq j \leq 4$ and for all $1 \leq t \leq T$ that
    \begin{equation*}
        \begin{aligned}
        \bigl\| \jxi^{\frac32} \calI_j^{\delta_0}(t,\xi) - \calI_{j,asympt}^{\delta_0}(t,\xi) \bigr\|_{L^\infty_\xi} &\lesssim \varepsilon^3 t^{-\frac98 + 6\delta},
        \end{aligned}
    \end{equation*}
    where 
    \begin{equation} \label{equ:pointwise_calIasympt_delta_definitions}
    \begin{aligned}
        \calI_{1,asympt}^{\delta_0}(t,\xi) &:= \frac{2\pi}{t} \frac{1}{\sqrt{3}} e^{i\frac{\pi}{2}} e^{it(-\jxi+3\jap{\frac{\xi}{3}})}  \frakm_{\ulell,+++}^{\delta_0}\bigl(\xi, \tstyfrakxithree, \tstyfrakxithree, \tstyfrakxithree \bigr) \jxi^{\frac32} \Bigl( \gulellsh\bigl(t, \tstyfrakxithree \bigr) \Bigr)^3, \\
        \calI_{2,asympt}^{\delta_0}(t,\xi) &:= \frac{1}{t} \jxi^{\frac32} \bigl|\gulellsh(t,\xi)\bigr|^2 \gulellsh(t,\xi), \\ 
        \calI_{3,asympt}^{\delta_0}(t,\xi) &:= \frac{2\pi}{t} e^{it(-2\jxi)} \frakm_{\ulell,+--}^{\delta_0}\bigl(\xi,-\xi,-\xi,-\xi\bigr) \jxi^{\frac32} \bigl| \gulellsh(t,-\xi) \bigr|^2 \overline{\gulellsh(t,-\xi)}, \\ 
        \calI_{4,asympt}^{\delta_0}(t,\xi) &:= \frac{2\pi}{t} \frac{1}{\sqrt{3}} e^{-i\frac{\pi}{2}} e^{it(-\jxi-3\jap{\frac{\xi}{3}})} \frakm_{\ulell,---}^{\delta_0}\bigl(\xi, -\tstyfrakxithree, -\tstyfrakxithree, -\tstyfrakxithree \bigr) \jxi^{\frac32} \Bigl( \overline{\gulellsh\bigl(t, -\tstyfrakxithree \bigr)} \Bigr)^3.
    \end{aligned}
    \end{equation}
\end{lemma}

For the proof of Lemma~\ref{lem:cubic_dirac_stationary_phase} we recall the list~\eqref{equ:structure_cubic_nonlinearities_dirac_kernel} of cubic interaction terms with a Dirac kernel
\begin{equation} 
    \begin{aligned}
        \calI_1^{\delta_0}(t,\xi) &= \iint e^{it \Psi_{1}(\xi,\xi_1,\xi_2)} \gulellsh(t,\xi_1) \, \gulellsh(t,\xi_2) \, \gulellsh(t,\xi_3) \, \jxione^{-1} \jxitwo^{-1} \jxithree^{-1} \\ 
        &\qquad \qquad \qquad \qquad \qquad \qquad \qquad \qquad \qquad \times \frakm_{\ulell, +++}^{\delta_0}(\xi,\xi_1,\xi_2,\xi_3) \, \ud \xi_1 \, \ud \xi_2 \\ 
        \calI_2^{\delta_0}(t,\xi) &= \iint e^{it \Psi_{2}(\xi,\xi_1,\xi_2)} \gulellsh(t,\xi_1) \, \overline{\gulellsh(t,\xi_2)} \, \gulellsh(t,\xi_3) \, \jxione^{-1} \jxitwo^{-1} \jxithree^{-1} \\ 
        &\qquad \qquad \qquad \qquad \qquad \qquad \qquad \qquad \qquad \times \frakm_{\ulell, +-+}^{\delta_0}(\xi,\xi_1,\xi_2,\xi_3) \, \ud \xi_1 \, \ud \xi_2 \\ 
        \calI_3^{\delta_0}(t,\xi) &= \iint e^{it \Psi_{3}(\xi,\xi_1,\xi_2)} \gulellsh(t,\xi_1) \, \overline{\gulellsh(t,\xi_2)} \, \overline{\gulellsh(t,\xi_3)} \, \jxione^{-1} \jxitwo^{-1} \jxithree^{-1} \\ 
        &\qquad \qquad \qquad \qquad \qquad \qquad \qquad \qquad \qquad \times \frakm_{\ulell, +--}^{\delta_0}(\xi,\xi_1,\xi_2,\xi_3) \, \ud \xi_1 \, \ud \xi_2 \\         
        \calI_4^{\delta_0}(t,\xi) &= \iint e^{it \Psi_{4}(\xi,\xi_1,\xi_2)} \overline{\gulellsh(t,\xi_1)} \, \overline{\gulellsh(t,\xi_2)} \, \overline{\gulellsh(t,\xi_3)} \, \jxione^{-1} \jxitwo^{-1} \jxithree^{-1} \\ 
        &\qquad \qquad \qquad \qquad \qquad \qquad \qquad \qquad \qquad \times \frakm_{\ulell, ---}^{\delta_0}(\xi,\xi_1,\xi_2,\xi_3) \, \ud \xi_1 \, \ud \xi_2,                 
    \end{aligned}
\end{equation}
where
\begin{equation*}
    \xi_3 := \left\{ \begin{aligned}
                        &\xi-\xi_1-\xi_2, \quad \quad \, j = 1, \\
                        &\xi-\xi_1+\xi_2, \quad \quad \, j = 2, \\
                        &-\xi+\xi_1-\xi_2, \quad j = 3, \\
                        &-\xi-\xi_1-\xi_2, \quad j = 4,
                     \end{aligned} \right.
\end{equation*}
and the phase functions are given by
\begin{equation}
    \begin{aligned}
        \Psi_{1}(\xi,\xi_1,\xi_2) &:= -\jxi + \jxione + \jxitwo + \jap{\xi-\xi_1-\xi_2}, \\
        \Psi_{2}(\xi,\xi_1,\xi_2) &:= -\jxi + \jxione - \jxitwo + \jap{\xi-\xi_1+\xi_2}, \\
        \Psi_{3}(\xi,\xi_1,\xi_2) &:= -\jxi + \jxione - \jxitwo - \jap{\xi-\xi_1+\xi_2}, \\
        \Psi_{4}(\xi,\xi_1,\xi_2) &:= -\jxi - \jxione - \jxitwo - \jap{\xi+\xi_1+\xi_2}.
    \end{aligned}
\end{equation}

For each fixed frequency $\xi \in \bbR$, the phases $\Psi_j(\xi,\cdot,\cdot)$, $1 \leq j \leq 4$, of the oscillatory integrals $\calI_j^{\delta_0}(t,\xi)$ have a unique non-degenerate critical point. 
Heuristically, the leading order behavior as $t \to \infty$ of each oscillatory integral $\calI_j^{\delta_0}(t,\xi)$, $1 \leq j \leq 4$, listed in \eqref{equ:pointwise_calIasympt_delta_definitions} in the statement of Lemma~\ref{lem:cubic_dirac_stationary_phase} can be determined by a formal application of the stationary phase formula 
\begin{equation} \label{equ:pointwise_stationary_phase_formula}
    \int_{\bbR^n} e^{i\lambda \Psi(x)} a(x) \, \ud x = \biggl( \frac{2\pi}{\lambda} \biggr)^{\frac{n}{2}} \frac{e^{i \lambda \Psi(x_0)} e^{i \frac{\pi}{4} \mathrm{sign} \, \mathrm{Hess} \, \Psi(x_0)}}{\sqrt{\bigl| \det \mathrm{Hess} \, \Psi(x_0)\bigr|}} a(x_0) + \ldots
\end{equation}
Indeed, by direct computation we find for every fixed frequency $\xi \in \bbR$ that
\begin{equation*}
    \begin{aligned}
    \nabla_{\xi_1,\xi_2} \Psi_1(\xi,\xi_1,\xi_2) = 0 \quad &\Leftrightarrow \quad (\xi_1, \xi_2) = \Bigl(\frac{\xi}{3}, \frac{\xi}{3}\Bigr), \\    
    \nabla_{\xi_1,\xi_2} \Psi_2(\xi,\xi_1,\xi_2) = 0 \quad &\Leftrightarrow \quad (\xi_1, \xi_2) = (\xi,\xi), \\
    \nabla_{\xi_1,\xi_2} \Psi_3(\xi,\xi_1,\xi_2) = 0 \quad &\Leftrightarrow \quad (\xi_1, \xi_2) = (-\xi, -\xi), \\ 
    \nabla_{\xi_1,\xi_2} \Psi_4(\xi,\xi_1,\xi_2) = 0 \quad &\Leftrightarrow \quad (\xi_1, \xi_2) = \Bigl(-\frac{\xi}{3}, -\frac{\xi}{3}\Bigr).
    \end{aligned}
\end{equation*}
The Hessians of the phases at the critical points are given by
\begin{equation} \label{equ:pointwise_hessians_delta_phases}
    \begin{aligned}
        \mathrm{Hess}_{\xi_1,\xi_2} \Psi_1\Bigl( \xi, \frac{\xi}{3}, \frac{\xi}{3} \Bigr) &= \jap{ {\textstyle \frac{\xi}{3}} }^{-3} \begin{bmatrix} 2 & 1 \\ 1 & 2 \end{bmatrix}, \\
        \mathrm{Hess}_{\xi_1,\xi_2} \Psi_2\bigl( \xi, \xi, \xi \bigr) &= \jxi^{-3} \begin{bmatrix} 2 & -1 \\ -1 & 0 \end{bmatrix}, \\ 
        \mathrm{Hess}_{\xi_1,\xi_2} \Psi_3\bigl( \xi, -\xi, -\xi \bigr) &= \jxi^{-3} \begin{bmatrix} 0 & 1 \\ 1 & -2 \end{bmatrix}, \\
        \mathrm{Hess}_{\xi_1,\xi_2} \Psi_4\Bigl( \xi, -\frac{\xi}{3}, -\frac{\xi}{3} \Bigr) &= - \jap{ {\textstyle \frac{\xi}{3}} }^{-3} \begin{bmatrix} 2 & 1 \\ 1 & 2 \end{bmatrix}.
    \end{aligned}
\end{equation}
Correspondingly, their signatures at the critical points are
\begin{equation*}
    \begin{aligned}
        \mathrm{sign} \, \mathrm{Hess}_{\xi_1,\xi_2} \Psi_1\Bigl( \xi, \frac{\xi}{3}, \frac{\xi}{3} \Bigr) &= 2, \\
        \mathrm{sign} \, \mathrm{Hess}_{\xi_1,\xi_2} \Psi_2\bigl( \xi, \xi, \xi \bigr) &= 0, \\
        \mathrm{sign} \, \mathrm{Hess}_{\xi_1,\xi_2} \Psi_3\bigl( \xi, -\xi, -\xi \bigr) &= 0, \\
        \mathrm{sign} \, \mathrm{Hess}_{\xi_1,\xi_2} \Psi_4\Bigl( \xi, -\frac{\xi}{3}, -\frac{\xi}{3} \Bigr) &= -2.
    \end{aligned}
\end{equation*}
Moreover, the values of the phases at the critical points are
\begin{equation*}
    \begin{aligned}
        \Psi_1\Bigl( \xi, \frac{\xi}{3}, \frac{\xi}{3} \Bigr) &= -\jxi + 3 \jap{ {\textstyle \frac{\xi}{3}} }, \\
        \Psi_2\bigl( \xi, \xi, \xi \bigr) &= 0, \\
        \Psi_3\bigl( \xi, -\xi, -\xi \bigr) &= - 2 \jxi, \\
        \Psi_4\Bigl( \xi, -\frac{\xi}{3}, -\frac{\xi}{3} \Bigr) &= -\jxi - 3 \jap{ {\textstyle \frac{\xi}{3}} }.
    \end{aligned}
\end{equation*}
We obtained the expression $\calI_{2,asympt}^{\delta_0}(t,\xi)$ using the additional observation from Subsection~\ref{subsec:cubic_spectral_distributions} that 
\begin{equation*}
    \frakm_{\ulell,+-+}^{\delta_0}\bigl(\xi,\xi,\xi,\xi\bigr) = \frac{1}{2\pi}.
\end{equation*}

However, the effective profiles are not smooth and fast decaying functions, instead they only enjoy the limited regularity and decay encoded in the bootstrap quantity \eqref{equ:prop_profile_bounds_assumption2}.
A simple application of the stationary phase formula \eqref{equ:pointwise_stationary_phase_formula} is therefore not possible, and we have to carry out the stationary phase analysis under the limited regularity and decay assumptions on the effective profiles.
We provide the details for the rigorous derivation of the leading order asymptotics for the (resonant) term $\jxi^{\frac32} \calI_2^{\delta_0}(t,\xi)$, and we leave the analogous details for the other (non-resonant) terms $\jxi^{\frac32} \calI_j^{\delta_0}(t,\xi)$, $j = 1, 3, 4$, to the reader. 

Recall from the analysis of the cubic spectral distributions in Subsection~\ref{subsec:cubic_spectral_distributions} that $\frakm_{\ulell,+-+}^{\delta_0}(\xi,\xi_1,\xi_2,\xi_3)$ is a linear combination of products of the form $\overline{\frakb(\xi)} \frakb_1(\xi_1) \overline{\frakb_2(\xi_2)} \frakb_3(\xi_3)$, where $\frakb(\xi)$ is given by 
\begin{equation*}
    \frac{1}{|\ulg (\xi+\ulell\jxi)| - i} \quad \text{or} \quad \frac{\ulg(\xi+\ulell\jxi)}{|\ulg (\xi+\ulell\jxi)| - i},
\end{equation*}
and analogously for $\frakb_i(\xi_i)$, $1 \leq i \leq 3$. 
It therefore suffices to analyze the leading order behavior of the oscillatory integral
\begin{equation}
    \begin{aligned}
        \calI(t,\xi) &:= \jxi^{\frac32} \overline{\frakb(\xi)} \iint e^{it\Psi_2(\xi,\xi_1,\xi_2)} \frakb_1(\xi_1) \gulellsh(t,\xi_1) \overline{\frakb_2(\xi_2) \gulellsh(t,\xi_2)} \frakb_3(\xi_3) \gulellsh(t,\xi_3) \\ 
        &\qquad \qquad \qquad \qquad \qquad \qquad \qquad \qquad \qquad \qquad \times \jxione^{-1} \jxitwo^{-1} \jxithree^{-1} \, \ud \xi_1 \, \ud \xi_2,
    \end{aligned}
\end{equation}
where we recall that $\xi_3 := \xi - \xi_1 + \xi_2$.
In what follows, we use the short-hand notation
\begin{equation*}
    h(t,\xi_i) := \frakb_i(\xi_i) \gulellsh(t,\xi_i).
\end{equation*}
Proceeding as in \cite[Subsection 5.3]{LS1}, we introduce a Littlewood-Paley decomposition
\begin{equation*}
    \begin{aligned}
        \calI(t,\xi) = \sum_{k,m \geq 0} \calI_{k,m}(t,\xi)
    \end{aligned}
\end{equation*}
with 
\begin{equation*}
    \begin{aligned}
        \calI_{k,m}(t,\xi) &:= \jxi^{\frac32} \overline{\frakb(\xi)} \iint e^{it\Psi_2} \jxione^{-1} \rho_k(\xi_1) h(t,\xi_1) \, \overline{\jxitwo^{-1} \rho_m(\xi_2) h(t,\xi_2)} \, \jxithree^{-1} h(t,\xi_3) \, \ud \xi_1 \, \ud \xi_2, 
    \end{aligned}
\end{equation*}
where $\rho_k(\xi_1)$ is supported on $\{ |\xi_1| \simeq 2^k \}$ for $k \geq 1$ and on $\{ |\xi_1| \lesssim 1 \}$ for $k=0$. 
We only consider the case of large output frequencies $|\xi| \simeq 2^j$ for some $j \gg 1$, and leave the easier case $ |\xi| \lesssim 1$ to the reader. 
Our starting point is to decompose $\calI(t,\xi)$ into high-low, low-high, high-high, and critical contributions,
\begin{equation*}
    \calI(t,\xi) = \calI_{hl}(t,\xi) + \calI_{lh}(t,\xi) + \calI_{hh}(t,\xi) + \calI_{crit}(t,\xi),
\end{equation*}
where 
\begin{equation*}
    \begin{aligned}
        \calI_{hl}(t,\xi) &:= \sum_{k \geq j+10} \sum_{0 \leq m \leq k-5} \calI_{km}(t,\xi), \\
        \calI_{lh}(t,\xi) &:= \sum_{m \geq j+10} \sum_{0 \leq k \leq m-5} \calI_{km}(t,\xi), \\
        \calI_{hh}(t,\xi) &:= \sum_{k \geq j+10} \sum_{|m-k| < 5} \calI_{km}(t,\xi), \\
        \calI_{crit}(t,\xi) &:= \sum_{0 \leq k,m < j+10} \calI_{km}(t,\xi).
    \end{aligned}
\end{equation*}
As in earlier sections, we work with the auxiliary linear Klein-Gordon evolutions given by
    \begin{align}
        v_i(t,y) &:= \widehat{\calF}^{-1}\Bigl[ e^{it(\jap{\xi_i}+\ulell\xi_i)} \jap{\xi_i}^{-1} \frakb(\xi_i) \gulellsh(t,\xi_i) \Bigr](y), \quad 1 \leq i \leq 3, \label{equ:pointwise_aux_KG_v} \\
        w_i(t,y) &:= \widehat{\calF}^{-1}\Bigl[ e^{it(\jap{\xi_i}+\ulell\xi_i)} \partial_{\xi_i}\bigl( \jap{\xi_i}^{-1} \frakb(\xi_i) \gulellsh(t,\xi_i) \bigr) \Bigr](y), \quad 1 \leq i \leq 3. \label{equ:pointwise_aux_KG_w}
    \end{align}
We note that the Klein-Gordon evolutions $v_i(t,y)$, $1 \leq i \leq 3$, satisfy all decay estimates and energy bounds listed in item (11) of the statement of Corollary~\ref{cor:consequences_bootstrap_assumptions}.
Additionally, by Plancherel's theorem and by the bootstrap assumption \eqref{equ:prop_profile_bounds_assumption2}, the evolutions $w_i(t,y)$, $1 \leq i \leq 3$, satisfy the energy bound
\begin{equation} \label{equ:pointwise_aux_KG_w_energy_bound}
    \sup_{0 \leq t \leq T} \, \jt^{-2\delta} \bigl\|w_i(t)\bigr\|_{H^3_y} \lesssim \varepsilon.
\end{equation}

The standard trilinear estimate below will be a frequent tool in what follows.
\begin{lemma} \label{lem:pointwise_standard_trilinear}
    Suppose $\frakm \colon \bbR^2 \to \bbC$ is a Schwartz function. Then we have for any exponents $p, q, r \in [1,\infty]$ with $\frac{1}{p} + \frac{1}{q} + \frac{1}{r} = 1$ that
    \begin{equation} 
        \biggl| \iint \frakm(\xi_1,\xi_2) f(\xi_1) \overline{g(\xi_2)} h(-\xi_1+\xi_2) \, \ud \xi_1 \, \ud \xi_2 \biggr| \lesssim \bigl\| \widehat{\calF}^{-1}[\frakm] \bigr\|_{L^1(\bbR^2)} \bigl\|\widehat{\calF}^{-1}[f]\bigr\|_{L^p} \bigl\|\widehat{\calF}^{-1}[g]\bigr\|_{L^q} \bigl\|\widehat{\calF}^{-1}[h]\bigr\|_{L^r}.
    \end{equation}
\end{lemma}
\begin{proof}
    By direct computation, we find 
    \begin{equation}
        \begin{aligned}
            &\iint \frakm(\xi_1,\xi_2) f(\xi_1) \overline{g(\xi_2)} h(-\xi_1+\xi_2) \, \ud \xi_1 \, \ud \xi_2 \\
            &= \frac{1}{(2\pi)^{\frac32}} \iiint \biggl( \iint \frakm(\xi_1,\xi_2) e^{ix_1\xi_1} e^{ix_2\xi_2} \, \ud \xi_1 \, \ud \xi_2 \biggr) \\ 
            &\qquad \qquad \qquad \qquad \qquad \times \widehat{\calF}^{-1}[f](z-x_1) \overline{\widehat{\calF}^{-1}[g](z+x_2)} \widehat{\calF}^{-1}[h](z) \, \ud x_1 \, \ud x_2 \, \ud z.
        \end{aligned}
    \end{equation}
    The asserted bound then follows by H\"older's inequality.
\end{proof}

We begin by showing that $\calI_{hl}(t,\xi)$, $\calI_{lh}(t,\xi)$, and $\calI_{hh}(t,\xi)$ are error terms.
\begin{lemma} \label{lem:pointwise_error_terms1}
    Suppose the assumptions in the statement of Proposition~\ref{prop:profile_bounds} are in place.
    Assume $T \geq 1$. Then we have uniformly for all $\xi \in \bbR$ and for all times $1 \leq t \leq T$ that 
    \begin{equation}
        |\calI_{hl}(t,\xi)| + |\calI_{lh}(t,\xi)| + |\calI_{hh}(t,\xi)| \lesssim \varepsilon^3 t^{-\frac32+3\delta}.
    \end{equation}
\end{lemma}
\begin{proof}
Recall that we consider output frequencies $|\xi| \simeq 2^j \gg 1$.
In the high-low regime, we integrate by parts in the frequency variable $\xi_2$. This gives
\begin{equation} \label{equ:pointwise_high_low_intbyparts}
    \begin{aligned}
        \calI_{k,m}(t,\xi) &= - \frac{\jxi^{\frac32} \overline{\frakb(\xi)}}{it} \iint e^{it\Psi_2} \pxitwo \biggl( \frac{1}{\pxitwo \Psi_2}  \jxione^{-1} \rho_k(\xi_1) h(t,\xi_1) \\ 
        &\qquad \qquad \qquad \qquad \qquad \qquad \quad \times \overline{\jxitwo^{-1} \rho_m(\xi_2) h(t,\xi_2)} \, \jxithree^{-1} h(t,\xi_3) \biggr) \, \ud \xi_1 \, \ud \xi_2.
    \end{aligned}
\end{equation}
We now want to use Lemma~\ref{lem:pointwise_standard_trilinear} to estimate \eqref{equ:pointwise_high_low_intbyparts} uniformly for all $|\xi| \simeq 2^j$.
Since in the high-low regime $|\xi_1| \simeq |\xi_3| \simeq 2^k$ and $|\xi_2| \simeq 2^m$, we deduce from 
\begin{equation}
    \frac{1}{\pxitwo \Psi_2} = \frac{\jxitwo \jxithree (\xi_2 \jxithree + \xi_3 \jxitwo)}{\xi_3^2 - \xi_2^2} 
\end{equation}
that
\begin{equation}
    \biggl| \pxione^{n_1} \pxitwo^{n_2} \frac{1}{\pxitwo \Psi_2} \biggr| \lesssim_{n_1, n_2} 2^{2m} 2^{-n_1 k - n_2 m} \quad \text{ for all } n_1, n_2 \geq 0.
\end{equation}
Hence, 
\begin{equation}
    \frakm(\xi_1,\xi_2) := \frac{1}{\pxitwo \Psi_2} \rho_k(\xi_1) \rho_m(\xi_2)
\end{equation}
satisfies 
\begin{equation} \label{equ:pointwise_trilinear_bound1}
    \bigl\| \widehat{\calF}^{-1}\bigl[ \frakm \bigr] \bigr\|_{L^1(\bbR^2)} \lesssim 2^{2m}, \quad \bigl\| \widehat{\calF}^{-1}\bigl[ \pxitwo \frakm \bigr] \bigr\|_{L^1(\bbR^2)} \lesssim 2^{m}.
\end{equation}
Using Lemma~\ref{lem:pointwise_standard_trilinear} along with the bound \eqref{equ:pointwise_trilinear_bound1}, we conclude from 
\eqref{equ:pointwise_high_low_intbyparts} using the notation \eqref{equ:pointwise_aux_KG_v}, \eqref{equ:pointwise_aux_KG_w} and the bounds \eqref{equ:consequences_aux_KG_disp_decay}, \eqref{equ:consequences_aux_KG_energy_bounds}, \eqref{equ:pointwise_aux_KG_w_energy_bound} that
\begin{equation}
    \begin{aligned}
        \bigl| \calI_{k,m}(t,\xi) \bigr| &\lesssim \frac{2^{\frac{3}{2} j}}{t} \biggl( \bigl\| \widehat{\calF}^{-1}\bigl[ \pxitwo \frakm \bigr] \bigr\|_{L^1(\bbR^2)} 2^{-3k} \|v_1(t)\|_{H^3_y} \|v_2(t)\|_{L^\infty_y} 2^{-3k} \|v_3(t)\|_{H^3_y} \\ 
        &\qquad \qquad + \bigl\| \widehat{\calF}^{-1}\bigl[ \frakm \bigr] \bigr\|_{L^1(\bbR^2)} 2^{-3k} \|v_1(t)\|_{H^3_y} 2^{-3m} \|w_2(t)\|_{H^3_y} \|v_3(t)\|_{L^\infty_y} \\ 
        &\qquad \qquad + \bigl\| \widehat{\calF}^{-1}\bigl[ \frakm \bigr] \bigr\|_{L^1(\bbR^2)} 2^{-3k} \|v_1(t)\|_{H^3_y} \|v_2(t)\|_{L^\infty_y} 2^{-3k} \|w_3(t)\|_{H^3_y} \biggr) \\        
        &\lesssim 2^{\frac{3}{2} j} \Bigl( 2^m 2^{-6k} + 2^{2m} 2^{-3k} 2^{-3m} + 2^{2m} 2^{-3k} 2^{-3k} \Bigr) \varepsilon^3 t^{-\frac32+3\delta}.
    \end{aligned}
\end{equation}
Summing over the parameters of the high-low case then gives $\bigl|\calI_{hl}(t,\xi)\bigr| \lesssim \varepsilon^3 t^{-\frac32+3\delta}$.
Analogously, we estimate the low-high term by $\bigl|\calI_{lh}(t,\xi)\bigr| \lesssim \varepsilon^3 t^{-\frac32+3\delta}$.

In the high-high regime, we infer as in \cite[(5.45)]{LS1} that
\begin{equation}
    \bigl( \pxione \Psi_2 \bigr)^2 + \bigl( \pxitwo \Psi_2 \bigr)^2 \gtrsim \jap{\xi_3}^{-4}.
\end{equation}
We set
\begin{equation} \label{equ:pointwise_definition_calL}
    \calL := \bigl( ( \pxione \Psi_2 )^2 + ( \pxitwo \Psi_2 )^2 \bigr)^{-1} \bigl( (\pxione \Psi_2) \pxione + (\pxitwo \Psi_2) \pxitwo \bigr).
\end{equation}
Using that $\calL \bigl( e^{it\Psi_2} \bigr) = it e^{it\Psi_2}$, we then find 
\begin{equation} \label{equ:pointwise_high_high_intbyparts}
    \begin{aligned}
        \calI_{k,m}(t,\xi) &= \frac{\jxi^{\frac32} \overline{\frakb(\xi)}}{it} \sum_{0 \leq n \leq k+10} \iint e^{it\Psi_2} \calL^\ast \biggl( \jxione^{-1} \rho_k(\xi_1) h(t,\xi_1) \\ 
        &\qquad \qquad \qquad \qquad \qquad \times \overline{\jxitwo^{-1} \rho_m(\xi_2) h(t,\xi_2)} \, \jxithree^{-1} \rho_n(\xi_3) h(t,\xi_3) \biggr) \, \ud \xi_1 \, \ud \xi_2 \\ 
        &=: \sum_{0 \leq n \leq k+10} \calI_{k,m,n}(t,\xi),
    \end{aligned}
\end{equation}
where in the last line we inserted another dyadic frequency decomposition relative to the frequency~$\xi_3$.
It is straightforward to check that for any integers $n_1, n_2 \geq 0$ with $n_1 + n_2 \geq 2$,
\begin{equation}
    \bigl| \pxione^{n_1} \pxitwo^{n_2} \Psi_2 \bigr| \lesssim \jxithree^{-1-n_1-n_2}.
\end{equation}
In order to bound \eqref{equ:pointwise_high_high_intbyparts} in $L^\infty_\xi$ via Lemma~\ref{lem:pointwise_standard_trilinear}, we introduce
\begin{equation} \label{equ:pointwise_definitions_frakms}
    \begin{aligned}
        \frakm_1(\xi_1,\xi_2) &:= - \bigl( ( \pxione \Psi_2 )^2 + ( \pxitwo \Psi_2 )^2 \bigr)^{-1} (\pxione \Psi_2)(\xi_1,\xi_2) \rho_k(\xi_1) \rho_m(\xi_2) \rho_n(\xi_3), \\
        \frakm_2(\xi_1,\xi_2) &:= - \bigl( ( \pxione \Psi_2 )^2 + ( \pxitwo \Psi_2 )^2 \bigr)^{-1} (\pxitwo \Psi_2)(\xi_1,\xi_2) \rho_k(\xi_1) \rho_m(\xi_2) \rho_n(\xi_3), \\   
        \frakm_3 &:= \pxione \frakm_1 + \pxitwo \frakm_2.
    \end{aligned}
\end{equation}
We claim that for $1 \leq i \leq 3$
\begin{equation} \label{equ:pointwise_trilinear_bound2}
    \bigl\| \widehat{\calF}^{-1}\bigl[ \frakm_i \bigr] \bigr\|_{L^1(\bbR^2)} \lesssim 2^{2n} \bigl( 2^k 2^m 2^{-2n} \bigr)^{\frac98}.
\end{equation}
To see this in the case of $\frakm_1(\xi_1, \xi_2)$, we observe that after a change of variables and integration by parts,
\begin{equation}
    \begin{aligned}
        &\biggl\| \iint_{\bbR^2} e^{i(x_1 \xi_1 + x_2 \xi_2)} \frakm_1(\xi_1,\xi_2) \, \ud \xi_1 \, \ud \xi_2 \biggr\|_{L^1_{x_1,x_2}(\bbR^2)} \\ 
        &= \biggl\| \iint_{\bbR^2} e^{i (x_1 \eta_1 + x_2 \eta_2)} \bigl( ( \pxione \Psi_2 )^2 + ( \pxitwo \Psi_2 )^2 \bigr)^{-1} (\pxione \Psi_2)(2^k \eta_1, 2^m \eta_2) \\ 
        &\qquad \qquad \qquad \qquad \qquad \qquad \qquad \times \rho(\eta_1) \rho(\eta_2) \rho\bigl( 2^{-n} (\xi - 2^{k} \eta_1 + 2^m \eta_2) \bigr) \biggr\|_{L^1_{x_1,x_2}(\bbR^2)} \\
        &\lesssim 2^{2n} \bigl( 2^k 2^m 2^{-2n} \bigr)^{\frac98}.
    \end{aligned}
\end{equation}
The $\frac98$ power here produces integrable pointwise decay of the form $(\jap{x_1} \jap{x_2})^{-\frac98}$.
The other functions $\frakm_2$ and $\frakm_3$ satisfy the same bound.
Thus, using Lemma~\ref{lem:pointwise_standard_trilinear} with \eqref{equ:pointwise_trilinear_bound2} and the bounds \eqref{equ:consequences_aux_KG_disp_decay}, \eqref{equ:consequences_aux_KG_energy_bounds}, \eqref{equ:pointwise_aux_KG_w_energy_bound}, we conclude from \eqref{equ:pointwise_high_high_intbyparts} that
\begin{equation}
    \begin{aligned}
        &\bigl| \calI_{k,m,n}(t,\xi) \bigr| \\ 
        &\lesssim \frac{2^{\frac32 j}}{t} \bigl\| \widehat{\calF}^{-1}\bigl[ \frakm_1 \bigr] \bigr\|_{L^1(\bbR^2)} 2^{-m} \|v_2(t)\|_{W^{1,\infty}_y} 2^{-3k} 2^{-3n} \Bigl( \|w_1(t)\|_{H^3_y} \|v_3(t)\|_{H^3_y} + \|v_1(t)\|_{H^3_y} \|w_3(t)\|_{H^3_y} \Bigr) \\ 
        &\quad + \frac{2^{\frac32 j}}{t} \bigl\| \widehat{\calF}^{-1}\bigl[ \frakm_2 \bigr] \bigr\|_{L^1(\bbR^2)} 2^{-k} \|v_1(t)\|_{W^{1,\infty}_y} 2^{-3m} 2^{-3n} \Bigl( \|w_2(t)\|_{H^3_y} \|v_3(t)\|_{H^3_y} + \|v_2(t)\|_{H^3_y} \|w_3(t)\|_{H^3_y} \Bigr) \\ 
        &\quad + \frac{2^{\frac32 j}}{t} \bigl\| \widehat{\calF}^{-1}\bigl[ \frakm_3 \bigr] \bigr\|_{L^1(\bbR^2)} 2^{-3k} \|v_1(t)\|_{H^3_y} 2^{-3m} \|v_2(t)\|_{H^3_y} \|v_3(t)\|_{W^{1,\infty}_y} \\
        &\lesssim 2^{\frac32 j} 2^{\frac98 (k+m)} \Bigl( 2^{-m} 2^{-3k} 2^{-3n} + 2^{-k} 2^{-3m} 2^{-3n} + 2^{-3k} 2^{-3m} 2^{-n} \Bigr) \varepsilon^3 t^{-\frac32+3\delta}.
    \end{aligned}
\end{equation}
Summing over the high-high parameter regime yields the asserted bound $\bigl|\calI_{hh}(t,\xi)\bigr| \lesssim \varepsilon^3 t^{-\frac32+3\delta}$.
\end{proof}

Next, we consider $\calI_{crit}(t,\xi)$, which contains the critical point.
The integration region is of the form $|\xi_1| + |\xi_2| \lesssim |\xi|$ and we have $\sum_{i=1}^3 |\xi_i| \simeq |\xi|$.
The unique critical point is at $(\xi_1, \xi_2) = (\xi,\xi)$. We consider the case $|\xi| \gg 1$ and we denote
\begin{equation*}
    U_\ast := \Bigl\{ (\xi_1, \xi_2) \in \bbR^3 \, \Big| \, \max \bigl\{ |\xi_1-\xi|, |\xi_2-\xi| \bigr\} \leq c_\ast |\xi| \Bigr\},
\end{equation*}
where $0 < c_\ast \ll 1$ is a small constant that will be specified in Lemma~\ref{lem:pointwise_critical_term} below.
From \cite[Lemma 5.14]{LS1} we have the following characterization of the neighborhood $U_\ast$.

\begin{lemma} \label{lem:pointwise_Uast_characterization}
    The neighborhood $U_\ast$ is characterized by the property 
    \begin{equation*}
        |\nabla_{\xi_1, \xi_2} \Psi_2(\xi_1, \xi_2)| \ll |\xi|^{-2},
    \end{equation*}
    and in $U_\ast$ it holds that
    \begin{equation*}
        |\nabla_{\xi_1,\xi_2} \Psi_2(\xi+\eta_1,\xi+\eta_2)| \simeq |\xi|^{-3} \bigl( |\eta_1| + |\eta_2| \bigr).
    \end{equation*}
\end{lemma}

Correspondingly, we decompose
\begin{equation*}
    \calI_{crit}(t,\xi) = \calI_{U_\ast}(t,\xi) + \calI_{U_\ast}^c(t,\xi)
\end{equation*}
with
\begin{equation} \label{equ:pointwise_definition_IUast}
    \calI_{U_\ast}(t,\xi) := \jxi^{\frac32} \overline{\frakb(\xi)} \iint e^{it\Psi_2} \jxione^{-1} h(t,\xi_1) \, \overline{\jxitwo^{-1} h(t,\xi_2)} \, \jxithree^{-1} h(t,\xi_3) \chi_{U_\ast}(\xi_1,\xi_2) \, \ud \xi_1 \, \ud \xi_2, 
\end{equation}
where $\chi_{U_\ast}(\xi_1, \xi_2)$ is a smooth bump function adapted to $U_\ast$. 

Next, we show that $\calI_{U_\ast}^c(t,\xi)$ is also an error term. 
\begin{lemma} \label{lem:pointwise_error_terms2}
    Suppose the assumptions in the statement of Proposition~\ref{prop:profile_bounds} are in place.
    Assume $T \geq 1$. Then we have uniformly for all $\xi \in \bbR$ and all times $1 \leq t \leq T$ that
    \begin{equation}
        \bigl| \calI_{U_\ast}^c(t,\xi) \bigr| \lesssim \varepsilon^3 t^{-\frac32+3\delta}.
    \end{equation}
\end{lemma}
\begin{proof}
Again, we only consider the more difficult case of large output frequencies $|\xi| \simeq 2^j \gg 1$.
We further decompose the term $\calI_{U_\ast}^c(t,\xi)$ into high-low, low-high, and high-high contributions relative to the frequency variables $\xi_1$ and $\xi_2$,
\begin{equation}
    \calI_{U_\ast}^c(t,\xi) = \calI_{U_\ast, hl}^c(t,\xi) + \calI_{U_\ast, lh}^c(t,\xi) + \calI_{U_\ast, hh}^c(t,\xi),
\end{equation}
where 
\begin{equation}
    \begin{aligned}
        \calI_{U_\ast, hl}^c(t,\xi) &:= \sum_{0 \leq k \leq j+10} \sum_{0 \leq m \leq k-5} \calI_{U_\ast, km}^c(t,\xi), \\ 
        \calI_{U_\ast, lh}^c(t,\xi) &:= \sum_{0 \leq m \leq j+10} \sum_{0 \leq k \leq m-5} \calI_{U_\ast, km}^c(t,\xi), \\
        \calI_{U_\ast, hh}^c(t,\xi) &:= \sum_{0 \leq k \leq j+10} \sum_{|m-k|<5} \sum_{0 \leq n \leq j+10} \calI_{U_\ast, kmn}^c(t,\xi),
    \end{aligned}
\end{equation}
and where 
\begin{equation}
    \begin{aligned}
        \calI_{U_\ast, km}^c(t,\xi) &:= \jxi^{\frac32} \overline{\frakb(\xi)} \iint e^{it\Psi_2(\xi,\xi_1,\xi_2)} \jxione^{-1} \rho_k(\xi_1) h(t,\xi_1) \, \overline{\jxitwo^{-1} \rho_m(\xi_2) h(t,\xi_2)} \\
        &\qquad \qquad \qquad \qquad \qquad \qquad \qquad \quad \times \, \jxithree^{-1} h(t,\xi_3) \bigl( 1 - \chi_{U_\ast}(\xi_1,\xi_2) \bigr) \, \ud \xi_1 \, \ud \xi_2, \\
        \calI_{U_\ast, kmn}^c(t,\xi) &:= \jxi^{\frac32} \overline{\frakb(\xi)} \iint e^{it\Psi_2(\xi,\xi_1,\xi_2)} \jxione^{-1} \rho_k(\xi_1) h(t,\xi_1) \, \overline{\jxitwo^{-1} \rho_m(\xi_2) h(t,\xi_2)} \\
        &\qquad \qquad \qquad \qquad \qquad \qquad \qquad \quad \times \, \jxithree^{-1} \rho_n(\xi_3) h(t,\xi_3) \bigl( 1 - \chi_{U_\ast}(\xi_1,\xi_2) \bigr) \, \ud \xi_1 \, \ud \xi_2.
    \end{aligned}
\end{equation}
The high-low contributions $\calI_{U_\ast, hl}^c(t,\xi)$ and the low-high contributions $\calI_{U_\ast, lh}^c(t,\xi)$ can be treated in the exact same manner as the corresponding high-low and low-high terms in the proof of Lemma~\ref{lem:pointwise_error_terms1}. We omit the details. 

Instead, the high-high contributions $\calI_{U_\ast, hh}^c(t,\xi)$, now much closer to the critical point, require a more careful treatment.
Integrating by parts using $\calL$ defined in \eqref{lem:pointwise_error_terms1}, 
we find that
\begin{equation}
    \calI_{U_\ast,hh}^c(t,\xi) = \sum_{0 \leq k \leq j+10} \sum_{|m-k|<5} \sum_{0 \leq n \leq j+10} J_{kmn}(t,\xi)
\end{equation}
with 
\begin{equation} \label{equ:pointwise_definition_Jkmn}
\begin{aligned}
    J_{kmn}(t,\xi) &:= \frac{\jxi^{\frac32} \overline{\frakb(\xi)}}{it} \iint e^{it\Psi_2} \calL^\ast \biggl( \jxione^{-1} \rho_k(\xi_1) h(t,\xi_1) \overline{\jxitwo^{-1} \rho_m(\xi_2) h(t,\xi_2)} \\ 
    &\qquad \qquad \qquad \qquad \qquad \qquad \quad \times \, \jxithree^{-1} \rho_n(\xi_3) h(t,\xi_3) \biggr) \bigl( 1 - \chi_{U_\ast}(\xi_1,\xi_2) \bigr) \, \ud \xi_1 \, \ud \xi_2.
\end{aligned}
\end{equation}
In order to bound \eqref{equ:pointwise_definition_Jkmn} in $L^\infty_\xi$ using Lemma~\ref{lem:pointwise_standard_trilinear}, in analogy to \eqref{equ:pointwise_definitions_frakms} we introduce
\begin{equation} 
    \begin{aligned}
        \frakn_1(\xi_1,\xi_2) &:= - \bigl( ( \pxione \Psi_2 )^2 + ( \pxitwo \Psi_2 )^2 \bigr)^{-1} (\pxione \Psi_2)(\xi_1,\xi_2) \rho_k(\xi_1) \rho_m(\xi_2) \rho_n(\xi_3) \bigl( 1 - \chi_{U_\ast}(\xi_1,\xi_2) \bigr), \\
        \frakn_2(\xi_1,\xi_2) &:= - \bigl( ( \pxione \Psi_2 )^2 + ( \pxitwo \Psi_2 )^2 \bigr)^{-1} (\pxitwo \Psi_2)(\xi_1,\xi_2) \rho_k(\xi_1) \rho_m(\xi_2) \rho_n(\xi_3) \bigl( 1 - \chi_{U_\ast}(\xi_1,\xi_2) \bigr), \\   
        \frakn_3 &:= \pxione \frakn_1 + \pxitwo \frakn_2.
    \end{aligned}
\end{equation}
By \cite[(5.69)--(5.73)]{LS1} we have that
\begin{equation} \label{equ:pointwise_trilinear_bound3}
    \bigl\| \widehat{\calF}^{-1}\bigl[ \frakn_1 \bigr] \bigr\|_{L^1(\bbR^2)} + \bigl\| \widehat{\calF}^{-1}\bigl[ \frakn_2 \bigr] \bigr\|_{L^1(\bbR^2)} + \bigl\| \widehat{\calF}^{-1}\bigl[ \frakn_3 \bigr] \bigr\|_{L^1(\bbR^2)} \lesssim 2^{k+m}.
\end{equation}
Then using Lemma~\ref{lem:pointwise_standard_trilinear} with \eqref{equ:pointwise_trilinear_bound3} along with the bounds \eqref{equ:consequences_aux_KG_disp_decay}, \eqref{equ:consequences_aux_KG_energy_bounds}, \eqref{equ:pointwise_aux_KG_w_energy_bound}, we obtain from \eqref{equ:pointwise_definition_Jkmn} that
\begin{equation}
    \begin{aligned}
        &\bigl| J_{kmn}(t,\xi) \bigr| \\
        &\lesssim \frac{2^{\frac32 j}}{t} \bigl\| \widehat{\calF}^{-1}\bigl[ \frakn_1 \bigr] \bigr\|_{L^1(\bbR^2)} 2^{-m} \|v_2(t)\|_{W^{1,\infty}_y} 2^{-3k} 2^{-3n} \Bigl( \|w_1(t)\|_{H^3_y} \|v_3(t)\|_{H^3_y} + \|v_1(t)\|_{H^3_y} \|w_3(t)\|_{H^3_y} \Bigr) \\
        &\quad + \frac{2^{\frac32 j}}{t} \bigl\| \widehat{\calF}^{-1}\bigl[ \frakn_2 \bigr] \bigr\|_{L^1(\bbR^2)} 2^{-k} \|v_1(t)\|_{W^{1,\infty}_y} 2^{-3m} 2^{-3n} \Bigl( \|w_2(t)\|_{H^3_y} \|v_3(t)\|_{H^3_y} + \|v_2(t)\|_{H^3_y} \|w_3(t)\|_{H^3_y} \Bigr) \\
        &\quad + \frac{2^{\frac32 j}}{t} \bigl\| \widehat{\calF}^{-1}\bigl[ \frakn_3 \bigr] \bigr\|_{L^1(\bbR^2)} 2^{-k} \|v_1(t)\|_{W^{1,\infty}_y} 2^{-3m} \|v_2(t)\|_{H^3_y} 2^{-3n} \|v_3(t)\|_{H^3_y} \\
        &\lesssim 2^{\frac32 j} 2^{k+m} \Bigl( 2^{-3k} 2^{-m} 2^{-3n} + 2^{-k} 2^{-3m} 2^{-3n} \Bigr) \varepsilon^3 t^{-\frac32+3\delta}.
    \end{aligned}
\end{equation}
Keeping in mind that $\max \, \{2^k, 2^m, 2^n\} \simeq 2^j$ and summing over the frequency parameters for the high-high regime, we conclude $\bigl|\calI_{U_\ast, hh}^c(t,\xi)\bigr| \lesssim \varepsilon^3 t^{-\frac32+3\delta}$. This finishes the proof of the lemma.
\end{proof}

Finally, we extract the leading order behavior of $\calI_{U_\ast}(t,\xi)$.

\begin{lemma} \label{lem:pointwise_critical_term}
    Suppose the assumptions in the statement of Proposition~\ref{prop:profile_bounds} are in place.
    Assume $T \geq 1$. For any $0 < \alpha < \frac14$ we have uniformly for all times $1 \leq t \leq T$ that
    \begin{equation} \label{equ:pointwise_critical_term_asymptotics}
        \begin{aligned}
            \biggl\| \calI_{U_\ast}(t,\xi) - \frac{2\pi}{t} \jxi^{\frac32} \overline{\frakb(\xi)} \frakb_1(\xi) \overline{\frakb_2(\xi)} \frakb_3(\xi) \bigl| \gulellsh(t,\xi) \bigr|^2 \gulellsh(t,\xi) \biggr\|_{L^\infty_\xi} \lesssim \varepsilon^3 t^{-1-\alpha+6\delta}.
        \end{aligned}
    \end{equation}
\end{lemma}
\begin{proof}
We proceed along the lines of the proof of \cite[Lemma 5.15]{LS1}. 
Recall that we consider a large output frequency $|\xi| \simeq 2^j \gg 1$ and that $|\xi_1| \simeq |\xi_2| \simeq |\xi_3| \simeq |\xi|$ on the support of $U_\ast$.
First, we make the substitution
\begin{equation}
    \xi_1 = \xi + \jxi \zeta_1, \quad \xi_2 = \xi+\jxi \zeta_2
\end{equation}
in \eqref{equ:pointwise_definition_IUast} and we rescale the phase as
\begin{equation}
    \Psi_2(\xi,\xi_1,\xi_2) = \jxi^{-1} \Psi(\zeta_1, \zeta_2).
\end{equation}
Then we have $\partial_{\zeta_1} \Psi(0,0) = \partial_{\zeta_2} \Psi(0,0) = 0$ as well as $\Psi(0,0)=0$, and by \eqref{equ:pointwise_hessians_delta_phases} it holds that
\begin{equation}
    \mathrm{Hess} \, \Psi(0,0) = \begin{bmatrix} 2 & - 1 \\ -1 & 0 \end{bmatrix}.
\end{equation}
Furthermore, Lemma~\ref{lem:pointwise_Uast_characterization} implies that for all $|\zeta_1| + |\zeta_2| \lesssim 1$ and for all integers $n_1, n_2 \geq 0$ with $n_1 + n_2 \geq 1$,
\begin{equation}
    \bigl| \partial^{n_1}_{\zeta_1} \partial^{n_2}_{\zeta_2} \Psi(\zeta_1, \zeta_2) \bigr| \lesssim_{n_1,n_2} 1.
\end{equation}
Next, we set 
\begin{equation}
    F(\zeta_1, \zeta_2) := f_j(t,\xi_1) \, \overline{f_j(t,\xi_2)} \, f_j(t,\xi_3)
\end{equation}
with 
\begin{equation}
    f_j(t,\xi_k) := \jap{\xi_k}^{-1} \widetilde{\rho}_j(\xi_k) h(t,\xi_k), \quad 1 \leq k \leq 3,
\end{equation}
where on the support of $U_\ast$ we could freely insert a fattened cut-off  $\tilde{\rho}_j(\xi_k)$ localized around $|\xi_k| \simeq |\xi| \simeq 2^j$.
Moreover, we define $\lambda := t \jxi^{-1}$ and $\chi_{U_\ast}(\xi_1,\xi_2) := \chi_0(\zeta_1, \zeta_2)$ with $\chi_0$ being a smooth cut-off to a neighborhood of $(0,0)$ of size $0 < c_\ast \ll 1$.
Then we have
\begin{equation}
    \begin{aligned}
         \calI_{U_\ast}(t,\xi) &= \jxi^{\frac72} \overline{\frakb(\xi)} \iint e^{i\lambda \Psi(\zeta_1 \zeta_2)} \chi_0(\zeta_1, \zeta_2) F(\zeta_1, \zeta_2) \, \ud \zeta_1 \, \ud \zeta_2 \\
         &= (2\pi)^{-1} \jxi^{\frac72} \overline{\frakb(\xi)} \iint G_\lambda(z_1,z_2) \widehat{F}(z_1, z_2) \, \ud z_1 \, \ud z_2,
    \end{aligned}
\end{equation}
where we introduce the function
\begin{equation}
    G_\lambda(z_1,z_2) := \iint e^{i(z_1 \zeta_1 + z_2 \zeta_2)} e^{i\lambda\Psi(\zeta_1,\zeta_2)} \chi_0(\zeta_1, \zeta_2) \, \ud z_1 \, \ud z_2.
\end{equation}
It follows that
\begin{equation}
    \begin{aligned}
        \calI_{U_\ast}(t,\xi) &= \jxi^{\frac72} \overline{\frakb(\xi)} G_\lambda(0,0) F(0,0) + \calO \Bigl( \jxi^{\frac72} \bigl\| \bigl( G_\lambda - G_\lambda(0,0) \bigr) \widehat{F} \bigr\|_{L^1_{z_1,z_2}} \Bigr).
    \end{aligned}
\end{equation}
Carrying out exactly the same stationary phase analysis as in \cite[(5.65)--(5.67)]{LS1}, we then conclude for any $0 \leq \alpha \leq 1$ that
\begin{equation} \label{equ:pointwise_asymptotics_deriv1}
    \begin{aligned}
        \calI_{U_\ast}(t,\xi) &= \frac{2\pi}{t} \jxi^{\frac92} \overline{\frakb(\xi)} \bigl| h_j(t,\xi) \bigr|^2 h_j(t,\xi) + \calO \Bigl( \jxi^{\frac92+\alpha} \bigl\| \jap{|z_1| + |z_2|}^{2\alpha} \widehat{F} \bigr\|_{L^1_{z_1, z_2}} t^{-1-\alpha} \Bigr),
    \end{aligned}
\end{equation}
where 
\begin{equation}
    \begin{aligned}
        \widehat{F}(z_1, z_2) &= \frac{1}{2\pi} \iint e^{-i(z_1 \zeta_1 + z_2 \zeta_2)} f_j\bigl(t,\xi+\jxi\zeta_1\bigr) \overline{f_j\bigl(t,\xi+\jxi\zeta_2\bigr)} f_j\bigl(t,\xi-\jxi\zeta_1+\jxi\zeta_2\bigr) \, \ud \zeta_1 \, \ud \zeta_2 \\
        &= \frac{1}{\sqrt{2\pi}} \jxi^{-2} e^{i(z_1+z_2)\frac{\xi}{\jxi}} \int e^{-i\xi y_3} \widehat{\calF}^{-1}\bigl[f_j(t)\bigr]\bigl(y_3-\jxi^{-1}z_1\bigr) \\ 
        &\qquad \qquad \qquad \qquad \qquad \qquad \qquad \qquad \times \overline{\widehat{\calF}^{-1}\bigl[f_j(t)\bigr]\bigl(y_3+\jxi^{-1}z_2\bigr)}  \widehat{\calF}^{-1}\bigl[f_j(t)\bigr](y_3) \, \ud y_3. 
    \end{aligned}
\end{equation}
Using elementary estimates such as $\jap{z_1} \lesssim \jxi \bigl( \jap{y_3-\jxi^{-1}z_1} + \jap{y_3} \bigr)$ and keeping in mind that $|\xi| \simeq 2^j$, we infer from the preceding identity by H\"older's inequality that for any $0 < \alpha < \frac14$,
\begin{equation} \label{equ:pointwise_asymptotics_deriv2}
    \begin{aligned}
        \bigl\| \jap{ |z_1| + |z_2| }^{2\alpha} \widehat{F}(z_1,z_2) \bigr\|_{L^1_{z_1,z_2}} &\lesssim 2^{2\alpha j} \bigl\| \jy^{2\alpha} \widehat{\calF}^{-1}\bigl[f_j(t)\bigr] \bigr\|_{L^1_y}^2 \bigl\| \widehat{\calF}^{-1}\bigl[f_j(t)\bigr] \bigr\|_{L^1_y} \\ 
        &\lesssim 2^{2\alpha j} \bigl\| \jy \widehat{\calF}^{-1}\bigl[f_j(t)\bigr] \bigr\|_{L^2_y}^3.
    \end{aligned}
\end{equation}
Finally, by standard properties of the flat Fourier transform and by the bootstrap assumption \eqref{equ:prop_profile_bounds_assumption2} we have 
\begin{equation} \label{equ:pointwise_asymptotics_deriv3}
    \begin{aligned}
        \bigl\| \jy \widehat{\calF}^{-1}\bigl[f_j(t)\bigr] \bigr\|_{L^2_y} &\lesssim 2^{-3j} \Bigl( \bigl\| \jxi^{-1} \gulellsh(t,\xi) \bigr\|_{L^2_\xi} + \bigl\| \jxi^2 \pxi \gulellsh(t,\xi) \bigr\|_{L^2_\xi} \Bigr) \lesssim 2^{-3j} \varepsilon \jt^{2\delta}.
    \end{aligned}
\end{equation}
Combining \eqref{equ:pointwise_asymptotics_deriv1}, \eqref{equ:pointwise_asymptotics_deriv2}, and \eqref{equ:pointwise_asymptotics_deriv3}, we arrive at the asserted asymptotics \eqref{equ:pointwise_critical_term_asymptotics}.
\end{proof}

Using Lemma~\ref{lem:pointwise_critical_term} with $\alpha = \frac18$, summing over all terms constituting $\frakm_{\ulell,+-+}^{\delta_0}\bigl(\xi,\xi,\xi,\xi\bigr)$, and invoking \eqref{equ:cubic_spectral_distributions_pmp_delta_identity}, we arrive at the expression $\calI_{2,asympt}^{\delta_0}(t,\xi)$ defined in \eqref{equ:pointwise_calIasympt_delta_definitions} for the leading order behavior of $\jxi^{\frac32} \calI_2^{\delta_0}(t,\xi)$.
This completes the proof of Lemma~\ref{lem:cubic_dirac_stationary_phase} for the case $\jxi^{\frac32} \calI_2^{\delta_0}(t,\xi)$.

\subsubsection{Contributions of cubic interactions with a Hilbert-type kernel} \label{subsubsec:stationary_phase_Hilbert}

In the next lemma we determine the critically decaying leading order behavior of all cubic interactions with a Hilbert-type kernel.

\begin{lemma} \label{lem:cubic_Hilbert_stationary_phase}
    Suppose the assumptions in the statement of Proposition~\ref{prop:profile_bounds} are in place.
    Assume $T \geq 1$.
    Let $\chi_0(\cdot)$ be a smooth even cut-off supported on $[-1,1]$ and set $\chi_1(\cdot) := 1 - \chi_0(\cdot)$.
    Then we have for all $1 \leq t \leq T$ that
    \begin{equation} \label{equ:pointwise_calIpv_asymptotics}
    \begin{aligned}
        \biggl\| \jxi^{\frac32} \calI_1^{\pvdots}(t,\xi) - \calI^{\pvdots}_{1, asympt}(t,\xi) \biggr\|_{L^\infty_\xi} &\lesssim \varepsilon^3 t^{-1-2\delta}, \\
        \bigl\| \jxi^{\frac32} \calI_2^{\pvdots}(t,\xi) \bigr\|_{L^\infty_\xi} &\lesssim \varepsilon^3 t^{-1-2\delta}, \\
        \biggl\| \jxi^{\frac32} \calI_3^{\pvdots}(t,\xi) - \calI^{\pvdots}_{3, asympt}(t,\xi) \biggr\|_{L^\infty_\xi} &\lesssim \varepsilon^3 t^{-1-2\delta}, \\
        \biggl\| \jxi^{\frac32} \calI_4^{\pvdots}(t,\xi) - \calI^{\pvdots}_{4, asympt}(t,\xi) \biggr\|_{L^\infty_\xi} &\lesssim \varepsilon^3 t^{-1-2\delta},
    \end{aligned}
    \end{equation}
    where 
    \begin{equation} \label{equ:pointwise_calIasympt_pv_definitions}
        \begin{aligned}
            \calI^{\pvdots}_{1, asympt}(t,\xi) &:= - \frac{4\pi i \ulg}{t} \frac{1}{\sqrt{3}} e^{i\frac{\pi}{2}} e^{it(-\jxi+3\jap{\frac{\xi}{3}})} \frakm_{\ulell,+++}^{\pvdots}\bigl(\xi, \tstyfrakxithree, \tstyfrakxithree, \tstyfrakxithree \bigr) \jxi^{\frac32} \Bigl( \gulellsh\bigl(t, \tstyfrakxithree \bigr) \Bigr)^3 \qquad \qquad \\ 
            &\qquad \qquad \qquad \qquad \qquad \qquad \qquad \times \tanh\bigl( \ulg (\partial_\xi \phi_1)(\xi) t \bigr) \chi_1\bigl( (\tstyfrakxithree + \ulg\ulell) t^{5\delta} \bigr) \\
            \calI^{\pvdots}_{3, asympt}(t,\xi) &:= - \frac{4\pi i \ulg}{t} e^{it(-2\jxi))} \frakm_{\ulell,+--}^{\pvdots}\bigl(\xi, -\xi, -\xi, -\xi \bigr) \jxi^{\frac32} \bigl| \gulellsh(t,-\xi) \bigr|^2 \overline{\gulellsh(t,-\xi)} \\ 
            &\qquad \qquad \qquad \qquad \qquad \qquad \qquad \times \tanh\bigl( \ulg (\partial_\xi \phi_3)(\xi) t \bigr) \chi_1\bigl( (-\xi+ \ulg\ulell) t^{5\delta} \bigr) \\
            \calI^{\pvdots}_{4, asympt}(t,\xi) &:= - \frac{4\pi i \ulg}{t} \frac{1}{\sqrt{3}} e^{-i\frac{\pi}{2}} e^{it(-\jxi-3\jap{\frac{\xi}{3}})} \frakm_{\ulell,---}^{\pvdots}\bigl(\xi, -\tstyfrakxithree, -\tstyfrakxithree, -\tstyfrakxithree \bigr) \jxi^{\frac32} \Bigl( \overline{\gulellsh\bigl(t, -\tstyfrakxithree \bigr)} \Bigr)^3 \\ 
            &\qquad \qquad \qquad \qquad \qquad \qquad \qquad \times \tanh\bigl( \ulg (\partial_\xi \phi_4)(\xi) t \bigr) \chi_1\bigl( (-\tstyfrakxithree + \ulg\ulell) t^{5\delta} \bigr)
        \end{aligned}
    \end{equation}
    with the phases $\phi_j(\xi)$, $j = 1, 3, 4$, given by
    \begin{equation*}
        \begin{aligned}
            \phi_1(\xi) := 3 \jap{ \tstyfracthird \xi } + \ulell \xi, \quad 
            \phi_3(\xi) := -\jap{\xi} + \ulell \xi, \quad 
            \phi_4(\xi) := -3 \jap{ \tstyfracthird \xi } + \ulell \xi.
        \end{aligned}
    \end{equation*}
\end{lemma}

For the proof of Lemma~\ref{lem:cubic_Hilbert_stationary_phase} we recall the list \eqref{equ:structure_cubic_nonlinearities_PV} of cubic interaction terms with a Hilbert-type kernel,
\begin{equation} \label{equ:structure_cubic_nonlinearities_PV_with_cutoff}
    \begin{aligned}
        \calI_1^{\pvdots}(t,\xi) &= \iiint e^{it \Omega_{1,\ulell}(\xi,\xi_1,\xi_2,\xi_4)} \gulellsh(t,\xi_1) \, \gulellsh(t,\xi_2) \, \gulellsh(t,\xi_3) \, \jxione^{-1} \jxitwo^{-1} \jxithree^{-1} \\ 
        &\qquad \qquad \qquad \qquad \times \frakm_{\ulell, +++}^{\pvdots}(\xi, \xi_1, \xi_2, \xi_3) \widetilde{\chi}_0(\xi_4) \, \pvdots \cosech\Bigl( \frac{\pi}{2\ulg} \xi_4 \Bigr) \, \ud \xi_1 \, \ud \xi_2 \, \ud \xi_4 \\
        \calI_2^{\pvdots}(t,\xi) &= \iiint e^{it \Omega_{2,\ulell}(\xi,\xi_1,\xi_2,\xi_4)} \gulellsh(t,\xi_1) \, \overline{\gulellsh(t,\xi_2)} \, \gulellsh(t,\xi_3) \, \jxione^{-1} \jxitwo^{-1} \jxithree^{-1} \\ 
        &\qquad \qquad \qquad \qquad \times \frakm_{\ulell, +-+}^{\pvdots}(\xi, \xi_1, \xi_2, \xi_3) \widetilde{\chi}_0(\xi_4) \, \pvdots \cosech\Bigl( \frac{\pi}{2\ulg} \xi_4 \Bigr) \, \ud \xi_1 \, \ud \xi_2 \, \ud \xi_4 \\        
        \calI_3^{\pvdots}(t,\xi) &= \iiint e^{it \Omega_{3,\ulell}(\xi,\xi_1,\xi_2,\xi_4)} \gulellsh(t,\xi_1) \, \overline{\gulellsh(t,\xi_2)} \, \overline{\gulellsh(t,\xi_3)} \, \jxione^{-1} \jxitwo^{-1} \jxithree^{-1} \\ 
        &\qquad \qquad \qquad \qquad \times \frakm_{\ulell, +--}^{\pvdots}(\xi, \xi_1, \xi_2, \xi_3) \widetilde{\chi}_0(\xi_4) \, \pvdots \cosech\Bigl( \frac{\pi}{2\ulg} \xi_4 \Bigr) \, \ud \xi_1 \, \ud \xi_2 \, \ud \xi_4 \\                
        \calI_4^{\pvdots}(t,\xi) &= \iiint e^{it \Omega_{4,\ulell}(\xi,\xi_1,\xi_2,\xi_4)} \overline{\gulellsh(t,\xi_1)} \, \overline{\gulellsh(t,\xi_2)} \, \overline{\gulellsh(t,\xi_3)} \, \jxione^{-1} \jxitwo^{-1} \jxithree^{-1} \\ 
        &\qquad \qquad \qquad \qquad \times \frakm_{\ulell, ---}^{\pvdots}(\xi, \xi_1, \xi_2, \xi_3) \widetilde{\chi}_0(\xi_4) \, \pvdots \cosech\Bigl( \frac{\pi}{2\ulg} \xi_4 \Bigr) \, \ud \xi_1 \, \ud \xi_2 \, \ud \xi_4                       
    \end{aligned}
\end{equation}
with 
\begin{equation*}
    \xi_3 := \left\{ \begin{aligned}
                        &\xi-\xi_1-\xi_2+\xi_4, \quad \quad \, j = 1, \\
                        &\xi-\xi_1+\xi_2+\xi_4, \quad \quad \, j = 2, \\
                        &-\xi+\xi_1-\xi_2-\xi_4, \quad j = 3, \\
                        &-\xi-\xi_1-\xi_2-\xi_4, \quad j = 4,
                     \end{aligned} \right.
\end{equation*}
and
\begin{equation}
    \begin{aligned}
        \Omega_{1,\ulell}(\xi,\xi_1,\xi_2,\xi_4) &:= -\jxi + \jxione + \jxitwo + \jxithree + \ulell \xi_4, \\
        \Omega_{2,\ulell}(\xi,\xi_1,\xi_2,\xi_4) &:= -\jxi + \jxione - \jxitwo + \jxithree + \ulell \xi_4, \\
        \Omega_{3,\ulell}(\xi,\xi_1,\xi_2,\xi_4) &:= -\jxi + \jxione - \jxitwo - \jxithree + \ulell \xi_4, \\
        \Omega_{4,\ulell}(\xi,\xi_1,\xi_2,\xi_4) &:= -\jxi - \jxione - \jxitwo - \jxithree + \ulell \xi_4.
    \end{aligned}
\end{equation}
Here we freely inserted a cut-off $\widetilde{\chi}_0(\xi_4)$ to frequencies $|\xi_4| \ll 1$, because for frequencies $|\xi_4| \gtrsim 1$ the interactions \eqref{equ:structure_cubic_nonlinearities_PV} are just regular interactions due to the rapid decay of the $\cosech$ function. In the next subsection all regular interactions are shown to be integrable error terms.

We extract the asserted asymptotics \eqref{equ:pointwise_calIpv_asymptotics} of the oscillatory integrals with a Hilbert-type kernel $\jxi^{\frac32} \calI_j^{\pvdots}(t,\xi)$, $1 \leq j \leq 4$, in a two step procedure. 
First, for every fixed output frequency $\xi \in \bbR$ and for every fixed frequency $\xi_4 \in \bbR$, we carry out a stationary phase analysis for the two-dimensional oscillatory integrals (with respect to the frequency variables $\xi_1$, $\xi_2$) in the integrand of the one-dimensional integral with respect to the frequency variable $\xi_4$. 
To leading order, this yields (one-dimensional) oscillatory integrals in the frequency variable $\xi_4$ with a Hilbert-type kernel. 
In a second step we then extract their asymptotic behavior.

We present the details at the example of the resonant term $\jxi^{\frac32} \calI_2^{\pvdots}(t,\xi)$. The other non-resonant terms can be treated analogously.
For every fixed output frequency $\xi \in \bbR$ and for every fixed $\xi_4 \in \bbR$, the phases $\Omega(\xi, \cdot, \cdot, \xi_4)$ have a unique non-degenerate critical point. By direct computation we find 
\begin{equation*}
    \begin{aligned}
    \nabla_{\xi_1,\xi_2} \Omega_{1,\ulell}(\xi,\xi_1,\xi_2,\xi_4) = 0 \quad &\Leftrightarrow \quad (\xi_1, \xi_2) = \bigl( \tstyfracthird (\xi+\xi_4), \tstyfracthird (\xi+\xi_4) \bigr), \\    
    \nabla_{\xi_1,\xi_2} \Omega_{2,\ulell}(\xi,\xi_1,\xi_2,\xi_4) = 0 \quad &\Leftrightarrow \quad (\xi_1, \xi_2) = (\xi+\xi_4,\xi+\xi_4), \\
    \nabla_{\xi_1,\xi_2} \Omega_{3,\ulell}(\xi,\xi_1,\xi_2,\xi_4) = 0 \quad &\Leftrightarrow \quad (\xi_1, \xi_2) = (-\xi-\xi_4, -\xi-\xi_4), \\ 
    \nabla_{\xi_1,\xi_2} \Omega_{4,\ulell}(\xi,\xi_1,\xi_2,\xi_4) = 0 \quad &\Leftrightarrow \quad (\xi_1, \xi_2) = \bigl( -\tstyfracthird (\xi+\xi_4), -\tstyfracthird (\xi+\xi_4) \bigr).
    \end{aligned}
\end{equation*}
The Hessians of the phases $\Omega_{j,\ulell}(\xi,\cdot,\cdot,\xi_4)$, $1 \leq j \leq 4$, at the critical points are given by
\begin{equation*}
    \begin{aligned}
        \mathrm{Hess}_{\xi_1,\xi_2} \Omega_{1,\ulell}\bigl( \xi, \tstyfracthird (\xi+\xi_4) , \tstyfracthird (\xi+\xi_4), \xi_4 \bigr) &= \jap{ \tstyfracthird (\xi+\xi_4) }^{-3} \begin{bmatrix} 2 & 1 \\ 1 & 2 \end{bmatrix}, \\
        \mathrm{Hess}_{\xi_1,\xi_2} \Omega_{2,\ulell}\bigl( \xi, \xi+\xi_4, \xi+\xi_4, \xi_4 \bigr) &= \jap{\xi+\xi_4}^{-3} \begin{bmatrix} 2 & -1 \\ -1 & 0 \end{bmatrix}, \\ 
        \mathrm{Hess}_{\xi_1,\xi_2} \Omega_{3,\ulell}\bigl( \xi, -\xi-\xi_4, -\xi-\xi_4, \xi_4 \bigr) &= \jap{\xi+\xi_4}^{-3} \begin{bmatrix} 0 & 1 \\ 1 & -2 \end{bmatrix}, \\
        \mathrm{Hess}_{\xi_1,\xi_2} \Omega_{4,\ulell}\bigl( \xi, -\tstyfracthird (\xi+\xi_4), -\tstyfracthird (\xi+\xi_4), \xi_4 \bigr) &= - \jap{\tstyfracthird (\xi+\xi_4)}^{-3} \begin{bmatrix} 2 & 1 \\ 1 & 2 \end{bmatrix}.
    \end{aligned}
\end{equation*}
Correspondingly, their signatures at the critical points are
\begin{equation*}
    \begin{aligned}
        \mathrm{sign} \, \mathrm{Hess}_{\xi_1,\xi_2} \Omega_{1,\ulell}\bigl( \xi, \tstyfracthird (\xi+\xi_4) , \tstyfracthird (\xi+\xi_4), \xi_4 \bigr) &= 2, \\
        \mathrm{sign} \, \mathrm{Hess}_{\xi_1,\xi_2} \Omega_{2,\ulell}\bigl( \xi, \xi+\xi_4, \xi+\xi_4, \xi_4 \bigr) &= 0, \\
        \mathrm{sign} \, \mathrm{Hess}_{\xi_1,\xi_2} \Omega_{3,\ulell}\bigl( \xi, -\xi-\xi_4, -\xi-\xi_4, \xi_4 \bigr) &= 0, \\
        \mathrm{sign} \, \mathrm{Hess}_{\xi_1,\xi_2} \Omega_{4,\ulell}\bigl( \xi, -\tstyfracthird (\xi+\xi_4), -\tstyfracthird (\xi+\xi_4), \xi_4 \bigr) &= -2.
    \end{aligned}
\end{equation*}
Moreover, the values of the phases at the critical points are
\begin{equation*}
    \begin{aligned}
        \Omega_{1,\ulell}\bigl( \xi, \tstyfracthird (\xi+\xi_4) , \tstyfracthird (\xi+\xi_4), \xi_4 \bigr) &= -\jxi + 3 \jap{ \tstyfracthird (\xi+\xi_4) } + \ulell \xi_4, \\
        \Omega_{2,\ulell}\bigl( \xi, \xi+\xi_4, \xi+\xi_4, \xi_4 \bigr) &= -\jxi + \jap{\xi+\xi_4} + \ulell \xi_4, \\
        \Omega_{3,\ulell}\bigl( \xi, -\xi-\xi_4, -\xi-\xi_4, \xi_4 \bigr) &= -\jxi -\jap{\xi+\xi_4} + \ulell \xi_4, \\
        \Omega_{4,\ulell}\bigl( \xi, -\tstyfracthird (\xi+\xi_4), -\tstyfracthird (\xi+\xi_4), \xi_4 \bigr) &= -\jxi - 3 \jap{ \tstyfracthird (\xi+\xi_4) } + \ulell \xi_4.
    \end{aligned}
\end{equation*}
Now for every fixed $\xi \in \bbR$ and for every fixed $\xi_4 \in \bbR$, we formally apply the stationary phase formula for the (two-dimensional) oscillatory integrals with respect to the frequency variables $\xi_1$ and $\xi_2$. 
We obtain that to leading order as $t \to \infty$,
\begin{equation} \label{equ:pointwise_pv_osc_integrals_two_dimensional_extracted}
\begin{aligned}
    \jxi^{\frac32} \calI_1^{\pvdots}(t,\xi) &= \frac{2\pi}{t} \frac{1}{\sqrt{3}} e^{i\frac{\pi}{2}} \jxi^{\frac32} \int e^{it(-\jxi+3\jap{\frac13(\xi+\xi_4)}+\ulell \xi_4)} \Bigl( \gulellsh\bigl(t, \tstyfracthird(\xi+\xi_4) \bigr) \Bigr)^3 \\
    &\qquad \times \frakm_{\ulell,+++}^{\pvdots}\bigl(\xi, \tstyfracthird(\xi+\xi_4), \tstyfracthird(\xi+\xi_4), \tstyfracthird(\xi+\xi_4) \bigr) \widetilde{\chi}_0(\xi_4) \, \pvdots \cosech\Bigl( \frac{\pi}{2\ulg} \xi_4 \Bigr) \, \ud \xi_4 + \ldots \\
    \jxi^{\frac32} \calI_2^{\pvdots}(t,\xi) &= \frac{2\pi}{t} \jxi^{\frac32} \int e^{it(-\jxi + \jap{\xi+\xi_4} + \ulell \xi_4)} \bigl| \gulellsh(t,\xi+\xi_4) \bigr|^2 \gulellsh(t,\xi+\xi_4) \\
    &\qquad \times \frakm_{\ulell,+-+}^{\pvdots}\bigl(\xi, \xi+\xi_4, \xi+\xi_4, \xi+\xi_4 \bigr) \widetilde{\chi}_0(\xi_4) \, \pvdots \cosech\Bigl( \frac{\pi}{2\ulg} \xi_4 \Bigr) \, \ud \xi_4 + \ldots \\ 
    \jxi^{\frac32} \calI_3^{\pvdots}(t,\xi) &= \frac{2\pi}{t} \jxi^{\frac32} \int e^{it(-\jxi -\jap{\xi+\xi_4} + \ulell \xi_4)} \bigl| \gulellsh(t,-\xi-\xi_4) \bigr|^2 \overline{\gulellsh(t,-\xi-\xi_4)} \\ 
    &\qquad \times \frakm_{\ulell,+--}^{\pvdots}\bigl(\xi,-\xi-\xi_4,-\xi-\xi_4,-\xi-\xi_4\bigr) \widetilde{\chi}_0(\xi_4) \, \pvdots \cosech\Bigl( \frac{\pi}{2\ulg} \xi_4 \Bigr) \, \ud \xi_4 + \ldots \\ 
    \jxi^{\frac32} \calI_4^{\pvdots}(t,\xi) &= \frac{2\pi}{t} \frac{1}{\sqrt{3}} e^{-i\frac{\pi}{2}} \jxi^{\frac32} \int e^{it(-\jxi -3\jap{\frac13 (\xi+\xi_4)} + \ulell \xi_4)} \Bigl( \overline{\gulellsh\bigl(t, -\tstyfracthird(\xi+\xi_4) \bigr)} \Bigr)^3 \\
    &\qquad \times \frakm_{\ulell,---}^{\pvdots}\bigl(\xi, -\tstyfracthird(\xi+\xi_4), -\tstyfracthird(\xi+\xi_4), -\tstyfracthird(\xi+\xi_4) \bigr) \widetilde{\chi}_0(\xi_4) \, \pvdots \cosech\Bigl( \frac{\pi}{2\ulg} \xi_4 \Bigr) \, \ud \xi_4 + \ldots 
\end{aligned}
\end{equation}
We describe at the example of $\jxi^{\frac32} \calI_2^{\pvdots}(t,\xi)$ how to make the preceding analysis rigorous.
In order to deal with the singularity of the Hilbert-type kernel, we remove a small neighborhood around $\xi_4 = 0$ as follows. We denote the two-dimensional oscillatory integral with respect to the frequency variables $\xi_1$ and $\xi_2$ by 
\begin{equation} \label{equ:pointwise_pv_stat_phase_rig1}
    \begin{aligned}
        H(t,\xi, \xi_4) &:= \iint e^{it \Omega_{2,\ulell}(\xi,\xi_1,\xi_2,\xi_4)} \gulellsh(t,\xi_1) \, \overline{\gulellsh(t,\xi_2)} \, \gulellsh(t,\xi_3) \, \jxione^{-1} \jxitwo^{-1} \jxithree^{-1} \\ 
        &\qquad \qquad \qquad \qquad \qquad \qquad \qquad \qquad \qquad \qquad \times \frakm_{\ulell, +-+}^{\pvdots}(\xi, \xi_1, \xi_2, \xi_3) \, \ud \xi_1 \, \ud \xi_2.
    \end{aligned}
\end{equation}
Then we write 
\begin{equation} \label{equ:pointwise_pv_stat_phase_rig2}
    \begin{aligned}
        \jxi^{\frac32} \calI_2^{\pvdots}(t,\xi) &= \int H(t,\xi,\xi_4) \chi_1\bigl(t^{10}\xi_4\bigr)  \widetilde{\chi}_0(\xi_4) \, \frac{\varphi(\xi_4)}{\xi_4} \, \ud \xi_4 \\
        &\quad + \int H(t,\xi,\xi_4) \chi_0\bigl(t^{10}\xi_4\bigr)  \widetilde{\chi}_0(\xi_4) \, \pvdots \frac{\varphi(\xi_4)}{\xi_4} \, \ud \xi_4,
    \end{aligned}
\end{equation}
where we introduced the Schwartz function 
\begin{equation} \label{equ:pointwise_schwartz_varphi_def}
    \varphi(\xi_4) := \xi_4 \cosech\Bigl( \frac{\pi}{2\ulg} \xi_4 \Bigr).
\end{equation}
The second integral on the right-hand side of \eqref{equ:pointwise_pv_stat_phase_rig2} is a fast decaying remainder term. To see this, we can use the definition of the principal value and the fact that $\varphi$ is even to write 
\begin{equation}
    \begin{aligned}
        \int H(t,\xi,\xi_4) \chi_0\bigl(t^{10}\xi_4\bigr)  \widetilde{\chi}_0(\xi_4) \, \pvdots \frac{\varphi(\xi_4)}{\xi_4} \, \ud \xi_4 = \int \bigl( H(t,\xi,\xi_4) - H(t,\xi,0) \bigr) \chi_0\bigl(t^{10}\xi_4\bigr)  \widetilde{\chi}_0(\xi_4) \, \pvdots \frac{\varphi(\xi_4)}{\xi_4} \, \ud \xi_4.
    \end{aligned}
\end{equation}
Using the fundamental theorem of calculus and the bootstrap assumptions \eqref{equ:prop_profile_bounds_assumption2} about the effective profile, we can show a crude bound such as $\bigl|H(t,\xi,\xi_4) - H(t,\xi,0) \bigr| \lesssim \varepsilon^3 t^2 |\xi_4|^{\frac12}$. While the latter grows in time quite strongly, it provides the factor $|\xi_4|^{\frac12}$ to remove the singularity from the Hilbert-type kernel, and then we can gain back much more time decay from the small frequency support of $\chi_0\bigl(t^{10}\xi_4\bigr)$. Overall the second integral on the right-hand side of \eqref{equ:pointwise_pv_stat_phase_rig2} then becomes a remainder term with fast time decay.

In the first integral on the right-hand side of \eqref{equ:pointwise_pv_stat_phase_rig2} we have removed the singularity at $\xi_4 = 0$, so that we can first carry out a stationary phase analysis for the (two-dimensional) oscillatory integral $H(t,\xi,\xi_4)$ by proceeding analogously as in the preceding Subsection~\ref{subsubsec:stationary_phase_Dirac}. We omit the lengthy but similar arguments. Subsequently we still have to carry out the integration with respect to $\xi_4$ over the region $t^{-10} \lesssim |\xi_4| \lesssim 1$, but this only costs an additional $\log(t)$ factor, which can be easily absorbed into the faster decaying remainder terms stemming from the stationary phase analysis of $H(t,\xi,\xi_4)$. 

This leaves us with the one-dimensional oscillatory integrals with a Hilbert-type kernel on the right-hand sides of \eqref{equ:pointwise_pv_osc_integrals_two_dimensional_extracted}.
In order to determine their asymptotics as $t \to \infty$, we continue with the example of $\jxi^{\frac32} \calI_2^{\pvdots}(t,\xi)$. The other terms can be treated analogously.
Concretely, we are thus left to determine the asymptotics as $t\to\infty$ of
\begin{equation} \label{equ:pointwise_calI2pvdots_leading_order_recalled}
    \begin{aligned}
        &\frac{2\pi}{t} \jxi^{\frac32} e^{-it(\jxi+\ulell\xi)} \int e^{it(\jap{\xi+\xi_4} + \ulell (\xi + \xi_4)} \bigl| \gulellsh(t,\xi+\xi_4) \bigr|^2 \gulellsh(t,\xi+\xi_4) \\
        &\qquad \qquad \qquad \qquad \qquad \times \frakm_{\ulell,+-+}^{\pvdots}\bigl(\xi, \xi+\xi_4, \xi+\xi_4, \xi+\xi_4 \bigr) \widetilde{\chi}_0(\xi_4) \, \pvdots \cosech\Bigl( \frac{\pi}{2\ulg} \xi_4 \Bigr) \, \ud \xi_4.
    \end{aligned}
\end{equation}
To this end recall from Subsection~\ref{subsec:cubic_spectral_distributions} that $\frakm_{\ulell,+-+}^{\pvdots}\bigl(\xi, \xi_1, \xi_2, \xi_3\bigr)$ is a linear combination of tensorized products $\overline{\frakb(\xi)} \frakb_1(\xi_1) \overline{\frakb_2(\xi_2)} \frakb_3(\xi_3)$, where $\frakb(\xi)$ is given by
\begin{equation*}
    \frac{1}{|\ulg(\xi+\ulell\jxi)|-i} \quad \text{or} \quad \frac{\ulg(\xi+\ulell\jxi)}{|\ulg(\xi+\ulell\jxi)|-i},
\end{equation*}
and where $\frakb_j(\xi_j)$, $1 \leq j \leq 3$, are given by
\begin{equation*}
    \frac{1}{|\ulg(\xi_j+\ulell\jap{\xi_j})|-i} \quad \text{or} \quad \frac{\ulg(\xi_j+\ulell\jap{\xi_j})}{|\ulg(\xi_j+\ulell\jap{\xi_j})|-i}.
\end{equation*}
Importantly, at least one of the factors in each product $\overline{\frakb(\xi)} \frakb_1(\xi_1) \overline{\frakb_2(\xi_2)} \frakb_3(\xi_3)$ must be of the type
\begin{equation} \label{equ:pointwise_calI2pvdots_good_coefficients}
    \frac{\ulg(\xi+\ulell\jxi)}{|\ulg(\xi+\ulell\jxi)|-i}, \quad \text{respectively} \quad \frac{\ulg(\xi_j+\ulell\jap{\xi_j})}{|\ulg(\xi_j+\ulell\jap{\xi_j})|-i}, \quad 1 \leq j \leq 3.
\end{equation}
For the analysis below, it will be useful to rewrite the terms \eqref{equ:pointwise_calI2pvdots_good_coefficients} as
\begin{equation*}
    \frac{\ulg(\xi+\ulell\jxi)}{|\ulg(\xi+\ulell\jxi)|-i} = \bigl(\xi+\ulell\jxi\bigr) \jxi^{-1} \frac{\ulg \jxi}{|\ulg(\xi+\ulell\jxi)|-i} =: \bigl(\xi+\ulell\jxi\bigr) \jxi^{-1} \frakc(\xi),
\end{equation*}
where we trivially observe that $\frakc \in W^{1,\infty}$. We introduce analogous coefficients $\frakc_j(\xi_j)$ for rewriting the coefficients $\frakb_j(\xi_j)$, $1 \leq j \leq 3$, when they are of the form \eqref{equ:pointwise_calI2pvdots_good_coefficients}.

Hence, the leading order expression \eqref{equ:pointwise_calI2pvdots_leading_order_recalled} for $\calI_2^{\pvdots}(t,\xi)$ can be written as a linear combination of two types of terms. If one of the coefficients $\frakb_j(\xi_j)$, $1 \leq j \leq 3$, for the input frequencies is of the type \eqref{equ:pointwise_calI2pvdots_good_coefficients}, then the corresponding contribution to \eqref{equ:pointwise_calI2pvdots_leading_order_recalled} (up to permutations of the inputs) is given by 
\begin{equation} \label{equ:pointwise_calI2pvdots_leading_order_type1}
    \begin{aligned}
        &\frac{2\pi}{t} \jxi^{\frac32} e^{-it\phi(\xi,0)} \overline{\frakb(\xi)} \int e^{it\phi(\xi,\sigma)} (\partial_\sigma \phi)(\xi,\sigma) F(t,\xi+\sigma) \widetilde{\chi}_0(\sigma) \, \pvdots \cosech\Bigl( \frac{\pi}{2\ulg} \sigma \Bigr) \, \ud \sigma,
    \end{aligned}
\end{equation}
where we have introduced the phase 
\begin{equation} \label{equ:pointwise_phi_xi_sigma_definition}
    \phi(\xi,\sigma) := \jap{\xi+\sigma} + \ulell (\xi+\sigma)
\end{equation}
and where 
\begin{equation*}
    F(t,\xi+\sigma) := \bigl( \frakc_1(\xi+\sigma) \gulellsh(t,\xi+\sigma) \bigr) \overline{\bigl( \frakb_2(\xi+\sigma) \gulellsh(t,\xi+\sigma) \bigr)} \bigl( \frakb_3(\xi+\sigma) \gulellsh(t,\xi+\sigma) \bigr). 
\end{equation*}
Here we used that $\frakb_1(\xi+\sigma) = (\partial_\sigma \phi)(\xi,\sigma) \frakc_1(\xi+\sigma)$ in view of the preceding definitions and the fact that
\begin{equation*}
    (\partial_\sigma \phi)(\xi,\sigma) = \bigl( \xi + \sigma + \ulell \jap{\xi+\sigma} \bigr) \jap{\xi+\sigma}^{-1}.
\end{equation*}

Instead, if the coefficient $\frakb(\xi)$ for the output frequency is of the type \eqref{equ:pointwise_calI2pvdots_good_coefficients}, then the corresponding contribution to \eqref{equ:pointwise_calI2pvdots_leading_order_recalled} is given by 
\begin{equation} \label{equ:pointwise_calI2pvdots_leading_order_type2}
    \begin{aligned}
        &\frac{2\pi}{t} \jxi^{\frac32} e^{-it\phi(\xi,0)} (\partial_\sigma \phi)(\xi,0) \overline{\frakc(\xi)} \int e^{it\phi(\xi,\sigma)} G(t,\xi+\sigma) \widetilde{\chi}_0(\sigma) \, \pvdots \cosech\Bigl( \frac{\pi}{2\ulg} \sigma \Bigr) \, \ud \sigma,
    \end{aligned}
\end{equation}
where $\phi(\xi,\sigma)$ is defined in \eqref{equ:pointwise_phi_xi_sigma_definition} and where
\begin{equation*}
    G(t,\xi+\sigma) := \bigl( \frakb_1(\xi+\sigma) \gulellsh(t,\xi+\sigma) \bigr) \overline{\bigl( \frakb_2(\xi+\sigma) \gulellsh(t,\xi+\sigma) \bigr)} \bigl( \frakb_3(\xi+\sigma) \gulellsh(t,\xi+\sigma) \bigr). 
\end{equation*}

The next lemma extracts asymptotics for the contributions to \eqref{equ:pointwise_calI2pvdots_leading_order_recalled} of the first type \eqref{equ:pointwise_calI2pvdots_leading_order_type1}.
Its proof is an adaptation of the proof of \cite[Lemma 7.6]{CP22} to the Klein-Gordon case.
\begin{lemma} \label{lem:pointwise_calI2pvdots_final_asymptotics_extract_type1}
    Suppose the assumptions in the statement of Proposition~\ref{prop:profile_bounds} are in place.
    Assume $T \geq 1$.
    Let $\chi_0(\cdot)$ be a smooth even cut-off supported on $[-1,1]$ and set $\chi_1(\cdot) := 1 - \chi_0(\cdot)$.
    Let $\frakc_1, \frakb_2, \frakb_3 \in W^{1,\infty}(\bbR)$. 
    Consider
    \begin{equation*}
        \calJ(t,\xi) := \int_\bbR e^{i\phi(\xi,\sigma) t} (\partial_\sigma \phi)(\xi, \sigma) F(t,\xi+\sigma) \widetilde{\chi}_0(\sigma) \, \pvdots \cosech\Bigl( \frac{\pi}{2\ulg} \sigma \Bigr) \, \ud \sigma 
    \end{equation*}
    with
    \begin{equation*}
        \phi(\xi,\sigma) := \jap{\xi+\sigma} + \ulell (\xi + \sigma)
    \end{equation*}
    and 
    \begin{equation*}
        F(t,\xi+\sigma) = \bigl( \frakc_1(\xi+\sigma) \gulellsh(t,\xi+\sigma) \bigr) \overline{\bigl( \frakb_2(\xi+\sigma) \gulellsh(t,\xi+\sigma) \bigr)} \bigl( \frakb_3(\xi+\sigma) \gulellsh(t,\xi+\sigma) \bigr).
    \end{equation*}
    Then uniformly for all $1 \leq t \leq T$,
    \begin{equation*}
        \biggl\| \jxi^{\frac32} \calJ(t,\xi) + 2i \ulg e^{i \phi(\xi,0) t} (\partial_\sigma \phi)(\xi,0) \jxi^{\frac32} F(t,\xi) \tanh\bigl( \ulg (\partial_\sigma \phi)(\xi,0) t \bigr) \chi_1\bigl( (\xi+\ulg\ulell) t^{3\delta} \bigr) \biggr\|_{L^\infty_\xi} \lesssim \varepsilon^3 t^{-1-2\delta}.
    \end{equation*}
\end{lemma}
\begin{proof}
    We again use the short-hand notation \eqref{equ:pointwise_schwartz_varphi_def}.
    Observe that by the bootstrap assumptions \eqref{equ:prop_profile_bounds_assumption2} and the localization $|\sigma| \lesssim 1$, we have uniformly for all $\xi \in \bbR$ and all $0 \leq t \leq T$,
    \begin{equation} \label{equ:pointwise_calI2pvdots_final_asymptotics_extract_type1_Fbounds}
        \begin{aligned}
            \bigl\| \jxi^{\frac32} F(t,\xi+\sigma) \widetilde{\chi}_0(\sigma) \bigr\|_{L^\infty_\sigma} &\lesssim \varepsilon^3, \\ \bigl\| \jxi^{\frac32} \partial_\sigma F(t,\xi+\sigma) \chi_0(\sigma) \bigr\|_{L^2_\sigma} &\lesssim \varepsilon^3 \jt^{2\delta}.
        \end{aligned}
    \end{equation}    
    Introducing suitable frequency cut-offs, we decompose 
    \begin{equation} \label{equ:pointwise_calI2pvdots_final_asymptotics_extract_type1_decomposition_calJ}
        \calJ(t,\xi) = \calJ_1(t,\xi) + \calJ_2(t,\xi) + \calJ_3(t,\xi) + \calJ_4(t,\xi),
    \end{equation}
    where 
    \begin{equation*}
        \begin{aligned}
            \calJ_1(t,\xi) &:= \jxi^{\frac32} \int_\bbR e^{i\phi(\xi,\sigma) t} (\partial_\sigma \phi)(\xi, \sigma) F(t,\xi+\sigma) \chi_0\bigl(\sigma t^3\bigr) \chi_0\bigl( (\sigma+\xi+\ulg\ulell) t^{5\delta} \bigr) \widetilde{\chi}_0(\sigma) \, \pvdots \, \frac{\varphi(\sigma)}{\sigma} \, \ud \sigma, \\ 
            \calJ_2(t,\xi) &:= \jxi^{\frac32} \int_\bbR e^{i\phi(\xi,\sigma) t} (\partial_\sigma \phi)(\xi, \sigma) F(t,\xi+\sigma) \chi_1\bigl(\sigma t^3\bigr) \chi_0\bigl( (\sigma+\xi+\ulg\ulell) t^{5\delta} \bigr) \widetilde{\chi}_0(\sigma) \, \frac{\varphi(\sigma)}{\sigma} \, \ud \sigma, \\ 
            \calJ_3(t,\xi) &:= \jxi^{\frac32} \int_\bbR e^{i\phi(\xi,\sigma) t} (\partial_\sigma \phi)(\xi, \sigma) F(t,\xi+\sigma) \chi_1\bigl(\sigma t^{1-10\delta}\bigr) \chi_1\bigl( (\sigma+\xi+\ulg\ulell) t^{5\delta} \bigr) \widetilde{\chi}_0(\sigma) \, \frac{\varphi(\sigma)}{\sigma} \, \ud \sigma, \\ 
            \calJ_4(t,\xi) &:= \jxi^{\frac32} \int_\bbR e^{i\phi(\xi,\sigma) t} (\partial_\sigma \phi)(\xi, \sigma) F(t,\xi+\sigma) \chi_0\bigl(\sigma t^{1-10\delta}\bigr) \chi_1\bigl( (\sigma+\xi+\ulg\ulell) t^{5\delta} \bigr) \widetilde{\chi}_0(\sigma) \, \pvdots \, \frac{\varphi(\sigma)}{\sigma} \, \ud \sigma.
        \end{aligned}
    \end{equation*}
    The leading order asymptotics of $\calJ(t,\xi)$ are contained in the term $\calJ_4(t,\xi)$. We first bound the other remainder terms $\calJ_j(t,\xi)$, $1 \leq j \leq 3$. 

    \medskip 

    \noindent \underline{Estimate for $\calJ_1(t,\xi)$.} 
    On the support of the integrand of $\calJ_1(t,\xi)$ we have $|\sigma| \lesssim t^{-3}$. We set
    \begin{equation*}
     H(t,\xi,\sigma) := \jxi^{\frac32} e^{i\phi(\xi,\sigma) t} (\partial_\sigma \phi)(\xi, \sigma) F(t,\xi+\sigma) \chi_0\bigl( (\sigma+\xi+\ulg\ulell) t^{5\delta} \bigr).
    \end{equation*} 
    Then by \eqref{equ:pointwise_calI2pvdots_final_asymptotics_extract_type1_Fbounds}, we have $\|\partial_\sigma H(t,\xi,\sigma)\|_{L^2_\sigma} \lesssim \varepsilon^3 t$. Using the principal value and that $\chi_0(\cdot)$, $\varphi(\cdot)$ are even, we can then bound
    \begin{equation*}
        \begin{aligned}
            \bigl| \calJ_1(t,\xi) \bigr| &= \biggl| \int_\bbR \bigl( H(t,\xi,\sigma) - H(t,\xi,0) \bigr) \chi_0\bigl(\sigma t^3\bigr) \widetilde{\chi}_0(\sigma) \, \pvdots \, \frac{\varphi(\sigma)}{\sigma} \, \ud \sigma \biggr| \\ 
            &\lesssim \int_\bbR \varepsilon^3 t |\sigma|^{\frac12} \chi_0\bigl(\sigma t^3\bigr) \frac{|\varphi(\sigma)|}{|\sigma|} \, \ud \sigma \lesssim \varepsilon^3 t^{-\frac12}.
        \end{aligned}
    \end{equation*}

    \medskip 

    \noindent \underline{Estimate for $\calJ_2(t,\xi)$.}     
    Here we observe that for every fixed $\xi \in \bbR$, we have $(\partial_\sigma \phi)(\xi,\sigma) = 0$ if and only if $\sigma = -\xi-\ulg \ulell$. Hence, we can freely insert $(\partial_\sigma \phi)(\xi,-\xi-\ulg \ulell)$ in the integrand and then use the localization $|\sigma+\xi+\ulg \ulell| \lesssim t^{5\delta}$ to bound 
    \begin{equation*}
        \begin{aligned}
            \bigl| \calJ_2(t,\xi) \bigr| &= \biggl| \jxi^{\frac32} \int_\bbR e^{i\phi(\xi,\sigma) t} \Bigl( (\partial_\sigma \phi)(\xi, \sigma) - (\partial_\sigma \phi)(\xi,-\xi-\ulg) \Bigr) F(t,\xi+\sigma) \\ 
            &\qquad \qquad \qquad \qquad \qquad \qquad \times \chi_1\bigl(\sigma t^3\bigr) \chi_0\bigl( (\sigma+\xi+\ulg\ulell) t^{5\delta} \bigr) \widetilde{\chi}_0(\sigma) \, \frac{\varphi(\sigma)}{\sigma} \, \ud \sigma \biggr| \\ 
            &\lesssim \bigl\| \jxi^{\frac32} F(t,\xi+\sigma) \chi_0(\sigma) \bigr\|_{L^\infty_\sigma} \int_\bbR \bigl| \sigma + \xi + \ulg\bigr| \chi_1\bigl(\sigma t^3\bigr) \chi_0\bigl( (\sigma+\xi+\ulg\ulell) t^{5\delta} \bigr) \widetilde{\chi}_0(\sigma) \, \frac{|\varphi(\sigma)|}{|\sigma|} \, \ud \sigma \\
            &\lesssim \varepsilon^3 t^{2\delta} \cdot t^{-5\delta} \cdot \log(\jt) \lesssim \varepsilon^3 t^{-2\delta}.
        \end{aligned}
    \end{equation*}

    \medskip 

    \noindent \underline{Estimate for $\calJ_3(t,\xi)$.}     
    Using that $|\sigma| \gtrsim t^{-1+10\delta}$ on the support of the integrand in $\calJ_3(t,\xi)$, we integrate by parts and use \eqref{equ:pointwise_calI2pvdots_final_asymptotics_extract_type1_Fbounds} to obtain
    \begin{equation*}
        \begin{aligned}
            \bigl| \calJ_3(t,\xi) \bigr| &\lesssim t^{-1} \int_\bbR \, \biggl| \partial_\sigma \biggl( \jxi^{\frac32} F(t,\xi+\sigma) \chi_1\bigl(\sigma t^{1-10\delta}\bigr) \chi_1\bigl( (\sigma+\xi+\ulg\ulell) t^{5\delta} \bigr) \widetilde{\chi}_0(\sigma) \, \frac{\varphi(\sigma)}{\sigma} \biggr) \biggr| \, \ud \sigma \lesssim \varepsilon^3 t^{-10\delta}.            
        \end{aligned}
    \end{equation*}   
    
    \medskip 

    \noindent \underline{Asymptotics for the term $\calJ_4(t,\xi)$.}
    In $\calJ_4(t,\xi)$ we can replace 
    \begin{equation*}
        (\partial_\sigma \phi)(\xi, \sigma) F(t,\xi+\sigma) \chi_1\bigl( (\sigma+\xi+\ulg\ulell) t^{5\delta} \bigr) \widetilde{\chi}_0(\sigma)
    \end{equation*}
    by the same expression with $\sigma = 0$, i.e., by $(\partial_\sigma \phi)(\xi, 0) F(t,\xi) \chi_1\bigl( (\xi+\ulg\ulell) t^{5\delta} \bigr)$, at the expense of picking up a mild error term
    \begin{equation*}
        \begin{aligned}
            \biggl\| \calJ_4(t,\xi) - \jxi^{\frac32} (\partial_\sigma \phi)(\xi, 0) F(t,\xi) \chi_1\bigl( (\xi+\ulg\ulell) t^{5\delta} \bigr) \int_\bbR e^{i\phi(\xi,\sigma) t} \chi_0\bigl(\sigma t^{1-10\delta}\bigr) \, \pvdots \, \frac{\varphi(\sigma)}{\sigma} \, \ud \sigma \biggr\|_{L^\infty_\xi} \lesssim \varepsilon^3 t^{-\frac12+10\delta}.
        \end{aligned}
    \end{equation*}
    Next, we carry out a Taylor expansion of the phase around $\sigma = 0$,
    \begin{equation*}
        e^{i \phi(\xi,\sigma) t} = e^{i \phi(\xi,0) t} e^{i \partial_\sigma \phi(\xi,0) \sigma t} e^{i \calO(\sigma^2) t},
    \end{equation*}
    and we replace $e^{i \calO(\sigma^2) t}$ by $1$ at the expense of a mild error term bounded by $\varepsilon^3 t^{-1+20\delta}$. 
    At this point we have 
    \begin{equation*}
        \begin{aligned}
            &\biggl\| \calJ_4(t,\xi) - \jxi^{\frac32} (\partial_\sigma \phi)(\xi, 0) F(t,\xi) \chi_1\bigl( (\xi+\ulg\ulell) t^{5\delta} \bigr) e^{i\phi(\xi,0) t} \\
            &\qquad \qquad \qquad \qquad \qquad \times \biggl( \int_\bbR e^{i \partial_\sigma \phi(\xi,0) \sigma t} \chi_0\bigl(\sigma t^{1-10\delta}\bigr) \, \pvdots \, \frac{\varphi(\sigma)}{\sigma} \, \ud \sigma \biggr) \biggr\|_{L^\infty_\xi} \lesssim \varepsilon^3 t^{-\frac12+10\delta}.
        \end{aligned}
    \end{equation*}
    We can now freely add $\chi_1(\sigma t^{1-10\delta})$ back into the integrand. 
    Indeed, since $\partial_\sigma \phi(\xi,0) = 0$ if and only if $\xi = -\ulg \ulell$, and since $|\xi+\ulg \ulell| \gtrsim t^{-5\delta}$, we have $|\partial_\xi \phi(\xi,0)|^{-1} \lesssim t^{5\delta}$, so we can integrate by parts with respect to $\sigma$ in the corresponding additional error term and bound it by $\varepsilon^3 t^{-5\delta}$. 
    Recalling that
    \begin{equation*}
        \pvdots \, \frac{\varphi(\sigma)}{\sigma} = \pvdots \cosech\Bigl( \frac{\pi}{2\ulg} \sigma \Bigr) = -i \sqrt{\frac{2}{\pi}} \whatcalF\bigl[ \tanh(\cdot) \bigr]\Bigl(\frac{\sigma}{\gamma}\Bigr) = -i \ulg \sqrt{\frac{2}{\pi}} \whatcalF\bigl[ \tanh(\ulg \cdot) \bigr](\sigma),
    \end{equation*}
    we obtain in the sense of distributions 
    \begin{equation*}
        \int_\bbR e^{i \partial_\sigma \phi(\xi,0) \sigma t} \, \pvdots \, \frac{\varphi(\sigma)}{\sigma} \, \ud \sigma = \sqrt{2\pi} \whatcalF^{-1}\Bigl[ \pvdots \cosech\Bigl( \frac{\pi}{2\ulg} \cdot \Bigr) \Bigr]\bigl( (\partial_\sigma \phi)(\xi,0) t \bigr) = -2i\ulg \tanh\bigl( \ulg (\partial_\sigma \phi)(\xi,0) t \bigr),
    \end{equation*}
    whence 
    \begin{equation*}
        \begin{aligned}
            \biggl\| \calJ_4(t,\xi) + \jxi^{\frac32} (\partial_\sigma \phi)(\xi, 0) F(t,\xi) \chi_1\bigl( (\xi+\ulg\ulell) t^{5\delta} \bigr) e^{i\phi(\xi,0) t}  2i \ulg \tanh\bigl( \ulg (\partial_\sigma \phi)(\xi,0) t \bigr) \biggr\|_{L^\infty_\xi} \lesssim \varepsilon^3 t^{-5\delta}.
        \end{aligned}
    \end{equation*}

    Combining the preceding estimates finishes the proof of the lemma.
\end{proof}

The next lemma extracts asymptotics for the contributions to \eqref{equ:pointwise_calI2pvdots_leading_order_recalled} of the second type \eqref{equ:pointwise_calI2pvdots_leading_order_type2}.
\begin{lemma} \label{lem:pointwise_calI2pvdots_final_asymptotics_extract_type2}
    Suppose the assumptions in the statement of Proposition~\ref{prop:profile_bounds} are in place.
    Assume $T \geq 1$.
    Let $\chi_0(\cdot)$ be a smooth even cut-off supported on $[-1,1]$ and set $\chi_1(\cdot) := 1 - \chi_0(\cdot)$.
    Let $\frakb_1, \frakb_2, \frakb_3, \frakc \in W^{1,\infty}(\bbR)$. 
    Consider
    \begin{equation*}
        \calK(t,\xi) := (\partial_\sigma \phi)(\xi,0) \int_\bbR e^{i\phi(\xi,\sigma) t} G(t,\xi+\sigma) \widetilde{\chi}_0(\sigma) \, \pvdots \cosech\Bigl( \frac{\pi}{2\ulg} \sigma \Bigr) \, \ud \sigma 
    \end{equation*}
    with
    \begin{equation*}
        \phi(\xi,\sigma) := \jap{\xi+\sigma} + \ulell (\xi + \sigma)
    \end{equation*}
    and 
    \begin{equation*}
        G(t,\xi+\sigma) = \bigl( \frakb_1(\xi+\sigma) \gulellsh(t,\xi+\sigma) \bigr) \overline{\bigl( \frakb_2(\xi+\sigma) \gulellsh(t,\xi+\sigma) \bigr)} \bigl( \frakb_3(\xi+\sigma) \gulellsh(t,\xi+\sigma) \bigr).
    \end{equation*}
    Then we have uniformly for all $1 \leq t \leq T$,
    \begin{equation*}
        \biggl\| \jxi^{\frac32} \calK(t,\xi) + 2i \ulg e^{i \phi(\xi,0) t} (\partial_\sigma \phi)(\xi,0) \jxi^{\frac32} G(t,\xi) \tanh\bigl( \ulg (\partial_\sigma \phi)(\xi,0) t \bigr) \chi_1\bigl( (\xi+\ulg\ulell) t^{3\delta} \bigr) \biggr\|_{L^\infty_\xi} \lesssim \varepsilon^3 t^{-1-2\delta}.
    \end{equation*}
\end{lemma}
\begin{proof}
    We can largely proceed as in the proof of Lemma~\ref{lem:pointwise_calI2pvdots_final_asymptotics_extract_type1} and first enact a frequency decomposition of $\calK(t,\xi)$ analogous to \eqref{equ:pointwise_calI2pvdots_final_asymptotics_extract_type1_decomposition_calJ}. 
    In the proof of Lemma~\ref{lem:pointwise_calI2pvdots_final_asymptotics_extract_type1} for the terms $\calJ_j(t,\xi)$, $1 \leq j \leq 4$, the presence of the factor $(\partial_\sigma \phi)(\xi,\sigma)$ in the integrand was only relevant for bounding the term $\calJ_2(t,\xi)$.
    The analogous term for $\calK(t,\xi)$ reads
    \begin{equation*}
        \begin{aligned}
            \calK_2(t,\xi) = \jxi^{\frac32} (\partial_\sigma \phi)(\xi, 0) \int_\bbR e^{i\phi(\xi,\sigma) t} G(t,\xi+\sigma) \chi_1\bigl(\sigma t^3\bigr) \chi_0\bigl( (\sigma+\xi+\ulg\ulell) t^{5\delta} \bigr) \widetilde{\chi}_0(\sigma) \, \frac{\varphi(\sigma)}{\sigma} \, \ud \sigma.  
        \end{aligned}
    \end{equation*}
    But here we can still obtain an acceptable error bound by distinguishing the cases $|\xi+\ulg \ulell| \lesssim t^{-3\delta}$ and $|\xi+\ulg\ulell| \gg t^{-3\delta}$. In the former case, the localization $|\xi+\ulg \ulell| \lesssim t^{-3\delta}$ implies $|(\partial_\sigma \phi)(\xi,\sigma)| \lesssim t^{-3\delta}$, whence
    \begin{equation*}
        \bigl| \calK_2(t,\xi) \bigr| \lesssim t^{-3\delta} \cdot \bigl\| \jxi^{\frac32} G(t,\xi+\sigma) \chi_0(\sigma) \bigr\|_{L^\infty_\sigma} \cdot \int_\bbR \chi_1\bigl(\sigma t^3\bigr) \frac{|\varphi(\sigma)|}{|\sigma|} \, \ud \sigma \lesssim t^{-3\delta} \cdot \varepsilon^3 \cdot \log(\jt) \lesssim \varepsilon^3 t^{-2\delta}.
    \end{equation*}
    In the latter case $|\xi+\ulg\ulell| \gg t^{-3\delta} \gg t^{-5\delta}$ (for $t \gg 1$) the localization $|\sigma+\xi+\ulg\ulell| \lesssim t^{-5\delta}$ implies that we must then actually have $|\sigma| \gtrsim t^{-3\delta}$, which is much better than the bound away from zero furnished by the cut-off $\chi_1(\sigma t^3)$. Thus, in the latter case we can estimate by
    \begin{equation*}
        \begin{aligned}
        \bigl| \calK_2(t,\xi) \bigr| &\lesssim \bigl\| \jxi^{\frac32} G(t,\xi+\sigma) \chi_0(\sigma) \bigr\|_{L^\infty_\sigma} \int_{\{|\sigma| \gtrsim t^{-3\delta}\}} \chi_0\bigl( (\sigma+\xi+\ulg\ulell) t^{5\delta} \bigr) \widetilde{\chi}_0(\sigma) \, \frac{\varphi(\sigma)}{\sigma} \, \ud \sigma \\ 
        &\lesssim \varepsilon^3 \cdot t^{-5\delta} \cdot t^{3\delta} \lesssim \varepsilon^3 t^{-2\delta}.
        \end{aligned}
    \end{equation*}
    This suffices and concludes the proof.
\end{proof}

Finally, writing $\frakm_{\ulell,+-+}^{\pvdots}\bigl(\xi, \xi_1, \xi_2, \xi_3)$ as a linear combination of tensorized products  using Lemma~\ref{lem:pointwise_calI2pvdots_final_asymptotics_extract_type1},  Lemma~\ref{lem:pointwise_calI2pvdots_final_asymptotics_extract_type2}, and reassembling all terms, we extract the following leading order behavior uniformly for $1 \leq t \leq T$,
\begin{equation}
    \begin{aligned}
        &\biggl\| \jxi^{\frac32} \calI_2^{\pvdots}(t,\xi) + \frac{4\pi i\ulg}{t} \frakm_{\ulell,+-+}^{\pvdots}\bigl(\xi, \xi, \xi, \xi) \jxi^{\frac32} \bigl| \gulellsh(t,\xi) \bigr|^2 \gulellsh(t,\xi) \\ 
        &\qquad \qquad \qquad \qquad \qquad \times \tanh\bigl( \ulg (\partial_\sigma \phi)(\xi,0) t \bigr) \chi_1\bigl( (\xi+\ulg\ulell) t^{5\delta} \bigr) \biggr\|_{L^\infty_\xi} \lesssim \varepsilon^3 t^{-1-2\delta}.
    \end{aligned}
\end{equation}
Owing to the remarkable vanishing property $\frakm_{\ulell,+-+}^{\pvdots}\bigl(\xi, \xi, \xi, \xi) = 0$ for all $\xi \in \bbR$ observed in \eqref{equ:cubic_spectral_distributions_pv_vanishing},
we arrive at the leading order behavior for $\jxi^{\frac32} \calI_2^{\pvdots}(t,\xi)$ asserted in the statement of Lemma~\ref{lem:cubic_Hilbert_stationary_phase}. 
Similarly, we deduce the asserted asymptotics for the other cubic interaction terms with Hilbert-type kernels

\subsubsection{Contributions of regular cubic interactions}

Finally, we show that all regular cubic interaction terms are integrable remainder terms.

\begin{lemma}
    Suppose the assumptions in the statement of Proposition~\ref{prop:profile_bounds} are in place.
    Assume $T \geq 1$.
    For each $1 \leq j \leq 4$ we then have for all $1 \leq t \leq T$ that
    \begin{equation}
        \bigl\| \jxi^{\frac32} \calI_j^{\mathrm{reg}}(t,\xi) \bigr\|_{L^\infty_\xi} \lesssim \varepsilon^3 t^{-\frac32}.
    \end{equation}
\end{lemma}
\begin{proof}
We provide the details for the term $\calI_2^{\reg}(t,\xi)$, the other terms being analogous. 
In view of the fine structure of the cubic spectral distributions analyzed in Subsection~\ref{subsec:cubic_spectral_distributions}, the term $\calI_2^{\mathrm{reg}}(t,\xi)$ is a linear combination of terms of the schematic form
\begin{equation*}
    \begin{aligned}
        \calI_{2, \mathrm{schem}}^{\mathrm{reg}}(t,\xi) &:= \frakb(\xi) \iiint e^{i t \Phi_{2,\ulell}(\xi,\xi_1,\xi_2,\xi_3)} h_1(t,\xi_1) \, \overline{h_2(t,\xi_2)} h_3(s,\xi_3) \\ 
        &\quad \quad \quad \quad \quad \quad \quad \quad \times \ulg^{-1} \widehat{\calF}^{-1}\bigl[\varphi\bigr]\bigl( \ulg^{-1} (-\xi+\xi_1-\xi_2+\xi_3) \bigr) \, \ud \xi_1 \, \ud \xi_2 \, \ud \xi_3,
    \end{aligned}
\end{equation*}
with the phase
\begin{equation*}
    \Phi_{2,\ulell}(\xi,\xi_1,\xi_2,\xi_3) := -(\jxi + \ulell \xi) + (\jxione + \ulell \xi_1) - (\jxitwo + \ulell \xi_2) + (\jxithree + \ulell \xi_3),
\end{equation*}
and for some coefficients $\frakb, \frakb_1, \frakb_2, \frakb_3 \in W^{1,\infty}(\bbR)$ such that the inputs are given by 
\begin{equation*}
    h_j(t,\xi_j) := \jap{\xi_j}^{-1} \frakb_j(\xi_j) \gulellsh(t,\xi_j), \quad 1 \leq j \leq 3,
\end{equation*}
and where $\varphi \in \calS(\bbR)$ is some Schwartz function.
In terms of the auxiliary linear Klein-Gordon evolutions defined in \eqref{equ:consequences_aux_KG_evolutions_def}, we thus have
\begin{equation*}
    \begin{aligned}
        \jxi^{\frac32} \calI_{2, \mathrm{schem}}^{\mathrm{reg}}(t,\xi) &= (2\pi)^{\frac32} e^{-it(\jxi+\ulell\xi)} \frakb(\xi) \jxi^{\frac32} \widehat{\calF}\bigl[ v_1(t,\cdot) \overline{v_2(t,\cdot)} v_3(t,\cdot) \varphi(\ulg \cdot) \bigr](\xi). 
    \end{aligned}
\end{equation*}
By standard Fourier theory, Sobolev estimates, H\"older's inequality, and the local decay estimate~\eqref{equ:consequences_aux_KG_local_H3y_decay}, we conclude 
\begin{equation}
    \begin{aligned}
        \bigl\| \jxi^{\frac32} \calI_{2, \mathrm{schem}}^{\mathrm{reg}}(t,\xi) \bigr\|_{L^\infty_\xi} &\lesssim \bigl\| v_1(t) \overline{v_2(t)} v_3(t) \varphi(\ulg \cdot) \bigr\|_{W^{2,1}_y} \\ 
        &\lesssim \bigl\|\jy^3 \varphi(\ulg \cdot) \bigr\|_{W^{2,\infty}_y} \bigl\| \jy^{-1} v_1(t) \bigr\|_{H^2_y} \bigl\| \jy^{-1} v_2(t) \bigr\|_{H^2_y} \bigl\| \jy^{-1} v_3(t) \bigr\|_{H^2_y} \\ 
        &\lesssim \varepsilon^3 \jt^{-\frac32},
    \end{aligned}
\end{equation}
as desired.
\end{proof}

\subsection{Auxiliary pointwise bounds} \label{subsec:auxiliary_pointwise_bounds}

Next, we deduce auxiliary pointwise bounds that suffice to treat the second and the third term on the right-hand side of the renormalized evolution equation for the effective profile \eqref{equ:pointwise_evol_equ_g_multiplied_by_jxi} as integrable remainder terms.

\begin{lemma} \label{lem:pointwise_remainder_terms_bounds}
    Suppose the assumptions in the statement of Proposition~\ref{prop:profile_bounds} are in place.
    Then we have for all $0 \leq t \leq T$,
    \begin{align}
        \bigl\| \jxi^3 B\bigl[ \gulellsh \bigr](t,\xi) \bigr\|_{L^\infty_\xi} &\lesssim \varepsilon^2 \jt^{-1}, \label{equ:pointwise_weighted_B_Linfty_bound} \\
        \bigl\| \jxi^{\frac32} \calR(t,\xi) \bigr\|_{L^\infty_\xi} + \bigl\| \jxi^{\frac32} \widetilde{\calR}(t,\xi) \bigr\|_{L^\infty_\xi} &\lesssim \varepsilon^2 \jt^{-\frac32+2\delta}. \label{equ:pointwise_weighted_calR_Linfty_bound}
    \end{align}
\end{lemma}
\begin{proof}
    The first asserted bound \eqref{equ:pointwise_weighted_B_Linfty_bound} is immediate from the definition of $B\bigl[ \gulellsh \bigr](t,\xi)$ in \eqref{equ:setting_up_definition_B} in view of Remark~\ref{rem:setting_up_rapid_decay_qjs} and the decay estimate \eqref{equ:consequences_hulell_decay}. 
    So we now turn to the proof of the second asserted bound \eqref{equ:pointwise_weighted_calR_Linfty_bound}.
    We verify \eqref{equ:pointwise_weighted_calR_Linfty_bound} separately for all terms in the definition of $\widetilde{\calR}(t,\xi)$ in \eqref{equ:setting_up_definition_wtilR} and in the definition of $\calR(t,\xi)$ in \eqref{equ:setting_up_definition_calR}.

    \medskip 

    \noindent \underline{Contribution of $(\dot{q}(t)-\ulell) \calL_\ulell(\bmu(t))$.}
    Using the mapping property \eqref{eq:calLellLinfty} along with the decay estimates \eqref{equ:consequences_qdot_minus_ulell_decay}, \eqref{equ:consequences_local_decay_uonetwo}, we obtain
    \begin{equation*}
        \begin{aligned}
            \bigl\| \jxi^{\frac32} (\dot{q}(t) - \ulell) \calL_\ulell\bigl( \bmu(t) \bigr)(\xi) \bigr\|_{L^\infty_\xi} &\lesssim |\dot{q}(t)-\ulell| \bigl\| \jxi^{\frac32} \calL_\ulell\bigl( \bmu(t) \bigr)(\xi) \bigr\|_{L^\infty_\xi} \\ 
            &\lesssim |\dot{q}(t)-\ulell| \Bigl( \bigl\| \sech^2(\ulg \cdot) u_1(t) \bigr\|_{W^{2,1}_y} + \bigl\| \sech^2(\ulg \cdot) u_2(t) \bigr\|_{W^{1,1}_y} \Bigr) \\
            &\lesssim \varepsilon \jt^{-1+\delta} \cdot \varepsilon \jt^{-\frac12} \lesssim \varepsilon^2 \jt^{-\frac32+\delta}.
        \end{aligned}
    \end{equation*}

    \medskip 

    \noindent \underline{Contribution of $\calMod(t)$.} 
    By Corollary~\ref{cor:TellP} we have
    \begin{equation*}
        \begin{aligned}
            \calF_{\ulell,D}^{\#}\bigl[ \bigl(\calMod\bigr)_1(t) \bigr](\xi) - \calF_{\ulell}^{\#}\bigl[ \bigl(\calMod\bigr)_2(t) \bigr](\xi) = \calF_{\ulell,D}^{\#}\bigl[ \bigl(\ulPe \calMod\bigr)_1(t) \bigr](\xi) - \calF_{\ulell}^{\#}\bigl[ \bigl(\ulPe \calMod\bigr)_2(t) \bigr](\xi).
        \end{aligned}
    \end{equation*}
    Moreover, since $\ulPe \bmY_{j,\ulell} = 0$ for $j=0,1$, we can write
    \begin{equation*}
        \begin{aligned}
            \ulPe \calMod(t,y) = - \bigl(\dot{q}(t)-\ell(t)\bigr) \ulPe \bigl( \bmY_{0,\ell(t)}(y) - \bmY_{0,\ulell}(y) \bigr) - \dot{\ell}(t) \ulPe \bigl( \bmY_{1,\ell(t)}(y) - \bmY_{1,\ulell}(y) \bigr).
        \end{aligned}
    \end{equation*}
    Then using the mapping properties \eqref{equ:mapping_property_calFulellsh_jxi32_Linfty}, \eqref{equ:mapping_property_calFulellDsh_jxi32_Linfty} and the decay estimates \eqref{equ:prop_profile_bounds_assumption1}, \eqref{equ:consequences_crude_decay_modulation_parameters}, we find that 
    \begin{equation*}
        \begin{aligned}
            &\Bigl\| \jxi^{\frac32} \Bigl( \calF_{\ulell,D}^{\#}\bigl[ \bigl(\calMod(t)\bigr)_1 \bigr](\xi) - \calF_{\ulell}^{\#}\bigl[ \bigl(\calMod(t)\bigr)_2 \bigr](\xi) \Bigr) \Bigr\|_{L^\infty_\xi} \\ 
            &\lesssim |\dot{q}(t)-\ell(t)| \bigl\| \ulPe \bigl( \bmY_{0,\ell(t)}(y) - \bmY_{0,\ulell}(y) \bigr) \bigr\|_{W^{3,1}_y} + |\dot{\ell}(t)| \bigl\| \ulPe \bigl( \bmY_{1,\ell(t)}(y) - \bmY_{1,\ulell}(y) \bigr) \bigr\|_{W^{3,1}_y} \\ 
            &\lesssim \bigl( |\dot{q}(t)-\ell(t)| + |\dot{\ell}(t)| \bigr) |\ell(t) - \ulell| \lesssim \varepsilon \jt^{-1} \cdot \varepsilon \jt^{-1+\delta} \lesssim \varepsilon^2 \jt^{-2+\delta}.
        \end{aligned}
    \end{equation*}

    \medskip 
    
    \noindent \underline{Contribution of $\calC_{\ell(t)}\bigl(\usubeone(t)\bigr)$.}    
    Similarly, we obtain by \eqref{equ:mapping_property_calFulellsh_jxi32_Linfty}, Sobolev embedding, and \eqref{equ:consequences_local_decay_usubeonetwo} that
    \begin{equation*}
        \begin{aligned}
            \bigl\| \jxi^{\frac32} \calF_{\ulell}^{\#}\bigl[ \calC_{\ell(t)}\bigl( \usubeone(t) \bigr) \bigr](\xi) \bigr\|_{L^\infty_\xi} \lesssim \bigl\| \sech^2(\gamma(t) \cdot) \usubeone(t)^3 \bigr\|_{W^{2,1}_y} 
            &\lesssim \bigl\| \jy^{-1} \usubeone(t) \bigr\|_{H^2_y}^3 \lesssim \varepsilon^3 \jt^{-\frac32}.
        \end{aligned}
    \end{equation*}

    \medskip 
    
    \noindent \underline{Contributions of $\calR_1\bigl(\usubeone(t)\bigr)$ and $\calE_j(t)$, $1 \leq j \leq 3$.}      
    The contributions of the other spatially localized terms with at least cubic-type time decay $\calR_1\bigl(\usubeone(t)\bigr)$ and $\calE_j(t)$, $1 \leq j \leq 3$, can be estimated analogously to the preceding term $\calC_{\ell(t)}\bigl(\usubeone(t)\bigr)$.

    \medskip 
    
    \noindent \underline{Contribution of $\calR_2\bigl(\usubeone(t)\bigr)$.}
    The contributions of the quintic and higher-order nonlinearities in $\calR_2\bigl(\usubeone(t)\bigr)$ can be bounded using the mapping property \eqref{equ:mapping_property_calFulellsh_jxi32_Linfty} and \eqref{equ:consequences_dispersive_decay_ulPe_radiation}, \eqref{equ:consequences_local_decay_usubeonetwo} as follows
    \begin{equation*}
        \begin{aligned}
            \bigl\| \jxi^{\frac32} \calF_{\ulell}^{\#}\bigl[ \calR_2\bigl(\usubeone(t)\bigr) \bigr](\xi) \bigr\|_{L^\infty_\xi} &\lesssim \bigl\| \calR_2\bigl(\usubeone(t)\bigr)\bigr\|_{W^{2,1}_y} \lesssim \|\usubeone(t)\|_{H^2_y}^2 \|\usubeone(t)\|_{L^\infty_y}^3 \lesssim \varepsilon^5 \jt^{-\frac32 + 2\delta}.
        \end{aligned}
    \end{equation*}

    \noindent \underline{Contribution of $\calQ_{\ulell,r}\bigl(\usubeone(t)\bigr)$.}    
    Moreover, by the mapping property \eqref{equ:mapping_property_calFulellsh_jxi32_Linfty} and \eqref{equ:consequences_hulell_decay}, \eqref{equ:consequences_Remusubeone_H3y_local_decay}, we have  
    \begin{equation*}
        \begin{aligned}
            &\bigl\| \jxi^{\frac32} \calF_{\ulell}^{\#}\bigl[ \calQ_{\ulell,r}\bigl(\usubeone(t)\bigr)\bigr](\xi) \bigr\|_{L^\infty_\xi} 
            \lesssim \bigl\| \calQ_{\ulell,r}\bigl(\usubeone(t)\bigr) \bigr\|_{W^{2,1}_y} \\ 
            &\lesssim \bigl\| \jy^5 \alpha(\ulg y) \bigr\|_{L^\infty_y} \Bigl( |h_\ulell(t)| + \bigl\| \jy^{-3} R_{\usubeone}(t,y) \bigr\|_{H^2_y} \Bigr) \bigl\| \jy^{-3} R_{\usubeone}(t,y) \bigr\|_{H^2_y} 
            \lesssim \varepsilon^2 \jt^{-\frac32+2\delta}.
        \end{aligned}
    \end{equation*}

    \noindent \underline{Contribution of $\calR_q(t)$.}
    Finally, in view of Remark~\ref{rem:setting_up_rapid_decay_qjs} and the decay estimates \eqref{equ:consequences_hulell_decay}, \eqref{equ:consequences_hulell_phase_filtered_decay}, we obtain straightaway from the definition of $\calR_q(t,\xi)$ in \eqref{equ:setting_up_definition_calRq} that
    \begin{equation*}
        \begin{aligned}
            \bigl\| \jxi^{\frac32} \calR_q(t,\xi) \bigr\|_{L^\infty_\xi} \lesssim \varepsilon^2 \jt^{-\frac32+\delta}.
        \end{aligned}
    \end{equation*}

    Combining the preceding estimates implies the second asserted estimate \eqref{equ:pointwise_weighted_calR_Linfty_bound}.    
\end{proof}

\subsection{Proof of Proposition~\ref{prop:pointwise_profile_bounds}} \label{subsec:proof_pointwise_profile_bounds}

Thanks to the preparations in the preceding two subsections, we are now in the position to present the proof of Proposition~\ref{prop:pointwise_profile_bounds}.

\begin{proof}[Proof of Proposition~\ref{prop:pointwise_profile_bounds}]
In view of Lemma~\ref{lem:pointwise_cubic_terms_asymptotics} and Lemma~\ref{lem:pointwise_remainder_terms_bounds},
the renormalized evolution equation \eqref{equ:pointwise_evol_equ_g_multiplied_by_jxi} for the (weighted) effective profile reads
\begin{equation} \label{equ:pointwise_evol_equ_g_renorm_leading_order}
    \begin{aligned}
            &\pt \biggl( e^{-i\xi\theta(t)} \Bigl( \jxi^{\frac32} g_\ulell^{\#}(t,\xi) + \jxi^{\frac32} B\bigl[ \gulellsh \bigr](t,\xi) \Bigr) \biggr) \\
            &= \frac{i}{16 t} \bigl| \gulellsh(t,\xi) \bigr|^2 e^{-i\xi\theta(t)} \jxi^{\frac32} \gulellsh(t,\xi) - \frac16 e^{-i\xi\theta(t)} \calT_{osc}(t,\xi) + \calO_{L^\infty_\xi}\bigl( \varepsilon^2 t^{-1-2\delta} \bigr),
    \end{aligned}
\end{equation}
where the term $\calT_{osc}(t,\xi)$ defined in \eqref{equ:pointwise_Tosc_definition} contains all critically decaying (non-resonant) cubic terms that oscillate in time.
We remove the critically decaying resonant cubic term on the right-hand side of \eqref{equ:pointwise_evol_equ_g_renorm_leading_order} using an integrating factor,
\begin{equation} \label{equ:pointwise_evol_equ_g_with_integrating_factor}
    \begin{aligned}
        &\pt \biggl( e^{-i\Lambda(t,\xi)} e^{-i\xi\theta(t)} \Bigl( \jxi^{\frac32} g_\ulell^{\#}(t,\xi) + \jxi^{\frac32} B\bigl[ \gulellsh \bigr](t,\xi) \Bigr) \biggr) \\
        &= e^{-i\Lambda(t,\xi)} \biggl( - \frac16 e^{-i\xi\theta(t)} \calT_{osc}(t,\xi) - \frac{i}{16} \frac{1}{t} \bigl| \gulellsh(t,\xi) \bigr|^2 e^{-i\xi\theta(t)} \jxi^{\frac32} B\bigl[ \gulellsh \bigr](t,\xi) + \calO_{L^\infty_\xi}\bigl( \varepsilon^2 t^{-1-2\delta} \bigr) \biggr) 
    \end{aligned}
\end{equation}
with 
\begin{equation*}
    \Lambda(t,\xi) := \frac{1}{16} \jxi^{-3} \int_1^t \frac{1}{s} \bigl| \jxi^{\frac32} \gulellsh(s,\xi) \bigr|^2 \, \ud s.
\end{equation*}
Below we show that uniformly for all $1 \leq t_1 \leq t_2 \leq T$,
\begin{equation} \label{equ:pointwise_integrated_bound_oscillating_terms}
    \biggl| \int_{t_1}^{t_2} e^{-i\Lambda(s,\xi)} e^{-i\xi\theta(s)} \calT_{osc}(s,\xi) \, \ud s \biggr| \lesssim \varepsilon^3 t_1^{-\frac32}.
\end{equation}
The asserted estimates \eqref{equ:pointwise_proposition_asserted_bound} and \eqref{equ:pointwise_proposition_asserted_difference_bound} in the statement of Proposition~\ref{prop:pointwise_profile_bounds} then follow from integrating \eqref{equ:pointwise_evol_equ_g_with_integrating_factor} in time, taking the $L^\infty_\xi$ norm, and using \eqref{equ:pointwise_weighted_B_Linfty_bound}.

It remains to prove \eqref{equ:pointwise_integrated_bound_oscillating_terms}. We provide the details for the contributions of the terms $\calI_{1,asympt}^{\delta_0}(t,\xi)$ and $\calI_{1,asympt}^{\pvdots}(t,\xi)$ in the definition \eqref{equ:pointwise_Tosc_definition} of $\calT_{osc}(t,\xi)$, and we leave the analogous details for the treatment of the other terms to the reader.
We begin with $\calI_{1,asympt}^{\delta_0}(t,\xi)$. 
Exploiting the oscillations of the phase $e^{it(-\jxi+3\jap{\frac{\xi}{3}})}$ in $\calI_{1,asympt}^{\delta_0}(t,\xi)$, we write
\begin{equation}
    \begin{aligned}
        &e^{-i\Lambda(t,\xi)} e^{-i\xi\theta(t)} \calI_{1,asympt}^{\delta_0}(t,\xi) \\
        &\quad = \pt \Bigl( I_{1,(a)}^{\delta_0}(t,\xi) \Bigr) + I_{1,(b)}^{\delta_0}(t,\xi) + I_{1,(c)}^{\delta_0}(t,\xi) + I_{1,(d)}^{\delta_0}(t,\xi) + I_{1,(e)}^{\delta_0}(t,\xi),
    \end{aligned}
\end{equation}
where 
\begin{align*}
    I_{1,(a)}^{\delta_0}(t,\xi) &:= e^{-i\Lambda(t,\xi)} e^{-i\xi\theta(t)} \frac{2\pi}{t} \frac{1}{\sqrt{3}} e^{i\frac{\pi}{2}} e^{it(-\jxi+3\jap{\frac{\xi}{3}})}  (-i) \bigl(-\jxi+3\jap{\tstyfrakxithree}\bigr)^{-1} \\ 
    &\qquad \qquad \qquad \qquad \qquad \qquad \qquad \times \frakm_{\ulell,+++}^{\delta_0}\bigl(\xi, \tstyfrakxithree, \tstyfrakxithree, \tstyfrakxithree \bigr) \jxi^{\frac32} \Bigl( \gulellsh\bigl(t, \tstyfrakxithree \bigr) \Bigr)^3, \\ 
    I_{1,(b)}^{\delta_0}(t,\xi) &:= e^{-i\Lambda(t,\xi)} e^{-i\xi\theta(t)} \frac{2\pi}{t^2} \frac{1}{\sqrt{3}} e^{i\frac{\pi}{2}} e^{it(-\jxi+3\jap{\frac{\xi}{3}})}  (-i) \bigl(-\jxi+3\jap{\tstyfrakxithree}\bigr)^{-1} \\
    &\qquad \qquad \qquad \qquad \qquad \qquad \qquad \times \frakm_{\ulell,+++}^{\delta_0}\bigl(\xi, \tstyfrakxithree, \tstyfrakxithree, \tstyfrakxithree \bigr) \jxi^{\frac32} \Bigl( \gulellsh\bigl(t, \tstyfrakxithree \bigr) \Bigr)^3 \\ 
    I_{1,(c)}^{\delta_0}(t,\xi) &:= - e^{-i\Lambda(t,\xi)} e^{-i\xi\theta(t)} \frac{2\pi}{t} \frac{1}{\sqrt{3}} e^{i\frac{\pi}{2}} e^{it(-\jxi+3\jap{\frac{\xi}{3}})}  (-i) \bigl(-\jxi+3\jap{\tstyfrakxithree}\bigr)^{-1} \\ 
    &\qquad \qquad \qquad \qquad \qquad \qquad \qquad \times \frakm_{\ulell,+++}^{\delta_0}\bigl(\xi, \tstyfrakxithree, \tstyfrakxithree, \tstyfrakxithree \bigr) \jxi^{\frac32} 3 \Bigl( \gulellsh\bigl(t, \tstyfrakxithree \bigr) \Bigr)^2 (\pt \gulellsh)\bigl(t, \tstyfrakxithree \bigr) \\ 
    I_{1,(d)}^{\delta_0}(t,\xi) &:= i \pt \Lambda(t,\xi) e^{-i\Lambda(t,\xi)} e^{-i\xi\theta(t)} \frac{2\pi}{t} \frac{1}{\sqrt{3}} e^{i\frac{\pi}{2}} e^{it(-\jxi+3\jap{\frac{\xi}{3}})}  (-i) \bigl(-\jxi+3\jap{\tstyfrakxithree}\bigr)^{-1} \\
    &\qquad \qquad \qquad \qquad \qquad \qquad \qquad \times \frakm_{\ulell,+++}^{\delta_0}\bigl(\xi, \tstyfrakxithree, \tstyfrakxithree, \tstyfrakxithree \bigr) \jxi^{\frac32} \Bigl( \gulellsh\bigl(t, \tstyfrakxithree \bigr) \Bigr)^3 \\
    I_{1,(e)}^{\delta_0}(t,\xi) &:= i \xi (\pt \theta)(t) e^{-i\Lambda(t,\xi)} e^{-i\xi\theta(t)} \frac{2\pi}{t} \frac{1}{\sqrt{3}} e^{i\frac{\pi}{2}} e^{it(-\jxi+3\jap{\frac{\xi}{3}})}  (-i) \bigl(-\jxi+3\jap{\tstyfrakxithree}\bigr)^{-1} \\ 
    &\qquad \qquad \qquad \qquad \qquad \qquad \qquad \times \frakm_{\ulell,+++}^{\delta_0}\bigl(\xi, \tstyfrakxithree, \tstyfrakxithree, \tstyfrakxithree \bigr) \jxi^{\frac32} \Bigl( \gulellsh\bigl(t, \tstyfrakxithree \bigr) \Bigr)^3.
\end{align*}
Next, we show that
\begin{equation} \label{equ:pointwise_I1a_bound}
    \sup_{1 \leq t \leq T} \, t \cdot \bigl\| I_{1,(a)}^{\delta_0}(t,\xi) \bigr\|_{L^\infty_\xi} \lesssim \varepsilon^3
\end{equation}
as well as 
\begin{equation} \label{equ:pointwise_I1bcde_bound}
    \sup_{1 \leq t \leq T} \, t^{\frac32} \Bigl( \bigl\| I_{1,(b)}^{\delta_0}(t,\xi) \bigr\|_{L^\infty_\xi} + \bigl\| I_{1,(c)}^{\delta_0}(t,\xi) \bigr\|_{L^\infty_\xi} + \bigl\| I_{1,(d)}^{\delta_0}(t,\xi) \bigr\|_{L^\infty_\xi}  + \bigl\| I_{1,(e)}^{\delta_0}(t,\xi) \bigr\|_{L^\infty_\xi} \Bigr) \lesssim \varepsilon^3.
\end{equation}
Then the bound \eqref{equ:pointwise_integrated_bound_oscillating_terms} for the contribution of $\calI_{1,asympt}^{\delta_0}(t,\xi)$ follows immediately from \eqref{equ:pointwise_I1a_bound} and \eqref{equ:pointwise_I1bcde_bound} upon integrating in time.

For the proofs of \eqref{equ:pointwise_I1a_bound} and of \eqref{equ:pointwise_I1bcde_bound}, we observe that $\bigl(-\jxi+3\jap{\tstyfrakxithree}\bigr)^{-1} \simeq \jxi$ and that by inspection from Subsection~\ref{subsec:cubic_spectral_distributions}, we have 
\begin{equation*}
    \sup_{\xi \in \bbR} \, \bigl| \frakm_{\ulell,+++}^{\delta_0}\bigl(\xi, \tstyfrakxithree, \tstyfrakxithree, \tstyfrakxithree \bigr) \bigr| \lesssim 1.
\end{equation*}
The bound \eqref{equ:pointwise_I1a_bound} follows easily using \eqref{equ:prop_profile_bounds_assumption2},
\begin{equation*}
    \bigl\| I_{1,(a)}^{\delta_0}(t,\xi) \bigr\|_{L^\infty_\xi} \lesssim t^{-1} \Bigl\| \jxi^{\frac52} \Bigl( \gulellsh\bigl(t, \tstyfrakxithree \bigr) \Bigr)^3 \Bigr\|_{L^\infty_\xi} \lesssim t^{-1} \varepsilon^3.
\end{equation*}
Next we prove the bound \eqref{equ:pointwise_I1bcde_bound}. Using \eqref{equ:prop_profile_bounds_assumption2} we infer straightaway that
\begin{equation*}
    \bigl\| I_{1,(b)}^{\delta_0}(t,\xi) \bigr\|_{L^\infty_\xi} \lesssim t^{-2} \Bigl\| \jxi^{\frac52} \Bigl( \gulellsh\bigl(t, \tstyfrakxithree \bigr) \Bigr)^3 \Bigr\|_{L^\infty_\xi} \lesssim t^{-2} \varepsilon^3.
\end{equation*}
In order to estimate the term $I_{1,(c)}^{\delta_0}(t,\xi)$, we insert the evolution equation \eqref{equ:setting_up_g_evol_equ3} for the effective profile, which gives
\begin{equation*}
    \begin{aligned}
        \bigl\| I_{1,(c)}^{\delta_0}(t,\xi) \bigr\|_{L^\infty_\xi} &\lesssim t^{-1} \Bigl\| \jxi^{\frac52} \Bigl( \gulellsh\bigl(t, \tstyfrakxithree \bigr) \Bigr)^2 \Bigr\|_{L^\infty_\xi} \Bigl( |\dot{q}(t) - \ulell| \bigl\| \tstyfrakxithree \gulellsh(t,\tstyfrakxithree) \bigr\|_{L^\infty_\xi} + \bigl\| \calFulellsh\bigl[ \calQ_{\ulell}\bigl(\usubeone(t)\bigr) \bigr](\tstyfrakxithree) \bigr\|_{L^\infty_\xi} \\ 
        &\qquad \qquad \qquad \qquad \qquad \qquad \qquad \qquad \qquad + \bigl\| \calFulellsh\bigl[ {\textstyle \frac16} \usubeone(t)^3 \bigr](\tstyfrakxithree) \bigr\|_{L^\infty_\xi} + \bigl\| \widetilde{\calR}(t,\tstyfrakxithree) \bigr\|_{L^\infty_\xi} \Bigr).
    \end{aligned}
\end{equation*}
By \eqref{equ:consequences_dispersive_decay_radiation} and \eqref{equ:consequences_energy_bound_radiation}, we can bound the second and the third term in the parentheses as follows
\begin{align*}
    \bigl\| \calFulellsh\bigl[ \calQ_{\ulell}\bigl(\usubeone(t)\bigr) \bigr](\tstyfrakxithree) \bigr\|_{L^\infty_\xi} &\lesssim \bigl\| \calQ_{\ulell}\bigl(\usubeone(t)\bigr) \bigr\|_{L^1_y} \lesssim \bigl\| \alpha(\ulg y) \bigr\|_{L^1_y} \|\usubeone(t)\|_{L^\infty_y}^2 \lesssim \varepsilon^2 t^{-1}, \\
    \bigl\| \calFulellsh\bigl[ {\textstyle \frac16} \usubeone(t)^3 \bigr](\tstyfrakxithree) \bigr\|_{L^\infty_\xi} &\lesssim \bigl\| \usubeone(t)^3 \bigr\|_{L^1_y} \lesssim \|\usubeone(t)\|_{L^2_y}^2 \|\usubeone(t)\|_{L^\infty_y} \lesssim \varepsilon^3 t^{-\frac12}.
\end{align*}
By additionally invoking \eqref{equ:prop_profile_bounds_assumption2}, \eqref{equ:consequences_qdot_minus_ulell_decay}, and \eqref{equ:pointwise_weighted_calR_Linfty_bound}, we thus obtain
\begin{equation*}
    \begin{aligned}
        \bigl\| I_{1,(c)}^{\delta_0}(t,\xi) \bigr\|_{L^\infty_\xi} &\lesssim t^{-1} \varepsilon^2 \bigl( \varepsilon t^{-1+\delta} + \varepsilon t^{-1} + \varepsilon^3 t^{-\frac12} + \varepsilon^2 t^{-\frac32+2\delta} \bigr) \lesssim \varepsilon^3 t^{-\frac32}.
    \end{aligned}
\end{equation*}
The estimates for $I_{1,(d)}^{\delta_0}(t,\xi)$ and $I_{1,(e)}^{\delta_0}(t,\xi)$ are again straightforward. Using \eqref{equ:prop_profile_bounds_assumption2} and \eqref{equ:consequences_qdot_minus_ulell_decay}, we find 
\begin{equation*}
    \bigl\| I_{1,(d)}^{\delta_0}(t,\xi) \bigr\|_{L^\infty_\xi} \lesssim t^{-1} \bigl\| \gulellsh(t,\xi) \bigr\|_{L^\infty_\xi}^2 \cdot t^{-1} \Bigl\| \jxi^{\frac52} \Bigl( \gulellsh\bigl(t, \tstyfrakxithree \bigr) \Bigr)^3 \Bigr\|_{L^\infty_\xi} \lesssim t^{-2} \varepsilon^5,
\end{equation*}
and 
\begin{equation*}
    \bigl\| I_{1,(e)}^{\delta_0}(t,\xi) \bigr\|_{L^\infty_\xi} \lesssim |\dot{q}(t)-\ulell| \bigl\| \jxi \gulellsh(t,\xi) \bigr\|_{L^\infty_\xi} \cdot t^{-1} \Bigl\| \jxi^{\frac52} \Bigl( \gulellsh\bigl(t, \tstyfrakxithree \bigr) \Bigr)^3 \Bigr\|_{L^\infty_\xi} \lesssim t^{-2+\delta} \varepsilon^4.
\end{equation*}
Combining the preceding estimates finishes the proof of the bound \eqref{equ:pointwise_I1bcde_bound}, and thus concludes the estimates for the contribution of the term $\calI_{1,asympt}^{\delta_0}(t,\xi)$ to the bound \eqref{equ:pointwise_integrated_bound_oscillating_terms}. 

The treatment of the contributions of the term $\calI_{1,asympt}^{\pvdots}(t,\xi)$ defined in \eqref{equ:pointwise_calIasympt_pv_definitions} is mostly identical to the preceding discussion of the term $\calI_{1,asympt}^{\delta_0}(t,\xi)$, and proceeds by integration by parts in time exploiting the oscillations of the same phase function $e^{it(-\jxi+3\jap{\frac{\xi}{3}})}$. The only significant difference is that the term $\calI_{1,asympt}^{\pvdots}(t,\xi)$ additionally features the factors $\tanh\bigl( \ulg (\partial_\xi \phi_1)(\xi) t \bigr) \chi_1\bigl( (\tstyfrakxithree + \ulg\ulell) t^{5\delta} \bigr)$. When the time derivative falls onto the cut-off $\chi_1\bigl( (\tstyfrakxithree + \ulg\ulell) t^{5\delta} \bigr)$, this obviously produces additional time decay. Instead when the time derivative falls onto $\tanh\bigl( \ulg (\partial_\xi \phi_1)(\xi) t \bigr)$, it produces the term
\begin{equation*}
    \ulg (\pxi \phi_1)(\xi) \sech^2\bigl( \ulg (\partial_\xi \phi_1)(\xi) t \bigr).
\end{equation*}
Thanks to the cut-off $\chi_1\bigl( (\tstyfrakxithree + \ulg\ulell) t^{5\delta} \bigr)$, we have $|(\tstyfrakxithree + \ulg\ulell)| \gtrsim t^{-5\delta}$, whence $|(\pxi \phi_1)(\xi)| \gtrsim t^{-5\delta}$ and the rapid decay of the hyperbolic secant function therefore furnishes more than enough additional time decay.

This finishes the proof of Proposition~\ref{prop:pointwise_profile_bounds}.
\end{proof}

\appendix

\section{Fourier transforms of hyperbolic functions}

In this appendix we collect several Fourier transform identities for hyperbolic functions that are used in the computation of the nonlinear spectral distributions in Section~\ref{sec:nonlinear_spectral_distributions} and in Lemma~\ref{lem:null_structure1}. 
We recall the well-known Fourier transforms of $\sech(x)$ and $\sech^2(x)$ given by
\begin{align}
      \int_\bbR e^{ix\zeta} \sech(x) \, \ud x &= \pi \cosech\Bigl( \frac{\pi}{2} \zeta \Bigr), \label{equ:appendix_FT_sech} \\ 
      \int_\bbR e^{ix\zeta} \sech^2(x) \, \ud x &= \pi \zeta \cosech\Bigl( \frac{\pi}{2} \zeta \Bigr).\label{equ:appendix_FT_sech2}
\end{align}
See \cite[Example 3.3]{SteinShakarchi_Complex} for a proof of \eqref{equ:appendix_FT_sech} and \cite[Corollary 5.7]{LS1} for a proof of \eqref{equ:appendix_FT_sech2}.
We refer to \cite[Appendix A]{LL2} for a systematic derivation of the following identities from \eqref{equ:appendix_FT_sech} and \eqref{equ:appendix_FT_sech2},
\begin{align}
  \int_\bbR e^{ix\zeta} \sech(x) \tanh(x) \, \ud x &= i \pi \zeta \sech\Bigl(\frac{\pi}{2} \zeta\Bigr), \label{equ:appendix_FT_sech_tanh1} \\
  \int_\bbR e^{ix\zeta} \sech(x) \tanh^2(x) \, \ud x &= \frac{\pi}{2} \bigl(1-\zeta^2\bigr) \sech\Bigl(\frac{\pi}{2} \zeta\Bigr), \label{equ:appendix_FT_sech_tanh2} \\
  \int_\bbR e^{ix\zeta} \sech(x) \tanh^3(x) \, \ud x &= -i \frac{\pi}{6} \zeta \bigl(-5+\zeta^2\bigr) \sech\Bigl(\frac{\pi}{2} \zeta\Bigr), \label{equ:appendix_FT_sech_tanh3} \\
  \int_\bbR e^{ix\zeta} \sech(x) \tanh^4(x) \, \ud x &= \frac{\pi}{24} \bigl(9-14\zeta^2+\zeta^4\bigr) \sech\Bigl(\frac{\pi}{2} \zeta\Bigr), \label{equ:appendix_FT_sech_tanh4} \\
    \int_\bbR e^{ix\zeta} \sech^2(x) \tanh(x) \, \ud x &= \frac{i}{4} \zeta^2 \cosech\Bigl( \frac{\pi}{2} \zeta \Bigr), \label{equ:appendix_FT_sech2_tanh1} \\
    \int_\bbR e^{ix\zeta} \sech^2(x) \tanh^2(x) \, \ud x &= -\frac{1}{12} \zeta (-2+\zeta^2) \cosech\Bigl( \frac{\pi}{2} \zeta \Bigr), \label{equ:appendix_FT_sech2_tanh2} \\
    \int_\bbR e^{ix\zeta} \sech^2(x) \tanh^3(x) \, \ud x &= -\frac{i}{48} \zeta^2 (-8 + \zeta^2) \cosech\Bigl( \frac{\pi}{2} \zeta \Bigr). \label{equ:appendix_FT_sech2_tanh3}
\end{align}

\section{Deeper quadratic null structures for the radiation term} \label{appendix:deeper_null_structure}

In the course of this work, we discovered deeper null structures in all quadratic nonlinearities of the evolution equation for the radiation term. Although the proof of Theorem~\ref{thm:main} relies solely on the null structure revealed in Lemma~\ref{lem:null_structure1}, we record these findings in this appendix for the (algebraically simpler) case of odd perturbations of the static sine-Gordon kink \eqref{eq:staticK}.
Decomposing the corresponding solution $\phi(t,x)$ to the sine-Gordon equation \eqref{equ:intro_sG} as
\begin{equation}
    \phi(t,x) = K(x) + u(t,x),
\end{equation}
we compute that the radiation term $u(t,x)$ must solve the nonlinear Klein-Gordon equation
\begin{equation}\label{eq:oddperbeq}
 \bigl( \pt^2 - \px^2 - 2 \sech^2(x) + 1 \bigr) u = \alpha(x) u^2 + \beta_0 u^3 + \beta(x) u^3 + \bigl\{\text{higher order}\bigr\}
\end{equation}
with
\begin{equation*}
 \begin{aligned}
    \alpha(x) := - \sech(x) \tanh(x), \quad \beta_0 := \frac16, \quad \beta(x) := -\frac13 \sech^2(x).
 \end{aligned}
\end{equation*}
Setting 
\begin{equation}
 v(t) := \frac12 \bigl( u(t) - i \jwtilD^{-1} \pt u(t) \bigr), \quad \jap{\wtilD} := \sqrt{1 -\px^2 -2\sech^2(x)},
\end{equation}
we can pass from \eqref{eq:oddperbeq} to the following (complexified) first-order equation
\begin{equation} \label{equ:appendix_odd_perturbation_1sorder}
 \bigl( \pt - i \jwtilD \bigr) v = (2i\jwtilD)^{-1} \Bigl( \alpha(x) \bigl(v+\bar{v}\bigr)^2 + \beta_0 \bigl(v+\bar{v}\bigr)^3 + \beta(x) \bigl(v+\bar{v}\bigr)^3  + \{\text{higher order}\} \Bigr).
\end{equation}
Taking a space-time resonances approach to determine the long-time behavior of the solution $v(t,x)$ to \eqref{equ:appendix_odd_perturbation_1sorder}, 
we consider its profile $f(t) := e^{-it\jap{\widetilde{D}}} v(t)$ and study the evolution equation of the (modified) distorted Fourier transform of the profile
\begin{equation*}
    \tilde{f}(t,\xi) := \widetilde{\calF}^{\#}\bigl[ f(t) \bigr](\xi) = \int_\bbR \overline{e^{\#}(x,\xi)} f(t,x) \, \ud x,
\end{equation*}
where 
\begin{equation} \label{equ:definition_dist_FT_basis}
 e^{\#}(x,\xi) := \frac{1}{\sqrt{2\pi}} \frac{\xi + i \tanh(x)}{|\xi|-i} e^{i x \xi}.
\end{equation}
Here, $\widetilde{\calF}^{\#}$ denotes the (modified) distorted Fourier transform associated with the (selfadjoint) scalar Schr\"odinger operator $-\partial_x^2 - 2 \sech^2(x)$ from \eqref{eq:oddperbeq}.
Expressing the solution $v(t,x)$ in terms of the (modified) distorted Fourier transform of its profile
\begin{equation}
    v(t,x) = \widetilde{\calF}^{\#,\ast}\bigl[ e^{it\jxi} \tilde{f}(t,\xi) \bigr](x) = \int_\bbR e^{\#}(x,\xi) e^{it\jxi} \tilde{f}(t,\xi) \, \ud \xi,
\end{equation}
we find that the contribution of the spatially localized quadratic nonlinearity in \eqref{equ:appendix_odd_perturbation_1sorder} to the evolution equation for the (modified) distorted Fourier transform of the profile is given by
 \begin{align}
    &e^{-it\jxi} \wtilcalF^{\#}\Bigl[ \alpha(x) \bigl( v(t,x) + \barv(t,x) \bigr)^2 \Bigr](\xi) \\
    &= e^{-it\jxi} \wtilcalF^{\#}\Bigl[ \alpha(x) v(t,x)^2 \Bigr](\xi) 
    + 2 e^{-it\jxi} \wtilcalF^{\#}\Bigl[ \alpha(x) v(t,x) \barv(t,x) \Bigr](\xi) + e^{-it\jxi} \wtilcalF^{\#}\Bigl[ \alpha(x) \barv(t,x)^2 \Bigr](\xi) \\
    &= \iint_{\bbR^2} e^{it(-\jxi+\jeta+\jsigma)} \tilf(t,\eta) \tilf(t,\sigma) \frakq_{++}(\xi,\eta,\sigma) \, \ud \eta \, \ud \sigma \label{eq:deepnullphase1}\\
    &\quad + 2 \iint_{\bbR^2} e^{it(-\jxi-\jeta+\jsigma)} \overline{\tilf}(t,\eta) \tilf(t,\sigma) \frakq_{-+}(\xi,\eta,\sigma) \, \ud \eta \, \ud \sigma  \label{eq:deepnullphase2}\\
    &\quad + \iint_{\bbR^2} e^{it(-\jxi-\jeta-\jsigma)} \overline{\tilf}(t,\eta) \overline{\tilf}(t,\sigma) \frakq_{--}(\xi,\eta,\sigma) \, \ud \eta \, \ud \sigma \label{eq:deepnullphase3}
 \end{align}
with the quadratic spectral distributions 
\begin{equation*}
 \begin{aligned}
  \frakq_{++}(\xi,\eta,\sigma) &:= \int_\bbR \alpha(x) \overline{e^{\#}(x,\xi)} \, e^{\#}(x,\eta) \, e^{\#}(x,\sigma) \, \ud x, \\
  \frakq_{-+}(\xi,\eta,\sigma) &:= \int_\bbR \alpha(x) \overline{e^{\#}(x,\xi)} \, \overline{e^{\#}(x,\eta)} \, e^{\#}(x,\sigma) \, \ud x, \\
  \frakq_{--}(\xi,\eta,\sigma) &:= \int_\bbR \alpha(x) \overline{e^{\#}(x,\xi)} \, \overline{e^{\#}(x,\eta)} \, \overline{e^{\#}(x,\sigma)} \, \ud x.
 \end{aligned}
\end{equation*}
In the next lemma we compute these quadratic spectral distributions explicitly.

\begin{lemma} \label{lem:deepnull}
It holds that 
    \begin{align}
\frakq_{++}(\xi,\eta,\sigma) &= \bigl(-\jxi+\jeta+\jsigma\bigr) \, \frakp_{++}(\xi,\eta,\sigma), \label{eq:deepnull1}\\
        \frakq_{-+}(\xi,\eta,\sigma) &= \bigl(-\jxi-\jeta+\jsigma\bigr) \, \frakp_{-+}(\xi,\eta,\sigma), \label{eq:deepnull2}\\
        \frakq_{--}(\xi,\eta,\sigma) &= \bigl(-\jxi-\jeta-\jsigma\bigr) \, \frakp_{--}(\xi,\eta,\sigma)\label{eq:deepnull3},
\end{align}
with 
\begin{equation}
    \begin{aligned}
    \frakp_{++}(\xi,\eta,\sigma)&=i \frac{\pi}{8}\bigl( \jxi + \jeta + \jsigma \bigr) \bigl( \jxi^2 - (\jeta-\jsigma)^2 \bigr) \sech\Bigl(\frac{\pi}{2} (\xi-\eta-\sigma) \Bigr), \\
        \frakq_{-+}(\xi,\eta,\sigma) &=  i \frac{\pi}{8} \bigl(\jxi+\jeta+\jsigma\bigr) \Bigl( - (\jxi-\jeta)^2 + \jsigma^2 \Bigr) \sech\Bigl(\frac{\pi}{2} (\xi+\eta-\sigma) \Bigr), \\
        \frakq_{--}(\xi,\eta,\sigma) &=-i \frac{\pi}{8}\bigl( -\jxi + \jeta + \jsigma \bigr) \bigl( \jxi^2 - (\jeta-\jsigma)^2 \bigr) \sech\Bigl(\frac{\pi}{2} (\xi+\eta+\sigma) \Bigr).
    \end{aligned}
    \end{equation}
\end{lemma}

We observe that the expressions inside the parentheses in \eqref{eq:deepnull1}, \eqref{eq:deepnull2}, \eqref{eq:deepnull3} are exactly the phase functions in the corresponding integrals \eqref{eq:deepnullphase1}, \eqref{eq:deepnullphase2}, \eqref{eq:deepnullphase3}.
Thus, via suitable normal forms one can in fact completely transform the quadratic nonlinearities on the right-hand side of the evolution equation \eqref{equ:appendix_odd_perturbation_1sorder} for the radiation term into spatially localized cubic nonlinearities.

In the proof of \cite[Theorem 1.1]{LS1} or in the proof of Theorem~\ref{thm:main} (for $\ell_0 = 0$), only the vanishing of $\frakq_{++}(\xi,\eta,\sigma)$ at $(\xi,\eta,\sigma) = (\pm \sqrt{3}, 0, 0)$ is crucially exploited, while for all other frequency interactions it suffices to exploit the improved local decay of the radiation term at higher frequencies and the fact that the phase functions in the integrals \eqref{eq:deepnull2}, \eqref{eq:deepnull3} do not have time resonances at $\eta = \sigma = 0$ for arbitrary $\xi \in \bbR$.

\begin{proof}[Proof of Lemma~\ref{lem:deepnull}]
We start with the derivation of \eqref{eq:deepnull1}. 
Inserting \eqref{equ:definition_dist_FT_basis} into $\frakq_{++}(\xi,\eta,\sigma)$, we obtain
\begin{equation*}
 \begin{aligned}
  \frakq_{++}(\xi,\eta,\sigma) &= - \frac{1}{(2\pi)^{\frac32}} \frac{1}{|\xi|+i} \frac{1}{|\eta|-i} \frac{1}{|\sigma|-i} \\
  &\quad \quad \times \int_{\bbR} e^{-ix(\xi-\eta-\sigma)} \sech(x) \tanh(x) \bigl(\xi-i\tanh(x)\bigr) \bigl(\eta+i\tanh(x)\bigr) \bigl(\sigma+i\tanh(x)\bigr) \, \ud x.
 \end{aligned}
\end{equation*}
Expanding the integrand in the preceding line gives
\begin{equation*}
 \begin{aligned}
  &\int_{\bbR} e^{-ix(\xi-\eta-\sigma)} \sech(x) \tanh(x) \bigl(\xi-i\tanh(x)\bigr) \bigl(\eta+i\tanh(x)\bigr) \bigl(\sigma+i\tanh(x)\bigr) \, \ud x \\
  &= \xi \eta \sigma \int_\bbR e^{-ix(\xi-\eta-\sigma)} \sech(x) \tanh(x) \, \ud x \\
  &\quad + i \bigl( \xi \eta + \xi \sigma - \eta \sigma \bigr) \int_\bbR e^{-ix(\xi-\eta-\sigma)} \sech(x) \tanh^2(x) \, \ud x \\
  &\quad + \bigl(-\xi+\eta+\sigma\bigr) \int_\bbR e^{-ix(\xi-\eta-\sigma)} \sech(x) \tanh^3(x) \, \ud x \\
  &\quad + i \int_\bbR e^{-ix(\xi-\eta-\sigma)} \sech(x) \tanh^4(x) \, \ud x.
 \end{aligned}
\end{equation*}
Using the Fourier transform identities \eqref{equ:appendix_FT_sech_tanh1}, \eqref{equ:appendix_FT_sech_tanh2}, \eqref{equ:appendix_FT_sech_tanh3}, \eqref{equ:appendix_FT_sech_tanh4},
we find
\begin{equation*}
 \begin{aligned}
  &\int_{\bbR} e^{-ix(\xi-\eta-\sigma)} \sech(x) \tanh(x) \bigl(\xi-i\tanh(x)\bigr) \bigl(\eta+i\tanh(x)\bigr) \bigl(\sigma+i\tanh(x)\bigr) \, \ud x \\
  &= \xi \eta \sigma (-i) \pi (\xi-\eta-\sigma) \sech\Bigl(\frac{\pi}{2} (\xi-\eta-\sigma) \Bigr) \\
  &\quad + i \bigl( \xi \eta + \xi \sigma - \eta \sigma \bigr) \frac{\pi}{2} \bigl( 1 - (\xi-\eta-\sigma)^2 \bigr) \sech\Bigl(\frac{\pi}{2} (\xi-\eta-\sigma) \Bigr) \\
  &\quad + (-\xi+\eta+\sigma) \frac{\pi}{6} i (\xi-\eta-\sigma) \bigl( -5 + (\xi-\eta-\sigma)^2 \bigr) \sech\Bigl(\frac{\pi}{2} (\xi-\eta-\sigma) \Bigr) \\
  &\quad + i \frac{\pi}{24} \bigl( 9 - 14 (\xi-\eta-\sigma)^2 + (\xi-\eta-\sigma)^4 \bigr) \sech\Bigl(\frac{\pi}{2} (\xi-\eta-\sigma) \Bigr) \\
  &= i \frac{\pi}{8} \Bigl( 3 - \xi^4 - \eta^4 + 2\sigma^2 - \sigma^4 + 2 \eta^2 (1+\sigma^2) + 2 \xi^2 (1+\eta^2+\sigma^2) \Bigr) \sech\Bigl(\frac{\pi}{2} (\xi-\eta-\sigma) \Bigr) \\
  &= i \frac{\pi}{8}\bigl( -\jxi + \jeta + \jsigma \bigr) \bigl( \jxi + \jeta + \jsigma \bigr) \bigl( \jxi^2 - (\jeta-\jsigma)^2 \bigr) \sech\Bigl(\frac{\pi}{2} (\xi-\eta-\sigma) \Bigr)\\
  &=: \bigl(-\jxi+\jeta+\jsigma\bigr) \, \frakp_{++}(\xi,\eta,\sigma).
 \end{aligned}
\end{equation*}

Next, we prove \eqref{eq:deepnull2}. Inserting \eqref{equ:definition_dist_FT_basis} into $\frakq_{-+}(\xi,\eta,\sigma)$, we obtain
\begin{equation*}
 \begin{aligned}
  \frakq_{-+}(\xi,\eta,\sigma) &= - \frac{1}{(2\pi)^{\frac32}} \frac{1}{|\xi|+i} \frac{1}{|\eta|+i} \frac{1}{|\sigma|-i} \\
  &\quad \quad \times \int_{\bbR} e^{-ix(\xi+\eta-\sigma)} \sech(x) \tanh(x) \bigl(\xi-i\tanh(x)\bigr) \bigl(\eta-i\tanh(x)\bigr) \bigl(\sigma+i\tanh(x)\bigr) \, \ud x.
 \end{aligned}
\end{equation*}
Expanding the integrand in the preceding line gives
\begin{equation*}
 \begin{aligned}
  &\int_{\bbR} e^{-ix(\xi+\eta-\sigma)} \sech(x) \tanh(x) \bigl(\xi-i\tanh(x)\bigr) \bigl(\eta-i\tanh(x)\bigr) \bigl(\sigma+i\tanh(x)\bigr) \, \ud x \\
  &= \xi \eta \sigma \int_{\bbR} e^{-ix(\xi+\eta-\sigma)} \sech(x) \tanh(x) \, \ud x \\
  &\quad + i ( \xi \eta - \xi \sigma - \eta \sigma ) \int_\bbR  e^{-ix(\xi+\eta-\sigma)} \sech(x) \tanh^2(x) \, \ud x \\
  &\quad + (\xi+\eta-\sigma) \int_\bbR e^{-ix(\xi+\eta-\sigma)} \sech(x) \tanh^3(x) \, \ud x \\
  &\quad - i \int_\bbR e^{-ix(\xi+\eta-\sigma)} \sech(x) \tanh^4(x) \, \ud x.
 \end{aligned}
\end{equation*}
Using again the Fourier transform identities  \eqref{equ:appendix_FT_sech_tanh1}, \eqref{equ:appendix_FT_sech_tanh2}, \eqref{equ:appendix_FT_sech_tanh3}, \eqref{equ:appendix_FT_sech_tanh4}, we find 
\begin{equation*}
 \begin{aligned}
  &\int_{\bbR} e^{-ix(\xi+\eta-\sigma)} \sech(x) \tanh(x) \bigl(\xi-i\tanh(x)\bigr) \bigl(\eta-i\tanh(x)\bigr) \bigl(\sigma+i\tanh(x)\bigr) \, \ud x \\
  &= \xi \eta \sigma (-i) \pi (\xi+\eta-\sigma) \sech\Bigl(\frac{\pi}{2} (\xi+\eta-\sigma) \Bigr) \\
  &\quad + i \bigl( \xi \eta - \xi \sigma - \eta \sigma \bigr) \frac{\pi}{2} \bigl( 1 - (\xi+\eta-\sigma)^2 \bigr) \sech\Bigl(\frac{\pi}{2} (\xi+\eta-\sigma) \Bigr) \\
  &\quad + (\xi+\eta-\sigma) \frac{\pi}{6} i (\xi+\eta-\sigma) \bigl( -5 + (\xi+\eta-\sigma)^2 \bigr) \sech\Bigl(\frac{\pi}{2} (\xi+\eta-\sigma) \Bigr) \\
  &\quad - i \frac{\pi}{24} \bigl( 9 - 14 (\xi+\eta-\sigma)^2 + (\xi+\eta-\sigma)^4 \bigr) \sech\Bigl(\frac{\pi}{2} (\xi+\eta-\sigma) \Bigr) \\
  &= i \frac{\pi}{8} \Bigl( -3 -2\xi^2 + \xi^4 - 2 \eta^2 - 2 \xi^2 \eta^2 + \eta^4 - 2 \sigma^2 - 2 \xi^2 \sigma^2 - 2 \eta^2 \sigma^2 + \sigma^4 \Bigr) \sech\Bigl(\frac{\pi}{2} (\xi+\eta-\sigma) \Bigr) \\
  &= i \frac{\pi}{8} \bigl(-\jxi-\jeta+\jsigma\bigr)\bigl(\jxi+\jeta+\jsigma\bigr) \Bigl( - (\jxi-\jeta)^2 + \jsigma^2 \Bigr) \sech\Bigl(\frac{\pi}{2} (\xi+\eta-\sigma) \Bigr)\\
  &=:\bigl(-\jxi-\jeta+\jsigma\bigr) \, \frakp_{-+}(\xi,\eta,\sigma).
 \end{aligned}
\end{equation*}

Lastly, we establish \eqref{eq:deepnull3}. Inserting \eqref{equ:definition_dist_FT_basis} into $\frakq_{--}(\xi,\eta,\sigma)$ gives
\begin{equation*}
 \begin{aligned}
  \frakq_{--}(\xi,\eta,\sigma) &= - \frac{1}{(2\pi)^{\frac32}} \frac{1}{|\xi|+i} \frac{1}{|\eta|+i} \frac{1}{|\sigma|+i} \\
  &\quad \quad \times \int_{\bbR} e^{-ix(\xi+\eta+\sigma)} \sech(x) \tanh(x) \bigl(\xi-i\tanh(x)\bigr) \bigl(\eta-i\tanh(x)\bigr) \bigl(\sigma-i\tanh(x)\bigr) \, \ud x.
 \end{aligned}
\end{equation*}
Expanding the integrand in the preceding line yields
\begin{equation*}
 \begin{aligned}
  &\int_{\bbR} e^{-ix(\xi+\eta+\sigma)} \sech(x) \tanh(x) \bigl(\xi-i\tanh(x)\bigr) \bigl(\eta-i\tanh(x)\bigr) \bigl(\sigma-i\tanh(x)\bigr) \, \ud x \\
  &= \xi \eta \sigma \int_\bbR e^{-ix(\xi+\eta+\sigma)} \sech(x) \tanh(x) \, \ud x \\
  &\quad - i \bigl( \xi \eta + \xi \sigma + \eta \sigma \bigr) \int_\bbR e^{-ix(\xi+\eta+\sigma)} \sech(x) \tanh^2(x) \, \ud x \\
  &\quad - \bigl(\xi+\eta+\sigma\bigr) \int_\bbR e^{-ix(\xi+\eta+\sigma)} \sech(x) \tanh^3(x) \, \ud x \\
  &\quad + i \int_\bbR e^{-ix(\xi+\eta+\sigma)} \sech(x) \tanh^4(x) \, \ud x.
 \end{aligned}
\end{equation*}
Using the Fourier transform identities \eqref{equ:appendix_FT_sech_tanh1}, \eqref{equ:appendix_FT_sech_tanh2}, \eqref{equ:appendix_FT_sech_tanh3}, \eqref{equ:appendix_FT_sech_tanh4},
we find
\begin{equation*}
 \begin{aligned}
  &\int_{\bbR} e^{-ix(\xi+\eta+\sigma)} \sech(x) \tanh(x) \bigl(\xi-i\tanh(x)\bigr) \bigl(\eta-i\tanh(x)\bigr) \bigl(\sigma-i\tanh(x)\bigr) \, \ud x \\
  &= \xi \eta \sigma (-i) \pi (\xi+\eta+\sigma) \sech\Bigl(\frac{\pi}{2} (\xi+\eta+\sigma) \Bigr) \\
  &\quad -i \bigl( \xi \eta + \xi \sigma + \eta \sigma \bigr) \frac{\pi}{2} \bigl( 1 - (\xi+\eta+\sigma)^2 \bigr) \sech\Bigl(\frac{\pi}{2} (\xi+\eta+\sigma) \Bigr) \\
  &\quad - (\xi+\eta+\sigma) \frac{\pi}{6} i (\xi+\eta+\sigma) \bigl( -5 + (\xi+\eta+\sigma)^2 \bigr) \sech\Bigl(\frac{\pi}{2} (\xi+\eta+\sigma) \Bigr) \\
  &\quad + i \frac{\pi}{24} \bigl( 9 - 14 (\xi+\eta+\sigma)^2 + (\xi+\eta+\sigma)^4 \bigr) \sech\Bigl(\frac{\pi}{2} (\xi+\eta+\sigma) \Bigr) \\
  &= i \frac{\pi}{8}\bigl( \jxi + \jeta + \jsigma \bigr) \bigl( -\jxi + \jeta + \jsigma \bigr) \bigl( \jxi^2 - (\jeta-\jsigma)^2 \bigr) \sech\Bigl(\frac{\pi}{2} (\xi+\eta+\sigma) \Bigr)\\
  &=: \bigl(-\jxi-\jeta-\jsigma\bigr) \, \frakp_{--}(\xi,\eta,\sigma).
 \end{aligned}
\end{equation*}
This completes the proof of Lemma~\ref{lem:deepnull}.
\end{proof}

\bibliographystyle{amsplain}
\bibliography{references}

\end{document}